\documentclass[11pt]{amsart}

\newif\ifpreprintcolor
\preprintcolortrue

\usepackage[T1]{fontenc}
\usepackage{lmodern}
\usepackage{amssymb}
\usepackage{mathtools}
\usepackage{microtype}
\usepackage{tikz}
\usetikzlibrary{arrows.meta,positioning,fit,calc}
\usepackage{xcolor}
\usepackage{hyperref}

\newtheorem{theorem}{Theorem}[subsection]
\newtheorem{lemma}[theorem]{Lemma}
\newtheorem{proposition}[theorem]{Proposition}
\newtheorem{corollary}[theorem]{Corollary}
\newtheorem{introtheorem}{Theorem}[section]
\theoremstyle{definition}

\theoremstyle{remark}
\newtheorem{remark}[theorem]{Remark}
\numberwithin{equation}{subsection}

\newcommand{\Z}{\mathbb{Z}}
\newcommand{\R}{\mathbb{R}}
\newcommand{\E}{\mathbb{E}}
\newcommand{\cP}{\mathcal{P}}
\newcommand{\Psep}{P_{\mathrm{sep}}}

\DeclarePairedDelimiter{\abs}{\lvert}{\rvert}
\DeclarePairedDelimiter{\floor}{\lfloor}{\rfloor}
\DeclarePairedDelimiter{\ceil}{\lceil}{\rceil}

\newcommand{\PaperAuthor}{Ho Dang Dung}
\newcommand{\PaperAffiliation}{HSE University}
\newcommand{\PaperEmail}{dkho@edu.hse.ru}
\title{Prime Labelings of Trees}
\author{\PaperAuthor}
\address{\PaperAffiliation}
\email{\PaperEmail}
\date{}

\makeatletter
\renewcommand{\@setauthors}{%
  \begin{center}
    \normalfont\scshape\PaperAuthor\par
    \vspace{.6em}
    \normalfont\small\PaperAffiliation\par
    \href{mailto:\PaperEmail}{\texttt{\PaperEmail}}\par
    \vspace{.8em}
    \normalfont\small 30 September 2026\par
  \end{center}%
}
\let\@setaddresses\relax
\makeatother
\makeatletter
\renewcommand\subsection{\@startsection{subsection}{2}{\z@}%
  {.9\linespacing\@plus.5\linespacing}{.4\linespacing}%
  {\normalfont\bfseries}}
\makeatother
\hypersetup{pdftitle={Prime Labelings of Trees},pdfauthor={Ho Dang Dung}}
\definecolor{linkblue}{RGB}{0,80,180}
\ifpreprintcolor
  \hypersetup{colorlinks=true,linkcolor=linkblue,citecolor=linkblue,
    urlcolor=linkblue,filecolor=linkblue,menucolor=linkblue,
    runcolor=linkblue,anchorcolor=linkblue}
\else
  \hypersetup{hidelinks}
\fi
\makeatletter
\def\indentlabel#1{\setbox\z@\hbox{#1}%
  \ifdim\wd\z@=\z@ \ifnum\@tocthislevel=\@ne \else
    \hb@xt@\@tempdima{\hfil}\fi
  \else\ifdim\wd\z@<\@tempdima\hb@xt@\@tempdima{\unhbox\z@\hfil}%
  \else\box\z@\fi\fi}
\def\@tocthislevel{0}
\def\@tocline#1#2#3#4#5#6#7{\relax
  \ifnum #1>\c@tocdepth
  \else
    \par \addpenalty\@secpenalty\addvspace{#2}%
    \begingroup \hyphenpenalty\@M \def\@tocthislevel{#1}%
    \@ifempty{#4}{%
      \@tempdima\csname r@tocindent\number#1\endcsname\relax
    }{%
      \@tempdima#4\relax
    }%
    \parindent\z@ \leftskip#3\relax \advance\leftskip\@tempdima\relax
    \rightskip\@pnumwidth plus4em \parfillskip-\@pnumwidth
    #5\leavevmode\hskip-\@tempdima #6\nobreak\relax
    \ifnum #1>\@ne
      \leaders\hbox{$\m@th\mkern4.5mu\hbox{.}\mkern4.5mu$}\hfil
    \else \hfil \fi
    \nobreak\hbox to\@pnumwidth{\@tocpagenum{#7}}\par
    \nobreak
    \endgroup
  \fi}
\def\@tocpagenum#1{\hss #1}
\def\l@section{\@tocline{1}{2.5pt plus1pt}{0pt}{2.3em}{\bfseries}}
\def\l@subsection{\@tocline{2}{0pt}{2.3em}{3.2em}{}}
\makeatother

\begin{document}

\begin{abstract}
A prime labeling of a graph labels its vertices bijectively with the
integers from one to the number of vertices, so that adjacent vertices
receive coprime labels. We prove the Entringer--Tout conjecture: every
tree has a prime labeling. Haxell, Pikhurko, and Taraz proved this for
all sufficiently large trees; our proof covers every order. The proof
is computer-assisted. Trees with at most 1359 vertices are handled by
finite certificates whose rules are proved sound. For larger trees,
primes above half the order label a small separator, the even labels
are assigned to one side of the remaining forest with prescribed
numbers in each divisibility class, and the odd labels are placed by a
matching obtained from Hall's theorem. For these trees, the
inequalities needed for this matching are verified by exact rational
computation when the natural logarithm of the order is at most 256, and
by an analytic argument otherwise.
Most of the proof development, computational verification, and drafting
of this manuscript were carried out by generative AI systems,
principally GPT Astra and Claude Opus, under the author's direction; see
the section on use of generative AI for details.
\end{abstract}

\maketitle
\begingroup
\small
\linespread{0.84}\selectfont
\tableofcontents
\endgroup
\clearpage

% Editorial dependency map: manuscript_plan.md. Writing conventions: ../style.md.
\section{Introduction}\label{sec:introduction}

For a positive integer $n$, write $[n]=\{1,\ldots,n\}$. A \emph{prime
labeling} of a graph $G$ of order $n$ is a bijection $f\colon V(G)\to[n]$
such that
\[
    \gcd(f(u),f(v))=1 \qquad\text{whenever }uv\in E(G).
\]
A graph with a prime labeling is called \emph{prime}. Entringer
conjectured around 1980 that every tree is prime. Haxell, Pikhurko, and
Taraz~\cite[Section~1]{HPT11} trace the conjecture to that time, find its
earliest written statement in a 1982 paper of Tout, Dabboucy, and
Howalla~\cite{TDH82}, and call it the Entringer--Tout conjecture. They proved it for
all sufficiently large trees~\cite[Theorem~1]{HPT11}. We prove it for
every order.

\begin{introtheorem}\label{thm:main}
For every integer $n\geq1$ and every tree $T$ on $n$ vertices, there is a
bijection $f\colon V(T)\to[n]$ such that $\gcd(f(u),f(v))=1$ for every
edge $uv$ of $T$.
\end{introtheorem}

The proof is computer-assisted. The structural reductions, the matching
criteria, and the soundness of every certificate check are proved in the
text. The finite certificate conditions are verified by exact integer
and rational computation. Section~\ref{sec:part-completion} completes the
proof of Theorem~\ref{thm:main}, and Section~\ref{sec:part-computation}
lists the computed conditions, the programs that verify them, and the
commands that reproduce the verification.

\subsection{Main ideas and contributions}\label{sec:main-ideas}

The theorem of Haxell, Pikhurko, and Taraz is stated for all
sufficiently large trees, without an explicit threshold.
They report an estimate of Spiess, based on an earlier version of their
argument, that the threshold $10^{10^{100}}$ is a suitable choice in
their Theorem~1, and they note that closing the remaining gap is beyond
any small-order approach~\cite[Section~1]{HPT11}. The published finite
verifications reach orders $15$, $34$, and $50$
\cite{FH94,Pik02,Pik07}. Direct order-by-order verification is not a
viable way to close the gap between these ranges. The intermediate
orders are therefore treated here by arguments that are uniform over
long ranges of orders and whose constants are explicit enough to be
checked.

We keep the special labels of Haxell, Pikhurko, and Taraz: primes
above $n/2$ label a small separator, and the labels $2p$ with
$n/3<p\leq n/2$ are coprime to all odd labels except $p$, so a final
exchange removes the few pairs $p,2p$
(Section~\ref{sec:repair}). The quantitative work is organized
differently, in four ways.

\begin{enumerate}
\item \emph{Explicit preparation.} Under explicit counting conditions
  on the order, a tree is reduced to a bipartite residual forest, a set
  of even labels for one side, a set of odd labels for the other, and a
  reserve of positions that accept every remaining odd label
  (Proposition~\ref{prop:residual-problem}). Special even labels are
  placed on vertices of largest degree, which bounds the residual
  degree. The prime counts in these conditions are computed exactly for
  the orders through $1982759$ and bounded by Dusart's explicit
  estimates~\cite{Dus10} for $\log n\geq29/2$.
\item \emph{Exact quotas.} Labels are grouped by their divisibility by
  small primes. The even labels are assigned to the vertices of the
  residual forest with exactly prescribed class sizes, and the
  resulting errors in degree sums and position tags are bounded
  deterministically
  (Proposition~\ref{prop:two-measure-assignment} and
  Theorem~\ref{thm:rare-type-assignment}). Probabilistic arguments are
  confined to permutations within exact classes. Larger prime divisors
  are controlled by first-moment or exponential bounds for random
  permutations, and, from the fifth logarithmic bridge on and in the
  unbounded range, by sorting labels against degrees within these
  classes (Theorem~\ref{thm:large-prime-sorting}).
\item \emph{Hall criteria with finitely many rational conditions.}
  The odd labels are placed by a perfect matching. Several sufficient
  conditions for Hall's condition are proved, each depending on a
  small number of normalized parameters
  (Theorems~\ref{thm:simultaneous-matching},
  \ref{thm:three-range-hall}, \ref{thm:direct-reserve-hall},
  and~\ref{thm:asymmetric-product-hall}, and
  Lemma~\ref{lem:two-sided-integrated-hall}).
\item \emph{Certificates for ranges of orders.} Orders through $1359$
  are covered by finite proof graphs for coprime labelings by general
  label sets, whose soundness is proved
  (Theorem~\ref{thm:finite-certificate-soundness}). Integer intervals
  through order $1982759$ are covered by bounds that hold uniformly on
  each interval (Theorem~\ref{thm:integer-interval-certificate} and its
  variants). The orders with $29/2\leq\log n\leq256$ are covered by
  finitely many closed cells in the logarithm of the order and, for
  most bridges, a normalized degree, and the orders with
  $\log n\geq256$ by an analytic argument
  (Theorem~\ref{thm:analytic-tail}).
\end{enumerate}

Some tools may be of use elsewhere: the rounding of fractional type
assignments in forests with exact class sizes and bounded degree error
(Section~\ref{sec:forest-types}), Hall criteria for type-compatibility
graphs with a reserve of universal vertices
(Sections~\ref{sec:refined-hall}, \ref{sec:direct-reserve-hall},
and~\ref{sec:adaptive-refinement-hall}), the certificate rules for
coprime labelings by arbitrary finite label sets
(Section~\ref{sec:finite-certificates}), and the directed rational
evaluation of probability bounds
(Section~\ref{sec:directed-probability}).

\subsection{Related work}\label{sec:related-work}

We recall the results most directly connected with the two parts of
our argument: reductions for finite label sets and embeddings for large
orders. For results on particular families and on other graph classes
we refer to the account in~\cite[Section~1]{HPT11}.

\subsubsection{Finite orders and particular families}

Fu and Huang~\cite[Theorem~2.8]{FH94} proved that every tree of order at
most $15$ is prime. Pikhurko extended the verified range to $34$
vertices~\cite{Pik02} and subsequently to $50$ vertices~\cite[Theorem~5]{Pik07}.
The latter two proofs proceed by explicit reductions that can be checked
by hand. Haxell, Pikhurko, and Taraz~\cite[Section~1]{HPT11} also mention
the bound $206$, claimed in an unpublished manuscript of Kuo and Fu. We
have not seen that manuscript, and the present proof does not use it.
Removing a suitable part of a tree may leave a nonconsecutive
label set, so the induction must allow general finite sets of integers.
This is also the reason for using the predicate $\cP(L)$ here.

Results for particular families follow a complementary approach.
Robertson and Small~\cite{RS09} extended coprime mapping results and used
them to label palm trees, banana trees, binomial trees, and certain spider
colonies. Such constructions exploit the structure of the chosen family;
the reductions below must instead apply to an arbitrary input tree.

\subsubsection{Asymptotic embeddings}

Let $m(n)$ be the least integer $m$ such that every tree of order $n$
has an injective coprime labeling using labels in $[m]$. Pikhurko proved
that $m(n)=(1+o(1))n$~\cite[Theorem~6]{Pik07}. Allowing unused labels
removes part of the difficulty of a spanning embedding. The prime
labeling conjecture requires the exact equality $m(n)=n$.

Haxell, Pikhurko, and Taraz~\cite[Theorems~1 and~2]{HPT11} proved that
all sufficiently large trees are prime, as a consequence of a more
general theorem for bipartite degenerate graphs with suitable separators.
Their proof separates a structural partition, an arithmetic partition,
and an embedding step. Primes above $n/2$ and the labels $2p$ with
$n/3<p\leq n/2$ provide the exceptional labels needed at the separators;
a final exchange removes the possible conflicts between $p$ and $2p$.

\subsubsection{Relation to the present argument}

Our finite reductions continue the approach of keeping track of the
entire remaining label set. The large-order construction retains the
special-label and separator method just described. Its quantitative
requirements are expressed through exact type counts, prescribed
class sizes, degree sums, and matching inequalities. These quantities
must be controlled for every tree in each claimed order range.

For rounding fractional type assignments, we use the sorted-slot
matching construction of Shmoys and Tardos~\cite{ST93}. Its application
here preserves exact class sizes and controls the change in degree
sums; Section~\ref{sec:forest-types} gives the required form with proof.
Our local sampling bounds use Hoeffding's convex comparison between
sampling with and without replacement~\cite[Theorem~4]{Hoe63}.
Section~\ref{sec:local-probability} includes a proof and derives the
finite-population and forest-degree forms used here.

The arithmetic and structural statements proved below are the general
ingredients for that task. To deduce the all-orders statement from them,
we supply the finite reduction certificates through order $1359$,
verify the parameter ranges for every larger order with $\log n\leq256$,
and prove the estimates used in the unbounded range $\log n\geq256$.
Section~\ref{sec:overview} summarizes these dependencies, and
Sections~\ref{sec:part-completion} and~\ref{sec:part-computation} make
them explicit.

\subsection{Organization}\label{sec:organization}

Section~\ref{sec:overview} describes the proof and lists the ranges of
orders with the results that cover them. Sections~\ref{sec:part-structural}
and~\ref{sec:part-arithmetic} contain the general structural,
arithmetic, and matching results. Sections~\ref{sec:part-finite},
\ref{sec:part-intervals}, \ref{sec:part-integer-bridges},
and~\ref{sec:part-log-bridges} cover the orders with $\log n\leq256$ in
increasing order: certificates through order $1359$, integer intervals
through $34999$, integer bridges through $1982759$, and logarithmic
bridges from $\log n=29/2$ to $256$. Section~\ref{sec:part-analytic}
contains uniform estimates and the argument for $\log n\geq256$; it is
placed before the bridges because several logarithmic bridges use its
tools.
Section~\ref{sec:part-completion} proves Theorem~\ref{thm:main}, and
Section~\ref{sec:part-computation} describes the computations and their
reproduction. Appendix~\ref{sec:notation-index} lists the notation, and
Appendix~\ref{sec:cover-tables} contains the degree covers of the
retained-position certificates.

\subsection{Notation}\label{sec:notation}

All graphs are finite and simple. The order of a graph is its number of
vertices. We write $V(G)$ and $E(G)$ for its vertex and edge sets,
$e(G)=\abs{E(G)}$, $N_G(v)$ for the neighborhood of $v$, and
$d_G(v)=\abs{N_G(v)}$. For a vertex set $X$, put
$N_G(X)=\bigcup_{x\in X}N_G(x)$. The notation $G[X]$ means the induced
subgraph on $X$, and $G-X=G[V(G)\setminus X]$. A forest may have isolated
vertices or be empty. Set $\Delta(\varnothing)=0$ and, for
$X\subseteq V(G)$, $\Delta_G(X)=\max_{x\in X}d_G(x)$, with value zero
when $X$ is empty. A graph is $C_4$-free if it contains no cycle of length
four, whether induced or not.

For a finite set $L$ of positive integers, a \emph{coprime labeling by $L$}
is a bijection onto $L$ with coprime labels on every edge. Write $\cP(L)$
if every forest on $\abs{L}$ vertices admits such a labeling, and
$\cP_d(L)$ if this holds for every such forest of maximum degree at most
$d$. In particular $\cP(\varnothing)$ and $\cP_d(\varnothing)$ hold.
A label $q\in L$ is \emph{universal in $L$} if it is coprime to every
element of $L\setminus\{q\}$. Universality always refers to the specified
label set.

We use $\log$ for the natural logarithm and $\log_2$ for the logarithm to
base two. For real $x\geq0$, $\pi(x)$ denotes the number of primes at most
$x$. All label sets and vertex sets in disjoint unions are understood to
be pairwise disjoint. Further notation is introduced where it is first
used; Appendix~\ref{sec:notation-index} lists the symbols used in more
than one section.

\section{Overview of the proof}\label{sec:overview}

The proof uses one reduction at every order, several ways of solving the
reduced problem, and a chain of ranges of orders, each covered by one of
these ways. This section describes the steps in the order of the
argument and lists the ranges.

\subsection{The residual problem}\label{sec:overview-residual}

Since $\floor{n/2}$ labels are even, the vertices receiving them must
form an independent set. We remove a set $K$ of at most
$b(n)=\floor{\log_2(n-2)}+2$ vertices and split the rest of the tree
into independent sets $S$ and $Z$ with $\abs S=\floor{n/2}$
(Proposition~\ref{prop:prescribed-side}). If $n\geq8$ and $(n/2,n]$ contains at
least $b(n)$ primes, the vertices of $K$ receive such primes
(Corollary~\ref{cor:even-side}); each of them is coprime to every
other label in $[n]$.

Two further families of labels have few conflicts
(Lemma~\ref{lem:special-labels}). The special even labels are the powers
of $2$ and the labels $2p$ with $p$ prime and $n/3<p\leq n/2$; each is
coprime to every odd label except possibly $p$. The special odd labels
are $1$ and the primes above $n/3$ not used on $K$; each is coprime to
every even label except possibly its double. The pairs $p,2p$ are allowed
during the construction and removed at the end by an exchange
(Lemma~\ref{lem:sparse-repair} and Corollary~\ref{cor:prime-pair-repair}).

The special even labels are placed first, on a set $U\subseteq S$. It
contains the neighbors of a reserve $W\subseteq Z$ of positions, so that
every position of $W$ accepts every odd label, and it may contain
vertices of largest degree, which bounds the degrees that remain
(Lemma~\ref{lem:cyclic-neighborhood} and
Proposition~\ref{prop:residual-problem}). What remains is a bipartite
residual forest $G=(S_0,Z_0)$. Its even side must receive the remaining
even labels bijectively, and its positions must receive distinct odd
labels, coprime along the edges of $G$. Once the even labels are fixed,
the odd labels are placed by a matching, and the reserve reduces
Hall's condition by the number of reserved positions
(Corollary~\ref{cor:residual-hall}).

\subsection{Types, quotas, and matchings}\label{sec:overview-types}

Fix a set $P$ of small odd primes. The \emph{type} of an odd label is
the set of primes of $P$ dividing it, and the type of an even label is
that of its half. The numbers of labels of each type in an interval are
compared with the product law $\mu=\bigotimes_{p\in P}
\operatorname{Bernoulli}(1/p)$ (Section~\ref{sec:arithmetic-types}).
The even labels are assigned to the vertices of $S_0$ type by type, with
exactly prescribed numbers of vertices of each type. Every position of
$Z_0$ receives a \emph{tag}, a set of primes containing every prime of
$P$ that divides a neighboring even label; in some constructions a
position is instead marked, and it later receives a special odd label.
The assignment controls the distribution of the tags and the degree sums
of the classes (Sections~\ref{sec:forest-types}
and~\ref{sec:rare-type-assignment}).

A label of type $I$ can be placed at a position with tag $J$ only if
$I\cap J=\varnothing$; primes outside $P$ may delete further pairs.
Within each exact class the even labels are then permuted, either at
random or by sorting them against the degrees of their vertices. The
number of deleted pairs is bounded in both cases
(Sections~\ref{sec:matching-estimates} and~\ref{sec:large-prime-cost}).
Hall criteria reduce the existence of the matching to lower bounds on
degrees in the compatibility graph and an upper bound on the number of
deleted pairs (Sections~\ref{sec:matching-estimates},
\ref{sec:refined-hall}, \ref{sec:direct-reserve-hall},
\ref{sec:asymmetric-reserve-hall}, and~\ref{sec:adaptive-refinement-hall}).

\subsection{How the ranges are covered}\label{sec:overview-ranges}

\emph{Orders through $1359$.} These orders are small enough for
finite certificates. The assertion $\cP([n])$ is proved by finite proof graphs whose nodes
are statements $\cP(L)$ about general label sets $L$. The rules of these
graphs are proved sound in Theorem~\ref{thm:finite-certificate-soundness},
and a program checks that every rule is applied correctly and every case
list is complete. Branch peeling
(Theorem~\ref{thm:peeled-prefix-step}) and certified rooted blocks
reduce the number of cases.

\emph{Integer intervals from $1360$ through $1982759$.} For an interval
$[\ell,u]$ of orders, the preparation constants are bounded uniformly,
and the largest residual degree $D$ takes one of finitely many values
(Lemma~\ref{lem:interval-degree-cover}). Each value of $D$ has a row
specifying a construction and its parameters. The class permutations are
random. A finite family of events covers every failure of a local degree
condition, and the deleted pairs are bounded by a first-moment or an
exponential estimate. If the sum of the failure probabilities is less
than one, some permutation succeeds
(Theorem~\ref{thm:integer-interval-certificate}, its variants with
the Hall estimates of Proposition~\ref{prop:two-sided-hall} and
Lemma~\ref{lem:single-range-hall}, and
Theorem~\ref{thm:retained-interval-criterion}). All bounds are evaluated by
directed rational arithmetic (Section~\ref{sec:directed-probability}).
From $35000$ the same scheme is used with a reserve on the label side
(Section~\ref{sec:asymmetric-reserve-hall}).

\emph{Logarithmic bridges for $29/2\leq\log n\leq256$.} Here the order
enters only through $\log n$ and, in the first fifteen bridges, a
normalized residual degree. All quantities are bounded monotonically on
closed cells: rectangles in these two coordinates for $\log n\leq48$,
and unit intervals of $\log n$ beyond. Each cell carries a finite list
of rational inequalities that imply the Hall conditions, and the cells
cover $29/2\leq\log n\leq256$. The nineteen bridges differ in the
preparation, the type assignment, the Hall criterion, or the arithmetic
bounds; Table~\ref{tab:log-bridges} lists them.

\emph{The unbounded range $\log n\geq256$.} The cutoff for the small
primes grows with $\log n$, and every estimate is proved analytically,
using monotonicity and explicit comparisons at finitely many rational
points, mostly at $\log n=256$ (Theorem~\ref{thm:analytic-tail}).

\subsection{The chain of ranges}\label{sec:overview-chain}

Table~\ref{tab:range-chain} lists the ranges. The \emph{direct} result
of a range treats the orders of that range: it is a theorem, or a
certificate proposition which, together with its criterion, proves the
assertion for these orders. A \emph{cumulative} result states the
assertion for all orders up to an endpoint; its proof applies the direct
result to the new orders and the preceding cumulative result to the
others. For the logarithmic range, Theorem~\ref{thm:log-bridges-all}
instead collects the nineteen bridges. Adjacent integer ranges share no order and omit none; the
passage to the logarithmic range uses
$\log1982759<29/2<\log1982760$.

\begin{table}[ht]
\centering
\caption{The ranges of orders and the results that cover them.}
\label{tab:range-chain}
\par\vspace{.5\baselineskip}
\renewcommand{\arraystretch}{1.15}
\small
\begin{tabular}{@{}lp{2.3cm}p{3.6cm}p{2.3cm}@{}}
\hline
Orders&
  \raggedright Method&
  \raggedright Direct result&
  \raggedright Cumulative result\tabularnewline
\hline
$1\leq n\leq1359$&
  \raggedright finite proof graphs&
  \raggedright Theorem~\ref{thm:initial-prefix-1359}&
  \raggedright Theorem~\ref{thm:initial-prefix-1359}\tabularnewline
$1360\leq n\leq1521$&
  \raggedright interval certificate&
  \raggedright Theorem~\ref{thm:interval-1360-1521}&
  \raggedright Corollary~\ref{cor:initial-prefix-1521}\tabularnewline
$1522\leq n\leq3499$&
  \raggedright interval certificates&
  \raggedright Propositions~\ref{prop:verified-forest-continuation}, \ref{prop:verified-exponential-continuation}, \ref{prop:verified-single-range-continuation}&
  \raggedright Theorem~\ref{thm:initial-prefix-3499}\tabularnewline
$3500\leq n\leq34999$&
  \raggedright retained positions&
  \raggedright Theorem~\ref{thm:retained-interval-criterion} and Proposition~\ref{prop:verified-retained-continuations}&
  \raggedright Theorem~\ref{thm:initial-prefix-34999}\tabularnewline
$35000\leq n\leq1982759$&
  \raggedright integer bridges&
  \raggedright Propositions~\ref{prop:integer-bridge-certificate}, \ref{prop:profile-block-certificate}, \ref{prop:full-bridge-certificate}, \ref{prop:upper-integer-certificates}&
  \raggedright Theorem~\ref{thm:initial-prefix-1982759}\tabularnewline
$29/2\leq\log n\leq256$&
  \raggedright logarithmic bridges&
  \raggedright Table~\ref{tab:log-bridges}&
  \raggedright Theorem~\ref{thm:log-bridges-all}\tabularnewline
$\log n\geq256$&
  \raggedright analytic estimates&
  \raggedright Theorem~\ref{thm:analytic-tail}&
  \raggedright --\tabularnewline
\hline
\end{tabular}
\end{table}

The cumulative results are joined in Section~\ref{sec:part-completion}.
Except for Theorems~\ref{thm:initial-prefix-1359}
and~\ref{thm:analytic-tail}, each direct result above either is a
certificate proposition, which
states that the finite conditions of a criterion hold, or combines such
a proposition with the criterion. The computations behind the
certificate propositions are described in
Section~\ref{sec:part-computation}.

\section{Structural reductions and special labels}
\label{sec:part-structural}

This part reduces the construction of a prime labeling to a matching
problem. Section~\ref{sec:finite-reductions} collects local
reductions for coprime labelings by general finite label sets; the
finite certificates of Section~\ref{sec:part-finite} use them.
Section~\ref{sec:repair} introduces the special labels
and the exchange that removes the pairs $p,2p$
(Lemma~\ref{lem:sparse-repair}). Sections~\ref{sec:balance}
and~\ref{sec:matching} choose the balanced independent set, the
separator, and the reserve, and state the residual problem
(Proposition~\ref{prop:residual-problem}) with its Hall condition
(Corollary~\ref{cor:residual-hall}).

\subsection{Local reductions for finite label sets}\label{sec:finite-reductions}

Removing vertices from a tree usually leaves a forest. We therefore state
the extension rules for forests and record first why this does not enlarge
the terminal family that needs to be checked.

\begin{lemma}\label{lem:forest-completion}
Every nonempty forest can be extended to a tree on the same vertex set by
adding edges. If its maximum degree is at most $d$, where $d\geq2$, the
extension can be chosen to have maximum degree at most $d$.
\end{lemma}

\begin{proof}
Order the components as $F_1,\ldots,F_j$. In each nontrivial component,
choose two distinct leaves, one as a left attachment vertex and one as a
right attachment vertex. In an isolated-vertex component use its sole
vertex for both roles. For $i<j$, join the right attachment vertex of
$F_i$ to the left attachment vertex of $F_{i+1}$.

The resulting graph is connected and has one fewer edge than vertices,
so it is a tree. A selected leaf gains at most one edge, and an isolated
vertex gains at most two. All other degrees are unchanged. This proves
both assertions.
\end{proof}

Consequently, for a nonempty label set, $\cP(L)$ is equivalent to the
corresponding assertion for trees. The same equivalence holds for
$\cP_d(L)$ when $d\geq2$. A degree-restricted assertion is nevertheless
distinct from $\cP(L)$: a labeling table omitting high-degree trees does
not establish the latter.

\subsubsection{Deleting a leaf and its neighbor}

The next lemma makes a two-vertex deletion useful even when the available
label set supports only a restricted degree bound.

% Source: BOUNDED_EDGE_DEGREE_REDUCTION.md, Section 1.
\begin{lemma}\label{lem:leaf-degree}
Let $T$ be a tree on $m\geq4$ vertices. There is a leaf $x$ with neighbor
$y$ such that
\[
    \Delta(T-\{x,y\})
    \leq \min\left\{\Delta(T),\floor*{\frac{m-1}{2}}\right\}.
\]
\end{lemma}

\begin{proof}
Put $h=\floor{(m-1)/2}$. If $\Delta(T)\leq h$, any leaf and its neighbor
satisfy the conclusion. Otherwise choose a vertex $y$ of degree $s>h$.

We first show that $y$ has a leaf neighbor. Suppose otherwise. Each of
its $s$ neighbors has a further neighbor distinct from $y$. Since $T$ is
a tree, these further vertices are distinct and lie outside
$\{y\}\cup N_T(y)$. Thus $T$ has at least $1+2s>m$ vertices, a
contradiction. Choose a leaf $x$ adjacent to $y$.

Let $v\in V(T)\setminus\{x,y\}$. If $v$ is adjacent to $y$, their
neighborhoods are disjoint: a common neighbor would form a triangle.
Hence $d_T(v)+s\leq m$, and deleting $y$ gives
\[
    d_{T-\{x,y\}}(v)\leq m-1-s.
\]
If $v$ is not adjacent to $y$, both neighborhoods avoid $v$ and $y$ and
intersect in at most one vertex, since two common neighbors would form
a four-cycle. Thus $d_T(v)+s\leq m-1$, giving the same bound.

Finally, $s\geq h+1$ implies $m-1-s\leq m-2-h\leq h$.
Deletion also preserves the bound by $\Delta(T)$.
\end{proof}

\begin{corollary}\label{cor:leaf-label-extension}
Let $L$ have $m\geq4$ elements, and let $q,a\in L$ be distinct, with
$q$ universal in $L$. Let $d\geq2$ and let $d'$ be a nonnegative integer
such that
\[
 d'\geq\min\left\{d,\floor*{\frac{m-1}{2}}\right\}.
\]
Then $\cP_{d'}(L\setminus\{q,a\})$ implies $\cP_d(L)$.
\end{corollary}

\begin{proof}
Complete the given forest to a tree of maximum degree at most $d$ by
Lemma~\ref{lem:forest-completion}, and choose $x,y$ by
Lemma~\ref{lem:leaf-degree}. Label $x$ by $a$ and $y$ by $q$. The
remaining forest has maximum degree at most $d'$, so the hypothesis
labels it by $L\setminus\{q,a\}$. Every edge incident with $y$ has
coprime labels, and $x$ has no other neighbor. Restricting to the original
forest proves the assertion.
\end{proof}

\subsubsection{Removing a block through one vertex}

A larger deletion must control both the size and the internal degree of
the removed forest. For an integer $b\geq3$, define
\[
 \mathcal K_b=\{b\}\cup\{k\in2\Z:b<k\leq2b-2\}.
\]
The even sizes in this set will allow arithmetic blocks to be assembled
from pairs of consecutive integers.

% Source: DEGREE_BOUNDED_PACKET_REDUCTION.md, Section 1.
\begin{lemma}\label{lem:detached-block}
Let $b\geq3$, and let $T$ be a tree of order $n\geq b+2$. There are a
vertex $u$ and a union $A$ of components of $T-u$ such that
\[
 \abs{A}\in\mathcal K_b,
 \qquad \Delta(T[A])\leq b-1.
\]
\end{lemma}

\begin{proof}
Root $T$ at a leaf. Choose a vertex $v$ of greatest depth whose rooted
subtree has at least $b+1$ vertices. The root is not such a deepest
vertex, since its unique child has a subtree of order $n-1\geq b+1$.
Thus $v$ has a parent, and every child subtree of $v$ has at most $b$
vertices.

If some union of these child subtrees has order in $\mathcal K_b$, take
it as $A$ and put $u=v$. Each component of $T[A]$ then has at most $b$
vertices, giving the degree bound.

Suppose no such union exists. No child subtree has order $b$. Pair the
child subtrees of odd order, leaving at most one unpaired, and treat each
even-order child subtree as a separate group. Every group has even total
order at most $2b-2$. Its order must be less than $b$, or it would be an
allowed union. If the groups together had total order at least $b$,
accumulating them until first reaching $b$ would give an even total
between $b$ and $2b-2$, again an allowed union. Hence the total order
$g$ of all groups satisfies $g\leq b-1$.

The child subtrees together have at least $b$ vertices. There must
therefore be an unpaired odd-order subtree, of order $a\leq b-1$.
The odd integer $g+a$ lies between $b$ and $2b-2$, and so is at most
$2b-3$. It follows that the whole rooted subtree of $v$ has even order
$g+a+1\in\mathcal K_b$. Take this subtree as $A$ and let $u$ be the
parent of $v$.

It remains to check the degree at $v$. If $r$ is its number of children,
the grouped subtrees contain $r-1$ children, so $r-1\leq g\leq b-1$.
Equality $r=b$ would force all grouped child subtrees to be singletons.
The unpaired subtree, together with $b-a$ of these $b-1$ singletons,
would then have order $b$, a contradiction. Thus $r\leq b-1$.
Every other vertex lies in a child subtree of order at most $b-1$;
even its attachment vertex has degree at most $b-1$ after adding the
edge to $v$. This completes the proof.
\end{proof}

\begin{proposition}[Extension through a universal label]
\label{prop:block-extension}
Let $\abs{L}=n\geq b+2$, where $b\geq3$. For every
$k\in\mathcal K_b$ with $k\leq n-1$, suppose there is a partition
\[
 L=\{q_k\}\sqcup B_k\sqcup C_k,
 \qquad \abs{B_k}=k,
\]
such that $q_k$ is universal in $L$, $\cP_{b-1}(B_k)$ holds, and
$\cP(C_k)$ holds. Then $\cP(L)$ holds.
\end{proposition}

\begin{proof}
By Lemma~\ref{lem:forest-completion}, it suffices to label a tree $T$
of order $n$. Apply Lemma~\ref{lem:detached-block} and put $k=\abs{A}$.
Give $u$ the label $q_k$, label $T[A]$ by $B_k$, and label
$T-(A\cup\{u\})$ by $C_k$. These assignments use every label exactly
once. Since $A$ is a union of components of $T-u$, every edge between
one of the two forests and its exterior is incident with $u$.
Universality of $q_k$ checks all such edges.
\end{proof}

The partitions may depend on $k$. Thus a finite certificate needs a
valid arithmetic partition for each permitted size, while the structural
lemma applies to every input tree. No enumeration of all trees of order
$n$ is required at this step.

\subsection{Exceptional pairs and special labels}\label{sec:repair}

We temporarily relax coprimality for a small family of pairs. The next lemma
shows that a labeling obtained under this relaxation can be repaired without
changing which vertices carry even labels. It is a finite $C_4$-free variant
of the switching argument in \cite[Section~4]{HPT11}.

For a finite set $L$ of distinct positive integers, define
\[
    \mathcal Q(L)=\left\{q\in L:
    \begin{array}{l}
        q\geq3\text{ is odd},\quad 2q\in L,\\
        \gcd(q,t)=1\text{ for every }t\in L\setminus\{q,2q\}
    \end{array}\right\}.
\]
Given a bijection $f\colon V(G)\to L$, call an edge $uv$ \emph{bad} if
$\gcd(f(u),f(v))>1$.

% Mathematical source: ../SPARSE_LABEL_REPAIR.md, Sections 1--2 and 6.
\begin{lemma}\label{lem:sparse-repair}
Let $G$ be a $C_4$-free graph, and let $L$ be a set of $\abs{V(G)}$ distinct
positive integers containing $1$ and two distinct powers of $2$ greater
than $1$.
Suppose that $f\colon V(G)\to L$ is a bijection and that every bad edge has
labels $q,2q$ for some $q\in \mathcal Q(L)$. Then $G$ has a bijective labeling by $L$
in which every edge has coprime labels. Moreover, the parity of the label at
each vertex can be preserved.
\end{lemma}

\begin{proof}
We show how to remove one bad edge without creating another. Repeating this
operation will prove the lemma.

First observe that the vertices carrying even labels form an independent
set. An edge between two such vertices would be bad and would not have the
allowed form. This remains true after any exchanges that preserve parity
at each vertex.

Fix a bad edge $xy$, with $f(x)=q$ and $f(y)=2q$. Let $u$ be the vertex
carrying $1$, and let $v_r,v_s$ carry two distinct powers of $2$ greater
than $1$, denoted
by $r$ and $s$. These five vertices are distinct. By the definition of
$\mathcal Q(L)$, the only other label in $L$ that is not coprime to $q$ is $2q$.
Among odd labels, the only one that is not coprime to $2q$ is $q$.

\emph{Case 1: $x$ is not adjacent to at least one of $v_r,v_s$.}
Suppose that $xv_r\notin E(G)$. Exchange the labels $2q$ and $r$.
All neighbors of $y$ and $v_r$ carry odd labels. The new label at $y$ is
therefore coprime to every neighboring label. The new label at $v_r$ can
conflict only with $q$, which is at the nonadjacent vertex $x$. Thus every
affected edge is good.

\emph{Case 2: $x$ is adjacent to both $v_r$ and $v_s$.}
The vertex $u$ cannot be adjacent to both: otherwise
\[
    x,v_r,u,v_s,x
\]
would be a cycle of length four. After interchanging the names $r,s$ if
necessary, assume that $uv_r\notin E(G)$. Exchange $q$ with $1$, and $2q$
with $r$.
The new labels at $x$ and $y$ are $1$ and a power of $2$, respectively,
so their incident edges are good. The only possible conflict involving the
new label $q$ at $u$ is with the new label $2q$ at $v_r$; these vertices
are not adjacent. Conversely, all neighbors of $v_r$ have odd labels, so
$q$ is its only possible conflicting label. All affected edges are again
good.

In both cases the exchanges preserve the label set and the parity at every
vertex. The chosen bad edge disappears, and no new bad edge is introduced.
The number of bad edges is a nonnegative integer and strictly decreases.
The procedure therefore terminates with the required labeling.
\end{proof}

\begin{figure}[tbp]
\centering
\begin{tikzpicture}[
 vertex/.style={circle,draw,fill=white,minimum size=7mm,inner sep=1pt},
 every node/.style={font=\small},x=9mm,y=11mm]
\begin{scope}
 \node[vertex,label=above:$x$] (x) at (0,1.4) {$q$};
 \node[vertex,label=below:$y$] (y) at (-1.3,0) {$2q$};
 \node[vertex,label=below:$v_s$] (s) at (0,0) {$s$};
 \node[vertex,label=below:$v_r$] (r) at (1.3,0) {$r$};
 \node[vertex,label=above:$u$] (u) at (2.6,1.4) {$1$};
 \draw (x)--(y) (x)--(s) (x)--(r);
 \draw[densely dashed] (u)--(r);
\end{scope}
\draw[-{Stealth[length=2.5mm]},thick] (3.4,.7)--(4.5,.7);
\begin{scope}[xshift=5.8cm]
 \node[vertex,label=above:$x$] (x2) at (0,1.4) {$1$};
 \node[vertex,label=below:$y$] (y2) at (-1.3,0) {$r$};
 \node[vertex,label=below:$v_s$] (s2) at (0,0) {$s$};
 \node[vertex,label=below:$v_r$] (r2) at (1.3,0) {$2q$};
 \node[vertex,label=above:$u$] (u2) at (2.6,1.4) {$q$};
 \draw (x2)--(y2) (x2)--(s2) (x2)--(r2);
 \draw[densely dashed] (u2)--(r2);
\end{scope}
\end{tikzpicture}
\caption{The second case of Lemma~\ref{lem:sparse-repair}. Vertex names
are unchanged; the entries inside vertices are labels. Solid segments are
edges required in this case, and the dashed segment marks the chosen
nonedge $uv_r$. Other adjacencies are not specified. The exchanges move
$q$ and $2q$ to nonadjacent vertices.}
\label{fig:repair}
\end{figure}
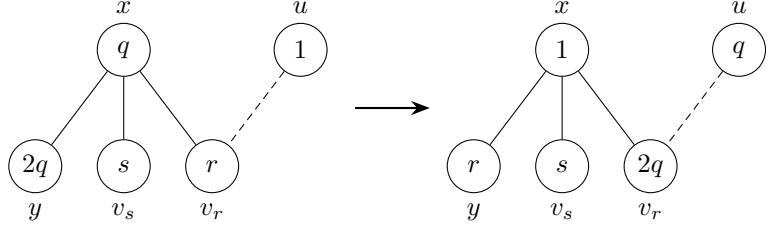

\begin{corollary}\label{cor:prime-pair-repair}
Under the assumptions on $G$ and $L$ in Lemma~\ref{lem:sparse-repair}, suppose
that every bad edge has labels $p,2p$, where $p$ is an odd prime and
$3p>\max L$. Then the labeling can be repaired while preserving parity at
every vertex.
\end{corollary}

\begin{proof}
For each such pair, the only multiples of $p$ in $L$ are $p$ and $2p$.
Since $p$ is prime, it follows that $p\in \mathcal Q(L)$. Apply
Lemma~\ref{lem:sparse-repair}.
\end{proof}

\subsubsection{The arithmetic families}

Fix $n\geq8$ and put
\[
 \begin{aligned}
 E_n&=\{2,4,\ldots,2\floor{n/2}\},\qquad O_n=[n]\setminus E_n,\\
 P_n^+&=\{p:n/2<p\leq n,\ p\text{ prime}\},\\
 P_n^0&=\{p:n/3<p\leq n/2,\ p\text{ prime}\}.
 \end{aligned}
\]
For distinct $a,b\in[n]$, write $a\sim_n b$ if either
$\gcd(a,b)=1$ or $\{a,b\}=\{p,2p\}$ for some $p\in P_n^0$.
We call this relation \emph{relaxed compatibility}. Corollary~
\ref{cor:prime-pair-repair} shows that a bijective labeling of a tree
by $[n]$ satisfying this relation on every edge can be repaired to a
prime labeling. The labels $1,2,4$ required for the repair are present.

\begin{lemma}\label{lem:special-labels}
Let $n\geq8$, let $\Psep\subseteq P_n^+$, and put $O=O_n\setminus\Psep$.
Define
\begin{align}
 E_*&=\{2^j:j\geq1,\ 2^j\leq n\}
          \cup\{2p:p\in P_n^0\},\label{eq:special-even}\\
 O_*&=\{1\}\cup
       \bigl(\{p:n/3<p\leq n,\ p\text{ prime}\}\setminus\Psep\bigr).
       \label{eq:special-odd}
\end{align}
Then the following hold.
\begin{enumerate}
 \item Every label in $\Psep$ is universal in $[n]$.
 \item Every label in $E_*$ is relaxed-compatible with every label in $O$.
 \item Every label in $O_*$ is relaxed-compatible with every label in $E_n$.
\end{enumerate}
Moreover, the two families in \eqref{eq:special-even} are disjoint, and
\begin{align}
 w:=\abs{E_*}&=\floor{\log_2 n}+\pi(n/2)-\pi(n/3),
       \label{eq:even-supply}\\
 u:=\abs{O_*}&=1+\pi(n)-\pi(n/3)-\abs{\Psep}.
       \label{eq:odd-supply}
\end{align}
\end{lemma}

\begin{proof}
A prime $p>n/2$ has no other positive multiple in $[n]$, proving the
first assertion. A power of two is coprime to every odd number. If
$p\in P_n^0$, the only multiples of $p$ in $[n]$ are $p$ and $2p$;
thus $p$ is the only odd label that can conflict with $2p$, and this
pair is allowed. This proves the second assertion. The same observation,
together with the universality of $1$ and of primes greater than $n/2$,
proves the third.

Since $n\geq8$, every prime greater than $n/3$ is odd. Therefore none
of the labels $2p$ in \eqref{eq:special-even} is a power of two, and
all labels in \eqref{eq:special-odd} are odd. Counting the indicated
primes gives \eqref{eq:even-supply} and \eqref{eq:odd-supply}.
\end{proof}

The even family is universal only against odd labels, and only for the
relaxed relation. We will assign it to an independent vertex set and
invoke the repair lemma after every label has been placed.

\subsection{Balanced independent sets and a reserve}
\label{sec:balance}

We need an independent set of exactly $\floor{n/2}$ vertices, with an
independent complement after deleting a small separator. A centroid
ordering lets us change a two-coloring one vertex at a time while keeping
all conflicts inside a set of logarithmic size.

\subsubsection{An ordering with a small boundary}

For an ordering $v_1,\ldots,v_N$ of a graph $G$, define its boundary
after position $j$ by
\[
 \partial_j=\{v_i:i\leq j,\ N_G(v_i)\cap\{v_{j+1},\ldots,v_N\}
                         \ne\varnothing\}.
\]
In particular $\partial_0=\partial_N=\varnothing$.

\begin{lemma}\label{lem:centroid-order}
Every forest on $N\geq1$ vertices has an ordering for which
\[
 \abs{\partial_j}\leq\floor{\log_2 N}
 \qquad(0\leq j\leq N).
\]
\end{lemma}

\begin{proof}
First consider a tree. We recall why it has a vertex whose deletion
leaves components of order at most $N/2$. Starting at any vertex, move
to the neighbor in a component of order greater than $N/2$, if such a
component exists. After crossing an edge in this way, the component on
the side just left has order less than $N/2$, so that edge cannot be
crossed back. A walk in a tree that revisits a vertex must backtrack
along an edge. The procedure therefore terminates, at a vertex with
the required property. Call it a centroid.

Put the centroid first, and then list the components of its deletion
consecutively, applying the same rule within each component. We prove
the bound by induction on $N$. For $N=1$ it is immediate. During the
processing of a component $C$, only the centroid can contribute to
the boundary from outside $C$. Since $\abs{C}\leq\floor{N/2}$, the
inductive bound is at most
\[
 1+\floor*{\log_2\floor{N/2}}\leq\floor{\log_2 N}.
\]
Cuts between components have boundary of size at most one. Finally,
for a forest process its components consecutively. Completed components
have no neighbors in those still unprocessed, so the same bound applies.
\end{proof}

% Source: ORIGINAL_INTERVAL_RESERVE.md, Section 1.
\begin{proposition}[Prescribing an independent side]
\label{prop:prescribed-side}
Let $T$ be a tree of order $n\geq3$, let $r$ be a specified leaf, and
let $a,n-a$ be the sizes of its bipartition classes. For every integer
$t$ between $a$ and $n-a$, inclusive, there is a partition
\[
 V(T)=S\sqcup Z\sqcup K
\]
such that $S$ and $Z$ are independent, $r\in S$, $\abs{S}=t$, and
\begin{equation}\label{eq:separator-size}
 1\leq\abs{K}\leq b(n):=\floor{\log_2(n-2)}+2.
\end{equation}
The neighbor of $r$ belongs to $K$.
\end{proposition}

\begin{proof}
Let $A,B$ be the bipartition classes with $r\in A$, and write
$\abs{A}=a$. Let $u$ be the neighbor of $r$. Keep $r$ in $S$ and $u$
in $K$ throughout. Initially put $S=A$, $Z=B\setminus\{u\}$, and
$K=\{u\}$.

Order the vertices of $F=T-\{r,u\}$ by Lemma~\ref{lem:centroid-order}.
Process them in this order. To process a vertex, first move it from its
current class to $K$. Then, one at a time, release every processed vertex
having no unprocessed neighbor in $F$, placing it in the class opposite
its original class. A released vertex is not moved again.

At each stage, unprocessed vertices retain their original classes and
released vertices have exchanged classes. There is no edge between these
two sets: such an edge would have prevented the processed endpoint from
being released. Edges within either set remain properly two-colored,
and the other vertices are in $K$. Edges incident with $u$ cause no
problem, and $r$ has no neighbor other than $u$. Thus $S$ and $Z$
remain independent throughout.

After all releases at a given position, $K\setminus\{u\}$ is exactly
the processed boundary in $F$. Immediately after adding the next vertex,
its size has increased by at most one. Lemma~\ref{lem:centroid-order}
therefore gives \eqref{eq:separator-size} at every intermediate stage.

Initially $\abs{S}=a$. At the end every vertex of $F$ has changed
class, so
\[
 S=\{r\}\cup(B\setminus\{u\}),\qquad \abs{S}=n-a.
\]
Each individual move changes $\abs{S}$ by at most one. Hence it attains
every integer between these endpoints. Stop when it first attains $t$.
\end{proof}

\begin{corollary}\label{cor:even-side}
Let $n\geq8$ and suppose
\begin{equation}\label{eq:prime-separator-budget}
 \pi(n)-\pi(n/2)\geq b(n).
\end{equation}
Every tree $T$ of order $n$, with a specified leaf $r$, admits a partition
$V(T)=S\sqcup Z\sqcup K$ satisfying Proposition~\ref{prop:prescribed-side}
with $\abs{S}=m:=\floor{n/2}$. The vertices of $K$ can be labeled
bijectively by a set $\Psep\subseteq P_n^+$ so that every incident edge will
have coprime labels, whatever distinct remaining labels are subsequently
assigned. The remaining label sets are $E_n$ for $S$ and
$O=O_n\setminus\Psep$ for $Z$, and
\begin{equation}\label{eq:remaining-side-sizes}
 \abs{Z}=q=\ceil{n/2}-\abs{K},
 \qquad m-b(n)\leq q\leq m.
\end{equation}
Moreover, $r$ is isolated in $T-K$.
\end{corollary}

\begin{proof}
The integer $m$ lies between the two bipartition sizes, so apply
Proposition~\ref{prop:prescribed-side} with $t=m$. Condition
\eqref{eq:prime-separator-budget} supplies $\abs{K}$ distinct primes
in $P_n^+$. They are universal by Lemma~\ref{lem:special-labels}.
They are also odd, so all $m$ even labels remain available for $S$.
Counting the unused odd labels gives \eqref{eq:remaining-side-sizes},
using $1\leq\abs{K}\leq b(n)$. The final assertion follows because
the neighbor of $r$ lies in $K$.
\end{proof}

\subsubsection{Choosing positions with few neighbors}

The same ordering gives more than a separator. It also finds many
positions whose combined neighborhood is small. This is more economical
than bounding the degree of each position separately.

% Source: CYCLIC_NEIGHBORHOOD_AND_CORE.md, Section 1.
\begin{lemma}[Cyclic neighborhood bound]\label{lem:cyclic-neighborhood}
Let $G$ be a bipartite graph with parts $A,B$, where $\abs{A}=M$ and
$\abs{B}=Q>0$. Suppose $G$ has a vertex ordering with
$\abs{\partial_j}\leq t$ at every cut, where $t$ is a nonnegative
integer. For every integer $h$ with $1\leq h\leq Q-t$, there is a set
$W\subseteq B$ such that
\begin{equation}\label{eq:cyclic-bound}
 \abs{W}=h,
 \qquad \abs{N_G(W)}\leq
          \floor*{\frac{M(h+t)}{Q}}+t.
\end{equation}
\end{lemma}

\begin{proof}
Read the given ordering cyclically. Divide it into $Q$ blocks, each
ending at a vertex of $B$ and containing all vertices after the preceding
vertex of $B$. Each block contains exactly one vertex of $B$; the blocks
together contain all $M$ vertices of $A$.

Put $\ell=h+t$. Among the $Q$ cyclic intervals consisting of $\ell$
consecutive blocks, each vertex of $A$ occurs exactly $\ell$ times.
Thus some interval $I$ satisfies
\[
 \abs{A\cap I}\leq\floor*{\frac{M\ell}{Q}}.
\]
Let $a,b$ denote its initial and final cuts in the original linear
ordering. Remove from $B\cap I$ every vertex in $\partial_b$. Since
$\abs{B\cap I}=\ell$ and $\abs{\partial_b}\leq t$, at least $h$
vertices remain; choose any $h$ of them as $W$.

We claim that
\begin{equation}\label{eq:cyclic-containment}
 N_G(W)\subseteq (A\cap I)\cup(A\cap\partial_a).
\end{equation}
For an interval not crossing the end of the ordering, an external
neighbor before $I$ lies in $\partial_a$. An external neighbor after
$I$ would put its endpoint in $W$ into $\partial_b$, contrary to the
choice of $W$.

If $I$ crosses the end, it consists of a terminal segment and an initial
segment. A vertex of $W$ in the initial segment cannot have a neighbor
in the intervening complement, because it was excluded from
$\partial_b$. An external neighbor of a vertex in the terminal segment
belongs to $\partial_a$. Finally, if $\ell=Q$, the interval is the
whole vertex set and has no external neighbors. This proves
\eqref{eq:cyclic-containment} in every case, and hence
\eqref{eq:cyclic-bound}.
\end{proof}

For a forest one may take $t=\floor{\log_2(M+Q)}$ by
Lemma~\ref{lem:centroid-order}. When $M\leq Q$, the bound becomes
$\abs{N_G(W)}\leq h+2t$. Thus reserving $h$ positions costs only an
additional logarithmic number of vertices on the other side.

\subsection{Reduction to a matching problem}\label{sec:matching}

We now combine the structural and arithmetic preparations. The conclusion
will be conditional only on labeling a smaller bipartite forest. In
particular, all edges outside that forest and every unused label will be
accounted for explicitly.

\subsubsection{The graph of available positions}

Suppose a partial labeling of a graph leaves an independent set $X$ of
unlabeled vertices and a set $L$ of unused labels, with
$\abs{X}=\abs{L}=q$. Assume every edge with both endpoints already
labeled satisfies the required compatibility relation. Form a bipartite
graph $H$ with parts $X,L$ by putting $x\lambda\in E(H)$ precisely when
$\lambda$ is compatible with the label of every neighbor of $x$.
A perfect matching in $H$ is exactly a completion of the partial labeling.
There are no edges inside $X$, so no additional conditions are needed.

We use the following form of the matching criterion, including its proof
to make the subsequent reduction self-contained.

\begin{lemma}[Hall's criterion]\label{lem:hall}
A bipartite graph $H$ with parts $X,L$ of the same finite size has a
perfect matching if and only if
\begin{equation}\label{eq:hall}
 \abs{N_H(Q)}\geq\abs{Q}\qquad\text{for every }Q\subseteq L.
\end{equation}
\end{lemma}

\begin{proof}
Necessity follows because the matched partners of $Q$ are distinct
members of $N_H(Q)$. For sufficiency, let $M$ be a matching of maximum
size. Suppose some label $\lambda\in L$ is unmatched. Starting at
$\lambda$, follow nonmatching edges from $L$ to $X$ and matching edges
from $X$ to $L$. If an unmatched vertex of $X$ is reached, a shortest
such directed path alternates between nonmatching and matching edges.
Exchanging the two kinds of edge along the path increases $\abs{M}$,
which is impossible.

Let $Q\subseteq L$ and $Y\subseteq X$ be the reachable sets. Every
vertex of $Y$ is matched, and its matched partner belongs to
$Q\setminus\{\lambda\}$. Conversely, every vertex of
$Q\setminus\{\lambda\}$ is reached along its matching edge. These
edges give a bijection between $Y$ and $Q\setminus\{\lambda\}$.
Every neighbor of $Q$ lies in $Y$: nonmatching edges can be followed,
and the matching neighbor of a reachable label was already reached.
Thus $N_H(Q)=Y$ and $\abs{N_H(Q)}=\abs{Q}-1$, contradicting
\eqref{eq:hall}.
\end{proof}

For later use, a Hall obstruction can equivalently be written as a pair
$Q\subseteq L$, $Y\subseteq X$ with no edges between them and
$\abs{Q}+\abs{Y}>q$. Indeed, a deficient $Q$ gives
$Y=X\setminus N_H(Q)$; conversely, any such $Y$ forces
$\abs{N_H(Q)}\leq q-\abs{Y}<\abs{Q}$. We will refer to the absent
pairs $Q\times Y$ as an \emph{empty rectangle}.

\begin{lemma}[Matching with a reserve]\label{lem:reserve-matching}
Let $H$ have parts $X,L$ with $\abs{X}=\abs{L}=q$. Suppose
$W\subseteq X$ and $R\subseteq L$, of sizes $h$ and $u$, satisfy
\[
 W\times L\subseteq E(H),\qquad X\times R\subseteq E(H).
\]
Put $H_0=H[X\setminus W,L\setminus R]$. The following are equivalent:
\begin{enumerate}
 \item $H$ has a perfect matching;
 \item $H_0$ has a matching of size at least $\max\{0,q-h-u\}$;
 \item $\abs{N_{H_0}(Q)}\geq\abs{Q}-h$ for every
       $Q\subseteq L\setminus R$.
\end{enumerate}
\end{lemma}

\begin{proof}
A perfect matching has at most $h+u$ edges incident with $W\cup R$,
which proves the second assertion from the first.

Conversely, let a matching in $H_0$ have size $v\geq q-h-u$.
It leaves $q-h-v\leq u$ unmatched vertices of $X\setminus W$;
match them to distinct labels in $R$. It also leaves
$q-u-v\leq h$ unmatched labels of $L\setminus R$; match them
to distinct vertices of $W$. The numbers still unused in $W$ and
$R$ are both $h+u+v-q$, and they can be matched arbitrarily.
This produces a perfect matching in $H$.

Finally, any set of labels meeting $R$ has neighborhood $X$ and
satisfies Hall's condition automatically. For nonempty
$Q\subseteq L\setminus R$, its neighborhood in $H$ is the disjoint
union of $W$ and $N_{H_0}(Q)$. Lemma~\ref{lem:hall} therefore proves
the equivalence of the first and third assertions. The empty set
satisfies the third assertion as well.
\end{proof}

\subsubsection{Preparing the residual forest}

Use the partition and labels supplied by Corollary~\ref{cor:even-side},
and let $F=T-K$, with parts $S,Z$. Retain the notation
\[
 \abs{S}=m,\quad \abs{Z}=q,\quad
 E=E_n,\quad O=O_n\setminus\Psep,
\]
and let $E_*,O_*,w,u$ be as in Lemma~\ref{lem:special-labels}.
For this construction $q>0$: indeed $Z$ has as many vertices as there
are odd labels outside $\Psep$, and the label $1$ is outside $\Psep$.

\begin{proposition}[Residual labeling problem]\label{prop:residual-problem}
Choose integers $a,h$ such that $0\leq a\leq w$, $h\geq1$, and,
on putting $t=\floor{\log_2(m-a+q)}$,
\begin{equation}\label{eq:reserve-budget}
 h+t\leq q,
 \qquad
 \floor*{\frac{(m-a)(h+t)}q}+t\leq w-a.
\end{equation}
Choose also an integer $d$ with $0\leq d\leq\min\{u,q-h\}$.
There are sets $U\subseteq S$ and disjoint sets $W,D\subseteq Z$,
and a set $R\subseteq O_*$, with
\[
 \abs{U}=w,\qquad\abs{W}=h,\qquad\abs{D}=\abs{R}=d,
 \qquad N_F(W)\subseteq U,
\]
such that, on setting
\begin{equation}\label{eq:residual-definitions}
 \begin{aligned}
 S_0&=S\setminus U,& Z_0&=Z\setminus(W\cup D),\\
 E_0&=E\setminus E_*,& O_0&=O\setminus R,
 \end{aligned}
 \qquad G=F[S_0\cup Z_0],
\end{equation}
one has
\begin{align}
 \abs{E_0}&=\abs{S_0}=m-w,
 &\abs{O_0}&=\abs{Z_0}+h=q-d,\label{eq:residual-sizes}\\
 \Delta_G(S_0)&\leq\floor*{\frac{m+q-1}{a+1}},
 &\Delta_G(Z_0)&\leq\floor*{\frac{m+q-1}{d+1}}.
 \label{eq:residual-degrees}
\end{align}
Every bijection $f_0\colon S_0\to E_0$ and injection
$g_0\colon Z_0\to O_0$ satisfying
\begin{equation}\label{eq:residual-coprimality}
 \gcd(f_0(s),g_0(z))=1\qquad\text{for every }sz\in E(G)
\end{equation}
can be converted into a prime labeling of $T$.
\end{proposition}

\begin{figure}[tbp]
\centering
\begin{tikzpicture}[
 box/.style={draw,rounded corners=1pt,minimum width=31mm,
             minimum height=10mm,align=center,fill=white},
 every node/.style={font=\small},x=1cm,y=1cm]
\node at (0,1.65) {even positions $S$};
\node at (6,1.65) {odd positions $Z$};
\node[box] (u) at (0,0.9) {$U\longleftrightarrow E_*$};
\node[box] (s) at (0,-1.4) {$S_0\longleftrightarrow E_0$};
\node[box] (w) at (6,0.9) {$W$: $h$ reserved positions};
\node[box] (d) at (6,-0.25) {$D\longleftrightarrow R\subseteq O_*$};
\node[box] (z) at (6,-1.4) {$Z_0\lhook\joinrel\longrightarrow O_0$};
\draw (u.east)--(w.west);
\draw (u.east)--(d.west);
\draw (u.east)--(z.west);
\draw (s.east)--(d.west);
\draw[thick] (s.east)--node[below,fill=white,inner sep=2pt]
 {$G$: residual edges}(z.west);
\end{tikzpicture}
\caption{The construction in Proposition~\ref{prop:residual-problem}.
The boxes partition the two sides of $F=T-K$. Lines indicate possible
classes of edges, not complete bipartite graphs. There are no edges from
$S_0$ to $W$. Only the edges between $S_0$ and $Z_0$ remain to be checked
by the residual labeling. Edges incident with $K$ are omitted.}
\label{fig:residual}
\end{figure}
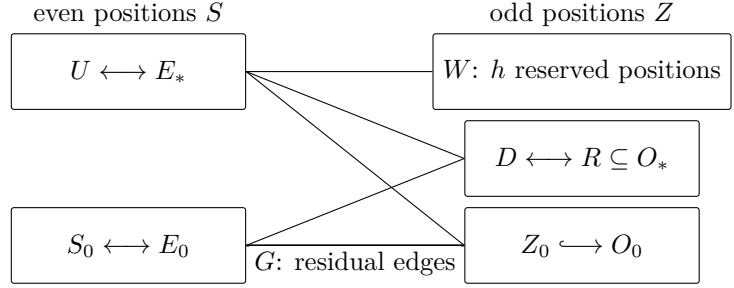

\begin{proof}
\emph{Step 1: Reserve the positions and the even labels.}
Choose $A\subseteq S$ consisting of $a$ vertices of largest degree in
$F$, breaking ties arbitrarily. Apply Lemmas~\ref{lem:centroid-order}
and~\ref{lem:cyclic-neighborhood} to $F-A$, whose parts have sizes
$m-a$ and $q$. Condition~\eqref{eq:reserve-budget} gives a set
$W\subseteq Z$ with $\abs{W}=h$ and
$\abs{N_{F-A}(W)}\leq w-a$. Thus $A\cup N_F(W)$ has at most
$w$ vertices. Enlarge it to a set $U\subseteq S$ of size $w$, and
assign the labels $E_*$ bijectively to $U$.

\emph{Step 2: Place some of the universal odd labels.}
Among the vertices of $Z\setminus W$, choose $d$ of largest degree
in $F$ and call the resulting set $D$. Choose any $d$ labels
$R\subseteq O_*$, and label $D$ bijectively by $R$. The imposed bounds
on $d$ make both choices possible. This construction is illustrated in
Figure~\ref{fig:residual}.

\emph{Step 3: Check the residual parameters.}
The cardinalities in \eqref{eq:residual-sizes} follow by subtraction.
For $s\in S_0$, the vertices of $A$ all have degree at least $d_F(s)$.
Since $s\notin A$ and each edge of $F$ has exactly one endpoint in $S$,
\[
 (a+1)d_F(s)\leq\sum_{v\in S}d_F(v)=e(F)\leq m+q-1.
\]
Likewise, for $z\in Z_0$, the $d$ vertices of $D$ all have degree at
least $d_F(z)$, and
\[
 (d+1)d_F(z)\leq\sum_{v\in Z}d_F(v)=e(F)\leq m+q-1.
\]
Restricting to $G$ cannot increase degrees. These inequalities prove
\eqref{eq:residual-degrees}; if a side is empty, its maximum degree
is zero by convention.

\emph{Step 4: Restore all labels and edges.}
Suppose $f_0,g_0$ satisfy \eqref{eq:residual-coprimality}. Exactly $h$
labels of $O_0$ remain unused; assign them bijectively to $W$. Together
with the labels on $K,U,D$, this uses every element of $[n]$ exactly once.

Edges incident with $K$ have coprime labels by
Corollary~\ref{cor:even-side}. In $F$, edges incident with $U$ have
relaxed-compatible labels by Lemma~\ref{lem:special-labels}(2), and
edges incident with $D$ have relaxed-compatible labels by part (3)
of that lemma. Every edge incident with $W$ has its other endpoint in
$U$. All remaining edges lie in $G$ and satisfy
\eqref{eq:residual-coprimality}. Thus every edge of $T$ satisfies the
relaxed relation. Corollary~\ref{cor:prime-pair-repair} repairs the
labeling, proving the assertion.
\end{proof}

\begin{corollary}\label{cor:residual-hall}
In Proposition~\ref{prop:residual-problem}, fix a bijection
$f_0\colon S_0\to E_0$. Form the bipartite graph $H_0$ with parts
$Z_0,O_0$, where $z\lambda$ is an edge if
\[
 \gcd(f_0(s),\lambda)=1
 \quad\text{for every }s\in N_G(z).
\]
There is an injection $g_0$ satisfying
\eqref{eq:residual-coprimality} if and only if
\begin{equation}\label{eq:residual-hall}
 \abs{N_{H_0}(Q)}\geq\abs{Q}-h
 \qquad\text{for every }Q\subseteq O_0.
\end{equation}
Whenever this condition holds, $T$ has a prime labeling.
\end{corollary}

\begin{proof}
Add $h$ new vertices to the $Z_0$ side of $H_0$, each adjacent to all
labels in $O_0$. The resulting graph has equal side sizes by
\eqref{eq:residual-sizes}. An injection labeling all of $Z_0$ extends
to a perfect matching by assigning the $h$ unused labels to the new
vertices; conversely, a perfect matching restricts to such an injection.
For a nonempty set $Q$ of labels, its neighborhood in the enlarged
graph has size $h+\abs{N_{H_0}(Q)}$. Lemma~\ref{lem:hall} therefore
gives \eqref{eq:residual-hall}. The last assertion follows from
Proposition~\ref{prop:residual-problem}.
\end{proof}

This identifies the remaining task for the matching part of the proof:
choose the even-label bijection so that every set of odd labels satisfies
\eqref{eq:residual-hall}. The prime supply
\eqref{eq:prime-separator-budget}, the reserve budget
\eqref{eq:reserve-budget}, and these simultaneous matching inequalities
are separate obligations. In later sections they must be established
for the parameter ranges in which this reduction is used.

\section{Arithmetic types and Hall estimates}
\label{sec:part-arithmetic}

This part develops the arithmetic and matching tools used for all
orders beyond the finite certificates.
Sections~\ref{sec:arithmetic-types} and~\ref{sec:compatible-measures}
define the arithmetic types over a set of small primes, compare the type
counts of an interval with the product law, and bound the compatible
type laws. Section~\ref{sec:forest-types} assigns the even labels to
the residual forest with exact class sizes
(Proposition~\ref{prop:two-measure-assignment}).
Sections~\ref{sec:matching-estimates} and~\ref{sec:refined-hall}
give Hall criteria for the type-restricted matching
(Theorems~\ref{thm:simultaneous-matching}
and~\ref{thm:three-range-hall}).

\subsection{Arithmetic types and their distributions}\label{sec:arithmetic-types}

The matching graph is controlled in two stages. We first exclude common
small prime divisors by assigning types to vertices. Common larger prime
divisors are then counted as deletions from that initial compatibility
graph. This section supplies the arithmetic information for both stages.

\subsubsection{Exact counts and the product law}

Let $P$ be a finite set of odd primes. For a positive integer $x$, define
its \emph{$P$-type} by
\[
 J_P(x)=\{p\in P:p\mid x\}.
\]
For $I\subseteq P$, put $d_I=\prod_{p\in I}p$, with $d_\varnothing=1$.
The associated product probability measure on $2^P$ is
\begin{equation}\label{eq:prime-product-law}
 \mu(I)=\prod_{p\in I}\frac1p
              \prod_{p\in P\setminus I}\left(1-\frac1p\right).
\end{equation}
For an integer $M\geq0$, let $C_I(M)$ count the integers $x\in[M]$
of type $I$, and let $C_I^{\mathrm{odd}}(M)$ count those of type $I$
among $1,3,\ldots,2M-1$. Both counts are zero when $M=0$.

\begin{lemma}\label{lem:exact-type-counts}
For every $I\subseteq P$ and integer $M\geq0$,
\begin{align}
 C_I(M)&=\sum_{J\subseteq P\setminus I}(-1)^{\abs{J}}
                     \floor*{\frac{M}{d_{I\cup J}}},
                     \label{eq:exact-type-count}\\
 C_I^{\mathrm{odd}}(M)&=C_I(2M)-C_I(M).
                     \label{eq:odd-type-count}
\end{align}
If $\Pi=\prod_{p\in P}p$, each of the two sequences of types has
period $\Pi$, and the number of occurrences of type $I$ in a period
is $\Pi\mu(I)$.
\end{lemma}

\begin{proof}
The indicator of the event $J_P(x)=I$ is
\[
 \prod_{p\in I}\boldsymbol{1}_{p\mid x}
 \prod_{p\in P\setminus I}(1-\boldsymbol{1}_{p\mid x}).
\]
Expanding the second product and summing over $x\in[M]$ proves
\eqref{eq:exact-type-count}. Multiplication by two preserves the
$P$-type, since all primes in $P$ are odd. The even integers in
$[2M]$ therefore have the same type counts as $[M]$, proving
\eqref{eq:odd-type-count}.

For completeness, the residues modulo the distinct primes of $P$ vary
independently over a full period. Indeed the map from residues modulo
$\Pi$ to their tuple of residues modulo the primes is injective: if two
residues have the same tuple, their difference is divisible by each
prime and hence by $\Pi$. Domain and codomain both have $\Pi$
elements, so the map is a bijection. There is one zero residue and
$p-1$ nonzero residues modulo $p$, giving \eqref{eq:prime-product-law}.
The map $x\mapsto2x-1$ is also a permutation modulo $\Pi$. This proves
the same period count for the odd sequence.
\end{proof}

For probability measures $\eta,\mu$ on a finite set, write
\[
 \|\eta-\mu\|_{\mathrm{TV}}
       =\frac12\sum_I\abs{\eta(I)-\mu(I)}.
\]
Since the coordinate differences sum to zero, their positive and
negative sums agree. It follows that this distance is also
$\sup_{0\leq f\leq1}\abs{\E_\eta f-\E_\mu f}$; the supremum is
attained by the indicator of the positive coordinates. Let $\mu_M$
and $\mu_M^{\mathrm{odd}}$ denote the two empirical type laws for
$M\geq1$.

\begin{corollary}[Periodic discrepancy]\label{cor:periodic-discrepancy}
For $0\leq s<\Pi$, define
\[
 D(s)=\frac12\sum_I\abs{C_I(s)-s\mu(I)},\qquad
 D^{\mathrm{odd}}(s)=\frac12\sum_I
                  \abs{C_I^{\mathrm{odd}}(s)-s\mu(I)}.
\]
If $M\equiv s\pmod\Pi$, then
\begin{equation}\label{eq:periodic-tv}
 \|\mu_M-\mu\|_{\mathrm{TV}}=\frac{D(s)}M,
 \qquad
 \|\mu_M^{\mathrm{odd}}-\mu\|_{\mathrm{TV}}
                                  =\frac{D^{\mathrm{odd}}(s)}M.
\end{equation}
In particular, both distances are bounded by the corresponding maximum
discrepancy over one period divided by $M$.
\end{corollary}

\begin{proof}
Write $M=v\Pi+s$. Lemma~\ref{lem:exact-type-counts} gives
$C_I(M)=v\Pi\mu(I)+C_I(s)$, and the same identity holds for the odd
counts. Subtract $M\mu(I)$, sum the absolute values, and divide by $2M$.
\end{proof}

The period can be too large to enumerate. The next estimate controls
the entire distribution using only divisibility moments of bounded order.

% Source: MOMENT_TYPE_CONTROL.md, Section 1.
\begin{proposition}[Truncated moment bound]\label{prop:moment-tv}
Put $k=\abs{P}$, and let $\sigma$ be the sum of the reciprocals of
the primes in $P$. For an integer $r\geq0$, set $t=r+1$ and
\[
 B_r=\sum_{j=0}^{\min(r,k)}2^j\binom{k}{j},
 \qquad
 R_r=3^t\sum_{\substack{I\subseteq P\\\abs{I}=t}}\frac1{d_I}.
\]
Then, for every $M\geq1$,
\begin{equation}\label{eq:moment-tv}
 \|\mu_M-\mu\|_{\mathrm{TV}}\leq\frac{B_r}{M}+2R_r,
 \qquad
 \|\mu_M^{\mathrm{odd}}-\mu\|_{\mathrm{TV}}
                          \leq\frac{B_r}{M}+3R_r,
\end{equation}
and $R_r\leq(3\sigma)^t/t!$. The sum defining $R_r$ is zero if $t>k$.
\end{proposition}

\begin{proof}
Fix $f\colon2^P\to[0,1]$. Define its finite differences and their
truncated expansion by
\[
 \Delta_I f=\sum_{J\subseteq I}(-1)^{\abs{I}-\abs{J}}f(J),
 \qquad
 T_rf(S)=\sum_{\substack{I\subseteq S\\\abs{I}\leq r}}\Delta_I f.
\]
We have $\abs{\Delta_I f}\leq2^{\abs{I}}$. The full expansion equals
$f(S)$: the coefficient of $f(J)$ in it is
$\sum_{J\subseteq I\subseteq S}(-1)^{\abs{I}-\abs{J}}$, which is
one for $J=S$ and zero otherwise. Thus $T_rf(S)=f(S)$ when
$\abs{S}\leq r$.

If $s=\abs{S}>r$, then
\begin{equation}\label{eq:moment-remainder}
 \abs{f(S)-T_rf(S)}
 \leq1+\sum_{j=0}^r2^j\binom{s}{j}
 \leq3^t\binom{s}{t}.
\end{equation}
To see the last inequality, use
$\binom{s}{j}\leq\binom{s}{t}\binom{t}{j}$ for $j<t\leq s$.
The resulting sum is $(3^t-2^t)\binom{s}{t}$, and the missing term
$2^t\binom{s}{t}$ is at least one. With the convention
$\binom{s}{t}=0$ for $s<t$, \eqref{eq:moment-remainder} holds for all
$S$.

Under $\mu$, the divisibility moment $\Pr(I\subseteq S)$ is $1/d_I$.
Under $\mu_M$, it is $\floor{M/d_I}/M$, so its error is at most
$1/M$ and its value is at most $1/d_I$. Under the odd law, the same
event picks one residue class of indices modulo $d_I$. Its error is
again at most $1/M$. Its value is at most $2/d_I$, since the number of
odd multiples of $d_I$ is no greater than the number of all multiples
up to $2M-1$.

It follows that the error in the expectation of $T_rf$ is at most
$B_r/M$ for either empirical law. Also
\[
 \E\binom{\abs{S}}{t}
       =\sum_{\substack{I\subseteq P\\\abs{I}=t}}\Pr(I\subseteq S).
\]
Applying \eqref{eq:moment-remainder} to the empirical and product laws
therefore adds at most $2R_r$ in the ordinary case and $3R_r$ in the
odd case. Taking the supremum over $f$ proves \eqref{eq:moment-tv}.
Finally, the expansion of $\sigma^t$ contains each product with $t$
distinct primes exactly $t!$ times, proving the bound for $R_r$.
\end{proof}

\subsubsection{Grouping primes into blocks}

We may use fewer types by recording only whether a label is divisible
by some prime in each of several disjoint blocks. Fix a partition
\[
 P=B_1\sqcup\cdots\sqcup B_s
\]
into nonempty blocks, and define
\[
 T(x)=\{j\in[s]:p\mid x\text{ for some }p\in B_j\},
 \qquad q_j=1-\prod_{p\in B_j}\left(1-\frac1p\right).
\]
We call $T(x)$ the \emph{block type} of $x$. Write $\phi(J)$ for the
block type obtained from a full type $J\subseteq P$.

% Source: PRIME_BLOCK_GROUPING_LEMMA.md, Sections 1 and 3.
\begin{lemma}\label{lem:block-types}
The image of $\mu$ under $\phi$ is the product law
\begin{equation}\label{eq:block-law}
 \bar\mu(I)=\prod_{j\in I}q_j
                   \prod_{j\notin I}(1-q_j),\qquad I\subseteq[s].
\end{equation}
Its empty mass is
\begin{equation}\label{eq:empty-mass}
 \alpha:=\bar\mu(\varnothing)
        =\prod_{j=1}^s(1-q_j)=\prod_{p\in P}(1-1/p).
\end{equation}
For any two probability laws $\eta,\mu$ on the full types,
\[
 \|\phi_*\eta-\phi_*\mu\|_{\mathrm{TV}}
                       \leq\|\eta-\mu\|_{\mathrm{TV}}.
\]
Moreover, $T(x)=\varnothing$ exactly when $J_P(x)=\varnothing$,
and $T(x)\cap T(y)=\varnothing$ implies that $x,y$ share no prime in $P$.
\end{lemma}

\begin{proof}
Under \eqref{eq:prime-product-law}, divisibility indicators for distinct
primes are independent. The event defining coordinate $j$ uses only
primes in $B_j$, so these events are independent across blocks and have
probabilities $q_j$. This proves \eqref{eq:block-law} and
\eqref{eq:empty-mass}.

For each block type, its discrepancy is the sum of the discrepancies
over its fiber under $\phi$. Apply the triangle inequality within each
fiber and then sum over the fibers to obtain the total-variation bound.
The last two assertions follow directly from the definitions.
\end{proof}

Let $\bar C_I(M)=\sum_{\phi(J)=I}C_J(M)$ be the exact block-type
count, and define $\bar C_I^{\mathrm{odd}}$ similarly. All preceding
counting and distribution bounds pass to these counts. In particular,
the empirical block law has exact distance
\begin{equation}\label{eq:rational-tv}
 \frac{1}{2M\Pi}\sum_I
       \abs{\Pi\bar C_I(M)-M\Pi\bar\mu(I)}
\end{equation}
from $\bar\mu$. The numerator contains only integers, since
$\Pi\bar\mu(I)$ counts residues in a full period. Odd-prefix distances
have the same form with the odd counts. Thus maxima over finite ranges
can be compared by integer cross-multiplication.

Disjoint block types give a sufficient compatibility condition. The
converse need not hold: if $5$ and $7$ lie in the same block, their block
types intersect although the two integers are coprime. Both sides of the
matching graph must use the coarser relation, as in
Figure~\ref{fig:block-compatibility}.

\begin{figure}[tbp]
\centering
\begin{tikzpicture}[x=13mm,y=10mm,every node/.style={font=\small}]
\node at (1.5,1.9) {position type};
\node[rotate=90] at (-2.3,-1.5) {label type};
\foreach \i/\lab in {0/\varnothing,1/\{1\},2/\{2\},3/{\{1,2\}}} {
 \node[anchor=south,inner sep=2pt] at (\i,.85) {$\lab$};
 \node[anchor=east,inner sep=2pt] at (-.85,-\i) {$\lab$};
}
\draw[thin] (-.5,.5) rectangle (3.5,-3.5);
\foreach \i in {0.5,1.5,2.5} {
 \draw[gray!45] (\i,.5)--(\i,-3.5);
 \draw[gray!45] (-.5,-\i)--(3.5,-\i);
}
\foreach \x/\y in {0/0,1/0,2/0,3/0,0/1,2/1,0/2,1/2,0/3}
 \fill (\x,-\y) circle (1.5pt);
\end{tikzpicture}
\caption{Compatibility for the blocks $B_1=\{3\}$ and
$B_2=\{5,7\}$. A dot marks a pair of disjoint block types. The empty
type is compatible with every type; the full type is compatible only
with the empty type. The same relation is used for labels and positions.}
\label{fig:block-compatibility}
\end{figure}

\subsubsection{Removing special labels}

Deleting labels changes the reference law. The following exact identity
keeps this change separate from arithmetic discrepancy.

\begin{lemma}[Depletion of the empty type]\label{lem:empty-depletion}
Let a population of $M$ objects have type law $\eta$, and remove $w<M$
objects, all of empty type. Put $\xi=w/M$. If $\mu$ is a reference law
with empty mass $\alpha\geq\xi$, define
\[
 \eta'=\frac{\eta-\xi\delta_\varnothing}{1-\xi},
 \qquad
 \nu_\xi=\frac{\mu-\xi\delta_\varnothing}{1-\xi}.
\]
Here $\delta_\varnothing$ is the unit mass at the empty type. Then
\begin{equation}\label{eq:depleted-tv}
 \|\eta'-\nu_\xi\|_{\mathrm{TV}}
                 =\frac{\|\eta-\mu\|_{\mathrm{TV}}}{1-\xi}.
\end{equation}
\end{lemma}

\begin{proof}
The removed mass is $\xi\delta_\varnothing$ in each numerator.
Subtracting the two laws leaves $(\eta-\mu)/(1-\xi)$; now take the
sum of absolute values. Both displayed laws are probability measures
by the assumptions on the removed objects and on $\alpha$.
\end{proof}

\begin{corollary}[Exact residual quotas]\label{cor:residual-type-quotas}
Use the construction of Proposition~\ref{prop:residual-problem}. Let
$3\leq z\leq n/3$, and let $P$ in the present section be the set of
odd primes at most $z$, partitioned into blocks as above. To distinguish
it from the separator label set, denote that label set here by
$P_{\mathrm{sep}}$ and put $k_0=\abs{P_{\mathrm{sep}}}$.

Writing an even label as $2j$, the numbers of available even half-labels
and odd labels of block type $I$ are respectively
\begin{align}
 a_I&=\bar C_I(m)-w\boldsymbol{1}_{I=\varnothing},
                      \label{eq:exact-even-quota}\\
 o_I&=\bar C_I(n)-\bar C_I(m)
                   -(k_0+d)\boldsymbol{1}_{I=\varnothing},
                      \label{eq:exact-odd-quota}
\end{align}
where $m=\floor{n/2}$ and $d$ is the number of placed special odd labels.
In particular, $\sum_Ia_I=m-w$ and
$\sum_Io_I=\ceil{n/2}-k_0-d$.
\end{corollary}

\begin{proof}
A special even half-label is a power of two or a prime greater than
$n/3$. Neither has an odd prime divisor at most $z$, so each has empty
type. Subtracting them from $[m]$ proves \eqref{eq:exact-even-quota}.

Among $[n]$, the even labels have exactly the type counts of $[m]$.
Hence the initial odd counts are $\bar C_I(n)-\bar C_I(m)$. The
separator labels are primes greater than $n/2$, and the $d$ other
placed special odd labels belong to $O_*$: each is $1$ or a prime
greater than $n/3$. They are distinct and all have empty type. Their
removal proves \eqref{eq:exact-odd-quota} and the two total counts.
\end{proof}

For a positive residual even population, its empirical law is
$a_I/(m-w)$. Lemma~\ref{lem:empty-depletion} compares it to
$(\bar\mu-\xi\delta_\varnothing)/(1-\xi)$ when
$\xi=w/m\leq\alpha$. Alternatively, the exact quotas can be used
directly, avoiding a product-law approximation. An empty population
requires no such normalization.

\subsubsection{Counting the remaining prime divisors}

Continue to assume that $P$ consists of all odd primes at most $z$,
where $z\geq3$. For an odd label $y$, let
\[
 P_z(y)=\{p>z:p\text{ is prime and }p\mid y\},
 \qquad \tau_z(y)=\sum_{p\in P_z(y)}\frac1p.
\]

\begin{lemma}[Exact large-prime counts]\label{lem:large-prime-counts}
For $M\geq0$, the number $B_I(M,y)$ of integers $j\in[M]$ of block
type $I$ sharing a prime greater than $z$ with $y$ is
\begin{equation}\label{eq:large-prime-union}
 B_I(M,y)=\sum_{\varnothing\ne A\subseteq P_z(y)}
             (-1)^{\abs{A}+1}\bar C_I\!\left(\floor*{\frac{M}{d_A}}\right),
\end{equation}
where $d_A=\prod_{p\in A}p$. In particular
$\sum_I B_I(M,y)\leq M\tau_z(y)$. For $m=\floor{n/2}$,
define
\[
 T_n(j)=\sum_{\substack{p>z,\ p\text{ prime}\\p\mid j}}
                     \ceil*{\frac{\floor{n/p}}2}.
\]
This is an upper bound for the number of odd labels in $[n]$ sharing
a prime greater than $z$ with $2j$, and
\begin{equation}\label{eq:total-prime-cost}
 \sum_{j=1}^mT_n(j)
    =\sum_{\substack{z<p\leq m\\p\text{ prime}}}
           \floor*{\frac mp}\ceil*{\frac{\floor{n/p}}2}.
\end{equation}
\end{lemma}

\begin{proof}
Multiplication by $d_A$ preserves every small-prime block type. Thus
the number of integers of type $I$ divisible by all primes in $A$ is
$\bar C_I(\floor{M/d_A})$. Inclusion--exclusion proves
\eqref{eq:large-prime-union}. Counting the union by the sum of its
parts gives $\sum_I B_I(M,y)\leq\sum_{p\in P_z(y)}\floor{M/p}
\leq M\tau_z(y)$.

Since $p>z\geq3$ is odd, its odd multiples in $[n]$ correspond to
the odd integers at most $\floor{n/p}$, of which there are
$\ceil{\floor{n/p}/2}$. Summing over $p\mid j$ can only overcount
the conflicting labels. Summing again over $j\leq m$ and interchanging
the finite sums proves \eqref{eq:total-prime-cost}.
\end{proof}

The block type also limits how many large primes can divide a label.
This gives a useful bound even when the labels are not enumerated.

\begin{lemma}[A bound from the type]\label{lem:type-tail-bound}
Let $r_1<r_2<\cdots$ be the primes greater than $z$. For a block type
$I$, put $b_I=\prod_{j\in I}\min B_j$. If $b_I\leq n$, let
\[
 h_I=\max\{h\in\Z_{\geq0}:b_Ir_1\cdots r_h\leq n\},
\]
where the empty product is one. Every integer $y\leq n$ of type $I$
satisfies
\begin{equation}\label{eq:type-tail-bound}
 \abs{P_z(y)}\leq h_I,
 \qquad \tau_z(y)\leq\sum_{j=1}^{h_I}\frac1{r_j}.
\end{equation}
If $b_I>n$, there is no such integer $y$.
\end{lemma}

\begin{proof}
Choose one prime divisor of $y$ from each block indexed by $I$.
These primes are distinct and at most $z$, and their product is at
least $b_I$. If the distinct larger prime divisors of $y$ are
$p_1<\cdots<p_h$, then $p_j\geq r_j$ and
\[
 b_Ir_1\cdots r_h\leq
       \left(\prod_{j\in I}\text{chosen prime in }B_j\right)
          p_1\cdots p_h\leq y\leq n.
\]
This proves the bound on $h$ and excludes $b_I>n$. Summing
$1/p_j\leq1/r_j$ proves the second inequality.
\end{proof}

All bounds in this section concern label populations alone. Their use in
the matching problem also requires control of which types are assigned to
the forest vertices; that geometric step is separate.

\subsection{Compatible type laws and rectangle bounds}
\label{sec:compatible-measures}

The arithmetic estimates compare the label populations with a product
law. We now show how that law supplies compatible pairs. Two distinct
tools will be useful: a coupling supported on disjoint types, and an
exact bound using only the compatibility degrees. Both concern weighted
type populations and are independent of the forest construction.

\subsubsection{A coupling supported on disjoint types}

In this subsection let $\mu$ denote a product law on $2^{[s]}$ with
coordinate probabilities $q_1,\ldots,q_s$, where $s\geq1$ and
$0<q_j<1/2$. This is the block law of \eqref{eq:block-law}, with the
bar suppressed. Define
\begin{equation}\label{eq:coupling-parameters}
 \alpha=\prod_{j=1}^s(1-q_j),\qquad
 \gamma=\prod_{j=1}^s\frac{1-2q_j}{1-q_j},\qquad
 \rho=\max_j\frac{q_j}{1-q_j}.
\end{equation}
For $0\leq\xi<\gamma$, put
\begin{equation}\label{eq:coupling-depletion}
 \nu_\xi=\frac{\mu-\xi\delta_\varnothing}{1-\xi},\qquad
 b_\xi=\frac{\alpha-\xi}{1-\xi},\qquad
 \theta_\xi=\rho\sqrt{\frac{\alpha}{\alpha-\xi}}.
\end{equation}
These quantities are well defined: each factor of $\gamma$ is positive
and at most the corresponding factor of $\alpha$, so
$0<\gamma\leq\alpha<1$.

% Source: PRIME_BLOCK_GROUPING_LEMMA.md, Section 2.
\begin{lemma}[Depleted product coupling]\label{lem:depleted-coupling}
There is a probability measure $\kappa_\xi$ on pairs of types with
marginals $\mu,\nu_\xi$, supported on $I\cap J=\varnothing$, such
that
\begin{equation}\label{eq:coupling-density}
 \kappa_\xi(I,J)\leq\frac{\mu(I)\nu_\xi(J)}{b_\xi}.
\end{equation}
For all real functions $F,G$ on the type space,
\begin{equation}\label{eq:coupling-correlation}
 \abs{\operatorname{Cov}_{\kappa_\xi}(F(I),G(J))}
 \leq\theta_\xi
       \sqrt{\operatorname{Var}_\mu(F)\operatorname{Var}_{\nu_\xi}(G)}.
\end{equation}
\end{lemma}

\begin{proof}
We first construct the coupling without depletion. In coordinate $j$,
give the pairs $(0,0),(1,0),(0,1),(1,1)$ respective probabilities
$1-2q_j,q_j,q_j,0$, and take the product over coordinates. Call the
result $\kappa_0$. Its marginals are both $\mu$, and its support
consists of disjoint types.

Relative to the independent product of the marginals, the nonzero
one-coordinate density ratios are
\[
 \frac{1-2q_j}{(1-q_j)^2},\qquad\frac1{1-q_j},\qquad\frac1{1-q_j}.
\]
All are at most $1/(1-q_j)$. Hence
$\kappa_0(I,J)\leq\mu(I)\mu(J)/\alpha$. Also
\begin{equation}\label{eq:empty-column}
 \frac{\kappa_0(I,\varnothing)}{\mu(I)}
      =\prod_{j\notin I}\frac{1-2q_j}{1-q_j}\in[\gamma,1].
\end{equation}
We may therefore define a nonnegative measure
\begin{equation}\label{eq:depleted-coupling-definition}
 \kappa_\xi=
      \frac{\kappa_0-\xi\mu\otimes\delta_\varnothing}{1-\xi}.
\end{equation}
Its total mass is one and its marginals are $\mu,\nu_\xi$.
For $J\ne\varnothing$, the preceding density bound gives
\[
 \kappa_\xi(I,J)\leq\frac{\mu(I)\nu_\xi(J)}\alpha
                         \leq\frac{\mu(I)\nu_\xi(J)}{b_\xi}.
\]
For the empty column, \eqref{eq:empty-column} instead gives
$\kappa_\xi(I,\varnothing)\leq\mu(I)$, which is exactly
\eqref{eq:coupling-density}, since $\nu_\xi(\varnothing)=b_\xi$.

To bound covariance, first work under $\kappa_0$. The functions
\[
 \chi_A(I)=\prod_{j\in A}
       \frac{\boldsymbol{1}_{j\in I}-q_j}{\sqrt{q_j(1-q_j)}},
       \qquad A\subseteq[s],
\]
form an orthonormal basis for functions under $\mu$. Independence of
the coordinate pairs shows that
\[
 \E_{\kappa_0}[\chi_A(I)\chi_B(J)]
 =\begin{cases}
   \displaystyle\prod_{j\in A}\frac{-q_j}{1-q_j},&A=B,\\
   0,&A\ne B.
  \end{cases}
\]
For nonempty $A$, the absolute value of the product is at most $\rho$,
because every factor has absolute value less than one. Expanding the
centered functions $F,G$ in this basis and applying Cauchy--Schwarz
therefore gives
\begin{equation}\label{eq:undepleted-covariance}
 \abs{\operatorname{Cov}_{\kappa_0}(F,G)}
 \leq\rho\sqrt{\operatorname{Var}_\mu(F)\operatorname{Var}_\mu(G)}.
\end{equation}
Equation~\eqref{eq:depleted-coupling-definition} implies, by subtracting
the product of the new means, that
\[
 \operatorname{Cov}_{\kappa_\xi}(F,G)
            =\frac{\operatorname{Cov}_{\kappa_0}(F,G)}{1-\xi}.
\]

It remains to compare the two variances of $G$. Put
$g_0=G(\varnothing)$, and let $m_1,v_1$ be the mean and variance of
$G$ under $\mu$ conditioned on a nonempty type. Then
\begin{align*}
 \operatorname{Var}_\mu(G)
   &=(1-\alpha)v_1+\alpha(1-\alpha)(g_0-m_1)^2,\\
 \operatorname{Var}_{\nu_\xi}(G)
   &=\frac{1-\alpha}{1-\xi}v_1
     +\frac{(\alpha-\xi)(1-\alpha)}{(1-\xi)^2}(g_0-m_1)^2.
\end{align*}
Since $\alpha(1-\xi)\geq\alpha-\xi$, comparison of the two
coefficients gives
\[
 \operatorname{Var}_{\nu_\xi}(G)
 \geq\frac{\alpha-\xi}{\alpha(1-\xi)^2}
                    \operatorname{Var}_\mu(G).
\]
Combine this with \eqref{eq:undepleted-covariance} to obtain
\eqref{eq:coupling-correlation}.
\end{proof}

A \emph{submeasure} of a probability vector $\mu$ is a vector $f$ with
$0\leq f_I\leq\mu(I)$; its mass is $\sum_I f_I$. For two submeasures
on the type spaces, write
\[
 \mathcal C(f,g)=\sum_{I\cap J=\varnothing}f_Ig_J
\]
for their compatible pair mass. This counts pairs using the independent
product of the selected populations; the coupling is a device for
bounding that quantity.

\begin{corollary}[Correlation bound for rectangles]
\label{cor:correlation-rectangle}
If $f\leq\mu$ and $g\leq\nu_\xi$ have masses $x,y$, then
\begin{equation}\label{eq:correlation-rectangle}
 \mathcal C(f,g)\geq b_\xi
   \left(xy-\theta_\xi\sqrt{x(1-x)y(1-y)}\right)_+,
\end{equation}
where $a_+=\max\{a,0\}$. If $x+y\geq1$ and $\theta_\xi<1$, this
implies $\mathcal C(f,g)\geq b_\xi(1-\theta_\xi)xy$.
\end{corollary}

\begin{proof}
Put $F(I)=f_I/\mu(I)$ and $G(J)=g_J/\nu_\xi(J)$, so that
$0\leq F,G\leq1$ and their respective means are $x,y$.
Since $F^2\leq F$ and $G^2\leq G$, their variances are at most
$x(1-x)$ and $y(1-y)$. The support and density statements of
Lemma~\ref{lem:depleted-coupling} give
\[
 \mathcal C(f,g)\geq b_\xi\E_{\kappa_\xi}[FG].
\]
Apply its covariance bound and nonnegativity to obtain
\eqref{eq:correlation-rectangle}. Finally $x+y\geq1$ implies
$(1-x)(1-y)\leq xy$, which proves the last assertion.
\end{proof}

The compatibility degrees of these two laws have simple exact forms.
For a type $I$, write
\[
 c(I)=\prod_{j\in I}(1-q_j).
\]
Independence gives
\begin{equation}\label{eq:two-model-degrees}
 \sum_{J\cap I=\varnothing}\mu(J)=c(I),\qquad
 \sum_{J\cap I=\varnothing}\nu_\xi(J)
                                  =\frac{c(I)-\xi}{1-\xi}.
\end{equation}
The two expressions generally differ. Grouping primes changes $q_j$
and hence both these degree profiles and the factor $\rho$ in
\eqref{eq:coupling-parameters}; these quantities must be recomputed.

\subsubsection{Using an upper bound on depletion}

A numerical interval may use a common upper bound $\omega$ on the
actual removed proportion $\xi$. Replacing $\xi$ by $\omega$ is
justified by monotonicity of compatibility, without adding a discrepancy
term for this replacement.

\begin{lemma}\label{lem:depletion-transport}
Let $0\leq\xi\leq\omega<\alpha$. For every submeasure
$g\leq\nu_\xi$ there is a submeasure $\widetilde g\leq\nu_\omega$
of the same mass such that, for every nonnegative vector $f$,
\begin{equation}\label{eq:transport-compatible}
 \mathcal C(f,\widetilde g)\leq\mathcal C(f,g).
\end{equation}
\end{lemma}

\begin{proof}
Put $t=(\omega-\xi)/(1-\xi)$. Direct substitution gives
\[
 \nu_\xi=(1-t)\nu_\omega+t\delta_\varnothing.
\]
Thus a type of law $\nu_\omega$, independently replaced by the empty
type with probability $t$, has law $\nu_\xi$. Let $K(I,J)$ be this
transition probability. It vanishes unless $J=I$ or $J=\varnothing$.
Put $h_J=g_J/\nu_\xi(J)$, which belongs to $[0,1]$, and define
\[
 \widetilde g_I=\nu_\omega(I)\sum_J K(I,J)h_J.
\]
This is bounded by $\nu_\omega$, and summing over $I$ gives
$\sum_I\widetilde g_I=\sum_Jg_J$.

Whenever $K(I,J)>0$, disjointness from $I$ implies disjointness from
$J$, since $J\subseteq I$. Multiply this indicator inequality by the
nonnegative quantities $f_A\nu_\omega(I)K(I,J)h_J$ and sum over
$A,I,J$. The resulting inequality is \eqref{eq:transport-compatible}.
\end{proof}

This lemma transfers a lower bound proved for the more depleted model
to every less depleted model. If a correlation bound is used, its
additional hypothesis $\omega<\gamma$ must still be checked.

\subsubsection{Bounds from the compatibility degrees}

We next give a bound that does not require a product structure. Let
$\mathcal A,\mathcal B$ be finite type sets, let
$\mathcal R\subseteq\mathcal A\times\mathcal B$ be a compatibility
relation, and let $\mu,\nu$ be probability vectors on the two sets.
In this subsection write
\[
 \mathcal C(f,g)=\sum_{(i,j)\in\mathcal R}f_ig_j,
 \qquad
 a_i=\sum_{j:(i,j)\in\mathcal R}\nu_j,
 \qquad
 b_j=\sum_{i:(i,j)\in\mathcal R}\mu_i.
\]
Fix $0\leq V<\tau\leq1/2$. For $0\leq s\leq1/2-V$, define
\begin{align}
 F_{\mathrm L}(s)&=
    \min_{\substack{0\leq f_i\leq\mu_i\\\sum_i f_i=s}}
                \sum_i f_i(a_i-s-2V)_+,\label{eq:left-profile}\\
 F_{\mathrm R}(s)&=
    \min_{\substack{0\leq g_j\leq\nu_j\\\sum_j g_j=s}}
                \sum_j g_j(b_j-s-2V)_+.
                \label{eq:right-profile}
\end{align}
The feasible sets are nonempty and compact, so these minima exist. Put
\begin{equation}\label{eq:profile-budget}
 W(V,\tau)=\min\left\{
     \min_{\tau-V\leq s\leq1/2-V}F_{\mathrm L}(s),\quad
     \min_{\tau-V\leq s\leq1/2-V}F_{\mathrm R}(s)\right\}.
\end{equation}

% Source: PROFILE_RECTANGLE_LEMMA.md.
\begin{proposition}[Degree-profile rectangle bound]
\label{prop:profile-rectangle}
If $f\leq\mu$ and $g\leq\nu$ have masses $x,y$ with
$x+y=1-2V$ and $x,y\geq\tau-V$, then
\[
 \mathcal C(f,g)\geq W(V,\tau).
\]
More generally, suppose probability laws $\widehat\mu,\widehat\nu$
satisfy
\[
 \|\widehat\mu-\mu\|_{\mathrm{TV}}\leq V,
 \qquad \|\widehat\nu-\nu\|_{\mathrm{TV}}\leq V.
\]
If $\widehat f\leq\widehat\mu$ and
$\widehat g\leq\widehat\nu$ have masses $u,v\geq\tau$ and
$u+v>1$, then $\mathcal C(\widehat f,\widehat g)\geq W(V,\tau)$.
\end{proposition}

\begin{proof}
For a left type $i$, at most $1-y$ of its compatible $\nu$-mass was
omitted from $g$. Hence
\[
 \sum_{j:(i,j)\in\mathcal R}g_j
               \geq(a_i-1+y)_+=(a_i-x-2V)_+.
\]
If $x\leq1/2-V$, multiply by $f_i$ and sum to obtain
$\mathcal C(f,g)\geq F_{\mathrm L}(x)\geq W(V,\tau)$.
If $x>1/2-V$, then $y<1/2-V$ and the symmetric argument with
$F_{\mathrm R}(y)$ proves the same result.

For the second assertion, first decrease the two masses to $u',v'$
with $u'+v'=1$, $u',v'\geq\tau$, $u'\leq u$, and $v'\leq v$.
Such a choice exists because the interval
\[
 [\max\{\tau,1-v\},\,\min\{u,1-\tau\}]
\]
is nonempty; choose $u'$ there and put $v'=1-u'$.
Scale down each original submeasure to its chosen mass. Next truncate
the coordinates of the left submeasure at $\mu_i$, and those on the
right at $\nu_j$. The mass removed on the left is at most
\[
 \sum_i(\widehat\mu_i-\mu_i)_+
                         =\|\widehat\mu-\mu\|_{\mathrm{TV}}\leq V,
\]
and similarly on the right. The resulting submeasures have masses at
least $u'-V$ and $v'-V$; scale them down to exactly these masses.
Every operation only decreases compatible pair mass. The first
assertion now applies with $x=u'-V$ and $y=v'-V$.
\end{proof}

\begin{lemma}[Finite evaluation of the profile bound]
\label{lem:profile-breakpoints}
Order the positive-mass left types so that $a_1\leq\cdots\leq a_r$,
and put $P_i=\sum_{j<i}\mu_j$ for $1\leq i\leq r+1$.
The minimum of $F_{\mathrm L}$ in \eqref{eq:profile-budget} is attained
at one of the points
\begin{equation}\label{eq:profile-breakpoints}
 \tau-V,\quad 1/2-V,\quad P_i,\quad a_i-2V
\end{equation}
that lie in $[\tau-V,1/2-V]$, with $1\leq i\leq r+1$ for $P_i$
and $1\leq i\leq r$ for $a_i-2V$. The analogous statement holds on the
right. Consequently $W(V,\tau)$ can be computed exactly by finitely
many rational operations whenever the input data are rational.
\end{lemma}

\begin{proof}
For fixed $s$, the coefficients $(a_i-s-2V)_+$ are nondecreasing in
$i$. If a feasible vector assigns positive mass to a later type while
an earlier type is not full, moving the smaller of the available
amounts to the earlier type cannot increase the objective. Repeating
this operation produces a minimizer that fills the types in order,
with at most one partially filled type. Thus, when
$P_i\leq s\leq P_{i+1}$,
\[
 F_{\mathrm L}(s)
    =\sum_{j<i}\mu_j(a_j-s-2V)_+
                    +(s-P_i)(a_i-s-2V)_+.
\]
Between consecutive points in \eqref{eq:profile-breakpoints}, both
the partial index and the signs inside the positive parts are fixed.
Every full-type term is affine or zero. The partial term is either
zero or a quadratic with coefficient $-1$ on $s^2$. Therefore
$F_{\mathrm L}$ is concave on each such interval, and its minimum
there occurs at an endpoint. The formula is continuous across the
endpoints. This proves the claim, including ties among the $a_i$ and
the case when the prescribed interval consists of one point.
\end{proof}

To apply this to the block model, take the label law to be $\mu$ and
the position law to be $\nu_\omega$, with $\omega<\alpha$. The
degrees are then $a_I=(c(I)-\omega)/(1-\omega)$ and $b_I=c(I)$ by
\eqref{eq:two-model-degrees}. Lemma~\ref{lem:depletion-transport}
justifies a common depletion upper bound. If the coupling estimate is
also applicable, both lower bounds are available and the larger may
be used.

For a matching graph with $\ell$ labels and $\ell$ positions, normalize
type counts by $\ell$. Then $\ell^2\mathcal C(f,g)$ is the number
of compatible pairs in the selected rectangle before any large-prime
deletions. In particular, a rectangle covered by
Proposition~\ref{prop:profile-rectangle} cannot become empty if the
number of deleted pairs is strictly less than $\ell^2W(V,\tau)$.
This addresses the range in which both selected sides have mass at
least $\tau$. Smaller sides require the separate degree estimates
in the later matching argument.

\subsection{Assigning types to a forest}\label{sec:forest-types}

The arithmetic counts prescribe how many even half-labels belong to each
type. We must place these types on the even side of the residual forest
and control which odd labels are allowed at the opposite vertices.
We give a construction in which every separating vertex lies on the
opposite side. This preserves the even quotas exactly during separation.
Two additive quantities suffice: the number of even positions and the
sum of their degrees.

\subsubsection{Completion and fractional assignments}

It is sometimes useful to add edges before constructing the assignment.
The following criterion records the cost of making the forest connected.

\begin{lemma}[Bipartite completion]\label{lem:bipartite-completion}
Let $G=(A,B)$ be a bipartite forest with both parts nonempty, and let
$D\geq1$ be an integer. If $\Delta_G(A)\leq D$ and
\[
 D\abs{A}\geq\abs{A}+\abs{B}-1,
\]
then $G$ is contained in a tree $H$ on the same vertices with the same
bipartition and $\Delta_H(A)\leq D$.
\end{lemma}

\begin{proof}
Join each isolated vertex of $A$ to any vertex of $B$. These edges
create no cycle and respect the degree bound. Every resulting component
contains a vertex of $B$. If there are $c>1$ components, their total
number of edges is $\abs{A}+\abs{B}-c$. The unused degree capacity on
$A$ is therefore at least $c-1>0$. Choose an $A$-vertex with unused
capacity and join it to a $B$-vertex in a different component. This
reduces both the unused capacity and $c-1$ by one. Repetition joins
all components without a cycle or a violation of the degree bound.
\end{proof}

We will temporarily give a connected piece several colors with
fractional weights. A piece is \emph{pure} if one color has weight one,
and \emph{mixed} otherwise. The next elementary reduction bounds the
number of pieces that still need to be split.

\begin{lemma}[Fractional pruning]\label{lem:fractional-pruning}
Let $d\geq0$ and $g\geq1$ be integers, and let finitely many objects
have vectors $v_j\in\R^d$. Suppose
$w_{ji}\geq0$, $1\leq i\leq g$, and $\sum_iw_{ji}=1$ for each
object $j$. The weights can be changed so as to preserve every vector
sum $\sum_jw_{ji}v_j$ and leave at most $d(g-1)$ mixed objects.
Objects that are already pure may be kept fixed. If the data are
rational, the changed weights may be taken rational.
\end{lemma}

\begin{proof}
Freeze the pure objects and suppose that $f>d(g-1)$ mixed objects
remain. Each has at least two positive weights. Vary only these positive
entries, subject to preserving the $f$ row sums and the $d$ coordinates
of the vector sum for each of the first $g-1$ colors. There are more
than $f+d(g-1)$ variables and at most that many homogeneous equations,
so a nonzero direction exists. The vector sum of the last color is
also preserved, because all row sums are fixed.

A nonzero direction preserving row sums has a negative coordinate.
Move in that direction until the first positive weight becomes zero.
All weights remain nonnegative, and their row sums remain one. No
previously zero entry becomes positive. Freeze any newly pure objects
and repeat if more than $d(g-1)$ mixed objects remain. Each step
strictly reduces the number of positive entries, so the process ends.
For rational data, Gaussian elimination supplies a rational direction,
and the step to the first zero is rational as well.
\end{proof}

\subsubsection{Rounding the last assignments}

Once each mixed piece contains only one even position, its color can be
rounded. Exact class sizes alone do not control the sums of degrees;
we use the ordering of those degrees to obtain that control. The
construction uses the sorted slots and bipartite matching of Shmoys
and Tardos~\cite[Section~2]{ST93}. We give the exact-quota form and
the two-sided error estimate needed here.

\begin{lemma}[Sorted-slot rounding]\label{lem:sorted-slot-rounding}
Let $f\geq0$ and $g\geq1$ be integers. Let $d_1,\ldots,d_f$ lie
in a real interval $[a,b]$, and let $w_{ji}\geq0$, $1\leq j\leq f$,
$1\leq i\leq g$, satisfy $\sum_iw_{ji}=1$. Suppose every quota
$k_i=\sum_jw_{ji}$ is an integer. There is an integral assignment
$x_{ji}\in\{0,1\}$ with the same row sums and quotas for which,
writing $R=b-a$ and
\[
 \epsilon_i=\sum_{j=1}^fd_j(x_{ji}-w_{ji}),
\]
one has
\begin{equation}\label{eq:slot-rounding-error}
 \abs{\epsilon_i}\leq R,
 \qquad \frac12\sum_{i=1}^g\abs{\epsilon_i}
                       \leq R\floor{g/2}.
\end{equation}
For any one specified color $i_*$, the assignment can also satisfy
$\epsilon_{i_*}\leq0$. Zero quotas and $f=0$ are allowed.
\end{lemma}

\begin{proof}
For each color $i$, list its positive fractional pieces $w_{ji}$ in
nonincreasing order of $d_j$. Their total length is the integer $k_i$.
Divide the list into $k_i$ slots of length one, splitting pieces at
slot boundaries. Let $y_{js}$ be the amount of item $j$ in slot $s$.
There are $f$ slots; each item and each slot has total mass one.

Join an item to a slot when $y_{js}>0$. For any set of items, its
total outgoing mass is its cardinality, while each neighboring slot
receives at most one. Hall's criterion, Lemma~\ref{lem:hall}, gives a
perfect matching. Assign each item the color of its matched slot.
The class sizes are exactly $k_i$.

For the $h$th slot of a fixed color with $k_i>0$, let $U_h,L_h$ be
the largest and smallest weights $d_j$ in its positive support.
The sorted order implies $U_h\geq L_h\geq U_{h+1}$. Consequently
\[
 \sum_{h=1}^{k_i}(U_h-L_h)
   =U_1-L_{k_i}+\sum_{h=1}^{k_i-1}(U_{h+1}-L_h)\leq R.
\]
Both the fractional and the matched degree sums lie between
$\sum_hL_h$ and $\sum_hU_h$. Their difference therefore has absolute
value at most $R$. A zero quota has zero error. Moreover
$\sum_i\epsilon_i=0$, since every item is assigned once. If $p$
errors are positive, their sum is at most both $pR$ and $(g-p)R$.
This proves \eqref{eq:slot-rounding-error} for every such matching.

To obtain the last assertion, choose a perfect matching minimizing
the sum of $d_j$ assigned to slots of color $i_*$. Its cost is at
most the fractional cost. Here is a direct justification. Starting
with $y$, choose a perfect matching in its positive support and
subtract the least entry on that matching from all its edges. The
row and column sums decrease by the same amount, and at least one
positive entry vanishes. If the common sum is still positive, the
same mass argument gives another perfect matching. Repetition ends
after finitely many steps and expresses $y$ as a convex combination
of perfect matching matrices. A minimum-cost matching has cost at
most their weighted average. Thus $\epsilon_{i_*}\leq0$, while the
already proved bounds remain valid. The empty case is immediate.
\end{proof}

\subsubsection{Separating from one side}

We measure the size of a piece by its number of $A$-vertices. A piece
with just one such vertex is already suitable for the preceding
rounding lemma, regardless of the number of $B$-vertices it contains.

\begin{lemma}[A weighted one-sided cut]\label{lem:weighted-side-cut}
Let $U=(A_U,B_U)$ be a bipartite tree with $m=\abs{A_U}\geq2$
and $\Delta_U(A_U)\leq D$, where $D\geq2$. Put
$\rho=(D-1)/D$. There is a vertex $b\in B_U$ such that every
component $V$ of $U-b$ satisfies
\begin{equation}\label{eq:weighted-side-contraction}
 \abs{V\cap A_U}-1\leq\rho(m-1).
\end{equation}
\end{lemma}

\begin{proof}
Choose a vertex $c$ minimizing $\sum_{a\in A_U}\operatorname{dist}_U(c,a)$.
Each branch of $U-c$ contains at most $m/2$ vertices of $A_U$:
otherwise a move into a branch with $s>m/2$ such vertices changes
the sum by $m-2s<0$. If $c\in B_U$, take $b=c$. Every resulting component
then satisfies
$\abs{V\cap A_U}-1\leq m/2-1\leq\rho(m-1)$.

Suppose instead that $c\in A_U$. Its branches have total $A_U$-count
$m-1$ and there are at most $D$ of them. Let $s$ be the largest
branch count, and let $b$ be the neighbor of $c$ in such a branch.
Then $(m-1)/D\leq s\leq m/2$. The component of $U-b$ containing
$c$ has $A_U$-count $m-s$, so its count minus one is at most
$\rho(m-1)$. Every other component lies in the chosen branch and
has count minus one at most $s-1\leq m/2-1$. This proves
\eqref{eq:weighted-side-contraction} in both cases.
\end{proof}

For integers $M\geq1$, $d\geq0$, and $0<\rho<1$, define
\begin{equation}\label{eq:weighted-separator-budget}
 \mathcal S(M,d;\rho)=
   \sum_{j\geq0}\min\left\{d,
          \floor*{\frac{M}{\ceil{\rho^{-j}}+1}}\right\}.
\end{equation}
Only finitely many summands are nonzero. This expression counts
separations at different component sizes.

\begin{proposition}[Two-measure assignment]\label{prop:two-measure-assignment}
Let $H=(A,B)$ be a tree with both parts nonempty, $M=\abs{A}$,
$e=e(H)$, and $\Delta_H(A)\leq D$, where $D\geq2$ is an integer.
Prescribe nonnegative integers $a_1,\ldots,a_g$ with sum $M$, and put
$\pi_i=a_i/M$ and $\rho=(D-1)/D$. There are a coloring
$c\colon A\to[g]$ and a set $C\subseteq B$ such that
\begin{enumerate}
 \item exactly $a_i$ vertices have color $i$;
 \item every vertex of $B\setminus C$ has neighbors of a single color;
 \item $\abs{C}\leq\mathcal S(M,2(g-1);\rho)$;
 \item on writing
 \[
  E_i=\sum_{a:\,c(a)=i}d_H(a),\qquad h_i=E_i-e\pi_i,
 \]
 one has $\sum_i h_i=0$, $\abs{h_i}\leq D$, and
 \begin{equation}\label{eq:two-measure-degree-error}
  H_{\mathrm{deg}}:=\frac12\sum_i\abs{h_i}
                                      \leq D\floor{g/2}.
 \end{equation}
\end{enumerate}
For one specified color $i_*$, we may also require $h_{i_*}\leq0$.
\end{proposition}

\begin{proof}
\emph{Step 1: Preserve the two measures.}
Initially regard the whole tree as one object with color weights
$\pi_i$. For a connected object $U$, use the measure vector
\[
 v(U)=\left(\abs{U\cap A},\ \sum_{a\in U\cap A}d_H(a)\right).
\]
The degrees are always those of the fixed tree $H$, including edges
to vertices removed later. The total vector assigned to color $i$
is initially $(a_i,e\pi_i)$.

Apply Lemma~\ref{lem:fractional-pruning} and freeze every pure object.
While a mixed object contains at least two $A$-vertices, choose one
with largest $\abs{U\cap A}-1$ and cut it at a vertex of $B$ given
by Lemma~\ref{lem:weighted-side-cut}. Put that vertex in $C$, replace
the object by the components of its deletion, and give each child
the parent's color weights. The cut vertex contributes zero to both
measures, so all target vectors remain unchanged. Prune again, keeping
frozen objects fixed. After every pruning there are at most
$d=2(g-1)$ mixed objects.

Each cut removes a new vertex. The process therefore terminates.
Every remaining component contains an $A$-vertex: an undeleted
$B$-vertex retains all its $A$-neighbors. At termination each mixed
object has exactly one $A$-vertex. The objects, together with $C$,
partition the vertices. Every edge whose endpoints are outside $C$
lies within one object.

\emph{Step 2: Count the cuts.}
For a nonterminal object put $L(U)=\abs{U\cap A}-1$. Its bucket
is the unique integer $j\geq0$ for which
\[
 \rho^{-j}\leq L(U)<\rho^{-(j+1)}.
\]
Each child has $L$ at most $\rho L(U)$, so it lies below its parent's
bucket if it is nonterminal. When the first cut in bucket $j$ is
made, no mixed object remains in a higher bucket, by the choice of
a largest object. There are at most $d$ mixed objects in bucket $j$
at that moment, and no later cut can create another. The objects
cut in that bucket are disjoint, and each contains at least
$\ceil{\rho^{-j}}+1$ vertices of $A$. Their number is at most
both $d$ and $\floor{M/(\ceil{\rho^{-j}}+1)}$. Summing gives
the asserted bound on $\abs{C}$.

\emph{Step 3: Round and check the edges.}
Keep the colors of the frozen objects. Their contributions to each
$A$-quota are integers. The fractional quotas on the mixed terminal
objects are therefore nonnegative integers as well. Apply
Lemma~\ref{lem:sorted-slot-rounding} to their single $A$-vertices,
with weights $d_H(a)\in[0,D]$ and the requested preferred color.
The quotas remain exact. Since the degree targets before rounding
were $e\pi_i$, its error bounds give all assertions in part (4),
including the preferred sign.

A frozen component is monochromatic on $A$. A terminal mixed
component contains just one $A$-vertex, so rounding its color also
makes it monochromatic on $A$. Every $b\in B\setminus C$ and all
its neighbors lie in one such component. Thus part (2) survives
rounding. If $M=1$ or $g=1$, the same construction ends without
cuts and the preceding conclusions still apply.
\end{proof}

\subsubsection{Recovering the position tags}

The construction did not preserve $B$-counts. Connectedness allows us
to recover them from the degree sums. For vectors of counts with equal
total mass, use the unnormalized variation
\[
 \operatorname{tv}(x,y)=\frac12\sum_i\abs{x_i-y_i}.
\]
Dividing by their positive common mass gives the total variation
between the corresponding probability laws.

\begin{lemma}[Integral recovery of tags]\label{lem:integral-tag-recovery}
Let $H=(A,B)$ be a tree with both parts nonempty. Let $A$ be colored
with quotas $a_i$, put $M=\abs{A}$, $N=\abs{B}$, $\pi_i=a_i/M$,
and use the quantities $E_i,h_i,H_{\mathrm{deg}}$ above. Suppose
$C\subseteq B$ and every vertex of $B\setminus C$ has neighbors
of a single color. All vertices of $B$ can be given color tags,
agreeing with that forced color outside $C$, such that their count
vector $b$ satisfies
\begin{equation}\label{eq:integral-tag-error}
 \operatorname{tv}(b,N\pi)
       \leq H_{\mathrm{deg}}+1-\max_i\pi_i.
\end{equation}
\end{lemma}

\begin{proof}
Every $A$-degree is positive, so $x_i=E_i-a_i$ is a nonnegative
integer. Since $e(H)=M+N-1$, we have $\sum_i x_i=N-1$.
Choose $j$ maximizing $a_j$ and set $b_i=x_i+\boldsymbol1_{i=j}$.
These are nonnegative integers with total $N$.

Suppose $C\ne\varnothing$. Every component of $H-C$ meets $C$
by an edge. Its $A$-vertices have one color, because the path
between any two of them alternates through undeleted $B$-vertices.
For color $i$, let $f_i$ be the number of $B$-vertices in those
components, $t_i$ their number of components, and $e_{iC}$ their
number of edges to $C$. The tree edge identity gives
\[
 E_i=a_i+f_i-t_i+e_{iC}.
\]
Each component contributes at least one edge to $C$, and these edges
have their other endpoint in $A$. Thus $e_{iC}\geq t_i$ and
$f_i\leq E_i-a_i=x_i\leq b_i$. All forced tags fit below the
chosen counts, so the vertices of $C$ can fill the remaining deficits.

If $C=\varnothing$, connectedness implies that all $A$-vertices
have one color. The chosen $j$ is that color, and $b$ assigns it to
every vertex of $B$, as required.

Finally, in either case,
\[
 b_i-N\pi_i=h_i+\boldsymbol1_{i=j}-\pi_i.
\]
The half sum of the absolute values of the last two terms is
$1-\pi_j$. The triangle inequality proves
\eqref{eq:integral-tag-error}.
\end{proof}

At separating vertices the tag can instead record all neighboring
types. This retains the positions for matching and pays for them
through the tag distribution.

\begin{lemma}[Retaining the separating positions]\label{lem:retained-tags}
Let $G=(A,B)$ be a bipartite graph. Give each $a\in A$ a type
$I(a)\subseteq[s]$ and each $b\in B$ a tag $J(b)\subseteq[s]$.
Suppose $Z\subseteq B$ and, on each edge $ab$ with $b\notin Z$,
either $I(a)=\varnothing$ or $I(a)=J(b)$. Set
\[
 J'(b)=
 \begin{cases}
  J(b),&b\notin Z,\\
  \displaystyle\bigcup_{a\in N_G(b)}I(a),&b\in Z.
 \end{cases}
\]
Then $I(a)\subseteq J'(b)$ on every edge. If $\beta,\beta'$ are the count
vectors of the two tag assignments and $q$ is any reference vector
of the same total mass, then
\begin{equation}\label{eq:retagging-variation}
 \operatorname{tv}(\beta',q)\leq\operatorname{tv}(\beta,q)+\abs{Z}.
\end{equation}

Suppose these are block types for all odd primes at most $z$, as in
Section~\ref{sec:arithmetic-types}. Assign distinct even labels $2j(a)$
to $A$, with $T(j(a))=I(a)$, and let $Y$ be any set of odd labels.
Write $j\sim_z y$ when $j,y$ share a prime greater than $z$.
Initially allow $(b,y)$ whenever $T(y)\cap J'(b)=\varnothing$.
For $y\in Y$, put
\begin{equation}\label{eq:retained-local-deletion}
 K_y=\sum_{\substack{a\in A\\I(a)\cap T(y)=\varnothing}}
        d_G(a)\,\boldsymbol1_{j(a)\sim_z y}.
\end{equation}
The number of initially allowed positions at which $y$ fails actual
coprimality is at most $K_y$. More generally, if
\[
 T_I(j)=\abs{\{y\in Y:T(y)\cap I=\varnothing,\ j\sim_z y\}},
\]
the total number of pairs deleted for failed coprimality is at most
\begin{equation}\label{eq:retained-global-deletion}
 \sum_{a\in A}d_G(a)T_{I(a)}(j(a)).
\end{equation}
\end{lemma}

\begin{proof}
The inclusion of types follows from the assumption outside $Z$ and
from the definition on $Z$. Changing one tag subtracts one unit from
one coordinate of its count vector and adds one to another, so its
variation cost is at most one. Summing over $Z$ and using the
triangle inequality proves \eqref{eq:retagging-variation}.

If $(b,y)$ is initially allowed, then $T(y)$ is disjoint from the
type of every neighbor of $b$. Thus no odd prime at most $z$ divides
both $y$ and an adjacent label. Since $y$ is odd, failure of
coprimality must be witnessed by a neighbor $a$ for which $j(a)$
and $y$ share a prime greater than $z$. Such a witness also satisfies
$I(a)\cap T(y)=\varnothing$. For fixed $a$, there are at most
$d_G(a)$ positions to charge to it. Counting all witnesses can only
overcount the failed positions, proving
\eqref{eq:retained-local-deletion}. Summing over $y\in Y$ and
interchanging the finite sums gives \eqref{eq:retained-global-deletion}.
\end{proof}

This statement is deterministic. Permuting actual labels inside their
prescribed classes preserves its assumptions. In particular, retaining
the separating positions introduces no uncounted large-prime deletions.

\begin{corollary}[Type assignment with retained positions]
\label{cor:forest-type-assignment}
Let $G=(A,B)$ be a bipartite forest with $M=\abs{A}>0$ and
$N=\abs{B}>0$. Suppose $D\geq2$ is an integer,
$\Delta_G(A)\leq D$, and $DM\geq M+N-1$. Prescribe an integer
quota $a_I\geq0$ for every type $I\subseteq[s]$, with
$\sum_Ia_I=M$. Put $g=2^s$, $\pi_I=a_I/M$, and
$\rho=(D-1)/D$.

There are a completing tree $H$, a type assignment $I(a)$ with these
exact quotas, and tags $J'(b)$ at every vertex of $B$ such that
$I(a)\subseteq J'(b)$ on every edge of $G$. Their count vector $b'$
satisfies
\begin{equation}\label{eq:forest-tag-budget}
 \operatorname{tv}(b',N\pi)
 \leq D\floor{g/2}+1-\max_I\pi_I
                         +\mathcal S(M,2(g-1);\rho).
\end{equation}
The degree errors in $H$ satisfy
\eqref{eq:two-measure-degree-error} and $\abs{h_I}\leq D$;
one specified error may be required to be nonpositive.
\end{corollary}

\begin{proof}
Complete $G$ by Lemma~\ref{lem:bipartite-completion} and apply
Proposition~\ref{prop:two-measure-assignment} to $H$. Interpret its
colors as the $g$ subsets of $[s]$. Lemma~\ref{lem:integral-tag-recovery}
assigns tags at every $B$-vertex with error at most
$H_{\mathrm{deg}}+1-\max_I\pi_I$. Outside the resulting separator
$C$, a tag agrees with all neighboring types in $H$, hence also
in $G$. Apply Lemma~\ref{lem:retained-tags} to $G$ with $Z=C$.
Its variation cost is at most $\abs{C}$, giving
\eqref{eq:forest-tag-budget}. The degree assertions are those of
Proposition~\ref{prop:two-measure-assignment} on $H$.
\end{proof}

The completion and the original forest have different degree sums.
The distributional error in this corollary is measured using $H$,
where $e(H)=M+N-1$. The deterministic deletion bounds of
Lemma~\ref{lem:retained-tags} may use the original degrees in $G$.
Any later estimate that replaces these by completion degrees is an
upper bound and must retain the completion's larger edge count.
Sharper separator constructions and alternative treatments of rare
types require their own bounds. The one-sided construction proved
here supplies the underlying quota and tag arguments.

\subsection{Matching estimates}\label{sec:matching-estimates}

After the types have been placed, the remaining choice is a permutation
of the labels within each type. The initial matching graph is fixed by
the types and tags. Only conflicts involving larger primes depend on
these permutations. We first bound these deletions, then compare them
with the rectangles required by Hall's criterion.

Throughout this section, let $G=(A,B)$ be a bipartite graph and let $Y$
be a set of odd labels, where $\abs{B}=\abs{Y}=R>0$. Partition all odd
primes at most $z\geq3$ into blocks as in Section~\ref{sec:arithmetic-types}.
Each vertex $a\in A$ has a type $I(a)$, and each $b\in B$ has a tag
$J(b)$ satisfying
\begin{equation}\label{eq:matching-tag-inclusion}
 I(a)\subseteq J(b)\qquad\text{whenever }ab\in E(G).
\end{equation}
For each type $I$, let $A_I$ be the set of vertices of that type and
let $\Lambda_I$ be a set of available even half-labels of the same
type, with
\[
 \abs{A_I}=\abs{\Lambda_I}=a_I.
\]
The sets $\Lambda_I$ are disjoint and contain positive integers.
Choose a uniform bijection $j\colon A_I\to\Lambda_I$ for each
nonempty class, independently across classes, and give $a$ the even
label $2j(a)$. All types, tags, and vertex sets are fixed before this
choice. Classes of size zero are omitted from sums involving averages.
Write
\[
 d_a=d_G(a),\qquad E_I=\sum_{a\in A_I}d_a,
 \qquad e=e(G)=\sum_I E_I,
\]
and fix a real bound $D>0$ with $d_a\leq D$ for every $a$.

Let $\Gamma_0$ have parts $Y,B$, with $yb$ an edge when
$T(y)\cap J(b)=\varnothing$. Obtain $\Gamma$ by deleting the edges
for which $y$ is not coprime to every even label adjacent to $b$ in
$G$. Thus every edge of $\Gamma$ is an available assignment, although
some available assignments may have been excluded by the block types.
A perfect matching in $\Gamma$ is sufficient for the labeling problem.
Let
\[
 Z=\abs{E(\Gamma_0)\setminus E(\Gamma)}.
\]
When applying this setup to Corollary~\ref{cor:residual-hall}, include
the $h$ additional positions as isolated vertices on the $B$ side of
$G$, with empty tags. They are adjacent to every label in both
matching graphs. If $\beta_I$ counts the tags on the original
positions, the enlarged position law is
$(\beta_I+h\boldsymbol1_{I=\varnothing})/R$. Both type laws and
their variation bounds must refer to this enlarged side size $R$.

\subsubsection{Class permutations}

Sampling without replacement gives the moment bound needed for both
individual labels and the total number of deletions. We prove it for
arbitrary nonnegative costs.

\begin{lemma}[A permutation moment bound]\label{lem:permutation-moment}
Let $x_1,\ldots,x_m\geq0$, where $m\geq1$, and let
$0\leq d_1,\ldots,d_m\leq D$, where $D>0$.
Put $E=\sum_jd_j$. If $\pi$ is a uniform permutation of $[m]$,
then, for every $t\geq0$,
\begin{equation}\label{eq:permutation-moment}
 \log\E\exp\left(\frac{t}{D}\sum_{j=1}^m d_jx_{\pi(j)}\right)
 \leq\frac ED\log\left(\frac1m\sum_{j=1}^m e^{t x_j}\right).
\end{equation}
\end{lemma}

\begin{proof}
For $u\in[0,1]^m$, let
\[
 F(u)=\log\E_\pi\exp\left(t\sum_j u_jx_{\pi(j)}\right).
\]
This function is symmetric and convex. Indeed its second derivative
along any line is the variance of the corresponding linear statistic
under the positive exponential weights, and is therefore nonnegative.
On the cube slice $\sum_j u_j=E/D$, a convex function attains a
maximum at an extreme point. Such a point has $q$ coordinates equal
to one, at most one coordinate $\theta\in[0,1)$, and all others
zero, where $q+\theta=E/D$. If two coordinates were strictly between
zero and one, they could be perturbed in opposite directions, so the
point would not be extreme.

First suppose $\theta=0$. Put $v_j=e^{t x_j}$. The moment at this
point is $e_q(v)/\binom mq$, where $e_q$ is the elementary symmetric
polynomial of degree $q$, with $e_0=1$. We claim that
\[
 \frac{e_q(v)}{\binom mq}\leq
                  \left(\frac1m\sum_jv_j\right)^q.
\]
For a fixed sum of the nonnegative $v_j$, choose a maximizer of
$e_q$ on the resulting compact simplex, and among those choose one
minimizing $\sum_jv_j^2$. Replacing two coordinates by their average
cannot decrease $e_q$: with the other coordinates fixed, its only
quadratic term in this pair is their product times a nonnegative
coefficient, and its linear terms depend only on their sum.
If the two coordinates differed, averaging would reduce the sum of
squares. Hence all coordinates of the chosen maximizer are equal,
which proves the claim.

If $0<\theta<1$, convexity in that coordinate bounds $F$ by
$(1-\theta)$ times its value with $q$ ones plus $\theta$ times its
value with $q+1$ ones. The preceding bound gives
$(q+\theta)\log(m^{-1}\sum_jv_j)$. The case $q=m$ has no
fractional coordinate. Applying this bound to $u_j=d_j/D$ proves
\eqref{eq:permutation-moment}.
\end{proof}

For $y\in Y$ and $a_I>0$, define
\[
 p_I(y)=\frac{\boldsymbol1_{I\cap T(y)=\varnothing}}{a_I}
       \abs{\{j\in\Lambda_I:j\sim_z y\}},
 \qquad \mu_y=\sum_I E_I p_I(y),
\]
where $j\sim_z y$ means that $j,y$ share a prime greater than $z$.
The count in $p_I(y)$ may be obtained from
Lemma~\ref{lem:large-prime-counts} after removing the unavailable labels.
Put
\[
 K_y=\sum_{a\in A}
    d_a\boldsymbol1_{I(a)\cap T(y)=\varnothing}
        \boldsymbol1_{j(a)\sim_z y}.
\]
By Lemma~\ref{lem:retained-tags} and
\eqref{eq:matching-tag-inclusion}, at most $K_y$ edges incident with
$y$ are deleted. Also $\E K_y=\mu_y$.

\begin{corollary}[Simultaneous local bounds]\label{cor:local-conflict-tails}
For every $t\geq0$,
\begin{equation}\label{eq:local-moment}
 \E e^{tK_y}\leq
       \exp\left(\frac{e^{tD}-1}{D}\mu_y\right).
\end{equation}
Fix $c_y>0$ and $0<\delta_y<1$, before choosing the permutations,
and put
\begin{equation}\label{eq:local-conflict-budget}
 L_y=c_y(e^{1/c_y}-1)\mu_y+c_yD\log(1/\delta_y).
\end{equation}
Then
\[
 \Pr\{K_y>L_y\}\leq\delta_y,
 \qquad
 \Pr\{K_y\leq L_y\text{ for every }y\in Y\}
                              \geq1-\sum_{y\in Y}\delta_y.
\]
If $\mu_y=0$, then $K_y=0$ surely; this label requires no failure
budget and we may instead set $L_y=\delta_y=0$.
\end{corollary}

\begin{proof}
Within class $I$, apply Lemma~\ref{lem:permutation-moment} to the
zero--one cost of a conflicting half-label, using parameter $tD$.
The logarithm of its moment is at most
\[
 \frac{E_I}{D}\log\bigl(1+p_I(y)(e^{tD}-1)\bigr)
 \leq\frac{E_I p_I(y)}D(e^{tD}-1),
\]
where $\log(1+x)\leq x$ for $x\geq0$ follows by integration of
$(1+x)^{-1}\leq1$. The class permutations are independent, so their
moments multiply. Summing the logarithms proves
\eqref{eq:local-moment}.

For $t>0$, Markov's inequality gives
$\Pr\{K_y>L\}\leq e^{-tL}\E e^{tK_y}$. Substituting
$t=1/(c_yD)$ and $L=L_y$ proves the first probability bound.
The union bound gives the simultaneous assertion; no independence
between the different $K_y$ is needed. Finally, a nonnegative random
variable of mean zero is zero surely.
\end{proof}

Thus a common choice $c_y=c$ and $\delta_y=\delta/R$ pays the
additive local term $cD\log(R/\delta)$. Any predetermined upper
bound for $\mu_y$ can replace its exact value in
\eqref{eq:local-conflict-budget}.

\subsubsection{Degree sums and arithmetic costs}

The preceding expectations depend on class degree sums. The next
identity records exactly what the assignment errors cost, including
the edges added in a connected completion.

\begin{lemma}[Comparison by degree variation]\label{lem:degree-cost-comparison}
Let $g\geq1$, $a_1,\ldots,a_g>0$, $M=\sum_i a_i$, and
$E_i\geq0$, with $e=\sum_iE_i$.
Put $h_i=E_i-ea_i/M$ and $H=\frac12\sum_i\abs{h_i}$.
For real numbers $f_i$,
\begin{equation}\label{eq:degree-cost-comparison}
 \sum_i E_i f_i
 \leq\frac eM\sum_i a_i f_i
                  +H(\max_i f_i-\min_i f_i).
\end{equation}
If $f_i\geq0$, $0\leq q_i\leq E_i$, and $q=\sum_iq_i$, then
\begin{equation}\label{eq:completion-cost-credit}
 \sum_i(E_i-q_i)f_i
 \leq\frac eM\sum_i a_i f_i
       +H(\max_i f_i-\min_i f_i)-q\min_i f_i.
\end{equation}
\end{lemma}

\begin{proof}
The numbers $h_i$ sum to zero, and their positive and negative sums
both have absolute value $H$. Therefore
\[
 \sum_i h_i f_i\leq H\max_i f_i-H\min_i f_i.
\]
Adding $(e/M)\sum_i a_i f_i$ proves the first assertion. Subtract
$\sum_iq_i f_i\geq q\min_i f_i$ to prove the second.
\end{proof}

In a completion $G\subseteq H$, take $E_i$ from $H$ and
$q_i=E_i^H-E_i^G$. Then $q=e(H)-e(G)$, since every added edge
has one endpoint on the even side. The first term in
\eqref{eq:completion-cost-credit} still uses $e(H)$. Replacing it
by $e(G)$ without the displayed correction would not be justified.
The lemma applies to the local frequencies $p_I(y)$, to class means
of conflict costs, and to their nonnegative log moments.

For $j\in\Lambda_I$, let
\[
 r_I(j)=\abs{\{y\in Y:T(y)\cap I=\varnothing,\ j\sim_z y\}}.
\]
Choose any deterministic costs $q_I(j)\geq r_I(j)$. The exact
union count $r_I$ is one choice. A sum over prime divisors is another;
the latter admits the following finite formula.

\begin{lemma}[Filtered prime costs]\label{lem:filtered-prime-costs}
Let $n\geq1$ be an integer. Suppose $Y$ consists of odd integers in
$[n]$ and all half-labels
lie in $[m]$, where $m=\floor{n/2}$. For $v\geq0$, put
\begin{align*}
 O_I(v)&=\sum_{J:\,J\cap I=\varnothing}
       \left(\bar C_J(\floor v)
                    -\bar C_J(\floor{v/2})\right),\\
 Q_{n,I}(j)&=\sum_{\substack{p>z,\ p\text{ prime}\\p\mid j}}
                               O_I(\floor{n/p}).
\end{align*}
Then $q_I(j)=Q_{n,I}(j)$ is an admissible cost, and its complete
class total satisfies
\begin{equation}\label{eq:filtered-prime-total}
 \sum_{\substack{j\leq m\\T(j)=I}}Q_{n,I}(j)
  =\sum_{\substack{z<p\leq m\\p\text{ prime}}}
      \bar C_I(\floor{m/p})\,O_I(\floor{n/p}).
\end{equation}
The costs for a fixed $j$ and the complete class totals are
nondecreasing in the integer $n$.
\end{lemma}

\begin{proof}
Multiplication by two preserves the block type. Thus $O_I(v)$
counts the odd positive integers at most $v$ whose types are disjoint
from $I$. For a prime $p>z$, multiplication by $p$ preserves both
the block type and oddness. The number of such odd labels in $[n]$
divisible by $p$ is therefore $O_I(\floor{n/p})$.
Every label counted by $r_I(j)$ shares with $j$ at least one prime greater than $z$;
summing over these primes can only overcount it. This proves
$Q_{n,I}(j)\geq r_I(j)$.

Interchange the sums over $j$ and its prime divisors. The number of
$j\leq m$ of type $I$ divisible by $p>z$ is
$\bar C_I(\floor{m/p})$. This gives
\eqref{eq:filtered-prime-total}. All counts increase when $n$
increases, proving the last assertion.
\end{proof}

The formula concerns the complete prefix. Means for the random
permutations must use the actual sets $\Lambda_I$. If $3\leq z\leq n/3$,
a removed special half-label has empty type: powers of two have cost
zero, while a prime $p$ with $n/3<p\leq n/2$ has cost one.
Indeed the only odd multiple of such a prime in $[n]$ is $p$ itself.
This permits exact subtraction for those removals. Other removals
must be accounted for separately.

\subsubsection{The total deletion bound}

For the chosen costs define
\[
 S=\sum_{a\in A}d_a q_{I(a)}(j(a)),\qquad
 f_I=\frac1{a_I}\sum_{j\in\Lambda_I}q_I(j),
\]
and
\[
 \psi_I(t)=\log\left(\frac1{a_I}
                         \sum_{j\in\Lambda_I}e^{tq_I(j)}\right).
\]
The deterministic charging argument of Lemma~\ref{lem:retained-tags}
gives $Z\leq S$. Notice that $S$ may count the same deleted pair
several times; this only strengthens its use as an upper bound.

\begin{proposition}[First and exponential moments]
\label{prop:global-conflict-tails}
In the preceding setup,
\begin{equation}\label{eq:global-conflict-moments}
 \E S=\sum_I E_I f_I,
 \qquad
 \log\E e^{tS/D}\leq\frac1D\sum_I E_I\psi_I(t)
                          \quad(t\geq0).
\end{equation}
For every $b>0$,
\begin{equation}\label{eq:global-conflict-tail}
 \Pr\{Z\geq b\}\leq p_{\mathrm g}(b):=
 \min\left\{1,\ \frac{\sum_I E_I f_I}{b},\
   \inf_{t>0}\exp\left(\frac{\sum_I E_I\psi_I(t)-tb}{D}\right)
                                                   \right\}.
\end{equation}
Any proved upper bounds for the two sums in this formula may be used
in their place.
\end{proposition}

\begin{proof}
Each half-label in a class is equally likely at each vertex of that
class, so linearity of expectation gives the first identity.
Apply Lemma~\ref{lem:permutation-moment} in each class and multiply
the resulting moments, using independence across classes, to obtain
the second inequality. Since $0\leq Z\leq S$, Markov's inequality
applied to $S$ and to $e^{tS/D}$ gives the two bounds in
\eqref{eq:global-conflict-tail}. Taking their minimum with one and
then the infimum over $t$ preserves the bound on this same event.
No independence between these different estimates is asserted.
\end{proof}

Lemma~\ref{lem:degree-cost-comparison} can be applied separately to
$f_I$ and $\psi_I(t)$. Both are nonnegative. In particular, a
connected completion gives the correction
\eqref{eq:completion-cost-credit} for each sum. For a numerical
certificate it is enough to choose one positive $t$ in advance and
bound its expression; the infimum need not be computed. The same
argument bounds deletions incident with a specified subset of $Y$
when the costs are defined using that subset.

\subsubsection{Excluding Hall obstructions}

The local estimates exclude small sides of an empty rectangle. The
remaining rectangles have a uniform positive supply of compatible
pairs, provided by Section~\ref{sec:compatible-measures}.

\begin{proposition}[Matching from a profile bound]
\label{prop:profile-hall}
Let $\Gamma_0$ be the graph of a type compatibility relation on
parts $Y,B$ of size $R>0$, and let $\Gamma$ be a spanning subgraph
of $\Gamma_0$.
Let $\widehat\mu,\widehat\nu$ be the empirical type laws on these
parts. Suppose reference probability laws $\mu,\nu$ satisfy
\[
 \|\widehat\mu-\mu\|_{\mathrm{TV}}\leq V,
 \qquad \|\widehat\nu-\nu\|_{\mathrm{TV}}\leq V,
 \qquad 0\leq V<\tau\leq1/2.
\]
Let $W(V,\tau)$ be the profile bound of
\eqref{eq:profile-budget} for this relation and these reference laws.
If every vertex of $\Gamma$ has degree at least $\tau R$ and
\begin{equation}\label{eq:profile-hall-deletion}
 \abs{E(\Gamma_0)\setminus E(\Gamma)}<R^2W(V,\tau),
\end{equation}
then $\Gamma$ has a perfect matching.
\end{proposition}

\begin{figure}[tbp]
\centering
\begin{tikzpicture}[x=1cm,y=1cm,font=\small]
 \fill[black!8] (2.35,1.5) rectangle (5.5,3.2);
 \draw (0,0) rectangle (5.5,3.2);
 \draw (0,1.5)--(5.5,1.5);
 \draw (2.35,0)--(2.35,3.2);
 \node[anchor=south] at (1.175,3.35) {$N_\Gamma(U)$};
 \node[anchor=south] at (3.925,3.35) {$X$};
 \node[anchor=east] at (-.15,2.35) {$U$};
 \node[anchor=east] at (-.15,.75) {$Y\setminus U$};
 \node[align=center] at (3.925,2.35) {no edges\\in $\Gamma$};
 \node[anchor=north] at (2.75,-.18) {$\abs U+\abs X>R$};
\end{tikzpicture}
\caption{Rows represent labels and columns represent positions.
For a Hall-deficient $U$, the set $X=B\setminus N_\Gamma(U)$
satisfies $\abs U+\abs X>R$. Every initially compatible pair in
the shaded rectangle has been deleted. The other three rectangles
carry no edge assertion.}
\label{fig:hall-rectangle}
\end{figure}
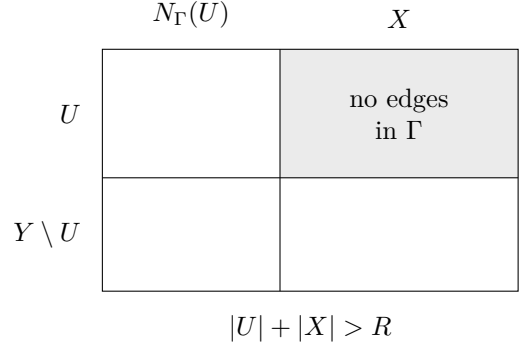

\begin{proof}
Suppose otherwise. By Lemma~\ref{lem:hall}, there is a nonempty
$U\subseteq Y$ with $\abs{N_\Gamma(U)}<\abs U$. Put
$X=B\setminus N_\Gamma(U)$. Then $U\times X$ has no edges of
$\Gamma$ and $\abs U+\abs X>R$.

If $\abs U\leq\tau R$, any vertex in $U$ has at least
$\tau R\geq\abs U$ neighbors, a contradiction. If
$\abs X\leq\tau R$, any vertex in $X$ has all its neighbors
outside $U$, and hence has degree at most
$R-\abs U<\abs X\leq\tau R$, again a contradiction.
Thus both $U$ and $X$ have size greater than $\tau R$.

Let $\widehat f,\widehat g$ be their type counts divided by $R$.
They are submeasures of $\widehat\mu,\widehat\nu$, their masses
sum to more than one, and both masses exceed $\tau$. By
Proposition~\ref{prop:profile-rectangle}, the number of
$\Gamma_0$-edges between $U$ and $X$ is at least $R^2W(V,\tau)$.
All of those edges were deleted, as illustrated in
Figure~\ref{fig:hall-rectangle}, contradicting
\eqref{eq:profile-hall-deletion}.
\end{proof}

We can now give a sufficient condition expressed entirely in terms of
the prepared graph, its class populations, and the preceding bounds.
For each $b\in B$, put
\begin{equation}\label{eq:position-deletion-budget}
 P_b=\sum_{a\in N_G(b)}\max_{j\in\Lambda_{I(a)}}q_{I(a)}(j).
\end{equation}
This is a deterministic upper bound for the deletions incident with
$b$. Indeed every initially allowed label at $b$ has type disjoint
from $I(a)$ for each neighbor $a$, so the same charging argument
applies with $b$ fixed. A sharper deterministic bound may replace
$P_b$ whenever available.

\begin{theorem}[A simultaneous matching criterion]
\label{thm:simultaneous-matching}
Use the random assignment and costs defined in this section, and
suppose the two empirical type laws satisfy the hypotheses of
Proposition~\ref{prop:profile-hall}, with $W(V,\tau)>0$.
Choose local budgets $L_y,\delta_y$ as in
Corollary~\ref{cor:local-conflict-tails}. Suppose
\begin{align}
 d_{\Gamma_0}(y)-L_y&\geq\tau R&&\text{for every }y\in Y,
                                      \label{eq:label-degree-condition}\\
 d_{\Gamma_0}(b)-P_b&\geq\tau R&&\text{for every }b\in B,
                                      \label{eq:position-degree-condition}\\
 \sum_{y\in Y}\delta_y+p_{\mathrm g}(R^2W(V,\tau))&<1.
                                      \label{eq:joint-success-condition}
\end{align}
Then some choice of the class bijections produces a graph $\Gamma$
with a perfect matching.
\end{theorem}

\begin{proof}
The event that some $K_y$ exceeds $L_y$ has probability at most
$\sum_y\delta_y$. The event $Z\geq R^2W(V,\tau)$ has probability
at most $p_{\mathrm g}(R^2W(V,\tau))$, by
Proposition~\ref{prop:global-conflict-tails}. Their union has
probability less than one by \eqref{eq:joint-success-condition}.
Hence there is an assignment for which all local budgets hold and
$Z<R^2W(V,\tau)$ simultaneously.

For that assignment, \eqref{eq:label-degree-condition} gives the
required label degrees. The position bounds are deterministic, so
\eqref{eq:position-degree-condition} gives the required degrees on
the other side as well. Proposition~\ref{prop:profile-hall} supplies
the perfect matching. The argument uses no independence between
the local and global failure events.
\end{proof}

The initial degrees in this criterion are finite type counts:
\begin{align*}
 d_{\Gamma_0}(y)&=\sum_{J:\,J\cap T(y)=\varnothing}
                    \abs{\{b:J(b)=J\}},\\
 d_{\Gamma_0}(b)&=\sum_{I:\,I\cap J(b)=\varnothing}
                    \abs{\{y:T(y)=I\}}.
\end{align*}
Consequently the criterion can use exact populations, or model lower
bounds with their stated variation errors. When its setup comes from
Proposition~\ref{prop:residual-problem}, the matching restricts to
the required injection on $Z_0$. That proposition then restores the
prime labeling of the original tree. The reserve and prime-supply
hypotheses of that proposition remain necessary.

\subsubsection{A deterministic bound for selected positions}

Bounding the total number of deleted pairs can be wasteful for a
rectangle with few positions. The forest structure gives a bound
depending on the mass of those positions. It holds for every
permutation, so it introduces no further failure event.

\begin{lemma}[Weighted forest bound]\label{lem:weighted-forest-cost}
Let $G=(A,B)$ be a bipartite forest with $\Delta_G(A)\leq D$,
where $D\geq1$ is an integer. Give each $b\in B$ a weight
$w_b\in[0,1]$, and put
\[
 v=\sum_{b\in B}w_b,\qquad
 d_w(a)=\sum_{b\in N_G(a)}w_b.
\]
Let $T_a\geq0$. Form the multiset with $D-1$ copies of each $T_a$,
and let $\mathcal L(t)$ be the sum of its largest $t$ entries, allowing
a fractional last entry. Set $\mathcal L(0)=0$ and keep it constant
after all entries have been used. Then
\begin{equation}\label{eq:weighted-forest-cost}
 \sum_{a\in A}d_w(a)T_a
 \leq\sum_{a\in A}T_a+\mathcal L\bigl((v-1)_+\bigr).
\end{equation}
\end{lemma}

\begin{proof}
We first prove
\begin{equation}\label{eq:weighted-forest-excess}
 \sum_{a\in A}(d_w(a)-1)_+\leq(v-1)_+.
\end{equation}
If $w$ is the indicator of a nonempty set $S\subseteq B$, the
forest induced by $S$ and its nonisolated neighbors has at most
$\abs S+\abs{N_G(S)}-1$ edges. Subtracting one for each neighbor
gives
\[
 \sum_a(d_w(a)-1)_+=e(G[S\cup N_G(S)])-\abs{N_G(S)}
                                                   \leq\abs S-1.
\]
For $S=\varnothing$ both sides of
\eqref{eq:weighted-forest-excess} are zero.

For general weights with fixed sum $v=q+\theta$, where
$q=\floor v$ and $0\leq\theta<1$, the left side of
\eqref{eq:weighted-forest-excess} is convex on the cube slice.
It has a maximum at a point with $q$ ones and at most one remaining
coordinate $\theta$. Convex interpolation between that point's two
adjacent integral vectors bounds its value by
\[
 (1-\theta)(q-1)_++\theta q=(q+\theta-1)_+.
\]
Here for $q=0$ the left expression is zero, and for $q\geq1$
it is $q+\theta-1$. If $q=\abs B$, the point is already integral.
This proves \eqref{eq:weighted-forest-excess}.

Write $e_a=(d_w(a)-1)_+$. Each $e_a$ is at most $D-1$, and
their sum is at most $(v-1)_+$. For nonnegative costs, the largest
possible value of $\sum_a e_aT_a$ under these constraints is
$\mathcal L((v-1)_+)$: move available excess to larger costs
until each earlier capacity is filled. Finally,
$d_w(a)\leq1+e_a$. Multiplying by $T_a$ and summing proves
\eqref{eq:weighted-forest-cost}. If $D=1$, the multiset is empty
and the same proof has $e_a=0$ for all $a$.
\end{proof}

\begin{corollary}[A uniform rectangle deletion bound]
\label{cor:forest-rectangle-deletions}
In the matching setup, suppose $G$ is a forest and $D\geq1$ is an
integer degree bound on $A$. Define
\[
 Q_0=\sum_I\sum_{j\in\Lambda_I}q_I(j),
\]
and let $\mathcal L$ be the largest-entry sum of the multiset with
$D-1$ copies of every cost $q_I(j)$ in these label pools. For any
weights $u_y,w_b\in[0,1]$, put $v=\sum_b w_b$. Then, for every
choice of the class bijections,
\begin{equation}\label{eq:forest-rectangle-deletions}
 \sum_{yb\in E(\Gamma_0)\setminus E(\Gamma)}u_yw_b
                  \leq Q_0+\mathcal L\bigl((v-1)_+\bigr).
\end{equation}
The bound holds simultaneously for all such weighted rectangles.
\end{corollary}

\begin{proof}
For each deleted pair, there is a neighbor $a$ of $b$ whose even
label shares a prime greater than $z$ with $y$. Moreover
$I(a)\cap T(y)=\varnothing$, by
\eqref{eq:matching-tag-inclusion}. Charge the pair to all such
neighbors and use $u_y\leq1$. For a fixed $a$, summing over the
eligible $y$ gives at most $r_{I(a)}(j(a))\leq q_{I(a)}(j(a))$.
The weighted deletion count is therefore at most
\[
 \sum_a d_w(a)q_{I(a)}(j(a)).
\]
Apply Lemma~\ref{lem:weighted-forest-cost} with
$T_a=q_{I(a)}(j(a))$. Since each class bijection uses every
half-label once, the multiset of these costs is fixed, independently
of the permutation. Its base sum is $Q_0$, and its repeated
largest-entry sum is the stated $\mathcal L$. This proves the
uniform assertion.
\end{proof}

This last bound can replace a global probability estimate for a
specified range of Hall rectangles. To do so, its dependence on $v$
must be compared with a compatible-pair lower bound at the same
rectangle masses. It is not a single bound on $Z$ at an arbitrarily
chosen value of $v$. The sharper treatment of exceptional label
classes and mass-restricted rectangles will use this distinction.

\subsection{Hall inequalities with type restrictions}\label{sec:refined-hall}

A deficient set in the matching graph need not have an arbitrary type
distribution. Degree bounds exclude some label types from small sets,
and a position type present in the complementary set imposes an
arithmetic restriction on every selected label. We use these restrictions
to replace a single rectangle bound by several smaller minimizations.

We retain the setting of Section~\ref{sec:matching-estimates}. In
particular, $\Gamma\subseteq\Gamma_0$ has parts $Y,B$ of size $R>0$,
and adjacency in $\Gamma_0$ means disjointness of the label type and
the position tag. Write $\mathcal T=2^{[s]}$, where $s\geq1$, and
abbreviate the empty and full types by $0=\varnothing$ and $F=[s]$.
Throughout this section,
\[
 C(I,J)=\boldsymbol1_{I\cap J=\varnothing}.
\]
A universal label means a vertex of $Y$ adjacent to every position in
$\Gamma$, after the actual coprimality deletions.

\subsubsection{Degree targets and deficient sets}

The least numerous types may be given weaker degree targets. What
matters is their population and the degree of each label, not the
number of events used to prove those degree bounds.

\begin{lemma}[Three ranges of deficient sets]\label{lem:three-hall-ranges}
Let
\[
 0\leq\varepsilon<\tau\leq\eta\leq\gamma,
 \qquad u\geq0,\qquad \gamma+u\leq1.
\]
Suppose $E\subseteq Y$ consists of full-type labels,
$\abs E\leq\varepsilon R$, and the following strict bounds hold:
\[
 d_\Gamma(y)>
 \begin{cases}
  \varepsilon R,&y\in E,\\
  \tau R,&T(y)=F,\ y\notin E,\\
  \eta R,&\varnothing\ne T(y)\ne F,\\
  \gamma R,&T(y)=\varnothing.
 \end{cases}
\]
Assume also that at least $uR$ labels are universal and have empty
type. If $U\subseteq Y$ is deficient, meaning
$\abs{N_\Gamma(U)}<\abs U$, and $V=B\setminus N_\Gamma(U)$, then
\begin{equation}\label{eq:deficient-range-masses}
 \abs U>\tau R,\qquad \abs V>uR,\qquad \abs U+\abs V>R.
\end{equation}
No universal label lies in $U$. Moreover, $U$ consists only of full-type
labels if $\abs U\leq\eta R$, and contains no empty-type label if
$\abs U\leq\gamma R$.
\end{lemma}

\begin{proof}
A deficient set is nonempty. If $U\subseteq E$, every member of $U$
has more than $\varepsilon R\geq\abs U$ neighbors, a contradiction.
Otherwise $U$ contains a label outside $E$, of degree greater than
$\tau R$. Its neighbors lie in $N_\Gamma(U)$, so $\abs U>\tau R$.

A universal member would give $N_\Gamma(U)=B$, also impossible.
Thus $\abs U\leq(1-u)R$, and
$\abs V=R-\abs{N_\Gamma(U)}>R-\abs U\geq uR$.
The same identity gives the last inequality in
\eqref{eq:deficient-range-masses}. Finally, a member of degree greater
than $\eta R$ cannot belong to a deficient set of size at most
$\eta R$. This excludes every type except $F$ in the first range.
The identical argument with $\gamma R$ excludes type $0$ in the
second range.
\end{proof}

Some full-type labels meet the stronger target without a probability
estimate. For an integer $h\geq1$ and a full-type odd label $y\leq h$,
put
\[
 b_0(h,y)=\abs{\{j\leq\floor{h/2}:T(j)=0,\ j\sim_z y\}}.
\]
Here, as before, $j\sim_z y$ means that $j,y$ share a prime greater
than the cutoff $z$.

\begin{lemma}[Protected full-type labels]\label{lem:protected-full-labels}
Suppose all labels lie in $[h]$, $\Delta_G(A)\leq D$, and
$R\geq R_0>0$. If $d_{\Gamma_0}(y)\geq cR$ for full-type labels,
where $c>\tau$, then
\begin{equation}\label{eq:protected-full-label}
 D b_0(h,y)<R_0(c-\tau)
\end{equation}
implies $d_\Gamma(y)>\tau R$ for every permutation within the even
type classes. Consequently the full-type labels failing
\eqref{eq:protected-full-label}, counted among all odd integers at
most $h$, give a uniform upper bound for the exception population.
\end{lemma}

\begin{proof}
A full-type label is initially compatible only with an empty-tag
position. Every neighbor of such a position has empty type, by
\eqref{eq:matching-tag-inclusion}. Each deletion is therefore witnessed
by an empty-type half-label counted in $b_0(h,y)$. Such a half-label
occupies one vertex of degree at most $D$, so it can witness at most
$D$ deleted positions. Hence
\[
 d_\Gamma(y)\geq cR-D b_0(h,y)
       >cR-R_0(c-\tau)\geq\tau R.
\]
The complete endpoint pool includes every usable half-label, and the
complete odd prefix includes every actual exception. This proves the
uniform assertion.
\end{proof}

For example, if the endpoint exception count is $e_0$, one may use
$\varepsilon=e_0/R_0$, provided $\varepsilon<\tau$. The weak degree
target must still be checked for each actual exception.

\subsubsection{Coordinate bounds that preserve mass}

Clipping a selected type distribution at a reference law loses mass.
Instead we move the excess to the full type, which is no more
compatible than any other type. This gives coordinate bounds without
changing the size of the selected set.

\begin{lemma}[Moving excess to the full type]\label{lem:full-type-transport}
Let $\rho,\pi$ be probability laws on $\mathcal T$ with
$\|\rho-\pi\|_{\mathrm{TV}}\leq v$, and let $0\leq f_I\leq\rho_I$.
There is a vector $\widehat f$ of the same total mass such that
\begin{equation}\label{eq:full-type-caps}
 0\leq\widehat f_I\leq\pi_I\quad(I\ne F),\qquad
 0\leq\widehat f_F\leq\pi_F+v.
\end{equation}
For every nonnegative vector $g$,
\[
 \sum_{I,J}C(I,J)\widehat f_I g_J
       \leq\sum_{I,J}C(I,J)f_Ig_J.
\]
The analogous assertion holds on the other side. If $f_0=0$, then
$\widehat f_0=0$. If $f_0\leq\rho_0-u$, then
\[
 \widehat f_0\leq\min\{\pi_0,\max(0,\pi_0+v-u)\}.
\]
For $I\ne F$, the retained coordinate can be realized by reducing
weights on the originally selected objects of type $I$.
\end{lemma}

\begin{proof}
Set $e_I=(f_I-\pi_I)_+$ and $e=\sum_I e_I$. Since $f\leq\rho$,
\[
 e\leq\sum_I(\rho_I-\pi_I)_+
       =\|\rho-\pi\|_{\mathrm{TV}}\leq v.
\]
Define $\widehat f_I=\min(f_I,\pi_I)$ and then add $e$ to its
full coordinate. This preserves the total mass and proves
\eqref{eq:full-type-caps}. The change moves mass only from a type
$I$ to $F$, and $C(F,J)\leq C(I,J)$ for every $J$. Thus compatibility
cannot increase. The matrix is symmetric, giving the other-side
assertion. Since $F\ne0$, the empty coordinate is only clipped.
Its additional bound follows from $\rho_0\leq\pi_0+v$.
Finally, multiply all originally selected weights of a nonfull type
by $\widehat f_I/f_I$ when $f_I>0$, and leave a zero coordinate empty.
\end{proof}

For finite order intervals, the label coordinates can be bounded by
exact counts. No transport on that side is then needed. Write
$\overline C_I(m)$ for the block-type prefix count of
Section~\ref{sec:arithmetic-types}, and put
\begin{equation}\label{eq:odd-endpoint-caps}
 c_I(h)=\overline C_I(h)-\overline C_I(\floor{h/2}).
\end{equation}
Multiplication by two preserves these types, so $c_I(h)$ counts
the odd integers of type $I$ in $[h]$.

\begin{corollary}[Uniform finite coordinate bounds]\label{cor:finite-coordinate-caps}
Suppose the order belongs to an integer interval $[l,h]$, all residual
odd labels lie in $[n]$, and their number satisfies $R\geq R_0>0$.
Every selected label submeasure, normalized by $R$, has coordinate
at most $c_I(h)/R_0$. If the selection excludes at least $uR$
universal empty-type labels, its empty coordinate is at most
\[
 \max\{0,c_0(h)/R_0-u\}.
\]
Suppose the actual position-tag law $\rho(n)$ satisfies
$\|\rho(n)-\pi(n)\|_{\mathrm{TV}}\leq v_*$ for a specified probability
law $\pi(n)$. After the transport of
Lemma~\ref{lem:full-type-transport}, every selected position submeasure
has coordinates bounded by
\begin{equation}\label{eq:finite-position-caps}
 \kappa_I=\max_{l\leq n\leq h}\pi_I(n)
                               +v_*\boldsymbol1_{I=F}.
\end{equation}
Its total mass is unchanged. The caps in these assertions need not
sum to one.
\end{corollary}

\begin{proof}
The selected population of a type is at most its complete endpoint
count. Divide by $R\geq R_0$. Excluding the universal labels subtracts
at least $u$ from the empty coordinate before this upper bound is
applied. The position assertion follows by applying
Lemma~\ref{lem:full-type-transport} at the actual order and then taking
the separate coordinate maxima.
\end{proof}

The reference law may be the exact usable even law
\[
 \pi_I(n)=\frac{a_I(n)}{\sum_Ja_J(n)},
\]
when this denominator is positive.
The comparison with the actual position law remains a hypothesis:
all rounding, reserve, and retagging errors must be included in $v_*$.
In particular, exact arithmetic quotas do not remove the separator
cost. For a nonempty label type $Q$, these hypotheses give the initial
degree fraction at least
\[
 \min_{l\leq n\leq h}\sum_{I\cap Q=\varnothing}\pi_I(n)-v_*.
\]
An empty-type label has initial degree exactly $R$, without a variation
loss, since it is compatible with every tag.

\subsubsection{Finite rectangle minima}

We next justify the finite minimizations used with these coordinate
bounds. The assertion is for a nonnegative matrix; it is not an
integrality assertion for arbitrary bilinear optimization.

\begin{proposition}[Capped bilinear minimum]\label{prop:capped-bilinear-minimum}
Let $C_{ij}\geq0$, let $p_i,q_j\geq0$ be finite coordinate caps, and
let $a,b,L\geq0$ with $a+b\leq L$. Consider
\begin{equation}\label{eq:capped-bilinear-domain}
 0\leq x_i\leq p_i,\quad 0\leq w_j\leq q_j,\quad
 \sum_i x_i+\sum_jw_j=L,\quad
 \sum_i x_i\geq a,\quad\sum_jw_j\geq b.
\end{equation}
The domain is nonempty precisely when
\begin{equation}\label{eq:capped-feasible-mass}
 \max\{a,L-\textstyle\sum_jq_j\}
       \leq\min\{\textstyle\sum_ip_i,L-b\}.
\end{equation}
On a nonempty domain the minimum of $\sum_{i,j}C_{ij}x_iw_j$ is
attained at a point of one of the following forms:
\begin{enumerate}
\item all but at most one coordinate, on the two sides together, are
      either zero or at their cap;
\item $\sum_i x_i=a$, with at most one coordinate on each side
      strictly between zero and its cap;
\item $\sum_j w_j=b$, with the same coordinate restriction.
\end{enumerate}
Thus rational data give a finite exact enumeration of the minimum.
\end{proposition}

\begin{proof}
Possible left masses form the interval in
\eqref{eq:capped-feasible-mass}; every mass between zero and the sum
of a side's caps can be filled on that side. This proves feasibility.
The domain is compact, so a minimum exists.

Start from a minimizer. If neither side-mass constraint is an equality
and two coordinates are strictly between their bounds, transfer mass
between them, leaving their sum fixed. If they are on the same side,
the objective is affine in the transfer. If they are on opposite
sides, say $x_i$ increases by $t$ and $w_j$ decreases by $t$, the
quadratic coefficient is $-C_{ij}\leq0$. In either case the objective
is concave on the permitted transfer interval. One endpoint has value
no greater than the starting value. Move to that endpoint. A coordinate
now reaches a bound, or a side-mass constraint becomes an equality.

Keep coordinates at their bounds fixed. If a side-mass constraint
has become an equality, both side masses are fixed, since their sum
is $L$. Transfer mass only within one side thereafter; the objective
is affine and an endpoint puts another coordinate at a bound. Each
step decreases the number of unbounded coordinates unless it first
fixes the side masses. The procedure terminates in one of the stated
forms. To enumerate them, choose the coordinates at their caps, the
possible remaining coordinate on each side, and the indicated mass
equality. The remaining coordinates are determined by these equalities.
All candidate values are rational when the data are rational.
\end{proof}

For caps $p,q$ on $\mathcal T$, let $\mathcal D(p,q;a,b)$ be the
domain in \eqref{eq:capped-bilinear-domain} with $L=1$, and write
\[
 \mathcal B(p,q;a,b)
   =\min_{(x,w)\in\mathcal D(p,q;a,b)}
                            \sum_{I,J}C(I,J)x_Iw_J.
\]
When the domain is empty we omit that case rather than assign it a
deletion event.

\begin{theorem}[Hall criterion with three type ranges]\label{thm:three-range-hall}
Assume the hypotheses of Lemma~\ref{lem:three-hall-ranges}. Let
$p_I\geq\abs{\{y\in Y:T(y)=I\}}/R$. Suppose every selected position
submeasure admits the transport of
Lemma~\ref{lem:full-type-transport} with caps $\kappa_I$.
Define three label cap vectors by
\[
 p_I^{(1)}=p_I\boldsymbol1_{I=F},\qquad
 p_I^{(2)}=p_I\boldsymbol1_{I\ne0},\qquad
 p_I^{(3)}=
 \begin{cases}\max(0,p_0-u),&I=0,\\p_I,&I\ne0.\end{cases}
\]
The ranges, their mass bounds, and their relevant label types are
given in Table~\ref{tab:hall-ranges}. For each nonempty range whose
capped domain is feasible, put
\[
 \beta_k=\mathcal B(p^{(k)},\kappa;a_k,b_k).
\]
If fewer than $R^2\beta_k$ edges of $\Gamma_0$ incident with the
types in the last column of the table are deleted, for every such
range, then $\Gamma$ has a perfect matching.
\end{theorem}

\begin{table}[htbp]
\centering
\renewcommand{\arraystretch}{1.3}
\begin{tabular}{c c c c c}
\hline
$k$ & Deficient-set mass & $a_k$ & $b_k$ & Label types counted\\
\hline
1 & $\tau<\abs U/R\leq\eta$ & $\tau$ & $1-\eta$ & $F$\\
2 & $\eta<\abs U/R\leq\gamma$ & $\eta$ & $1-\gamma$ & $I\ne0$\\
3 & $\gamma<\abs U/R\leq1-u$ & $\gamma$ & $u$ & all types\\
\hline
\end{tabular}
\par\vspace{.6em}
\caption{Mass bounds after trimming a deficient rectangle to total
mass one. A range with coincident endpoints is empty.}
\label{tab:hall-ranges}
\end{table}

\begin{proof}
Suppose there is a deficient set $U$, and put
$V=B\setminus N_\Gamma(U)$. By Lemma~\ref{lem:three-hall-ranges},
it belongs to one of the three ranges. Write $\xi=\abs U/R$ and
$\zeta=\abs V/R$. In range $k$, we have
$\xi\geq a_k$, $\zeta\geq b_k$, and $\xi+\zeta>1$.
Choose
\[
 t\in[\max\{a_k,1-\zeta\},\ \min\{\xi,1-b_k\}].
\]
The interval is nonempty: each term in its lower endpoint is at most
each term in its upper endpoint. Reduce the weights of the members
of $U,V$ to have respective normalized masses $t,1-t$. The label
coordinates obey $p^{(k)}$; in the third range all universal labels
were excluded from $U$.

Transport only the selected position mass to the caps $\kappa$.
Mass is preserved and compatibility can only decrease. The resulting
pair lies in the capped domain, so this domain is feasible and its
compatible mass is at least $\beta_k$. Before transport, all of these
compatible pairs belonged to $U\times V$, which has no edge of
$\Gamma$. Its weighted compatible mass times $R^2$ is therefore at
most the number of deleted pairs of the relevant types. This
contradicts the strict deletion bound. Hall's criterion now gives
the matching.
\end{proof}

\begin{corollary}[Simultaneous probability bound]\label{cor:three-range-probability}
Fix the type assignment, tags, caps, and degree targets before
permuting the even labels within their classes. Suppose the degree
conditions of Lemma~\ref{lem:three-hall-ranges} fail with probability
at most $p_{\mathrm{loc}}$, and the required universal labels exist
for every permutation. For each feasible nonempty range, suppose
the probability of at least $R^2\beta_k$ relevant deletions is at
most $p_k$. If
\begin{equation}\label{eq:three-range-probability}
 p_{\mathrm{loc}}+\sum_kp_k<1,
\end{equation}
then some choice of the class permutations admits a perfect matching.
\end{corollary}

\begin{proof}
The union bound gives positive probability that all degree and deletion
conditions hold simultaneously. Apply Theorem~\ref{thm:three-range-hall}
to that choice. No independence between these events is needed.
\end{proof}

Corollary~\ref{cor:local-conflict-tails} supplies $p_{\mathrm{loc}}$,
using a loss allowance that leaves each required degree strictly
above its target. Proposition~\ref{prop:global-conflict-tails} supplies
$p_k$ by restricting the odd labels counted in the cost to the
corresponding column of Table~\ref{tab:hall-ranges}. A range excluded
by a deterministic argument requires no probability term.

\subsubsection{Arithmetic restrictions from position types}

We now give such a deterministic argument for the middle range.
It needs a stronger property of the type assignment. Suppose
$G=(A,B)$ is the actual residual forest, $\Delta_G(A)\leq D$, where
$D\geq1$ is an integer, and $C_B\subseteq B$ has full tags. Require
\begin{equation}\label{eq:pure-neighbor-tags}
  b\notin C_B,\ ab\in E(G)\quad\Longrightarrow\quad I(a)=J(b).
\end{equation}
Reserve positions may be isolated, in which case the requirement is
vacuous. This property holds when components of a completed forest
after deleting $C_B$ receive a single type, and the separating
positions are assigned the full tag. Changing a tag at one position
changes its empirical law by total variation at most $1/R$.
Thus these full tags are permitted with the same retagging charge
$\abs{C_B}/R$; they must be included in the position error bound.
The inclusion of neighbor types remains valid. Condition
\eqref{eq:pure-neighbor-tags} does not follow from that inclusion
alone.

Fix an integer endpoint $h$ and use the complete half-label pools
\[
 \Lambda_J^+(h)=\{a\leq\floor{h/2}:T(a)=J\}.
\]
They include half-labels removed in preparing the residual problem.
For an odd $y\leq h$ and $J\cap T(y)=\varnothing$, define
$\operatorname{bad}_h(y,J)$ to mean that some member of
$\Lambda_J^+(h)$ shares a prime greater than $z$ with $y$.
It is an exact finite predicate:
\begin{equation}\label{eq:arithmetic-support-test}
 \operatorname{bad}_h(y,J)\quad\Longleftrightarrow\quad
 \sum_{\substack{p\mid y,\ p\text{ prime}\\z<p\leq h/2}}
     \overline C_J\left(\floor*{\frac{h}{2p}}\right)>0.
\end{equation}
Indeed multiplication by $p>z$ preserves the block type; positivity
asserts that at least one such multiple is in the pool. Multiplicity
in this sum is harmless because only positivity is used.

For a nonempty set $S\subseteq\mathcal T$ of position types, define
$c_Q(h,S)$ to be the number of odd $y\leq h$ of nonempty type $Q$
such that, for every $J\in S$,
\begin{equation}\label{eq:support-label-condition}
 Q\cap J\ne\varnothing\quad\hbox{or}\quad
                    \operatorname{bad}_h(y,J).
\end{equation}
Set $c_0(h,S)=0$. These are counts of labels, not counts of prime
divisor patterns.

\begin{lemma}[Support restrictions]\label{lem:arithmetic-slot-support}
Assume \eqref{eq:pure-neighbor-tags}, all actual labels lie in $[h]$,
and $R\geq R_0>0$. Let $U\times V$ contain no edge of $\Gamma$,
where every member of $U$ has nonempty type. Give these vertices
weights in $[0,1]$ and transport their selected position submeasure
by Lemma~\ref{lem:full-type-transport}. Write $x_Q,w_J$ for the
normalized label and transported position masses, and suppose the
position mass is positive. If $S=\{J:w_J>0\}$, then
\begin{equation}\label{eq:support-coordinate-cap}
 x_Q\leq c_Q(h,S)/R_0.
\end{equation}
Every nonfull position coordinate after transport can still be
represented by weights on the original positions of that tag.
\end{lemma}

\begin{proof}
By Lemma~\ref{lem:full-type-transport}, the nonfull coordinates
retain their original positions. Fix a positively weighted
$y\in U$ and $J\in S$. If $J=F$, its intersection with $T(y)$
is nonempty. If $J\ne F$, positive mass at $J$ is carried by an
actual position $b\in V\setminus C_B$. Suppose $J\cap T(y)$ is
empty and $\operatorname{bad}_h(y,J)$ is false. Every neighbor of
$b$ has a half-label in $\Lambda_J^+(h)$, by
\eqref{eq:pure-neighbor-tags}. Disjoint block types exclude common
primes at most $z$, and the failed badness predicate excludes all
larger ones. Since $y$ is odd, it is coprime to every even label
at a neighbor of $b$. Thus $yb\in E(\Gamma)$, a contradiction.

Every selected label therefore satisfies
\eqref{eq:support-label-condition}. Its weight is at most one, so
the total selected weight of type $Q$ is at most $c_Q(h,S)$. Divide
by $R\geq R_0$ to obtain \eqref{eq:support-coordinate-cap}.
\end{proof}

We also separate the deterministic forest cost by position type.
For $a\in\Lambda_J^+(h)$, let $q_J(a)\geq0$ dominate the number
of odd $y\leq h$ with
\[
 T(y)\ne0,\qquad T(y)\cap J=\varnothing,\qquad a\sim_z y.
\]
An exact union count or the filtered prime sum of
Lemma~\ref{lem:filtered-prime-costs}, restricted to nonempty odd
types, is suitable. Take $q_F=0$. Put
\[
 B_J=\sum_{a\in\Lambda_J^+(h)}q_J(a),
\]
and let $L_J(t)$ be the sum of the largest $t$ entries in the
multiset containing $D-1$ copies of each $q_J(a)$. A fractional
last entry is allowed, and the function is constant after the
multiset is exhausted. In particular, $L_J=0$ when $D=1$.
Each $L_J$ is nondecreasing, concave, and piecewise linear.

\begin{proposition}[Separate forest bounds]\label{prop:separate-forest-hall}
Assume the graph and endpoint hypotheses of
Lemma~\ref{lem:arithmetic-slot-support},
and let $N\geq R$ be a specified upper bound. Suppose the transported
coordinates of every selected position submeasure satisfy
$w_J\leq\kappa_J$.
Fix $0<a\leq b<1$. For every nonempty $S\subseteq\mathcal T$,
consider the domain
\begin{equation}\label{eq:separate-forest-domain}
 \begin{gathered}
 0\leq x_Q\leq c_Q(h,S)/R_0,\qquad
 0\leq w_J\leq\kappa_J\boldsymbol1_{J\in S},\\
 \sum_Qx_Q+\sum_Jw_J=1,\qquad
                  a\leq\sum_Qx_Q\leq b.
 \end{gathered}
\end{equation}
If, on every nonempty such domain,
\begin{equation}\label{eq:separate-forest-slack}
 \sum_{Q,J}C(Q,J)x_Qw_J
     -\frac1{R_0^2}\sum_{J\in S}\bigl(B_J+L_J(Nw_J)\bigr)>0,
\end{equation}
then there is no deficient set $U$ of nonempty-type labels with
$a<\abs U/R\leq b$. This exclusion holds for every permutation
within the even classes.
\end{proposition}

\begin{proof}
Suppose such a set exists, and put $V=B\setminus N_\Gamma(U)$.
As in the proof of Theorem~\ref{thm:three-range-hall}, trim the two
sides to combined normalized mass one, keeping the label mass in
$[a,b]$. Transport only the position mass. It remains positive
because $b<1$. Let $S$ be its support. By
Lemma~\ref{lem:arithmetic-slot-support}, the resulting vectors
belong to \eqref{eq:separate-forest-domain}.

For a nonfull $J\in S$, use the retained position weights of that
tag, whose unnormalized sum is $Rw_J$. By
\eqref{eq:pure-neighbor-tags}, all their neighbors lie in the
even class $J$. Charge each deleted pair to the neighboring
half-labels that witness a common large prime, as in the proof of
Corollary~\ref{cor:forest-rectangle-deletions}. Its weighted count
is at most $\sum_{a\in A_J}d_w(a)q_J(j(a))$, because the selected
odd-label weights are at most one and all have nonempty type.
The graph on these vertices is a subforest of $G$, with even
degree at most $D$. Apply Lemma~\ref{lem:weighted-forest-cost}
to that sum. Extending from the usable pool to its complete
endpoint superset can only increase both the base sum and the
largest-cost sums. Thus the number of weighted deletions on these
positions is at most
\[
 B_J+L_J((Rw_J-1)_+)\leq B_J+L_J(Nw_J).
\]
The full coordinate, including any artificial mass introduced by
transport, has no compatible pair with a nonempty label. Its cost
is zero.

All compatible pairs on the retained nonfull positions are deleted,
since their weights still lie in $U\times V$. Hence their total
normalized compatible mass is at most
\[
 \frac1{R^2}\sum_{J\in S}\bigl(B_J+L_J(Nw_J)\bigr)
 \leq\frac1{R_0^2}\sum_{J\in S}\bigl(B_J+L_J(Nw_J)\bigr).
\]
This contradicts \eqref{eq:separate-forest-slack}. Every step was
uniform over the class permutations.
\end{proof}

For an order interval $[l,h]$, one may take any proved lower bound
$R_0$ and $N=\ceil{h/2}$. To exclude the middle range of
Table~\ref{tab:hall-ranges}, use $a=\eta$, $b=\gamma<1$.
An infeasible support needs no check. Allowing zero coordinates for
types inside an enumerated $S$ enlarges the domain and is safe;
their base costs $B_J$ are still paid. Types outside $S$ contribute
neither base nor excess cost. The full-type and high ranges, and
the local degree conditions, retain their separate checks.

\subsubsection{Exact minimization with separate costs}

The preceding strict inequality is a finite rational test. We prove
this for a general nonnegative compatibility matrix, so that the
reason for the test does not depend on the arithmetic of a
particular interval.

\begin{lemma}[Finite minimization with convex column costs]\label{lem:separate-cost-minimum}
Let $C_{ij}\geq0$ be a finite rational matrix, $p_i,q_j$ nonnegative
rational caps, and $0\leq a\leq b\leq1$. For each $j$, let
$h_j$ be a continuous convex piecewise linear function on $[0,q_j]$
with rational breakpoints and rational affine pieces. The minimum of
\begin{equation}\label{eq:column-cost-objective}
 \Phi(x,w)=\sum_{i,j}C_{ij}x_iw_j+\sum_jh_j(w_j)
\end{equation}
subject to
\[
 0\leq x_i\leq p_i,\quad0\leq w_j\leq q_j,\quad
 \sum_i x_i+\sum_jw_j=1,\quad a\leq\sum_ix_i\leq b
\]
can be computed by finitely many rational evaluations, with
infeasibility detected by \eqref{eq:capped-feasible-mass} using
$L=1$ and lower right mass $1-b$.
\end{lemma}

\begin{proof}
Write $t=\sum_i x_i$. The feasible interval is
\[
 \max\{a,1-\textstyle\sum_jq_j\}
       \leq t\leq\min\{b,\textstyle\sum_ip_i\}.
\]
Assume it is nonempty. The domain is compact and the objective is
continuous, so a minimum exists. We describe an exhaustive finite
procedure.

\emph{Step 1: choose a left vertex.}
For fixed $t$ and $w$, the objective is affine in $x$. Its minimum
on the capped simplex $\sum_i x_i=t$ is attained at a point with
at most one coordinate strictly between its bounds. Indeed, if two
coordinates are interior, transfer mass between them to an endpoint
without increasing the affine objective, and repeat. Applied to a
global minimizer while keeping its $t,w$ fixed, this shows that some
global minimum has this form.

Enumerate a set $A_0$ of coordinates at their caps and a possible
remaining coordinate $i\notin A_0$. With $m_0=\sum_{k\in A_0}p_k$,
the corresponding family is
\[
 x_k=p_k\ (k\in A_0),\qquad x_i=t-m_0,
 \qquad x_k=0\ (k\notin A_0\cup\{i\}),
\]
restricted to $m_0\leq t\leq m_0+p_i$ and the feasible interval.
Include also all-bound choices, for which $t=m_0$ is a single
point. These families exhaust the required left vertices, including
degenerate caps.

\emph{Step 2: fill the right coordinates.}
For one such family, divide each interval $[0,q_j]$ into the linear
pieces of $h_j$. Their slopes are nondecreasing within a column.
After adding the bilinear term, each segment has a marginal cost
\[
 A_{jr}+B_jt,\qquad B_j=C_{ij}\geq0,
\]
where $r$ indexes the segment. Its length is a fixed rational number.
For fixed $t$, fill segments in increasing order of these marginal
costs until total length $1-t$ is reached. Ties may be ordered so
that earlier segments of a column come first. Convexity guarantees
that this order respects the within-column constraint.

This filling minimizes the cost: if positive length were used in a
more expensive available segment while a cheaper segment remained
unfilled, transferring the smaller available length would decrease
the objective. Equivalently, ignoring the within-column constraint
gives this same greedy optimum, and the ordered construction already
satisfies that constraint. Add the constant $\sum_jh_j(0)$ to the
segment costs. All-bound left choices are evaluated by this same
procedure at their fixed mass.

\emph{Step 3: find the masses that need evaluation.}
For a family with variable $t$, partition its closed feasibility
interval at every crossing of two distinct affine marginal costs.
On each open interval the segment order is fixed. Partition it
further wherever $1-t$ equals a cumulative segment length in that
order. Between consecutive points, every used segment is fixed in
full except possibly the last. If its preceding length is $d_0$,
the last contribution is
\[
 (1-t-d_0)(A+Bt),\qquad B\geq0.
\]
All full-segment contributions are affine in $t$. Thus the optimized
value on this interval is a quadratic with leading coefficient
$-B\leq0$, or an affine function if no segment is partial. It is
concave and attains its minimum at an endpoint.

At a cost crossing, tied segments have equal marginal cost, so
changing their order does not change the limiting value. When the
last segment becomes empty or full, the two formulas also agree.
The optimized value is therefore continuous at all subdivision
points. Evaluating the finite list of endpoints gives the minimum
for this family, and then for all families. Crossings, cumulative
lengths, and feasibility endpoints are rational. Hence every
evaluation is rational, as claimed.
\end{proof}

To apply the lemma to \eqref{eq:separate-forest-slack}, use the
label caps in \eqref{eq:separate-forest-domain} and, for $J\in S$,
\[
 h_J(w)=-\frac{B_J+L_J(Nw)}{R_0^2}.
\]
These functions are convex and piecewise linear. For $J\notin S$,
the cap and cost are zero. With rational costs and parameters the
test is therefore exact. The use of $L_J(Nw)$ is deliberate:
replacing its argument by $(Nw-1)_+$ would introduce a slope change
that need not preserve this convexity. The weaker forest bound
used here is sufficient and has the required finite minimization.

\section{Certificates for the orders through 1359}
\label{sec:part-finite}

This part proves Theorem~\ref{thm:initial-prefix-1359}: every forest of
order $n\leq1359$ has a coprime labeling by $[n]$.
Section~\ref{sec:finite-certificates} defines finite proof graphs for
the statements $\cP(L)$ and proves that their rules are sound
(Theorem~\ref{thm:finite-certificate-soundness}).
Sections~\ref{sec:branch-peeling}--\ref{sec:larger-blocks} add
branch peeling and certified rooted blocks, which reduce the number of
cases. Sections~\ref{sec:initial-prefix}--\ref{sec:finite-endpoint}
apply them to the orders $1$--$299$, $300$--$1352$, and
$1353$--$1359$.

\subsection{Finite labeling certificates}\label{sec:finite-certificates}

We return to the finite reductions of Section~\ref{sec:finite-reductions}.
Their arithmetic hypotheses concern every forest on a specified label
set, sometimes with a restriction on an attachment vertex. We describe
finite certificates for these assertions. The role of a certificate is
to provide the choices in a proved reduction, together with a complete
check of every terminal case.

It is useful to separate the integer labels from the graph of permitted
pairs. A single graph certificate can then be used for several label
sets, provided its edges and attachment restrictions are checked anew
for each set.

\subsubsection{Compatibility graphs and marked vertices}

Let $H$ be a simple graph on a set $S$ of $m\geq1$ positions, and
let $O\subseteq S$. For $0\leq k\leq m$, write
$\mathcal U_k(H,O)$ if, for every tree $T$ of order $m$ and every
$M\subseteq V(T)$ with $\abs M\leq k$, there is a bijection
$f\colon V(T)\to S$ such that
\begin{equation}\label{eq:marked-host-contract}
 f(v)f(w)\in E(H)\quad(vw\in E(T)),\qquad f(M)\subseteq O.
\end{equation}
Vertices of $M$ are called marked. The set $O$ is irrelevant when
$k=0$; we then write $\mathcal U_0(H)$.

The same assertions hold for forests. Complete the input forest to a
tree by Lemma~\ref{lem:forest-completion}, retaining its marked set, and
restrict the resulting labeling. Conversely trees are forests. It
suffices to check marked sets of size exactly $k$, since smaller ones
can be enlarged. If $O=S$, the condition on marks is automatic, so
$\mathcal U_0(H)$ implies $\mathcal U_k(H,S)$ for every $k\leq m$.

For a set $L$ of distinct positive integers, its \emph{coprimality
graph} $H_L$ has vertex set $L$ and joins distinct $x,y$ when
$\gcd(x,y)=1$. Thus $\cP(L)$ is equivalent to $\mathcal U_0(H_L)$
when $L\ne\varnothing$. If a branch is attached to a vertex labeled
$c\notin L$, its required assertion is
\begin{equation}\label{eq:arithmetic-root-contract}
 \mathcal U_1(H_L,O_c(L)),\qquad
 O_c(L)=\{x\in L:\gcd(c,x)=1\}.
\end{equation}
The marked vertex is the actual attachment vertex. It need not be
a leaf or a centroid.

\begin{lemma}[Transfer of a host certificate]\label{lem:host-transfer}
Let $H,H'$ have the same order $m$, let $O\subseteq V(H)$ and
$O'\subseteq V(H')$, and let $0\leq k\leq k'\leq m$. Suppose
a bijection $\sigma\colon V(H)\to V(H')$ maps every edge of $H$
to an edge of $H'$. If $k>0$, assume also $\sigma(O)\subseteq O'$.
Then
\[
 \mathcal U_{k'}(H,O)\quad\Longrightarrow\quad
                                  \mathcal U_k(H',O').
\]
In particular, a graph certificate can be bound to an integer label
set by an explicit bijection preserving all required edges and, when
needed, the allowed attachment labels.
\end{lemma}

\begin{proof}
Given $T,M$ with $\abs M\leq k$, take a labeling $f$ furnished by
the source assertion. The composition $\sigma\circ f$ is bijective,
maps every tree edge to an edge of $H'$, and sends every marked
vertex into $O'$ when marks are present.
\end{proof}

Equality of the two compatibility graphs is sufficient but unnecessary.
Additional permitted edges cannot invalidate a labeling. The inclusion
of allowed attachment positions must, however, hold for the entire
source set $O$, not only for a position used by one previously found
labeling.

\subsubsection{Complete terminal tables}

A terminal table is useful only if its index set is complete. We give
a recursive enumeration whose completeness follows directly from the
decomposition of a rooted tree at its root.

Encode the one-vertex rooted tree by the string \texttt{()}. In
general, its code is an opening parenthesis, the codes of the rooted
child trees in lexicographic order, and a closing parenthesis. The
order of a code is its number of opening parentheses. Such a code
also specifies a representative: create a vertex at each opening
parenthesis, joined to the vertex of the innermost still open pair.
Vertices are indexed in their order of creation, starting with the
root.

\begin{lemma}[Enumeration by rooted codes]\label{lem:tree-code-enumeration}
Two rooted trees have the same code if and only if they are isomorphic
by an isomorphism preserving the root. All rooted types of order $m$
are obtained, once each, by taking all multisets of previously
enumerated rooted types whose orders sum to $m-1$, sorting their
codes, and enclosing the result in parentheses.

For an unrooted tree $T$, let
\[
 c(T)=\min_{v\in V(T)}c(T,v),
\]
where $c(T,v)$ is its code rooted at $v$. Retaining from the rooted
list just those codes that are the minimum over all rootings of
their decoded tree gives exactly one representative of each
unrooted tree type of order $m$.
\end{lemma}

\begin{proof}
Proceed by induction on the order. Matching outer parentheses
uniquely splits a root code into its child codes. A root-preserving
isomorphism permutes the child components and restricts to rooted
isomorphisms within them. Conversely, rooted isomorphisms between
matched child components extend, with the roots identified, to an
isomorphism of the whole rooted trees. By induction, equality of
the sorted child-code lists is therefore equivalent to rooted
isomorphism.

Every child has smaller order, and the child orders sum to $m-1$. Thus
the multiset recursion includes the children of every rooted tree.
Different multisets give different sorted lists and hence different
codes. Enumerating multisets as nondecreasing lists of type indices
introduces no repetitions. This proves the first assertion, including
its one-vertex base.

An unrooted isomorphism identifies the sets of rooted codes obtained
by choosing all possible roots. Hence it preserves their minimum.
If two minima agree, their decoded rooted trees are isomorphic,
and therefore their unrooted trees are isomorphic. The minimum is
itself one of the enumerated rooted codes, so each unrooted type
has exactly one retained representative.
\end{proof}

\begin{proposition}[Terminal witness criterion]\label{prop:terminal-witnesses}
Let $H$ have order $m\geq1$, and let $0\leq k\leq m$.
For every unrooted representative in
Lemma~\ref{lem:tree-code-enumeration} and every $k$-subset of its
vertices, suppose a table contains a bijection satisfying
\eqref{eq:marked-host-contract}. Then $\mathcal U_k(H,O)$ holds.
For $k=1$, one may instead use one witness for each rooted code,
requiring that its root map into $O$.

For an integer label set $L$ and $d\geq2$, a complete unrooted table
restricted to trees of maximum degree at most $d$, with a coprime
labeling by $L$ in each row, proves $\cP_d(L)$.
\end{proposition}

\begin{proof}
Enlarge a given marked set to size $k$. Map the input tree to its
unrooted representative by an isomorphism. The image of the marked
set is among the checked subsets; compose the corresponding witness
with that isomorphism. This gives the desired bijection and preserves
every edge and mark. For one marked vertex, apply the rooted
enumeration instead. All choices of attachment vertex are included
in its rooted types.

In the degree-restricted case, first complete the input forest to a
tree of maximum degree at most $d$, using
Lemma~\ref{lem:forest-completion}. Its representative occurs in the
restricted table. Transfer the witness back and restrict to the
original forest.
\end{proof}

The checks in each row are finite: the listed images form a permutation
of the host positions, every tree edge maps to a host edge, and every
mark has an allowed image. For integer labels, the edge check is an
exact gcd computation. Coverage requires equality with the entire
expected index set, not merely the expected number of rows. Counting
rows is sufficient only if validity, distinctness, and membership in
that expected set have also been checked.

\begin{proposition}[Two terminal label sets]\label{prop:two-terminal-sets}
The following label sets satisfy $\cP(L_{13})$ and $\cP(L_{14})$:
\[
 \begin{split}
 L_{13}&=[10]\cup\{12,14,15\},\\
 L_{14}&=\{1,2,3,4,7,8,9,10,15,16,20,21,22,27\}.
 \end{split}
\]
\end{proposition}

\begin{proof}
This is a finite verification by complete terminal tables. The
recursion of Lemma~\ref{lem:tree-code-enumeration} gives $1301$
unrooted types of order $13$ and $3159$ of order $14$. The
accompanying terminal-table supplement contains a labeling for
each of these types with the indicated inventory.

The supplied checker reconstructs both complete index sets and
checks their equality with the table keys. Every row is a bijection
onto its stated inventory, and every edge has coprime endpoint
labels. There are $1301\cdot12=15612$ edge checks in the first
table and $3159\cdot13=41067$ in the second; all pass. The row
format and a standalone checker are included in the supplement.
Proposition~\ref{prop:terminal-witnesses} now gives both forest
assertions. No bound on the degrees of these trees is imposed.
\end{proof}

\subsubsection{Recursive host certificates}

Complete tree tables become large. The following terminal rule avoids
enumeration when all positions in the designated set can receive an
arbitrary neighbor.

\begin{lemma}[A split host]\label{lem:split-host}
Let $H$ have order $m$, and suppose every vertex of $O$ is adjacent
to every other vertex of $H$. Put $e=m-\abs O$. If
\begin{equation}\label{eq:split-host-condition}
 e\leq\ceil*{\frac{m-k}{2}},
\end{equation}
then $\mathcal U_k(H,O)$ holds.
\end{lemma}

\begin{proof}
Given a tree $T$ and at most $k$ marks, delete the marked vertices.
Each component of the resulting forest is bipartite. Taking the larger
part in each component gives an independent set of size at least
$\ceil{(m-k)/2}$. Assign the $e$ labels outside $O$ to any $e$
vertices of this independent set, and the labels in $O$ to all
remaining vertices. Every marked vertex receives a label in $O$.
No edge has both endpoints labeled outside $O$, and every other
edge is permitted by the hypothesis on $O$.
\end{proof}

The next rule replaces a large tree by the components at a centroid.
The marks must be distributed among those components explicitly;
a prescribed attachment vertex cannot be assumed to be the centroid.

For $m\geq2$, a \emph{centroid case} for $k$ marks consists of
$h\in\{0,1\}$ and a multiset of pairs $(j_i,t_i)$ satisfying
\begin{equation}\label{eq:centroid-mark-cases}
 1\leq j_i\leq\floor{m/2},\quad 0\leq t_i\leq j_i,
 \qquad \sum_i j_i=m-1,\quad h+\sum_i t_i=k.
\end{equation}
There are finitely many cases. They can be enumerated by first
listing the integer partitions of $m-1$ with parts at most
$\floor{m/2}$, then distributing $k-h$ marks among the parts
within their size bounds, and finally sorting the pairs to remove
duplicates.

\begin{proposition}[Marked centroid reduction]\label{prop:marked-centroid-host}
Let $H$ have vertex set $S$ of order $m\geq2$, and let
$O\subseteq S$, $0\leq k\leq m$. For every case in
\eqref{eq:centroid-mark-cases}, suppose there are a center $c$ and a
partition
\[
 S=\{c\}\sqcup\bigsqcup_i S_i,\qquad \abs{S_i}=j_i,
\]
with $c\in O$ if $h=1$. For each part, require either of the
following two assertions:
\begin{enumerate}
\item $c$ is adjacent to every member of $S_i$, and
      $\mathcal U_{t_i}(H[S_i],O\cap S_i)$ holds;
\item there is $O_i\subseteq N_H(c)\cap S_i$, with $O_i\subseteq O$
      when $t_i>0$, for which
      $\mathcal U_{\min(j_i,t_i+1)}(H[S_i],O_i)$ holds.
\end{enumerate}
Then $\mathcal U_k(H,O)$ holds. A child assertion with a larger
mark budget suffices; if its designated set is its entire inventory,
an unmarked assertion suffices.
\end{proposition}

\begin{proof}
Enlarge the marked set of the input tree to size $k$. Choose a
centroid $z$, whose existence was proved in
Lemma~\ref{lem:centroid-order}. Each component of $T-z$ has order
at most $\floor{m/2}$. Let $h$ indicate whether $z$ is marked, and
let $j_i,t_i$ be the orders and mark counts of these components.
They give one of the listed cases. Match components of the same
size and mark count with its parts in any order, and assign $c$ to $z$.

Every component has a unique vertex adjacent to $z$. In the first
alternative, apply the child assertion to its inherited marks;
the attachment edge is valid because $c$ is adjacent to the whole
child inventory. In the second, additionally mark the attachment
vertex. There are at most $\min(j_i,t_i+1)$ distinct marks. The
child assertion puts its attachment in $N_H(c)$, and, if there are
inherited marks, puts all of them in $O$ as well. If $z$ itself is
marked, its label $c$ lies in $O$ by hypothesis.

The child bijections give all internal edges, and the attachment
checks give all remaining edges. The inventories are disjoint and
exhaust $S$, so their union is the required bijection. The assertions
about larger budgets and an unrestricted designated set follow from
the definition of $\mathcal U_k$.
\end{proof}

The choice $O_i=N_H(c)\cap O\cap S_i$ handles both inherited marks
and the new attachment. When $t_i=0$, the intersection with $O$ is
unnecessary.

We next give extensions that avoid new tree tables. Additional host
edges are allowed by Lemma~\ref{lem:host-transfer}.

\begin{lemma}[Extensions at leaves and a universal vertex]
\label{lem:host-local-extensions}
Let $H$ have $m\geq1$ vertices and designated set $O$.
\begin{enumerate}
\item Add a vertex $u$ adjacent to every vertex of $H$, obtaining
      $H^+$. If $0\leq k\leq m+1$ and
      $\mathcal U_{\min(k,m)}(H,O)$ holds, then
      $\mathcal U_k(H^+,O\cup\{u\})$ holds.
\item If $\mathcal U_0(H)$ holds, then the same $H^+$ satisfies
      $\mathcal U_1(H^+,\{u\})$.
\item Add vertices $u,v$, with $u$ adjacent to all other vertices
      and $v$ adjacent only to $u$. For $k\in\{0,1\}$,
      $\mathcal U_k(H,O)$ implies
      $\mathcal U_k(H^+,O\cup\{u\})$ for this enlarged host.
\item If $\mathcal U_1(H,\{a\})$ holds for some $a\in V(H)$,
      adding a single new vertex adjacent to $a$ gives an
      unmarked universal host.
\end{enumerate}
\end{lemma}

\begin{proof}
For the first assertion, remove a leaf of the input tree, assign
it $u$, and use the child assertion on the remaining tree with its
remaining marks. There are at most $\min(k,m)$ such marks. If the
removed leaf was marked, its new label is allowed; its only edge
is valid because $u$ is universal.

For the second assertion, designate a root, choosing an arbitrary
vertex if no mark was specified. Delete it and label the remaining
forest using $H$. Put $u$ on the deleted vertex. Every restored
edge meets $u$, proving the required fixed-root assertion.

For the third assertion, the input has at least three vertices.
If there is a marked vertex, choose a leaf different from it; a
tree of order at least two has at least two leaves. Without a mark,
choose any leaf. Remove this leaf $x$ and its neighbor $y$, and
place them at $v$ and $u$, respectively. If $y$ was marked, its
label $u$ is allowed. Otherwise the mark remains in the residual
forest and the source assertion places it in $O$. Every edge
incident with the removed vertices meets $y$, so its image meets
the universal vertex $u$. This proves the assertion, also when
the residual forest has a single vertex.

For the last assertion, delete any leaf of the input tree. Root the
remaining tree at its neighbor, embed it in $H$ with that root at
$a$, and assign the new host vertex to the deleted leaf.
\end{proof}

The third assertion retains a set of possible root positions. It
does not force the input root to the new universal position, since
the root may remain in the residual forest.

\begin{lemma}[An attached triangle]\label{lem:host-attached-triangle}
Suppose $H$ has order $m\geq1$ and
$\mathcal U_1(H,\{a\})$ holds. Add two vertices $x,y$ and the
three edges $ax,ay,xy$. The resulting host $H^+$ satisfies
$\mathcal U_0(H^+)$.
\end{lemma}

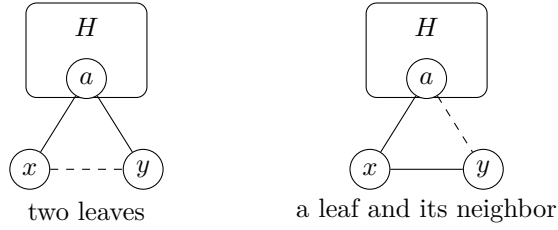
\begin{figure}[htbp]
\centering
\begin{tikzpicture}[scale=1,
 vertex/.style={circle,draw,inner sep=1.5pt,minimum size=15pt},
 every node/.style={font=\small}]
 \foreach \id/\shift/\title in {1/0/{two leaves},2/4.5/{a leaf and its neighbor}} {
   \begin{scope}[xshift=\shift cm]
   \draw[rounded corners] (-.8,-.55) rectangle (.8,.7);
   \node at (0,.38) {$H$};
   \node[vertex,fill=white] (a\id) at (0,-.3) {$a$};
   \node[vertex] (x\id) at (-.75,-1.5) {$x$};
   \node[vertex] (y\id) at (.75,-1.5) {$y$};
   \node[align=center] at (0,-2.05) {\title};
   \end{scope}
 }
 \draw (a1)--(x1) (a1)--(y1);
 \draw[dashed] (x1)--(y1);
 \draw (a2)--(x2) (x2)--(y2);
 \draw[dashed] (a2)--(y2);
\end{tikzpicture}
\caption{The attached triangle in Lemma~\ref{lem:host-attached-triangle}.
The remaining tree is embedded in $H$ with its attachment at $a$.
Solid edges restore the deleted vertices; the dashed edge is available
in the host but is unused in that case.}
\label{fig:attached-triangle}
\end{figure}

\begin{proof}
Every tree of order at least three has either two leaves with a
common neighbor or a leaf whose neighbor has degree two. To see
this, take an endpoint $x'$ of a longest path and its neighbor $y'$.
If $y'$ has degree two, the second alternative holds. Otherwise
$y'$ has a neighbor off that path. Every such neighbor is a leaf:
a path of length at least two into its component, followed by the
rest of the original longest path, would be longer. Thus $x'$
and one of these other leaves have a common neighbor.

In the first alternative, delete the two leaves and call their
common neighbor $a'$. In the second, delete the leaf and its
degree-two neighbor, and let $a'$ be the latter's other neighbor.
The remaining graph is a tree of order $m$. Embed it in $H$ with
$a'$ at $a$. In the first case use $x,y$ for the two leaves. In
the second put the degree-two vertex at $x$ and its leaf at $y$.
The edges $ax,ay$ or $ax,xy$, respectively, restore the input tree.
Figure~\ref{fig:attached-triangle} shows the two possibilities.
\end{proof}

This conclusion is unmarked: it gives no arbitrary prescribed-root
condition for $H^+$. Adding a universal vertex afterwards restores
such a condition by Lemma~\ref{lem:host-local-extensions}(2).

\subsubsection{From branch contracts to a label-set assertion}

The next reduction explains which rooted block contracts suffice for
an entire tree. Fix integers $b\geq2$ and $n\geq b+2$. A
\emph{short-branch profile} is a vector of nonnegative integers
$\mathbf k=(k_1,\ldots,k_{b-1})$ such that
\begin{equation}\label{eq:short-branch-profile}
 b\leq s=\sum_{j=1}^{b-1}jk_j\leq n-2,
 \qquad \sum_{j=1}^{b-1}j a_j\ne b
       \quad(0\leq a_j\leq k_j, a_j\in\Z).
\end{equation}
The second requirement says that no selection of branches has total
order $b$.

\begin{lemma}[A detached set or short branches]\label{lem:stump-dichotomy}
Every tree $T$ of order $n\geq b+2$ has one of the following
configurations:
\begin{enumerate}
\item a vertex $u$ and a union of components of $T-u$ of total
      order $b$;
\item adjacent vertices $u,v$ such that, on the side of the edge
      $uv$ containing $v$, the components attached to $v$ other
      than the one containing $u$ have orders below $b$ and
      multiplicities forming a short-branch profile.
\end{enumerate}
All components on the $v$ side are included in the second alternative.
\end{lemma}

\begin{proof}
Root $T$ at a leaf and choose a deepest vertex $v$ whose descendant
subtree has at least $b+1$ vertices. The child of the root has order
$n-1\geq b+1$, so $v$ is not the root and has a parent $u$.
Every child subtree of $v$ has at most $b$ vertices. If some
selection of those subtrees has total order $b$, it gives the
first alternative, with contact vertex $v$.

Otherwise none has order $b$, all have orders below $b$, and their
total is at least $b$. It is at most $n-2$, since $u,v$ are not
among their vertices. The multiplicities have no submultiset of
total order $b$, by the failure of the first choice. They therefore
satisfy \eqref{eq:short-branch-profile}, and the edge $uv$ gives
the second alternative.
\end{proof}

\begin{proposition}[Extension through rooted branch blocks]
\label{prop:rooted-branch-extension}
Let $L$ be a set of $n\geq b+2$ labels, where $b\geq2$.
Suppose there is a partition
\[
 L=\{q\}\sqcup B\sqcup C,\qquad \abs B=b,
\]
such that $q$ is universal in $L$ and $\cP(B),\cP(C)$ hold.
For each short-branch profile, suppose there is one of the following
partitions, with $k_j$ sets $B_{j,r}$ of size $j$:
\begin{align*}
 &L=\{c\}\sqcup C\sqcup\bigsqcup_{j,r}B_{j,r},
       &&\gcd(c,x)=1\quad(x\in C);\tag{retain}\label{eq:branch-retain}\\
 &L=\{q,c\}\sqcup C\sqcup\bigsqcup_{j,r}B_{j,r},
       &&q\text{ universal in }L.\tag{remove}\label{eq:branch-remove}
\end{align*}
In either case require $\cP(C)$ and
$\mathcal U_1(H_{B_{j,r}},O_c(B_{j,r}))$ for every block.
The partitions and their labels may depend on the profile.
Then $\cP(L)$ holds.
\end{proposition}

\begin{proof}
It suffices to label an arbitrary tree of order $n$, by
Lemma~\ref{lem:forest-completion}. Apply
Lemma~\ref{lem:stump-dichotomy}. In its first alternative, put $q$
on the contact vertex, label the detached forest by $B$, and label
the remaining forest by $C$. Every crossing edge meets the contact,
so its labels are coprime.

In the second alternative, put $c$ on $v$. For each size $j$, match
the $k_j$ blocks with the $k_j$ actual child subtrees. Root each
subtree at its unique neighbor of $v$ and use its rooted block
contract. This proves all internal branch edges and all edges from
the branches to $v$.

In the retain case, delete $v$ and all its child subtrees and label
the remaining tree by $C$. The parent $u$ remains there. The only
edge from that remainder to the deleted part is $uv$, and its labels
are coprime because the label of $u$ lies in $C$.

In the remove case, also delete $u$, put $q$ on it, and label the
remaining forest by $C$. Every further crossing edge meets $u$ and
is valid by universality of $q$. The empty remainder is allowed.
In both cases all inventories are disjoint and exhaust $L$. The
result is the required bijection, with every edge checked.
\end{proof}

The two rows distinguish whether the parent belongs to the recursively
labeled remainder. Their remainder orders are $n-s-1$ and $n-s-2$,
respectively. Neither may be replaced by a different initial interval
without checking the actual inventory.

For clarity, the profile index set itself has a finite exact test.
For a proposed vector $\mathbf k$, start with $A_0=\{0\}$ and
compute
\[
 A_j=\{a+t j:a\in A_{j-1},\ 0\leq t\leq k_j\}\cap[0,b]
             \qquad(1\leq j<b).
\]
Here $[0,b]$ denotes the integer set $\{0,1,\ldots,b\}$.
Inductively $A_j$ is precisely the set of sums at most $b$ formed
using the first $j$ branch sizes. Thus the last condition in
\eqref{eq:short-branch-profile} is exactly $b\notin A_{b-1}$.
Enumerating each $k_j$ between zero and $\floor{(n-2)/j}$ and
retaining just the vectors passing this test and the total-size
test gives every required profile. Smaller sufficient families of
retained cores will be treated separately; they require a selection
lemma before this full checklist can be shortened.

\subsubsection{Soundness of a finite proof graph}

We can now specify the logical content of a certificate, independently
of how its witnesses were found. A state records either a marked
host assertion $\mathcal U_k(H,O)$, an arithmetic forest assertion
$\cP(L)$, or a degree-restricted assertion $\cP_d(L)$ with $d\geq2$.
Its data include the exact graph or integer inventory and all mark
or degree restrictions. A proof node names a rule, supplies its
finite choices, and refers to the states used as hypotheses.

The rules in the present system are complete terminal tables,
the split-host rule, marked centroid reduction, the local host
extensions, the attached triangle, host transfer and arithmetic
binding, rooted branch extension, and the leaf-pair and detached-block
rules of Section~\ref{sec:finite-reductions}. Trivial empty arithmetic
inventories and singletons are terminal states. Forgetting a mark
restriction is allowed; forgetting a degree restriction is allowed
only when that degree cap already includes every tree of the stated
order. A host on one vertex meets a positive mark requirement
exactly when its position is allowed.

In particular, a terminal degree-restricted table is checked only
against its declared tree family and is never registered as an
unrestricted assertion by that check alone. In the leaf-pair rule,
the child inventory and degree cap must satisfy
Corollary~\ref{cor:leaf-label-extension}. In a detached-block rule,
every required size in Proposition~\ref{prop:block-extension} must
have its own checked partition and child assertions.

\begin{theorem}[Soundness of the certificate rules]
\label{thm:finite-certificate-soundness}
Let a finite directed graph of proof nodes use only the preceding
rules. Suppose every referenced node is present, there is no directed
cycle, every rule's arithmetic and graph hypotheses are checked, and
every finite case list is complete. Then the assertion at every node
is true. In particular, an accepted root state $\cP([n])$ proves a
prime labeling for every tree of order $n$.
\end{theorem}

\begin{proof}
Follow an order in which all dependencies of a node precede it;
such an order exists for every finite acyclic directed graph.
Terminal states hold by their direct checks,
Proposition~\ref{prop:terminal-witnesses}, or
Lemma~\ref{lem:split-host}. At a nonterminal node all referenced
assertions already hold. Its checked rule is an implication proved
above or in Section~\ref{sec:finite-reductions}, so its assertion
holds as well. Host binding is the composition of
Lemma~\ref{lem:host-transfer} with the exact coprimality and allowed-root
relations in \eqref{eq:arithmetic-root-contract}.

For a centroid node, completeness means that every size and mark
case in \eqref{eq:centroid-mark-cases} is present. For a branch
extension it means that every profile in
\eqref{eq:short-branch-profile} is present. The proofs of the
corresponding reductions show that every input tree selects one
of these checked cases. Terminal completeness has the same role
through the tree-code enumeration. Thus induction proves the full
universal assertion at each node, and not only a statement about
the representatives used during discovery.
\end{proof}

The centroid rule strictly decreases order, but a host transfer
preserves it. Such transfers must still be checked for cycles.
One sufficient convention is to permit an equal-order transfer only
from a host with fewer edges; then the pair consisting of order and
edge count decreases along those dependencies. Explicit cycle
detection also permits transfers between isomorphic states without
assuming that order alone makes the proof well founded.

This theorem concerns the stated rule system. Further host rules or
shortened profile lists require their own implications and coverage
proofs before they can be used. To establish a finite prefix, one
must also supply certificates whose accepted roots include every
initial interval in that prefix. Those inventories and coverage
checks are separate from the soundness result just proved.

\subsection{Peeling branches and finite-prefix induction}\label{sec:branch-peeling}

The branch extension in Section~\ref{sec:finite-certificates} asks for
every short-branch profile, whose total may be almost $n$. We now
retain only a bounded collection of branches. The others receive
consecutive intervals, and the exterior forest receives an initial
interval of its actual order. This reduces the number of arithmetic
assertions without changing the set of input trees covered.

Throughout this section $b\geq2$, $n\geq b+2$, and a profile has
branch sizes in $\{1,\ldots,b-1\}$. For a profile $\mathbf k$, put
\[
 |\mathbf k|=\sum_{j=1}^{b-1}jk_j,
 \qquad o(\mathbf k)=\sum_{\substack{1\leq j<b\\j\text{ odd}}}k_j.
\]
Thus $|\mathbf k|$ counts vertices and $o(\mathbf k)$ counts odd-order
branches. A subprofile $\mathbf h\leq\mathbf k$ is obtained by
retaining some whole branches. The short-branch condition remains
\eqref{eq:short-branch-profile}; in particular no subprofile has
total $b$. For integers $a,j\geq1$, write
$I(a,j)=\{a,a+1,\ldots,a+j-1\}$.
Intervals with reversed endpoints are empty. A \emph{branch
configuration} means the second alternative of
Lemma~\ref{lem:stump-dichotomy}, with adjacent vertices $u,v$
and all branches at $v$ on the side away from $u$.

\subsubsection{Selecting a bounded collection of branches}

\begin{lemma}[Minimal retained cores]\label{lem:minimal-branch-core}
Let $A\geq b$ and let $\mathbf k$ have total at least $A$.
One may require, in addition, that a retained subprofile contain
an odd-order branch, provided $o(\mathbf k)>0$. There is an
inclusion-minimal subprofile $\mathbf h$ meeting these requirements,
with total $t$ satisfying
\[
 A\leq t\leq A+b-2.
\]
If all branch orders are even, no odd branch is required, and
$E=2\floor{(b-1)/2}$, the upper bound improves to $A+E-1$.

Writing $\varepsilon=1$ when an odd branch is required and
$\varepsilon=0$ otherwise, minimality is equivalent to
\begin{equation}\label{eq:core-minimality}
 h_j>0\quad\Longrightarrow\quad
 t-j<A\ \text{ or }\
       (\varepsilon=1,\ j\text{ odd},\ o(\mathbf h)=1).
\end{equation}
Every retained subprofile of a short-branch profile still has no
subprofile of total $b$.
\end{lemma}

\begin{proof}
If an odd branch is required, choose one first. Its order is below
$b\leq A$. Continue choosing branches until their total first
reaches $A$. Without the extra requirement, start with no branches
and use the same procedure. Just before the last choice the total
is at most $A-1$; the last branch has order at most $b-1$, or at
most $E$ in the even case. These are the stated upper bounds.

Now delete a retained branch whenever the remaining total is at
least $A$ and any required odd branch is still present. Each
deletion decreases the number of retained branches, so the
procedure terminates. It preserves both requirements and cannot
increase the upper bound. A branch of order $j$ is removable
exactly when both alternatives on the right of
\eqref{eq:core-minimality} fail.

It remains to check inclusion-minimality, rather than merely the
failure of a single deletion. If a proper subprofile met the
requirements, choose a branch omitted from it. Deleting that
branch alone would leave a profile containing the proposed
subprofile, so it too would meet the requirements. This
contradicts termination. Conversely, an allowed single deletion
contradicts inclusion-minimality. Finally, a forbidden sum in a
subprofile would already have been a forbidden sum in the original
profile.
\end{proof}

At odd order we will apply this lemma with $A=b$ and no parity
requirement. At even order, a profile with an odd-order branch
will retain one. We call these the \emph{ordinary profiles}.
The remaining profiles occur at even order and have only even
branch sizes. Their threshold will depend on the chosen center
label. The distinction concerns branch orders; it does not
prescribe the parity of every label inside a branch.

\subsubsection{Arithmetic assertions for the retained branches}

We specify two assertions about a profile $\mathbf h$ of total
$t\leq n-2$. They concern only its retained branches and the
vertex or vertices at which they attach.

Write $\mathcal R_n(\mathbf h)$ if there are a label $q$ universal
in $[n]$ and a partition
\begin{equation}\label{eq:retained-core-inventory}
 \{n-t,\ldots,n\}
       =\{q\}\sqcup\bigsqcup_{j,r} B_{j,r},
 \qquad |B_{j,r}|=j\quad(1\leq r\leq h_j),
\end{equation}
such that $\cP(B_{j,r})$ holds for every block. The center will
receive $q$. No root restriction on these blocks is needed,
because $q$ is universal in the full interval.

For a specified center label $a$, write $\mathcal D_n(\mathbf h;a)$
if there are a universal label $q$ and a partition
\begin{equation}\label{eq:removed-core-inventory}
 \{n-t-1,\ldots,n\}
       =\{q,a\}\sqcup\bigsqcup_{j,r} B_{j,r},
 \qquad |B_{j,r}|=j\quad(1\leq r\leq h_j),
\end{equation}
such that every block satisfies
$\mathcal U_1(H_{B_{j,r}},O_a(B_{j,r}))$.
The center receives $a$, and its parent receives $q$.
The displayed partition requires $q\ne a$ and membership of both
labels in the stated top interval. Universality refers to $[n]$,
not just that interval.

These are finite arithmetic assertions. Their block conditions
may be supplied by the certificate rules of
Section~\ref{sec:finite-certificates}. Choosing a retained
profile does not itself establish either assertion.

\begin{proposition}[Peeling at a universal center]
\label{prop:ordinary-peeling}
Let $\mathbf k$ be a short-branch profile of total $s$, and let
$\mathbf h\leq\mathbf k$ have total $t\geq b$. Suppose
$\mathcal R_n(\mathbf h)$ holds. Put $m=n-s-1$ and order the
omitted branches arbitrarily. Give them consecutive intervals,
in that order, partitioning
\[
 \{m+1,\ldots,n-t-1\}.
\]
If each of these intervals satisfies $\cP$ and $\cP([m])$ holds,
then every branch configuration with profile $\mathbf k$ has a
coprime labeling by $[n]$, with its parent retained in the
exterior tree.
\end{proposition}

\begin{proof}
Use the notation $u,v$ of Lemma~\ref{lem:stump-dichotomy}.
Choose the actual branches described by $\mathbf h$, put $q$
on $v$, and label the retained branches by the blocks in
\eqref{eq:retained-core-inventory}. Their orders agree, so their
forest assertions apply to the actual branch trees.

The omitted branches have total $s-t$. The displayed interval
also has length $n-t-1-m=s-t$, so it is partitioned exactly by
their assigned intervals. Label each omitted branch using its
interval assertion. Delete $v$ and all its branches. The remaining
tree contains $u$ and has order $m$; label it by $[m]$.

The three inventories are disjoint and exhaust $[n]$. Every edge
inside a piece is valid. Every edge between pieces meets $v$,
whose label $q$ is universal. This proves the assertion.
\end{proof}

For even branches we use the stronger interval assertion
\begin{equation}\label{eq:odd-root-interval}
 \mathcal U_1(H_{I(a,j)},O^{\rm odd}(a,j)),
 \qquad O^{\rm odd}(a,j)=\{x\in I(a,j):x\text{ is odd}\}.
\end{equation}
It places the actual attachment vertex at an odd label. It is
not a claim about just one choice of root in an unrooted tree.

\begin{proposition}[Peeling even branches]
\label{prop:even-peeling}
Let $n$ be even, let $\mathbf k$ be a short-branch profile with
only even branch sizes, and let $\mathbf h\leq\mathbf k$ have
total $t\geq b$. Write $s=|\mathbf k|$ and $m=n-s-2$.
Suppose $c\geq2$ is a power of two and
$\mathcal D_n(\mathbf h;c)$ holds. Partition
\[
 \{m+1,\ldots,n-t-2\}
\]
into consecutive intervals for the omitted branches, in any
order. If each interval satisfies \eqref{eq:odd-root-interval}
and $\cP([m])$ holds, the full branch configuration has a
coprime labeling by $[n]$ with its parent removed.
Every interval used here has even length and odd starting label.
\end{proposition}

\begin{proof}
Put $c$ on $v$ and the universal label $q$ on $u$. Label the
retained branches by \eqref{eq:removed-core-inventory}, rooting
each at its neighbor of $v$. Label each omitted branch by its
assigned interval, with an odd label at that same attachment
vertex. Its edge to $v$ is coprime because $c$ is a power of two.

After deleting $u,v$ and all branches at $v$, the exterior is a
forest of order $m$. Label it by $[m]$. The intermediate interval
has length $s-t$, and the top interval has length $t+2$; these
inventories and $[m]$ partition $[n]$. Internal edges and the
edges at $v$ have been checked. Every other crossing edge meets
$u$ and is valid by universality of $q$.

Both $n$ and $s$ are even, so $m$ is even. The first interval
starts at $m+1$, and each later start is obtained by adding even
branch orders. Thus every start is odd, as claimed.
\end{proof}

\begin{figure}[htbp]
\centering
\begin{tikzpicture}[x=1cm,y=1cm,every node/.style={font=\normalfont\normalsize}]
 \node[anchor=west] at (0,1.15) {Parent retained: $m=n-s-1$};
 \draw (0,.25) rectangle (11,.85);
 \draw (3,.25)--(3,.85) (7,.25)--(7,.85);
 \node at (1.5,.55) {exterior: $m$};
 \node at (5,.55) {omitted: $s-t$};
 \node at (9,.55) {$v$ and core: $t+1$};
 \node[anchor=north] at (0,.15) {$1$};
 \node[anchor=north] at (3,.15) {$m\mid m+1$};
 \node[anchor=north] at (7,.15) {$n-t-1\mid n-t$};
 \node[anchor=north] at (11,.15) {$n$};
 \node[anchor=west] at (0,-1.15) {Parent removed: $m=n-s-2$};
 \draw (0,-2.05) rectangle (11,-1.45);
 \draw (3,-2.05)--(3,-1.45) (7,-2.05)--(7,-1.45);
 \node at (1.5,-1.75) {exterior: $m$};
 \node at (5,-1.75) {omitted: $s-t$};
 \node at (9,-1.75) {$u,v$ and core: $t+2$};
 \node[anchor=north] at (0,-2.15) {$1$};
 \node[anchor=north] at (3,-2.15) {$m\mid m+1$};
 \node[anchor=north] at (7,-2.15) {$n-t-2\mid n-t-1$};
 \node[anchor=north] at (11,-2.15) {$n$};
\end{tikzpicture}
\caption{The label inventories in the two peeling constructions.
The numbers inside the boxes are their sizes; the labels below
give the integer boundaries. Widths are schematic. An exterior
or omitted-branch interval may be empty.}
\label{fig:branch-peeling}
\end{figure}
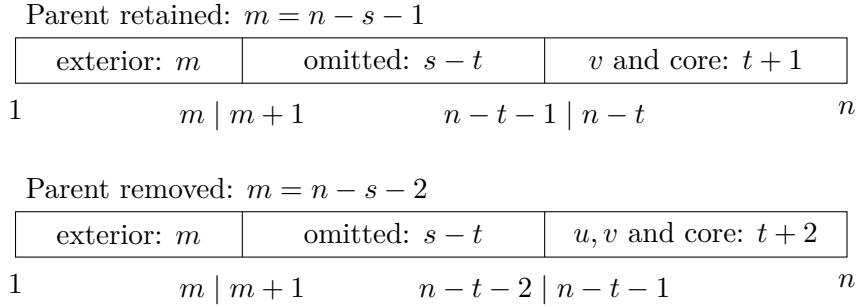

The exterior orders in these two propositions depend on the full
total $s$, not on the retained total $t$. The intermediate labels
are used by the omitted branches. They are not added to the
initial interval used for the exterior forest.

\subsubsection{Exchanging the center label}

The power of two in Proposition~\ref{prop:even-peeling} must
belong to the top interval. If it lies far below $n$, this may
require a large retained total. The next proposition permits
a label $2p$ near $n$ in that interval, and uses a power of two
only in the final exchange.

\begin{proposition}[Center exchange]\label{prop:peeling-center-exchange}
Let $n$ be even and $b\geq3$. Choose a power of two $c$ and an
odd prime $p$ such that
\begin{equation}\label{eq:peeling-transfer-arithmetic}
 c\leq n,\qquad 2p\leq n<3p,\qquad c-p\geq b.
\end{equation}
Let $\mathbf k$ be a short-branch profile with only even sizes,
and let $\mathbf h\leq\mathbf k$ have total $t\geq b$.
Suppose
\begin{equation}\label{eq:peeling-transfer-core}
 \mathcal D_n(\mathbf h;2p)\quad\text{and}\quad n-t-1>p.
\end{equation}
Put $s=|\mathbf k|$ and $m=n-s-2$. As in
Proposition~\ref{prop:even-peeling}, assign consecutive intervals
to the omitted branches. If these satisfy
\eqref{eq:odd-root-interval} and $\cP([m])$ holds, the full
configuration has a coprime labeling by $[n]$.
\end{proposition}

\begin{proof}
Put $2p$ on $v$ and the universal label $q$ on $u$. Use the
retained block labelings, the odd-root interval labelings, and
an exterior labeling by $[m]$. This gives a bijection with
coprime internal edges. Every branch root has an odd label:
for a retained branch this follows from coprimality with $2p$,
and for an omitted branch it is part of the interval assertion.
The only odd label in $[n]$ not coprime to $2p$ is $p$, since
$3p>n$. Thus an edge from $v$ to a root labeled $p$ is the only
possible failure so far.

If $m\geq c$, then $p<c\leq m$. Every branch label exceeds
$m$, so no branch root is labeled $p$. All edges at $v$ are
valid, and every remaining crossing edge meets $u$. This gives
the required labeling in this case.

Suppose now that $m<c$. The label $c$ occurs in a branch.
Indeed it is outside $[m]$, differs from $2p$, and differs
from $q$. The first distinction holds because an odd prime
cannot divide a power of two. For the second, every universal
label in $[n]$, with $n\geq4$, is odd: an even label other
than $2$ conflicts with $2$, while $2$ conflicts with $4$.
The vertex carrying $c$ is not a branch root, since all those
roots have odd labels.

Exchange $c$ and $2p$. The center now has a power-of-two label,
so all of its branch edges are coprime. We must check the
vertex in the branch whose label changes from $c$ to $2p$.
Its neighbors all lie inside that branch and have odd labels,
because its previous incident edges were coprime to $c$.
None has label $p$. If $c$ is in a retained block, every label
of that block is at least $n-t-1>p$. If $c$ is in an omitted
interval of length $j<b$, any two labels there differ by at
most $j-1\leq b-2$, whereas $c-p\geq b$.

Every odd neighbor is therefore coprime to $2p$. All other
internal edges are unchanged, and edges at $u$ remain valid
by universality of $q$. The exchange preserves the inventory
and bijectivity. This proves the assertion.
\end{proof}

The preliminary labeling in the second case need not be valid
at $v$. The proof checks those edges after the exchange; it
uses only the already valid internal branch edges to determine
the parities of the neighbors of $c$.

\subsubsection{Restricting the longest odd intervals}

An interval assertion for a large branch may be harder to
establish than the assertions for smaller branches. We can
restrict the starts needed for the largest odd branch size
without enlarging the retained-core family.

\begin{lemma}[Retaining the largest odd branches]
\label{lem:largest-odd-core}
Assume $b\geq4$, and let $j_*$ be the largest odd integer below
$b$. An ordinary short-branch profile has a minimal retained
core, of total between $b$ and $2b-2$, chosen as follows:
\begin{enumerate}
\item if at least two have order $j_*$, retain two;
\item if exactly one does, retain it and a branch of largest
      remaining order;
\item if none does, use Lemma~\ref{lem:minimal-branch-core}
      with threshold $b$ and the required odd branch, if any.
\end{enumerate}
In each case minimality includes the odd-branch requirement
when $n$ is even.
\end{lemma}

\begin{proof}
Put $d=b-j_*\in\{1,2\}$. In the first case the total is
$2j_*\in[b,2b-2]$, and deleting either branch leaves a total
below $b$. Both branches have odd order.

In the second case another branch exists, since $j_*<b$ and
the full total is at least $b$. Let $h$ be its largest possible
order. If $d=1$, a branch of order $1$ together with the
$j_*$-branch would give the forbidden sum $b$, so $h\geq2$.
If $d=2$, a branch of order $2$ is forbidden. If every other
branch had order $1$, there would be at least two of them,
again giving the forbidden sum with $j_*$. Thus $h\geq3$.
In both cases $h\geq d+1$, and hence
\[
 b+1\leq j_*+h\leq2b-2.
\]
Deleting either retained branch leaves total below $b$, and
the $j_*$-branch supplies the required odd order.
The last case is exactly Lemma~\ref{lem:minimal-branch-core}.
All selected profiles inherit the absence of a sum $b$.
\end{proof}

\begin{corollary}[The required interval starts]
\label{cor:largest-odd-intervals}
Use Lemma~\ref{lem:largest-odd-core} in ordinary peeling and
place the omitted $j_*$-branches last. Then the plain interval
assertions of length $j_*$ are needed only at starts in
\begin{equation}\label{eq:largest-odd-starts}
 \mathcal A_{n,j_*}
   =\{n-(2+r)j_*:r\geq1,\ n-(2+r)j_*>0\}.
\end{equation}
All other plain interval assertions may retain their full
range of starts. The odd-root assertions for even-length
intervals are unchanged.
\end{corollary}

\begin{proof}
If any $j_*$-branch is omitted, at least three were present
and the retained core consists of exactly two. The peeled
intervals end at $n-2j_*-1$. Since the $j_*$-intervals come
last, counting them backwards gives starts
$n-3j_*,n-4j_*,\ldots$. Their actual starts are positive,
so all occur in \eqref{eq:largest-odd-starts}.
If at most one $j_*$-branch was present, it was retained and
no interval of this length is used. Profiles with only even
branches contain no $j_*$-branch at all.
\end{proof}

For example, $n=1359$ and $b=32$ give $j_*=31$ and the $41$
starts $26,57,\ldots,1266$. This calculation restricts an
arithmetic checklist; it does not assert that the interval
or retained-block assertions have already been proved.

\subsubsection{The finite-prefix step}

We can now state the induction step. The case lists are defined
mathematically, so their completeness can be checked separately
from the arithmetic assertions assigned to them.

Let $\mathcal O_{n,b}$ consist of all short-branch profiles
$\mathbf h$ with
\[
 b\leq t=|\mathbf h|\leq\min(2b-2,n-2)
\]
that satisfy \eqref{eq:core-minimality} with $A=b$ and
$\varepsilon=0$ for odd $n$, or $\varepsilon=1$ for even $n$.
In the latter case also require $o(\mathbf h)>0$.

For even $n$, put $E=2\floor{(b-1)/2}$. When $E\geq2$,
choose a nominal center label $a$ and define
\begin{equation}\label{eq:even-core-threshold}
 B=\max(b,n-a-1),\qquad H=B+E-1.
\end{equation}
Let $\mathcal E^-_{n,b}(a)$ contain all short-branch profiles
with only even sizes and total $b\leq t<B$. Let
$\mathcal E^+_{n,b}(a)$ contain all such profiles with
$B\leq t\leq\min(H,n-2)$ and
\[
 h_j>0\quad\Longrightarrow\quad t-j<B.
\]
The restriction $t\leq n-2$ is part of the short-branch
definition in both families. The low family is not filtered
by minimality.

For interval assertions, use all pairs $(j,a_0)$ with
$1\leq j<b$ and $1\leq a_0\leq n-j+1$, except that, when
$b\geq4$ and $j=j_*$, use only $a_0\in\mathcal A_{n,j_*}$.
Denote this list by $\mathcal I_{n,b}$.

\begin{theorem}[A step in the finite prefix]
\label{thm:peeled-prefix-step}
Put $\ell=n-b-1$ and suppose $\cP([m])$ holds for every
$0\leq m\leq\ell$. Assume the following assertions.
\begin{enumerate}
\item There is a label $q_0$ universal in $[n]$ for which
\[
 \{\ell+1,\ldots,n\}=\{q_0\}\sqcup B_0,
 \qquad \cP(B_0).
\]
\item $\mathcal R_n(\mathbf h)$ holds for every
      $\mathbf h\in\mathcal O_{n,b}$.
\item $\cP(I(a_0,j))$ holds for every
      $(j,a_0)\in\mathcal I_{n,b}$. If $n$ is even,
      \eqref{eq:odd-root-interval} also holds for every even
      $j<b$ and every odd start $1\leq a_0\leq n-j+1$.
\end{enumerate}
If $n$ is even and $E\geq2$, make one of the following choices:
\begin{enumerate}
\item[(i)] $a=c$, where $2\leq c\leq n$ is a power of two;
\item[(ii)] $a=2p$, where a power of two $c$ and an odd prime
      $p$ satisfy \eqref{eq:peeling-transfer-arithmetic}, and
      $n-H-1>p$.
\end{enumerate}
For this choice require, in addition,
\[
 \begin{array}{ll}
 \mathcal D_n(\mathbf h;a(\mathbf h))
      \text{ for some }a(\mathbf h),
       &\mathbf h\in\mathcal E^-_{n,b}(a),\\[2pt]
 \mathcal D_n(\mathbf h;a),
       &\mathbf h\in\mathcal E^+_{n,b}(a).
 \end{array}
\]
There are no even-family requirements when $n$ is odd or $E=0$.
Then $\cP([n])$ holds.
\end{theorem}

\begin{proof}
By Lemma~\ref{lem:forest-completion}, it suffices to label an
arbitrary tree of order $n$. Apply
Lemma~\ref{lem:stump-dichotomy}. In its detached-set case,
put $q_0$ on the contact vertex, use $B_0$ for the detached
forest of order $b$, and use $[\ell]$ for the remaining
forest. The first assertion and universality of $q_0$ check
all edges.

In the branch case, let $\mathbf k$ be the full profile and
write $s=|\mathbf k|$. Suppose first that it is ordinary.
Choose a minimal core by Lemma~\ref{lem:largest-odd-core}
when $b\geq4$, or by Lemma~\ref{lem:minimal-branch-core}
otherwise. Its total is at most $2b-2$ and at most $s\leq n-2$,
so it belongs to $\mathcal O_{n,b}$. Place omitted branches
of order $j_*$ last when that size is defined. The required
plain intervals are all in $\mathcal I_{n,b}$, by
Corollary~\ref{cor:largest-odd-intervals}. The actual exterior
order is $m=n-s-1\leq\ell$. Thus
Proposition~\ref{prop:ordinary-peeling} applies.

It remains to treat an even $n$ and a profile consisting only
of even branches. Such a profile requires $E\geq2$. If $s<B$,
the full profile itself lies in $\mathcal E^-_{n,b}(a)$.
Use its asserted $\mathcal D_n$ partition and retain every
branch. Its center may depend on the profile: no branch is
peeled. Put its universal label on $u$ and label the exterior
forest by $[n-s-2]$. This is the remove-parent construction of
Proposition~\ref{prop:rooted-branch-extension}.

If $s\geq B$, apply Lemma~\ref{lem:minimal-branch-core} with
threshold $B$, no odd-branch requirement, and the bound $E$
on branch size. The selected total satisfies $B\leq t\leq H$;
its minimality puts the selected profile in
$\mathcal E^+_{n,b}(a)$. In particular
$t\geq n-a-1$, so the nominal center belongs to the required
top interval. The exterior has order $m=n-s-2\leq\ell$.
All omitted intervals have even length and odd start, as
proved in Proposition~\ref{prop:even-peeling}.

For choice (i), that proposition finishes the labeling.
For choice (ii), the bound $t\leq H$ gives
$n-t-1\geq n-H-1>p$, so
Proposition~\ref{prop:peeling-center-exchange} applies.
Every branch profile has now been treated, and the detached-set
case was handled at the start. This proves $\cP([n])$.
\end{proof}

\begin{corollary}[Induction over a finite range]
\label{cor:finite-prefix-induction}
Suppose $\cP([m])$ is established for $0\leq m\leq N_0$.
For each integer $N_0<n\leq N$, choose $b_n$ with
$2\leq b_n\leq n-2$ and supply all the arithmetic assertions
of Theorem~\ref{thm:peeled-prefix-step}, using its exact case
families. Then $\cP([n])$ holds for every $0\leq n\leq N$.
\end{corollary}

\begin{proof}
Proceed in increasing order of $n$. Its required initial
intervals have orders at most $n-b_n-1<n$, so they are covered
by the induction hypothesis. Theorem~\ref{thm:peeled-prefix-step}
establishes the next order. Induction gives the whole range.
\end{proof}

All case families in the theorem are finite and explicit.
To enumerate them, use the bounded subset-sum recurrence after
Proposition~\ref{prop:rooted-branch-extension}, restrict the
total to the displayed bounds, and test
\eqref{eq:core-minimality}. Retain every low even profile.
The interval list is given by integer endpoints and, for its
one restricted size, \eqref{eq:largest-odd-starts}.
A certificate must give exactly these cases, together with
the detached-set partition, and bind every block to its full
arithmetic assertion. Equality of the expected and supplied
case sets checks coverage; the soundness results check the
assertions themselves.

The remaining block theorems and their complete arithmetic
certificates are needed to apply this induction to a specified
numerical prefix.

\subsection{Bipartite hosts for finite blocks}\label{sec:bipartite-hosts}

The branch reductions require coprime labelings of intervals with a
prescribed root. Checking every rooted tree separately becomes costly
as the interval grows. We instead assign universal labels to the root
and a small set of further vertices. The remaining forest has a
bipartition with a fixed small side and a degree bound on that side.
Its nonleaf structure can then be enumerated independently of the
number of leaves.

We first prove the structural reduction. We then describe complete
certificates for the bipartite hosts and apply one to a rooted block
of order $28$. The certificate checks every residual tree described
by the reduction; it makes no assumption on the original tree's
maximum degree.

\subsubsection{A small side with bounded degrees}

Choose a bipartition $(C,I)$ of a forest $F$ by putting the smaller
side of each nontrivial component in $C$, and every isolated vertex
in $I$. Thus $C$ contains no isolated vertices and
$|C|\leq\floor{|V(F)|/2}$. The next lemma quantifies how many
vertices of $I$ become isolated when suitable vertices of $C$ are
deleted.

\begin{samepage}
\begin{lemma}[One-side deletion]\label{lem:one-side-isolation}
Let $(C,I)$ be a bipartition of a forest of order $N$, with no
isolated vertex in $C$, and put $c=|C|$. Let $1\leq k<c$ and
$D\geq1$. Suppose that
\[
 H=\{v\in C:d_F(v)\geq D+1\},\qquad h=|H|\leq k.
\]
There is a set $S\subset C$, of size $k$ and containing $H$,
such that $\Delta_{F-S}(C\setminus S)\leq D$ and the number $Z$
of isolated vertices of $I$ in $F-S$ satisfies
\begin{equation}\label{eq:one-side-isolation}
 Z\geq
 \frac{(k-h)(N+1-2c)+(c-k)\max\{0,Dh+2-c\}}{c-h}.
\end{equation}
\end{lemma}
\end{samepage}

\begin{proof}
Write $z,L,q$ for the numbers of vertices in $I$ of degrees
$0,1$, and at least $2$, respectively. Let $m_1$ be the number
of nontrivial components. The isolated vertices are precisely
the $z$ vertices just counted, so
\[
 e(F)=N-z-m_1,\qquad e(F)\geq L+2q,
 \qquad N=c+z+L+q.
\]
Consequently $q\leq c-m_1\leq c-1$ and
\begin{equation}\label{eq:isolation-leaf-count}
 z+L\geq N+1-2c.
\end{equation}
Here $m_1\geq1$ because $c>k\geq1$. Counting the isolated
vertices as separate components is essential to this inequality.

Let $a$ be the number of leaves in $I$ adjacent to $H$. If
$h>0$, the subgraph between $H$ and the $q$ nonleaf vertices
of $I$ is a forest, and has at most $h+q-1$ edges. It follows that
\[
 (D+1)h\leq a+h+q-1,\qquad
 a\geq Dh-q+1\geq Dh+2-c.
\]
For $h=0$ we have $a=0$ and $c\geq2$, so in either case
$a\geq\max\{0,Dh+2-c\}$.

After choosing $H$, add to $S$ the $k-h$ vertices of
$C\setminus H$ with the most leaf neighbors in $I$. Their
leaf counts sum to at least $(k-h)(L-a)/(c-h)$, by averaging
the largest $k-h$ of the $c-h$ counts. Thus
\begin{align*}
 Z&\geq z+a+\frac{(k-h)(L-a)}{c-h}\\
  &=\frac{(k-h)(z+L)+(c-k)(z+a)}{c-h}.
\end{align*}
Both coefficients in the numerator are nonnegative. Substituting
\eqref{eq:isolation-leaf-count} and the bound for $a$ proves
\eqref{eq:one-side-isolation}. Every vertex of $C\setminus S$
has degree at most $D$, since all vertices of larger degree
belong to $H$.
\end{proof}

\begin{theorem}[Uniform small sides]\label{thm:uniform-small-sides}
Let $(b,D,k_0)$ be one of
\[
 (8,5,3),\qquad (9,6,4),\qquad (10,6,4).
\]
For every $k\geq k_0$ and every forest $F$ of order
$2k+2b+1$, there is a set $S$ of $k$ vertices such that
$F-S$ has a bipartition $(B,A)$ satisfying
\[
 |B|=b,\qquad |A|=k+b+1,\qquad
 \Delta_{F-S}(B)\leq D.
\]
\end{theorem}

\begin{proof}
Choose $(C,I)$ as above and put $c=|C|$. If $c\leq k$,
delete all of $C$ and any $k-c$ other vertices. The remaining
graph is independent, and can be divided into parts of the
required sizes.

Suppose $c=k+r$, where $1\leq r\leq b$. At most $k$
vertices of $C$ have degree at least $D+1$. Otherwise
\[
 e(F)\geq(D+1)(k+1)>2k+2b=|V(F)|-1,
\]
where the strict inequality holds in each of the three
parameter ranges. Apply Lemma~\ref{lem:one-side-isolation}
and write $s=k-h$. It gives
\begin{equation}\label{eq:uniform-isolation-bound}
 Z\geq
 \frac{s(2b+2-2r)+r\max\{0,(D-1)k+2-r-Ds\}}{s+r}.
\end{equation}
There are $r$ vertices in $C\setminus S$. It suffices to have
$b-r$ isolated vertices in $I$, which can then be moved to
that side without changing any degree.

For $r=b$ nothing remains to show. For $1\leq r<b$, the
right side of \eqref{eq:uniform-isolation-bound} is greater
than $b-1-r$ if
\[
 f(s)=s(b+3-r)+r\max\{0,(D-1)k+2-r-Ds\}
                     -r(b-1-r)>0.
\]
We verify this for every real $s\geq0$. If $r=1$, both slopes
of $f$ are positive in the stated ranges, and
$f(0)=(D-1)k+3-b>0$. If $r\geq2$, the slope before the
breakpoint is $b+3-(D+1)r<0$, and that after the breakpoint
is $b+3-r>0$. The minimum is therefore attained at
\[
 s_0=\frac{(D-1)k+2-r}{D}>0.
\]
At this point
\[
 Df(s_0)=(b+3-r)((D-1)k+2-r)-Dr(b-1-r).
\]
This expression increases with $k$. Its values at $k_0$ are
the positive expressions in Table~\ref{tab:small-side-inequalities}.
Thus $Z>b-1-r$. Since $Z$ is an integer, $Z\geq b-r$.
Move that many isolated vertices to $C\setminus S$ and call
the resulting side $B$. The other side has
$2k+2b+1-k-b=k+b+1$ vertices, as required.
\end{proof}

\begin{table}[htbp]
\centering
\renewcommand{\arraystretch}{1.45}
\begin{tabular}{c c c c}
\hline
$b$ & $D$ & $k_0$ & Lower bound for $Df(s_0)$ \\
\hline
$8$ & $5$ & $3$ & $6(r-5)^2+4$ \\
$9$ & $6$ & $4$ & $7(r-41/7)^2+167/7$ \\
$10$ & $6$ & $4$ & $7(r-89/14)^2+87/28$ \\
\hline
\end{tabular}
\par\vspace{.5\baselineskip}
\caption{Exact positive bounds used in the small-side reduction,
for $2\leq r<b$.}
\label{tab:small-side-inequalities}
\end{table}

\subsubsection{From bipartite hosts to rooted blocks}

Let $X$ be a bipartite graph with specified sides of sizes
$b,a$. We say that $X$ \emph{covers $(b,a,D)$} if every forest
with specified bipartition $(B,A)$, $|B|=b$, $|A|=a$, and
$\Delta(B)\leq D$ has a side-preserving edge-preserving
bijection to $X$. In particular, the bipartition is part of the
input; it is not chosen anew when testing a witness.

\begin{proposition}[Adding universal vertices]
\label{prop:bipartite-universal-extension}
Choose $b,D,k$ as in Theorem~\ref{thm:uniform-small-sides},
and suppose $X$ covers $(b,k+b+1,D)$. Adjoin $k$ vertices,
each adjacent to every other vertex, to obtain a graph $H$.
Then every forest of order $2k+2b+1$ embeds bijectively in $H$.

If instead $k+1$ such vertices are adjoined, every forest of
order $2k+2b+2$ with a specified vertex $v$ has an embedding
that sends $v$ to any specified one of the new vertices.
\end{proposition}

\begin{proof}
In the first assertion use Theorem~\ref{thm:uniform-small-sides}
to choose $S$, and embed the residual forest in $X$. Send
$S$ bijectively to the new vertices. Edges inside the residual
forest are preserved by its embedding; every other edge has
at least one end sent to a universal vertex.

For the rooted assertion, delete $v$ first and apply the same
construction to the remaining forest of order $2k+2b+1$.
Send $v$ to the specified new vertex and $S$ to the other $k$.
The edges incident to $v$ are valid by universality. The two
constructions use every host vertex exactly once.
\end{proof}

It remains to prove the covering assertions for the hosts
actually used. The next reduction retains all vertices of the
small side, which we color black. Vertices on the other side
are colored white.

\begin{lemma}[Kernel and leaves]\label{lem:bipartite-kernel}
Let $T$ be a tree with a specified bipartition into $b\geq2$
black vertices and $a\geq1$ white vertices, and with black
degrees at most $D$. Delete
all white leaves. The resulting graph $K$ is a tree containing
all $b$ black vertices and $q$ white vertices, where
\[
 1\leq q\leq\min\{a,b-1\}.
\]
Every white vertex of $K$ has degree at least $2$, and every
black vertex has degree between $1$ and $D$. The tree $T$
is recovered by attaching $\ell_v$ white leaves to each black
vertex $v$, with
\begin{equation}\label{eq:kernel-leaf-vectors}
 \sum_{v\in B}\ell_v=a-q,
 \qquad 0\leq\ell_v\leq D-d_K(v).
\end{equation}
Conversely every such kernel and vector give a tree of the
specified orders and degree bound.
\end{lemma}

\begin{proof}
Removing leaves preserves connectivity among the remaining
vertices. A black vertex could become isolated only if all
its neighbors were white leaves, in which case its component
in $T$ would contain no other black vertex. This contradicts
$b\geq2$ and connectivity. Hence $K$ is a tree, every black
vertex still has a neighbor, and $q\geq1$. No edge at a
remaining white vertex was removed, so its degree is at least
$2$. Counting edges from the white side gives
\[
 2q\leq e(K)=b+q-1,
\]
and thus $q\leq b-1$. The deleted leaves give exactly
\eqref{eq:kernel-leaf-vectors}. Conversely, attaching those
leaves preserves the tree property, adds precisely $a-q$
white vertices, and raises no black degree above $D$.
\end{proof}

The embedding of $K$ is allowed to depend on its leaf vector.
An embedding that works for one vector need not have enough
unused neighbors for another. Our certificates therefore
supply a full tree embedding for each required vector.

\subsubsection{Complete enumeration of the kernels}

\begin{proposition}[Star partitions]\label{prop:kernel-star-partitions}
Every kernel in Lemma~\ref{lem:bipartite-kernel} is obtained
by the following procedure, up to a color-preserving
isomorphism. Take a tree $R$ on $b$ black vertices and partition
its edges into nonempty blocks, each having a common endpoint.
For each block introduce a white vertex adjacent to all
endpoints of its edges, and discard the edges of $R$.
Retain only the graphs with at most $a$ white vertices and
black degrees at most $D$.

Conversely, every retained graph is a tree with the stated
kernel conditions. Consequently it is enough to use one representative
of each unrooted tree $R$ and all its edge partitions.
\end{proposition}

\begin{proof}
For each white vertex of a kernel choose one of its black
neighbors, and contract the edge between them. These edges
form a forest, each of whose components contains exactly one
black vertex. The quotient is therefore a tree $R$ on the
$b$ black vertices. The uncontracted edges originally incident
to a given white vertex become a nonempty star centered at
its chosen black neighbor. These stars partition $E(R)$.
Their endpoint sets are exactly the original white
neighborhoods, so the procedure reconstructs the kernel.

Conversely, a block of $t$ edges with a common endpoint has
$t+1$ distinct endpoints. Replacing it by a new white vertex
and its incident edges preserves connectivity. If there are
$q$ blocks, the resulting graph has $b+q$ vertices and
\[
 \sum_{\text{blocks}}(t+1)=(b-1)+q
\]
edges, and is therefore a tree. Each white degree is at least
$2$; every black degree is positive. The two stated filters
give the remaining conditions. Finally an isomorphism of
$R$ carries all its edge partitions with it, so passing to
unrooted representatives loses no kernel.
\end{proof}

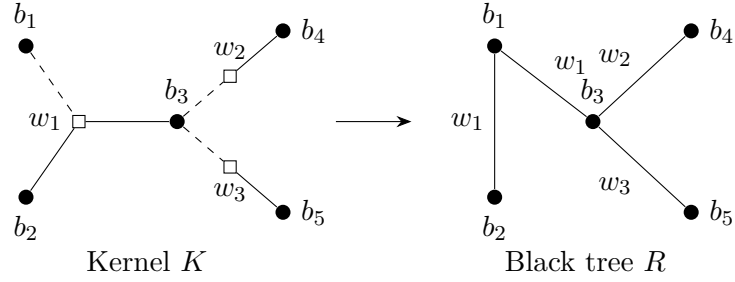
\begin{figure}[htbp]
\centering
\begin{tikzpicture}[x=1cm,y=1cm,
  every node/.style={font=\normalfont\normalsize},
  blackv/.style={circle,fill=black,inner sep=2pt},
  whitev/.style={rectangle,draw,fill=white,inner sep=2.3pt}]
\node[blackv,label=above:$b_1$] (b1) at (0,2) {};
\node[blackv,label=below:$b_2$] (b2) at (0,0) {};
\node[blackv,label=above:$b_3$] (b3) at (2,1) {};
\node[blackv,label=right:$b_4$] (b4) at (3.4,2.2) {};
\node[blackv,label=right:$b_5$] (b5) at (3.4,-.2) {};
\node[whitev,label=left:$w_1$] (w1) at (.7,1) {};
\node[whitev,label=above:$w_2$] (w2) at (2.7,1.6) {};
\node[whitev,label=below:$w_3$] (w3) at (2.7,.4) {};
\draw[dashed] (b1)--(w1) (b3)--(w2) (b3)--(w3);
\draw (b2)--(w1)--(b3) (w2)--(b4) (w3)--(b5);
\node at (1.6,-.85) {Kernel $K$};
\draw[-{Stealth[length=2mm]}] (4.1,1)--(5.1,1);
\begin{scope}[xshift=6.2cm]
\node[blackv,label=above:$b_1$] (c1) at (0,2) {};
\node[blackv,label=below:$b_2$] (c2) at (0,0) {};
\node[blackv,label=above:$b_3$] (c3) at (1.3,1) {};
\node[blackv,label=right:$b_4$] (c4) at (2.6,2.2) {};
\node[blackv,label=right:$b_5$] (c5) at (2.6,-.2) {};
\draw (c1)--node[left] {$w_1$}(c2)
      (c1)--node[above right] {$w_1$}(c3)
      (c3)--node[above left] {$w_2$}(c4)
      (c3)--node[below left] {$w_3$}(c5);
\node at (1.2,-.85) {Black tree $R$};
\end{scope}
\end{tikzpicture}
\caption{Contract the dashed edges, one at each white vertex.
The labels on the edges of $R$ identify their blocks.
Replacing each block by a white vertex adjacent to its
endpoints recovers $K$.}
\label{fig:kernel-stars}
\end{figure}

The enumeration in Lemma~\ref{lem:tree-code-enumeration}
supplies the trees $R$. An alternative, used by the checker
below, starts with the one-vertex tree and repeatedly attaches
a leaf at every possible vertex, retaining one representative
of each resulting isomorphism class. It is complete because
every tree of order at least two has a leaf, whose deletion
gives a tree of the preceding order.

For the colored kernels, put the symbol $B$ or $W$ immediately
after each opening parenthesis in the rooted code. Sort the
child codes as before, and take the least code over all roots.
The induction proving Lemma~\ref{lem:tree-code-enumeration}
then proves that equality of these codes is exactly
color-preserving isomorphism. Thus repeated kernels can be
removed without discarding any required case.

There is a further reduction when black vertices have equal
neighborhoods. It is useful for larger blocks, though the
specific table below retains every feasible leaf vector.

\begin{lemma}[Permuting equal neighborhoods]\label{lem:kernel-twin-vectors}
Partition the black vertices of a kernel into classes with
equal white neighborhoods. The vectors in
\eqref{eq:kernel-leaf-vectors} need only be checked when their
coordinates are nondecreasing in a fixed order within each
class.

For a class of size $t$ whose vertices have degree $d$, put
\[
 P_{t,D-d}(z)=
 \sum_{0\leq x_1\leq\cdots\leq x_t\leq D-d}
                          z^{x_1+\cdots+x_t}.
\]
The number of vectors that remain is the coefficient of
$z^{a-q}$ in the product of these polynomials over the
classes.
\end{lemma}

\begin{proof}
Every permutation within a class, fixing all other vertices,
is an automorphism of the kernel. All vertices in that class
have the same degree, and therefore the same upper bound
$D-d$ for their leaf coordinates. Sorting the coordinates
preserves admissibility. Permuting the black vertices and
their attached leaves gives an isomorphism of the enlarged
trees, so an embedding for the sorted vector gives one for
the original vector.

The polynomial for one class lists its allowed nondecreasing
coordinate sequences, with exponent equal to their total.
Multiplication makes the choices independently for the
classes, and extraction of the coefficient enforces the
total $a-q$. No quotient by other kernel automorphisms is
asserted or needed.
\end{proof}

\begin{proposition}[Bipartite witness criterion]
\label{prop:bipartite-witnesses}
Let $b\geq2$, $a\geq1$, and $Db\geq a+b-1$. Let $X$ be a
bipartite graph with sides of sizes $b,a$. Enumerate the
kernels by Proposition~\ref{prop:kernel-star-partitions}.
For each kernel and every vector in
\eqref{eq:kernel-leaf-vectors}, supply a side-preserving
bijection of the enlarged tree to $X$. If every edge maps
to an edge of $X$, then $X$ covers $(b,a,D)$.

It is enough to supply the sorted vectors of
Lemma~\ref{lem:kernel-twin-vectors}. In this case, completeness
of a strictly lexicographically increasing list of admissible
sorted vectors is certified by equality of its length with
the coefficient in that lemma.
\end{proposition}

\begin{proof}
Lemma~\ref{lem:bipartite-kernel} describes every tree with the
specified sides and degree bound. The star partitions contain
its kernel type, and the supplied vectors contain its leaf
data, up to the permitted permutations. Pulling back the
checked witness along this isomorphism embeds the tree in
$X$.

Now let $F$ be any forest with the specified sides and black
degree bound. Lemma~\ref{lem:bipartite-completion}, with the
black side as its degree-constrained side, completes $F$ to
a tree on the same bipartition: both sides are nonempty and
$Db\geq a+b-1$. An embedding of that tree also embeds $F$.
For the final assertion, strict ordering excludes duplicate
vectors, admissibility includes the listed vectors in the
required finite set, and equality of cardinalities proves
that none is omitted.
\end{proof}

The distinction between the two sides also gives a useful
extension that requires no new tree table.

\begin{lemma}[Adding a white vertex]\label{lem:white-host-extension}
Suppose $X$ covers $(b,a,D)$, where $b\geq2$, $a\geq b-1$,
and $Db\geq b+a$. Adjoin a white vertex adjacent to every
black vertex of $X$. The resulting host covers $(b,a+1,D)$.
\end{lemma}

\begin{proof}
Apply Lemma~\ref{lem:bipartite-completion} to complete the
input forest to a tree on the same sides. The capacity
condition is $Db\geq b+(a+1)-1$. This tree has a
white leaf. Otherwise its number of edges would be at least
$2(a+1)>b+a$, since $a\geq b-1$, whereas a tree of this
order has $b+a$ edges. Delete a white leaf and embed the
remaining tree in $X$. Send the leaf to the new white vertex.
Its unique neighbor is black, so its edge is preserved.
Restricting to the input forest proves the assertion.
\end{proof}

\subsubsection{A rooted interval of order twenty-eight}

We now give one explicit arithmetic application. Let $X$ be
the bipartite coprimality graph with black labels
\begin{equation}\label{eq:twentyeight-black-labels}
 129,133,135,141,143,145,147,153,155
\end{equation}
and white labels $130,132,\ldots,156$. An edge joins two
labels on opposite sides exactly when their gcd is $1$.
Table~\ref{tab:rooted-twentyeight-certificate} records the
complete witness collection supplied with this manuscript.

\begin{table}[htbp]
\centering
\renewcommand{\arraystretch}{1.25}
\begin{tabular}{l r}
\hline
Enumerated or checked objects & Number \\
\hline
Unrooted trees on nine black vertices & $47$ \\
Partitions of their eight edge positions & $4\,140$ \\
Distinct kernels with black degree at most six & $493$ \\
Feasible kernel--leaf-vector pairs & $8\,205\,524$ \\
Edges checked in their full embeddings & $180\,521\,528$ \\
\hline
\end{tabular}
\par\vspace{.5\baselineskip}
\caption{The exhaustive certificate for the $9$-by-$14$
bipartite host. The $4\,140$ partitions are tested for each
black tree; only partitions into stars contribute kernels.}
\label{tab:rooted-twentyeight-certificate}
\end{table}

\begin{theorem}[A rooted block of order twenty-eight]
\label{thm:rooted-twentyeight-block}
Every forest $F$ on $28$ vertices with a specified vertex
$v$ has a coprime labeling by $\{129,\ldots,156\}$ in which
$v$ receives $131$.
\end{theorem}

\begin{proof}
The host $X$ covers $(9,14,6)$ by the exhaustive certificate
in Table~\ref{tab:rooted-twentyeight-certificate} and
Proposition~\ref{prop:bipartite-witnesses}. We specify the
checks that establish this finite assertion.

The checker regenerates all $47$ black trees by leaf
augmentation. It generates each partition of eight positions
by successively putting the next position into an existing
block or into a new last block. Every partition appears
exactly once when its blocks are ordered by their least
elements. The star construction and the black degree filter
give exactly $493$ colored kernel codes. The supplied kernel
list is required to contain each of these codes exactly once.

For each kernel, with black vertices indexed first, the
checker enumerates all nonnegative leaf vectors of total
$14-q$, in lexicographic order, and retains precisely those
satisfying \eqref{eq:kernel-leaf-vectors} with $D=6$.
Each vector has a supplied record of $23$ host indices.
The record lists the images of the kernel vertices first,
then those of the attached leaves in black-vertex order.
The first nine entries must permute the black host indices,
and the remaining entries must permute the white indices.
The checker reconstructs the arithmetic host by exact gcd
calculations, then checks all kernel edges and all added
leaf edges. Each enlarged tree has $22$ edges. The record
counts and boundaries must agree with the independently
enumerated vectors, with no missing or trailing records.
All $8\,205\,524$ records pass these checks. Finally,
$6\cdot9\geq9+14-1$, so the witness criterion applies
also to forests.

The five labels
\[
 131,137,139,149,151
\]
are prime and exceed $156/2$. They are therefore universal
in $\{129,\ldots,156\}$. Together with the two sides of
$X$ they partition that interval. Apply
Proposition~\ref{prop:bipartite-universal-extension} with
$b=9$, $D=6$, and $k=4$. Its rooted assertion has order
$2k+2b+2=28$, and the specified universal label may be
chosen to be $131$. This gives the required labeling.
\end{proof}

The supplement contains the kernel list, every embedding
record, and a standalone checker using only exact integer
operations and the Python standard library. Compression
changes only the storage of the records. The checks above,
rather than a file hash or a search success rate, establish
the covering assertion. The universal small-side theorem
supplies the passage from arbitrary input forests to that
finite assertion.

\subsection{Reserved vertices in finite blocks}\label{sec:reserved-blocks}

A universal label can be placed on a deleted vertex without
restricting its neighbors. Some useful intervals have too few
such labels, but contain a further label with only one possible
conflict. We will place its forbidden partner at a vertex
isolated in the residual forest. Controlling that vertex's
degree before deletion leaves a nonneighbor among the deleted
vertices, where the nearly universal label can be placed.

This section gives three reductions for forests of order $25$.
Together with complete bipartite certificates, they yield
rooted blocks of order $26$. We also verify one family of
arithmetic transfers over an explicit finite range.

\subsubsection{The integer deletion bounds}

We refine the averaging step in
Lemma~\ref{lem:one-side-isolation}. The improvement uses only
that the leaf counts are integers.

\begin{lemma}[Largest integer weights]\label{lem:largest-integer-weights}
Let $r\geq1$, $0\leq t\leq r$, and let $r$ nonnegative
integer weights have total $T$. The sum of their $t$ largest
values is at least
\begin{equation}\label{eq:balanced-integer-weights}
 M(t,r,T)=t\floor{T/r}+\min\{t,T-r\floor{T/r}\}.
\end{equation}
\end{lemma}

\begin{proof}
If two weights differ by at least two, decrease a largest
weight by one and increase a smallest weight by one. The
sum of the $t$ largest cannot increase. Indeed, the only
$t$-element sets whose sums increase contain the smallest
entry and omit the largest; replacing the former by the
latter gave a still larger sum before the change.

Each change strictly decreases the sum of the squares, so
the process terminates with all weights differing by at most
one. Write $T=rQ+R$, where $0\leq R<r$. The final weights
are $R$ copies of $Q+1$ and $r-R$ copies of $Q$. Their $t$
largest sum to $tQ+\min\{t,R\}$, as claimed.
\end{proof}

The following finite calculation isolates the cases in which
leaf counts alone do not suffice. We give every exceptional
tuple; the next proof resolves them using degree-two vertices.

\begin{lemma}[Finite deletion calculation]\label{lem:finite-deletion-calculation}
Consider the following three choices, denoted by
$\mathrm I$, $\mathrm {II}$, and $\mathrm {III}$:
\begin{equation}\label{eq:twentyfive-parameters}
 (k,D,\varepsilon)=(4,4,1),\quad(3,6,0),\quad(3,8,1).
\end{equation}
Let $c,h,z,L,a,q$ be integers satisfying
\begin{gather}
 k+1\leq c\leq12,\qquad
 0\leq h\leq\min\{k,\floor{(24-c)/D}\},\label{eq:deletion-parameter-domain}\\
 z,L\geq0,\quad z+L\leq25-c,\quad
 q=25-c-z-L,\quad 0\leq a\leq L,\notag\\
 \max\{c+Dh,L+2q\}\leq24-z,\notag\\
 a=0\ \text{if }h=0,\qquad
 (D+1)h-a\leq h+q-1\ \text{if }h>0.\notag
\end{gather}
Put
\begin{equation}\label{eq:integer-isolation-bound}
 \Lambda=z+a+M(k-h,c-h,L-a).
\end{equation}
Then $\Lambda\geq12-c+\varepsilon$ except for the tuples
in Table~\ref{tab:deletion-exceptions}.
\end{lemma}

\begin{proof}
This is an exhaustive calculation on bounded integers. For
each choice in \eqref{eq:twentyfive-parameters}, enumerate
$c$ and $h$ in the displayed ranges. Enumerate
$0\leq z\leq25-c$ and $0\leq L\leq25-c-z$, set
$q=25-c-z-L$, and discard a tuple if
$\max\{c+Dh,L+2q\}>24-z$. Finally enumerate $a$ from
$0$ to $L$ and impose the last line of
\eqref{eq:deletion-parameter-domain}. These nested ranges
include every permitted tuple exactly once.

There remain respectively $9\,181$, $6\,142$, and $4\,010$
tuples. Evaluate \eqref{eq:integer-isolation-bound} by the
quotient and remainder of $L-a$ upon division by $c-h$.
The complete lists for which the required inequality fails
are precisely those in Table~\ref{tab:deletion-exceptions}.
The accompanying checker performs these integer loops and
requires equality with all eight displayed rows. No tree
enumeration or numerical approximation enters this calculation.
\end{proof}

\begin{table}[htbp]
\centering
\renewcommand{\arraystretch}{1.2}
\begin{tabular}{c r r r r r r r r}
\hline
Case & $c$ & $h$ & $z$ & $L$ & $a$ & $q$ & $\Lambda$ & Required \\
\hline
$\mathrm I$   & 9 & 1 & 0 & 8 & 0 & 8 & 3 & 4 \\
$\mathrm I$   & 9 & 2 & 0 & 8 & 1 & 8 & 3 & 4 \\
$\mathrm I$   &10 & 2 & 0 & 6 & 0 & 9 & 2 & 3 \\
$\mathrm {II}$& 9 & 1 & 0 & 8 & 0 & 8 & 2 & 3 \\
$\mathrm {III}$&9 & 0 & 0 & 8 & 0 & 8 & 3 & 4 \\
$\mathrm {III}$&9 & 0 & 0 & 9 & 0 & 7 & 3 & 4 \\
$\mathrm {III}$&9 & 1 & 0 & 8 & 1 & 8 & 3 & 4 \\
$\mathrm {III}$&10& 1 & 0 & 6 & 0 & 9 & 2 & 3 \\
\hline
\end{tabular}
\par\vspace{.5\baselineskip}
\caption{The exceptions to the integer leaf-count bound.
The required number is $12-c+\varepsilon$.}
\label{tab:deletion-exceptions}
\end{table}

\begin{theorem}[Reductions at order twenty-five]
\label{thm:twentyfive-reductions}
For each choice $(k,D,\varepsilon)$ in
\eqref{eq:twentyfive-parameters}, every forest $F$ of order
$25$ has a partition
\[
 V(F)=S\sqcup B\sqcup A,\qquad
 |S|=k,\quad |B|=12-k,\quad |A|=13,
\]
such that $(B,A)$ is a bipartition of $F-S$ and
$\Delta_{F-S}(B)\leq D$. When $\varepsilon=1$, the set
$A$ contains a vertex $r$ isolated in $F-S$ with
$d_F(r)\leq2$.
\end{theorem}

\begin{proof}
Choose $(C,I)$ by taking the smaller side of each nontrivial
component in $C$ and placing all isolated vertices in $I$.
Write $c=|C|\leq12$. If $c\leq k$, delete all of $C$ and
any $k-c$ vertices of $I$. All remaining vertices are
isolated. Their original degrees sum to at most
$e(F)\leq24$. There are $25-k\geq21$ of them, so one has
original degree at most one. Put that vertex in $A$, complete
$A$ to size $13$, and put the others in $B$.

Suppose $c>k$, and put
\[
 H=\{v\in C:d_F(v)\geq D+1\},\qquad h=|H|.
\]
Every vertex of $C$ has positive degree, so
$c+Dh\leq e(F)\leq24$. In all three parameter ranges
this gives $h\leq k$. Let $z,L,q$ count vertices of $I$
of degrees zero, one, and at least two, and let $a$ count
the leaves adjacent to $H$. As in the proof of
Lemma~\ref{lem:one-side-isolation},
\[
 \max\{c+Dh,L+2q\}\leq e(F)\leq24-z.
\]
If $h>0$, the forest between $H$ and the $q$ nonleaf
vertices has at most $h+q-1$ edges; thus
$(D+1)h-a\leq h+q-1$. All conditions of
Lemma~\ref{lem:finite-deletion-calculation} hold.

Give each vertex of $C$ its number of leaf neighbors in $I$
as weight. Delete $H$ and the $k-h$ vertices of $C\setminus H$
with largest weights. Lemma~\ref{lem:largest-integer-weights}
shows that these deletions isolate at least $\Lambda$
vertices of $I$, all of original degree zero or one.
Outside Table~\ref{tab:deletion-exceptions}, there are at
least $12-c+\varepsilon$ such vertices.

We now treat the table. Except for its sixth row,
$z=0$ and $L+2q=24$. These equalities force $F$ to be
connected and every nonleaf vertex of $I$ to have degree
two. Delete the $I$-leaves and contract one edge at each
remaining $I$-vertex. This gives a tree $R$ on $C$.
Whenever both endpoints of an edge of $R$ are deleted,
the intervening degree-two vertex of $I$ becomes isolated.
Such vertices are distinct from the leaves counted by their
weights.

\emph{The first and fourth rows.}
Here $H$ consists of one vertex of weight zero; the eight
other weights sum to eight. Put $t=k-1$, which is three
in the first row and two in the fourth. If the largest $t$
weights sum to at least $t+1$, their leaves suffice.
Otherwise all eight weights are one. Indeed, if a largest
weight were at least two, either the other $t-1$ largest
were positive, giving sum at least $t+1$, or fewer than
$t$ entries would carry the total eight, making the largest
$t$ sum to eight. Choose $t$ vertices
outside $H$, including a neighbor in $R$ of its unique
vertex. The deletion isolates $t$ leaves and the intervening
vertex, as required.

\emph{The second row.}
The seven weights outside the two vertices of $H$ sum to
seven. If the largest two sum to at least three, these leaves
and the one leaf at $H$ suffice. Otherwise all seven weights
are one. If the two vertices of $H$ are adjacent in $R$,
their connector is already isolated. If they are not, choose
one of the two additional vertices to be a neighbor of $H$
in $R$. In either case the deletion isolates at least
$1+2+1=4$ vertices.

\emph{The third row.}
The eight outside weights sum to six. If the largest two
sum to at least three, we are done. Otherwise the weights
are six ones and two zeros. Both vertices of $H$ have
degree at least five in $R$, since $a=0$. They must be
adjacent: if not, their neighborhoods lie in the other
eight vertices and have at most one common vertex, since
$R$ is a tree. Their degree sum would then be at most
nine, a contradiction. Delete two positive-weight vertices
in addition to $H$. Their two leaves and the connector
between the vertices of $H$ give three isolates.

\emph{The fifth row.}
There is no vertex of $H$. If the largest three weights
sum to at least four, we are done. Otherwise the nine
weights are eight ones and one zero. If $R$ has an edge
joining positive-weight vertices, delete its endpoints
and a third positive-weight vertex. This isolates three
leaves and their connector. If no such edge exists, every
edge of $R$ meets the zero-weight vertex, so $R$ is a star.
Delete its center and two leaves, isolating two original
leaves and two connectors.

\emph{The sixth row.}
The nine weights sum to nine. Unless the three largest sum
to at least four, all nine are one. Among the seven nonleaf
vertices of $I$ there is one of degree two: otherwise
$e(F)\geq9+3\cdot7>24$. Delete its two neighbors in $C$
and any third vertex of $C$. These deletions isolate that
vertex and three leaves. Connectivity of $F$ is not needed
in this case.

\emph{The seventh and eighth rows.}
In each row $H$ consists of one vertex. Its degree in $R$
is at least $D+1-a$, which is eight in the seventh row and
nine in the eighth. Hence $R$ is a star centered there.
The outside weights sum to seven or six, respectively;
their two largest have sum at least two. Delete those two
vertices and $H$. Their two connectors, their leaves, and
any leaf at $H$ give at least five isolates in the seventh
row and four in the eighth, exceeding what is required.

Thus in every case $S$ contains $H$ and isolates at least
$12-c+\varepsilon$ vertices of $I$ of original degree at
most two. Move $12-c$ of them to $C\setminus S$. The
resulting side $B$ has $c-k+12-c=12-k$ vertices, and
the other side $A$ has $13$. Moved vertices are isolated,
and every surviving vertex of $C$ has degree at most $D$.
If $\varepsilon=1$, retain one of the counted isolates
in $A$ and call it $r$. This proves every assertion.
\end{proof}

\subsubsection{Restoring the deleted vertices}

The degree condition on $r$ concerns the forest before
deletion. Isolation after deletion alone would not ensure
that a deleted vertex avoids $r$.

\begin{proposition}[A reserved label]\label{prop:reserved-label-host}
Let $X$ cover $(8,12,4)$ or $(9,12,8)$, in the sense of
Section~\ref{sec:bipartite-hosts}, and put $k=4$ or $3$,
respectively. Adjoin a vertex $e$, a vertex $\zeta$, and
a set $U$ of $k$ vertices. Suppose every vertex of $U$ is
universal in the resulting graph $H$, and $\zeta$ is adjacent
to every other vertex except possibly $e$. No further
edges at $e$ are required. Then every forest on $26$
vertices, with a specified vertex $v$, embeds bijectively
in $H$ with $v$ sent to any specified member of $U$.

If $X$ instead covers $(9,13,6)$, adjoining four universal
vertices gives the same rooted assertion, with the root
at any one of those vertices.
\end{proposition}

\begin{proof}
Write $F$ for the forest obtained by deleting $v$.
Apply case $\mathrm I$ or $\mathrm {III}$
of Theorem~\ref{thm:twentyfive-reductions}. Put $e$ on
the reserved vertex $r$. All neighbors of $r$ in this
$25$-vertex forest lie in $S$, and there are at most two.
Since $|S|=k\geq3$, choose $s\in S\setminus N_F(r)$
and put $\zeta$ on $s$. Assign the other $k-1$ vertices
of $S$ to $U$ other than the prescribed root image.
Embed the forest on $B\sqcup(A\setminus\{r\})$ in $X$,
and send $v$ to its prescribed member of $U$.

This is a bijection onto $H$. Edges within the residual
forest are preserved by its embedding. Every edge meeting
a vertex assigned to $U$ is allowed by universality.
Every remaining edge at $s$ is allowed by the condition
on $\zeta$: its only possible forbidden partner lies on
$r$, which was chosen nonadjacent to $s$. The vertex $r$
has no neighbor in the residual forest. These checks cover
all edges, including those at the restored root.

For the last assertion, delete $v$ and use case
$\mathrm {II}$ of Theorem~\ref{thm:twentyfive-reductions}.
Embed the remaining $9$-by-$13$ forest in $X$, assign the
three deleted vertices to three universal vertices, and
send $v$ to the fourth. Every unexamined edge meets a
universal vertex.
\end{proof}

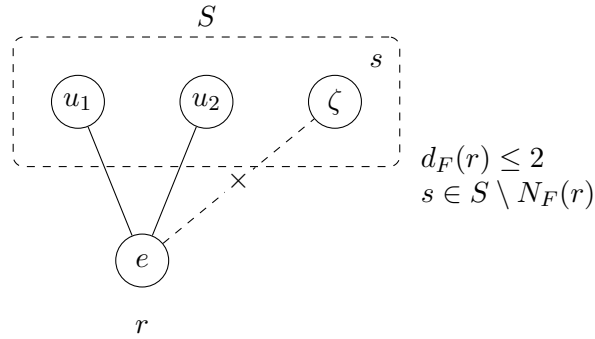
\begin{figure}[htbp]
\centering
\begin{tikzpicture}[x=1cm,y=1cm,
  every node/.style={font=\normalfont\normalsize},
  vertex/.style={circle,draw,fill=white,minimum size=7mm,inner sep=1pt}]
\node[vertex] (u1) at (0,1.1) {$u_1$};
\node[vertex] (u2) at (1.7,1.1) {$u_2$};
\node[vertex] (z) at (3.4,1.1) {$\zeta$};
\node[draw,dashed,rounded corners,fit=(u1)(u2)(z),inner sep=5mm,
      label=above:{$S$}] {};
\node[vertex] (r) at (.85,-1) {$e$};
\node[below=3mm of r] {$r$};
\draw (r)--(u1) (r)--(u2);
\draw[dashed] (r)--node[fill=white,inner sep=1pt] {$\times$}(z);
\node[align=left,anchor=west] at (4.4,.1)
  {$d_F(r)\leq2$\\$s\in S\setminus N_F(r)$};
\node[above right=1mm of z] {$s$};
\end{tikzpicture}
\caption{The reserved-label construction when $|S|=3$.
The circles display assigned labels: $u_1,u_2$ are universal,
and $\zeta$ may fail only at $e$. The crossed dashed segment
is absent from the input forest. Deleting $S$ isolates $r$;
the prescribed root, omitted from the drawing, receives
another universal label.}
\label{fig:reserved-label}
\end{figure}

\subsubsection{The finite hosts}

We specify a sparse host that will be transferred to many
label sets. Its black vertices have names
\[
 45,49,51,55,57,63,65,69.
\]
Its white vertices and all their black neighbors are given
in Table~\ref{tab:reserved-sparse-host}. The names are
integers, and every listed edge is coprime, but only the
listed edges are used. Denote this graph by $X_0$.

\begin{table}[htbp]
\centering
\renewcommand{\arraystretch}{1.15}
\begin{tabular}{r l}
\hline
White vertex & Black neighbors \\
\hline
$46$ & $57$ \\
$48$ & $55$ \\
$50$ & $51,57$ \\
$52$ & $51,55$ \\
$54$ & $55,65$ \\
$56$ & $45,51,65$ \\
$58$ & $45,49,51,55,57,63,65,69$ \\
$62$ & $45,49,55,57$ \\
$64$ & $45,51,55,57,63,65$ \\
$66$ & $49,65$ \\
$68$ & $45,63,65,69$ \\
$70$ & $51,57,69$ \\
\hline
\end{tabular}
\par\vspace{.5\baselineskip}
\caption{The $38$ edges of the host $X_0$.}
\label{tab:reserved-sparse-host}
\end{table}

There are four further hosts. For each row of
Table~\ref{tab:reserved-arithmetic-hosts}, remove $U$ and
the displayed labels $\zeta,e$, when present, from the
interval $L$. On the remaining labels take the bipartite
coprimality graph, with odd labels black and even labels
white. Denote the resulting graphs by $X_1,X_2,X_3,X_4$
in the order of the table. The first has sides of sizes
$9,13$; the other three have sides of sizes $9,12$.

\begin{table}[htbp]
\centering
\renewcommand{\arraystretch}{1.25}
\begin{tabular}{c l r r}
\hline
$L$ & $U$ & $\zeta$ & $e$ \\
\hline
$\{183,\ldots,208\}$ & $191,193,197,199$ & --- & --- \\
$\{243,\ldots,268\}$ & $251,257,263$ & $253$ & $264$ \\
$\{795,\ldots,820\}$ & $797,809,811$ & $799$ & $816$ \\
$\{825,\ldots,850\}$ & $827,829,839$ & $841$ & $840$ \\
\hline
\end{tabular}
\par\vspace{.5\baselineskip}
\caption{Arithmetic inventories for $X_1,X_2,X_3,X_4$.
Every label in $U$ is universal in $L$, and $\zeta$ is
coprime to every other label except possibly $e$.}
\label{tab:reserved-arithmetic-hosts}
\end{table}

\begin{proposition}[Five finite hosts]\label{prop:reserved-finite-hosts}
The graph $X_0$ covers $(8,12,4)$, the graph $X_1$ covers
$(9,13,6)$, and each of $X_2,X_3,X_4$ covers $(9,12,8)$.
\end{proposition}

\begin{proof}
We apply Proposition~\ref{prop:bipartite-witnesses}. The
capacity inequalities for the three parameter choices are
$4\cdot8\geq19$, $6\cdot9\geq21$, and
$8\cdot9\geq20$, respectively. Thus it suffices to check
every kernel and feasible leaf vector.

Leaf augmentation gives all $23$ black trees at order
eight and all $47$ at order nine. Enumerating their star
partitions and applying the degree bounds gives respectively
$154$, $493$, and $496$ colored kernel types. The first
enumeration tests all $877$ partitions of seven edge
positions for each black tree; the other two test all
$4\,140$ partitions of eight positions. The supplied
kernel lists have exactly these types, each once.

For each kernel enumerate all vectors in
\eqref{eq:kernel-leaf-vectors}, with the appropriate white
order and black degree bound. The resulting numbers of
kernel--vector pairs are
\[
 \begin{array}{c|r|r|r}
 \text{Host} & \text{Kernels} & \text{Vectors} & \text{Edge checks}\\
 \hline
 X_0 &154&258\,592&4\,913\,248\\
 X_1 &493&4\,651\,581&97\,683\,201\\
 X_2,X_3,X_4\text{ (each)}&496&2\,652\,306&53\,046\,120
 \end{array}
\]
The certificate for $X_0$ lists each vector explicitly,
together with its full embedding. The checker requires
each vector to occur exactly once. The other certificates
list full embeddings in the lexicographic vector order
used in Section~\ref{sec:bipartite-hosts}; their offsets,
counts, and end of stream are checked against that
enumeration. In every record the images must be a
side-preserving bijection, and every kernel and leaf edge
must be an edge of the specified host.

All $12\,867\,091$ records pass, with $261\,734\,809$
edge checks in total. The sparse host is reconstructed
from Table~\ref{tab:reserved-sparse-host}; the other hosts
are reconstructed by exact gcd calculations from
Table~\ref{tab:reserved-arithmetic-hosts}. Complete kernel
and vector coverage, together with the checked embeddings,
therefore proves all five assertions by the witness criterion.
\end{proof}

\begin{corollary}[Explicit rooted intervals]
\label{cor:reserved-rooted-intervals}
Every forest of order $26$ has a coprime labeling by
$\{45,\ldots,70\}$, with a specified vertex assigned
any member of $\{47,53,59,61,67\}$.
The same assertion holds for each interval in
Table~\ref{tab:reserved-arithmetic-hosts}, with its root
assigned any label in the corresponding set $U$.
\end{corollary}

\begin{proof}
For the first interval use $X_0$, reserve the label $60$,
and use $47,53,59,61$ as $U$ and $67$ as $\zeta$ in
Proposition~\ref{prop:reserved-label-host}. These five
labels are all universal in the interval. The proposition
gives, in particular, a labeling with root $47$. Exchanging
two universal labels preserves every edge and permits
any of the five root labels.

For the other rows, use Proposition~\ref{prop:reserved-finite-hosts}
and the corresponding case of
Proposition~\ref{prop:reserved-label-host}. Direct gcd
checks give the universal and nearly universal properties
in the table. For example, $253=11\cdot23$ conflicts only
with $264$ in its interval, and $799=17\cdot47$ conflicts
only with $816$ in its interval. The label $841=29^2$
is fully universal in $\{825,\ldots,850\}$, which also
meets the stated hypothesis. The inventories are disjoint
and exhaust their respective intervals, so the embeddings
are the required coprime labelings.
\end{proof}

\subsubsection{A finite family of arithmetic transfers}

The sparse host $X_0$ needs substantially fewer arithmetic
edges than the full coprimality graph on its vertex names.
An edge-preserving image of these $38$ edges is enough to
transfer its covering property.

\begin{theorem}[Intervals with five universal labels]
\label{thm:five-universal-intervals}
Let $a$ be odd, $1\leq a\leq7125$, and put
$L=\{a,a+1,\ldots,a+25\}$. If $L$ contains at least
five universal labels, then every forest of order $26$
has a coprime labeling by $L$. Any specified vertex may
receive any universal label of $L$.
\end{theorem}

\begin{proof}
For each of the $3\,563$ possible starts, compute the set
\[
 U(L)=\{x\in L:\gcd(x,y)=1\text{ for all }y\in L\setminus\{x\}\}.
\]
Exactly $1\,944$ starts have $|U(L)|\geq5$. For every
one of these starts the transfer table supplies a
bijection from $V(X_0)$ and six additional positions
onto $L$. The first five additional positions receive
distinct members of $U(L)$; the sixth receives an even
reserved label. The eight black vertices of $X_0$ receive
odd labels, the twelve white vertices receive even labels,
and all $38$ edges have coprime images.

The checker recomputes $U(L)$ for every permitted start,
requires the set of supplied starts to equal the qualifying
set, and checks every inventory and every displayed host
edge. All $1\,944$ transfers pass, giving $73\,872$
edge gcd checks. Thus each qualifying arithmetic graph
contains a copy of $X_0$, a reserve, and five universal
vertices. Choose any of the five as the root image and a
different one as $\zeta$. Proposition~\ref{prop:reserved-label-host}
then gives the required rooted labeling.

The prescribed universal label need not belong to the
five selected in the transfer. If it does not, exchange it
with the label at the root. Both labels are universal in
the entire inventory, so all edges incident to either
vertex remain coprime; every other edge and the bijection
are unchanged. This proves the full root assertion.
\end{proof}

The remaining $1\,619$ starts fail this theorem's hypothesis.
Of the qualifying intervals, $360$ end at or below $1354$.
These counts describe blocks of length $26$; they do not
assert prime labeling for every order up to $7150$.
Applying the finite-prefix induction still requires its
complete inventory of block assertions.

The supplement includes the five host certificates, all
$1\,944$ transfer bijections, and the exact deletion
calculation. Its checker regenerates the required kernel,
leaf-vector, parameter, and interval sets before checking
their witnesses. The structural cases in the proof of
Theorem~\ref{thm:twentyfive-reductions} account for every
tuple excluded from the numerical leaf-count bound.

\subsection{Larger finite blocks}\label{sec:larger-blocks}

The same small-side reductions give blocks of order $30$.
Three of the intervals below contain six universal labels;
their residual hosts have sides of sizes $9$ and $15$.
A fourth interval has only five universal labels, so we use
the ten-vertex side of Theorem~\ref{thm:uniform-small-sides}.
We verify these four hosts together, distinguishing the
complete leaf-vector lists from the list reduced by equal
neighborhoods. We then explain how one of the hosts supplies
blocks of two larger orders without another enumeration.

\subsubsection{Four arithmetic hosts}

For each row of Table~\ref{tab:thirty-inventories}, remove the
set $U$ from the interval $L$. On the remaining labels take
the bipartite coprimality graph, with odd labels black and
even labels white. Call the resulting graphs $X_1,X_2,X_3,X_4$,
in the order of the table. Each has $15$ white vertices.
The first three have $9$ black vertices, and the fourth has
$10$.

\begin{table}[htbp]
\centering
\renewcommand{\arraystretch}{1.35}
\begin{tabular}{c l c}
\hline
$L$ & $U$ & $b$ \\
\hline
$\{77,\ldots,106\}$ & $79,83,89,97,101,103$ & $9$ \\
$\{91,\ldots,120\}$ & $97,101,103,107,109,113$ & $9$ \\
$\{1285,\ldots,1314\}$ & $1289,1291,1297,1301,1303,1307$ & $9$ \\
$\{931,\ldots,960\}$ & $937,941,943,947,953$ & $10$ \\
\hline
\end{tabular}
\par\vspace{.5\baselineskip}
\caption{The four $30$-label inventories. The set $U$ is the
complete set of universal labels in $L$; $b$ is the number
of remaining odd labels.}
\label{tab:thirty-inventories}
\end{table}

All labels in $U$ are universal in their respective intervals.
This can be checked directly by gcds. The only composite
label in these sets is $943=23\cdot41$. No other multiple
of either factor lies between $931$ and $960$, so this label
is universal as well. The criterion concerns coprimality in
the interval, not primality of the label.

\begin{proposition}[Four hosts]\label{prop:larger-finite-hosts}
Each of $X_1,X_2,X_3$ covers $(9,15,6)$, and $X_4$ covers
$(10,15,6)$.
\end{proposition}

\begin{proof}
We use Proposition~\ref{prop:bipartite-witnesses}. The
completion inequalities are $6\cdot9\geq23$ and
$6\cdot10\geq24$. Thus it suffices to verify the complete
kernel and leaf-vector families with the stated side sizes.

For nine black vertices, leaf augmentation gives $47$
unrooted trees. Apply all $4\,140$ partitions of their eight
edge positions, retain the star partitions, and impose the
kernel degree conditions. Deduplication by colored codes
gives $493$ kernels. For each kernel with $q$ white vertices,
enumerate the vectors
\begin{equation}\label{eq:larger-leaf-vectors}
 \sum_{v\in B}\ell_v=15-q,
 \qquad 0\leq\ell_v\leq6-d_K(v).
\end{equation}
There are $13\,652\,198$ kernel--vector pairs. The same
complete list is checked separately against each of the
three arithmetic hosts.

For ten black vertices, the corresponding enumeration starts
with $106$ unrooted trees and all $21\,147$ partitions of
nine edge positions. It gives $1\,529$ colored kernels.
The number of vectors in \eqref{eq:larger-leaf-vectors} is
$63\,388\,489$. Here we use
Lemma~\ref{lem:kernel-twin-vectors}: only vectors sorted
within each class of equal black neighborhoods are required.
Their total number is $22\,617\,172$.

We specify how this last count is verified independently of
the vector list. For a class of size $t$ with coordinate cap
$c$, let $m_j$ be the number of coordinates equal to $j$.
Its factor in the counting polynomial is
\[
 \sum_{\substack{m_0,\ldots,m_c\geq0\\
                   m_0+\cdots+m_c=t}}
                  z^{\sum_{j=0}^c j m_j}.
\]
Multiply the factors over the classes and take the
coefficient of $z^{15-q}$. Each multiplicity choice determines
exactly one nondecreasing sequence in its class, so this is
the coefficient in Lemma~\ref{lem:kernel-twin-vectors}.
For every kernel the enumerated vectors are admissible,
strictly increasing in lexicographic order, and equal in
number to this coefficient. Hence the list is complete.
For comparison, the unreduced count is the coefficient of
$z^{15-q}$ in
\[
 \prod_{v\in B}(1+z+\cdots+z^{6-d_K(v)}).
\]
Both coefficient calculations use integer polynomial
arithmetic.

The certificates supply a full embedding for each required
pair. A record lists the images of the black kernel vertices,
the white kernel vertices, and then the attached white leaves
in black-vertex order. The checker requires a side-preserving
bijection and checks every kernel edge and every leaf edge
against the exact arithmetic host. It also checks the
complete kernel set, every vector count, all record offsets,
and the end of each stream. Table~\ref{tab:larger-host-counts}
gives the resulting counts.

\begin{table}[htbp]
\centering
\renewcommand{\arraystretch}{1.3}
\begin{tabular}{c r r r}
\hline
Host & Kernels & Checked embeddings & Checked edges \\
\hline
$X_1,X_2,X_3$ (each) & $493$ & $13\,652\,198$ & $314\,000\,554$ \\
$X_4$ & $1\,529$ & $22\,617\,172$ & $542\,812\,128$ \\
\hline
\end{tabular}
\par\vspace{.5\baselineskip}
\caption{The complete larger-host certificates. Equal-neighborhood
permutations reduce the vector list only for $X_4$. Every
record is a full embedding of a tree on $24$ or $25$ vertices.}
\label{tab:larger-host-counts}
\end{table}

All $63\,573\,766$ supplied embeddings pass, giving
$1\,484\,813\,790$ edge checks. For the first three hosts,
every feasible vector is checked. For the fourth,
Lemma~\ref{lem:kernel-twin-vectors} transfers an embedding
back from its sorted representative. The witness criterion
therefore proves all four covering assertions, including
their forest conclusions.
\end{proof}

\subsubsection{Rooted blocks and one-label deletions}

\begin{theorem}[Four rooted blocks of order thirty]
\label{thm:rooted-thirty-blocks}
For every row $(L,U)$ of Table~\ref{tab:thirty-inventories},
every forest on $30$ vertices with a specified vertex $v$
has a coprime labeling by $L$ assigning to $v$ any specified
label of $U$.
\end{theorem}

\begin{proof}
Delete $v$. In the first three rows, apply
Theorem~\ref{thm:uniform-small-sides} with $b=9$, $D=6$,
and $k=5$ to the remaining forest of order
$29=2\cdot5+2\cdot9+1$. It gives five deleted vertices and
a residual forest with sides $9$ and $15$, with black
degree at most six. Embed that forest in the appropriate
host by Proposition~\ref{prop:larger-finite-hosts}. Assign
the specified label to $v$ and the other five labels of $U$
to the deleted vertices.

For the fourth row use $b=10$, $D=6$, and $k=4$. Now
$29=2\cdot4+2\cdot10+1$, and the residual sides have sizes
$10$ and $15$. Embed in $X_4$, and use the five labels of
$U$ on the four deleted vertices and $v$, respecting the
specified root label.

In both cases the residual embedding uses precisely $L\setminus U$.
The restored labels are distinct and exhaust $U$, so the
full labeling is a bijection. Every edge not in the residual
forest meets a vertex carrying a universal label. This also
checks all edges incident to $v$ and completes the proof.
\end{proof}

\begin{corollary}[Deleting a universal label]
\label{cor:thirty-block-deletion}
For each row $(L,U)$ of Table~\ref{tab:thirty-inventories}
and every $u\in U$, the assertion $\cP(L\setminus\{u\})$ holds.
\end{corollary}

\begin{proof}
Adjoin an isolated vertex to any forest on $29$ vertices.
Use Theorem~\ref{thm:rooted-thirty-blocks} with that vertex
specified to receive $u$, and restrict the resulting labeling.
The restriction uses exactly $L\setminus\{u\}$ and preserves
all edges of the original forest.
\end{proof}

\subsubsection{Enlarging the source block}

An added white vertex and an added universal vertex have
different roles. The first need only meet the black side;
the second must meet every other vertex. Keeping the two
possibilities separate gives useful source graphs of orders
$31$ and $32$.

Write $J_t(X)$ for the graph obtained from $X$ by adjoining
$t$ vertices, each adjacent to every other vertex of the
resulting graph. In particular, the new vertices form a
clique. If $X$ is bipartite, write $X^+$ for the bipartite
graph obtained by adding one white vertex adjacent to all
black vertices and to no white vertices.

\begin{proposition}[Larger source blocks]
\label{prop:ten-side-enlargements}
Suppose $X$ covers $(10,15,6)$. Then every forest of order
$31$ embeds bijectively in each of $J_6(X)$ and $J_5(X^+)$.
Every forest of order $32$ with a specified vertex $v$
embeds bijectively in $J_6(X^+)$, sending $v$ to any
specified one of its six added universal vertices.
\end{proposition}

\begin{proof}
Let $F$ have order $31$. Theorem~\ref{thm:uniform-small-sides}
with $b=10$, $D=6$, and $k=5$ gives a set $S$ of five
vertices such that $F-S$ has sides $B,A$ of sizes $10,16$
and black degree at most six.

To embed in $J_6(X)$, choose any vertex $w\in A$ and delete
it as well. The remaining forest has sides $10,15$, and
deleting $w$ cannot raise a black degree. Embed it in $X$,
and send $S\cup\{w\}$ to the six added universal vertices.
Every edge outside the residual forest meets one of those
vertices, so it is preserved. The vertex $w$ need not be
a leaf; universality permits all of its edges.

For the second host, Lemma~\ref{lem:white-host-extension}
applies because $15\geq10-1$ and $6\cdot10\geq10+15$.
Thus $X^+$ covers $(10,16,6)$. Embed $F-S$ directly in
$X^+$, and assign $S$ to the five added universal vertices.
This again uses every host vertex once and preserves every
edge.

Finally delete the specified vertex of a forest of order
$32$. Embed the remaining forest in $J_5(X^+)$ by the
second construction, using five of the six universal
vertices. Send the specified vertex to the remaining,
prescribed one. Every edge restored at that vertex is
allowed. This proves the rooted assertion.
\end{proof}

These constructions use the same certified host $X_4$.
The two source graphs of order $31$ place the extra vertex
in different parts: $J_6(X_4)$ has one more universal
vertex, whereas $J_5(X_4^+)$ has one more white vertex.
For an arithmetic application, an edge-preserving bijection
from the chosen source graph to the actual label inventory
must still be supplied. In the rooted case it must also
carry a permitted universal vertex to the required root
label. Lemma~\ref{lem:host-transfer} then gives the labeling
assertion. No statement about every interval of length
$31$ or $32$ follows without those arithmetic transfers.

The supplement contains all four host headers, their complete
embedding streams, and a checker regenerating the kernel and
vector families used above. The first three streams contain
every feasible vector; the fourth contains exactly the
representatives justified by equal-neighborhood permutations.
These finite host assertions can now be used in the
finite-prefix certificate system of Sections~\ref{sec:finite-certificates}
and~\ref{sec:branch-peeling}. Verifying that system's complete
arithmetic inventories remains a separate step.

\subsection{The initial finite range}\label{sec:initial-prefix}

We next assemble a finite proof for the first $299$ orders. The
certificate uses the local reductions of
Section~\ref{sec:finite-reductions} and the branch rules of
Section~\ref{sec:finite-certificates}. We first record the few additional
forms needed to apply those rules to the saved arithmetic inventories.
The complete proof graph, including all terminal tables, accompanies
this section. No previously accepted interval is used as a terminal
assertion.

\subsubsection{Extensions near a prime}

\begin{lemma}\label{lem:prime-contact-extension}
Let $p\geq3$ be prime. If $\cP([p-1])$ holds, then
$\cP([p])$, $\cP([p+1])$, and $\cP([p+2])$ hold.
\end{lemma}

\begin{proof}
The label $p$ is universal in each of the three intervals, since
$p+2<2p$. For a tree of order $p$, delete a leaf and label the
remaining tree by $[p-1]$; give the leaf label $p$. For order $p+1$,
delete a leaf and its neighbor, label the remaining forest by $[p-1]$,
and restore the leaf and neighbor with labels $p+1$ and $p$.

For order $p+2$, use the configuration in the proof of
Lemma~\ref{lem:host-attached-triangle}. If two leaves have a common
neighbor, delete those three vertices, label the remaining forest
by $[p-1]$, and restore the neighbor with label $p$ and the leaves
with labels $p+1,p+2$. Otherwise a leaf has a neighbor of degree two.
Let $x,y,z$ be the leaf, that neighbor, and its other neighbor.
Delete all three, label the remainder by $[p-1]$, and put
$p+2,p+1,p$ on $x,y,z$, respectively. The edge $xy$ has consecutive
labels; every other restored edge meets $z$, whose label is universal.
Both constructions use the entire required inventory. Forest
completion proves the stated forest assertions.
\end{proof}

This implication advances from order $p-1$ by at most three vertices.
It requires the smaller assertion even when several of these intervals
share the same prime.

\subsubsection{One omitted block size}

The parity version of the detached-block lemma already gives most
of the required block choices. Two certificate nodes use a different
shortening of the range.

\begin{lemma}\label{lem:omit-penultimate-block}
Let $b\geq4$ and $n\geq b+2$. Every tree $T$ of order $n$ has a
vertex $u$ and a union $A$ of components of $T-u$ such that
\[
 b\leq\abs A\leq2b-2,\qquad \abs A\ne2b-3.
\]
Consequently, the conclusion of Proposition~\ref{prop:block-extension}
also holds if its size set is replaced by
\[
 \{b,b+1,\ldots,2b-4,2b-2\}\cap[1,n-1]
\]
and each block assertion is the unrestricted assertion $\cP(B_k)$.
\end{lemma}

\begin{proof}
Root $T$ at a leaf and choose a deepest vertex $v$ whose descendant
subtree has at least $b+1$ vertices. It has a parent, and its child
subtrees have orders at most $b$. Their total is at least $b$.
If a selection of them has a permitted total, take that selection
and put $u=v$.

Suppose there is no such selection. No child subtree has order $b$.
Choose an inclusion-minimal selection with total at least $b$.
Deleting its last selected member shows that its total is at most
$2b-2$; it must therefore be $2b-3$. Deleting any member leaves
at most $b-1$ vertices, so every selected member has order at least
$b-2$. Three such members would have total at least
$3b-6>2b-3$. There are exactly two, of orders $b-1,b-2$.

An additional child subtree of order $a$ would give a total
$b-1+a$ between $b$ and $2b-2$ together with the larger selected
subtree. Thus $a=b-2$. But the two subtrees of order $b-2$ would
then have permitted total $2b-4$, a contradiction. Hence there
are no additional children. Take the whole descendant subtree
of $v$, of order $2b-2$, and put $u$ at its parent.

For the consequence, label $u$ by the universal label $q_k$, the
detached forest by $B_k$, and the remaining forest by $C_k$, exactly
as in Proposition~\ref{prop:block-extension}. Here no degree bound
on the detached forest is needed, because its block assertion is
unrestricted.
\end{proof}

\subsubsection{Branches of orders at most four}

For small $b$, the full profile condition
\eqref{eq:short-branch-profile} has a simple description. The total
$s$ in the following statement is always at most $n-2$.

\begin{lemma}\label{lem:small-profile-families}
The short-branch profiles for $b=3,4,5$ are exactly the following
families, subject to $b\leq s\leq n-2$:
\[
\begin{array}{ccl}
 b&\text{profile}&\text{additional conditions}\\[2pt]
 3&(0,k)&k\geq2,\\
 4&(0,t,k)&t\in\{0,1\},\\
 5&(0,a,0,d)&a,d\geq0,\\
 5&(0,0,c,d)&c\geq1,\ d\geq0,\\
 5&(1,0,c,0)&c\geq2.
\end{array}
\]
In the row for $b=4$, the total condition says $k\geq2$ when
$t=0$, and $k\geq1$ when $t=1$.
\end{lemma}

\begin{proof}
For $b=3$, a branch of order one together with a branch of order
two gives the forbidden sum, as do three singletons. A profile
containing a singleton therefore has total below three. The remaining
profiles consist of at least two branches of order two.

For $b=4$, two branches of order two are forbidden. A singleton
cannot coexist with a branch of order three. If there is a singleton
and no branch of order three, avoiding four also permits at most
one branch of order two and at most one singleton when that branch
is present. Without that branch there are at most three singletons.
Thus any profile containing a singleton has total below four.
The stated family follows.

For $b=5$, first suppose there is no singleton. Branches of orders
two and three cannot coexist. This gives the first two rows for
$b=5$, with the second row requiring a branch of order three to
avoid overlap. If there is exactly one singleton, order four is
forbidden and there can be at most one branch of order two. Such
a branch excludes order three and gives total only three. Therefore
there is no branch of order two, leaving at least two branches of
order three. Finally, two or more singletons exclude orders three
and four. Two branches of order two together with a singleton are
also forbidden, and a single branch of order two permits at most
two singletons. In either case the total is below five. Conversely,
none of the listed profiles has a submultiset of total $b$.
\end{proof}

The one-vertex rooted block needs only a label coprime to the center.
For a two-vertex block, a coprime pair with at least one permitted
root label suffices. The next lemma gives the remaining two forms.
It is stated for graphs, so it applies equally to a coprimality
host and to a temporarily enlarged host.

\begin{lemma}\label{lem:small-rooted-hosts}
The following rooted host assertions hold.
\begin{enumerate}
\item If a graph on $\{a,b,d\}$ contains $ab,bd$, and
      $O$ contains $b$ and at least one of $a,d$, then
      $\mathcal U_1(H,O)$ holds.
\item If a graph on $\{a,b,e,f\}$ contains $ab,ae,af,be$, and
      $\{a,b\}\subseteq O$, then $\mathcal U_1(H,O)$ holds.
\end{enumerate}
\end{lemma}

\begin{proof}
Every three-vertex tree is a path. Its prescribed root is either
the middle vertex, which can receive $b$, or an endpoint, which
can receive the endpoint belonging to $O$.

A four-vertex tree is a star or a path. In the star, put its center
at $a$; a marked leaf can be put at $b$. In the path, a marked
endpoint can be put at $b$ using the order $b,e,a,f$, whereas a
marked internal vertex can be put at $b$ using the order $e,b,a,f$.
All displayed successive pairs are edges of the host. These choices
cover every possible marked position. Forest completion retains
the prescribed root.
\end{proof}

Thus, for a three-label branch $(a,b,d)$ attached to a center labeled
$c$, the conditions used in the certificate are
\[
 \gcd(a,b)=\gcd(b,d)=\gcd(c,b)=1,
 \qquad \gcd(c,a)=1\ \text{or}\ \gcd(c,d)=1.
\]
For a four-label branch $(a,b,e,f)$, it suffices to check the four
internal pairs $ab,ae,af,be$ and the two contact pairs $ca,cb$.
The parity conditions used to find some of these blocks are not
needed for this implication.

\subsubsection{Retaining the contact of a three-vertex block}

In a few cases, the contact vertex is retained in the smaller
labeling problem. The label it receives may vary within a certified
root set. We must therefore check the attachment for every label
in that set.

\begin{lemma}\label{lem:retained-three-block}
Let $L$ have $n\geq5$ elements. In the three-vertex detached
alternative of Lemma~\ref{lem:stump-dichotomy}, the following data
suffice to label the tree by $L$. For each row of
Table~\ref{tab:three-boundaries}, choose a set $A\subset L$ of size
three and a nonempty set $G\subseteq L\setminus A$ such that
\[
 \mathcal U_1(H_{L\setminus A},G)
\]
holds. For every $g\in G$, require an ordering $(a_0,a_1,a_2)$ of
$A$ such that the internal pairs in that row are coprime and
$\gcd(g,a_i)=1$ at every indicated contact position. The sets and
orderings may depend on the row, and the ordering may depend on $g$.
\end{lemma}

\begin{table}[bp]
\centering
\begin{tabular}{lll}
Detached forest&Internal pairs&Contact positions\\[3pt]
Three isolated vertices&none&$0,1,2$\\
An edge and an isolate&$01$&$0,2$\\
A path, contact at an endpoint&$01,12$&$0$\\
A path, contact at the middle&$01,12$&$1$
\end{tabular}
\par\vspace{.5\baselineskip}
\caption{All three-vertex detached forests with their attachments to
the retained contact. Positions are indexed before assigning labels.}
\label{tab:three-boundaries}
\end{table}

\begin{proof}
Let $u$ be the contact vertex and $F$ the detached forest. Each
component of $F$ has exactly one neighbor of $u$ in the original
tree: at least one by connectivity, and at most one by acyclicity.
The four rows list every possibility on three vertices, up to an
isomorphism preserving these attachments.

Delete $F$, retain $u$, and use the asserted rooted labeling with
$u$ marked. Let $g\in G$ be its resulting label. Choose the checked
ordering for this $g$ and put its labels on $F$ according to the
appropriate row. The internal edges are coprime, the contact-position
checks prove every edge to $u$, and there are no other crossing
edges. The two inventories partition $L$.
\end{proof}

This lemma replaces only the detached alternative. The short-branch
alternative must still be checked for every profile in
Lemma~\ref{lem:small-profile-families}. If the residual assertion is
unmarked, one may take $G=L\setminus A$, since it then allows every
possible label at the retained contact.

\subsubsection{A whole short-branch block}

It is sometimes more economical to label the center and all its
short branches as a single block. The appropriate internal degree
bound follows directly from the profile.

\begin{lemma}\label{lem:whole-short-branch-block}
In the second alternative of Lemma~\ref{lem:stump-dichotomy}, let
$\mathbf k$ be the branch profile and put
\[
 s=\sum_{j=1}^{b-1}jk_j,\qquad
 d(\mathbf k)=\max\left\{\sum_{j=1}^{b-1}k_j,
                         \max\{j:k_j>0\}\right\}.
\]
The subtree consisting of $v$ and all its child branches has order
$s+1$ and maximum degree at most $d(\mathbf k)$. In particular, if
\[
 L=\{q\}\sqcup B\sqcup C,\qquad \abs B=s+1,
\]
where $d=d(\mathbf k)$, $q$ is universal in $L$, $\cP_d(B)$ holds, and
$\cP(C)$ holds, then this alternative admits a labeling by $L$.
\end{lemma}

\begin{proof}
The degree of $v$ in this subtree is the number of child branches,
$\sum_j k_j$. A vertex in a branch of order $j$ has internal
branch degree at most $j-1$; its attachment vertex gains one edge
to $v$, giving degree at most $j$. This proves the degree bound.
Also $d(\mathbf k)\geq2$, since a single branch has order below $b$
and cannot have total $s\geq b$.

Put $q$ on the parent $u$, label the whole subtree at $v$ by $B$,
and label the forest outside that subtree and $u$ by $C$. Every
crossing edge meets $u$, so universality checks it. The three
inventories are disjoint and exhaust $L$.
\end{proof}

The bound concerns degrees after deleting the parent. Its use as a
terminal restriction is justified by forest completion at degree at
least two. Recursive deletion of a leaf and its neighbor within the
block uses precisely the degree decrease of
Corollary~\ref{cor:leaf-label-extension}.

\subsubsection{Temporary exceptional edges}

One node of the initial certificate uses the repair operation of
Section~\ref{sec:repair}. Its children still assert ordinary coprime
labelings. Only the graph used for joining these children is enlarged.

\begin{proposition}\label{prop:repaired-finite-node}
Let $L$ contain $1$ and two distinct powers of two greater than $1$. Form
$H_L^+$ from $H_L$ by adjoining the edges $q(2q)$ for
$q\in\mathcal Q(L)$. Suppose a checked finite construction embeds
every tree of order $\abs L$ into $H_L^+$. Then $\cP(L)$ holds.

In particular, the three-vertex branch construction may check its
contact labels and internal block edges in $H_L^+$, while using
ordinary coprime labelings for its smaller residual inventories.
Its conclusion is unmarked.
\end{proposition}

\begin{proof}
An embedding into $H_L^+$ is a bijective labeling whose only
possible noncoprime edges have labels $q,2q$ with
$q\in\mathcal Q(L)$. The input tree is $C_4$-free, so
Lemma~\ref{lem:sparse-repair} repairs it. Forest completion gives
$\cP(L)$.

For the stated application, the detached three-vertex block needs
a path on its three labels, since every three-vertex forest extends
to that path. Give its contact a label adjacent to all other labels
in $H_L^+$. In the other alternative, the profiles are $(0,k)$,
so every child branch is a single edge. For each branch, check its
internal edge and the edge from its attachment label to the chosen
center. In the retain case, check adjacency from the center label
to every residual label; in the remove case, use a parent label
adjacent to all other labels. These are exactly the edge checks in
the proof of Proposition~\ref{prop:rooted-branch-extension}, now
made in $H_L^+$. Ordinary child labelings satisfy them because
$H_L$ is a subgraph of $H_L^+$. The repair is applied only after
the whole labeling has been assembled. It is not asserted to
preserve a specified root label.
\end{proof}

\subsubsection{The complete certificate}

A node in the present arithmetic proof graph specifies a finite label
set, with an allowed root set when one is required. Its dependencies
are the exact smaller inventories in its chosen reduction. The only
same-order dependencies are explicit edge-preserving bijections of
label sets. These are subject to cycle detection.

The terminal tables use the codes of
Lemma~\ref{lem:tree-code-enumeration}. An unrestricted table includes
every unrooted tree of its order. A rooted table includes every
rooted tree, and its first image lies in the stated allowed set.
A degree-restricted table includes exactly those unrooted trees
within its declared degree cap. In every case the checker regenerates
the index set and checks equality with the supplied keys before
checking the witnesses. Duplicate keys are rejected when reading
the certificate.

\begin{theorem}\label{thm:initial-prefix-299}
Every forest of order $n$ with $1\leq n\leq299$ has a coprime
labeling by $[n]$.
\end{theorem}

\begin{proof}
We apply the soundness argument of
Theorem~\ref{thm:finite-certificate-soundness}, adding the implications
proved in this section. The supplied arithmetic proof graph has
$1195$ nodes. Its root list is exactly
\[
 \cP([1]),\cP([2]),\ldots,\cP([299]).
\]
The checker verifies this list, every dependency, every inventory,
and the absence of directed cycles. All nodes are reachable from
these roots. Empty inventories, singletons, and coprime pairs are
checked directly; no published initial interval is accepted as a base.

Table~\ref{tab:initial-proof-rules} lists the node types. The local
extensions and host transfers use
Lemma~\ref{lem:host-local-extensions} and
Lemma~\ref{lem:host-transfer}. Prime-contact extensions use
Lemma~\ref{lem:prime-contact-extension}. Detached-block nodes have
a checked partition for every size required by
Proposition~\ref{prop:block-extension} or
Lemma~\ref{lem:omit-penultimate-block}. A bounded block is proved
by a complete degree-restricted table, possibly preceded by
leaf-pair reductions with explicitly checked degree caps.

For each short-branch node the checker regenerates the entire family
\eqref{eq:short-branch-profile} by the subset-sum recursion of
Section~\ref{sec:finite-certificates}. It requires exactly the cases
in that family. Their simpler forms are justified by
Lemma~\ref{lem:small-profile-families}. The branch inventories and
contacts are checked by
Proposition~\ref{prop:rooted-branch-extension} and
Lemma~\ref{lem:small-rooted-hosts}; the alternatives retaining a
three-vertex contact or using a whole bounded block invoke
Lemmas~\ref{lem:retained-three-block} and
\ref{lem:whole-short-branch-block}, respectively. The sole enlarged
host uses Proposition~\ref{prop:repaired-finite-node}.

The check includes $317775$ witnesses in unrestricted tables,
$734195$ in rooted tables, and $169492$ in degree-restricted tables.
Thus all $1221462$ terminal witnesses are checked, with $18216342$
edge gcds. It also checks $215286$ short-branch cases and all $109$
detached-block partitions. The terminal enumeration and all arithmetic
checks pass. The supplement supplies the complete proof graph and a standalone
checker, with its precise counts recorded alongside the certificate.
Since every accepted rule is a proved implication, induction through
the acyclic graph proves every root assertion, including every
integer in the stated range.
\end{proof}

\begin{table}[tbp]
\centering
\begin{tabular}{lr}
Rule&Nodes\\[3pt]
Empty inventory, singleton, or coprime pair&$6$\\
Complete tree table&$150$\\
Transfer of a label set&$38$\\
Unmarked one- or two-vertex extension&$173$\\
Rooted one- or two-vertex extension&$26$\\
Extension near a prime&$164$\\
Detached blocks, unrestricted or degree-restricted&$27$\\
Three-vertex separation and short branches&$551$\\
Four-vertex separation and short branches&$19$\\
Five-vertex separation and short branches&$41$\\[2pt]
Total&$1195$
\end{tabular}
\par\vspace{.5\baselineskip}
\caption{The complete proof graph for orders $1$ through $299$.
The $150$ terminal nodes are distinct from degree-restricted tables
stored inside block reductions. One of the $551$ three-vertex nodes
uses the proved repair operation.}
\label{tab:initial-proof-rules}
\end{table}

\subsection{Continuation of the finite prefix}
\label{sec:prefix-continuation}

We now apply the branch-peeling theorem to a consecutive range of
orders. Each step uses only initial intervals already proved and
finite assertions about its retained blocks. We first specify the
case lists used in this range, then explain how the block assertions
are certified. This separates coverage of all input trees from the
arithmetic checks for a particular block.

\subsubsection{Complete bounded case lists}

Fix $3\leq b\leq n-2$. Let $\widehat{\mathcal O}_{n,b}$ consist
of all short-branch profiles $\mathbf h$ with
\[
 b\leq |\mathbf h|\leq\min(2b-2,n-2),
 \qquad o(\mathbf h)>0\quad\hbox{if $n$ is even}.
\]
Here and below the short-branch condition includes the absence of
a subprofile of total $b$. Unlike $\mathcal O_{n,b}$, this family
does not impose minimality. For even $n$, choose a nominal center
$a$ as in Theorem~\ref{thm:peeled-prefix-step}, and put
\[
 E=2\floor{(b-1)/2},\qquad B=\max(b,n-a-1),\qquad H=B+E-1.
\]
Let $\widehat{\mathcal E}_{n,b}(a)$ consist of all short-branch
profiles having only even sizes and satisfying
\[
 b\leq |\mathbf h|\leq\min(H,n-2).
\]
These enlarged lists allow the same bounded subset-sum enumeration
to be used throughout the range.

\begin{proposition}[Sufficient complete lists]
\label{prop:complete-peeling-lists}
Put $\ell=n-b-1$ and suppose $\cP([m])$ holds for
$0\leq m\leq\ell$. Suppose a detached block satisfies the first
condition of Theorem~\ref{thm:peeled-prefix-step}, and that
\begin{enumerate}
\item $\mathcal R_n(\mathbf h)$ holds for every
      $\mathbf h\in\widehat{\mathcal O}_{n,b}$;
\item $\cP(I(s,j))$ holds for every $1\leq j<b$ and
      $1\leq s\leq n-j+1$;
\item if $n$ is even, every interval in the preceding item of
      even length admits every prescribed root at an odd label.
\end{enumerate}
For even $n$, also require
$\mathcal D_n(\mathbf h;a(\mathbf h))$ for every
$\mathbf h\in\widehat{\mathcal E}_{n,b}(a)$ of total below $B$,
and $\mathcal D_n(\mathbf h;a)$ for every remaining profile
in that family. The choice of $a$ must satisfy all the arithmetic
hypotheses of Theorem~\ref{thm:peeled-prefix-step}, including
$n-H-1>p$ when $a=2p$. Then $\cP([n])$ holds.
\end{proposition}

\begin{proof}
The ordinary family contains $\mathcal O_{n,b}$: it has the same
order and parity bounds and merely omits the minimality test.
The even family contains both $\mathcal E^-_{n,b}(a)$ and
$\mathcal E^+_{n,b}(a)$. Indeed, $H\geq B$, all low totals are
below $B$, and the high family has upper bound $H$. The same
bound $n-2$ is imposed throughout. Its low profiles are checked
with their individual centers, and its high profiles with $a$,
as required. Finally, the interval list contains every member
of $\mathcal I_{n,b}$; its rooted part includes every even length
with odd start. Thus every hypothesis of
Theorem~\ref{thm:peeled-prefix-step} holds.
\end{proof}

For completeness, the profile enumeration can be described without
reference to any stored list. Process the possible sizes $j$ in
increasing order. Given the chosen multiplicities and the set $R$
of attainable sums at most $b$, adjoining $h_j$ copies of $j$
replaces $R$ by
\[
 \{r+uj:r\in R,\ 0\leq u\leq h_j,\ r+uj\leq b\}.
\]
Start with $R=\{0\}$ and reject a choice as soon as $b\in R$.
Bound each multiplicity by the remaining total divided by $j$.
At the final size, test the lower bound and the stated parity
condition. Every bounded profile has a unique sequence of
multiplicities, so this procedure enumerates exactly the displayed
families. The checker requires equality of the regenerated and
supplied case sets; it rejects missing and repeated cases.

\subsubsection{Certificates for the arithmetic blocks}

A \emph{block model} is a graph $G$ on $j$ indexed positions,
together with a nonempty allowed set $R$. Its assertion is
$\mathcal U_1(G,R)$. For an ordered label inventory
$L=\{x_1,\ldots,x_j\}$ and a contact value $d$, a block record
supplies a permutation $\rho$ of the positions. We check
\begin{equation}\label{eq:continuation-model-binding}
 \begin{split}
 uv\in E(G)&\quad\Longleftrightarrow\quad
                 \gcd(x_{\rho(u)},x_{\rho(v)})=1
                 \quad(u\ne v),\\
 u\in R&\quad\Longleftrightarrow\quad
                 \gcd(x_{\rho(u)},d)=1.
 \end{split}
\end{equation}
Thus one certified model can serve several arithmetic inventories,
but the entire graph and allowed set are checked for each use.
For a plain interval we take $d=1$; for an odd-root interval we
take $d=2$. For a retained branch, $d$ is its proposed center label.
The value $d=1$ serves only to specify an unrestricted root menu;
it need not denote a separate contact vertex in that application.

There are two ways of proving the model assertion. The first gives
a complete rooted-tree table, with a permutation in each row whose
root image lies in $R$. The second binds the model to a recursive
host certificate of Section~\ref{sec:finite-certificates}. Only
split hosts, complete marked tables, and centroid reductions are
needed here.

We specify the centroid choices to make the treatment of inherited
marks explicit. A node asserts $\mathcal U_k(H,O)$. For each case
$(h;(j_i,t_i)_i)$ in \eqref{eq:centroid-mark-cases}, its record
gives a center $c$, belonging to $O$ if $h=1$, and an ordered
partition of the other positions into sets $S_i$ of sizes $j_i$.
The child graph must be exactly $H[S_i]$ in the recorded order.
The following choices suffice:
\begin{enumerate}
\item If $c$ is universal in $H$, use $O\cap S_i$ with mark
      budget $t_i$.
\item If $c$ is adjacent to every member of $O\setminus\{c\}$,
      use $O\cap S_i$ with budget $\min(j_i,t_i+1)$.
\item If $k=0$, use $N_H(c)\cap S_i$ with budget $1$.
\end{enumerate}
The general choice $N_H(c)\cap O\cap S_i$ when $t_i>0$, and
$N_H(c)\cap S_i$ when $t_i=0$, is also valid with budget
$\min(j_i,t_i+1)$. It is accepted by the checker but is not
needed in the certificates for this range. A larger child budget
is permitted. If the child's designated set is all its positions,
an unmarked child assertion suffices.

The centroid case list is regenerated as all multisets of pairs
$1\leq j\leq\floor{m/2}$, $0\leq t\leq j$, with total size
$m-1$ and total mark count $k-h$, separately for $h=0,1$ with
$h\leq k$. Listing pairs in nondecreasing lexicographic order
removes repetitions. These are exactly the cases in
\eqref{eq:centroid-mark-cases}. In particular, equal component
orders do not justify discarding different distributions of marks.

\begin{proposition}[Soundness of the block checks]
\label{prop:continuation-block-soundness}
Suppose every terminal table has its complete index set and valid
witnesses, every split node meets Lemma~\ref{lem:split-host},
and every centroid node passes the checks above. Suppose the dependency
graph is acyclic. Then every node asserts its stated
$\mathcal U_k(H,O)$.

A model $(G,R)$ is certified either by a complete rooted table
or by an edge-preserving bijection $\sigma\colon V(H)\to V(G)$
from a proved node. If $R\ne V(G)$, require $\sigma(O)\subseteq R$
and either $k\geq1$ or $O=V(H)$. If $R=V(G)$, any proved mark
budget suffices. Every arithmetic block passing
\eqref{eq:continuation-model-binding} then satisfies
$\mathcal U_1(H_L,O_d(L))$.
\end{proposition}

\begin{proof}
Induct through the dependency graph. Complete marked tables
and split nodes are valid by
Proposition~\ref{prop:terminal-witnesses} and
Lemma~\ref{lem:split-host}. At a centroid node, the first choice
above is the first alternative of
Proposition~\ref{prop:marked-centroid-host}. In the second choice,
$O\cap S_i$ is contained in $N_H(c)$ because $c\notin S_i$.
It therefore handles both the inherited marks and the attachment.
The third choice has no inherited marks. The general intersection
choice likewise meets the second alternative of that proposition.
The regenerated case list is complete, so the proposition proves
the parent assertion. This establishes the induction.

A complete rooted table proves the model assertion directly.
For a transferred model with a proper allowed set, the source
assertion permits one mark: either $k\geq1$, or its designated
set is unrestricted and imposes no condition. The inclusion
$\sigma(O)\subseteq R$ then permits
Lemma~\ref{lem:host-transfer}. If $R$ is unrestricted, forget
any source marks, transfer $\mathcal U_0(H)$, and observe that
the root condition in $G$ is automatic. Thus in every case the
model satisfies $\mathcal U_1(G,R)$. Applying the same transfer
lemma to $u\mapsto x_{\rho(u)}$ proves the arithmetic assertion.
\end{proof}

The detached block needs only $\cP$. In the direct format it has
a complete unrooted table. In the recursive format its labels
are bijected to a proved host, with every gcd edge checked, and
any mark condition is forgotten. This use does not require a
rooted conclusion.

\subsubsection{Assembly through order one thousand}

\begin{theorem}\label{thm:initial-prefix-1000}
Every forest of order $n$ with $1\leq n\leq1000$ has a coprime
labeling by $[n]$.
\end{theorem}

\begin{proof}
We induct on $n$, from $300$ through $1000$, using
Theorem~\ref{thm:initial-prefix-299} as the base.
Whenever one of $n-2,n-1,n$ is prime, use
Lemma~\ref{lem:prime-contact-extension}; its required order
$p-1$ is smaller than $n$. There are $302$ such steps.
For the other $399$ orders use
Proposition~\ref{prop:complete-peeling-lists}. The supplement
contains one certificate for each of these orders, and its
checker reconstructs this exact division by trial-division
primality tests. No order is omitted or accepted from a saved
range report.

For each branch-peeling step, the checker reads $b$, the power
of two $c$, and, when used, the odd prime $p$. It recomputes
$\ell=n-b-1$, the nominal center $a$, and the bounds $B,H$.
In the exchange case it checks
\[
 2p\leq n<3p,\qquad c-p\geq b,\qquad n-H-1>p.
\]
It verifies that the detached inventory is exactly
$\{\ell+1,\ldots,n\}\setminus\{q_0\}$ and that $q_0$ is
universal in $[n]$. Its block assertion is proved by the direct
or recursive checks just described.

For every regenerated ordinary profile, the recorded center is
universal and the branch blocks give exactly
\eqref{eq:retained-core-inventory}. For every even profile, the
parent is universal and the blocks give exactly
\eqref{eq:removed-core-inventory}; the center must equal $a$
when the total is at least $B$. Each branch model is bound to
its actual block and contact by
\eqref{eq:continuation-model-binding}. Finally, every interval
of every required length and start is checked, with the odd-root
condition for every even length when $n$ is even. Hence all the
arithmetic assertions of Proposition~\ref{prop:complete-peeling-lists}
hold. Its remaining hypothesis concerns only orders at most
$\ell<n$, already proved by induction.

Table~\ref{tab:continuation-counts} records the complete finite
check. The block proofs include $6567024$ rooted witnesses and
$13274$ detached-block witnesses in the direct format. Recursive
hosts contribute $1335600$ marked witnesses. Thus $7915898$
terminal bijections are checked, with $81883548$ terminal edge
checks. The recursive proof graphs contain $1985$ split nodes,
$10045$ table nodes, and $17163$ centroid nodes. Their complete
case lists have $3211520$ cases: $3079520$ use a universal center,
$131687$ use the inherited allowed set, and $313$ have no
inherited marks and use the center's neighborhood.

All table, coverage, inventory, and dependency checks pass.
Proposition~\ref{prop:continuation-block-soundness} proves the
block assertions, and Proposition~\ref{prop:complete-peeling-lists}
proves each remaining order. Together with the prime-contact
steps, this establishes the induction through $1000$.
\end{proof}

\begin{table}[t]
\centering
\begin{tabular}{lrr}
Checked data&Orders $300$--$537$&Orders $538$--$1000$\\[3pt]
Prime-contact steps&$105$&$197$\\
Branch-peeling steps&$133$&$266$\\
Ordinary profiles&$6832$&$28763$\\
Even profiles&$17060$&$8792$\\
Interval assertions&$373297$&$1471354$\\
Retained-branch assertions&$201838$&$168542$\\
Terminal witnesses&$6580298$&$1335600$\\[2pt]
Block proof format&Rooted tables&Recursive hosts
\end{tabular}
\par\vspace{.5\baselineskip}
\caption{Complete checks for the continuation. Counts are summed
over the separate order certificates; shared models are checked once
within each certificate and bound anew to every arithmetic inventory.}
\label{tab:continuation-counts}
\end{table}

The last step gives a small example of the construction. For
$n=1000$, take $b=3$, $c=512$, and $p=499$, so $a=998$,
$B=3$, and $H=4$. The detached inventory is
$\{998,999,1000\}$, with universal contact $997$ and exterior
$[996]$. It labels the three-vertex path by placing $999$ at
its middle; forest completion gives its unmarked assertion.
There are no ordinary profiles: using a branch of order one
would create a subprofile of total three as soon as the total
reached three. The only even retained profile consists of two
branches of order two. Its exact partition is
\[
 \{995,996,997,998,999,1000\}
 =\{997,998\}\sqcup\{995,996\}\sqcup\{999,1000\}.
\]
Each pair is consecutive, and its attachment may receive the
odd label $995$ or $999$, respectively, both coprime to $998$.
All omitted intervals have length two and can put their attachment
at the odd endpoint. Finally, $512-499=13\geq3$ and
$1000-H-1=995>499$, so the center exchange applies whenever
it is needed. This checks the hypotheses at the endpoint without
any new tree table.

The supplement includes the exact case inventories, every model
and host proof, and all permutations used in their arithmetic
transfers. The checker uses integer arithmetic throughout.
Its default run checks the initial certificate through $299$
before checking this continuation. The separate continuation
run treats Theorem~\ref{thm:initial-prefix-299} as its explicit
base. We next extend the host rules while keeping the same
arithmetic case lists.

\subsubsection{Reusing host proofs and removing leaves}

A block need not have a new centroid decomposition if a proved
host of the same order already embeds in it. By
Lemma~\ref{lem:host-transfer}, a bijection from a source host
$H'$ to a target host $H$ suffices whenever it preserves every
source edge and, for a marked assertion, sends the entire source
allowed set into the target allowed set. The source mark budget
must be at least the target budget. In the certificates below,
we additionally require $e(H')<e(H)$. This last condition is not
needed for the labeling implication; it gives a decreasing
parameter for the proof graph.

The following reductions are the local extensions of
Lemma~\ref{lem:host-local-extensions}, applied inside a larger
host. We state the inherited mark conditions explicitly, since
a centroid reduction may have introduced several marks.

\begin{lemma}[Leaf reductions with inherited marks]
\label{lem:inherited-leaf-reductions}
Let $H$ be a graph of order $m\geq2$, let $O\subseteq V(H)$,
and let $q$ be universal in $H$.
\begin{enumerate}
\item For $0\leq k\leq m$, suppose $q\in O$ whenever $k>0$.
      If
      \[
       \mathcal U_{\min(k,m-1)}(H-q,O\setminus\{q\})
      \]
      holds, then $\mathcal U_k(H,O)$ holds.
\item Suppose $m\geq3$, $k\in\{0,1\}$, and $a\ne q$ is a
      vertex of $H$. Require $q\in O$ if $k=1$. If
      \[
       \mathcal U_k(H-\{a,q\},O\setminus\{a,q\})
      \]
      holds, then $\mathcal U_k(H,O)$ holds.
\end{enumerate}
In either assertion, a larger child mark budget suffices. If
the child designated set is its entire inventory, an unmarked
child assertion suffices.
\end{lemma}

\begin{proof}
For the first assertion, choose a leaf $x$ of the input tree
and give it the label $q$. The remaining tree has at most
$\min(k,m-1)$ marks, so apply the child assertion to it.
If $x$ was marked, its label is allowed. Its unique incident
edge is valid because $q$ is universal. The two inventories
are disjoint and exhaust the host.

For the second assertion, choose an unmarked leaf $x$. Such
a leaf exists: the endpoints of a longest path are leaves,
and at most one vertex is marked. Let $y$ be its neighbor.
Give $x$ the label $a$ and $y$ the label $q$. Label the
remaining forest by the child assertion, retaining all marks
that remain. If $y$ was marked, its label belongs to $O$;
every other marked vertex receives a label in
$O\setminus\{a,q\}$. The only edge incident to $x$ joins
it to $y$. Every other edge leaving the remaining forest
also meets $y$, whose label is universal. The child assertion
checks the internal edges, and the disjoint inventories
give bijectivity.

Forest completion gives the same assertions for forests.
The observations about larger budgets and unrestricted
designated sets follow from the definition of $\mathcal U_k$.
\end{proof}

The restriction to at most one mark in the second assertion
ensures that an unmarked leaf can be chosen. A path whose two
endpoints are marked has no such leaf. Thus this rule cannot
be used at a node with two marks merely because its child
certificate allows two marks.

\begin{samepage}
\begin{proposition}[Soundness of the enlarged host system]
\label{prop:extended-host-soundness}
Add the following rules to the block system of
Proposition~\ref{prop:continuation-block-soundness}: a same-order
transfer satisfying Lemma~\ref{lem:host-transfer} and having
strictly fewer source edges, and either leaf reduction of
Lemma~\ref{lem:inherited-leaf-reductions}. Require every child
inventory and designated set in a leaf reduction to be the
exact inherited ones. Then every completely checked node
asserts its stated $\mathcal U_k(H,O)$, and the arithmetic
binding rule \eqref{eq:continuation-model-binding} remains valid.
\end{proposition}
\end{samepage}

\begin{proof}
Order host assertions by the lexicographic pair
\[
 (\abs{V(H)},e(H)).
\]
A centroid child has at most half the parent order. A leaf
reduction decreases the order by one or two. A transfer
preserves the order and strictly decreases the number of
edges. Thus every dependency has smaller rank. Since both
coordinates are nonnegative integers, induction on this
rank is well founded.

The terminal, split, and centroid nodes are justified as in
Proposition~\ref{prop:continuation-block-soundness}. At a transfer
node apply Lemma~\ref{lem:host-transfer}, checking the whole
source allowed set when the target budget is positive.
At a leaf node apply Lemma~\ref{lem:inherited-leaf-reductions}
to the exact induced child graph. Its inherited designated
set and sufficient mark budget are part of the checked
hypotheses. These implications prove the induction. The
arithmetic binding is then the same composition of bijections
as before.
\end{proof}

The extra host rules change neither the retained profiles
nor the integer labels required by a branch-peeling step.
Every block still carries its full root assertion. In
particular, a transfer with no marks may ignore the designated
sets, but this does not make it a valid transfer of a rooted
block. The block's final use must separately pass the root
condition in \eqref{eq:continuation-model-binding}.

\subsubsection{Extension through order 1352}

\begin{theorem}\label{thm:initial-prefix-1352}
Every forest of order $n$ with $1\leq n\leq1352$ has a coprime
labeling by $[n]$.
\end{theorem}

\begin{proof}
Start from Theorem~\ref{thm:initial-prefix-1000} and proceed
through all $352$ integers from $1001$ to $1352$. Trial-division
primality checks identify $136$ orders of the form $p$, $p+1$,
or $p+2$. For each, Lemma~\ref{lem:prime-contact-extension}
applies to the already proved order $p-1$.

For every other order the supplement supplies one branch-peeling
certificate. Its acceptance conditions for the detached block,
retained profiles, interval assertions, initial-interval
dependencies, and center exchange are exactly those used in
the proof of Theorem~\ref{thm:initial-prefix-1000}. They are
rechecked for the new numerical inventories. Only the host
proof system is enlarged, as in
Proposition~\ref{prop:extended-host-soundness}. All case sets
are regenerated, every arithmetic transfer is checked, and
every referenced host node is verified. Hence
Proposition~\ref{prop:complete-peeling-lists} applies at each
of these $216$ orders. Its required initial intervals have
orders at most $n-b-1<n$, so induction supplies them.

Table~\ref{tab:extended-prefix-checks} gives the complete
division of the certificates by their host rules. The check
covers $163165$ retained profiles, $2266427$ interval assertions,
and $867973$ retained-branch assertions. There are $3902347$
terminal witnesses and $21462240$ centroid cases. All terminal
witnesses, all child inventories and mark budgets, all $11323$
same-order transfers, and all $63825$ leaf-pair reductions pass.
Each transfer preserves every required edge and decreases
the edge count; each leaf-pair node has at most one mark.
The one-vertex leaf rule is supported by the checker but is
not needed in these certificates.

The index contains exactly the non-prime-contact orders in
the range, and the induction visits every integer in order.
The fresh verification therefore establishes all its premises,
rather than importing an accepted range from an earlier report.
Proposition~\ref{prop:extended-host-soundness} supplies the block
assertions, Proposition~\ref{prop:complete-peeling-lists} gives
the next initial interval, and induction proves the theorem.
\end{proof}

\begin{table}[htbp]
\centering
\begin{tabular}{lrrr}
Host rules used&Orders&Terminal rows&Centroid cases\\[3pt]
Split, table, and centroid&$91$&$798609$&$4188144$\\
Also same-order transfers&$8$&$379985$&$2643819$\\
Also leaf-pair reductions&$117$&$2723753$&$14630277$\\[2pt]
Total&$216$&$3902347$&$21462240$
\end{tabular}
\par\vspace{.5\baselineskip}
\caption{Complete block checks for orders $1001$ through $1352$.
The remaining $136$ orders use the prime-contact lemma and need
no block certificate. Shared host states are counted separately
when they occur in different order certificates.}
\label{tab:extended-prefix-checks}
\end{table}

All calculations in this continuation use integers. The
supplement includes the $216$ order certificates and a checker
that shares the already specified arithmetic checks while
adding only the new host rules. Its default invocation also
rechecks the earlier initial range. Section~\ref{sec:finite-endpoint} completes the remaining
finite steps using minimal retained cores and the larger
bipartite hosts.

\subsection{Completing the finite prefix}\label{sec:finite-endpoint}

The last seven orders of the finite argument use the minimal
cores of Section~\ref{sec:branch-peeling}. Their larger branch
blocks are supplied by the bipartite hosts already proved.
Two further marked reductions allow these hosts to be used
inside the recursive certificates. We prove those reductions,
specify the source graphs, and then check the seven induction
steps together.

\subsubsection{Deleting marked vertices}

Deleting a prescribed vertex reduces the number of marks by
one, although it may leave several components. This is why
our host assertions are formulated for forests as well as trees.

\begin{lemma}\label{lem:endpoint-marked-deletions}
Let $H$ have order $m$, let $O\subseteq V(H)$, and let
$u\in O$ be universal in $H$.
\begin{enumerate}
\item If $m\geq2$, $1\leq k\leq m$, and
\[
 \mathcal U_{k-1}(H-u,O\setminus\{u\})
\]
holds, then $\mathcal U_k(H,O)$ holds.
\item Suppose $m\geq4$, choose $a\ne u$, and assume
\[
 \mathcal U_2(H-\{u,a\},O\setminus\{u,a\}).
\]
If $H$ has a spanning path with both endpoints in $O$,
then $\mathcal U_2(H,O)$ holds.
\end{enumerate}
\end{lemma}

\begin{proof}
For the first assertion, take a forest with at most $k$ marks.
Choose a marked vertex $v$ if one exists, and otherwise
choose any vertex. Delete $v$, embed the remaining forest
using the child assertion, and send $v$ to $u$. At most
$k-1$ marks remain, including when the original marked set
was empty. The inventories are disjoint and exhaustive.
Every restored edge meets $u$, and a mark on $v$ has an
allowed image. This proves the first assertion.

For the second, it suffices to consider a tree, by
Lemma~\ref{lem:forest-completion}. If it has an unmarked
leaf $x$, delete $x$ and its neighbor $y$. Embed the
remaining forest using the assumed two-mark assertion,
and send $x$ to $a$ and $y$ to $u$. All remaining marks
are respected; a mark on $y$ receives $u\in O$. Every
edge meeting a deleted vertex has an endpoint at $y$,
so all restored edges are allowed.

Otherwise every leaf is marked. A tree has at least two
leaves and there are at most two marks, so there are exactly
two leaves, both marked. The identity
\[
 \#\{v:d(v)=1\}
   =2+\sum_{v:d(v)\geq3}(d(v)-2)
\]
shows that no vertex has degree greater than two. Hence the
tree is a path, with precisely its endpoints marked.
Map it in order to the assumed spanning path of $H$.
The endpoint condition respects both marks and the path
edges respect every tree edge. This completes the proof.
\end{proof}

The path condition in the second assertion is checked by
one permutation of the entire host inventory. Checking only
the child with two marks would leave the marked-endpoint
path case untreated.

\subsubsection{Source graphs for the larger blocks}

We use the graph $J_t(X)$ defined in
Section~\ref{sec:larger-blocks}: it is obtained by adjoining
$t$ universal vertices to $X$. For a graph $X$ and $t\geq1$,
write $R_t(X)$ for the graph on
\[
 V(X)\sqcup U\sqcup\{\zeta,e\},\qquad |U|=t,
\]
whose edges are those of $X$, all edges meeting $U$, and
all edges from $\zeta$ to $V(X)$. Thus every member of
$U$ is universal, $\zeta$ is adjacent to every other
vertex except $e$, and $e$ is adjacent only to $U$.
Write $R_t^+(X)$ when the edge $\zeta e$ is also added.
In the latter graph $\zeta$ is universal as well.

The subscripts on the arithmetic hosts in the next statement
refer to their own defining sections.

\begin{proposition}\label{prop:endpoint-source-hosts}
Each of the following source graphs admits a bijective
embedding of every forest of its order, sending any specified
vertex to any specified added universal vertex:
\begin{enumerate}
\item $R_4(X_0)$ and $R_4^+(X_0)$, of order $26$, for the
      sparse host $X_0$ of Section~\ref{sec:reserved-blocks};
\item $J_4(X_1)$ and $R_3(X_i)$, $i=2,3,4$, of order $26$,
      for the other four hosts of that section;
\item $J_5(X)$, of order $28$, for the $9$-by-$14$ host
      $X$ in Theorem~\ref{thm:rooted-twentyeight-block};
\item $J_6(X_i)$, $i=1,2,3$, of order $30$, for the first
      three hosts of Section~\ref{sec:larger-blocks}.
\end{enumerate}
Consequently, an edge-preserving bijection from one of these
sources to a host $H$, taking one of its added universal
vertices into $O$, proves $\mathcal U_1(H,O)$.
\end{proposition}

\begin{proof}
Proposition~\ref{prop:reserved-finite-hosts} supplies the
covering assertions for the first two items.
Proposition~\ref{prop:reserved-label-host} then gives their
rooted embeddings. Adding the edge $\zeta e$ preserves
all these embeddings. In $R_4^+(X_0)$, exchanging $\zeta$
with a member of $U$ preserves every edge, since both are
universal; hence $\zeta$ is also an allowed prescribed image.

The host in the third item covers $(9,14,6)$, by the complete
verification in Theorem~\ref{thm:rooted-twentyeight-block}.
Apply Proposition~\ref{prop:bipartite-universal-extension}
with $b=9$, $k=4$, and $D=6$. The first three hosts in the
fourth item cover $(9,15,6)$, by
Proposition~\ref{prop:larger-finite-hosts}. Apply the same
extension with $b=9$, $k=5$, and $D=6$.
These constructions require only the stated bipartite edges
and the added universal vertices.

For the consequence, fix an added universal vertex whose
image lies in $O$. Embed the input forest with its marked
vertex there, and compose with the supplied bijection.
Every edge is preserved, and the mark lies in $O$.
If no vertex is marked, designate an arbitrary one.
\end{proof}

A source certificate records that entire bijection. The check
reconstructs the bipartite source from its stated inventory,
then checks every source edge and the distinguished root
image. Extra edges in the target cause no difficulty.
The covering assertions themselves are supplied by the
complete certificates of Sections~\ref{sec:bipartite-hosts}--\ref{sec:larger-blocks};
a graph name or a previously saved acceptance report does
not replace them.

\subsubsection{Soundness of the endpoint certificates}

\begin{proposition}\label{prop:endpoint-certificate-soundness}
Enlarge the system of Proposition~\ref{prop:extended-host-soundness}
by the two rules of Lemma~\ref{lem:endpoint-marked-deletions},
the source transfers of Proposition~\ref{prop:endpoint-source-hosts},
and the following previously proved extensions: a fixed-root
host with an added leaf, an unmarked host with an added
universal vertex, and a fixed-root host with an attached
triangle. Every accepted host assertion in this enlarged
system is valid. The arithmetic model binding of
\eqref{eq:continuation-model-binding} remains valid.
\end{proposition}

\begin{proof}
Use induction on $(|V(H)|,e(H))$, ordered lexicographically,
as in Proposition~\ref{prop:extended-host-soundness}.
A source transfer is a terminal implication from
Proposition~\ref{prop:endpoint-source-hosts}; it does not
create a dependency on another host of the same order.
Every other new rule has a child of smaller order.

For marked deletion, check the universal image in $O$,
the exact induced child graph and inherited designated set,
and the child budget $k-1$. For two-mark deletion, check
these conditions with two deleted positions and child
budget two. Also check that the extra path record uses
every position exactly once, has both endpoints in $O$,
and preserves every consecutive edge. The lemma then
applies. A larger child budget suffices, as does an
unrestricted child designated set.

The single-leaf and universal-vertex extensions are
Lemma~\ref{lem:host-local-extensions}(4) and (2), respectively.
The triangle extension is Lemma~\ref{lem:host-attached-triangle}.
In the leaf and triangle rules the child has a singleton
designated set and a budget of at least one, and the
conclusion is unmarked. In the universal-vertex rule the
child need only be unmarked; if the conclusion has a mark,
its new universal position must lie in $O$. Each rule
checks the full residual inventory and all required source
edges. These checks supply exactly the hypotheses of the
cited lemmas.

All recursive assertions therefore follow by induction.
The final arithmetic check is unchanged: compose the
host-to-model and model-to-label bijections, preserving
all edges and the entire required root menu.
\end{proof}

\subsubsection{The seven endpoint orders}

For each endpoint order, set
\begin{equation}\label{eq:endpoint-stumps}
 n\in\{1353,\ldots,1359\},\quad
 b=n-1327,\quad \ell=1326,\quad q_0=1327.
\end{equation}
Trial division proves that $1327$ is prime, and $2\cdot1327>1359$.
It is therefore universal in each $[n]$ under consideration.
The detached inventory is always $B_0=\{1328,\ldots,n\}$.
Its forest assertion is checked separately for each order.

\begin{theorem}\label{thm:initial-prefix-1359}
Every forest of order $n$, $1\leq n\leq1359$, has a
coprime labeling by $[n]$.
\end{theorem}

\begin{proof}
Theorem~\ref{thm:initial-prefix-1352} supplies the initial
range. For each of the seven remaining orders we apply
Theorem~\ref{thm:peeled-prefix-step} with
\eqref{eq:endpoint-stumps}. In fact its required smaller
intervals all have orders at most $1326$, so they are
already covered by the initial range.

First regenerate the minimal families
$\mathcal O_{n,b}$ and, for even $n$,
$\mathcal E^-_{n,b}(a)\cup\mathcal E^+_{n,b}(a)$.
Use the bounded subset-sum enumeration and the test
\eqref{eq:core-minimality}. Every low even profile is
retained without a minimality filter. Require equality
between these families and the supplied core records.
Each record partitions the entire inventory into its
initial interval, its contact label or labels, and its
branch blocks. Each branch assertion is bound to its
exact gcd graph and root menu.

At the even orders $1354$, $1356$, and $1358$, use the
center exchange with
\[
 c=1024,\qquad p=677,\qquad a=2p=1354.
\]
Trial division proves that $677$ is prime. We have
$2p\leq n<3p$ and $c-p=347\geq b$.
Here $B=b$ and $H=2b-2$, so the three values of
$n-H-1$ are $1301$, $1299$, and $1297$, all greater than
$p$. Thus all center-exchange hypotheses hold. There
are no low even profiles in these three steps. The odd
orders require no even-family assertion.

For orders $1353$ through $1358$, check every plain
interval of length $1\leq j<b$ and every permitted start.
At even order also check every even-length interval with
odd start, with the root required to receive an odd label.
For order $1359$, keep the same plain requirements except
at length $31$, where the starts are
\[
 26,57,\ldots,1266.
\]
These are exactly the $41$ members of
$\mathcal A_{1359,31}$. Lemma~\ref{lem:largest-odd-core}
retains two branches of order $31$ whenever any such branch
is omitted, and Corollary~\ref{cor:largest-odd-intervals}
therefore supplies precisely this shortened list. No even
root condition is needed at the odd order $1359$.

Table~\ref{tab:finite-endpoint-checks} records the complete
case counts. The supplement includes all seven certificates.
Its checker regenerates the case families just described,
checks every retained inventory and arithmetic transfer,
and verifies every referenced host assertion by
Proposition~\ref{prop:endpoint-certificate-soundness}.
All detached-block, core, and interval assertions pass.
They are exactly the hypotheses of
Theorem~\ref{thm:peeled-prefix-step}, or supersets of its
interval hypotheses. That theorem proves each of the
seven orders. Together with the initial range this proves
the assertion.
\end{proof}

\begin{table}[htbp]
\centering
\renewcommand{\arraystretch}{1.12}
\begin{tabular}{rrrrr}
$n$&$b$&Ordinary cores&Even cores&Interval assertions\\[3pt]
$1353$&$26$&$1179$&$0$&$33525$\\
$1354$&$27$&$1341$&$234$&$43602$\\
$1355$&$28$&$1689$&$0$&$36234$\\
$1356$&$29$&$1921$&$305$&$46991$\\
$1357$&$30$&$2396$&$0$&$38947$\\
$1358$&$31$&$2737$&$390$&$50385$\\
$1359$&$32$&$3381$&$0$&$40376$
\end{tabular}
\par\vspace{.5\baselineskip}
\caption{The seven endpoint certificates. Plain and odd-root
assertions on the same interval are counted separately.
Only the length-$31$ list at order $1359$ is shortened.}
\label{tab:finite-endpoint-checks}
\end{table}

The supplement's default command rechecks the nine covering
certificates and the earlier initial range before verifying
these seven steps. Its endpoint-only command checks the
implication from those established premises. Neither reads
a saved acceptance report. The finite prefix is now complete;
it remains to establish every order $n\geq1360$.

\section{Interval certificates through order 34999}
\label{sec:part-intervals}

This part covers the orders from $1360$ through $34999$ by certificates
on intervals of orders.
Sections~\ref{sec:local-probability}--\ref{sec:directed-probability}
provide the ingredients: local probability bounds, type assignments
controlled by decomposition height and in the original forest, finite
arithmetic bounds, and the directed rational evaluation of
probabilities. Section~\ref{sec:interval-construction} states the
certificate criterion on an integer interval
(Theorem~\ref{thm:integer-interval-certificate}), and
Sections~\ref{sec:local-certificates}
and~\ref{sec:hall-certificates} verify it for the orders $1360$--$1521$.
Sections~\ref{sec:forest-continuation}--\ref{sec:single-range-continuation}
continue through $3499$ with modified transport and Hall estimates, and
Section~\ref{sec:retained-continuation} continues through $34999$
with retained positions.

\subsection{Sharper local probability bounds}\label{sec:local-probability}

We return to the matching construction. The local estimate in
Corollary~\ref{cor:local-conflict-tails} uses only a degree bound and
an expectation. For the finite parameter ranges we can also retain
the number of conflicting labels, the degrees in the original forest,
and the fact that the deletion counts are integers. We establish these
refinements before applying them to numerical intervals.

The probability space is the one in Section~\ref{sec:matching-estimates}:
types and tags are fixed, and the bijections within the classes are
independent and uniform. A fixed set of $b_I$ half-labels in $\Lambda_I$
therefore occupies a uniform $b_I$-subset of $A_I$. This observation
allows us to regard the local conflict count as a sum sampled without
replacement from the class degree lists.

\subsubsection{Comparison with a binomial variable}

The comparison between sampling with and without replacement is due
to Hoeffding~\cite[Theorem~4]{Hoe63}. We include the argument, followed
by the two elementary convexity steps that give the form needed here.
For a real number $x$, write $x_+=\max(x,0)$.

\begin{lemma}[Convex comparison]\label{lem:class-convex-comparison}
Let $g\geq1$ be an integer. For each $i\in[g]$, let
$x_{i1},\ldots,x_{im_i}\in[0,d]$, where
$m_i\geq1$ and $d>0$, and fix an integer $0\leq b_i\leq m_i$.
Independently in each class, choose a uniform subset $S_i\subseteq[m_i]$
of size $b_i$, and put
\[
 X=\sum_{i=1}^g\sum_{j\in S_i}x_{ij},\qquad
 b=\sum_i b_i,\qquad m=\E X.
\]
If $b>0$ and $Z$ is binomial with parameters $b,m/(db)$, then
\begin{equation}\label{eq:class-convex-comparison}
 \E\phi(X/d)\leq\E\phi(Z)
\end{equation}
for every convex function $\phi\colon\R\to\R$. More generally, if
$B\geq\max(1,b)$ is an integer and $M\geq m$, put
$p=\min(1,M/(dB))$. Then \eqref{eq:class-convex-comparison} holds
with $Z$ binomial with parameters $B,p$ for every increasing convex
$\phi$. If $b=0$, then $X=0$ surely.
\end{lemma}

\begin{proof}
Divide all entries by $d$. We compare the resulting sum with three
successive sums of independent variables, preserving its mean until
the final possible enlargement of the Bernoulli parameters.

\emph{Step 1. Sampling with replacement.}
Consider one population $x_1,\ldots,x_a$ and a sample of size $r\leq a$.
For an integer vector $c\geq0$ with $\sum_jc_j=r$, define
\[
 F(c)=\frac1{a!}\sum_{\pi}\phi\left(\sum_jc_jx_{\pi(j)}\right),
\]
where the sum is over all permutations. Suppose two coordinates of
$c$ are $u\geq2$ and $0$. Pair the permutations by interchanging
the population entries $x,y$ in those two positions. If the contribution
from the other positions is $s$, the two arguments are $s+ux$ and
$s+uy$. After replacing $(u,0)$ by $(u-1,1)$, the arguments become
\[
 s+(u-1)x+y,\qquad s+x+(u-1)y.
\]
They are convex combinations of the original arguments, with
coefficients $1-1/u$ and $1/u$ in opposite order. Hence the sum of
their $\phi$-values does not increase.

As long as some coordinate is at least two, there is a zero coordinate,
since $r\leq a$. Repeating the replacement terminates: the integer
$\sum_jc_j^2$ strictly decreases. The final vector has $r$ ones and
all other coordinates zero. Its value of $F$ is the expectation for
sampling without replacement. In sampling with replacement, condition
on the unordered multiplicities of the chosen population positions.
Their assignment to positions is uniform, so the conditional expectation
is the corresponding $F(c)$. Averaging proves that replacement can
only increase the convex expectation. We may condition on the sum
from all other classes and apply the same argument to one class at
a time.

\emph{Step 2. Bernoulli variables.}
After Step 1, all $b$ draws are independent and belong to $[0,1]$.
For $0\leq w\leq1$ and any $s$,
\[
 \phi(s+w)\leq(1-w)\phi(s)+w\phi(s+1).
\]
Conditioning on the other draws and averaging over $w$ replaces
this draw by an independent Bernoulli variable of the same mean
without decreasing the expectation. Repeating gives a sum of $b$
independent Bernoulli variables whose parameters have sum $m/d$.

\emph{Step 3. Equal parameters.}
For two Bernoulli parameters $u,v$ and a fixed value $s$ of the
remaining sum, the expectation is
\[
 \phi(s)+(u+v)\bigl(\phi(s+1)-\phi(s)\bigr)
       +uv\bigl(\phi(s+2)-2\phi(s+1)+\phi(s)\bigr).
\]
The last coefficient is nonnegative. Replacing $u,v$ by their average
therefore cannot decrease this expression. On the compact set of
parameter vectors with sum $m/d$, choose a maximizer of the convex
expectation, and among these one minimizing the sum of squares of
the parameters. If two parameters differed, averaging them would
give another maximizer with a smaller sum of squares. Thus all
parameters are $m/(db)$, proving the first assertion.

For the second assertion, insert $B-b$ Bernoulli variables with
parameter zero before Step 3. This gives parameter $m/(dB)$, which
is at most $p$. Coupling Bernoulli variables by the same independent
uniform variables shows that increasing all parameters to $p$ cannot
decrease the expectation of an increasing function. This also covers
$b=0$ by starting with $B$ zero parameters.
\end{proof}

For an integer $B\geq1$, $0\leq p\leq1$, and $0<t\leq B$, define
\begin{equation}\label{eq:binomial-positive-part-bound}
 U(B,p,t)=\min\left\{1,
  \min_{\substack{h\in\{0,\ldots,B-1\}\\h<t}}
  \frac{\E(Z-h)_+}{t-h}\right\},
 \qquad Z\sim\operatorname{Bin}(B,p).
\end{equation}
The inner set is nonempty because it contains $h=0$. We set
$U(B,p,t)=1$ when $t\leq0$, and $U(B,p,t)=0$ when $t>B$.
The latter convention uses the support of $Z$ and of any sum being
bounded. For an empty sample, put $U(0,p,t)=0$ for $t>0$ and
$U(0,p,t)=1$ otherwise.

\begin{corollary}[A rational local tail bound]\label{cor:binomial-local-tail}
In Lemma~\ref{lem:class-convex-comparison}, suppose $b\leq B$ and
$\E X\leq M$, where $B\geq1$ is an integer. With
$p=\min(1,M/(dB))$, for every real $t$,
\begin{equation}\label{eq:binomial-local-tail}
 \Pr\{X/d\geq t\}\leq U(B,p,t).
\end{equation}
If $p=u/v$ with integers $0\leq u\leq v$ and $v>0$, then, for
each integer $0\leq h<B$,
\begin{equation}\label{eq:binomial-exact-positive-part}
 \E(Z-h)_+=\frac1{v^B}
  \sum_{k=h+1}^B(k-h)\binom Bk u^k(v-u)^{B-k}.
\end{equation}
Thus \eqref{eq:binomial-positive-part-bound} is an exact rational
calculation when $p,t$ are rational.
\end{corollary}

\begin{proof}
For $h<t$,
\[
 \boldsymbol1_{\{X/d\geq t\}}
       \leq\frac{(X/d-h)_+}{t-h}.
\]
The numerator is an increasing convex function of $X/d$.
Lemma~\ref{lem:class-convex-comparison} bounds its expectation by
$\E(Z-h)_+$. Take the minimum over the indicated integers $h$ and
the trivial bound one. If $t>B$, then $X/d\leq b\leq B$, so the
probability is zero. The binomial mass formula proves
\eqref{eq:binomial-exact-positive-part}, including $u=0$ and $u=v$
with the usual convention for a zeroth power.
\end{proof}

The function in \eqref{eq:binomial-positive-part-bound} is an upper
bound, not an identity for the tail probability. No assertion about
an optimal real choice of $h$ is needed. Any tested value $h<t$
already supplies a valid bound.

\subsubsection{Integer orders and connected completions}

\begin{proposition}[Uniform bounds on an integer interval]
\label{prop:integer-local-tail}
Let $1\leq\ell\leq u$ be integers and let $\theta>0$.
For each integer $n\in[\ell,u]$, let $K_n$ be an integer-valued
sum of the form in Lemma~\ref{lem:class-convex-comparison}, with
entries in $[0,d]$, at most $B\geq1$ selected entries, and
$\E K_n\leq n\mu$, where $\mu\geq0$. The populations and the
probability spaces may depend on $n$. Put
\begin{equation}\label{eq:integer-local-threshold}
 \theta_* =\min_{\ell\leq n\leq u}
                 \frac{\ceil{n\theta}}n,\qquad
 t=\frac{\ceil{\ell\theta_*}}d,\qquad
 p=\min\left(1,\frac{u\mu}{dB}\right).
\end{equation}
Then, for every such $n$,
\[
 \Pr\{K_n/n\geq\theta\}\leq U(B,p,t).
\]
One may replace $\theta_*$ by $\theta$ for a possibly weaker bound.
\end{proposition}

\begin{proof}
On the indicated event, integrality gives
$K_n\geq\ceil{n\theta}$, and hence $K_n/n\geq\theta_*$.
Since $n\geq\ell$, this implies $K_n\geq\ceil{\ell\theta_*}$.
Also $\E K_n\leq u\mu$. Apply
Corollary~\ref{cor:binomial-local-tail} with $M=u\mu$.
Finally, $\theta_*\geq\theta$, so the same argument works with
the latter threshold. This is a bound for each order, not a union
bound over the orders in the interval.
\end{proof}

Some type assignments use a connected completion of the forest.
The completion is useful for controlling the tags, but the local
conflicts need only charge the original edges. The following statement
keeps these two degree lists separate.

\begin{proposition}[Original degrees and a deterministic baseline]
\label{prop:original-local-degrees}
Let $G=(A,B)$ be a bipartite forest, with both parts nonempty, and
let $H$ be a tree on the same vertex set and bipartition containing
$G$. Suppose $d_H(a)\leq D$ for every $a\in A$, where $D\geq1$
is an integer. Partition $A$ into nonempty classes $A_i$ of sizes
$a_i$, and independently choose a uniform $b_i$-subset of each
class, where $0\leq b_i\leq a_i$. Put
\[
 \begin{gathered}
 f_i=\frac{b_i}{a_i},\qquad
 E_i^H=\sum_{a\in A_i}d_H(a),\qquad
 E_i^G=\sum_{a\in A_i}d_G(a),\\
 q_i=E_i^H-E_i^G,\qquad q=\sum_iq_i.
 \end{gathered}
\]
Let $K^H,K^G$ be the sums of the selected degrees in $H,G$,
respectively, using the same selected subsets. Then
\begin{equation}\label{eq:original-local-credit}
 K^G\leq K^H,\qquad
 \E K^G=\E K^H-\sum_iq_if_i
             \leq\E K^H-q\min_i f_i.
\end{equation}
Writing $b=\sum_i b_i$, there is also the identity
\begin{equation}\label{eq:completed-degree-baseline}
 K^H=b+X,\qquad \E X=\E K^H-b,
\end{equation}
where $X$ is sampled from the lists $d_H(a)-1\in[0,D-1]$.
In particular, for a positive integer $n$, if $D>1$ and
\[
 b_-\leq b/n\leq b_+,\qquad \E K^H/n\leq\mu_H,
\]
then
\begin{equation}\label{eq:baseline-local-event}
 \{K^G/n\geq\theta\}\subseteq\{X/n\geq\theta-b_+\},
 \qquad \E X/n\leq\mu_H-b_-.
\end{equation}
The comparison lemmas therefore apply to $X$ with degree bound
$D-1$. If $D=1$, then $X=0$ and $K^G\leq b$ deterministically.

Finally, suppose $G$ is obtained from a tree of order $n$ by first
deleting $k$ vertices, and then deleting an independent set of $w$
vertices, each of degree at least $D$ in the first remaining forest.
If the final forest is nonempty, every completion $H$ as above satisfies
\begin{equation}\label{eq:completion-edge-count}
 q\geq w(D-1).
\end{equation}
\end{proposition}

\begin{proof}
For each vertex, the probability of being selected is $f_i$ in
its class. Summing degrees proves the expectation identity in
\eqref{eq:original-local-credit}; the inequalities follow from
$q_i\geq0$ and $G\subseteq H$. Since $H$ is a nontrivial tree,
each of its degrees is at least one. Subtracting one from each of
the $b$ selected degrees proves \eqref{eq:completed-degree-baseline}.
On the event $K^G/n\geq\theta$, we have
$X/n=K^H/n-b/n\geq\theta-b_+$. Subtracting the lower bound
for $b/n$ from the expectation gives the other assertion in
\eqref{eq:baseline-local-event}.

For the last claim, the forest after the first deletion has at most
$n-k-1$ edges. The next $w$ vertices are independent, so the sets
of edges incident with them are disjoint; at least $wD$ edges are
removed. Consequently
\[
 e(G)\leq n-k-1-wD,\qquad e(H)=n-k-w-1.
\]
Subtracting gives \eqref{eq:completion-edge-count}.
\end{proof}

The mean in \eqref{eq:original-local-credit} may be bounded using
Lemma~\ref{lem:degree-cost-comparison}. The coarse gain
$q\min_i f_i$ is zero whenever any class has frequency zero.
In contrast, the baseline
in \eqref{eq:completed-degree-baseline} uses the completed degrees:
its mean is $\E K^H-b$. One must not subtract the completion credit
from this mean while retaining those same sampled lists.

\subsubsection{Using the excess degrees of a forest}

The previous baseline requires a connected completion. A different
refinement applies directly to the original forest. It uses the
limited total degree above one on a fixed side.

\begin{proposition}[A forest moment bound]\label{prop:forest-local-moment}
Let $G=(A,B)$ be a bipartite forest with both parts nonempty and
$d_G(a)\leq D$ for $a\in A$, where $D>1$ is an integer.
Use independent uniform $b_i$-subsets of nonempty classes $A_i$
as in Proposition~\ref{prop:original-local-degrees}, and let $K$
be the selected degree sum. Put
\[
 \begin{aligned}
 E_i&=\sum_{a\in A_i}d_G(a),&
 Q_i&=\sum_{a\in A_i}(d_G(a)-1)_+,\\
 f_i&=b_i/\abs{A_i},&
 \mu&=\sum_i f_iE_i=\E K,\qquad q=\sum_i f_iQ_i.
 \end{aligned}
\]
If $b=\sum_i b_i>0$, then for every $\lambda>0$,
\begin{equation}\label{eq:forest-local-moment}
 \begin{gathered}
 \log\E e^{\lambda K}
 \leq b\log\left(1+
       \frac{(e^\lambda-1)\mu+c_\lambda q}{b}\right),\\
 c_\lambda=\frac{e^{\lambda D}-D e^\lambda+D-1}{D-1}\geq0.
 \end{gathered}
\end{equation}
Moreover,
\begin{equation}\label{eq:local-excess-mean}
 q\leq\min\left\{\frac{D-1}{D}\mu,
                   (\abs{B}-1)\max_i f_i\right\}.
\end{equation}
The right side of \eqref{eq:forest-local-moment} remains an upper
bound if $b$ is replaced by any real $B_0\geq b$, and $\mu,q$
by upper bounds. If $b=0$, then $K=0$ surely.
\end{proposition}

\begin{figure}[ht]
\centering
\begin{tikzpicture}[x=12mm,y=5mm,every node/.style={font=\normalsize}]
\draw[->] (-.2,0)--(3.5,0) node[right] {$r$};
\draw[->] (0,-.2)--(0,9) node[above] {$e^{\lambda r}$};
\draw[gray,densely dotted,thick] (0,1)--(3,8);
\draw[thick] (0,1)--(1,2)--(3,8);
\draw[dashed,domain=0:3,samples=49,smooth] plot (\x,{2^\x});
\foreach \x/\y in {0/1,1/2,3/8}
  \fill (\x,\y) circle (1.6pt);
\node[anchor=east] at (-.05,1) {$1$};
\node[anchor=east] at (-.05,8) {$8$};
\node[anchor=north] at (0,-.15) {$0$};
\node[anchor=north] at (1,-.15) {$1$};
\node[anchor=north] at (3,-.15) {$D=3$};
\draw[gray,densely dotted,thick] (4.05,7.4)--(4.75,7.4);
\node[anchor=west] at (4.9,7.4) {single chord};
\draw[thick] (4.05,5.8)--(4.75,5.8);
\node[anchor=west] at (4.9,5.8) {two chords};
\draw[dashed] (4.05,4.2)--(4.75,4.2);
\node[anchor=west] at (4.9,4.2) {$e^{\lambda r}=2^r$};
\end{tikzpicture}
\caption{The bound \eqref{eq:two-degree-chords} for $D=3$ and
$\lambda=\log2$. The increase in slope is charged
only to degrees above one.}
\label{fig:degree-chords}
\end{figure}
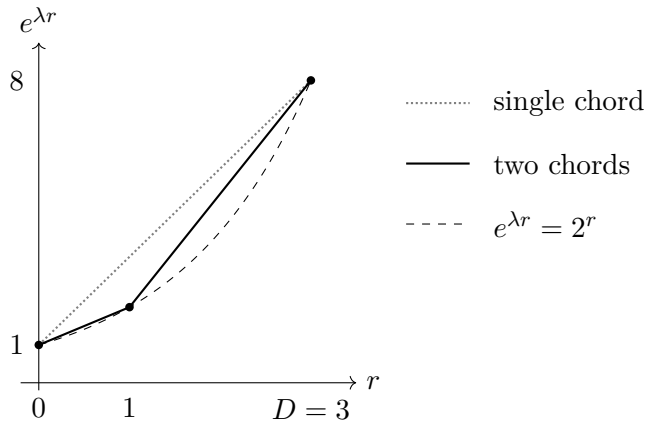

\begin{proof}
Convexity on $[0,1]$ and $[1,D]$ gives, for $0\leq r\leq D$,
\begin{equation}\label{eq:two-degree-chords}
 e^{\lambda r}\leq
  1+(e^\lambda-1)r+c_\lambda(r-1)_+.
\end{equation}
The expression on the right is the chord on each interval. The
slopes of successive chords of a convex function are increasing,
so $c_\lambda\geq0$; see Figure~\ref{fig:degree-chords}.

Step 1 of Lemma~\ref{lem:class-convex-comparison}, applied to
the exponential function in each class, gives
\[
 \E e^{\lambda K}\leq
 \prod_{i:b_i>0}\left(
    \frac1{\abs{A_i}}\sum_{a\in A_i}e^{\lambda d_G(a)}
                  \right)^{b_i}.
\]
Insert \eqref{eq:two-degree-chords}, take logarithms, and use
concavity of the logarithm with weights $b_i/b$. The resulting
bound is \eqref{eq:forest-local-moment}.

For $0\leq r\leq D$,
$(r-1)_+\leq (D-1)r/D$, so the first bound in
\eqref{eq:local-excess-mean} follows class by class.
Lemma~\ref{lem:weighted-forest-cost}, with all position weights
equal to one, gives
$\sum_i Q_i\leq\abs{B}-1$. Multiplication by $\max_i f_i$
gives the second bound.

Finally, for $x\geq0$ the function $v\log(1+x/v)$ is
increasing for $v>0$: its derivative is
\[
 \log(1+x/v)-\frac{x/v}{1+x/v}\geq0,
\]
because $\log(1+t)=\int_0^t(1+s)^{-1}\,ds\geq t/(1+t)$.
Both coefficients of $\mu,q$ are nonnegative. These observations
justify all stated replacements.
\end{proof}

\begin{corollary}[A uniform forest tail]\label{cor:uniform-forest-tail}
Suppose a family of sums $K_n$ as in
Proposition~\ref{prop:forest-local-moment}, with a common degree
bound $D$ and integer orders $1\leq\ell\leq n\leq u$, satisfies
\[
 b/n\leq\beta,\qquad \E K_n/n\leq\bar\mu,
 \qquad q/n\leq\bar q,
\]
where $\beta>0$ and $\bar\mu,\bar q\geq0$. For $\theta>0$
and $\lambda>0$, put
\begin{equation}\label{eq:uniform-forest-rate}
 a=\lambda\theta-\beta\log\left(1+
       \frac{(e^\lambda-1)\bar\mu+c_\lambda\bar q}{\beta}\right).
\end{equation}
For every such $n$,
\[
 \Pr\{K_n/n\geq\theta\}\leq
 \begin{cases}e^{-\ell a},&a>0,\\1,&a\leq0.\end{cases}
\]
Since $K_n$ is integer-valued, $\theta$ may first be replaced by
$\theta_*$ from \eqref{eq:integer-local-threshold}.
\end{corollary}

\begin{proof}
When $b>0$, Proposition~\ref{prop:forest-local-moment}, with
$B_0=n\beta$, bounds $\log\E e^{\lambda K_n}$ by
$n(\lambda\theta-a)$. Markov's inequality gives $e^{-na}$.
For $a>0$, this is at most $e^{-\ell a}$; otherwise use the
bound one. An empty sample has zero probability of exceeding the
positive threshold. The integer-threshold replacement follows by
the event inclusion proved in Proposition~\ref{prop:integer-local-tail}.
\end{proof}

There is no need to solve for an optimal $\lambda$: each fixed
positive choice gives a valid upper bound. Numerical evaluation
must bound the logarithmic moment upwards and the resulting
positive rate downwards. A numerical search for a parameter and
a verification of its inequality are separate operations.

\subsubsection{Which local events must be checked}

Different odd labels can give the same conflict event. More generally,
one event may contain another. We use this inclusion before taking
the union bound.

\begin{lemma}[Maximal prime supports]\label{lem:maximal-local-supports}
Fix the class permutations and the original degrees $d_a$.
For a block type $J$ and a finite set $P$ of primes greater than $z$,
put
\begin{equation}\label{eq:support-local-variable}
 K(J,P)=\sum_{a\in A}d_a
   \boldsymbol1_{I(a)\cap J=\varnothing}
   \boldsymbol1_{\{\text{some }p\in P\text{ divides }j(a)\}}.
\end{equation}
If $P\subseteq P'$, then $K(J,P)\leq K(J,P')$ for every outcome.
For a finite set $Y$ of odd labels, put
\[
 P(y)=\{p>z:p\mid y\}.
\]
Let $\mathcal M_J$ consist of the inclusion-maximal members of
$\{P(y):y\in Y,\ T(y)=J\}$.
For any real thresholds $L_J$,
\begin{equation}\label{eq:maximal-local-event-equality}
 \bigcup_{y\in Y}\{K(T(y),P(y))\geq L_{T(y)}\}
 =\bigcup_J\ \bigcup_{P\in\mathcal M_J}\{K(J,P)\geq L_J\}.
\end{equation}
One may first remove any primes that divide no available half-label.
In particular, if every half-label is at most $\floor{u/2}$,
all primes greater than $u/2$ may be removed.
\end{lemma}

\begin{proof}
The divisibility indicator in \eqref{eq:support-local-variable}
is increasing under inclusion of $P$, proving the first claim.
Every support in a finite family is contained in a maximal member
of that family, with the same $J$. This proves inclusion of the
left union in the right union. Each maximal support is itself the
support of a label in $Y$, proving the reverse inclusion.
Removing a prime dividing no available half-label leaves every
indicator unchanged. A positive integer at most $\floor{u/2}$
cannot have a prime divisor greater than $u/2$.
\end{proof}

\begin{corollary}[Simultaneous local estimates]
\label{cor:refined-simultaneous-local}
In the matching setup of Section~\ref{sec:matching-estimates},
fix positive thresholds $L_J$ for the label types being estimated
probabilistically. For each maximal support event in
Lemma~\ref{lem:maximal-local-supports}, let $\varepsilon_{J,P}$
be an upper bound for
$\Pr\{K(J,P)\geq L_J\}$. Then, with probability at least
\begin{equation}\label{eq:refined-local-success}
 1-\sum_J\sum_{P\in\mathcal M_J}\varepsilon_{J,P},
\end{equation}
every such label $y$ loses fewer than $L_{T(y)}$ incident edges.
The bounds in this section may be combined by taking their minimum
for the same event. If a further event has failure probability at
most $\varepsilon_0$, both conclusions hold simultaneously with
positive probability whenever
\[
 \varepsilon_0+
       \sum_J\sum_{P\in\mathcal M_J}\varepsilon_{J,P}<1.
\]
\end{corollary}

\begin{proof}
For fixed $J,P$, let $b_I$ be the number of half-labels in
$\Lambda_I$ divisible by a prime in $P$, setting $b_I=0$ when
$I\cap J\ne\varnothing$. Their assigned positions form independent
uniform $b_I$-subsets of the classes. Thus $K(J,P)$ has precisely
the sampling form used in this section. The local deletion count
for $y$ is at most $K(T(y),P(y))$, as proved in
Section~\ref{sec:matching-estimates}. Use
\eqref{eq:maximal-local-event-equality} and the union bound.
The minimum of several valid upper bounds is still an upper bound.
A final union bound with the additional failure event proves the
last assertion. No independence between these events is required.
\end{proof}

The thresholds must be common within each family to which
Lemma~\ref{lem:maximal-local-supports} is applied. With several
degree targets, apply it separately within each target class, or
use a common smaller threshold. Labels protected deterministically
by Lemma~\ref{lem:protected-full-labels} need not occur in these
families. Distinct maximal supports with the same numerical bound
still contribute separately to the sum; only actual event inclusion
or equality justifies deleting an event. Likewise, this reduction
does not change the multiplicities of labels in a global edge count.

Corollary~\ref{cor:refined-simultaneous-local} supplies the local
probability term in Corollary~\ref{cor:three-range-probability}.
For a numerical order interval, the remaining task is to construct
the type assignment, establish all mean and population bounds
uniformly, and verify the matching budgets. The present section
proves the probability estimates conditional on those quantities;
it does not assert a new interval of prime-labelable orders.

\subsection{Assignments controlled by decomposition height}
\label{sec:height-assignments}

The connected construction of Section~\ref{sec:forest-types} separates
mixed components until each contains one $A$-vertex. Its cost depends
on how quickly these components can be split by deleting vertices of
$B$. We now count the cuts by the height of such a decomposition.
Integer bounds for this height retain information that is lost in a
single geometric contraction factor.

We first sharpen the rounding and tag estimates used at the end of
the construction. Neither improvement changes the two preserved
measures. The resulting assignment theorem will keep the separator
positions and give all of them the full tag.

\subsubsection{Rounding without a prescribed sign}

\begin{lemma}[Averaging the slot matchings]\label{lem:averaged-slot-rounding}
Under the hypotheses of Lemma~\ref{lem:sorted-slot-rounding}, let
$q$ be the number of positive quotas and let $R=b-a$ be the span
of the item weights. There is an integral assignment with exact quotas
such that
\begin{equation}\label{eq:averaged-slot-error}
 \abs{\epsilon_i}\leq R,\qquad
 \frac12\sum_i\abs{\epsilon_i}\leq\frac{qR}{4}.
\end{equation}
No sign is prescribed for an individual error. The assignment may
be chosen from a finite convex decomposition of the fractional
item--slot matrix into perfect matchings.
\end{lemma}

\begin{proof}
Use the sorted slots and their convex decomposition from the proof
of Lemma~\ref{lem:sorted-slot-rounding}. Choose a matching with
the probabilities in that decomposition. For color $i$, write
$X_i$ for its matched degree sum and $T_i$ for its fractional target.
The decomposition gives $\E X_i=T_i$. Let $L_i,U_i$ be the sums
of the smallest and largest weights in its slots. Every matching
satisfies
\[
 L_i\leq X_i\leq U_i,\qquad U_i-L_i\leq R.
\]
If $U_i>L_i$, the chord bound for the convex function
$(x-T_i)_+$ gives
\[
 \E(X_i-T_i)_+
 \leq\frac{(U_i-T_i)(T_i-L_i)}{U_i-L_i}
 \leq\frac{U_i-L_i}{4}.
\]
If $U_i=L_i$, the error is zero. Since the sum of the errors is
zero for every matching, their half-$\ell^1$ norm equals the sum
of their positive parts. Its expectation is therefore at most
$qR/4$. At least one matching attains a value no larger than this
expectation. The individual error bound holds for every matching,
by Lemma~\ref{lem:sorted-slot-rounding}.
\end{proof}

In a completing tree the relevant degrees lie in $[1,D]$, so the
span in this lemma is $D-1$. The concentration bound on a single
degree remains $D$. A construction requiring a nonpositive error
in a distinguished class must instead use the preferred-sign
rounding lemma, or pay explicitly for a possible positive error.

\begin{lemma}[Choosing the last tag]\label{lem:optimal-last-tag}
Assume the hypotheses of Lemma~\ref{lem:integral-tag-recovery}
and $H_{\mathrm{deg}}\leq K$. For $g\geq2$, put
\begin{equation}\label{eq:optimal-tag-bound}
 \mathcal R_g(K)=\max\left\{K,\ 1-\frac1g,
                   \frac{g-2}{g-1}(K+1)\right\},
\end{equation}
and put $\mathcal R_1(K)=0$. There are integral tags respecting
every forced tag outside $C$ whose count vector $b$ satisfies
\begin{equation}\label{eq:optimal-tag-variation}
 \operatorname{tv}(b,\abs{B}\pi)
 \leq\min\left\{H_{\mathrm{deg}}+1-\max_i\pi_i,
                          \mathcal R_g(K)\right\}.
\end{equation}
\end{lemma}

\begin{proof}
If $C=\varnothing$, all $A$-vertices have the same color and
all tags can have that color. The variation is zero. Suppose
$C\ne\varnothing$. In the notation of
Lemma~\ref{lem:integral-tag-recovery}, every forced count is at
most $x_i=E_i-a_i$, and $\sum_i x_i=\abs{B}-1$. Thus each
choice $b=x+e_j$ is feasible, where $e_j$ is a coordinate unit
vector. Choose the one minimizing variation. The earlier choice
of a largest quota proves the first bound in
\eqref{eq:optimal-tag-variation}.

Write $q_i=\pi_i-h_i$, so $\sum_iq_i=1$, and put
$t=\sum_i(-q_i)_+$. Since $\pi_i\geq0$, we have
$t\leq\sum_i(h_i)_+=H_{\mathrm{deg}}\leq K$. The error
vector for the choice $j$ is $e_j-q$. Directly summing its
absolute coordinates gives
\begin{equation}\label{eq:last-tag-exact-minimum}
 \min_j\frac12\sum_i\abs{(e_j)_i-q_i}
       =t+1-\min\{1,\max_i q_i\}.
\end{equation}
For completeness, if $0\leq q_j\leq1$ the value is $t+1-q_j$;
if $q_j\geq1$ it is $t$; and if $q_j<0$ it is $t+1$.
Some coordinate is positive because their sum is one, which proves
\eqref{eq:last-tag-exact-minimum}.

If $\max_i q_i\geq1$, the minimum is $t\leq K$. If all
$q_i$ are nonnegative, it is at most $1-1/g$. In the remaining
case there are at most $g-1$ positive coordinates, with sum
$1+t$. Consequently $\max_iq_i\geq(1+t)/(g-1)$, and
\eqref{eq:last-tag-exact-minimum} is at most
$(g-2)(1+t)/(g-1)\leq(g-2)(1+K)/(g-1)$.
These three cases prove the claim for $g\geq2$. For one color
all degree and tag errors vanish.
\end{proof}

For example, with four colors and $K=D-1\geq2$, the bound is
$\mathcal R_4(K)=K$. There is then no additional constant beyond
the degree variation before changing the separator tags.

\subsubsection{Integer bounds for the decomposition height}

Let $U=(A_U,B_U)$ be a bipartite tree with $m=\abs{A_U}\geq1$.
Set $L(U)=m-1$. A component with $L(U)=0$ is terminal. In a
nonterminal component we may delete a vertex $b\in B_U$ having
at least two neighbors in that component. Such a vertex exists
on the path between two $A_U$-vertices. Each component of $U-b$
has a positive, strictly smaller number of $A_U$-vertices.

Define its minimum decomposition height by the finite recursion
\begin{equation}\label{eq:minimum-side-height}
 h(U)=
 \begin{cases}
  0,&L(U)=0,\\
  1+\displaystyle\min_b\ \max_{V\text{ a component of }U-b}h(V),
       &L(U)>0,
 \end{cases}
\end{equation}
where the minimum is over the admissible $B_U$-vertices just
specified. A minimizing cut reduces the height on every child
by at least one. Computing this minimum need not be efficient;
it is a finite construction for every input tree.

Fix an integer $D\geq3$. For a nonnegative integer $L$, set
\[
 a_D(0)=0,\qquad
 a_D(L)=\max\left\{0,\ L-\ceil{L/D},
                         \floor{(L-1)/2}\right\}\quad(L>0).
\]
For $L>0$ and $1\leq i\leq D$, put
\begin{equation}\label{eq:branch-height-remainders}
 \begin{split}
 b_i(L)&=\max\{0,\min(\floor{(L+1)/2},\floor{L/i})-1\},\\
 c_i(L)&=\max\{0,L-\ceil{iL/D}\},\\
 B_j(L)&=\max\left\{c_j(L),
               \max_{1\leq i\leq j}a_D^{\,j-i}(b_i(L))\right\}
                  \quad(1\leq j\leq D),
 \end{split}
\end{equation}
where the superscript denotes iteration and $a_D^0(x)=x$.
The function $a_D$ is nondecreasing and $a_D(L)<L$ for $L>0$.
All quantities $b_i(L),c_i(L),B_j(L)$ also lie in $[0,L-1]$.

\begin{lemma}[Peeling at a weighted centroid]\label{lem:integer-centroid-peeling}
Let $U$ be a bipartite tree with $L=L(U)>0$ and
$\Delta_U(A_U)\leq D$. A single admissible cut can make every
child's $L$-value at most $a_D(L)$. For each $1\leq j\leq D$,
there is a decomposition through at most $j$ levels, stopping
terminal components, for which every remaining $L$-value is at
most $B_j(L)$.

More precisely, if a weighted centroid lies in $A_U$, order the
positive $A_U$-counts in its branches as
$s_1\geq\cdots\geq s_d>0$, where $d\leq D$, and put
$S_i=\sum_{r\leq i}s_r$. Deleting its $B_U$-neighbors in that
order leaves a central component with $L$-value $L-S_i$ after
the $i$th cut. Every component detached at that cut has
$L$-value at most $s_i-1$.
\end{lemma}

\begin{proof}
Use the weighted centroid from Lemma~\ref{lem:weighted-side-cut},
with weight one on $A_U$ and zero on $B_U$. Every branch has
at most $\floor{(L+1)/2}$ vertices of $A_U$.
If the centroid belongs to $B_U$, deleting it leaves every child
with $L$-value at most $\floor{(L-1)/2}$. It is admissible,
since otherwise its unique branch would contain all $L+1$
vertices of $A_U$.

If the centroid belongs to $A_U$, its branch counts sum to $L$.
The largest is at least $\ceil{L/D}$. Delete its neighbor in
that branch. The central component has $L$-value at most
$L-\ceil{L/D}$, while every other child has $L$-value at most
$\floor{(L-1)/2}$. This proves the single-cut assertion.

For repeated peeling, ignore branches with zero $A_U$-count.
The neighbor in every positive branch has a further $A_U$-neighbor,
so its deletion is admissible. Distinct branches meet only at the
centroid. Thus cutting branch $i$ removes exactly $s_i$ vertices
of $A_U$ from the central component, and each detached component
has at most $s_i$ such vertices. This proves the last assertion.

Pad the branch counts by zeros to a list of length $D$.
The average of its first $i$ entries is at least the average of
all entries, giving $S_i\geq\ceil{iL/D}$. Also
$s_i\leq\floor{L/i}$ and $s_i\leq\floor{(L+1)/2}$.
Consequently the central remainder is at most $c_i(L)$, and
each component detached at step $i$ has remainder at most $b_i(L)$.
Apply ordinary single cuts to the detached components at subsequent
levels. Monotonicity of $a_D$ bounds their remainders at level
$j$ by $a_D^{\,j-i}(b_i(L))$. If fewer than $j$ positive
branches exist, the central component becomes terminal earlier.
This proves the bound $B_j(L)$.

When the fresh centroid is in $B_U$, all children are detached
at the first level and have remainder at most $b_1(L)$.
Their further ordinary cuts are covered by the $i=1$ term in
$B_j(L)$. Terminal components may always be left unchanged.
\end{proof}

Instead of making all detached components wait until level $j$,
we may decompose each of them as soon as it appears. The next
recursion takes the better of these two constructions.

Fix $1\leq k\leq D$. Starting with $F(0)=0$, define successively
for $L\geq1$
\begin{align}
 G_j(L)&=\max\left\{j+F(c_j(L)),
                   \max_{1\leq i\leq j}\bigl(i+F(b_i(L))\bigr)\right\},
       \label{eq:immediate-height-candidate}\\
 F(L)&=\max\left\{F(L-1),
       \min_{1\leq j\leq k}\min\{j+F(B_j(L)),G_j(L)\}\right\}.
       \label{eq:fixed-height-envelope}
\end{align}
Only previously defined values occur on the right.

\begin{proposition}[A height bound with a fixed set of stopping levels]
\label{prop:fixed-height-envelope}
The sequence in \eqref{eq:fixed-height-envelope} is nondecreasing,
satisfies $0\leq F(L)\leq L$, and bounds
$h(U)\leq F(L(U))$ for every bipartite tree whose $A$-degrees
are at most $D$.
\end{proposition}

\begin{proof}
Induct on $L$. Nondecreasingness is built into the recursion.
For a given $j$, apply the first $j$ levels of
Lemma~\ref{lem:integer-centroid-peeling}. Every remaining component
has smaller remainder, at most $B_j(L)$. By induction and
monotonicity its height is at most $F(B_j(L))$. This gives
the candidate $j+F(B_j(L))$.

Suppose next that the weighted centroid is in $A$. Peel $j$
branches, or stop earlier if they are exhausted, and immediately
decompose each detached component. A component detached at step
$i$ has depth at most $i+F(b_i(L))$. The surviving central
component has depth at most $j+F(c_j(L))$. Thus $G_j(L)$
is another valid bound. If the centroid is in $B$, its first
cut gives the bound $1+F(b_1(L))$, which is a term of every
$G_j(L)$. The same conclusion follows.

Every candidate bounds a full decomposition. Their minimum does
also, and increasing it to at least $F(L-1)$ preserves that
property. Finally the candidate with $j=1$ is at most
$1+F(a_D(L))\leq L$ by induction. Together with
$F(L-1)\leq L-1$, this proves $F(L)\leq L$.
\end{proof}

\subsubsection{Choosing the stopping level from the branch counts}

The preceding bound chooses among predetermined levels. We can do
better by allowing that choice to depend on the actual branch counts,
while still verifying the height bound by a short integer test.

Suppose a nondecreasing sequence $H$ has been defined through $L-1$,
with $H(0)=0$. For every integer $r$, set
\begin{equation}\label{eq:height-inverse}
 T_L(r)=\max\bigl(\{x\in\{0,\ldots,L-1\}:H(x)\leq r\}
                         \cup\{-1\}\bigr).
\end{equation}
To test a proposed integer height $h\geq H(L-1)$, first require
\begin{equation}\label{eq:adaptive-height-b-centroid}
 1+H(\floor{(L-1)/2})\leq h.
\end{equation}
Put $u=\floor{(L+1)/2}$. For each
$1\leq t\leq\min(D,h+1)$, define
\[
 a_t=T_L(h-t)+2,\qquad
 C_j=L-T_L(h-j)-1\quad(1\leq j<t).
\]
Require at least one of the following exclusions for each such $t$:
\begin{align}
 &a_t>u,\label{eq:adaptive-exclusion-one}\\
 &t a_t>L,\label{eq:adaptive-exclusion-two}\\
 &C_j<j a_t\quad\text{for some }1\leq j<t,
                       \label{eq:adaptive-exclusion-three}\\
 &\min\left\{Du,
       \min_{1\leq j<t}\bigl(C_j+(D-j)\floor{C_j/j}\bigr)\right\}<L.
                       \label{eq:adaptive-exclusion-four}
\end{align}
For $t=1$, the minimum in \eqref{eq:adaptive-exclusion-four}
has only its $Du$ term. Define $H(L)$ to be the least integer
$h\geq H(L-1)$ passing this test, beginning with $H(0)=0$.

\begin{proposition}[A height bound with an adaptive stopping level]
\label{prop:adaptive-height-envelope}
The preceding test defines a nondecreasing integer sequence with
$0\leq H(L)\leq L$. Every bipartite tree with $A$-degrees at
most $D$ satisfies $h(U)\leq H(L(U))$.
\end{proposition}

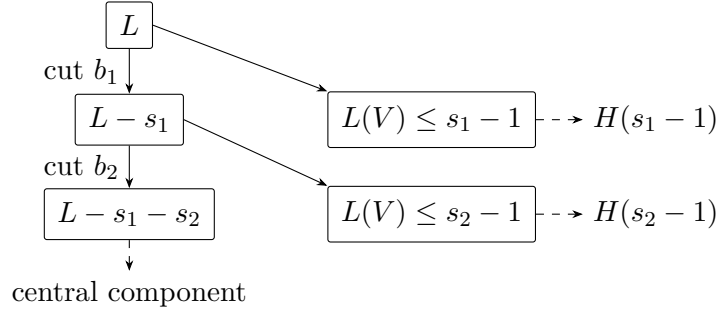
\begin{figure}[ht]
\centering
\begin{tikzpicture}[x=1cm,y=1cm,every node/.style={font=\normalsize},
 piece/.style={draw,rounded corners=1pt,inner sep=5pt,align=center},
 arrow/.style={-{Stealth[length=4pt]},thin}]
\node[piece] (u0) at (0,0) {$L$};
\node[piece] (u1) at (0,-1.25) {$L-s_1$};
\node[piece] (u2) at (0,-2.5) {$L-s_1-s_2$};
\node[piece] (v1) at (4,-1.25) {$L(V)\leq s_1-1$};
\node[piece] (v2) at (4,-2.5) {$L(V)\leq s_2-1$};
\draw[arrow] (u0)--node[left] {cut $b_1$}(u1);
\draw[arrow] (u0.east)--(v1.north west);
\draw[arrow] (u1)--node[left] {cut $b_2$}(u2);
\draw[arrow] (u1.east)--(v2.north west);
\draw[arrow,dashed] (u2)--(0,-3.25);
\draw[arrow,dashed] (v1.east)--(6,-1.25)
  node[anchor=west] {$H(s_1-1)$};
\draw[arrow,dashed] (v2.east)--(6,-2.5)
  node[anchor=west] {$H(s_2-1)$};
\node[anchor=north] at (0,-3.25) {central component};
\end{tikzpicture}
\caption{Peeling around a fixed $A$-centroid. A side box represents
all components detached at that step; they are decomposed immediately.
After $j$ central cuts, the depth is bounded by
\eqref{eq:actual-branch-stopping-bound}.}
\label{fig:height-peeling}
\end{figure}

\begin{proof}
Assume the assertions hold for all smaller remainders. If the
weighted centroid is in $B$, one cut and
\eqref{eq:adaptive-height-b-centroid} give height at most $h$.
Otherwise let $s_1\geq\cdots\geq s_d>0$ be its positive
branch counts, where $d\leq D$, $\sum_i s_i=L$, and $s_i\leq u$.
By Lemma~\ref{lem:integer-centroid-peeling}, stopping after
$j\leq d$ central cuts and decomposing each detached component
immediately gives height at most
\begin{equation}\label{eq:actual-branch-stopping-bound}
 \max\left\{j+H(L-S_j),
                \max_{1\leq i\leq j}\bigl(i+H(s_i-1)\bigr)\right\}.
\end{equation}
Here all arguments of $H$ are smaller than $L$. Figure~\ref{fig:height-peeling}
shows the separate contributions to this maximum.

Suppose no $j\in[d]$ makes \eqref{eq:actual-branch-stopping-bound}
at most $h$. There is then a first index $t$ with
$t+H(s_t-1)>h$. Indeed, otherwise the choice $j=d$ works:
the central remainder is zero and the last detached-branch inequality
gives $d\leq h$. Moreover $t\leq h+1$, since every positive
branch with index $h+1$ would violate its inequality.

By monotonicity of $H$, the failure at $t$ implies
$s_t\geq T_L(h-t)+2=a_t$. Hence $a_t\leq u$ and
$ta_t\leq L$, so neither \eqref{eq:adaptive-exclusion-one}
nor \eqref{eq:adaptive-exclusion-two} holds. For $j<t$, all
detached components fit within height $h$. The assumed failure
at that stopping point therefore gives
$j+H(L-S_j)>h$, equivalently $S_j\leq C_j$.
But $S_j\geq j a_t$, excluding \eqref{eq:adaptive-exclusion-three}.
Every later branch count is at most
$\floor{S_j/j}\leq\floor{C_j/j}$. Padding by zeros to $D$
terms yields
\[
 L\leq C_j+(D-j)\floor{C_j/j}\quad(1\leq j<t).
\]
Also $L\leq Du$. Thus \eqref{eq:adaptive-exclusion-four} fails
as well, contradicting the test. Some stopping point works, so
the proposed height bounds $h(U)$.

It remains to prove that a height is always accepted. Set $h=L$.
By induction $H(x)\leq x$, so
\eqref{eq:adaptive-height-b-centroid} holds. For $1\leq t\leq L$,
we have $T_L(L-t)\geq L-t$ and hence $a_t\geq L-t+2$.
Thus
\[
 t a_t-L\geq t(L-t+2)-L=(t-1)(L-t)+t>0.
\]
If $t=L+1$ is tested, then $a_t=1$ and again $ta_t>L$.
These are all possible tested indices, so every one satisfies
\eqref{eq:adaptive-exclusion-two}. The choice $h=L$ passes.
It follows that $H(L)\leq L$; its imposed lower bound makes
the sequence nondecreasing. This completes the induction.
\end{proof}

\subsubsection{Counting cuts and assigning types}

Let $R$ be either of the nondecreasing height bounds just proved.
For $M\geq1$ and $d\geq0$, define
\begin{equation}\label{eq:height-separator-budget}
 \mathcal C_R(M,d)=
 \sum_{j=1}^{R(M-1)}\min\left\{d,\floor{M/s_j}\right\},
 \qquad s_j=1+\min\{L\geq1:R(L)\geq j\}.
\end{equation}
The sum is empty when $M=1$. Every $s_j$ occurring in it is
defined using an $L\leq M-1$, so finite initial values of $R$
suffice to compute the budget.

\begin{theorem}[Type assignment with a height budget]
\label{thm:height-type-assignment}
Let $Q=(A,B)$ be a tree with both parts nonempty, $M=\abs{A}$,
$e=e(Q)$, and $\Delta_Q(A)\leq D$, where $D\geq3$ is an
integer. Prescribe nonnegative integer quotas $a_1,\ldots,a_g$
of total $M$, and put $\pi_i=a_i/M$. There are a coloring of
$A$ with these exact quotas and a set $C\subseteq B$ such that
every vertex of $B\setminus C$ has neighbors of a single color,
\begin{equation}\label{eq:height-assignment-errors}
 \abs{C}\leq\mathcal C_R(M,2(g-1)),\qquad
 \abs{h_i}\leq D-1,\qquad
 H_{\mathrm{deg}}\leq\frac{g(D-1)}4,
\end{equation}
where $h_i=\sum_{a:\,c(a)=i}d_Q(a)-e\pi_i$ and
$H_{\mathrm{deg}}=\frac12\sum_i\abs{h_i}$.
No preferred error sign is required or asserted.
\end{theorem}

\begin{proof}
Use the fractional construction of
Proposition~\ref{prop:two-measure-assignment}, preserving
$A$-count and the sum of degrees in the fixed tree $Q$.
After each pruning, at most $d=2(g-1)$ mixed objects remain.
Among the nonterminal mixed objects choose one of largest actual
height $h(U)$, and cut it at a vertex attaining the minimum in
\eqref{eq:minimum-side-height}. The cut belongs to $B$ and
contributes zero to both measures. Give the children their parent's
fractional weights, prune again, and freeze pure objects as before.

Consider the first cut at a particular positive height $j$.
No mixed object of larger height remains. At that moment there
are at most $d$ mixed objects of height $j$, and every later cut
produces children of strictly smaller height than its parent.
Thus at most $d$ objects of height $j$ are ever cut. They are
pairwise disjoint: no descendant of a cut object can have the
same height. Since $j=h(U)\leq R(L(U))$, each contains at
least $s_j$ vertices of $A$. There are at most
$\floor{M/s_j}$ such objects. Summing over heights proves
the separator bound.

Each cut removes a new $B$-vertex, so the construction terminates.
A terminal mixed component contains just one $A$-vertex. Its degree
in the fixed tree lies in $[1,D]$, even if all its current neighbors
have already been cut. Apply Lemma~\ref{lem:averaged-slot-rounding}
to these last fractional assignments. The residual quotas are
integers, since pure components contribute integer counts. There
are at most $g$ positive residual quotas, giving the degree bounds
in \eqref{eq:height-assignment-errors}. Finally, every edge outside
$C$ lies in one object, and rounding a one-$A$-vertex object leaves
it monochromatic on $A$. This proves the neighborhood assertion.
\end{proof}

\begin{corollary}[Retaining all positions with full separator tags]
\label{cor:height-retained-assignment}
Let $G=(A,B)$ be a bipartite forest with $M=\abs{A}>0$ and
$N=\abs{B}>0$. Suppose $D\geq3$ is an integer,
$\Delta_G(A)\leq D$, and $DM\geq M+N-1$.
For every type $I\subseteq[s]$, prescribe an integer quota $a_I\geq0$,
with $\sum_Ia_I=M$, and put $g=2^s$ and $\pi_I=a_I/M$.
There are a completing tree $Q$, an exact type assignment on $A$,
a set $C\subseteq B$, and tags $J(b)$ on all of $B$ such that
\begin{enumerate}
 \item $\Delta_Q(A)\leq D$, and the degree errors measured in $Q$
       satisfy \eqref{eq:height-assignment-errors};
 \item for $b\notin C$, every neighbor in $Q$ has type $J(b)$;
       for $b\in C$, one has $J(b)=[s]$;
 \item $\abs{C}\leq\mathcal C_R(M,2(g-1))$, and the tag count
       vector $\beta$ satisfies
 \begin{equation}\label{eq:height-retained-tag-budget}
  \operatorname{tv}(\beta,N\pi)
   \leq\mathcal R_g\left(\frac{g(D-1)}4\right)
                          +\mathcal C_R(M,2(g-1)).
 \end{equation}
\end{enumerate}
In particular $I(a)\subseteq J(b)$ on every edge of $G$.
No type is moved between $A$-vertices after the quota rounding.
\end{corollary}

\begin{proof}
Complete $G$ by Lemma~\ref{lem:bipartite-completion} and apply
Theorem~\ref{thm:height-type-assignment} to the resulting tree $Q$.
Recover integral tags by Lemma~\ref{lem:optimal-last-tag}, with
$K=g(D-1)/4$. Outside $C$ these tags already agree with every
neighboring type. Change the tags on $C$ to $[s]$. The variation
cost is at most $\abs{C}$, as in Lemma~\ref{lem:retained-tags},
which proves \eqref{eq:height-retained-tag-budget}. Both tag rules
give inclusion on every edge of $Q$, and therefore of $G$.
Tag recovery does not alter the rounded $A$-assignment.
\end{proof}

The full tag at each separating position has a useful additional
consequence: a label of nonempty type is incompatible with that
position. Any allowed pair involving such a label therefore uses
a position outside $C$, where all adjacent even half-labels have
exactly the position's tag. This is the hypothesis used in the
support-sensitive forest estimate of Section~\ref{sec:refined-hall}.
The local bounds of Section~\ref{sec:local-probability} may still
charge original degrees, with the completion credit proved there.
The degree errors in \eqref{eq:height-assignment-errors} themselves
refer to $Q$, whose edge count is $M+N-1$.

\begin{corollary}[A bound uniform over an order interval]
\label{cor:height-interval-budget}
Let $1\leq\ell\leq u$ be integers. Fix $D,g$ and one of the
height sequences above. Suppose the hypotheses of
Corollary~\ref{cor:height-retained-assignment} hold separately at
each order, with even-side sizes $M(n)\leq M_*$. Then one may use
the common normalized separator bound
\[
 \frac{\abs{C}}n\leq
          \frac{\mathcal C_R(M_*,2(g-1))}{\ell}.
\]
In the preparation where $M(n)=\floor{n/2}-w(n)$, a valid choice is
\begin{equation}\label{eq:height-interval-even-cap}
 M_* =\floor{u/2}-\min_{\ell\leq n\leq u}w(n),
\end{equation}
provided every actual $M(n)$ is positive. Here $w(n)$ is the
special-even count in \eqref{eq:even-supply}.
\end{corollary}

\begin{proof}
For fixed $R,d$, the budget $\mathcal C_R(M,d)$ is nondecreasing
in $M$: its upper summation limit cannot decrease and each existing
summand is nondecreasing. Thus $\abs{C}\leq\mathcal C_R(M_*,2(g-1))$.
Divide by $n\geq\ell$. The expression in
\eqref{eq:height-interval-even-cap} bounds every $M(n)$ by bounding
its first term above and the subtracted term below.
\end{proof}

All recursions and cut budgets in this section use integers; an
exact implementation is included in the supplement. The interval
calculation retains the completion-capacity hypothesis and does
not replace any arithmetic, probability, or Hall inequality.

For example, $w(1360)=w(1361)=46$, so the interval $[1360,1361]$
has $M_*=634$. Table~\ref{tab:height-construction-example} gives
the bounds for four types. At $D=3$ we use the fixed-level sequence
with $k=2$; at $4\leq D\leq8$ we use the adaptive sequence.
The tag variation is at most $D-1$ before retagging. The table
records only these construction quantities.

\begin{table}[ht]
\centering
\caption{Construction bounds for $M_*=634$ and four types.
Tag variation is measured before retagging the separator.}
\label{tab:height-construction-example}
\par\vspace{.5\baselineskip}
\begin{tabular}{lrrrrrr}
\hline
$D$&$3$&$4$&$5$&$6$&$7$&$8$\\
\hline
Number of cuts&$72$&$74$&$84$&$89$&$97$&$105$\\
Tag variation&$2$&$3$&$4$&$5$&$6$&$7$\\
\hline
\end{tabular}
\end{table}

\subsection{Assignments in the original forest}\label{sec:donor-construction}

A connected completion changes the degree sums used in the probability
estimates. We now give an alternative that uses the original forest
throughout. It preserves a third quantity, the number of opposite-side
positions of each type. The resulting separator may meet both sides.
We move its nonempty types to cheaper positions, then enlarge the tags
on the affected neighborhood. Every position remains available for
matching.

There are two costs to keep separate. Rounding creates signed errors
in the class degree sums. The subsequent exchanges transfer a
nonnegative amount of degree from each nonempty type to the empty
type. We retain both quantities until estimating the deletion means.

\subsubsection{Three measures and two choices of separator}

For nonnegative integers $v,e,d$, define the finite sum
\begin{equation}\label{eq:original-cut-budget}
 \mathcal C(v,e,d)=\sum_{j\geq1}
     \min\left\{d,\floor*{\frac v{2^j}},
                         \floor*{\frac e{2^j-1}}\right\}.
\end{equation}
It counts cuts against both the vertices and the edges available in
the forest. We also need a separator confined to $B$, particularly
when the degrees on $A$ are small.

\begin{lemma}[A cut confined to the opposite side]
\label{lem:original-side-cut}
Let $U=(A_U,B_U)$ be a bipartite tree with $v\geq2$ vertices and
$\Delta_U(A_U)\leq D$. If $D\geq2$, there is a vertex $b\in B_U$
such that every component $V$ of $U-b$ satisfies
\[
 \abs{V}-1\leq\frac{D-1}{D}(v-1).
\]
If $D=1$, a vertex of $B_U$ can be deleted so that every remaining
component is a singleton.
\end{lemma}

\begin{proof}
Choose an ordinary centroid $c$, as in the proof of
Lemma~\ref{lem:centroid-order}. If $c\in B_U$, all remaining
components have order at most $v/2$, which proves the first assertion.
Otherwise the branch orders at $c$ sum to $v-1$, and there are at
most $D$ branches. Let $t$ be the largest branch order and let $b$
be its neighbor of $c$. Then $(v-1)/D\leq t\leq v/2$.
After deleting $b$, the component containing $c$ has order $v-t$,
so its order minus one is at most $(D-1)(v-1)/D$. Every other
component has order at most $t-1$, and satisfies the same bound.
For $D=1$, every nontrivial component is a star with its center in
$B_U$; delete that center.
\end{proof}

\begin{theorem}[Assignment with the original degrees]
\label{thm:original-type-assignment}
Let $G=(A,B)$ be a bipartite forest, with $M=\abs A>0$,
$N=\abs B\geq0$, $e=e(G)$, and $\Delta_G(A)\leq D$, where
$D\geq1$ is an integer. Let $g\geq1$ be an integer. Prescribe nonnegative integer quotas
$a_0,\ldots,a_{g-1}$ of sum $M$, and put $\pi_i=a_i/M$.
There are colors $c(a)$ on $A$, tags $j(b)$ on $B$, and a set
$C\subseteq A\cup B$ with the following properties:
\begin{enumerate}
 \item the color quotas on $A$ are exactly $a_i$;
 \item every edge with both endpoints outside $C$ has $c(a)=j(b)$;
 \item writing $\beta_i=\abs{\{b:j(b)=i\}}$ and
 \[
   E_i=\sum_{a:c(a)=i}d_G(a),\qquad h_i=E_i-e\pi_i,
 \]
 one has
 \begin{equation}\label{eq:original-rounding-errors}
  \operatorname{tv}(\beta,N\pi)\leq g/4,\qquad
  \abs{h_i}\leq D,\qquad
  H:=\frac12\sum_i\abs{h_i}\leq gD/4;
 \end{equation}
 \item with $d=3(g-1)$, the ordinary construction has
 \[
  \abs C\leq\mathcal C(M+N,e,d).
 \]
 Alternatively, all cuts may be confined to $B$, with
 \begin{equation}\label{eq:original-separator-choices}
  \abs C\leq
  \begin{cases}
   \min(d,\floor{(M+N)/2}),&D=1,\\
   \mathcal S(M+N,d;(D-1)/D),&D\geq2.
  \end{cases}
 \end{equation}
\end{enumerate}
Here $\mathcal S$ is defined in \eqref{eq:weighted-separator-budget}.
No prescribed sign is asserted for any $h_i$.
\end{theorem}

\begin{proof}
\emph{Step 1. Preserve the three measures.}
Start with the connected components of $G$, each carrying color
weights $\pi_i$. Give an object $U$ the vector
\[
 v(U)=\left(\abs{U\cap A},\ \sum_{a\in U\cap A}d_G(a),\
                                      \abs{U\cap B}\right).
\]
All degrees are measured in the fixed original forest. Apply
Lemma~\ref{lem:fractional-pruning} with three coordinates. Freeze
pure objects; at most $d=3(g-1)$ mixed objects remain.

Choose a largest mixed object having more than one vertex. Cut it
at an ordinary centroid for the first construction, or at a vertex
of $B$ from Lemma~\ref{lem:original-side-cut} for the second.
Put the cut vertex in $C$. Replace the object by the cut vertex as
a singleton and the components of its deletion, each with the
parent's weights. Keeping the singleton is essential: a cut vertex
can contribute to the preserved measures. The vectors add to the
old vector, so all color targets remain $(a_i,e\pi_i,N\pi_i)$.
Prune again and continue. Each cut is at a new vertex; the process
ends with all mixed objects singletons.

\emph{Step 2. Count the cuts.}
In the ordinary construction, put an object in bucket $j\geq1$
when $2^j\leq\abs U<2^{j+1}$. Its children have order at most
$\floor{\abs U/2}<2^j$, and the cut atom has order one. When
bucket $j$ is first cut, no larger mixed bucket remains. There
are at most $d$ objects in that bucket, and no later cut creates
another. Objects cut in one bucket are pairwise disjoint. Each
contains at least $2^j$ vertices and, being a tree, at least
$2^j-1$ edges of $G$. Their number is therefore bounded by every
term of the minimum in \eqref{eq:original-cut-budget}.

For the second construction with $D\geq2$, use buckets for
$\abs U-1$ and contraction ratio $(D-1)/D$. The counting argument
in Proposition~\ref{prop:two-measure-assignment}, now counting all
vertices, gives the stated $\mathcal S$ bound. If $D=1$, each
cut leaves only singletons. Thus at most the initially mixed
nontrivial objects are cut; there are at most $d$, and their
disjoint orders are each at least two. This proves both choices.

\emph{Step 3. Round the singletons.}
The residual quotas on mixed $A$-singletons are nonnegative integers,
since pure objects contribute integer counts. Apply
Lemma~\ref{lem:averaged-slot-rounding} to their fixed degrees in
$[0,D]$. Its bounds give the last two inequalities in
\eqref{eq:original-rounding-errors}.

The residual targets on the $B$-singletons need not be integers.
Write them as $t_i=\floor{t_i}+r_i$, with $0\leq r_i<1$.
Their total is the integer number of these singletons, so
$k=\sum_i r_i$ is an integer. Round up the $k$ largest remainders
and round the others down. The resulting counts are nonnegative
integers with the correct total. Their variation from the targets is
\[
 k-\sum_{i\text{ rounded up}}r_i
       \leq k-\frac{k^2}{g}\leq\frac g4.
\]
For the middle inequality, the average of the $k$ largest remainders
is at least the average of all $g$; the case $k=0$ is immediate.
Assign the singletons to these counts. The two rounding operations
act on different sides and do not interfere.

Finally, an edge whose endpoints avoid $C$ stays within one final
object. If that object is pure, its endpoints have equal colors;
it cannot be a singleton. Hence the asserted edge property survives
both roundings. For $g=1$, color everything alike and take
$C=\varnothing$; all errors vanish.
\end{proof}

For an interval lower endpoint $\ell$ and an original forest with
at most $n$ vertices, $n\geq\ell$, the one-sided construction gives,
when $D\geq2$ and $\rho=(D-1)/D$,
\[
 \frac{\abs C}{n}\leq\sum_{j\geq0}
       \min\left\{\frac d\ell,\frac{\rho^j}{1+\rho^j}\right\}.
\]
Indeed, divide each summand of the $\mathcal S$ bound by $n$,
use $M+N\leq n$, and drop the ceiling and floor in the upper
bound. After the terms $0\leq j<J$, the remaining sum is at most
$\rho^J/(1-\rho)$. For $D=1$ the normalized bound is $d/\ell$.
Thus a finite evaluation with an explicit geometric remainder
also bounds the cuts throughout an order interval.

\subsubsection{Deleting the highest-degree even positions}

The following preparation uses all the special even labels on vertices
of largest degree. There is no reserve of odd positions in this variant.
It supplies an edge budget that also pays for isolated even positions.

\begin{proposition}[Preparation with highest-degree deletions]
\label{prop:highest-degree-preparation}
Assume the hypotheses of Corollary~\ref{cor:even-side}, and use its
partition $V(T)=S\sqcup Z\sqcup K$ and core labels $P$.
Let $m=\floor{n/2}$ and let $w$ be the count in
\eqref{eq:even-supply}. Suppose $w<m$. Choose $U\subseteq S$
to consist of $w$ vertices of largest degree in $T-K$, and put
\[
 A=S\setminus U,\quad B=Z,\quad G=T-(K\cup U),\quad
 M=m-w,\quad D=\Delta_G(A).
\]
Let $e=e(G)$, let $I$ be the number of isolated vertices in $A$,
and put $R=\sum_{a\in A}(d_G(a)-1)_+$. Then
\begin{align}
 \abs{V(G)}&\leq n-w,& e+I&\leq n-1-wD,
                    \label{eq:prepared-isolation-budget}\\
 R&\leq n-1-m-w(D-1)<n/2-w(D-1).
                    \label{eq:prepared-donor-excess}
\end{align}
Any bijective labeling of $A$ by $E_n\setminus E_*$ and of $B$
by $O_n\setminus\Psep$ that is coprime on every edge of $G$ extends,
after repair, to a prime labeling of $T$. If $D=0$, such a residual
labeling is automatic.
\end{proposition}

\begin{proof}
Since $S$ is independent, deleting $U$ does not change the degrees
of the remaining vertices of $S$ in $T-K$. Every selected vertex
therefore had degree at least $D$. Their incident edge sets are
disjoint, so at least $wD$ edges are deleted at this step.

If a vertex of $A$ is isolated in $G$, all its neighbors in $T$
are in $K$: it has no neighbor in $U\subseteq S$. It has at least
one such neighbor because $T$ is a nontrivial tree. Select one
edge to $K$ for each of these $I$ vertices. The selected edges are
distinct and disjoint from the $wD$ edges counted at $U$.
Together with the $e$ remaining edges, they prove
\eqref{eq:prepared-isolation-budget}.

All original degrees in $T$ are positive. Since $S$ is independent,
its degree sum is at most $n-1$. Hence
\[
 \sum_{s\in S}(d_T(s)-1)\leq n-1-m.
\]
Each of the $w$ removed vertices has original degree at least $D$,
so contributes at least $D-1$ to this sum. On $A$, residual degrees
cannot exceed original degrees, and
$(d_G(a)-1)_+\leq d_T(a)-1$. Subtracting proves the first bound
in \eqref{eq:prepared-donor-excess}. The strict inequality follows
from $n-1-\floor{n/2}<n/2$.

Assign $E_*$ bijectively to $U$. All labels are now used exactly
once, with odd labels on $B$ and the core, and even labels on $S$.
Edges meeting $K$ are coprime by Corollary~\ref{cor:even-side};
edges meeting $U$ are relaxed-compatible by
Lemma~\ref{lem:special-labels}. Every other edge is in $G$.
Thus Corollary~\ref{cor:prime-pair-repair} completes the labeling.
If $D=0$, there are no residual edges to check.
\end{proof}

\subsubsection{Counting the available empty positions}

We call a vertex of empty type a \emph{donor} when its type is
exchanged with a nonempty type on the separator. We need low-degree
donors, including isolated vertices. Their number depends on the
individual empty-class error, not just on its sign.

\begin{lemma}[Donor stock with a paid degree error]
\label{lem:paid-donor-stock}
Let a bipartite graph have $M>0$ vertices on side $A$, $e$ edges,
and $I$ isolated vertices in $A$. Suppose $e+I\leq E$.
A distinguished class has $a$ vertices, where $0<a\leq M$, and
degree sum $E_0\leq ea/M+h$, with $h\geq0$.
For an integer $s\geq1$, let $L_s$ be the number of its vertices
of degree at most $s$. Then
\begin{equation}\label{eq:paid-stock}
 L_s\geq\ell_s:=\max\left\{0,
   a-\floor*{\frac{Ea/M+h}{s+1}},
   \ceil*{\frac{a((s+1)M-E)-Mh}{(s+1)M-a}}\right\}.
\end{equation}
For data satisfying the hypotheses, $0\leq\ell_s\leq a$.
\end{lemma}

\begin{proof}
Put $t=s+1$. The high-degree vertices give
\[
 t(a-L_s)\leq E_0\leq ea/M+h,
 \qquad e\geq\frac Ma\bigl(t(a-L_s)-h\bigr).
\]
At least $L_s-I$ low-degree class vertices are nonisolated. Thus
\[
 E_0\geq t(a-L_s)+L_s-I.
\]
This remains a valid lower bound if $L_s-I<0$. Using $I\leq E-e$
and the upper bound on $E_0$ yields
\[
 ta-(t-1)L_s-E+e(1-a/M)-h\leq0.
\]
The coefficient of $e$ is nonnegative. Substituting its preceding
lower bound and simplifying gives
\[
 tM-E-\frac{tM-a}{a}L_s-\frac{Mh}{a}\leq0.
\]
The denominator $tM-a$ is positive, which proves the third term
in \eqref{eq:paid-stock}. Also $e\leq E$, so
$t(a-L_s)\leq Ea/M+h$. The integer $a-L_s$ is therefore at most
$\floor{(Ea/M+h)/t}$. This proves the second term. Finally
$0\leq L_s\leq a$ proves the asserted consistency of the bound.
\end{proof}

In Theorem~\ref{thm:original-type-assignment}, take the distinguished
class to be the empty type. Its individual error is at most $D$,
so the lemma applies with $h=D$, while the total degree variation
remains at most $gD/4$. We do not require the empty error to be
nonpositive. All these bounds concern the coloring before exchanges.
Also $e+I\geq M$: every nonisolated $A$-vertex contributes at least
one to the degree sum. In the preparation above, an order--degree
pair with $n-1-wD<M$ is therefore impossible. Such pairs may be
omitted from interval calculations.

\begin{lemma}[Choosing donors by degree]\label{lem:quantile-donors}
Let $G=(A,B)$ have maximum $A$-degree at most an integer $D\geq1$.
Give $A$ colors $0,\ldots,g-1$ and $B$ tags from the same set.
Suppose every edge with both endpoints outside a set $C$ joins
equal colors. For $1\leq s<D$, suppose $\ell_s$ is a nonnegative
integer lower bound for the number of color-$0$ vertices in $A$
of degree at most $s$. If $\abs C\leq c$, the coloring can be
changed, preserving all $A$-quotas and all $B$-tags, so that a set
$Z_C\subseteq B$ satisfies
\begin{equation}\label{eq:donor-boundary}
 \abs{Z_C}\leq\eta:=c+\sum_{s=1}^{D-1}(c-\ell_s)_+,
\end{equation}
and every edge to $B\setminus Z_C$ has $A$-color either zero or
equal to its $B$-tag. Each nonzero class loses a nonnegative
amount $\Delta_i$ of degree to class zero. If $\Delta=\sum_{i>0}\Delta_i$
and $R=\sum_{a\in A}(d_G(a)-1)_+$, then
\begin{equation}\label{eq:donor-transfer}
 \Delta\leq\min\{Dc,\ c+R\}.
\end{equation}
For $D=1$ or $C\subseteq B$, no exchange is needed and $\Delta=0$.
\end{lemma}

\begin{figure}[tbp]
\centering
\begin{tikzpicture}[x=1cm,y=1cm,>=Stealth,
  a/.style={circle,draw,minimum size=8mm,inner sep=1pt},
  b/.style={rectangle,draw,minimum size=7mm,inner sep=1pt}]
\foreach \off/\lefttype/\righttype/\lasttag in {0/I/\varnothing/\varnothing,6/\varnothing/I/I} {
 \node[a] (v\off) at (\off,1.1) {$\lefttype$};
 \node[a] (u\off) at (\off,-1.1) {$\righttype$};
 \node[left=2mm of v\off] {$v$};
 \node[left=2mm of u\off] {$u$};
 \foreach \j/\yy in {1/1.9,2/1.1,3/.3} {
  \node[b] (b\off\j) at (\off+1.5,\yy) {$\varnothing$};
  \draw (v\off)--(b\off\j);
 }
 \node[b] (c\off) at (\off+1.5,-1.1) {$\lasttag$};
 \draw (u\off)--(c\off);
}
\draw[->] (2.4,.4)--(4.8,.4);
\node at (3.6,.85) {exchange};
\node at (3.6,-.05) {and retag};
\end{tikzpicture}
\caption{An exchange with source degree three and donor degree one.
The entries at $v,u$ are $A$-types; the other entries are $B$-tags.
Initially $v\in C$, $u\notin C$, and the displayed tags are empty.
The exchange moves two units of degree to the empty type. Only the
donor's neighbor needs retagging in this example. Solid lines are
forest edges; the arrow indicates the change of assignment.}
\label{fig:donor-exchange}
\end{figure}
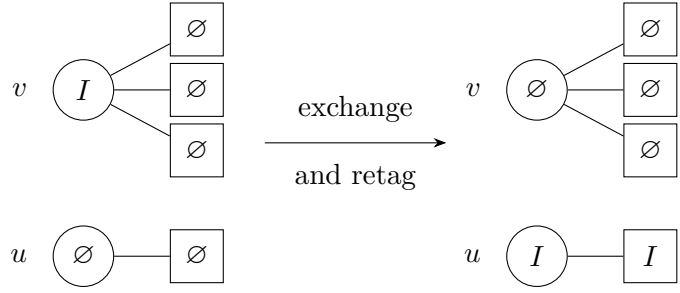

\begin{proof}
Let $X=C\cap A$, let $Z$ be its vertices of color zero, and
put $Q=X\setminus Z$, $q=\abs Q$. The $q$ nonzero colors can be
placed at vertices of
\[
 P=Q\ \sqcup\ \{a\in A\setminus C:c(a)=0\}.
\]
Assign a candidate the integer cost $\gamma(a)=\max\{1,d_G(a)\}$.
Choose the $q$ cheapest candidates, preferring members of $Q$
whenever costs tie, and call this set $W$. Members of $Q\cap W$
keep their colors. Pair $Q\setminus W$ bijectively with
$W\setminus Q$ and exchange their colors.

Every discarded member of $Q$ has strictly larger cost than every
newly selected member outside $Q$: an equal cost would have favored
the old member. It consequently has degree at least two and strictly
larger degree than its donor. Thus each exchange preserves the
quotas and transfers a positive degree difference from its nonzero
class to class zero. The cost $\max\{1,d_G(a)\}$ ensures that a
separator vertex of degree zero or one is never moved.

Set
\[
 Z_C=(C\cap B)\cup N_G(W).
\]
An unchanged nonzero vertex outside $W$ is outside $C$, and its
edges to $B\setminus Z_C$ retain their old equal tags. A moved
separator vertex is now zero. A newly nonzero donor is in $W$,
and all its neighbors are in $Z_C$. These cases check every edge.
Figure~\ref{fig:donor-exchange} illustrates one exchange.

Among the candidates there are at least $\ell_s-\abs Z$ vertices
of cost at most $s$, since only the empty cut vertices $Z$ may
have been excluded from the empty class. Therefore the cheapest
choice satisfies
\[
 \abs{\{a\in W:\gamma(a)>s\}}
        \leq(\abs X-\ell_s)_+.
\]
Indeed, when $\ell_s\geq\abs Z$ this follows by subtracting the
available cheap candidates from $q=\abs X-\abs Z$. Otherwise the
right side is at least $q$, which is already a bound. Summing the
integer cost layers gives
\begin{align*}
 \sum_{a\in W}d_G(a)
 &\leq\sum_{a\in W}\gamma(a)\\
 &=q+\sum_{s=1}^{D-1}\abs{\{a\in W:\gamma(a)>s\}}\\
 &\leq\abs X+\sum_{s=1}^{D-1}(\abs X-\ell_s)_+.
\end{align*}
Add $\abs{C\cap B}$ and use $\abs X\leq\abs C\leq c$ to obtain
\eqref{eq:donor-boundary}.

Let $V=Q\setminus W$ be the moved separator vertices and $u(v)$
the donor paired with $v$. Their total transfer is
\[
 \Delta=\sum_{v\in V}\bigl(d_G(v)-d_G(u(v))\bigr)
       \leq\sum_{v\in V}d_G(v)
       \leq\min\{D\abs V,\ \abs V+R\}.
\]
The last inequality uses $d_G(v)\geq2$ and $V\subseteq C\cap A$.
This proves \eqref{eq:donor-transfer}. If $D=1$, all candidate
costs are one and the tie rule chooses $W=Q$. If $C\subseteq B$,
then $Q=\varnothing$. In either case there are no exchanges.
\end{proof}

The boundary formula uses the entire sequence of degree counts.
For example, if $\ell_2\geq c$, it gives
$\eta\leq c+(c-\ell_1)_+$, provided the later bounds are also at
least $c$. More generally, the actual degree counts are nondecreasing;
any lower bounds may be replaced by their successive maxima.
No hypothesis that $c$ donors all have the same degree cap is needed.

\subsubsection{Paying for the transferred degree}

The next identity keeps the contribution removed from each nonempty
class. Replacing all transfers by an additional unsigned rounding error
would lose this cancellation.

\begin{lemma}[Weighted means after exchanges]\label{lem:donor-mean-cancellation}
Suppose the initial color quotas and degree sums are $a_i,E_i$,
with $M=\sum_i a_i>0$, $e=\sum_i E_i$, and
$h_i=E_i-ea_i/M$. Let $H=\frac12\sum_i\abs{h_i}$.
Suppose exchanges produce
\[
 E'_0=E_0+\Delta,\qquad E'_i=E_i-\Delta_i\quad(i>0),
 \qquad \Delta_i\geq0,\quad\Delta=\sum_{i>0}\Delta_i.
\]
For $g\geq2$ and real class means $f_i$, one has
\begin{align}
 \sum_i E'_i f_i
 &=\frac eM\sum_i a_i f_i+\sum_i h_i f_i
                      +\sum_{i>0}\Delta_i(f_0-f_i)
                       \label{eq:donor-mean-identity}\\
 &\leq\frac eM\sum_i a_i f_i
      +H(\max_i f_i-\min_i f_i)
      +\Delta(f_0-\min_{i>0}f_i)_+.
                       \label{eq:donor-mean-bound}
\end{align}
When $g=1$, there are no transfers and the identity holds with both
error terms zero.
\end{lemma}

\begin{proof}
Substitute the formulas for $E'_i$ and then for $E_i$ to obtain
\eqref{eq:donor-mean-identity}. Since $\sum_i h_i=0$, the positive
and negative parts each have mass $H$. Their weighted sum is at
most $H(\max_i f_i-\min_i f_i)$. Each $\Delta_i$ is nonnegative,
so the last sum is at most the last term in
\eqref{eq:donor-mean-bound}. No sign condition on $h_0$ is used.
\end{proof}

\begin{theorem}[Retained donor assignment]
\label{thm:retained-donor-assignment}
Let $G=(A,B)$, $M,N,e,D$ be as in
Theorem~\ref{thm:original-type-assignment}, with $N>0$.
Use the $g=2^s$ types $I\subseteq[s]$ as colors, with exact quotas
$a_I$ and $a_{\varnothing}>0$. Suppose $e+I_0\leq E$, where
$I_0$ counts isolated $A$-vertices, and let $R_*$ bound
$\sum_{a\in A}(d_G(a)-1)_+$ from above.
Choose either construction in that theorem and an integer bound
$c\geq\abs C$. Define $\ell_t$ by \eqref{eq:paid-stock} with
$a=a_{\varnothing}$, $h=D$, and $1\leq t<D$. Put
\[
 \eta=c+\sum_{t=1}^{D-1}(c-\ell_t)_+,
 \qquad \delta=\min\{Dc,c+R_*\}.
\]
There are an exact type assignment on $A$ and tags $J'(b)$ on all
of $B$ such that $I(a)\subseteq J'(b)$ on every edge and
\begin{equation}\label{eq:donor-retained-tag-bound}
 \operatorname{tv}(\beta',N\pi)\leq g/4+\eta.
\end{equation}
Their degree sums are obtained from initial errors
$\abs{h_I}\leq D$, $H\leq gD/4$, by nonnegative transfers toward
the empty type of total at most $\delta$.
If $C\subseteq B$, one may instead take $\eta=c$ and $\delta=0$;
if $D=1$, one may take $\delta=0$.

Fix finite half-label sets $\Lambda_I$ of cardinalities $a_I$, with
every member having type $I$, and permute them independently and
uniformly onto the final classes. For any real functions $F_I$ on
nonempty $\Lambda_I$, put
\[
 f_I=\frac1{a_I}\sum_{j\in\Lambda_I}F_I(j).
\]
Assign arbitrary means to classes with $a_I=0$. The expected cost
\[
 \E\sum_{a\in A}d_G(a)F_{I(a)}(j(a))
\]
is bounded by \eqref{eq:donor-mean-bound}, with $H\leq gD/4$ and
$\Delta\leq\delta$. For block types and nonnegative conflict
costs, these estimates bound the local and global deletions of
Lemma~\ref{lem:retained-tags}.
\end{theorem}

\begin{proof}
Apply Theorem~\ref{thm:original-type-assignment}. Its individual
empty error is at most $D$, so Lemma~\ref{lem:paid-donor-stock}
justifies every $\ell_t$. Apply Lemma~\ref{lem:quantile-donors};
it preserves the quotas and all the old tags, and supplies $Z_C$
of size at most $\eta$. Its transfer bound gives $\Delta\leq\delta$.
Now apply Lemma~\ref{lem:retained-tags}, replacing each tag on
$Z_C$ by the union of its neighboring types. This proves the
inclusion on every edge and \eqref{eq:donor-retained-tag-bound}.
One may alternatively use the full tag $[s]$ on $Z_C$, with the
same bound. Every position remains in $B$.

For the one-sided construction, $C\subseteq B$ already gives the
required edge property outside $C$, so no exchanges are needed and
only $C$ is retagged. The case $D=1$ follows from
Lemma~\ref{lem:quantile-donors}.

All colors, exchanges, and tags are fixed before the permutations.
At each vertex of class $I$, the assigned half-label is uniform on
$\Lambda_I$, so linearity of expectation gives the mean
$\sum_I E'_I f_I$. Empty classes contribute zero before and after
exchanges. Lemma~\ref{lem:donor-mean-cancellation} supplies the
stated bound. For a fixed odd label, use its conflict indicator,
setting it to zero in incompatible types. For the global cost,
use the number of compatible odd labels in conflict, or any
nonnegative upper bound for that number. The deterministic
estimates in Lemma~\ref{lem:retained-tags} then apply.
\end{proof}

The degree lists still consist of original degrees in $[0,D]$.
Consequently the class-sampling bounds of
Section~\ref{sec:local-probability} apply after the exchanges, with
the new means just proved. There is no completion credit in this
construction, and $H$ bounds the initial signed errors, not the
variation of the final degree sums.

\begin{corollary}[Donor bounds on an integer interval]
\label{cor:donor-interval}
Fix integers $1\leq\ell\leq u$, $D\geq1$, and $g\geq1$.
Suppose Proposition~\ref{prop:highest-degree-preparation} applies
at every order $n\in[\ell,u]$ under consideration, with residual
degree $D$. Write $w_n=w(n)$, $m_n=\floor{n/2}$,
$M_n=m_n-w_n$, and let $a_n>0$ be the exact empty quota.
Put $\mathcal E_n=n-1-Dw_n$ and retain only orders satisfying
$\mathcal E_n\geq M_n$, writing $\mathcal I$ for their set.
Assume $\mathcal I\ne\varnothing$; omitted orders cannot have
residual degree $D$.
For the ordinary construction define
\begin{align*}
 c_*&=\max_{n\in\mathcal I}
       \mathcal C(n-w_n,\mathcal E_n,3(g-1)),\\
 \eta_*&=\max_{n\in\mathcal I}
       \left(c_*+\sum_{t=1}^{D-1}(c_*-\ell_t(n))_+\right),\\
 \delta_*&=\max_{n\in\mathcal I}
       \min\{Dc_*,\ c_*+n-1-m_n-w_n(D-1)\},
\end{align*}
where $\ell_t(n)$ is \eqref{eq:paid-stock} with
$(M,a,E,h)=(M_n,a_n,\mathcal E_n,D)$. Then uniformly in the interval,
\[
 \frac{\abs C}{n}\leq\frac{c_*}{\ell},\qquad
 \frac{\abs{Z_C}}n\leq\frac{\eta_*}{\ell},\qquad
 \frac\Delta n\leq\frac{\delta_*}{\ell},\qquad
 \frac{H}{n}\leq\frac{gD}{4\ell}.
\]
In particular the normalized tag variation is at most
$(g/4+\eta_*)/\ell$. For $D=1$, set $\delta_*=0$.
\end{corollary}

\begin{proof}
The cut budget is nondecreasing in both its vertex and edge arguments.
Proposition~\ref{prop:highest-degree-preparation} therefore justifies
$c_*$ for every order. Its isolation bound justifies the paid stock
formula, and its excess bound justifies the expression for $\delta_*$.
Apply Theorem~\ref{thm:retained-donor-assignment} separately at each
order, and divide its bounds by $n\geq\ell$. The maxima are over
finite sets of integers, with exact rational calculations in the stock
formula. There is no union bound over orders.
\end{proof}

For example, at $n=1360$ and $D=9$, using four types formed from
$3,5,7,11,13$, the exact arithmetic gives
\[
 w=46,\qquad M=634,\qquad a_{\varnothing}=216,\qquad E=945.
\]
The ordinary cut bound is $c=65$. Paying $h=9$ gives
$\ell_1=61$, $\ell_2=120$, and all later stock bounds at least $120$.
Thus the boundary bound is $65+(65-61)=69$. These are construction
quantities only; the example does not check a labeling range.

The stock counts at different degree thresholds belong to the same
order. Taking the maximum of their combined cost, as above, retains
that dependence. These construction bounds still require the arithmetic
population estimates and the matching inequalities before they yield
a numerical range of prime labelings.

\subsection{Finite arithmetic for the matching estimates}\label{sec:finite-arithmetic}

The constructions in the preceding two sections fix the type at each
even position and the tag at each odd position. We now express the
remaining expectations and population bounds by finite arithmetic.
The distinction between the original forest and a connected completion
remains in force: the former allows degree transfers to the empty
class, while the latter allows a credit for the added edges.

Throughout this section, fix integers $12\leq\ell\leq u$ and a cutoff
$3\leq z<\ell/3$. Partition all odd primes at most $z$ into $s\geq1$
nonempty blocks. Write $\mathcal T=2^{[s]}$, $0=\varnothing$, and
$F=[s]$ for the type set, its empty type, and its full type.
The block type $T$ and the complete counts $\bar C_I$ are those of
Section~\ref{sec:arithmetic-types}. All extrema over $[\ell,u]$
in this section range over integers.

\subsubsection{The removed half-labels}

For an order $n$, put $m_n=\floor{n/2}$ and let
\[
 W_n=\{2^a:a\geq0,\ 2^a\leq m_n\}
       \mathbin{\dot\cup}\{p:n/3<p\leq m_n,\ p\text{ prime}\}.
\]
These are exactly the half-labels corresponding to the special even
labels. Set
\[
 r_n=\pi(m_n)-\pi(\floor{n/3}),\qquad
 w_n=\abs{W_n}=\floor{\log_2 n}+r_n.
\]
The usable half-labels and their quotas are
\begin{equation}\label{eq:finite-usable-classes}
 \Lambda_I(n)=\{j\in[m_n]\setminus W_n:T(j)=I\},\qquad
 a_I(n)=\bar C_I(m_n)-w_n\boldsymbol1_{I=0}.
\end{equation}
For a finite set $P$ of primes greater than $z$, let $B_I(m,P)$
count the integers $j\leq m$ of type $I$ divisible by some $p\in P$.
Lemma~\ref{lem:large-prime-counts} gives this count by
inclusion--exclusion, whether or not $P$ is the support of a label.

We will count deletions involving a specified family
$\mathcal J\subseteq\mathcal T$ of odd-label types. Define
\begin{align}
 O_I^{\mathcal J}(v)
  &=\sum_{\substack{J\in\mathcal J\\I\cap J=\varnothing}}
        \bigl(\bar C_J(v)-\bar C_J(\floor{v/2})\bigr),
                    \label{eq:restricted-odd-cost}\\
 Q_{n,I}^{\mathcal J}(j)
  &=\sum_{\substack{p>z,\ p\text{ prime}\\p\mid j}}
                   O_I^{\mathcal J}(\floor{n/p}),
                    \label{eq:restricted-prime-cost}\\
 S_I^{\mathcal J}(n)
  &=\sum_{\substack{j\leq m_n\\T(j)=I}}Q_{n,I}^{\mathcal J}(j).
                    \label{eq:restricted-class-total}
\end{align}
Here $v$ is a nonnegative integer. These are the filtered costs of
Lemma~\ref{lem:filtered-prime-costs}, restricted to $\mathcal J$.
They dominate the number of conflicts with usable odd labels in that
family. They may count a conflict more than once.

\begin{lemma}[Exact removal identities]\label{lem:finite-special-removal}
Every member of $W_n$ has empty type. Write
$\kappa=\boldsymbol1_{0\in\mathcal J}$ and
$r_n(P)=\abs{\{p\in P:n/3<p\leq m_n\}}$. Then
\begin{align}
 \sum_{j\in\Lambda_I(n)}Q_{n,I}^{\mathcal J}(j)
    &=S_I^{\mathcal J}(n)-\kappa r_n\boldsymbol1_{I=0},
                        \label{eq:removed-first-moment}\\
 \abs{\{j\in\Lambda_I(n):\text{some }p\in P\text{ divides }j\}}
    &=B_I(m_n,P)-r_n(P)\boldsymbol1_{I=0}.
                        \label{eq:removed-union-count}
\end{align}
For every $t\geq0$, one also has
\begin{equation}\label{eq:removed-exponential-moment}
 \begin{split}
 &\sum_{j\in\Lambda_I(n)}
           \bigl(e^{tQ_{n,I}^{\mathcal J}(j)/n}-1\bigr)\\
 &\quad=\sum_{\substack{j\leq m_n\\T(j)=I}}
           \bigl(e^{tQ_{n,I}^{\mathcal J}(j)/n}-1\bigr)
       -\kappa r_n(e^{t/n}-1)\boldsymbol1_{I=0}.
 \end{split}
\end{equation}
Moreover,
\begin{equation}\label{eq:restricted-total-prime-sum}
 S_I^{\mathcal J}(n)=
   \sum_{\substack{z<p\leq m_n\\p\text{ prime}}}
       \bar C_I(\floor{m_n/p})\,O_I^{\mathcal J}(\floor{n/p}).
\end{equation}
Each complete total is nondecreasing in $n$.
\end{lemma}

\begin{proof}
A power of two has no odd prime divisor. A prime in $(n/3,m_n]$
is greater than $z$; thus both sorts of removed half-label have empty
type. The former has cost zero. For the latter, the only odd multiple
of $p$ in $[n]$ is $p$ itself, whose type is empty. Its cost is
therefore $\kappa$. It belongs to the union counted by $B_0(m_n,P)$
exactly when $p\in P$. Subtracting these contributions proves the
three identities. There are $\floor{\log_2 n}$ powers of two in
$[m_n]$, which also justifies the formula for $w_n$.

Multiplication by a prime greater than $z$ preserves the block type.
Interchanging the sums over half-labels and their prime divisors gives
\eqref{eq:restricted-total-prime-sum}, just as in
Lemma~\ref{lem:filtered-prime-costs}. For each fixed half-label,
the cost cannot decrease when the odd prefix grows. Increasing $n$
also adds only nonnegative terms to the complete class sum, proving
monotonicity.
\end{proof}

\subsubsection{Uniform bounds for the class means}

Write $w_- =\min_nw_n$, $w_+=\max_nw_n$, and define $r_-,r_+$
similarly. Put
\begin{equation}\label{eq:finite-population-endpoints}
 \begin{split}
 c_I^-&=\bar C_I(m_\ell),\qquad c_I^+=\bar C_I(m_u),\\
 a_I^-&=c_I^- -w_+\boldsymbol1_{I=0},\qquad
 a_I^+=c_I^+ -w_-\boldsymbol1_{I=0},\qquad
 M^-=m_\ell-w_+.
 \end{split}
\end{equation}
Assume $a_I^->0$ for every $I$. Thus all actual quotas are positive,
$M(n)=m_n-w_n\geq M^->0$, and $a_I^-\leq a_I(n)\leq a_I^+$.
This positivity is a hypothesis to be checked, not a consequence of
the construction for arbitrary small orders.

\begin{proposition}[Arithmetic bounds for the global means]
\label{prop:finite-global-means}
For a fixed family $\mathcal J$, define
\begin{align*}
 U_I&=\frac{S_I^{\mathcal J}(u)-\kappa r_-\boldsymbol1_{I=0}}
                    {\ell a_I^-},&
 L_I&=\frac{(S_I^{\mathcal J}(\ell)-\kappa r_+\boldsymbol1_{I=0})_+}
                    {u a_I^+},\\
 A_*&=\frac{\sum_I S_I^{\mathcal J}(u)-\kappa r_-}{\ell M^-}.
\end{align*}
Also put
\[
 B_*=\max_I U_I-\min_I L_I,\qquad
 C_*=(U_0-\min_{I\ne0}L_I)_+.
\]
If $f_I(n)$ is the mean of $Q_{n,I}^{\mathcal J}$ on $\Lambda_I(n)$,
then
\begin{equation}\label{eq:finite-mean-brackets}
 L_I\leq f_I(n)/n\leq U_I,\qquad
 \frac{\sum_Ia_I(n)f_I(n)}{nM(n)}\leq A_*.
\end{equation}
Suppose an original-forest assignment has initial degree variation
$H$, edge count $e$, and donor transfer $\Delta$, with
$e/n\leq e_+$, $H/n\leq H_+$, and $\Delta/n\leq\delta$.
For the sum $S$ charging deletions in this family,
\begin{equation}\label{eq:finite-donor-global-mean}
 \E S/n^2\leq e_+A_*+H_+B_*+\delta C_*.
\end{equation}
For an assignment in a completing tree $Q$, with no donor transfers,
the analogous bound for original-forest costs is
\begin{equation}\label{eq:finite-completed-global-mean}
 \E S/n^2\leq e_+A_*+H_+B_*-q_-\min_I L_I,
\end{equation}
provided $e(Q)/n\leq e_+$, its degree variation is at most $nH_+$,
and $(e(Q)-e(G))/n\geq q_-\geq0$.
\end{proposition}

\begin{proof}
By Lemma~\ref{lem:finite-special-removal}, the numerator of the actual
class mean is $S_I^{\mathcal J}(n)-\kappa r_n\boldsymbol1_{I=0}$.
Bound its complete sum by the corresponding endpoint, its subtraction
by the opposite extremum, and its positive denominator by
\eqref{eq:finite-population-endpoints}. The numerator is nonnegative,
so the lower bound may be truncated at zero. This proves the first
part of \eqref{eq:finite-mean-brackets}. Sum the removal identity over
$I$ and use $nM(n)\geq\ell M^-$ to prove its second part.

Apply Lemma~\ref{lem:donor-mean-cancellation} to the class means $f_I(n)$.
After division by $n^2$, the three terms are bounded by
$e_+A_*$, $H_+B_*$, and $\delta C_*$. This proves
\eqref{eq:finite-donor-global-mean}. For the completed assignment,
apply \eqref{eq:completion-cost-credit} instead, retaining the edge
count and variation in $Q$ and subtracting the lower bound for the
added-edge cost. All class means are nonnegative, so this credit has
the asserted direction.
\end{proof}

It is sometimes convenient to use weaker endpoint bounds. For example,
one may omit the subtraction in $U_I$ and $A_*$, replace $r_+$ by
$w_+$ in $L_0$, or replace $a_I^+$ by $c_I^+$. Each change increases
an upper bound or decreases a lower bound. In particular, the removal
identities do not require the sets $W_n$ to be nested.

\subsubsection{Local frequencies and the number sampled}

For an actual odd label $y$, let $P(y)$ be the set of its prime
divisors in $(z,m_u]$, and put $J=T(y)$. Larger primes divide no
available half-label. More generally, fix any type $J$ and finite
set $P$ of primes in $(z,m_u]$.
Let $\mathcal A_J=\{I:I\cap J=\varnothing\}$. A type outside
$\mathcal A_J$ contributes zero to the local conflict count. Define
\begin{gather*}
 b_I^-=B_I(m_\ell,P),\qquad b_I^+=B_I(m_u,P),\\
 r^*(P)=\abs{\{p\in P:2p\leq u,\ 3p>\ell\}}.
\end{gather*}
The last number includes every possibly removed conflicting prime
at every order of the interval. It may include primes that are never
removed at the same order.

\begin{lemma}[Local arithmetic data]\label{lem:finite-local-data}
For $I\in\mathcal A_J$, put
\[
 U_I^P=\min\{1,b_I^+/a_I^-\},\qquad
 L_I^P=\frac{(b_I^- -r^*(P)\boldsymbol1_{I=0})_+}{a_I^+}.
\]
Set both quantities to zero for $I\notin\mathcal A_J$. Let $f_I$
be the fraction of usable half-labels in class $I$ contributing to
the event $(J,P)$, including this compatibility restriction. Then
$L_I^P\leq f_I\leq U_I^P$. The following are uniform bounds:
\begin{align}
 A_P&=\frac{\sum_{I\in\mathcal A_J}b_I^+}{M^-},&
 H_P&=\max_I U_I^P,&
 G_P&=(U_0^P-\min_{I\ne0}L_I^P)_+,
                                      \label{eq:finite-local-profiles}\\
 L_P&=\min_I L_I^P,&
 b_-&=\frac{(\sum_{I\in\mathcal A_J}b_I^- -r^*(P))_+}{u},&
 b_+&=\frac{\sum_{I\in\mathcal A_J}b_I^+}{\ell}.
                                      \label{eq:finite-local-support-bounds}
\end{align}
More precisely, $\sum_Ia_If_I/M\leq A_P$, $\max_I f_I\leq H_P$,
$(f_0-\min_{I\ne0}f_I)_+\leq G_P$, and $\min_I f_I\geq L_P$.
If $b=\sum_Ia_If_I$ is the number of selected half-labels, then
\begin{equation}\label{eq:finite-selected-population}
 b_-\leq b/n\leq b_+,\qquad
 b\leq B_P:=\sum_{I\in\mathcal A_J}b_I^+.
\end{equation}
If $J\ne0$, then $L_P=0$.
\end{lemma}

\begin{proof}
Apply \eqref{eq:removed-union-count}. Its complete count lies between
$b_I^-$ and $b_I^+$, and at most $r^*(P)$ members are removed, all
from class $0$. Divide by the appropriate endpoint population. The
upper bound may be truncated at one because it bounds a frequency.
For a forbidden class the frequency is identically zero.

Summing the exact usable counts over allowed classes gives $b$.
Bounding above by the complete upper prefix gives both $b\leq B_P$
and $b/M\leq A_P$. Bounding below by the complete lower prefix,
subtracting $r^*(P)$, and dividing by $n\leq u$ gives $b_-$.
The remaining bounds follow by taking maxima and minima of the class
brackets. When $J\ne0$, the full class is forbidden and its bracket
is zero, so $L_P=0$.
\end{proof}

\begin{proposition}[Means for the local probability bounds]
\label{prop:finite-local-means}
Use the arithmetic data of Lemma~\ref{lem:finite-local-data}.
For an original-forest donor assignment with normalized bounds
$e_+,H_+,\delta$ as in Proposition~\ref{prop:finite-global-means},
\begin{equation}\label{eq:finite-original-local-mean}
 \E K(J,P)/n\leq e_+A_P+H_+H_P+\delta G_P.
\end{equation}
For a connected assignment with $\Delta_Q(A)\leq D$ and no donor
transfers, put $\mu_Q=e_+A_P+H_+H_P$, using the edge count and variation in $Q$.
If the normalized number of added edges is at least $q_-$, then
\begin{equation}\label{eq:finite-connected-local-mean}
 \E K^G(J,P)/n\leq\mu_Q-q_-L_P.
\end{equation}
Alternatively, the completed count is $K^Q=b+X$, where $X$ samples
weights in $[0,D-1]$, and
\begin{equation}\label{eq:finite-baseline-data}
 \E X/n\leq\mu_Q-b_-,\qquad
 \{K^G/n\geq\theta\}\subseteq\{X/n\geq\theta-b_+\}.
\end{equation}
All these samples have size at most the integer $B_P$.
\end{proposition}

\begin{proof}
Apply the degree-transfer identity to the frequencies $f_I$.
The average term is at most $eA_P$, the oscillation is at most
$H_P$ because every frequency is nonnegative, and the transfer
coefficient is at most $G_P$. Divide by $n$ to prove
\eqref{eq:finite-original-local-mean}. In the connected case, use
Lemma~\ref{lem:degree-cost-comparison} and subtract the cost of the
added edges, whose total is at least $nq_-L_P$. This gives
\eqref{eq:finite-connected-local-mean}.

The baseline assertion is Proposition~\ref{prop:original-local-degrees},
with the bounds for $b/n$ just proved. It uses the uncredited mean
$\mu_Q$: the completion credit is not subtracted a second time.
The sample-size assertion follows from \eqref{eq:finite-selected-population}.
\end{proof}

These quantities supply the inputs to the integer and convex bounds
of Section~\ref{sec:local-probability}. For its forest moment estimate,
one may also use $\max_I f_I\leq H_P$ and any proved upper bound
for the total original excess degree. A local event is indexed by
its type $J$ and its prime support $P$. Inclusion-maximal supports
may replace all supports only with the common-threshold condition
of Lemma~\ref{lem:maximal-local-supports}; equality of numerical
profiles alone does not identify events.

\subsubsection{Endpoint bounds for exponential moments}

For the global probability estimate we need more than a first moment.
The same finite costs supply a histogram of integer values in each
complete class. In this subsection fix $\mathcal J$ and $t>0$, and
abbreviate $Q_{n,I}^{\mathcal J}$ to $Q_{n,I}$.

\begin{lemma}[Bounds for the logarithmic class moments]
\label{lem:finite-log-moments}
Define
\begin{align*}
 V_I^+&=\sum_{\substack{j\leq m_u\\T(j)=I}}
                (e^{tQ_{u,I}(j)/\ell}-1)
             -\kappa r_-(e^{t/\ell}-1)\boldsymbol1_{I=0},\\
 V_I^-&=\left(\sum_{\substack{j\leq m_\ell\\T(j)=I}}
                (e^{tQ_{\ell,I}(j)/u}-1)
             -\kappa r_+(e^{t/u}-1)\boldsymbol1_{I=0}\right)_+,\\
 v_I^-&=\log(1+V_I^-/a_I^+),\qquad
 v_I^+=\log(1+V_I^+/a_I^-).
\end{align*}
Then $V_I^+\geq0$, and the actual logarithmic mean satisfies
\begin{equation}\label{eq:finite-log-brackets}
 v_I^-\leq\psi_I(n,t):=
    \log\left(\frac1{a_I(n)}
        \sum_{j\in\Lambda_I(n)}e^{tQ_{n,I}(j)/n}\right)
                 \leq v_I^+.
\end{equation}
\end{lemma}

\begin{proof}
It suffices to bound the sum of exponential excesses over usable
half-labels, since the mean in \eqref{eq:finite-log-brackets} is one
plus that sum divided by the usable population.

For the upper bound, enlarge the prefix to $[m_u]$, replace each
cost by $Q_{u,I}$, and replace $n$ in the exponent by $\ell$.
Every usable term can only increase. If $0\in\mathcal J$, each of
the $r_n$ removed prime half-labels is still in this enlarged prefix
and has $Q_{u,0}(p)\geq1$: the odd label $p$ itself is counted.
Their omitted exponential excess is therefore at least
$r_-(e^{t/\ell}-1)$. If $0\notin\mathcal J$, omitting the
subtraction remains valid. Powers of two contribute zero in either
case. This proves the upper bound by $V_I^+$ and its nonnegativity.

For the lower bound, use only the half-labels in $[m_\ell]$ and
replace their costs by $Q_{\ell,I}$ and the denominator by $u$.
A removed prime present in this smaller prefix satisfies
$p>n/3\geq\ell/3$. Its only odd multiple in $[\ell]$ is $p$,
so its lower-endpoint cost is exactly $\kappa$. At most $r_+$
such primes are removed. Subtracting that many unit-cost excesses
and truncating at zero gives $V_I^-$. Finally divide by the
population bounds and use monotonicity of the logarithm.
\end{proof}

The proof uses inclusions of the usable sets into the complete
endpoint prefixes. Subtracting a removed term after independently
bounding $n$ in its exponent would give different, generally weaker,
bounds. The chosen subtractions are justified by the endpoint costs
of the removed primes themselves.

\begin{corollary}[A uniform global tail bound]\label{cor:finite-global-tail}
Write $\pi_I(n)=a_I(n)/M(n)$ and set
\begin{gather*}
 A_t=\max_{\ell\leq n\leq u}\sum_I\pi_I(n)v_I^+,\\
 B_t=\max_Iv_I^+-\min_Iv_I^-,\qquad
 C_t=(v_0^+-\min_{I\ne0}v_I^-)_+.
\end{gather*}
In the original-forest construction put
$\mathcal M=e_+A_t+H_+B_t+\delta C_t$.
In the connected construction put
$\mathcal M=e_+A_t+H_+B_t-q_-\min_Iv_I^-$,
using the corresponding bounds from
Proposition~\ref{prop:finite-global-means}.
Let $Z_{\mathcal J}$ count the deleted matching edges incident with
labels of types in $\mathcal J$. If $D>0$, $d_G(a)\leq D$, $b>0$, and
$t b>\mathcal M$, then, at every order in the interval,
\begin{equation}\label{eq:finite-global-tail}
 \Pr\{Z_{\mathcal J}\geq n^2b\}
       \leq\exp\left(-\frac{\ell}{D}(tb-\mathcal M)\right).
\end{equation}
Without the strict inequality one may use the bound one.
\end{corollary}

\begin{proof}
Apply Proposition~\ref{prop:global-conflict-tails} with parameter $t/n$.
The class log moments are exactly $\psi_I(n,t)$. Degree variation
and either donor cancellation or the completion credit give
\[
 \frac1n\sum_I E_I^G\psi_I(n,t)\leq\mathcal M.
\]
The average term uses the actual quota vector, bounded by the finite
maximum $A_t$. Markov's inequality now bounds the probability by
$\exp(-n(tb-\mathcal M)/D)$. For a positive rate, $n\geq\ell$
gives \eqref{eq:finite-global-tail}. All types, tags and transfers
were fixed before the class permutations, as required for that
moment calculation.
\end{proof}

The same inequalities hold if each $v_I^-$ is replaced by a proved
nonnegative lower bound and each $v_I^+$ by a proved upper bound.
In particular, an implementation must round the positive sums upwards
in $V_I^+$ and downwards in $V_I^-$. The subtracted exponential
excesses require the opposite directions. The logarithms are then
bounded in their respective directions. A search for $t$ supplies a
candidate; it does not certify any of these inequalities. Directed
rational evaluation will be specified with the probability certificate.

\subsubsection{From quota laws to degree targets}

The arithmetic law needed on the position side is the usable quota
law $\pi(n)$ itself. This avoids transporting it first to a product
law and then back to finite populations.

\begin{proposition}[Exact reference laws and local thresholds]
\label{prop:finite-quota-thresholds}
Suppose the retained position side has size $R(n)\geq n\rho>0$.
Let its tag law be $\nu(n)$ and assume that the construction gives
\[
 \operatorname{tv}(\beta(n),R(n)\pi(n))\leq n\varepsilon,
 \qquad \beta_I(n)=R(n)\nu_I(n).
\]
Put $v=\varepsilon/\rho$, and define
\begin{equation}\label{eq:finite-compatible-fractions}
 c_0=1,\qquad
 c_J=\min_{\ell\leq n\leq u}
            \sum_{I:\,I\cap J=\varnothing}\pi_I(n)-v
                    \quad(J\ne0).
\end{equation}
Then every odd label of type $J$ has at least $c_JR(n)$ neighbors
in the initial matching graph. For a target $t_J<c_J$, set
$\theta_J=\rho(c_J-t_J)>0$. The event
\begin{equation}\label{eq:finite-degree-failure-inclusion}
 \{d_\Gamma(y)\leq t_JR(n)\}
       \subseteq\{K(J,P(y))/n\geq\theta_J\}
\end{equation}
is valid for each label $y$ of that type.
The same construction supplies the position caps
\[
 \kappa_I=\max_{\ell\leq n\leq u}\pi_I(n)
                      +v\boldsymbol1_{I=F}
\]
for the Hall estimates after transport to the full type.
\end{proposition}

\begin{proof}
Divide the count variation bound by $R(n)\geq n\rho$ to get
$\|\nu(n)-\pi(n)\|_{\mathrm{TV}}\leq v$.
For $J\ne0$, apply this bound to the indicator of disjointness
from $J$ and then take the minimum over orders. An empty-type label
is compatible with every tag, giving the exact value $c_0=1$.
If its final degree fails the target, at least
$(c_J-t_J)R(n)\geq n\theta_J$ initially incident edges were
deleted. The retained-tag charging bound dominates these deletions
by $K(J,P(y))$, proving \eqref{eq:finite-degree-failure-inclusion}.
The position caps are Corollary~\ref{cor:finite-coordinate-caps}
with the same reference law and variation bound.
\end{proof}

For an original-forest assignment on all $R(n)$ positions,
Theorem~\ref{thm:retained-donor-assignment} permits any $\varepsilon$
with $n\varepsilon\geq g/4+\eta$, where $\eta$ bounds the changed tags.
For the connected construction the corresponding count bound is
\eqref{eq:height-retained-tag-budget}. If isolated reserve positions
are added afterwards, their empty tags change the count vector;
that contribution must also be included before applying the proposition.
The arithmetic law does not remove a combinatorial tag error.

The exact odd counts give the label caps in
Corollary~\ref{cor:finite-coordinate-caps}. For full-type labels,
Lemma~\ref{lem:protected-full-labels} applies with $R_0=\ell\rho$
and $c=c_F$. Thus the fixed endpoint set
\[
 \mathcal E=\{y\in[u]:y\text{ odd},\ T(y)=F,\
      D B_0(m_u,P(y))\geq\ell\rho(c_F-\tau)\}
\]
contains every label not protected at the degree target $\tau$,
provided $c_F>\tau$. Its size bounds the actual exception population;
the weaker degree target for those exceptions still requires a local
probability estimate. This separates the deterministic exception test
from the probabilistic estimate on the remaining events.

Every quantity introduced here is obtained from a finite prefix count,
a prime support, a class histogram, or an extremum over the stated
integer interval. These formulas connect the constructions to the
probability and Hall criteria. Proving a new interval of orders still
requires their simultaneous acceptance, the supply inequalities in the
preparation, and complete coverage of the residual degrees.

\subsection{Rational evaluation of the probability bounds}
\label{sec:directed-probability}

The estimates obtained so far involve exponentials and logarithms, but
their acceptance can be reduced to integer comparisons. We give the
rounding rules and prove their directions before assembling the failure
budget. The parameters used to choose a bound are distinct from the
arithmetic used to verify it. Any permitted parameter gives a valid
bound; a numerical search need not find the best one.

For a positive integer $Q$ and a rational number $x$, write
\[
 \lfloor x\rfloor_Q=\frac{\floor{Qx}}Q,\qquad
 \lceil x\rceil_Q=\frac{\ceil{Qx}}Q.
\]
All integer floors and ceilings are evaluated exactly. A nonnegative
mean may be rounded upwards when it is used as an upper bound.
A positive decay rate must instead be rounded downwards.

\subsubsection{Positive exponentials and logarithms}

Fix integers $Q\geq1$, $L\geq1$, and $N\geq0$. For an integer
$j\geq0$ with $j<L(N+2)$, start with $a_0^-=a_0^+=Q$ and put
\[
 a_k^-=\floor*{\frac{j a_{k-1}^-}{Lk}},\qquad
 a_k^+=\ceil*{\frac{j a_{k-1}^+}{Lk}}
                         \quad(1\leq k\leq N+1).
\]
Define
\begin{equation}\label{eq:directed-exp-grid}
 E_j^- =\frac1Q\sum_{k=0}^N a_k^-,\qquad
 E_j^+ =\frac1Q\left(\sum_{k=0}^N a_k^+
       +\ceil*{\frac{a_{N+1}^+L(N+2)}{L(N+2)-j}}\right).
\end{equation}

\begin{lemma}[Exponential enclosures]\label{lem:directed-exponential}
One has $E_j^-\leq e^{j/L}\leq E_j^+$. If $x\geq0$ is rational,
$j_- =\floor{Lx}$, $j_+=\ceil{Lx}$, and $j_+<L(N+2)$, then
\begin{equation}\label{eq:directed-exp-bracket}
 E_{j_-}^-\leq e^x\leq E_{j_+}^+.
\end{equation}
\end{lemma}

\begin{proof}
Induction gives
$a_k^-\leq Q(j/L)^k/k!\leq a_k^+$. Summing the lower bounds
through $N$ and omitting the positive tail gives the first lower bound.
For the upper tail, its first term is at most $a_{N+1}^+/Q$.
Every subsequent ratio is at most $j/(L(N+2))<1$.
Summing the resulting geometric series and rounding upwards proves
the upper bound in \eqref{eq:directed-exp-grid}. Finally,
$j_-/L\leq x\leq j_+/L$ and the exponential is increasing.
\end{proof}

For logarithms, use
\begin{equation}\label{eq:atanh-log-series}
 \log\frac{1+z}{1-z}=2\sum_{i=0}^{\infty}\frac{z^{2i+1}}{2i+1}
                                      \quad(0\leq z<1).
\end{equation}
This identity follows by integrating the geometric series for
$2/(1-z^2)$ on any compact subinterval of $[0,1)$.
We only need it for $z\leq1/3$.

For $Q\geq4$, $N\geq1$, and rational $0\leq z\leq1/3$, set
\[
 q^- =\floor{Qz},\quad q^+=\ceil{Qz},\qquad
 s^- =\floor{(q^-)^2/Q},\quad s^+=\ceil{(q^+)^2/Q}.
\]
Let $t_0^-=q^-$, $t_0^+=q^+$, and, for $0\leq i<N$, put
\[
 t_{i+1}^-=\floor{t_i^-s^-/Q},\qquad
 t_{i+1}^+=\ceil{t_i^+s^+/Q}.
\]
The two logarithmic series bounds are
\begin{align}
 A_Q^-(z)&=\frac2Q\sum_{i=0}^{N-1}
                           \floor*{\frac{t_i^-}{2i+1}},
                              \label{eq:directed-log-series-lower}\\
 A_Q^+(z)&=\frac1Q\left(2\sum_{i=0}^{N-1}
                           \ceil*{\frac{t_i^+}{2i+1}}
       +\ceil*{\frac{2t_N^+Q}{(2N+1)(Q-s^+)}}\right).
                              \label{eq:directed-log-series-upper}
\end{align}
The dependence on the fixed truncation length $N$ is suppressed.

\begin{lemma}[Logarithmic enclosures]\label{lem:directed-logarithm}
The denominator $Q-s^+$ is positive, and
\[
 A_Q^-(z)\leq\log\frac{1+z}{1-z}\leq A_Q^+(z).
\]
For rational $x\geq1$, choose the unique integer $k\geq0$ with
$2^k\leq x<2^{k+1}$ and set $z=(x-2^k)/(x+2^k)$. Then
\begin{equation}\label{eq:directed-log-brackets}
 \mathcal L^-(x):=kA_Q^-(1/3)+A_Q^-(z)
 \leq\log x\leq
 kA_Q^+(1/3)+A_Q^+(z)=:\mathcal L^+(x).
\end{equation}
\end{lemma}

\begin{proof}
The bounds $q^+\leq Q/3+1$ and
$s^+\leq(q^+)^2/Q+1$ give
\[
 s^+\leq Q/9+5/3+1/Q<Q\qquad(Q\geq4).
\]
Every recurrence is monotone on nonnegative inputs. Hence
$t_i^-\leq Qz^{2i+1}\leq t_i^+$ by induction.
The lower partial sum proves the lower estimate. The upper
partial sum bounds its first $N$ terms. For the remaining terms,
replace all denominators by $2N+1$ and all ratios by $s^+/Q$;
the geometric tail is at most
$2t_N^+/((2N+1)(Q-s^+))$. This proves the upper estimate.

For the last assertion, write $x=2^k y$, where $1\leq y<2$.
Then $z=(y-1)/(y+1)\in[0,1/3)$ and
$\log x=k\log2+\log y$. Apply the series bounds to both terms.
Since $k\geq0$, the directions are preserved. The integer $k$
is determined by exact comparisons with powers of two.
\end{proof}

Both lower bounds are nonnegative; the exponential lower bound is at
least one. These properties matter when they are subtracted from a
cost or used as logarithm arguments.

\subsubsection{From a decay rate to a rational probability}

For rational $r\geq0$, a positive integer $S$, and an integer $m\geq0$,
set $b_0=S$, $b_i=\floor{r b_{i-1}/i}$ for $1\leq i\leq m$, and
\begin{equation}\label{eq:inverse-exp-rational}
 B_{S,m}(r)=\frac1S\ceil*{\frac{S^2}{\sum_{i=0}^m b_i}}.
\end{equation}

\begin{lemma}[Rounding the rate downwards]\label{lem:directed-negative-exp}
For every such $r$, $e^{-r}\leq B_{S,m}(r)\leq1$.
Fix an integer $R\geq1$ and a rational cap $c>0$ with $Rc$ integral.
For a rational lower bound $a$ on a decay rate, define
\begin{equation}\label{eq:directed-probability-map}
 \mathcal P(a)=
 \begin{cases}
  1,&a\leq0,\\
  B_{S,m(r)}(r),&a>0,
 \end{cases}
 \quad
 r=\frac{\min\{Rc,\floor{Ra}\}}R,
\end{equation}
where $m(r)$ is any nonnegative integer.
If an event has probability at most $e^{-a}$ for $a>0$, then
$\mathcal P(a)$ bounds its probability in all cases.
\end{lemma}

\begin{proof}
Induction shows $b_i\leq Sr^i/i!$, so
$S\leq\sum_i b_i\leq Se^r$. Taking reciprocals and then rounding
up proves \eqref{eq:inverse-exp-rational} is at least $e^{-r}$
and at most one. In \eqref{eq:directed-probability-map},
$0\leq r\leq a$ whenever $a>0$. Thus $e^{-a}\leq e^{-r}$.
When $a\leq0$, the trivial probability bound is valid.
\end{proof}

Capping a large rate weakens the bound and is harmless. In contrast,
a small probability cannot be replaced by zero unless the event is
proved impossible. The only zero probabilities used below come from
an empty sample, a zero mean, or a threshold outside the support.

Table~\ref{tab:directed-constants} fixes the arithmetic constants.
Every positive exponential argument is checked to lie in $[0,10]$.
For the reciprocal calculation, rate capping instead gives a Taylor
argument in $[0,64]$.

\begin{table}[ht]
\centering
\caption{Integer constants for directed probability evaluation.
The inequalities above justify the rounding; the number of digits
alone does not do so.}
\label{tab:directed-constants}
\par\vspace{.5\baselineskip}
\renewcommand{\arraystretch}{1.2}
\begin{tabular}{ll}
\hline
Quantity&Choice\\
\hline
Scale for positive exponentials and logarithms&$Q=10^{80}$\\
Series length&$N=128$\\
Positive exponential grid&$L=2^{16}$\\
Scale for reciprocal exponentials&$S=10^{150}$\\
Rate grid and cap&$R=2^{20}$, $c=64$\\
Reciprocal Taylor length&$m(r)=4\ceil r+64$\\
Scale for local means and support bounds&$10^9$\\
Scale for final binomial bounds&$10^{80}$\\
\hline
\end{tabular}
\end{table}

\subsubsection{A moment bound with a finite sample cap}

The local bounds can use the number of half-labels that may conflict,
as well as their total expected degree. The next estimate allows a
real upper bound for that number, which is useful after normalization.

\begin{lemma}[A local moment with a support bound]\label{lem:bounded-local-moment}
Let $K$ be a sum of independent class samples as in
Lemma~\ref{lem:class-convex-comparison}, with entries in $[0,D]$,
where $D>0$. Write $b$ for the total number of selected entries
and $m=\E K$. If $b>0$, then for every $\lambda>0$,
\begin{equation}\label{eq:bounded-sample-moment}
 \log\E e^{\lambda K/D}
      \leq b\log\left(1+(e^\lambda-1)\frac{m}{Db}\right).
\end{equation}
Suppose $n\geq\ell>0$, $Db/n\leq s$, and $m/n\leq\mu$,
where $s>0$ and $0\leq\mu\leq s$. For a threshold $\theta>0$,
put
\begin{equation}\label{eq:bounded-local-rate}
 a=\frac{\ell}{D}\left(\lambda\theta
              -s\log\left(1+(e^\lambda-1)\frac\mu s\right)\right).
\end{equation}
If $a>0$, then $\Pr\{K/n\geq\theta\}\leq e^{-a}$.
One may replace any upper bound for $m/n$ by its minimum with $s$.
If $b=0$ or $m=0$, positive thresholds have probability zero.
\end{lemma}

\begin{proof}
Use Lemma~\ref{lem:class-convex-comparison} with
$\phi(x)=e^{\lambda x}$. Its comparison binomial variable has
parameters $b,m/(Db)$, giving \eqref{eq:bounded-sample-moment}.
For $x\geq0$, the function $v\log(1+x/v)$ is increasing in
$v>0$, as proved in Proposition~\ref{prop:forest-local-moment}.
It is also increasing in $x$. Since $Db/n\leq s$ and $m/n\leq\mu$,
\[
 \log\E e^{\lambda K/D}
       \leq\frac nD\,s\log\left(1+(e^\lambda-1)\frac\mu s\right).
\]
Apply Markov's inequality at the threshold $n\theta$, then use
$n\geq\ell$ when the parenthesized rate is positive.
Finally $0\leq K\leq Db\leq ns$, so $m/n\leq s$ automatically.
A nonnegative random variable of mean zero vanishes surely.
\end{proof}

A bound on the support of an arbitrary random variable would not by
itself justify \eqref{eq:bounded-sample-moment}. Here the bound $s$
comes from the number of entries sampled from fixed class lists.
For the arithmetic data of Section~\ref{sec:finite-arithmetic}, one
may take $s=DB_P/\ell$, with $B_P$ from
\eqref{eq:finite-selected-population}.

\subsubsection{Directed evaluation of the three moment estimates}

In the formulas below, write $\mathcal E^-(x),\mathcal E^+(x)$
for the exponential brackets in \eqref{eq:directed-exp-bracket},
and use $\mathcal L^-,\mathcal L^+$ from
\eqref{eq:directed-log-brackets}. All supplied means, quotas,
frequency bounds and degree bounds must already have been justified.
The following proposition concerns their numerical evaluation.

\begin{proposition}[Rational upper bounds for failure]
\label{prop:directed-moment-evaluation}
The following three constructions give rational upper bounds for
the indicated failure probabilities.

\emph{(i) Local samples.}
Under Lemma~\ref{lem:bounded-local-moment}, choose rational
$0<\lambda\leq10$. If $s>0$, put
\begin{equation}\label{eq:directed-bounded-rate}
 \widehat a=\frac\ell D\left(\lambda\theta
       -s\mathcal L^+\left(1+(\mathcal E^+(\lambda)-1)\frac\mu s\right)
                             \right).
\end{equation}
Then $\mathcal P(\widehat a)$ bounds
$\Pr\{K/n\geq\theta\}$. The bound zero is also valid if
$\mu=0$ or $\theta>s$, provided $\theta>0$.

\emph{(ii) Forest excess degrees.}
Under Corollary~\ref{cor:uniform-forest-tail}, assume $D>1$,
$\beta>0$, and choose rational $\lambda>0$ with $D\lambda\leq10$.
For the upper mean and excess bounds $\bar\mu,\bar q$, put
\begin{align}
 \widehat c_\lambda
   &=\max\left\{0,
       \frac{\mathcal E^+(D\lambda)-D\mathcal E^-(\lambda)+D-1}{D-1}
                     \right\},\label{eq:directed-forest-coefficient}\\
 \widehat a
   &=\ell\left(\lambda\theta
      -\beta\mathcal L^+\left(1+
          \frac{(\mathcal E^+(\lambda)-1)\bar\mu
                       +\widehat c_\lambda\bar q}{\beta}\right)\right).
                     \label{eq:directed-forest-rate}
\end{align}
Again $\mathcal P(\widehat a)$ is an upper bound for the local
failure probability.

\emph{(iii) Global deletions.}
In Corollary~\ref{cor:finite-global-tail}, suppose rational numbers
$0\leq\widehat v_I^-\leq v_I^-$ and
$\widehat v_I^+\geq v_I^+$ replace the logarithmic endpoints.
Form $\widehat{\mathcal M}$ by the same formula for the chosen
construction, and put
\begin{equation}\label{eq:directed-global-rate}
 \widehat a=\frac\ell D(tb-\widehat{\mathcal M}).
\end{equation}
For $b>0$, the value $\mathcal P(\widehat a)$ bounds
$\Pr\{Z_{\mathcal J}\geq n^2b\}$.
\end{proposition}

\begin{proof}
In (i), replacing the exponential and logarithm by upper bounds
increases the subtracted moment. Thus $\widehat a\leq a$ in
\eqref{eq:bounded-local-rate}. A positive $\widehat a$ permits
that lemma and Lemma~\ref{lem:directed-negative-exp}; otherwise
use one. The zero bounds follow from the stated deterministic
conditions.

In (ii), the coefficient of $e^\lambda$ in $c_\lambda$ is
negative. Using its lower bound and the upper bound for
$e^{D\lambda}$ therefore gives $\widehat c_\lambda\geq c_\lambda$.
Both coefficients multiplying $\bar\mu,\bar q$ are nonnegative.
Monotonicity of the logarithm shows that the displayed moment is
an upper bound for the one in \eqref{eq:uniform-forest-rate}.
Apply the same argument to its rate.

For (iii), increasing the upper class endpoints and decreasing the
lower ones can only increase the average, oscillation and donor
terms, and decrease the subtracted completion credit. Consequently
$\widehat{\mathcal M}$ bounds the normalized logarithmic moment
from above. Corollary~\ref{cor:finite-global-tail} and
Lemma~\ref{lem:directed-negative-exp} give the result.
\end{proof}

For clarity, the directed class endpoints in (iii) can themselves be
computed without any unspecified transcendental operation. Retain the
notation $Q_{n,I}$, $\kappa$, $r_-,r_+$ and $a_I^-,a_I^+$ of
Lemma~\ref{lem:finite-log-moments}, and form
\begin{align*}
 \widehat V_I^+
   &=\sum_{\substack{j\leq m_u\\T(j)=I}}
            (\mathcal E^+(tQ_{u,I}(j)/\ell)-1)
       -\kappa r_-(\mathcal E^-(t/\ell)-1)\boldsymbol1_{I=0},\\
 \widehat V_I^-
   &=\left(\sum_{\substack{j\leq m_\ell\\T(j)=I}}
            (\mathcal E^-(tQ_{\ell,I}(j)/u)-1)
       -\kappa r_+(\mathcal E^+(t/u)-1)\boldsymbol1_{I=0}\right)_+.
\end{align*}
Provided every displayed exponential argument is in $[0,10]$, one
has $\widehat V_I^+\geq V_I^+\geq0$ and
$0\leq\widehat V_I^-\leq V_I^-$. Hence the choices
\begin{equation}\label{eq:directed-class-log-endpoints}
 \widehat v_I^- =\mathcal L^-(1+\widehat V_I^-/a_I^+),\qquad
 \widehat v_I^+ =\mathcal L^+(1+\widehat V_I^+/a_I^-)
\end{equation}
satisfy the proposition. In particular, the subtracted terms use the
opposite rounding direction from the positive sums. Every denominator
is positive by the population hypotheses of Section~\ref{sec:finite-arithmetic}.

For (i), any positive grid value $\lambda=j/L\leq10$ suffices.
When $0<\mu<\theta<s$, differentiating the unrounded rate gives
the stationary value
\[
 \lambda_* =\log\frac{\theta(s-\mu)}{\mu(s-\theta)}.
\]
An approximation to this value may select a grid point; the final
bound always uses \eqref{eq:directed-bounded-rate} at that exact
point. At $\theta=s$ one may choose the largest grid point.
No optimality assertion is needed. Likewise, a floating-point search
among forest parameters only selects a rational $\lambda$; its
reported rate must be recomputed by \eqref{eq:directed-forest-rate}.

\subsubsection{Integer thresholds, binomial bounds, and baselines}

Let $K_n$ be integer-valued, with $\ell\leq n\leq u$, and let
$\theta>0$. Proposition~\ref{prop:integer-local-tail} permits the
replacement
\[
 \theta_* =\min_{\ell\leq n\leq u}\frac{\ceil{n\theta}}n
                   \quad\text{in the event }K_n/n\geq\theta.
\]
This finite minimum uses rational comparisons. For a mean bound
$\E K_n/n\leq\mu$, put $\widehat\mu=\lceil\mu\rceil_{10^9}$.
If at most $B\geq1$ entries from $[0,D]$ are sampled, define
\[
 p=\min\left\{1,\frac{u\widehat\mu}{DB}\right\},\qquad
 t_0=\frac{\ceil{\ell\theta_*}}D.
\]
An exact rational upper bound is
$U(B,p,t_0)$ from \eqref{eq:binomial-positive-part-bound}; its
upward rounding to denominator $10^{80}$ remains valid.
For $p=a/v$, each mass is the integer
$\binom Bk a^k(v-a)^{B-k}$ divided by $v^B$. Summing these masses
and first moments above an integer $h<t_0$ computes
$\E(Z-h)_+/(t_0-h)$ exactly. The minimum includes the bound one.
An empty sample or zero mean gives zero at a positive threshold;
a nonpositive threshold is assigned the bound one.

The binomial estimate and (i) of
Proposition~\ref{prop:directed-moment-evaluation} bound the same
local failure event, so their minimum is valid. For the connected
construction, a further estimate comes from
\eqref{eq:finite-baseline-data}. Use the completed mean, subtract
$b_-$, replace the degree bound by $D-1$, and the threshold by
$\theta_*-b_+$. The selected population bound remains $B$.
The event inclusion is
\[
 \{K_n^G/n\geq\theta\}
 \subseteq\{K_n^G/n\geq\theta_*\}
 \subseteq\{X_n/n\geq\theta_*-b_+\}.
\]
For $D=1$ the excess variable vanishes; for a nonpositive shifted
threshold this route supplies only the bound one. The mean used for
$X_n$ does not receive a completion credit. That credit belongs to
the separate estimate using the original degrees.

One may also take the minimum with (ii), or omit that calculation
when another bound is already small. Omitting an estimate cannot
invalidate the remaining ones. Replacing a small positive estimate
by zero would be a different operation and requires a separate proof.

\subsubsection{Exact acceptance of a simultaneous event}

The final step is an integer comparison, but it is meaningful only
after the events and their multiplicities have been specified.

\begin{proposition}[A rational failure budget]\label{prop:rational-failure-budget}
Fix an order interval and one construction whose deterministic
hypotheses have been proved uniformly on that interval. At each
order let $\{\mathcal F_e:e\in E\}$ be a finite family of local
failure events that covers all required local degree failures.
Suppose rational numbers $p_e\in[0,1]$ bound their probabilities,
and rational numbers $q_1,\ldots,q_h\in[0,1]$ bound the remaining
global failure events. If avoidance of these events implies the
chosen Hall criterion, then
\begin{equation}\label{eq:rational-failure-budget}
 \sum_{e\in E}p_e+\sum_{j=1}^h q_j<1
\end{equation}
guarantees a perfect matching at every order of the interval.

If events with the same certified bound are grouped, their contribution
is $m p$, where $m$ is the number of events in that group. Events may
be removed from the family only when their failure is already excluded,
or their inclusion in a retained event has been proved. If all bounds
have a common denominator $S$, condition
\eqref{eq:rational-failure-budget} is exactly the comparison of their
summed integer numerators with $S$.
\end{proposition}

\begin{proof}
Fix an actual order and an actual input forest. All construction choices
are made before the class permutations. The union bound makes the
probability of any listed failure at most the left side of
\eqref{eq:rational-failure-budget}. If this is less than one,
an outcome avoiding all failures exists. By hypothesis the matching
criterion applies to that outcome. The bounds hold separately for
each input and order, so this proves the assertion for all of them;
there is no union bound over input forests or over orders.

Grouping is an exact regrouping of a finite sum and preserves its
multiplicities. Removing an event requires one of the stated reasons,
so the remaining union still covers every relevant failure. Multiplying
a strict rational inequality by its positive common denominator gives
the stated integer comparison.
\end{proof}

In our application, the local family consists of maximal prime-support
events with a common threshold within each type and degree target,
as in Lemma~\ref{lem:maximal-local-supports}. A maximal event may
represent many labels; conversely, distinct maximal events with equal
arithmetic data still contribute separately. The three global families
are those in Table~\ref{tab:hall-ranges}. A Hall range is omitted only
when it is empty, its capped domain is infeasible, or a proved
deterministic estimate excludes all its obstructions. Otherwise its
positive deletion budget and failure estimate must both be checked.
For one global event, the minimum of its Markov and exponential bounds
may be taken before summing the different events.

The choices in Table~\ref{tab:directed-constants} make every exponential
failure bound an integer multiple of $1/S$. The rounded binomial bound
also has that property because $10^{80}$ divides $S=10^{150}$.
An exact Markov bound may instead be kept as a rational number, in
which case integer cross-multiplication or a common multiple of the
denominators gives the same strict test. Rounding any additional
probability term upwards is permitted; rounding it downwards is not.

This establishes the numerical soundness of the probability tests.
A range certificate must still supply the correct arithmetic and
construction inputs, verify the Hall budgets and the deterministic
preparation, and cover every possible residual degree and every integer
order claimed. The next assembly uses these requirements together;
the evaluation rules alone assert no new range of prime labelings.

\subsection{Preparation and certificates on integer intervals}
\label{sec:interval-construction}

We now connect the numerical bounds to an arbitrary input tree.
The order interval and its arithmetic data are fixed first. The tree
then determines a residual degree, and this degree selects a row of
the certificate. In particular, a complete table must account for
every degree that can occur after the preparation, not merely every
degree at which a numerical search succeeded.

We use the preparation of
Proposition~\ref{prop:highest-degree-preparation}. All special even
labels are assigned to vertices of largest degree; all remaining odd
positions are retained. Thus no position reserve is introduced in this
section. The final matching is between two sets of the same size.

\subsubsection{Uniform preparation and the universal odd labels}

Fix integers $12\leq\ell\leq u$. All extrema below are over integer
orders $n\in[\ell,u]$. Retain $m_n=\floor{n/2}$ and $w_n$ from
Section~\ref{sec:finite-arithmetic}, and put
\[
 \begin{gathered}
 b_* = b(u),\qquad s_n=1+\pi(n)-\pi(\floor{n/3}),\\
 W=\min_n\frac{w_n}{n},\qquad
 \Omega=\max_n\frac{w_n}{m_n},\qquad
 S_* =\min_n\frac{s_n}{n}.
 \end{gathered}
\]
The number $s_n$ counts the special odd labels before assigning the
core. Define
\begin{equation}\label{eq:interval-preparation-constants}
 \kappa=\frac{b_*}{\ell},\qquad
 \rho=\frac12-\kappa,\qquad
 o_+=\frac{\ell+1}{2\ell},\qquad
 \upsilon=\frac{S_*-\kappa}{o_+-\kappa}.
\end{equation}
We require
\begin{equation}\label{eq:interval-preparation-tests}
 \begin{gathered}
 0<w_n<m_n\quad\text{for every }n,\qquad
 \min_n\bigl(\pi(n)-\pi(m_n)\bigr)\geq b_*,\\
 \rho>0,\qquad S_*>\kappa.
 \end{gathered}
\end{equation}
These are finite integer or rational comparisons. They imply $W>0$
and $0\leq\Omega<1$.

\begin{proposition}[Uniform preparation]\label{prop:interval-preparation}
Assume \eqref{eq:interval-preparation-tests}. Every tree of order
$n\in[\ell,u]$ can be prepared as in
Proposition~\ref{prop:highest-degree-preparation}, with
\[
 k=\abs K\leq b_*,\qquad
 M=\abs A=m_n-w_n>0,\qquad
 N=\abs B=\ceil{n/2}-k>0.
\]
Its remaining even labels are $2j$ with $j\in[m_n]\setminus W_n$,
and its remaining odd labels form $Y=O_n\setminus\Psep$, with
$\abs Y=N$. Moreover,
\begin{equation}\label{eq:interval-side-bounds}
 \frac Nn\geq\rho,\qquad
 \frac Mn\geq\frac{(1-1/\ell)(1-\Omega)}2.
\end{equation}
At least $\upsilon N$ members of $Y$ are coprime to every remaining
even label. If the block cutoff satisfies $3z<\ell$, all these
labels have empty type. They are universal in the matching graph
for every permutation within the even classes.
\end{proposition}

\begin{proof}
Since $b(n)\leq b_*$, the prime count in
\eqref{eq:interval-preparation-tests} gives the hypothesis of
Corollary~\ref{cor:even-side}. Apply that corollary and then delete
the $w_n$ highest-degree vertices of its even side. Counting vertices
and unused labels gives the stated sizes. Also
\[
 N/n\geq1/2-b_*/\ell=\rho,
 \qquad
 M/n=(m_n/n)(1-w_n/m_n)
       \geq\frac{(1-1/\ell)(1-\Omega)}2.
\]

The special odd labels that remain are $1$ and the primes greater
than $n/3$, except for the $k$ primes used on the core. A prime
$p>n/2$ has no remaining multiple. If $n/3<p\leq n/2$, its only
possible even multiple is $2p$, which was removed with the special
even labels. Thus these $s_n-k$ labels are now coprime to every
remaining even label, without using the relaxed relation. Their
types are empty because each such prime exceeds $z$.

Write $t=k/n$, $s=s_n/n$, and $o=\ceil{n/2}/n$. Then
$0\leq t\leq\kappa$, $s\geq S_*$, and $o\leq o_+$.
Since $S_*\leq o_+$ and $S_*>\kappa$, the function
$(S_*-t)/(o_+-t)$ is nonincreasing for $0\leq t\leq\kappa$.
Consequently
\[
 \frac{s_n-k}{N}=\frac{s-t}{o-t}
       \geq\frac{S_*-t}{o_+-t}
       \geq\frac{S_*-\kappa}{o_+-\kappa}=\upsilon.
\]
An empty-type label is compatible with every tag, and the coprimality
just proved is independent of the assignment within each class.
This proves the last assertion.
\end{proof}

\subsubsection{Which residual degrees must be covered}

\begin{lemma}[A complete finite degree list]\label{lem:interval-degree-cover}
In this preparation, let $D=\Delta_G(A)$. If $D>0$, then
\begin{equation}\label{eq:interval-degree-list}
 D\leq D_*:=\ceil*{\frac1{2W}}.
\end{equation}
More precisely, an order $n$ can have residual degree $D$ only if
\begin{equation}\label{eq:interval-degree-feasibility}
 E_n(D):=n-1-Dw_n\geq M_n:=m_n-w_n,
\end{equation}
or, equivalently,
\[
 D\leq 1+\floor*{\frac{n-1-m_n}{w_n}}.
\]
The case $D=0$ already admits a residual labeling.
\end{lemma}

\begin{proof}
Let $I$ count the isolated vertices of $A$. Every nonisolated vertex
of $A$ contributes at least one to $e(G)$, so $e(G)+I\geq M_n$.
The isolation budget \eqref{eq:prepared-isolation-budget} proves
\eqref{eq:interval-degree-feasibility} and its equivalent integer
form. The nonnegative excess sum in
\eqref{eq:prepared-donor-excess} also gives
\[
 w_n(D-1)\leq n-1-m_n<n/2.
\]
Thus $D<1+n/(2w_n)\leq1+1/(2W)$.
An integer strictly smaller than $1+x$ is at most $\ceil x$,
which proves \eqref{eq:interval-degree-list}. If $D=0$, the residual
forest has no edges, as already observed in
Proposition~\ref{prop:highest-degree-preparation}.
\end{proof}

A table may retain the whole list $1,\ldots,D_*$ even when
\eqref{eq:interval-degree-feasibility} rules out some individual
order--degree pairs. Such pairs require no matching argument. They
must not be passed to a donor-stock formula whose hypotheses they
violate. In particular, infeasibility is established before evaluating
that formula; it is not inferred from a failed numerical estimate.

\begin{lemma}[Completion capacity and credit]\label{lem:interval-completion}
Fix a prepared forest of residual degree $D\geq3$. The inequality
\begin{equation}\label{eq:interval-completion-test}
 D\frac{(1-1/\ell)(1-\Omega)}2\geq1-W
\end{equation}
guarantees a bipartite completing tree $Q$ on $A\sqcup B$ with
$\Delta_Q(A)\leq D$. One then has
\begin{equation}\label{eq:interval-completion-costs}
 \frac{e(Q)}n\leq1-W,\qquad
 \frac{e(Q)-e(G)}n\geq W(D-1).
\end{equation}
For the original forest, $e(G)/n\leq1-DW$.
\end{lemma}

\begin{proof}
The desired completing tree has
$M+N-1=n-k-w_n-1\leq n(1-W)$ edges.
Equations \eqref{eq:interval-side-bounds} and
\eqref{eq:interval-completion-test} imply $DM\geq M+N-1$.
Apply Lemma~\ref{lem:bipartite-completion}.

Before removing $U$, the forest $T-K$ has at most $n-k-1$ edges.
Each of the $w_n$ removed vertices has degree at least $D$ there,
and their incident edge sets are disjoint. Hence
$e(G)\leq n-k-1-w_nD$. Subtraction from
$e(Q)=n-k-w_n-1$ gives
$e(Q)-e(G)\geq w_n(D-1)$. Divide by $n$ and use $w_n/n\geq W$.
The last assertion follows from
$e(G)\leq n-1-w_nD\leq n(1-DW)$.
\end{proof}

The degree variation for a connected construction is measured in
$Q$, while the deletion events concern $G$. The second bound in
\eqref{eq:interval-completion-costs} pays for this distinction.
It is a lower bound for an actual number of added edges, not a
subtraction from the degree list used by the completed baseline.

\subsubsection{Construction data for one row}

Fix a degree $1\leq D\leq D_*$, a cutoff $3z<\ell$, and a
partition into $s$ blocks. Write $g=2^s$. Assume that all endpoint
quota lower bounds $a_I^-$ in
\eqref{eq:finite-population-endpoints} are positive. The construction
row will supply numbers $e_+,H_+,\delta,q_-,\varepsilon_B$ with
\begin{equation}\label{eq:interval-construction-interface}
 \begin{gathered}
 e(\mathcal G)/n\leq e_+,\qquad H/n\leq H_+,\qquad
 \Delta/n\leq\delta,\\
 (e(Q)-e(G))/n\geq q_-\quad\text{in the connected case},\\
 \operatorname{tv}(\beta,N\pi(n))/n\leq\varepsilon_B.
 \end{gathered}
\end{equation}
Here $\mathcal G=G$ for an original-forest construction and
$\mathcal G=Q$ for a connected construction. In the former,
$H$ is the initial signed degree variation before donor exchanges
and $\Delta$ is the ensuing nonnegative transfer. In the latter,
$H$ is measured in $Q$ and $\Delta=0$. Set $q_-=0$ when no
completion is used.

For the original construction with cuts confined to $B$, put
$d=3(g-1)$ and define the rational bound
\begin{equation}\label{eq:interval-side-series}
 c_B=\begin{cases}
 d/\ell,&D=1,\\
 \displaystyle\sum_{j=0}^{127}
   \min\left\{d/\ell,\frac{r^j}{1+r^j}\right\}
       +\frac{r^{128}}{1-r},&D\geq2,\quad r=(D-1)/D.
 \end{cases}
\end{equation}
For the ordinary original-forest construction, let
$\mathcal I_D=\{n\in[\ell,u]:E_n(D)\geq M_n\}$.
Suppose this set is nonempty and $E_n(D)\geq0$ throughout the
interval. Put
\begin{align}
 c_*&=\max_{\ell\leq n\leq u}
       \mathcal C(n-w_n,E_n(D),3(g-1)),
                                      \label{eq:interval-donor-cuts}\\
 z_*&=\max_{n\in\mathcal I_D}
       \left(c_*+\sum_{t=1}^{D-1}(c_*-\ell_t(n))_+\right),
                                      \label{eq:interval-donor-boundary}\\
 \delta_*&=\min\left\{\frac{Dc_*}{\ell},\
                \frac{c_*}{\ell}+\frac12-W(D-1)\right\}.
                                      \label{eq:interval-donor-transfer}
\end{align}
Here $\ell_t(n)$ is the paid stock in \eqref{eq:paid-stock}, with
\[
 (M,a,E,h)=(M_n,a_0(n),E_n(D),D).
\]
For $D=1$, set $\delta_*=0$.
Taking the maximum in \eqref{eq:interval-donor-cuts} over extra
orders only weakens the cut bound; the stock maximum is restricted
to the orders where its premise holds.

For the connected construction choose either height sequence $R$
from Section~\ref{sec:height-assignments} and put
\begin{equation}\label{eq:interval-connected-cuts}
 M_*=\floor{u/2}-\min_nw_n,\qquad
 c_H=\mathcal C_R(M_*,2(g-1)).
\end{equation}
The choice of the sequence is part of the row, not an assumption
about the input tree.

\begin{proposition}[Three admissible constructions]
\label{prop:interval-construction-bounds}
The following choices satisfy
\eqref{eq:interval-construction-interface} for every feasible order
in the row. All positions are retained, and every original edge
has $I(a)\subseteq J(b)$.

\emph{(i) Original forest, cuts confined to $B$.}
One may take
\[
 e_+=1-DW,\quad H_+=\frac{gD}{4\ell},\quad
 \delta=q_-=0,\quad \varepsilon_B=\frac{g}{4\ell}+c_B.
\]

\emph{(ii) Original forest, with donor exchanges.}
Use the cut, boundary and transfer bounds in
\eqref{eq:interval-donor-cuts}--\eqref{eq:interval-donor-transfer}.
Then one may take
\[
 e_+=1-DW,\quad H_+=\frac{gD}{4\ell},\quad
 \delta=\delta_*,\quad q_-=0,\quad
 \varepsilon_B=\frac{g}{4\ell}+\frac{z_*}{\ell}.
\]

\emph{(iii) Connected completion.}
Assume $D\geq3$ and \eqref{eq:interval-completion-test}. One may take
\[
 \begin{gathered}
 e_+=1-W,\qquad H_+=\frac{g(D-1)}{4\ell},\qquad
 \delta=0,\qquad q_-=W(D-1),\\
 \varepsilon_B=
    \frac{\mathcal R_g(g(D-1)/4)+c_H}{\ell}.
 \end{gathered}
\]
Outside a set $C\subseteq B$ of at most $c_H$ vertices, every
neighbor in $Q$ has exactly its position's tag. The tags on $C$
are full.
\end{proposition}

\begin{proof}
For (i), use Theorem~\ref{thm:original-type-assignment} with cuts
confined to $B$, and retain those positions.
The normalized separator estimate following that theorem is at most
\eqref{eq:interval-side-series}: for the omitted terms use
$r^j/(1+r^j)\leq r^j$ and sum the geometric tail.
Theorem~\ref{thm:retained-donor-assignment} then gives the tag cost
and zero transfer. Its rounding gives $H\leq gD/4$.

For (ii), Proposition~\ref{prop:highest-degree-preparation} gives
at most $n-w_n$ vertices, edge bound $E_n(D)$, and the isolation
bound needed for every paid stock. Thus \eqref{eq:interval-donor-cuts}
bounds the cuts and \eqref{eq:interval-donor-boundary} bounds the
changed tags. The degree excess is less than
$n/2-w_n(D-1)$. The transfer estimate of
Theorem~\ref{thm:retained-donor-assignment}, divided by $n\geq\ell$,
is therefore bounded by \eqref{eq:interval-donor-transfer}.
The same theorem gives the remaining bounds.

For (iii), use Lemma~\ref{lem:interval-completion} and
Corollaries~\ref{cor:height-retained-assignment}
and~\ref{cor:height-interval-budget}. They give the stated height,
rounding and tag costs, including the stronger tag property.
Lemma~\ref{lem:interval-completion} gives both the edge bound and
the completion credit. All three constructions are fixed before
the class permutations, as required by the probability estimates.
\end{proof}

Each upper bound in the proposition may be increased, and the
completion credit may be decreased. If a row uses a larger degree
bound $H_+$ in the connected tag calculation, the value
$\mathcal R_g(\ell H_+)/\ell$ is still valid: the defining maximum
for $\mathcal R_g$ is nondecreasing. The rounding directions are
therefore the ones used in Section~\ref{sec:directed-probability}.
There is no donor credit in case (iii) and no completion credit in
cases (i) and (ii).

\subsubsection{Coordinate bounds and the first Hall range}

Set $v=\varepsilon_B/\rho$ and use the exact laws
$\pi_I(n)=a_I(n)/M_n$. Define
\begin{align}
 c_0&=1,\qquad
 c_J=\min_n\sum_{I\cap J=\varnothing}\pi_I(n)-v
                      \quad(J\ne0),\label{eq:interval-initial-degrees}\\
 \kappa_I&=\max_n\pi_I(n)+v\boldsymbol1_{I=F},\qquad
 p_I=\frac{\bar C_I(u)-\bar C_I(\floor{u/2})}{\ell\rho}.
                                      \label{eq:interval-hall-caps}
\end{align}

\begin{corollary}[Uniform caps and a vacuous range]
\label{cor:interval-caps}
The initial degree of a label of type $J$ is at least $c_JN$.
The position caps $\kappa_I$ are valid after full-type transport,
and $p_I$ bounds the normalized number of remaining odd labels of
type $I$. At least $\upsilon N$ labels have empty type and are
universal. If the local degree conditions of
Lemma~\ref{lem:three-hall-ranges} hold and
\begin{equation}\label{eq:interval-full-vacuity}
 p_F\leq\tau,
\end{equation}
then no deficient set belongs to the full-only range, for any
$\eta\geq\tau$.
\end{corollary}

\begin{proof}
The tag variation, divided by $N\geq n\rho$, is at most $v$.
Apply Proposition~\ref{prop:finite-quota-thresholds} and
Corollary~\ref{cor:finite-coordinate-caps}. The universal labels
were established in Proposition~\ref{prop:interval-preparation}.
Finally, every deficient set has size greater than $\tau N$,
whereas a set containing only full-type labels has size at most
$p_FN\leq\tau N$. These statements contradict each other.
\end{proof}

The strict numerical margin $p_F<\tau$ is convenient in a table,
but the proof only needs the weak inequality in
\eqref{eq:interval-full-vacuity}. This removes one global deletion
event. It does not remove the local degree estimates for full-type
labels.

\subsubsection{Soundness of an interval certificate}

We can now state the remaining acceptance conditions without referring
to how the parameters were found. All type assignments and tags in the
next theorem are constructed by the preceding propositions.

\begin{theorem}[Certificate criterion on an integer interval]
\label{thm:integer-interval-certificate}
Suppose \eqref{eq:interval-preparation-tests} holds. For each
$1\leq D\leq D_*$ that can occur, choose an admissible construction
row as in Proposition~\ref{prop:interval-construction-bounds}, with
positive quota lower bounds. Form $c_J,\kappa_I,p_I$ as above.
For each such row suppose the following conditions are verified
uniformly for its feasible orders.

\emph{(a) Degree targets and exceptions.}
Choose $0\leq\epsilon<\tau\leq\eta\leq\gamma$ with
$\gamma+\upsilon\leq1$. Specify a set $\mathcal E$ of full-type odd labels
in $[u]$, of size at most $\ell\rho\epsilon$. Every remaining
full-type label outside $\mathcal E$ has degree greater than $\tau N$
for every class permutation. This last assertion may be proved by
Lemma~\ref{lem:protected-full-labels}.

\emph{(b) Simultaneous local estimates.}
Set $t_0=\gamma$, $t_F=\epsilon$, and $t_J=\eta$ for the other
types, and require $\theta_J=\rho(c_J-t_J)>0$ for every local
family used. A finite collection of prime-support events covers
all local failures
\[
 K(J,P(y))\geq n\theta_J.
\]
Its rational probability bounds, with their event multiplicities,
have sum at most $p_{\rm loc}$.

\emph{(c) Every Hall range.}
For each of the three ranges in Table~\ref{tab:hall-ranges}, either
its emptiness is proved, its capped domain is infeasible, or all
its obstructions are excluded by a proved deterministic estimate.
Otherwise its capped minimum $\beta_k$ is positive, and a rational
number $p_k$ bounds the probability
\begin{equation}\label{eq:interval-global-failure}
 \Pr\{Z_k\geq n^2\rho^2\beta_k\},
\end{equation}
where $Z_k$ dominates the relevant number of deleted pairs.

If, for every row,
\begin{equation}\label{eq:interval-final-acceptance}
 p_{\rm loc}+\sum_{k\ \text{not excluded}}p_k<1,
\end{equation}
then every tree of every order in $[\ell,u]$ has a prime labeling.
\end{theorem}

\begin{proof}
Fix an arbitrary tree $T$ and order $n$ in the interval. Prepare it
by Proposition~\ref{prop:interval-preparation}. If its residual
degree is zero, Proposition~\ref{prop:highest-degree-preparation}
finishes the proof. Otherwise
Lemma~\ref{lem:interval-degree-cover} supplies a row, and the actual
order satisfies that row's feasibility condition.

Use its construction and permute the usable half-labels independently
within the prescribed classes. Every pair allowed by disjoint types
satisfies the small-prime conditions, and the conflict variables bound
its possible large-prime deletions. By
Proposition~\ref{prop:finite-quota-thresholds}, failure of a required
strict label degree implies one of the events in (b). Avoiding these
events therefore gives degrees greater than $\epsilon N$ for the
exceptional labels, greater than $\eta N$ for the proper nonempty
types, and greater than $\gamma N$ for the empty type. Condition (a)
gives the stronger target for the other full-type labels. Also
\[
 \abs{\mathcal E\cap Y}\leq\ell\rho\epsilon\leq N\epsilon.
\]
The hypotheses of Lemma~\ref{lem:three-hall-ranges} are now satisfied.

An unexcluded range has, outside its event in
\eqref{eq:interval-global-failure}, fewer than
$n^2\rho^2\beta_k\leq N^2\beta_k$ relevant deletions. Theorem
\ref{thm:three-range-hall} excludes its obstructions. The other
ranges have already been excluded by (c). The exact union bound
of Proposition~\ref{prop:rational-failure-budget}, using
\eqref{eq:interval-final-acceptance}, supplies a permutation for
which all these conclusions hold simultaneously. Hence the matching
graph has a perfect matching.

Assign each matched odd label to its position. This labels the
residual forest bijectively and coprimely on every edge. Proposition
\ref{prop:highest-degree-preparation} restores the special labels
and repairs any allowed exceptional pairs, giving a prime labeling
of $T$. The argument applies separately to every input tree and
order. No simultaneous choice of permutations for different trees
is required.
\end{proof}

A collection of adjacent integer intervals is accepted by applying
this theorem on each interval. Adjacency means that each new lower
endpoint is exactly one greater than the preceding upper endpoint.
Within an interval the full degree list is
$1,\ldots,\ceil{1/(2W)}$, with only the proved feasibility exclusions
allowed. These two checks settle the order and degree quantifiers.
They do not replace any of (a)--(c) or the strict probability test.

\subsubsection{The construction checks for orders 1360 through 1521}

Table~\ref{tab:interval-degree-cover} specifies the order intervals
and the largest degree in each list. The accompanying finite data
contain the prime blocks, construction choice, height choice, and
rational construction bounds for every row. They contain no claimed
acceptance of conditions (a)--(c) as a whole.

\begin{proposition}[Verified preparation and construction data]
\label{prop:finite-interval-construction-check}
For the intervals in Table~\ref{tab:interval-degree-cover}, the
preparation tests \eqref{eq:interval-preparation-tests} hold, and
the displayed degree lists cover every possible positive residual
degree. The supplied 509 construction rows satisfy
Proposition~\ref{prop:interval-construction-bounds}, with positive
quota lower bounds and the coordinate bounds
\eqref{eq:interval-initial-degrees}--\eqref{eq:interval-hall-caps}.
There are 64 original-forest rows with cuts confined to $B$,
231 original-forest rows with donors, and 214 connected rows.
In every row $p_F<\tau$, with minimum difference $10^{-6}$.
\end{proposition}

\begin{proof}
All quantities are finite counts or rational expressions, except
for the height budgets, which are finite integer recursions.
The supplement provides their input data and a standalone checker.
We describe the checks to specify the content of the assertion.

The 32 intervals are adjacent, their union is $[1360,1521]$, and
the degree list in each is exactly
$1,\ldots,\ceil{1/(2W)}$. Direct prime enumeration computes every
$s_n,w_n$ and core supply, hence all preparation constants.
The lists contain 2,598 order--degree pairs. Of these, exactly 18
fail \eqref{eq:interval-degree-feasibility}; the remaining 2,580
pairs are used whenever a stock hypothesis is needed.

For every row, exact class counts verify $3z<\ell$ and $a_I^->0$.
The checker recomputes the one-sided sum, the ordinary cut and donor
stock bounds, or the chosen connected height budget, as appropriate.
Every supplied upper bound dominates the corresponding bound of
Proposition~\ref{prop:interval-construction-bounds}; every supplied
completion credit and lower side or universal bound has the required
opposite direction. All 214 connected rows pass the completion
capacity test. No row combines donors with a completion credit.

Finally the checker recomputes every actual quota law, all coordinate
extrema, the total tag variation, and the odd endpoint counts.
These reproduce the supplied slot caps and compatibility bounds
exactly. Subtraction of the full-type population cap from $\tau$
gives a positive rational in every row; the least is $1/10^6$.
All comparisons use integers and rational arithmetic. This proves
the stated preparation and construction assertions.
\end{proof}

\begin{table}[ht]
\centering
\caption{The complete degree cover for the next finite interval.
Each entry represents all degrees $1\leq D\leq D_*$. This table
and Proposition~\ref{prop:finite-interval-construction-check}
certify the preparation and construction data; the remaining local
and Hall acceptance conditions must still be checked.}
\label{tab:interval-degree-cover}
\par\vspace{.5\baselineskip}
\renewcommand{\arraystretch}{1.15}
\begin{tabular}{rrrr}
\hline
Order interval&$D_*$&Order interval&$D_*$\\
\hline
$1360$--$1361$&15&$1416$--$1417$&16\\
$1362$--$1363$&15&$1418$--$1421$&16\\
$1364$--$1365$&15&$1422$--$1425$&16\\
$1366$--$1369$&15&$1426$--$1429$&16\\
$1370$--$1373$&15&$1430$--$1433$&16\\
$1374$--$1377$&15&$1434$--$1441$&16\\
$1378$--$1381$&16&$1442$--$1449$&16\\
$1382$--$1385$&16&$1450$--$1457$&16\\
$1386$--$1389$&16&$1458$--$1465$&16\\
$1390$--$1393$&16&$1466$--$1473$&17\\
$1394$--$1397$&16&$1474$--$1481$&17\\
$1398$--$1401$&16&$1482$--$1489$&16\\
$1402$--$1405$&16&$1490$--$1497$&16\\
$1406$--$1409$&16&$1498$--$1505$&16\\
$1410$--$1411$&16&$1506$--$1513$&17\\
$1412$--$1415$&16&$1514$--$1521$&16\\
\hline
\end{tabular}
\end{table}

The first Hall range is therefore absent as soon as the local degree
conditions hold. The remaining application of
Theorem~\ref{thm:integer-interval-certificate} must verify the
exception populations, the covering local events, and the middle and
high Hall ranges with their simultaneous probability bound. The
present finite assertion does not extend the proved order interval
beyond Theorem~\ref{thm:initial-prefix-1359}.

\subsection{Local degree certificates}
\label{sec:local-certificates}

We verify the local conditions of the interval criterion. The possible
odd labels are first divided by type and prime support. We then specify
which full-type labels require the weaker degree target, bound one event
for each maximal support, and add the resulting probabilities. The
arithmetic bounds and the event cover are both checked on the complete
endpoint population.

Throughout this section, fix one of the construction rows of
Section~\ref{sec:interval-construction}, on an interval $[\ell,u]$.
Retain its degree $D$, cutoff $z$, blocks, and normalized construction
bounds. Write $F$ for the full type. We may replace an initial degree
bound $c_J$ by a smaller bound $\widetilde c_J$, with
$\widetilde c_0=1$. This permits several types to use a common lower
bound; it does not identify their random conflict events.

\subsubsection{The complete endpoint event cover}

For every odd $y\leq u$, define
\[
 P_u(y)=\{p:z<p\leq\floor{u/2},\ p\mid y\}.
\]
For each type $J$, let $\mathcal M_J(u)$ be the inclusion-maximal
members of
\[
 \{P_u(y):1\leq y\leq u,\ y\text{ odd},\ T(y)=J\}.
\]
An empty type family contributes no event. If its only support is
empty, the empty support is its unique maximal member.
To obtain a disjoint record of the covered labels, assign each $y$
to the first $P\in\mathcal M_{T(y)}(u)$ containing $P_u(y)$,
in lexicographic order. Denote the resulting label sets by $Y_{J,P}$.

\begin{lemma}[Endpoint coverage and multiplicities]
\label{lem:endpoint-local-cover}
The sets $Y_{J,P}$ partition all odd integers in $[u]$. For arbitrary
real thresholds $L_J$ and every actual order $n\in[\ell,u]$,
\begin{equation}\label{eq:endpoint-local-cover}
 \bigcup_{y\in Y}\{K(T(y),P_u(y))\geq L_{T(y)}\}
 \subseteq
 \bigcup_J\ \bigcup_{P\in\mathcal M_J(u)}
                 \{K(J,P)\geq L_J\},
\end{equation}
where $Y$ is the actual set of remaining odd labels.
If each event on the right has probability at most $q_{J,P}$,
the probability of the union is at most
\[
 \sum_J\sum_{P\in\mathcal M_J(u)}q_{J,P}.
\]
Events with the same numerical bound may be evaluated together,
but that bound is multiplied by their number of distinct pairs
$(J,P)$, not by the number of labels in their sets $Y_{J,P}$.
\end{lemma}

\begin{proof}
Every member of a finite family of sets is contained in a maximal
member. Thus the specified first member exists and assigns each odd
label exactly once. Since every usable half-label is at most
$\floor{u/2}$, larger primes contribute nothing to its conflict
indicator. Lemma~\ref{lem:maximal-local-supports} now gives
\eqref{eq:endpoint-local-cover}: each actual label belongs to the
endpoint population and its support is contained in its assigned
maximal support, with the same type. Apply the union bound once to
each pair $(J,P)$. Several labels may be covered by the same event;
conversely, equality of probability bounds does not make distinct
events equal. Grouping equal summands proves the last assertion.
\end{proof}

The cover is formed before deleting the odd core labels. This can
include events unnecessary for a particular tree, but cannot omit
an event needed by that tree. In particular, no choice of the core
primes is hidden in the event enumeration.

\subsubsection{Choosing the full-type exceptions}

Set
\[
 b_F(u,y)=\abs{\{j\leq\floor{u/2}:T(j)=0,
                  \text{ some }p\in P_u(y)\text{ divides }j\}}.
\]
Fix the desired full-type target $\tau$. Two choices of an endpoint
exception set will be used:
\begin{align}
 \mathcal E_{\rm all}
   &=\{y\leq u:y\text{ odd},\ T(y)=F\},
                                      \label{eq:all-local-exceptions}\\
 \mathcal E_{\rm prot}
   &=\{y\in\mathcal E_{\rm all}:
        D b_F(u,y)\geq\ell\rho(\widetilde c_F-\tau)\}.
                                      \label{eq:protected-local-exceptions}
\end{align}
The second set uses a complete half-label prefix. It therefore
includes every exception that might remain at any order in the interval.

\begin{proposition}[Uniform exceptional degree targets]
\label{prop:local-exception-pool}
Choose either set $\mathcal E$ above and put
\begin{equation}\label{eq:local-exception-fraction}
 \epsilon=\frac{\abs{\mathcal E}}{\ell\rho}.
\end{equation}
Then $\abs{\mathcal E\cap Y}\leq\epsilon N$. Every actual full-type
label outside $\mathcal E$ has degree greater than $\tau N$ for every
class permutation. For the protected choice this follows from the
strict complementary inequality in
\eqref{eq:protected-local-exceptions}.

Suppose $0\leq\epsilon<\tau\leq\eta\leq\gamma$ and
$\gamma+\upsilon\leq1$. Set
\begin{equation}\label{eq:local-certificate-targets}
 t_J=\begin{cases}
 \gamma,&J=0,\\
 \epsilon,&J=F,\\
 \eta,&0\ne J\ne F,
 \end{cases}
 \qquad \theta_J=\rho(\widetilde c_J-t_J).
\end{equation}
If all $\theta_J>0$, avoiding all events
$K(J,P)\geq n\theta_J$ with $P\in\mathcal M_J(u)$ gives every
strict degree condition of Lemma~\ref{lem:three-hall-ranges}.
\end{proposition}

\begin{proof}
The population bound follows from $N\geq n\rho\geq\ell\rho$.
For $\mathcal E_{\rm all}$ there is no full-type label outside the
set. For $\mathcal E_{\rm prot}$, any such label satisfies
\[
 D b_F(u,y)<\ell\rho(\widetilde c_F-\tau).
\]
Apply Lemma~\ref{lem:protected-full-labels} with
$R_0=\ell\rho$ and $c=\widetilde c_F$.
If a label exists outside this set, the strict inequality itself
implies $\widetilde c_F>\tau$, as required there.

By Proposition~\ref{prop:finite-quota-thresholds}, failure of a
strict degree bound $d(y)>t_{T(y)}N$ implies
$K(T(y),P_u(y))\geq n\theta_{T(y)}$. The endpoint cover contains
this failure. Avoiding its events gives the target $\epsilon N$
for all full-type labels, hence for the exceptions. The preceding
deterministic argument gives the stronger target $\tau N$ for the
others. The proper nonempty and empty types have the required targets
$\eta N$ and $\gamma N$. These are exactly the four degree bounds
in Lemma~\ref{lem:three-hall-ranges}.
\end{proof}

We use the weaker target $\epsilon$ on the whole full-type event
family. Thus the support cover is independent of the choice of the
exception set. The local estimates establish a common target; the
protection test then supplies the stronger degree bound for each
nonexceptional full-type label.

\subsubsection{One probability bound for each maximal support}

Fix $(J,P)$ from the endpoint cover. Use the data
\[
 A_P,\ H_P,\ G_P,\ L_P,\ b_-,\ b_+,\ B_P
\]
from Lemma~\ref{lem:finite-local-data}. For consistency with the finite
data below, the lower empty-class frequency in $L_P$ may use the
complete upper-prefix population $c_0^+$ in place of $a_0^+$.
This only decreases a lower bound. The other lower class frequencies
are unchanged.

Define
\begin{equation}\label{eq:local-certificate-raw-mean}
 \mu_0=e_+A_P+H_+H_P+\delta G_P-q_-L_P.
\end{equation}
The construction has either $\delta=0$ or $q_-=0$; the two
corrections are never combined from different assignments.
With $Q_0=10^9$, put
\begin{equation}\label{eq:local-certificate-rounded-data}
 s=\frac{\ceil{Q_0 D b_+}}{Q_0},\qquad
 \mu=\min\left\{\frac{\ceil{Q_0\mu_0}}{Q_0},s\right\},\qquad
 \theta_* =\min_{\ell\leq n\leq u}
                         \frac{\ceil{n\theta_J}}n.
\end{equation}
All extrema here and below are over integer orders.
The minimum with $s$ is valid because the selected degree sum is
at most $Db$, where $b/n\leq b_+$. In particular it remains an
upper bound for the normalized mean.

We use the following notation for already proved rational bounds.
Let $\mathcal B_d(\theta;\mu,s)$ be the bounded-sample estimate from
\eqref{eq:directed-bounded-rate} and
\eqref{eq:directed-probability-map}, including the deterministic
support cases. Let $\mathcal C_d(\theta;\mu,B)$ be the convex
binomial estimate of Proposition~\ref{prop:integer-local-tail},
with the upward rounding of Section~\ref{sec:directed-probability}.
The arguments $[\ell,u]$ are fixed and suppressed. Thus both are
upper bounds for the same event when their stated mean and sample
bounds hold.

For a connected construction, set
\begin{equation}\label{eq:local-certificate-baseline}
 \mu_X=(e_+A_P+H_+H_P-b_-)_+.
\end{equation}
Its additional bound is
\begin{equation}\label{eq:local-certificate-baseline-bound}
 \begin{split}
 q^{\rm base}_{J,P}=\min\{&
 \mathcal B_{D-1}(\theta_*-b_+;\mu_X,(D-1)b_+),\\
 &\mathcal C_{D-1}(\theta_*-b_+;\mu_X,B_P)\}.
 \end{split}
\end{equation}
In an original-forest row this bound is omitted.

A row may also specify an integer $1\leq a_{J,P}\leq160$ and use
\[
 \lambda=\frac{a_{J,P}}{16D},\qquad
 \bar q=\min\left\{\frac{D-1}{D}\mu,
         \frac{\ceil{u/2}-1}{\ell}H_P\right\}
\]
in the directed forest estimate
\eqref{eq:directed-forest-rate}. Denote its probability bound by
$q^{\rm for}_{J,P}$. This option is used only when $D\geq2$.
The parameter zero in the finite data means that the option is
omitted; it does not represent a zero failure probability.

Finally, define
\begin{equation}\label{eq:local-certificate-event-bound}
 q_{J,P}=\min\{
   \mathcal B_D(\theta_*;\mu,s),\
   \mathcal C_D(\theta_*;\mu,B_P),\
   q^{\rm base}_{J,P},\ q^{\rm for}_{J,P}\},
\end{equation}
where omitted bounds are removed from the minimum.
Both sample estimates take the value one for a nonpositive threshold.
If $B_P=0$, the original positive-threshold event is impossible and
$q_{J,P}=0$; no forest moment is needed in this case.

\begin{theorem}[Simultaneous local degree certificate]
\label{thm:finite-local-degree-certificate}
Assume the preparation and construction conditions of
Section~\ref{sec:interval-construction}, the exceptional targets in
Proposition~\ref{prop:local-exception-pool}, and nonnegative mean
bounds in \eqref{eq:local-certificate-raw-mean}. Then, for each
actual tree and order in the row, the probability that any required
local degree condition fails is at most
\begin{equation}\label{eq:finite-local-failure}
 p_{\rm loc}=\sum_J\sum_{P\in\mathcal M_J(u)}q_{J,P}.
\end{equation}
Every term is obtained by rational arithmetic. No additional factor
for the number of labels, orders, or input trees is required.
\end{theorem}

\begin{proof}
Proposition~\ref{prop:finite-local-means} gives
$\E K(J,P)/n\leq\mu_0$. The rounded data in
\eqref{eq:local-certificate-rounded-data} retain the mean and
sample-degree upper bounds. Since $K$ is integer-valued,
\[
 \{K/n\geq\theta_J\}\subseteq\{K/n\geq\theta_*\}.
\]
The two unconditional terms in
\eqref{eq:local-certificate-event-bound} therefore bound the latter
event by the results of Sections~\ref{sec:local-probability}
and~\ref{sec:directed-probability}.

In the connected case write the completed sum as $K^Q=b+X$.
Equation~\eqref{eq:finite-baseline-data} gives
$\E X/n\leq\mu_X$ and the event inclusion
\[
 \{K^G/n\geq\theta_*\}
       \subseteq\{X/n\geq\theta_*-b_+\}.
\]
The sampled weights of $X$ lie in $[0,D-1]$, and their number is
at most $B_P$. This proves the baseline bound. Notice that
\eqref{eq:local-certificate-baseline} uses the uncredited completed
mean. The edge credit in \eqref{eq:local-certificate-raw-mean}
is not subtracted again.

For the forest option, Proposition~\ref{prop:forest-local-moment}
bounds the normalized sampled excess by $\bar q$. Indeed
$N-1\leq\ceil{u/2}-1$ and the largest class sampling fraction is
at most $H_P$. Moreover $0<D\lambda\leq10$, so the directed
exponential bounds apply. Corollary~\ref{cor:uniform-forest-tail}
and Proposition~\ref{prop:directed-moment-evaluation} now prove
$q^{\rm for}_{J,P}$. Every included term bounds the same original
event, so taking their minimum is valid.

Apply Lemma~\ref{lem:endpoint-local-cover} and then
Proposition~\ref{prop:local-exception-pool}. The permutations are
chosen separately for each fixed tree and order. The bounds are
uniform over these choices of input, which is why there is no union
bound over the inputs themselves.
\end{proof}

Numerical grouping can reduce the number of probability evaluations.
It keeps $J$, the degree target, all sample and mean bounds, the
baseline choice, and the chosen forest parameter fixed. The grouped
quantity is multiplied by the number of maximal events represented.
A search may select the forest parameter, but acceptance uses only
its integer index and the rational bounds above.

\subsubsection{Finite verification of the local conditions}

The supplement supplies the exceptional-set choice and all forest
parameter indices for each of the 509 construction rows in
Table~\ref{tab:interval-degree-cover}. No list of odd labels or
preaccepted support cover is needed: these are reconstructed from
the order interval and prime blocks.

\begin{proposition}[Local certificates for orders 1360 through 1521]
\label{prop:verified-local-certificates}
For every one of those 509 rows, the exceptional population satisfies
$0\leq\epsilon<\tau$, all margins $\theta_J$ are positive, and
\begin{equation}\label{eq:verified-local-upper-bound}
 p_{\rm loc}<\frac{997}{1000}.
\end{equation}
The 110 rows using $\mathcal E_{\rm all}$ and the 399 rows using
$\mathcal E_{\rm prot}$ together cover 363,264 endpoint label
occurrences by 152,825 maximal-support events. Among these events,
134,901 have a deterministic zero bound. The largest local failure
bound occurs on $[1386,1389]$ at $D=10$, and is less than
$996081/10^6$.
\end{proposition}

\begin{proof}
The standalone checker first accepts the complete construction cover
of Proposition~\ref{prop:finite-interval-construction-check}.
For each row it enumerates every odd integer at most $u$, computes
its type and prime support, and forms the maximal families. It assigns
each label to exactly one event and checks that every retained
maximal support has a representing label. Thus both the event list
and all multiplicities are complete.

It then constructs the indicated exceptional set using exact
empty-class union counts. The strict protection test is checked for
every label excluded from the exception set. Division by $\ell\rho$
gives $\epsilon$; the checker verifies $\epsilon<\tau$ and every
positive margin in \eqref{eq:local-certificate-targets}.
For each maximal event it recomputes the endpoint arithmetic,
rounded mean, integer threshold, and every selected probability
bound. The sum in \eqref{eq:finite-local-failure} is rational and
is compared to the displayed bounds by integer cross-multiplication.

An independent reconstruction uses a sieve for primes, division for
prime supports, and gcd tests for types and union counts. It checks
the class-frequency and sample-population bounds at every integer
order in each interval. The two computations agree on all exceptional
populations and local failure bounds. The finite counts and the
largest value stated above follow from this complete row list.
\end{proof}

\begin{table}[ht]
\centering
\caption{Local events in the complete construction cover. Counts are
summed over the rows; an event occurring in two different rows is
counted in each. A zero bound excludes that event deterministically.}
\label{tab:local-certificate-counts}
\par\vspace{.5\baselineskip}
\renewcommand{\arraystretch}{1.15}
\begin{tabular}{lrrr}
\hline
Construction&Rows&Maximal events&Zero bounds\\
\hline
Cuts confined to $B$&64&18,731&18,731\\
Donor exchanges&231&68,210&58,417\\
Connected height&214&65,884&57,753\\
\hline
Total&509&152,825&134,901\\
\hline
\end{tabular}
\end{table}

This establishes the local part of
Theorem~\ref{thm:integer-interval-certificate}. In particular, the
first Hall range is absent whenever these local events are avoided,
by Corollary~\ref{cor:interval-caps}. The middle and high ranges
still need their own deterministic or probabilistic bounds. Their
failure probabilities must be added to the actual value of
$p_{\rm loc}$ in each row, with total strictly below one. The
uniform bound \eqref{eq:verified-local-upper-bound} does not replace
that final row-by-row comparison.

\subsection{Hall certificates for orders 1360 through 1521}
\label{sec:hall-certificates}

We complete the interval argument by excluding its remaining deficient
sets. The local degree conditions are those of
Section~\ref{sec:local-certificates}. They already exclude the
full-only range. For the middle range, a connected assignment permits
an arithmetic restriction on the types of positions in the deficient
rectangle. For the high range, the forest deletion bound is often
sufficient for every permutation. The other cases use the global
exponential estimate. We compare the sum of the required probabilities
with one separately in each row.

Retain an integer interval $[\ell,u]$, its residual matching side of
size $N$, and the bounds of Section~\ref{sec:interval-construction}.
Put
\[
 N_0=\ell\rho,\qquad H=\ceil{u/2}.
\]
Thus $N_0\leq N\leq H$. Write $p_I,\kappa_I$ for the caps in
\eqref{eq:interval-hall-caps}. The three label cap vectors are those
of Theorem~\ref{thm:three-range-hall}; in particular,
\[
 p_I^{(2)}=p_I\boldsymbol1_{I\ne0},\qquad
 p_I^{(3)}=
 \begin{cases}(p_0-\upsilon)_+,&I=0,\\p_I,&I\ne0.\end{cases}
\]
All caps are relative to the actual side size $N$. They are uniform
upper bounds and need not sum to one.

\subsubsection{Minimization at fixed mass}

The high-range deletion bound depends on the number of selected
positions. We therefore keep that mass fixed before taking a minimum
of the compatible pairs. The next lemma also gives a second way to
compute the rectangle minima of Section~\ref{sec:refined-hall}.

\begin{lemma}[Concavity between subset masses]
\label{lem:fixed-mass-concavity}
Let $C_{ij}\geq0$ be a finite matrix and let $p_i,q_j\geq0$ be
coordinate caps. Suppose
\[
 I=\left[\max\{0,1-\textstyle\sum_jq_j\},\
          \min\{1,\textstyle\sum_ip_i\}\right]
\]
is nonempty. For $t\in I$, define
\begin{equation}\label{eq:fixed-mass-function}
 G(t)=\min\sum_{i,j}C_{ij}x_iw_j,
\end{equation}
where the minimum is over
\[
 0\leq x\leq p,\qquad0\leq w\leq q,\qquad
 \sum_i x_i=t,\qquad\sum_jw_j=1-t.
\]
The function $G$ is continuous. Partition $I$ at all points of the
forms
\begin{equation}\label{eq:subset-mass-breakpoints}
 \sum_{i\in A}p_i,\qquad 1-\sum_{j\in B}q_j
\end{equation}
that lie in $I$. On each resulting closed interval $G$ is concave.
For fixed $t$, it is enough to enumerate the vertices of the left
capped simplex and fill the right coordinates in increasing order
of their costs. Rational data therefore give exact rational values.
\end{lemma}

\begin{proof}
For fixed $w$, the objective is affine in $x$, so some minimizing
$x$ has at most one coordinate strictly between its bounds. With
this $x$ fixed, the same assertion holds for $w$. Thus some joint
minimum is attained at vertices of both capped simplices.

A left vertex family has a set $A$ of coordinates at their caps,
one possible remaining coordinate $i\notin A$, and
\[
 x_i=t-\sum_{k\in A}p_k.
\]
Its interval of validity has endpoints that are subset sums of $p$.
The analogous right family has
$w_j=1-t-\sum_{k\in B}q_k$, and its endpoints are complements of
subset sums of $q$. All-bound vertices are included at their fixed
masses. This also covers zero caps and one-point feasible intervals.

On an open interval between consecutive points in
\eqref{eq:subset-mass-breakpoints}, the available vertex families
are fixed. For one pair of families the objective is a quadratic
in $t$ with leading coefficient $-C_{ij}\leq0$. It is therefore
concave. The minimum of finitely many concave functions is concave:
apply the defining concavity inequality to every function and then
take their minimum. This proves the assertion on the open intervals.

To check continuity at an endpoint $t_0$, a convergent subsequence
of minimizers shows that any limiting value is at least $G(t_0)$,
by compactness and continuity of the objective. In the other direction,
choose a minimizer at $t_0$ and any feasible pair at a mass on the
specified side of $t_0$. Convex interpolation between the two pairs
is feasible at every intermediate mass, and its objective converges
to $G(t_0)$. This gives the reverse limiting inequality. Concavity
therefore extends to the closed intervals.

Finally, for a fixed left vector $x$, the cost of the right
coordinate $j$ is $\sum_iC_{ij}x_i$. Filling cheaper coordinates
first minimizes this linear expression under their caps and the
prescribed total mass. An exchange from a more expensive occupied
coordinate to a cheaper unfilled one proves the assertion. All
vertices, costs, and partial fills are rational for rational data.
\end{proof}

In particular, to minimize $G(t)$ on a closed subinterval of $I$,
it is enough to check its endpoints and the points in
\eqref{eq:subset-mass-breakpoints}. More generally, the same is true
of $G(t)-L(t)$ after adding every breakpoint of a continuous piecewise
linear function $L$. On each remaining interval $L$ is affine, so
$G-L$ is concave. This use of concavity does not assert that $G$ is
concave across the subset masses.

\subsubsection{A deterministic test for the high range}

Use the complete endpoint half-label pools, including the special
labels removed in the preparation. For a half-label $a$ of type $I$,
put
\[
 q(a)=Q_{u,I}^{\mathcal T}(a),\qquad
 Q_0=\sum_{a\leq\floor{u/2}}q(a),
\]
where the cost is defined in \eqref{eq:restricted-prime-cost}.
Let $L(t)$ be the largest-entry sum of the multiset containing
$D-1$ copies of each $q(a)$, with a fractional last entry and constant
extension after exhaustion. If $D=1$, set $L=0$.
The use of the complete pools only increases the forest bound.

For caps $p^{(3)},\kappa$, let $G$ be the function in
\eqref{eq:fixed-mass-function}, with the usual disjointness matrix.
The feasible position masses form
\begin{equation}\label{eq:high-mass-interval}
 J=\left[\max\{\upsilon,1-\textstyle\sum_Ip_I^{(3)},0\},\
          \min\{1-\gamma,\textstyle\sum_I\kappa_I,1\}\right].
\end{equation}
For $v\in J$, define
\begin{equation}\label{eq:high-forest-slack}
 \Psi(v)=\rho^2G(1-v)-\frac{Q_0+L((Hv-1)_+)}{\ell^2}.
\end{equation}

\begin{proposition}[High-range forest certificate]
\label{prop:high-forest-certificate}
Assume the preparation, type assignment, and local degree conditions
of Theorem~\ref{thm:integer-interval-certificate}. If $J$ is empty,
or if $\Psi(v)>0$ for every $v\in J$, then no deficient set belongs
to the high range. The assertion holds for every permutation within
the even classes satisfying the local conditions.

For rational data, the strict positivity test requires only the
endpoints of $J$, its position subset masses and complementary label
subset masses, the point $1/H$, and the points
\begin{equation}\label{eq:high-cost-breakpoints}
 v=(1+m)/H
\end{equation}
where $m$ is a cumulative multiplicity in the decreasing cost list
of $L$. Only points lying in $J$ are retained.
\end{proposition}

\begin{proof}
Suppose a high-range deficient set exists. Write $U$ for its labels
and $V=B\setminus N_\Gamma(U)$ for its complementary positions.
As in Theorem~\ref{thm:three-range-hall}, trim their weights to
combined normalized mass one, preserving label mass at least
$\gamma$ and position mass at least $\upsilon$. Let $v$ be the
latter mass. Universal labels do not lie in $U$, so its label
coordinates satisfy the caps $p^{(3)}$.

Apply full-type transport to the selected position coordinates.
This preserves $v$ and can only decrease compatible mass. The
transported caps are $\kappa$, so $v\in J$ and the original
compatible mass, normalized by $N^2$, is at least $G(1-v)$.
Every such original compatible pair is deleted, since the selected
weights are supported on $U\times V$.

The deletion count is bounded on the original positions, before
transport. Corollary~\ref{cor:forest-rectangle-deletions}, enlarged
to the complete endpoint pools, gives the bound
$Q_0+L((Nv-1)_+)$. Since $N\leq H$, this is at most
$Q_0+L((Hv-1)_+)$. Consequently
\[
 \rho^2G(1-v)\leq\frac{N^2}{n^2}G(1-v)
       \leq\frac{Q_0+L((Hv-1)_+)}{n^2}
       \leq\frac{Q_0+L((Hv-1)_+)}{\ell^2}.
\]
Here the first inequality uses $N/n\geq\rho$ and $G\geq0$.
This contradicts \eqref{eq:high-forest-slack}. The argument uses
the forest on the original remaining vertices; added completion
edges are not counted as actual conflicts.

For the finite test, apply Lemma~\ref{lem:fixed-mass-concavity}
with $t=1-v$. The second term of $\Psi$ is continuous and piecewise
linear. Its only additional breakpoints are $1/H$ and
\eqref{eq:high-cost-breakpoints}. Between all these points $\Psi$
is concave, so its minimum occurs at an endpoint. When $D=1$ the
cost function is constant, and the redundant point $1/H$ is harmless.
A one-point interval is checked directly. All evaluations are rational.
\end{proof}

The distinction between the transported and original positions in
this proof is useful. Artificial full-type mass need not represent
a set of actual full-tag positions. We use transport only to lower
the compatible mass; the forest deletion bound is applied to the
original weighted set, whose total mass is unchanged.

\subsubsection{The middle range and the remaining probability terms}

In each connected row, the positions in the separating set have full
tags and every other position has the same type as all its neighbors
in the completed tree. Proposition~\ref{prop:interval-construction-bounds}
therefore implies \eqref{eq:pure-neighbor-tags} also in the original
residual forest. This is the hypothesis needed for the support
restriction; it is not inferred merely from inclusion of neighbor
types in a tag.

For every nonempty $S\subseteq\mathcal T$, compute the endpoint
counts $c_Q(u,S)$ by \eqref{eq:arithmetic-support-test} and
\eqref{eq:support-label-condition}. In
Proposition~\ref{prop:separate-forest-hall}, take
\[
 R_0=N_0,\qquad N_{\rm upper}=H,\qquad
 a=\eta,\qquad b=\gamma.
\]
Here $N_{\rm upper}$ denotes the upper side bound called $N$ in
that proposition. Use the costs
$Q_{u,J}^{\mathcal T\setminus\{0\}}(a)$ in the complete half-label
pool of type $J$. Write $B_J$ for their sum and $L_J$ for their
largest-entry sum with multiplicity $D-1$. The required minimum is
\begin{equation}\label{eq:finite-middle-slack}
 \min\left\{
   \sum_{Q,J}C(Q,J)x_Qw_J
     -\frac{\sum_{J\in S}(B_J+L_J(Hw_J))}{N_0^2}
       \right\},
\end{equation}
where
\[
 \begin{gathered}
 0\leq x_Q\leq c_Q(u,S)/N_0,\qquad
 0\leq w_J\leq\kappa_J\boldsymbol1_{J\in S},\\
 \sum_Qx_Q+\sum_Jw_J=1,\qquad
 \eta\leq\sum_Qx_Q\leq\gamma.
 \end{gathered}
\]
An empty domain contributes no obstruction. Positive minima on all
other domains exclude the middle range for every class permutation.
Lemma~\ref{lem:separate-cost-minimum} evaluates these minima exactly.
The argument uses $L_J(Hw_J)$, without subtracting one inside its
argument, so the column cost retains the convexity required by that
lemma.

For completeness, an independent evaluation can use the fact that
$C(Q,J)$ is either zero or one. Fix a left vertex family with variable
mass $t$. The marginal costs of the right segments have slopes zero
or one in $t$. Within each slope class their order is fixed. Enumerate
the two prefixes filled completely and the possible last partial
segment. Requiring the full segments to be no more expensive than
the last, and the unused segments to be no cheaper, gives linear
inequalities in $t$. Its remaining length is $1-t-d$, where $d$
is the total full-segment length. The capacity inequalities give
another closed interval. On its intersection with the family domain
the objective is concave, so its two endpoints suffice. Include
all-bound left choices and the case in which every right segment
is full. This enumerates the same feasible optima without locating
all pairwise crossings of marginal costs.

For a range not excluded deterministically, put
\begin{align}
 b_2&=\rho^2\mathcal B(p^{(2)},\kappa;\eta,1-\gamma),
                  \label{eq:finite-middle-budget}\\
 b_3&=\rho^2\mathcal B(p^{(3)},\kappa;\gamma,\upsilon).
                  \label{eq:finite-high-budget}
\end{align}
An empty domain is omitted. Every remaining budget must be positive.
For $k=2,3$, use the odd-type families
\[
 \mathcal J_2=\mathcal T\setminus\{0\},\qquad
 \mathcal J_3=\mathcal T,
\]
and the row's positive rational parameter $t_k$. Let
$\widehat{\mathcal M}_k$ be the directed global moment for that
family and the chosen construction. Proposition
\ref{prop:directed-moment-evaluation} gives the failure bound
\begin{equation}\label{eq:finite-hall-failure}
 p_k=\mathcal P\left(\frac\ell D
                    (t_kb_k-\widehat{\mathcal M}_k)\right).
\end{equation}
Set $p_k=0$ for a deterministically excluded range. The original
construction uses its donor term; the connected construction uses
its completion credit. They are never combined in one moment.

There is no independence assertion between these events or between
them and the local events. The acceptance test is the single strict
inequality
\begin{equation}\label{eq:finite-hall-joint}
 p_{\rm loc}+p_2+p_3<1.
\end{equation}
Theorem~\ref{thm:integer-interval-certificate} then applies. The
full-type term is zero because $p_F<\tau$ in every row, as already
proved in Proposition~\ref{prop:finite-interval-construction-check}.

\subsubsection{Verification and the interval theorem}

The finite data specify the two deterministic choices and the two
global moment parameters for each row of
Table~\ref{tab:interval-degree-cover}. The construction and local
inputs remain those of Sections~\ref{sec:interval-construction}
and~\ref{sec:local-certificates}. The following assertion is the
complete acceptance test, rather than a comparison of conditional
moment values.

\begin{proposition}[Verified Hall and joint bounds]
\label{prop:verified-hall-certificates}
All 509 rows in Table~\ref{tab:interval-degree-cover} satisfy the
middle and high Hall conditions above and
\eqref{eq:finite-hall-joint}. More precisely,
\begin{equation}\label{eq:verified-joint-upper-bound}
 p_{\rm loc}+p_2+p_3<\frac{996242}{10^6}
\end{equation}
for every row. The largest sum occurs on $[1386,1389]$ at degree
$10$. All 214 connected rows pass the support test for the middle
range, and 438 rows pass the forest test for the high range. The
remaining 295 middle and 71 high cases pass the exponential checks.
\end{proposition}

\begin{proof}
The supplement supplies the complete rational input table and a
standalone checker. We describe its finite computations and why they
establish the assertion.

First the checker repeats the construction and local acceptance for
all 509 rows. The interval endpoints and degree lists must be exactly
those of Table~\ref{tab:interval-degree-cover}, and each local and
Hall row must have the same interval and degree as its construction
row. Thus no selection of a successful subcollection is permitted.
The caps, exception populations, and local failure sums are
recomputed from these inputs.

Next the checker forms the three capped domains. The full-only
domain is empty in every row. It evaluates both other rectangle
minima using Lemma~\ref{lem:fixed-mass-concavity}. A separate
enumeration of the three vertex forms in
Proposition~\ref{prop:capped-bilinear-minimum} gives the same minima.
All budgets used in a probability estimate are positive.

For the 438 high-range forest cases, the complete endpoint cost
histograms give $Q_0,L$ and every breakpoint in
Proposition~\ref{prop:high-forest-certificate}. There are 2,403
point evaluations, counted separately for each row. Their least
slack is the positive rational
\[
 \min\Psi=\frac{1637}{400207500000000}.
\]
This small margin is compared to zero exactly, without rounding to
a displayed decimal. Direct expansion of the repeated cost list and
a separate fixed-mass vertex enumeration reproduce these minima.

In each connected row there are four types, hence $2^4-1=15$
nonempty position supports. The checker treats all 3,210 supports;
2,116 have feasible domains. Each feasible minimum in
\eqref{eq:finite-middle-slack} is positive. The primary computation
uses the segment procedure of Lemma~\ref{lem:separate-cost-minimum};
the two-prefix procedure described above gives the same rational
minimum on every support. The endpoint label counts are also
reconstructed by testing gcds against the complete half-label pools.
The smallest support slack is greater than $1378/10^7$; its decimal
value is approximately $0.0001378211442$.

There remain 366 global moment evaluations. For each, the checker
recomputes the exact populations, endpoint cost histograms, and
directed logarithmic moments. It checks the domain of every
exponential evaluation and then applies
\eqref{eq:finite-hall-failure}. The sum with the row's actual local
bound is rational. Every sum satisfies
\eqref{eq:verified-joint-upper-bound} by integer cross-multiplication.
The largest sum, on the stated row, is approximately
$0.9962410558965252$. Its three components are approximately
\[
 p_{\rm loc}=0.9960809297572989,\qquad
 p_2=0.0001601261392263,\qquad p_3=0.
\]
These decimal values only describe the exact comparisons recorded
in the supplement. Table~\ref{tab:hall-certificate-counts} gives
the counts by construction. This verifies every condition asserted.
\end{proof}

\begin{table}[ht]
\centering
\caption{Deterministic exclusions in the accepted interval table.
Each remaining middle or high range uses its global exponential
bound. All 509 full-only ranges are excluded by their population.}
\label{tab:hall-certificate-counts}
\par\vspace{.5\baselineskip}
\renewcommand{\arraystretch}{1.15}
\begin{tabular}{lrrr}
\hline
Construction&Rows&Middle support&High forest\\
\hline
Cuts confined to $B$&64&0&64\\
Donor exchanges&231&0&185\\
Connected height&214&214&189\\
\hline
Total&509&214&438\\
\hline
\end{tabular}
\end{table}

\begin{theorem}[The interval from 1360 through 1521]
\label{thm:interval-1360-1521}
Every tree of order $n$ with $1360\leq n\leq1521$ has a prime
labeling.
\end{theorem}

\begin{proof}
Fix such an integer $n$ and a tree $T$ of order $n$. Proposition
\ref{prop:finite-interval-construction-check} supplies the uniform
preparation and a row covering its residual degree. If the residual
forest is edgeless, the immediate case in
Theorem~\ref{thm:integer-interval-certificate} completes the labeling.
Otherwise choose the row for its positive degree $D$. The necessary
degree exclusions cannot remove this actual order--degree pair, by
Lemma~\ref{lem:interval-degree-cover}.

Fix the type assignment supplied by that row. Proposition
\ref{prop:verified-local-certificates} provides the exception bounds,
strict local targets, and their joint local failure estimate.
Proposition~\ref{prop:verified-hall-certificates} supplies the Hall
conditions and makes the full failure sum strictly less than one.
Theorem~\ref{thm:integer-interval-certificate} therefore gives a
permutation for which the residual matching exists, and restores
a prime labeling of $T$.

The preparation and probability bounds hold separately for every
input tree at every integer order in its interval. The 32 adjacent
intervals cover all of $[1360,1521]$. Thus the conclusion holds for
every $T$ and $n$ in the statement.
\end{proof}

\begin{samepage}
\begin{corollary}[The established initial interval]
\label{cor:initial-prefix-1521}
Every tree of order at most $1521$ has a prime labeling.
\end{corollary}

\begin{proof}
Theorem~\ref{thm:initial-prefix-1359} treats $1\leq n\leq1359$.
Theorem~\ref{thm:interval-1360-1521} treats the next integer interval,
whose first order is $1360=1359+1$.
\end{proof}
\end{samepage}

The argument for subsequent orders must provide its own construction
bounds and accepted range data, followed by the estimates for an
unbounded interval. The finite certificate just proved establishes
no assertion about those remaining orders.

\subsection{Continuation through order 1753}
\label{sec:forest-continuation}

The interval criterion also treats the next 232 orders. We use the
same three retained-position constructions, the same local event
cover, and the same rational evaluation rules. The differences lie
in the finite data: some degree thresholds coincide, and the middle
range is handled by its global exponential bound whenever it is
nonempty. The arithmetic support test is not needed in this interval.

\subsubsection{The degree cover and empty ranges}

Table~\ref{tab:forest-continuation-cover} gives 52 adjacent intervals
with union $[1522,1753]$. For each interval $[\ell,u]$, the last
column of its entry is
\[
 D_* = \ceil{\frac1{2W}},\qquad
 W=\min_{\ell\leq n\leq u}\frac{w_n}{n}.
\]
Every degree $1\leq D\leq D_*$ has a row. As before, an actual
order--degree pair must satisfy
\eqref{eq:interval-degree-feasibility}; no other degree exclusion
is used. All rows have four types obtained from two nonempty
blocks of the odd primes at most their cutoff.

Each row specifies the construction, the prime blocks and, when
needed, the height choice. It supplies the bounds in
\eqref{eq:interval-construction-interface}: edge and degree errors,
transfer and retagging costs, residual side and universal fractions,
and any completion credit. Proposition~\ref{prop:interval-construction-bounds}
verifies these bounds with its stated hypotheses.

The degree targets obey
\[
 0\leq\epsilon<\tau\leq\eta\leq\gamma,
 \qquad \gamma+\upsilon\leq1.
\]
In 391 rows $\tau=\eta$. The first Hall range is then empty,
since its defining inequalities would be
$\tau<\abs U/N\leq\eta$. In fact the stronger population condition
$p_F<\tau$ holds in every row, so it excludes the first range even
when these thresholds differ.

In 24 rows $\eta=\gamma$, and the middle range is empty for the
same reason. These are all original-forest rows with cuts confined
to $B$; in them all three thresholds coincide. We set the omitted
range's failure probability to zero. This is an exclusion of an
empty family of deficient sets, not a zero estimate for a positive
probability event. The exceptional-label and local degree conditions
are still checked in each such row.

The remaining 831 middle cases use
\eqref{eq:finite-middle-budget} and the family
$\mathcal T\setminus\{0\}$ in
\eqref{eq:finite-hall-failure}. In the high range, 800 rows use
Proposition~\ref{prop:high-forest-certificate}; the other 55 use
\eqref{eq:finite-high-budget} and the family $\mathcal T$ in the
same exponential estimate. Thus each row has the complete failure
sum
\[
 p_{\rm loc}+p_2+p_3,
\]
with zero terms only for the deterministic exclusions just specified.

\begin{table}[!htbp]
\centering
\caption{The continuation cover. Each entry includes all degrees
$1\leq D\leq D_*$. A singleton denotes one integer order.}
\label{tab:forest-continuation-cover}
\par\vspace{.5\baselineskip}
\renewcommand{\arraystretch}{1.15}
\begin{tabular}{rrrr}
\hline
Order interval&$D_*$&Order interval&$D_*$\\
\hline
$1522$&16&$1582$--$1585$&17\\
$1523$--$1524$&16&$1586$--$1589$&17\\
$1525$--$1526$&16&$1590$--$1593$&17\\
$1527$&16&$1594$--$1597$&16\\
$1528$&16&$1598$--$1599$&16\\
$1529$&16&$1600$--$1603$&17\\
$1530$--$1531$&16&$1604$--$1611$&17\\
$1532$--$1533$&16&$1612$--$1615$&17\\
$1534$--$1537$&17&$1616$--$1623$&17\\
$1538$--$1539$&16&$1624$--$1631$&16\\
$1540$--$1543$&16&$1632$--$1639$&17\\
$1544$--$1547$&16&$1640$--$1643$&17\\
$1548$--$1551$&16&$1644$--$1647$&17\\
$1552$--$1555$&16&$1648$--$1655$&16\\
$1556$--$1559$&16&$1656$--$1663$&16\\
$1560$--$1563$&16&$1664$--$1671$&16\\
$1564$&16&$1672$--$1679$&16\\
$1565$&16&$1680$--$1687$&16\\
$1566$&16&$1688$--$1695$&16\\
$1567$--$1568$&16&$1696$--$1703$&17\\
$1569$&17&$1704$--$1711$&17\\
$1570$&17&$1712$--$1719$&17\\
$1571$--$1572$&17&$1720$--$1727$&16\\
$1573$&17&$1728$--$1735$&17\\
$1574$--$1577$&17&$1736$--$1743$&17\\
$1578$--$1581$&17&$1744$--$1753$&17\\
\hline
\end{tabular}
\end{table}

\subsubsection{Acceptance of the finite data}

The local input has one exception choice and one permitted forest
parameter index for each maximal event. The global input has the
high-range forest choice and two positive rational moment parameters.
These are finite choices in the estimates already proved; their
selection is not a premise of the acceptance test. The accompanying
checker reconstructs the populations and every bound before making
the final comparison.

\begin{proposition}[Verified continuation data]
\label{prop:verified-forest-continuation}
The 855 rows in Table~\ref{tab:forest-continuation-cover} satisfy
all hypotheses of Theorem~\ref{thm:integer-interval-certificate}.
Every row has positive local margins, an exceptional fraction
$\epsilon<\tau$, and
\begin{equation}\label{eq:forest-continuation-joint}
 p_{\rm loc}+p_2+p_3<\frac{999049}{10^6}<1.
\end{equation}
The largest sum occurs at the singleton order $1569$, with residual
degree $10$. The construction and Hall choices have the counts in
Table~\ref{tab:forest-continuation-modes}.
\end{proposition}

\begin{proof}
We apply the exact tests of Sections~\ref{sec:interval-construction}
through~\ref{sec:hall-certificates}. We give the finite checks in
logical order, including the conditions needed to pass from the
arithmetic data to arbitrary input trees.

\emph{Preparation and construction.}
Direct enumeration of primes computes the special even population,
the upper-prime supply and the universal odd population at every
integer order in the table. It verifies
\eqref{eq:interval-preparation-tests} and the stated $D_*$.
The intervals meet at consecutive integers, and each degree list is
exactly $1,\ldots,D_*$. The checker applies the proved edge-and-isolation
condition to every order--degree pair.

For each row it recomputes the one-sided separator sum, the paid
donor stock and boundary bound, or the connected height budget,
as specified. Every supplied upper bound dominates its proved
construction bound. The side and universal lower bounds and the
completion credit have the opposite required direction. All
260 connected rows satisfy the completion capacity test. No
original-forest row receives a completion credit, and no connected
row receives a donor term.

It next enumerates the exact usable quota laws at every order,
then computes the total tag error, coordinate caps and compatible
fractions. All class population lower bounds are positive. The
full-type odd population cap is smaller than $\tau$ in all 855
rows. Independently reconstructed prime sets and gcd type counts
give the same supplies, laws and caps.

\emph{Local degrees.}
There are 274 rows with the full endpoint exception pool and
581 with the protected exception pool. The strict test in
Proposition~\ref{prop:local-exception-pool} determines the latter.
Every excluded full-type label has its deterministic degree
strictly above $\tau N$. The weak common target remains in the
full-type event cover, as in Section~\ref{sec:local-certificates}.

The maximal-support covers contain 288,977 event occurrences across
the rows, covering 686,580 endpoint label occurrences. For each
event the checker recomputes the class-frequency and sample bounds,
the integer threshold, and all selected tail estimates of
Theorem~\ref{thm:finite-local-degree-certificate}. Exactly 251,803
events have a deterministic zero bound. Summing once per maximal
event gives the row's $p_{\rm loc}$. The sum retains the multiplicity
of distinct events that happen to have the same numerical bound.

The independent arithmetic check reconstructs the maximal supports
by division and their label classes by gcd tests. It verifies the
frequency and sample enclosures at every integer order of the
corresponding interval, not just at its endpoints. It also repeats
the exceptional-label test and the integer-threshold comparison
for every applicable label or event. All local margins are positive,
and all exception fractions are strictly less than $\tau$.

\emph{Hall conditions and the joint comparison.}
The full-only range is absent in every row. The 24 middle ranges
with equal endpoints are omitted; the other middle budgets are
positive. The high budgets are positive as well. The checker
computes each required capped minimum by
Lemma~\ref{lem:fixed-mass-concavity}; a separate enumeration of
Proposition~\ref{prop:capped-bilinear-minimum} gives the same result.

For all 800 high-range forest choices, every breakpoint comparison
in Proposition~\ref{prop:high-forest-certificate} is strictly
positive. The smallest slack is
\[
 \frac{1885043}{19012500000000000}>0,
\]
attained in the row $[1560,1563]$, $D=15$.
The endpoint cost histograms are independently reconstructed by
counting eligible odd multiples of each prime. Expanding their
repeated cost lists gives the same forest minima. The comparison
uses the complete endpoint pools, so removing the special even
labels cannot increase the actual deletion cost.

For the remaining ranges there are 831 middle and 55 high
exponential evaluations. Each uses the appropriate odd-type family
and the degree, edge, variation and transfer data of its own
construction. The directed rules of
Proposition~\ref{prop:directed-moment-evaluation} give rational
upper bounds for all these failure probabilities. The checker
adds them to the actual local sum in the row. Integer
cross-multiplication verifies \eqref{eq:forest-continuation-joint}
for all 855 rows.

The largest sum is approximately $0.9990485710199962$ on the stated
row. Its exact rational value, together with the individual local
and global bounds, is supplied in the verification record. The
strict test uses that rational value rather than its decimal
representation. This proves all the asserted conditions.
\end{proof}

\begin{table}[!htbp]
\centering
\caption{Construction choices and deterministic exclusions in the
continuation. Every nonempty middle range uses an exponential bound;
so does each high range not excluded by the forest test.}
\label{tab:forest-continuation-modes}
\par\vspace{.5\baselineskip}
\renewcommand{\arraystretch}{1.15}
\begin{tabular}{lrrr}
\hline
Construction&Rows&Empty middle&High forest\\
\hline
Cuts confined to $B$&104&24&104\\
Donor exchanges&491&0&456\\
Connected height&260&0&240\\
\hline
Total&855&24&800\\
\hline
\end{tabular}
\end{table}

\subsubsection{Joining the two intervals}

\begin{theorem}[The initial interval through 1753]
\label{thm:initial-prefix-1753}
Every tree of order at most $1753$ has a prime labeling.
\end{theorem}

\begin{proof}
Corollary~\ref{cor:initial-prefix-1521} treats every order at most
$1521$. It remains to treat a tree $T$ of an integer order
$n\in[1522,1753]$.

Select the unique interval in
Table~\ref{tab:forest-continuation-cover} containing $n$, and
prepare $T$ as in Proposition~\ref{prop:highest-degree-preparation}.
The prime supplies and side bounds hold by
Proposition~\ref{prop:verified-forest-continuation}. If the residual
forest is edgeless, its immediate labeling and restoration are
those of Theorem~\ref{thm:integer-interval-certificate}. Otherwise
its positive maximum even-side degree belongs to the complete
row list, by Lemma~\ref{lem:interval-degree-cover}. The necessary
feasibility test cannot exclude this actual order--degree pair.

Fix the type assignment before taking the class permutations.
Proposition~\ref{prop:verified-forest-continuation} gives the exception
bounds and the local and Hall estimates. The strict joint inequality
gives positive probability of avoiding every required event. Apply
Theorem~\ref{thm:integer-interval-certificate} to match the remaining
odd labels to positions and restore a prime labeling of $T$.

This argument applies to every tree at every integer order in the
52 intervals. Their first order is $1522=1521+1$, so they adjoin
the preceding theorem without a gap. This proves the assertion.
\end{proof}

Section~\ref{sec:exponential-continuation} continues from order
$1754$, using two-sided variation bounds for the arithmetic type
laws.

\subsection{Two-sided type transport through order 2209}
\label{sec:exponential-continuation}

Both remaining Hall ranges will use exponential estimates. To sharpen
their coordinate bounds, we compare the odd-type
distribution with a reference law and move excess mass to the full
type. We apply the same operation to the position tags. Both selected
masses are preserved, and their compatibility can only decrease.

\subsubsection{Transport on both sides of a deficient rectangle}

Retain the matching graphs $\Gamma\subseteq\Gamma_0$ on
$Y\sqcup B$, with $\abs Y=\abs B=N$. Write $\xi$ for the type
law on $Y$ and $\zeta$ for the tag law on $B$. Let $\lambda,\nu$
be probability laws on $\mathcal T$, and suppose
\begin{equation}\label{eq:two-sided-variation}
 \|\xi-\lambda\|_{\mathrm{TV}}\leq v_X,
 \qquad
 \|\zeta-\nu\|_{\mathrm{TV}}\leq v_B.
\end{equation}
For the two ranges with $\eta=\tau<\gamma$, define
\begin{align}
 \widehat\kappa_I&=\nu_I+v_B\boldsymbol1_{I=F},
                                      \label{eq:transport-position-caps}\\
 \widehat p_I^{(2)}&=
 \begin{cases}
 0,&I=0,\\
 \lambda_I+v_X\boldsymbol1_{I=F},&I\ne0,
 \end{cases}                           \label{eq:transport-middle-caps}\\
 \widehat p_I^{(3)}&=
 \begin{cases}
 \min\{\lambda_0,\max(0,\lambda_0+v_X-\upsilon)\},&I=0,\\
 \lambda_I+v_X\boldsymbol1_{I=F},&I\ne0.
 \end{cases}                           \label{eq:transport-high-caps}
\end{align}
Here $\upsilon N$ is a lower bound for the number of empty-type
universal labels. These vectors bound transported submeasures;
they need not bound the original type coordinates individually.

\begin{proposition}[Hall bounds after two-sided transport]
\label{prop:two-sided-hall}
Assume the degree and universal-label hypotheses of
Lemma~\ref{lem:three-hall-ranges}, with $\eta=\tau<\gamma$ and
$u=\upsilon>0$, and assume \eqref{eq:two-sided-variation}.
For every feasible domain put
\[
 \widehat\beta_2=
   \mathcal B(\widehat p^{(2)},\widehat\kappa;\tau,1-\gamma),
 \qquad
 \widehat\beta_3=
   \mathcal B(\widehat p^{(3)},\widehat\kappa;\gamma,\upsilon).
\]
Suppose each such minimum is positive. If fewer than
$N^2\widehat\beta_2$ compatible pairs incident with nonempty
label types, and fewer than $N^2\widehat\beta_3$ compatible pairs
in total, are deleted in their respective feasible ranges, then
$\Gamma$ has a perfect matching. An infeasible domain needs no
deletion bound.
\end{proposition}

\begin{proof}
Suppose $U\subseteq Y$ is deficient and put
$V=B\setminus N_\Gamma(U)$. The first range is empty because
$\eta=\tau$. By Lemma~\ref{lem:three-hall-ranges}, either
\[
 \tau<\abs U/N\leq\gamma
 \quad\text{or}\quad
 \gamma<\abs U/N\leq1-\upsilon.
\]
In the former $U$ has no empty-type label; in the latter it still
has no universal label. In both cases $\abs U+\abs V>N$.

As in the proof of Theorem~\ref{thm:three-range-hall}, reduce the
weights of the selected labels and positions to have normalized
masses $t$ and $1-t$. In range 2 choose
$\tau\leq t\leq\gamma$; in range 3 choose
$\gamma\leq t\leq1-\upsilon$. Such a choice exists because
the untrimmed masses satisfy the corresponding lower bounds and
have sum greater than one. Write $x,w$ for the resulting type
submeasures.

Apply Lemma~\ref{lem:full-type-transport} first to $x\leq\xi$
with reference $\lambda$, and then to $w\leq\zeta$ with
reference $\nu$. This produces $\widehat x,\widehat w$ of the
same respective masses. Compatibility does not increase in either
step, so
\[
 \sum_{I,J}C(I,J)\widehat x_I\widehat w_J
 \leq\sum_{I,J}C(I,J)\widehat x_Iw_J
 \leq\sum_{I,J}C(I,J)x_Iw_J.
\]
The position coordinates satisfy
\eqref{eq:transport-position-caps}. In range 2 the empty label
coordinate remains zero. In range 3, before transport it is at
most $\xi_0-\upsilon$, since all universal labels are outside
$U$. The last assertion of Lemma~\ref{lem:full-type-transport}
therefore gives \eqref{eq:transport-high-caps}. Thus the transported
pair lies in the appropriate capped domain.

If that domain is infeasible, we already have a contradiction.
Otherwise its compatible mass is at least $\widehat\beta_k$.
Every compatible pair in the original weighted rectangle belongs
to $U\times V$ and is absent from $\Gamma$. Its weight is at
most one. Hence at least $N^2\widehat\beta_k$ relevant pairs
were deleted. In range 2 these pairs have nonempty label types;
transport has introduced no empty label mass. This contradicts
the assumed deletion bound. There is no deficient set, and Hall's
criterion gives a perfect matching.
\end{proof}

The two transports are comparisons between measures, not changes
to the labeling or to its random permutations. In particular, the
deletion variables still count the original compatible pairs.
If the local degree failures and these two deletion events have
total probability less than one, the union bound gives a labeling
to which Proposition~\ref{prop:two-sided-hall} applies.

\subsubsection{Exact finite errors for the reference laws}

We now compute \eqref{eq:two-sided-variation} on an integer interval
$[\ell,u]$. Use the preparation and construction of
Section~\ref{sec:interval-construction}, with $3z<\ell$ and
$v=\varepsilon_B/\rho$. The actual usable quota law remains
$\pi_I(n)=a_I(n)/M_n$. For blocks $P_1,\ldots,P_s$, put
\[
 \alpha_j=\prod_{p\in P_j}(1-1/p),\qquad
 \lambda_I=\prod_{j\in I}(1-\alpha_j)
                 \prod_{j\notin I}\alpha_j.
\]
This is the block product law. Suppose $\Omega<\lambda_0$,
where $\Omega=\max_n w_n/m_n$, and set
\begin{equation}\label{eq:finite-reference-position-law}
 \nu_I=\frac{\lambda_I-\Omega\boldsymbol1_{I=0}}{1-\Omega}.
\end{equation}
It is a probability law with positive empty coordinate.

Let $c_I(n)$ count the odd integers of type $I$ in $[n]$.
For $0\leq k\leq b(n)$ define
\[
 \xi_I(n,k)=
 \frac{c_I(n)-k\boldsymbol1_{I=0}}{\ceil{n/2}-k}.
\]
The preparation supplies at least $b(n)$ empty-type odd primes,
so these are probability laws with positive denominators. All
maxima in the following definition are over integer indices:
\begin{align}
 v_X&=\max_{\ell\leq n\leq u}\ 
       \max_{0\leq k\leq b(n)}
          \|\xi(n,k)-\lambda\|_{\mathrm{TV}},
                                      \label{eq:finite-label-transport-error}\\
 v_A&=\max_{\ell\leq n\leq u}
          \|\pi(n)-\nu\|_{\mathrm{TV}},\qquad
 v_B=v_A+v.                            \label{eq:finite-position-transport-error}
\end{align}

\begin{lemma}[Finite transport errors]
\label{lem:finite-transport-errors}
For every prepared tree of every order in $[\ell,u]$, its odd
type law and retained position-tag law satisfy
\eqref{eq:two-sided-variation} with the errors above. All reference
coordinates and errors are exactly computable rational numbers.
\end{lemma}

\begin{proof}
The odd labels removed for the core are $k\leq b(n)$ primes
greater than $n/2>z$. They all have empty type. Thus the actual
odd law is exactly $\xi(n,k)$, irrespective of which of these
primes were chosen. The first maximum gives its required bound.
Including $k=0$ only enlarges the family being bounded.

The construction gives
\[
 \|\zeta-\pi(n)\|_{\mathrm{TV}}
 \leq\frac{n\varepsilon_B}{N}\leq v.
\]
The triangle inequality and
\eqref{eq:finite-position-transport-error} imply
$\|\zeta-\nu\|_{\mathrm{TV}}\leq v_B$.

Finally, divisibility counts determine $a_I(n)$ and $c_I(n)$
exactly. Every other operation is a rational product, absolute
value, or maximum over the displayed finite integer sets. This
proves exact computability as well as the uniform bounds.
\end{proof}

The reference law $\nu$ is used only for the coordinate comparison.
The type assignment still has the actual quotas $a_I(n)$.
Likewise, the local degree bounds use the exact compatible fractions
in \eqref{eq:interval-initial-degrees}. No product-law approximation
is substituted for either of these data.

For each row define the two normalized deletion budgets
\begin{equation}\label{eq:transported-finite-budgets}
 \widehat b_2=\rho^2\widehat\beta_2,
 \qquad
 \widehat b_3=\rho^2\widehat\beta_3.
\end{equation}
Since $N\geq n\rho$, a strict deletion bound of
$n^2\widehat b_k$ suffices in
Proposition~\ref{prop:two-sided-hall}. Formula
\eqref{eq:finite-hall-failure}, with $\widehat b_k$ in place of
$b_k$, bounds its failure probability. The middle event counts
only nonempty odd types, whereas the high event counts all types.
Both use the moment of the actual construction, including its
donor term or completion credit as appropriate.

\subsubsection{The complete finite certificate}

Table~\ref{tab:exponential-continuation-cover} gives the complete
order cover. Every entry contains one row for each integer degree
$1\leq D\leq D_*$, where $D_*=\ceil{1/(2W)}$.
There are 35 adjacent intervals and 622 rows. Each row uses two
blocks of small primes, hence four types, and has
$\eta=\tau<\gamma$. The first Hall range is empty; both remaining
ranges use their exponential bounds.

\begin{table}[!htbp]
\centering
\caption{The degree cover for the two-sided transport certificate.
Every degree from $1$ through the displayed $D_*$ is included.
A singleton denotes one integer order.}
\label{tab:exponential-continuation-cover}
\par\vspace{.5\baselineskip}
\renewcommand{\arraystretch}{1.15}
\begin{tabular}{rrrr}
\hline
Order interval&$D_*$&Order interval&$D_*$\\
\hline
$1754$--$1759$&16&$2008$&18\\
$1760$--$1765$&17&$2009$--$2011$&18\\
$1766$--$1771$&16&$2012$--$2014$&18\\
$1772$--$1796$&17&$2015$--$2017$&18\\
$1797$--$1808$&17&$2018$&18\\
$1809$--$1833$&17&$2019$&18\\
$1834$--$1858$&18&$2020$&18\\
$1859$--$1883$&18&$2021$--$2023$&18\\
$1884$--$1895$&18&$2024$--$2029$&18\\
$1896$--$1920$&18&$2030$--$2035$&18\\
$1921$--$1945$&18&$2036$--$2047$&18\\
$1946$--$1970$&18&$2048$--$2059$&18\\
$1971$--$1995$&18&$2060$--$2084$&18\\
$1996$&18&$2085$--$2109$&18\\
$1997$&18&$2110$--$2134$&18\\
$1998$--$2000$&18&$2135$--$2184$&18\\
$2001$--$2006$&18&$2185$--$2209$&18\\
$2007$&18&&\\
\hline
\end{tabular}
\end{table}

\begin{proposition}[Verified two-sided transport certificate]
\label{prop:verified-exponential-continuation}
Every row in Table~\ref{tab:exponential-continuation-cover}
supplies valid preparation and construction bounds, the local failure
estimates of Section~\ref{sec:local-certificates}, and exponential
failure bounds for both deletion events required by
Proposition~\ref{prop:two-sided-hall}. Its exceptional fraction
is strictly less than $\tau$, both budgets in
\eqref{eq:transported-finite-budgets} are positive, and
\begin{equation}\label{eq:exponential-continuation-joint}
 p_{\rm loc}+p_2+p_3<1.
\end{equation}
The largest sum is attained on $[1834,1858]$, at residual degree $8$.
\end{proposition}

\begin{proof}
The accompanying finite checker recomputes the complete construction,
local and Hall data. We specify how these checks establish each
hypothesis.

\emph{Preparation and construction.}
At all 456 integer orders, direct prime counts verify the core
supply, special-label counts and universal odd-label supply.
They also give the degree bound $D_*$. The interval endpoints are
adjacent and every degree list is complete. Each actual order--degree
pair is tested against \eqref{eq:interval-degree-feasibility}
before using the paid donor stock.

The rows use 70 original-forest constructions with cuts confined
to $B$, 399 with donor exchanges, and 153 connected height
constructions. For each row the checker recomputes the separator
or height budget, degree variation, transfer, tag error and side
bounds of Proposition~\ref{prop:interval-construction-bounds}.
It checks completion capacity wherever a completing tree is used.
All supplied upper and lower bounds have the required directions,
and all endpoint quota lower bounds are positive.

\emph{Reference laws and local degrees.}
The exact usable quota laws give the initial degree bounds.
Equations \eqref{eq:finite-label-transport-error} and
\eqref{eq:finite-position-transport-error} give the two transport
errors, including every core size in their finite maxima.
The comparison $\Omega<\lambda_0$ holds in every row.
An independent computation obtains the types by gcd tests, forms
the odd laws after removing $k$ empty-type labels, and sums the
positive coordinate differences. It gives the same maxima as
one half of the sum of absolute differences.

There are 531 rows using the complete full-type exception pool
and 91 using the protected pool of
\eqref{eq:protected-local-exceptions}. The latter include 76 donor
rows and 15 connected rows. Each omitted full-type label passes
the strict deterministic degree test. In every row the exception
fraction is less than $\tau$, all local margins are positive, and
$\gamma+\upsilon\leq1$.

The maximal-support cover has 256,870 event occurrences and
covers 616,062 endpoint label occurrences across the rows.
For every event, the checker reconstructs its class frequencies,
sample bounds, integer threshold and selected probability estimates
from Theorem~\ref{thm:finite-local-degree-certificate}.
Exactly 224,210 events have a deterministic zero bound. Summing
over the maximal events gives $p_{\rm loc}$, retaining the
multiplicity of different events with equal numerical bounds.
The independent arithmetic calculation also checks each profile
against the actual usable labels at every order of its interval.

\emph{The two Hall events.}
The checker forms the transported caps in
\eqref{eq:transport-position-caps}--\eqref{eq:transport-high-caps}.
It computes both minima using
Lemma~\ref{lem:fixed-mass-concavity}; a separate enumeration of
the vertices in Proposition~\ref{prop:capped-bilinear-minimum}
gives the same 1,244 positive minima. This comparison includes
the zero empty coordinate in the middle range and the reduced
empty coordinate in the high range.

For each range the endpoint cost histograms count precisely its
allowed odd-type family. They are independently reconstructed by
counting odd multiples of each larger prime. The directed moment
calculation then applies the bounds of
Proposition~\ref{prop:directed-moment-evaluation}, including the
special-label subtractions and the actual quota-weighted average.
There are two exponential evaluations in every row. The local
estimates and these two bounds concern one and the same random
permutation of the even classes.

All 622 comparisons in \eqref{eq:exponential-continuation-joint}
use exact rational arithmetic. The largest sum is approximately
$0.9270146969171777$. Its local, middle and high contributions
are about $0.007607$, $0.002053$ and $0.917355$, respectively.
The finite record supplies their exact rational values and every
accepted input. The rounded displays play no role in acceptance.
This proves the asserted conditions.
\end{proof}

\begin{theorem}[The initial interval through 2209]
\label{thm:initial-prefix-2209}
Every tree of order at most $2209$ has a prime labeling.
\end{theorem}

\begin{proof}
Theorem~\ref{thm:initial-prefix-1753} covers the orders at most
$1753$. Fix a tree $T$ of order $n\in[1754,2209]$, and select
its interval in Table~\ref{tab:exponential-continuation-cover}.
Prepare it by Proposition~\ref{prop:highest-degree-preparation}.
An edgeless residual forest is labeled immediately. Otherwise
Lemma~\ref{lem:interval-degree-cover} selects a row from the
complete degree list, and its actual order satisfies the feasibility
condition used by that row.

Fix the type assignment and tags before permuting the labels within
the even classes. Proposition~\ref{prop:verified-exponential-continuation}
and \eqref{eq:exponential-continuation-joint} give positive
probability that all local degree conditions hold and both deletion
counts are strictly below $n^2\widehat b_k$. Lemma
\ref{lem:finite-transport-errors} supplies the two variation bounds.
Since $n^2\widehat b_k\leq N^2\widehat\beta_k$,
Proposition~\ref{prop:two-sided-hall} gives a perfect matching
between the remaining odd labels and positions.

Place the matched labels, restore the special labels, and apply
the repair in Proposition~\ref{prop:highest-degree-preparation}.
This gives a prime labeling of $T$. The intervals cover every
integer from $1754=1753+1$ through $2209$, so the conclusion
joins the preceding theorem without a gap.
\end{proof}

Section~\ref{sec:single-range-continuation} continues from order
$2210$ by choosing equal degree targets, leaving a single Hall
range to check.

\subsection{A single Hall range through order 3499}
\label{sec:single-range-continuation}

The next intervals admit a common degree target for all labels outside
the exceptional full-type pool. The first two Hall ranges are then
empty. We retain the two-sided transport, but need only one bound
for the total number of deleted pairs.

\subsubsection{Coincident degree targets}

Use the notation of Section~\ref{sec:exponential-continuation}.
Thus $\Gamma\subseteq\Gamma_0$ has parts $Y,B$ of size $N$;
$\xi,\zeta$ are their actual type and tag laws; and $\lambda,\nu$
are reference probability laws satisfying
\eqref{eq:two-sided-variation}. Let
\[
 \widehat p=\widehat p^{(3)},\qquad
 \widehat\kappa_I=\nu_I+v_B\boldsymbol1_{I=F},
\]
where $\widehat p^{(3)}$ is given by
\eqref{eq:transport-high-caps}. These caps depend on the omitted
universal fraction $\upsilon$, but not on the choice of degree targets.

\begin{lemma}[A single transported Hall bound]
\label{lem:single-range-hall}
Suppose $0\leq\varepsilon<\tau$, $\upsilon>0$, and
$\tau+\upsilon\leq1$. Let $E\subseteq Y$ consist of full-type
labels, with $\abs E\leq\varepsilon N$. Assume
\[
 d_\Gamma(y)>\varepsilon N\quad(y\in E),\qquad
 d_\Gamma(y)>\tau N\quad(y\notin E),
\]
and suppose at least $\upsilon N$ empty-type labels are universal.
Assume also \eqref{eq:two-sided-variation}. If
\begin{equation}\label{eq:single-range-minimum}
 \widehat\beta=\mathcal B(\widehat p,\widehat\kappa;
                                      \tau,\upsilon)>0
\end{equation}
on a nonempty domain, and fewer than $N^2\widehat\beta$ pairs
of $\Gamma_0$ are deleted, then $\Gamma$ has a perfect matching.
If the domain is empty, no deletion bound is needed.
\end{lemma}

\begin{proof}
Suppose $U\subseteq Y$ is deficient and put
$V=B\setminus N_\Gamma(U)$. Lemma~\ref{lem:three-hall-ranges},
with $\eta=\gamma=\tau$, gives
\[
 \tau<\abs U/N\leq1-\upsilon,\qquad
 \abs V/N>\upsilon,\qquad \abs U+\abs V>N.
\]
Also $U$ contains no universal label. Choose
\[
 t\in[\tau,1-\upsilon]\cap
           [1-\abs V/N,\abs U/N].
\]
The displayed inequalities show that this intersection is nonempty.
Reduce the selected weights to masses $t$ and $1-t$, obtaining
submeasures $x\leq\xi$ and $w\leq\zeta$.

Apply Lemma~\ref{lem:full-type-transport} to each submeasure.
The resulting $\widehat x,\widehat w$ preserve both masses,
satisfy $\widehat x\leq\widehat p$ and
$\widehat w\leq\widehat\kappa$, and have no greater compatibility.
The smaller empty label cap is valid because $U$ omits all the
universal labels. Thus the domain in
\eqref{eq:single-range-minimum} is nonempty and
\[
 \widehat\beta
 \leq\sum_{I,J}C(I,J)\widehat x_I\widehat w_J
 \leq\sum_{I,J}C(I,J)x_Iw_J.
\]
Every compatible pair between $U$ and $V$ has been deleted.
Each selected object has weight at most one before normalization
by $N$, so the last sum is at most the number of deleted pairs
divided by $N^2$. This contradicts the assumed strict bound.
Hall's theorem completes the proof.
\end{proof}

The equality of the targets removes two ranges of possible deficient
sets. It does not remove the local degree conditions. In particular,
each member of the exception pool must still have degree strictly
greater than $\varepsilon N$.

\subsubsection{Two estimates for the same deletion event}

Fix one construction row on $[\ell,u]$, with $N\geq n\rho$,
and put
\begin{equation}\label{eq:single-range-budget}
 \widehat b=\rho^2\widehat\beta.
\end{equation}
Use all odd types in the global cost, so $\mathcal J=\mathcal T$.
Let $S$ be the nonnegative sum charging deleted pairs along the
original forest, as in Section~\ref{sec:matching-estimates}.
The number of deleted pairs is at most $S$. The event
$S<n^2\widehat b$ therefore supplies the strict bound in
Lemma~\ref{lem:single-range-hall}.

The exact first-moment data of
Proposition~\ref{prop:finite-global-means} give
\begin{equation}\label{eq:single-range-mean}
 \E S/n^2\leq M_*
 :=e_+A_*+H_+B_*+\delta C_*-q_-\min_I L_I.
\end{equation}
For the original-forest constructions $q_-=0$. For a connected
height construction $\delta=0$, and the last term is the credit
for the edges added in the completion. These are the two separate
bounds \eqref{eq:finite-donor-global-mean} and
\eqref{eq:finite-completed-global-mean}; the notation combines
them only by inserting a zero term.

Let $t>0$ be the row's rational exponential parameter, and let
$\widehat{\mathcal M}$ be its directed global moment for
$\mathcal T$. Define
\begin{equation}\label{eq:single-range-probability}
 p_{\rm glob}=
 \min\left\{1,\frac{M_*}{\widehat b},
 \mathcal P\!\left(\frac\ell D
              (t\widehat b-\widehat{\mathcal M})\right)\right\}.
\end{equation}
Here $M_*\geq0$ and $\widehat b>0$ are checked in each row.
Markov's inequality gives the second bound for
$\Pr(S\geq n^2\widehat b)$. Proposition
\ref{prop:directed-moment-evaluation} gives the third.
Their minimum bounds the same event, and it is added only once
to the local failure sum. The acceptance condition is
\begin{equation}\label{eq:single-range-joint}
 p_{\rm loc}+p_{\rm glob}<1.
\end{equation}
No independence between the two events is required. Parameter
selection precedes this exact comparison and makes no claim
that the chosen parameter is optimal.

\subsubsection{The complete order and degree cover}

Table~\ref{tab:single-range-cover} gives 89 adjacent intervals
covering the 1,290 integer orders from $2210$ through $3499$.
Each interval contains one row for every degree
$1\leq D\leq D_*=\ceil{1/(2W)}$, for a total of 1,631 rows.
Every row uses two blocks of small primes and has
$\tau=\eta=\gamma$. The actual quota laws still determine the
initial compatible fractions. The reference laws and their
finite errors are used only in the transported Hall comparison.

\begin{table}[!htbp]
\centering
\caption{The complete single-range degree cover. Read each pair of
columns downwards. Every degree from $1$ through $D_*$ is included.}
\label{tab:single-range-cover}
\par\vspace{.5\baselineskip}
\renewcommand{\arraystretch}{1.15}
\begin{tabular}{rrrrrr}
\hline
Order interval&$D_*$&Order interval&$D_*$&Order interval&$D_*$\\
\hline
$2210$&17&$2331$--$2334$&18&$2591$&19\\
$2211$--$2213$&17&$2335$--$2338$&18&$2592$--$2593$&19\\
$2214$&17&$2339$&18&$2594$--$2597$&19\\
$2215$--$2216$&17&$2340$&18&$2598$--$2601$&19\\
$2217$--$2218$&18&$2341$&18&$2602$--$2605$&18\\
$2219$&17&$2342$--$2343$&18&$2606$--$2623$&18\\
$2220$&17&$2344$--$2352$&18&$2624$--$2641$&18\\
$2221$--$2222$&17&$2353$--$2356$&18&$2642$--$2659$&18\\
$2223$--$2225$&17&$2357$--$2360$&18&$2660$--$2696$&19\\
$2226$--$2228$&17&$2361$&18&$2697$--$2733$&19\\
$2229$&18&$2362$--$2363$&18&$2734$--$2770$&19\\
$2230$&18&$2364$--$2367$&18&$2771$--$2807$&19\\
$2231$--$2232$&18&$2368$--$2369$&18&$2808$--$2844$&20\\
$2233$--$2234$&18&$2370$&18&$2845$&20\\
$2235$--$2236$&17&$2371$--$2372$&18&$2846$&20\\
$2237$--$2241$&17&$2373$--$2376$&18&$2847$&20\\
$2242$--$2247$&18&$2377$--$2378$&18&$2848$--$2850$&20\\
$2248$--$2253$&18&$2379$--$2382$&18&$2851$--$2856$&20\\
$2254$--$2260$&18&$2383$--$2391$&18&$2857$--$2862$&20\\
$2261$--$2264$&17&$2392$--$2400$&18&$2863$--$2874$&20\\
$2265$&17&$2401$--$2409$&18&$2875$--$2899$&19\\
$2266$--$2269$&17&$2410$--$2427$&18&$2900$--$2949$&19\\
$2270$--$2273$&18&$2428$--$2445$&18&$2950$--$2999$&19\\
$2274$--$2277$&18&$2446$--$2463$&18&$3000$--$3049$&19\\
$2278$--$2281$&18&$2464$--$2481$&18&$3050$--$3099$&20\\
$2282$--$2285$&18&$2482$--$2499$&19&$3100$--$3199$&20\\
$2286$--$2294$&18&$2500$--$2517$&19&$3200$--$3299$&20\\
$2295$--$2312$&18&$2518$--$2535$&19&$3300$--$3399$&20\\
$2313$--$2321$&18&$2536$--$2572$&19&$3400$--$3499$&20\\
$2322$--$2330$&18&$2573$--$2590$&19&&\\
\hline
\end{tabular}
\end{table}

\begin{proposition}[Verified single-range certificate]
\label{prop:verified-single-range-continuation}
Every row in Table~\ref{tab:single-range-cover} satisfies the
preparation, construction and local degree conditions, with
$\varepsilon<\tau$ and $\tau+\upsilon\leq1$. Its transported
budget in \eqref{eq:single-range-budget} is positive, and
\eqref{eq:single-range-joint} holds.
\end{proposition}

\begin{proof}
The finite checker recomputes every construction, local and global
bound from the supplied rational parameters. It checks interval
adjacency and complete degree lists. All comparisons are exact;
the decimals below only describe their sizes.

\emph{Preparation and assignment.}
Direct prime counts give the core supply, removed even family,
remaining side sizes and universal odd-label supply for every
integer order. The bound $D_*$ and the necessary feasibility test
\eqref{eq:interval-degree-feasibility} are evaluated before any
row is applied. For connected rows, completion capacity and the
added-edge credit are checked for every permitted core size.

The table uses 178 constructions with cuts confined to $B$,
1,146 original-forest donor constructions and 307 connected height
constructions. Proposition~\ref{prop:interval-construction-bounds}
supplies their degree variation, transfer and tag error bounds.
Every supplied upper or lower bound is compared in the required
direction. The exact usable class quotas are positive.

\emph{Type laws and local bounds.}
Lemma~\ref{lem:finite-transport-errors} gives the reference laws
and their errors. The finite maxima include every order of the
interval and every permitted number of core labels. A separate
calculation obtains types by gcd tests and total variation by
summing positive coordinate differences. Its maxima agree with
the half-sums of absolute differences used in the primary evaluator.

All rows use the complete full-type exception pool. Its normalized
size is less than $\tau$. The other labels have the common target
$\tau$, while the exception target is $\varepsilon$. The checker
forms every maximal prime-support event, including those with no
possible conflicts. It verifies each local profile against the
actual usable populations at every order, rounds the integer
threshold in the direction proved in
Section~\ref{sec:directed-probability}, and sums the selected
bounds of Theorem~\ref{thm:finite-local-degree-certificate}.
Distinct maximal events retain their multiplicities even when
they have the same numerical profile. The resulting 841,306 event
occurrences cover 2,052,537 endpoint label occurrences. Exactly
795,133 events have a deterministic zero bound.

\emph{The global bound and simultaneous success.}
The checker forms $\widehat p,\widehat\kappa$ and evaluates
\eqref{eq:single-range-minimum} using
Lemma~\ref{lem:fixed-mass-concavity}. Independent vertex enumeration
from Proposition~\ref{prop:capped-bilinear-minimum} gives the
same 1,631 positive minima. Thus no Hall range is omitted by
assuming that its domain is empty.

The endpoint cost histograms are independently reconstructed by
counting compatible odd multiples of each larger prime. Their
sums, together with the exact removed-label counts, reproduce
$L_I,U_I,A_*,B_*,C_*$ in
Proposition~\ref{prop:finite-global-means}. These give $M_*$ in
\eqref{eq:single-range-mean}. The same cost histograms enter the
directed exponential calculation, which uses the actual
quota-weighted mean of the class logarithmic moments.
Both calculations include the donor term or completion credit
appropriate to the construction.

The checker evaluates both global bounds in every row. The first-moment
bound is smaller in 13 rows and the exponential bound in 1,618 rows.
All 1,631 rational comparisons in \eqref{eq:single-range-joint}
are strict. The largest sum occurs on $[2383,2391]$ at degree $5$;
it is approximately $0.9866043684687599$. Its local and global terms
are about $3.88024\cdot10^{-9}$ and $0.9866043646$.

\end{proof}

\begin{theorem}[The initial interval through 3499]
\label{thm:initial-prefix-3499}
Every tree of order at most $3499$ has a prime labeling.
\end{theorem}

\begin{proof}
Theorem~\ref{thm:initial-prefix-2209} proves the assertion through
$2209$. Let $T$ have order $n\in[2210,3499]$, and choose its
interval in Table~\ref{tab:single-range-cover}.
Proposition~\ref{prop:highest-degree-preparation} prepares the
residual bipartite forest. If it is edgeless, label it directly.
Otherwise Lemma~\ref{lem:interval-degree-cover} places its
maximum degree in the complete row list, and its order satisfies
the row's necessary feasibility condition.

Fix the type assignment and tags furnished by that row, then
permute the even labels within their actual classes.
Proposition~\ref{prop:verified-single-range-continuation} and
\eqref{eq:single-range-joint} give positive probability that all
local degree conditions hold and $S<n^2\widehat b$.
Lemma~\ref{lem:finite-transport-errors} supplies the two variation
bounds. Since $n^2\widehat b\leq N^2\widehat\beta$,
Lemma~\ref{lem:single-range-hall} gives a perfect matching
between the remaining odd labels and positions.

Place those labels, restore the removed special labels, and
apply the repair of Proposition~\ref{prop:highest-degree-preparation}.
This gives a prime labeling of $T$. The intervals join at
$2210=2209+1$ and cover every subsequent integer through $3499$.
\end{proof}

Section~\ref{sec:retained-continuation} continues from order
$3500$ with the same Hall criterion and an explicit adjustment
of the common degree target.

\subsection{Retained positions and the orders through 34999}
\label{sec:retained-continuation}

From order $3500$ through $34999$ we use the three constructions of
Proposition~\ref{prop:interval-construction-bounds} and the single Hall
range of Lemma~\ref{lem:single-range-hall}, and we retain all remaining
positions. At an exceptional position, its tag contains every
neighboring even type; the resulting variation is paid before any
probability bound is applied. Seven certificates cover this range. They
use the same criterion and differ only in their interval covers, prime
blocks, and parameters. We state the criterion once
(Theorem~\ref{thm:retained-interval-criterion}) and then record the
certificates (Proposition~\ref{prop:verified-retained-continuations}).

\subsubsection{Degree targets and retained tags}

For one row on $[\ell,u]$, take the complete full-type pool
\[
 \mathcal E=\{y\leq u:y\text{ is odd and }T(y)=F\},
 \qquad \varepsilon=\frac{\abs{\mathcal E}}{\ell\rho}.
\]
Since the residual odd side has size $N\geq n\rho\geq\ell\rho$,
this bounds the actual exceptional fraction at every order.
Using the endpoint count $c_F(u)$ from
\eqref{eq:odd-endpoint-caps}, choose a rational parameter $\tau_0$
and set
\begin{equation}\label{eq:retained-continuation-target}
 \tau=\eta=\gamma
 =\max\left\{\tau_0,\frac{c_F(u)}{\ell\rho}+10^{-6}\right\}.
\end{equation}
Thus $\varepsilon<\tau$. Increasing this common target reduces
the domain of the remaining Hall minimum, but also strengthens
the required degrees outside $\mathcal E$. Both effects are
included in the acceptance test; the target choice alone does
not establish either degree or matching conditions.

For each row, take $\rho,\upsilon$ from
\eqref{eq:interval-preparation-constants}, and take
$e_+,H_+,\delta,q_-,\varepsilon_B$ to be the expressions in the chosen
case of Proposition~\ref{prop:interval-construction-bounds}. In a
connected row, the completion credit is $W(D-1)$ and the degree
variation is measured in the completing tree.

The construction bounds retain their separate meanings.
For cuts confined to $B$, the changed tags cost at most $c_Bn$.
For donor exchanges, their number is bounded by $z_*$ in
\eqref{eq:interval-donor-boundary}. This bound uses the full
sequence of paid stocks $\ell_t(n)$, not a single common degree
cap for all donors. In particular, it also applies when the
available degree-one donors do not suffice. For connected height
assignments, the boundary positions receive full tags and their
number is bounded by $c_H$.

In all three cases $I(a)\subseteq J(b)$ on every original edge.
Consequently Lemma~\ref{lem:retained-tags} bounds the actual
large-prime deletions by the same local and global conflict costs.
Retagging changes neither the label quotas nor the graph edges.
Original-forest rows pay the donor transfer in their means; only
connected rows subtract the completion credit. All assignments
and tags are fixed before the within-class permutations.

Use the exact quota laws for the local compatible fractions.
For the Hall comparison, take the reference laws and variation
bounds of Lemma~\ref{lem:finite-transport-errors}. With the caps
$\widehat p,\widehat\kappa$ from
Lemma~\ref{lem:single-range-hall}, put
\begin{equation}\label{eq:retained-continuation-budget}
 \widehat b=\rho^2\mathcal B(
             \widehat p,\widehat\kappa;\tau,\upsilon).
\end{equation}
The global bound $p_{\rm glob}$ is \eqref{eq:single-range-probability},
with parameter $t=10$ and all odd types. The local bound $p_{\rm loc}$
is the sum \eqref{eq:finite-local-failure} of
Theorem~\ref{thm:finite-local-degree-certificate}. As in
Section~\ref{sec:single-range-continuation}, the first-moment and
exponential estimates in \eqref{eq:single-range-probability} bound the
same global event; their minimum is added once to the local sum.

\subsubsection{The criterion}

\begin{theorem}[A criterion with retained positions]
\label{thm:retained-interval-criterion}
Let $[\ell,u]$ be an interval of orders satisfying
\eqref{eq:interval-preparation-tests}, and put $D_*=\ceil{1/(2W)}$.
Suppose that for every integer $1\leq D\leq D_*$ a row is specified:
one of the three constructions of
Proposition~\ref{prop:interval-construction-bounds} with its rational
parameters, a partition of the odd primes up to a cutoff into two prime
blocks, and a rational parameter $\tau_0$. Define the target by
\eqref{eq:retained-continuation-target}, $\widehat b$ by
\eqref{eq:retained-continuation-budget}, and $p_{\rm loc}$ and
$p_{\rm glob}$ as above. Suppose that every row satisfies, for all
orders $n\in[\ell,u]$ permitted by
\eqref{eq:interval-degree-feasibility}:
\begin{enumerate}
\item the preparation and construction conditions of
  Section~\ref{sec:interval-construction};
\item positive local degree margins at the target
  \eqref{eq:retained-continuation-target};
\item $\tau+\upsilon\leq1$ and $\widehat b>0$;
\item \begin{equation}\label{eq:retained-continuation-joint}
   p_{\rm loc}+p_{\rm glob}<1.
  \end{equation}
\end{enumerate}
Then every tree of every order in $[\ell,u]$ has a prime labeling.
\end{theorem}

\begin{proof}
Let $T$ have order $n\in[\ell,u]$ and prepare it by
Proposition~\ref{prop:highest-degree-preparation}. If the residual
forest is edgeless, assign its labels arbitrarily and finish with the
restoration and repair of that proposition. Otherwise
Lemma~\ref{lem:interval-degree-cover} places the residual even-side
degree $D$ in $\{1,\ldots,D_*\}$, and $n$ satisfies the feasibility
condition \eqref{eq:interval-degree-feasibility} of the row for $D$.

Fix the row's type assignment and retained tags before taking the
class permutations. The target \eqref{eq:retained-continuation-target}
gives $\varepsilon<\tau=\eta=\gamma$, and $\tau+\upsilon\leq1$ by the
third condition. Hence Proposition~\ref{prop:local-exception-pool}
applies to the full-type pool $\mathcal E$, and
$\abs{\mathcal E\cap Y}\leq\varepsilon N$. By the first three
conditions and Theorem~\ref{thm:finite-local-degree-certificate}, some
local degree condition fails with probability at most $p_{\rm loc}$.
By the construction of \eqref{eq:single-range-probability}, the global
conflict cost is at least $n^2\widehat b$ with probability at most
$p_{\rm glob}$. By \eqref{eq:retained-continuation-joint}, some choice
of the class permutations satisfies every local degree condition and
makes the global conflict cost less than $n^2\widehat b$.

For this choice we verify the hypotheses of
Lemma~\ref{lem:single-range-hall}. The avoided local events give its
degree conditions (Proposition~\ref{prop:local-exception-pool}).
Proposition~\ref{prop:interval-preparation} supplies at least
$\upsilon N$ universal labels of empty type, since the cutoff satisfies
$3z<\ell$. Lemma~\ref{lem:finite-transport-errors} gives the variation
bounds, and $\widehat b>0$ is the positivity of the rectangle minimum.
By Lemma~\ref{lem:retained-tags}, the number of deleted pairs is at
most the global conflict cost, which is less than
$n^2\widehat b=(n\rho)^2\widehat\beta\leq N^2\widehat\beta$, with
$\widehat\beta$ as in \eqref{eq:single-range-minimum}, because
$N\geq n\rho$. The lemma therefore supplies a perfect matching between
the remaining odd labels and positions.
Assign the matched labels, restore the special labels, and apply the
repair of Proposition~\ref{prop:highest-degree-preparation}. The
resulting labeling of $T$ is prime.
\end{proof}

\subsubsection{The certified intervals}

Tables~\ref{tab:retained-certificates}
and~\ref{tab:retained-largest-sums} list the seven certificates.
Their interval covers are given in Appendix~\ref{sec:cover-tables};
each interval carries one row for every degree $1,\ldots,D_*$.
The rows are of three kinds: original-forest constructions with cuts
confined to $B$, constructions with donor exchanges, and connected
height constructions.

\begin{table}[ht]
\centering
\caption{The retained-position certificates: the number of intervals
(Int.), the number of rows, as the total and the numbers of cut, donor,
and connected rows, and the number of maximal-support events of the
local sum (Events). The last columns give the number of these events
with a deterministic zero bound (Zero) and the number of endpoint label
occurrences covered by all events (Occ.).}
\label{tab:retained-certificates}
\par\vspace{.5\baselineskip}
\renewcommand{\arraystretch}{1.15}
\begin{tabular}{lrrrrr}
\hline
Orders&Int.&Rows&Events&Zero&Occ.\\
\hline
$3500$--$5219$&41&$846=82+641+123$&683727&646125&1671387\\
$5220$--$6771$&16&$359=32+294+33$&439818&412915&1090138\\
$6772$--$12499$&23&$555=46+463+46$&1104988&1030489&2724580\\
$12500$--$14999$&5&$129=10+109+10$&370855&342583&903500\\
$15000$--$19999$&5&$133=10+113+10$&496222&455137&1198500\\
$20000$--$24999$&2&$55=4+47+4$&272849&246848&653750\\
$25000$--$34999$&2&$58=4+50+4$&395270&350518&942500\\
\hline
\end{tabular}
\end{table}

\begin{table}[ht]
\centering
\caption{The largest joint sum $p_{\rm loc}+p_{\rm glob}$ of each
retained-position certificate, with the interval and degree where it
occurs. The decimals describe exact rational comparisons.}
\label{tab:retained-largest-sums}
\par\vspace{.5\baselineskip}
\renewcommand{\arraystretch}{1.15}
\begin{tabular}{lllr}
\hline
Orders&Interval&Degree&Largest joint sum, approximately\\
\hline
$3500$--$5219$&$[4220,4319]$&$7$&$0.9490722798465094$\\
$5220$--$6771$&$[5220,5319]$&$16$&$0.04302850038173173$\\
$6772$--$12499$&$[6772,7021]$&$15$&$0.005192492891846186$\\
$12500$--$14999$&$[12500,12999]$&$19$&$0.0005138912$\\
$15000$--$19999$&$[15000,15999]$&$20$&$0.0005116258$\\
$20000$--$24999$&$[20000,22499]$&$4$&$0.006724511$\\
$25000$--$34999$&$[25000,29999]$&$23$&$0.002633058$\\
\hline
\end{tabular}
\end{table}

\begin{proposition}[The retained-position certificates]
\label{prop:verified-retained-continuations}
For each certificate in Table~\ref{tab:retained-certificates}, each
interval of its cover in Appendix~\ref{sec:cover-tables}, and each
degree $1\leq D\leq D_*$, the specified row satisfies the four
conditions of Theorem~\ref{thm:retained-interval-criterion}.
\end{proposition}

\begin{proof}
The finite checker of each certificate reconstructs all hypotheses from
the supplied rational parameters. It accepts the whole order and degree
cover, not a subcollection of successful rows. The same steps are
carried out for every certificate.

\emph{Cover and preparation.} The intervals of each cover are adjacent,
and the degree list of every interval is exactly
$1,\ldots,\ceil{1/(2W)}$. Direct prime counts verify the preparation
at every integer order. Only \eqref{eq:interval-degree-feasibility}
excludes an order--degree pair, so every possible residual degree is
represented.

\emph{Constructions.} For each row, the checker recomputes the cut or
height budget, rounding variation, transfer, changed-tag bound, and
residual side bounds of
Proposition~\ref{prop:interval-construction-bounds} on the whole
interval and checks every supplied bound against its defining
expression. The
donor boundary uses the complete sequence of paid stocks at each
feasible order, so the test does not assume that the chosen donors have
degree at most one. In a connected row, the degree variation is
measured in the completing tree, and the completion capacity and
added-edge credit are checked for every permitted core size.

\emph{Types and targets.} Exact usable quotas give the compatible
fractions for the local degree estimates. The two finite type errors of
Lemma~\ref{lem:finite-transport-errors} are maxima over every order of
the interval and every permitted number of removed core labels. The
endpoint full-type count gives
$\varepsilon$, and the checker verifies $\varepsilon<\tau=\eta=\gamma$
and $\tau+\upsilon\leq1$. Thus only the last Hall range remains, while
the full-type labels keep their weak local degree condition. When
\eqref{eq:retained-continuation-target} raises the base target, the
local bounds are evaluated at the raised value.

\emph{Local sum.} The endpoint event cover is formed anew on each
interval, and Lemma~\ref{lem:endpoint-local-cover} covers every smaller
order: an actual label has the same type as its assigned endpoint event
and a support contained in that event's support. The representing
endpoint label need not itself occur at the actual order. For each
maximal-support event, the checker reconstructs its arithmetic profile.
For each type threshold
$\theta$, the checker computes $\min_{\ell\leq n\leq u}\ceil{n\theta}/n$
over all integer orders. The probability bounds are evaluated as in
Sections~\ref{sec:directed-probability}
and~\ref{sec:local-certificates}. Distinct events with the same
numerical bound keep their multiplicity in the local sum, and all local
margins are positive.

\emph{Global bound.} The checker reconstructs the endpoint cost
histograms and computes both bounds in
\eqref{eq:single-range-probability}. The directed exponential
calculation uses the same cost family and the actual quota-weighted
average; both estimates retain the donor term or the completion credit.
The transported rectangle minimum is positive in every row. The smaller
of the two global bounds is added once to the local sum, and all joint
comparisons \eqref{eq:retained-continuation-joint} are strict. All
comparisons use exact rational arithmetic.

\emph{Independent checks.} Each checker also runs reference
calculations that test its implementation; the proposition does not
depend on them. For every certificate, they compare each arithmetic
profile with the actual usable labels at every order of the interval,
recompute the endpoint cost histograms by direct counts of compatible
odd multiples, and obtain every rectangle minimum by a separate
enumeration of the rectangle vertices, using
Proposition~\ref{prop:capped-bilinear-minimum}. From order $6772$ on,
the population count is organized as follows. For a tail-prime support
$P$, list the half-labels $j$ of each type with
$\gcd\bigl(j,\prod_{p\in P}p\bigr)>1$; restricting to
$j\leq\floor{n/2}$ gives the complete population at order $n$. Powers
of two contribute no such half-label, and a removed odd prime
contributes precisely when it belongs to $P$. Subtracting these
contributions gives the usable population. From order $20000$ on,
independent factorization also reconstructs the endpoint covers with
their maximal supports and assigned label sets. For the certificates
through order $12499$, independent gcd-based counts reproduce the two
type errors, the odd reference law, and the usable even quota laws, and
the direct cost counts, together with the removed-label counts, give
independent values for the first-moment data of
Proposition~\ref{prop:finite-global-means}. The supplement lists these
checks for each certificate.

The seven certificates also involve the following particular features.
\begin{enumerate}
\item On $3500\leq n\leq5219$, the blocks partition the odd primes
  through one of $13$, $17$, $19$ into two sets, giving four types. The
  first-moment bound is smaller in $13$ rows and the exponential bound
  in the other $833$. In the largest joint sum, the local term is
  about $1.30731\cdot10^{-11}$ and the global term about
  $0.9490722798334362$.
\item On $5220\leq n\leq6771$, all rows have two prime blocks and hence
  four types. Universality of the special odd labels is checked by
  listing their possible even multiples: for a prime $p>n/3$ the only
  possible one in $[n]$ is $2p$, and when it is present it belongs to
  the removed special family; the label $1$ is universal directly.
\item On $6772\leq n\leq12499$, two prime blocks suffice throughout. In
  the degree-$13$ row on $[7272,7521]$ the target is raised by
  \eqref{eq:retained-continuation-target}.
\item On $12500\leq n\leq24999$, degree $4$ uses the blocks
  $\{3,13\}$ and $\{5,7,11\}$ with cutoff $13$; the other degrees use
  $\{3,13,19\}$ and $\{5,7,11,17\}$ with cutoff $19$. On
  $25000\leq n\leq34999$ the latter blocks are used for every degree. The base target is raised in
  no row on $12500\leq n\leq14999$, in $4$ rows on
  $15000\leq n\leq19999$, in $16$ rows on $20000\leq n\leq24999$, and
  in $27$ rows on $25000\leq n\leq34999$.
\item The core bound $b(n)=\floor{\log_2(n-2)}+2$ changes inside some
  intervals, for example $b(16385)=15$ and $b(16386)=16$, and
  $b(32769)=16$ and $b(32770)=17$. Taking $b_*=b(u)$ accounts for all
  values in the core supply and the residual side bounds. The label
  $32768$ enters the special even family at order $32768$; the removed
  label set is reconstructed at every order.
\item At order $20677=23\cdot29\cdot31$ the first support of size three
  above cutoff $19$ appears, and the support $\{23,29,37\}$ also occurs
  by order $24999$. These supports are included before maximality is
  tested, and no bound on the cardinality of a support is imposed.
\item On every interval from $12500$ on, the exponential global bound
  is the smaller one. All joint sums are less than $1/1000$ on
  $12500\leq n\leq19999$ and less than $1/100$ on
  $20000\leq n\leq34999$.
\end{enumerate}
On $5220\leq n\leq12499$, too, the exponential global bound is the
smaller one in every row. This proves the proposition.
\end{proof}

\subsubsection{The orders through 34999}

\begin{theorem}[The initial interval through 34999]
\label{thm:initial-prefix-34999}
Every tree of order at most $34999$ has a prime labeling.
\end{theorem}

\begin{proof}
Theorem~\ref{thm:initial-prefix-3499} covers the orders at most $3499$.
The seven covers of Table~\ref{tab:retained-certificates} partition the
orders from $3500=3499+1$ through $34999$ into adjacent intervals: the
consecutive ranges end at $5219$, $6771$, $12499$, $14999$, $19999$,
and $24999$, and each next range begins at the following integer. By
Proposition~\ref{prop:verified-retained-continuations}, every interval
satisfies the hypotheses of Theorem~\ref{thm:retained-interval-criterion}.
That theorem therefore gives a prime labeling for every tree of every
order from $3500$ through $34999$.
\end{proof}

The next finite range begins at $35000$. It requires the Hall
criterion with a label-side reserve of
Section~\ref{sec:asymmetric-reserve-hall} and is treated in
Section~\ref{sec:integer-bridge}. First,
Section~\ref{sec:part-analytic} develops the assignment needed for the
separate unbounded-range argument.

\section{Uniform estimates and the unbounded range}
\label{sec:part-analytic}

This part proves Theorem~\ref{thm:analytic-tail}: every tree of order
$n$ with $\log n\geq256$ has a prime labeling.
Sections~\ref{sec:rare-type-assignment}--\ref{sec:large-prime-cost}
give an assignment with rare arithmetic types
(Theorem~\ref{thm:rare-type-assignment}), a Hall estimate with a
reserve of universal positions (Theorem~\ref{thm:direct-reserve-hall}),
and a bound on the conflicts caused by prime divisors above the cutoff
(Theorem~\ref{thm:large-prime-sorting}); several logarithmic bridges
of Section~\ref{sec:part-log-bridges} also use them. Sections~\ref{sec:prime-estimates}--\ref{sec:analytic-tail}
choose the parameters as functions of $\log n$ and verify all
conditions analytically.

\subsection{Assignments with rare arithmetic types}
\label{sec:rare-type-assignment}

For a growing set of primes, keeping a separate rounding coordinate for
every type is expensive. We shall keep the types with small support and
combine the others during rounding. The exact types of the even labels
are still recorded. On the opposite side, a position whose tag becomes
too large is marked $\dagger$ and later receives a universal label.
The distinction matters: combining label classes is a temporary step in
meeting quotas, whereas marking a position creates an obligation to
assign it a compatible label.

Throughout this section, $G=(A,B)$ is a bipartite forest,
$N=\abs{A}+\abs{B}\geq1$, and $D\geq\max\{1,\Delta_G(A)\}$.
For two measures of the same total mass on a finite set, we use
$\|\xi-\zeta\|_{\mathrm{TV}}=\frac12\sum_x\abs{\xi(x)-\zeta(x)}$.
Thus errors in counting measures are measured in numbers of positions,
without normalization by $\abs{B}$.

\subsubsection{Component rounding with prescribed quotas}

We first cut the forest into small components. Unlike the one-sided
cuts used earlier, this cut may meet either part of the forest.

\begin{lemma}\label{lem:small-component-separator}
For every integer $t\geq1$, a forest on $N$ vertices has a set $C$ with
\[
 \abs{C}\leq\floor*{\frac{N}{t+1}}
\]
such that every component after deleting $C$ has at most $t$ vertices.
\end{lemma}

\begin{proof}
Root each component and process its vertices from the leaves upwards.
At a vertex $v$, join $v$ to the residual sets returned by its children.
If their union has more than $t$ vertices, put $v$ in $C$ and return
the empty set; otherwise return that union. Every returned nonempty set
is connected and has at most $t$ vertices.

When a vertex is put in $C$, charge it together with the residual sets
just received. These charged sets are disjoint and each contains at
least $t+1$ vertices. This proves the bound on $\abs{C}$. Each component
of the remaining forest is a residual set detached from a deleted
parent, or a residual set at an original root. Its size is at most $t$.
\end{proof}

Let $A_1,\ldots,A_q$ and $B_1,\ldots,B_b$ partition $A$ and $B$
into nonempty classes. An empty part has no classes. Put $d=q+b$,
so $1\leq d\leq N$. Fix integers $0\leq k_i\leq\abs{A_i}$ and
a number $0\leq\rho\leq1$, and write
\begin{equation}\label{eq:quota-positive-defect}
 \delta=\sum_{i=1}^q(k_i-\rho\abs{A_i})_+.
\end{equation}

\begin{lemma}[Quotas and additional tags]\label{lem:quota-tag-rounding}
There are sets $W\subseteq A$ and $T\subseteq B$ such that
\[
 \abs{W\cap A_i}=k_i,\qquad N_G(W)\subseteq T,
\]
and
\begin{equation}\label{eq:quota-tag-error}
 \sum_{j=1}^b\abs{\abs{T\cap B_j}-\rho\abs{B_j}}
 \leq2(D+1)\sqrt{Nd}+D\delta.
\end{equation}
More precisely, for any integer $t\geq1$ and the separator $C$ of
Lemma~\ref{lem:small-component-separator}, the construction gives the
upper bound $(D+1)(\abs{C}+dt)+D\delta$.
\end{lemma}

\begin{proof}
Represent each component of $G-C$ by the $d$ counts of its vertices
in the specified classes. Give each component selection weight $\rho$.
Apply Lemma~\ref{lem:fractional-pruning} with two colors, selected and
unselected. All class sums are preserved, and at most $d$ components
have fractional weights. Round these weights down to zero, and let
$V_*$ be the union of the components with weight one. For every class
$X$ in either part,
\[
 \abs{V_*\cap X}\leq\rho\abs{X\setminus C}.
\]
The sum of the deficits in this display, over all classes, is at most
$dt$: only the at most $d$ fractional components were discarded.

Start with $W_0=V_*\cap A$. In each $A_i$, delete excess vertices or
add missing vertices to obtain exactly $k_i$ selected vertices. Call
the resulting set $W$. The number of additions is at most
\begin{equation}\label{eq:quota-additions}
 L\leq\rho\abs{C\cap A}+dt+\delta.
\end{equation}
Indeed, if $a_i=\abs{W_0\cap A_i}$, then
\[
 (k_i-a_i)_+
 \leq(k_i-\rho\abs{A_i})_+
       +\rho\abs{A_i\cap C}
       +\rho\abs{A_i\setminus C}-a_i.
\]
Sum this inequality over $i$.

Keep the selected positions even when quota corrections remove their
neighbors, by defining
\begin{equation}\label{eq:quota-added-tag-set}
 T=(V_*\cap B)\cup N_G(W).
\end{equation}
Write $v_j=\abs{V_*\cap B_j}$ and $u_j=\abs{T\cap B_j}$.
We have $v_j\leq\rho\abs{B_j}$, and
\[
 \sum_j(\rho\abs{B_j}-v_j)\leq\rho\abs{C\cap B}+dt.
\]
The neighbors of $W_0$ outside $V_*$ lie in $C$. Each added vertex
has at most $D$ neighbors. Consequently
$\sum_j(u_j-v_j)\leq\abs{C\cap B}+DL$. It follows that
\begin{align*}
 \sum_j\abs{u_j-\rho\abs{B_j}}
 &\leq(1+\rho)\abs{C\cap B}+dt+DL\\
 &\leq(D+1)(\abs{C}+dt)+D\delta.
\end{align*}
Here $1+\rho\leq D+1$ and $D\rho\leq D+1$.
Finally choose $t=\floor{\sqrt{N/d}}\geq1$. Then
$\abs{C}\leq N/(t+1)<\sqrt{Nd}$ and $dt\leq\sqrt{Nd}$,
which gives \eqref{eq:quota-tag-error}.
\end{proof}

The extra positions in \eqref{eq:quota-added-tag-set} may create
unnecessary prohibitions. They preserve the lower count needed for
two-sided control of the tag distribution. Using only $N_G(W)$ would
not give that control.

\subsubsection{One refinement with rare classes combined}

Fix distinct odd primes $p_1,\ldots,p_k$ and an integer $r\geq0$.
The vertices of $A$ carry a bijection onto a fixed set of even labels.
After the first $i-1$ primes have been treated, let $A_I$ be the set
of vertices whose labels have exact type $I\subseteq\{p_1,\ldots,p_{i-1}\}$.
Write $k_I$ for the number of labels in this class divisible by $p_i$.
These numbers depend on the label set, not on its placement within $A_I$.

A position in $B$ has either a tag $J$ of size at most $r$ or the
mark $\dagger$. An ordinary tag contains every treated prime dividing
a neighboring label. For a marked position no compatibility is asserted.
Put
\[
 \mathcal R_i=\sum_{\abs{I}>r}\abs{A_I},\qquad
 \delta_i=\sum_I(k_I-\rho_i\abs{A_I})_+,
 \qquad 0\leq\rho_i\leq1.
\]

\begin{proposition}\label{prop:rare-refinement-step}
The labels can be permuted within their old exact classes so that the
new divisibility sets have all the prescribed sizes $k_I$. Tags can
be updated by adding $p_i$, leaving them unchanged, or replacing them
by $\dagger$ on overflow. An existing $\dagger$ remains $\dagger$.
Ordinary tags still contain all neighboring prime divisors.

Let $d_i$ be the number of nonempty classes used for rounding: one
for each $A_I$ with $\abs{I}\leq r$, one for their remaining union,
and one for each current tag class, including $\dagger$.
The total-variation error of the updated counting measure, relative
to independently adding $p_i$ with probability $\rho_i$ and then
marking overflow, is at most
\begin{equation}\label{eq:rare-step-error}
 E_i=2(D+1)\sqrt{Nd_i}+D(\delta_i+\mathcal R_i),\qquad
 d_i\leq2\sum_{j=0}^{\min(r,i-1)}\binom{i-1}{j}+2.
\end{equation}
\end{proposition}

\begin{proof}
Combine the classes $A_I$ with $\abs{I}>r$ into one class, prescribing
its total quota $\sum_{\abs{I}>r}k_I$. The positive quota defect of
this class is at most the sum of the positive defects of its parts.
Apply the construction of Lemma~\ref{lem:quota-tag-rounding} to this
partition and the current tag classes of $B$.

Within each rare exact class, now correct the selected set to size
$k_I$. These corrections require at most $\mathcal R_i$ additional vertices
in total. Let $W$ be the final selected set, and retain the original
union $V_*$ of selected components. Use
$T=(V_*\cap B)\cup N_G(W)$ as before. In the proof of
\eqref{eq:quota-additions}, the total number of additions is now at
most $\rho_i\abs{C\cap A}+dt+\delta_i+\mathcal R_i$.
The same two-sided tag estimate therefore gains at most $D\mathcal R_i$.

In each old exact class, place the $k_I$ labels divisible by $p_i$
on $W\cap A_I$ and place its other labels on the remaining vertices.
This is a bijection and preserves the type at every vertex with
respect to the first $i-1$ primes. Add $p_i$ to every ordinary tag
whose position lies in $T$, marking it $\dagger$ if its size exceeds
$r$. This gives the required containment of actual neighboring divisors.

For an ordinary old class of size $s_J$, let $u_J$ of its positions
receive the new prime. The two daughter counts differ from
$(1-\rho_i)s_J,\rho_i s_J$ by opposite errors. Their total-variation
contribution is exactly $\abs{u_J-\rho_i s_J}$. Combining the
overflow classes cannot increase this error. The old $\dagger$
class has no error because it is left fixed by both transitions.
The preceding rounding bound proves \eqref{eq:rare-step-error}.
There are at most the displayed number of small classes on each side,
and at most one additional class on each side. This proves the bound
on $d_i$.
\end{proof}

\subsubsection{Accumulating the errors}

Write $P=\{p_1,\ldots,p_k\}$. Let $c_r$ send a subset of $P$ to itself if its
size is at most $r$, and to $\dagger$ otherwise. Write
\[
 \nu=\bigotimes_{i=1}^k\operatorname{Bernoulli}(\rho_i),\qquad
 \beta=\nu\{J:\abs{J}>r\},\qquad
 \mathcal E=\sum_{i=1}^k E_i.
\]
Here the product law is identified with a law on subsets of $P$.

\begin{theorem}[Assignment with rare tags]\label{thm:rare-type-assignment}
Starting with empty tags, the successive refinements produce a
placement of the given even labels and a tag counting measure $\Xi$
on $B$, of total mass $s=\abs{B}$, such that
\begin{equation}\label{eq:rare-accumulated-error}
 \|\Xi-s(c_r)_*\nu\|_{\mathrm{TV}}\leq\mathcal E,
 \qquad
 C:=\Xi(\dagger)\leq s\beta+\mathcal E.
\end{equation}
Every ordinary tag contains all primes in $P$ dividing neighboring
labels. The construction preserves the exact arithmetic quotas at
each refinement, including those inside the combined rare class.
\end{theorem}

\begin{figure}[tbp]
\centering
\begin{tikzpicture}[
  >=Stealth,
  state/.style={draw,rounded corners=2pt,minimum height=8mm,
    minimum width=27mm,align=center,font=\normalsize},
  edge label/.style={font=\normalsize,fill=white,inner sep=2pt}]
\node[state] (old) at (0,1.2) {$J$, $\abs{J}\leq r$};
\node[state] (dag) at (0,-1.2) {$\dagger$};
\node[state] (same) at (6.3,2.1) {$J$};
\node[state] (grown) at (6.3,.6) {$J\cup\{p_i\}$};
\node[state] (overflow) at (6.3,-1.2) {$\dagger$};
\draw[->] (old) -- node[edge label,above,sloped] {$1-\rho_i$} (same);
\draw[->] (old) -- node[edge label,above,sloped] {$\rho_i$, $\abs{J}<r$} (grown);
\draw[->] (old.south east) -- node[edge label,above,sloped]
  {$\rho_i$, $\abs{J}=r$} (overflow.north west);
\draw[->] (dag) -- node[edge label,above] {$1$} (overflow);
\end{tikzpicture}
\caption{The model transition at one prime. An ordinary tag either
stays fixed or gains the new prime. Overflow and previously marked
positions enter the same absorbing class. Even-label classes remain
exact throughout this operation.}
\label{fig:rare-tag-transition}
\end{figure}
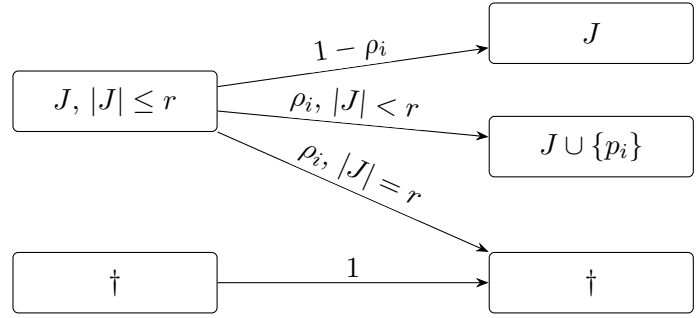

\begin{proof}
The independent transition at step $i$ leaves $\dagger$ fixed, and
either adds $p_i$ to an ordinary tag or leaves it unchanged. It sends
the compressed product law on the first $i-1$ primes to that on the
first $i$ primes: a set already larger than $r$ cannot become small
when a new coordinate is added. Figure~\ref{fig:rare-tag-transition}
shows this transition.

Every Markov transition contracts total variation. Indeed, if its
matrix is $K$, then
\[
 \frac12\sum_y\left|\sum_x(\xi(x)-\zeta(x))K(x,y)\right|
 \leq\frac12\sum_x\abs{\xi(x)-\zeta(x)}\sum_yK(x,y)
 =\|\xi-\zeta\|_{\mathrm{TV}}.
\]
Proposition~\ref{prop:rare-refinement-step} therefore adds at most
$E_i$ to the previous error. The initial error is zero, so induction
gives the first inequality in \eqref{eq:rare-accumulated-error}.
Apply it to the event $\{\dagger\}$ to obtain the second.
Tag containment and exact quotas follow at every step from the same
proposition. A marked position is counted once in $C$, even if it
would overflow at several later steps.
\end{proof}

\subsubsection{Arithmetic bounds for the errors}

For the remainder of the section, the even labels are obtained from
$\{2,4,\ldots,2m\}$ by deleting $w$ labels of empty $P$-type,
where $m\geq1$. Let $\mu$ be the product law in
\eqref{eq:prime-product-law}, and suppose
\[
 \|\mu_m-\mu\|_{\mathrm{TV}}\leq\varepsilon.
\]
Multiplication by two preserves these types. Suppose also that the
old empty class is nonempty at each step, and choose
\begin{equation}\label{eq:rare-empty-adjustment}
 \rho_i=\max\left\{\frac1{p_i},
                 \frac{k_\varnothing}{\abs{A_\varnothing}}\right\}.
\end{equation}
In particular $1/p_i\leq\rho_i\leq1$.

\begin{lemma}[The quota defect]\label{lem:rare-quota-defect}
For every step $i$,
\begin{equation}\label{eq:rare-arithmetic-defect}
 \delta_i\leq m(1+1/p_i)\varepsilon\leq2m\varepsilon.
\end{equation}
If $a_0=\abs{A_\varnothing}>0$ at that step, then also
\begin{equation}\label{eq:rare-inflated-coordinate}
 p_i\rho_i-1\leq\frac{w+(p_i+1)m\varepsilon}{a_0}.
\end{equation}
\end{lemma}

\begin{proof}
The empty class has no positive defect by
\eqref{eq:rare-empty-adjustment}. Let $\mathcal J$ be the nonempty
old types with positive defect, and put $q=1/p_i$. No deleted label
has any of these types. Replacing $\rho_i$ by $q$ can only increase
their sum of defects. Let $S$ be the full random type and $S_0$ its
restriction to the old primes. Dividing the sum of defects by $m$
therefore gives at most
\[
 \Pr_{\mu_m}(S_0\in\mathcal J,\ p_i\in S)
       -q\Pr_{\mu_m}(S_0\in\mathcal J).
\]
For the product law these two terms cancel, by independence.
Each probability changes by at most $\varepsilon$, because it is
an event in the full type space. This proves
\eqref{eq:rare-arithmetic-defect}.

Let $C_0$ be the size of the old empty class in the full interval.
Then $a_0=C_0-w$, and deleting empty full types removes no label
divisible by $p_i$. The same comparison of events gives
$k_\varnothing-qC_0\leq m(1+q)\varepsilon$. Hence
\[
 k_\varnothing-qa_0\leq qw+m(1+q)\varepsilon.
\]
Take the positive part, divide by $a_0$, and multiply by $p_i$.
This is \eqref{eq:rare-inflated-coordinate}.
\end{proof}

Define
\begin{equation}\label{eq:rare-tail-parameters}
 t=r+1,\qquad \sigma=\sum_{p\in P}\frac1p,\qquad
 Q_r=\frac{\sigma^t}{t!},\qquad
 C_r(k)=\sum_{j=0}^{\min(r,k)}\binom{k}{j}.
\end{equation}

\begin{corollary}\label{cor:rare-arithmetic-budget}
With these hypotheses, we have
\begin{align}
 \mathcal R_i&\leq mQ_r,\qquad
 \beta\leq\frac{(\sum_i\rho_i)^t}{t!},
                         \label{eq:rare-factorial-tails}\\
 \mathcal E&\leq
 2k(D+1)\sqrt{N(2C_r(k)+2)}
       +Dkm(2\varepsilon+Q_r).
                         \label{eq:rare-total-budget}
\end{align}
One may take $\varepsilon=B_r/m+2R_r$ from
Proposition~\ref{prop:moment-tv}.
If $\rho_i\leq(1+\eta)/p_i$ for all $i$, then
$\beta\leq((1+\eta)\sigma)^t/t!$.
\end{corollary}

\begin{proof}
For every integer $a\geq0$,
$\boldsymbol{1}_{a\geq t}\leq\binom at$. If $S$ is the type of a
uniform integer in $[m]$, then
\[
 \E\binom{\abs{S}}t
   =\sum_{\substack{I\subseteq P\\\abs{I}=t}}
       \frac{\floor{m/d_I}}m
   \leq\sum_{\substack{I\subseteq P\\\abs{I}=t}}\frac1{d_I}
   \leq\frac{\sigma^t}{t!}.
\]
The last inequality follows by expanding $\sigma^t$: each product
of $t$ distinct reciprocals occurs $t!$ times. An old rare type has
at least $t$ primes in its full type. Deleting labels cannot increase
their number, so $\mathcal R_i\leq mQ_r$.

Apply the same argument to independent Bernoulli coordinates of
probabilities $\rho_i$. Their joint moments are the corresponding
products, proving the bound for $\beta$. These arguments remain
valid when $t>k$, in which case the two exact tails vanish.
Sum \eqref{eq:rare-step-error}, use
$d_i\leq2C_r(k)+2$, and apply
Lemma~\ref{lem:rare-quota-defect}. This gives
\eqref{eq:rare-total-budget}. The final assertions follow from
Proposition~\ref{prop:moment-tv} and the bound on each $\rho_i$.
\end{proof}

\subsubsection{Assigning labels to marked positions}

The bound on $C$ must be paid with actual labels before the compressed
law can enter a matching argument. We record both distribution errors
after this payment.

\begin{lemma}[Removing the marked pairs]\label{lem:rare-universal-payment}
Let $\ell=s+h>0$, where $s=\abs{B}$ and $h\geq0$ is an integer.
Adjoin $h$ universal positions of empty tag. Suppose that the
$\ell$ available odd labels have type law $\widehat\mu$ with
$\|\widehat\mu-\mu\|_{\mathrm{TV}}\leq a$, and include $H$
distinct universal labels of empty $P$-type. Universality here means
compatibility with every opposite vertex under all arithmetic
constraints, including primes outside $P$.

Suppose \eqref{eq:rare-accumulated-error} holds, and there is an
integer $0\leq c<\ell$ such that
\begin{equation}\label{eq:rare-reserve-payment}
 s\beta+\mathcal E\leq c,\qquad H\geq h+c.
\end{equation}
Then the $C$ marked positions can receive distinct universal labels.
After deleting these matched pairs, there are $\ell-C$ vertices on
each side, at least $h$ universal vertices in each part, and only
ordinary tags. Their two type laws $\widehat\mu',\widehat\nu'$ satisfy
\begin{align}
 \|\widehat\mu'-\mu\|_{\mathrm{TV}}
       &\leq\frac{\ell a+c}{\ell-c},
                       \label{eq:rare-paid-label-error}\\
 \|\widehat\nu'-\nu\|_{\mathrm{TV}}
       &\leq\frac{\mathcal E+h+\ell\beta+c}{\ell-c}.
                       \label{eq:rare-paid-tag-error}
\end{align}
Any matching of the remaining vertices that satisfies all arithmetic
constraints extends by the $C$ fixed pairs.
\end{lemma}

\begin{proof}
By \eqref{eq:rare-accumulated-error}, $C\leq c$. Choose $C$ distinct
labels from the specified universal family and assign them bijectively
to the marked positions. Every chosen pair is admissible. At least
$H-C\geq h$ universal labels remain, and none of the $h$ adjoined
positions was used.

For the distribution estimate only, replace each $\dagger$ by the
empty type. The resulting image of $(c_r)_*\nu$ differs from $\nu$
by total variation $\beta$: all mass on types larger than $r$ has
been moved to the empty type. Contraction and the triangle inequality
give error at most $\mathcal E+s\beta$ for the $s$ positions.
Adding $h$ empty types changes the comparison with $\ell\nu$ by
at most $h$. Deleting the $C$ marked positions then changes it by
at most $C$, since
\[
 \|C\delta_\varnothing-C\nu\|_{\mathrm{TV}}\leq C.
\]
After this deletion the temporary empty representations are gone;
all remaining tags are the ordinary tags of the construction.
Divide by $\ell-C$ to obtain
\[
 \|\widehat\nu'-\nu\|_{\mathrm{TV}}
 \leq\frac{\mathcal E+s\beta+h+C}{\ell-C}
 \leq\frac{\mathcal E+\ell\beta+h+c}{\ell-c}.
\]
Likewise the original label counting measure differs from $\ell\mu$
by at most $\ell a$. Deleting $C$ empty-type labels adds at most
$C$ to this error relative to $(\ell-C)\mu$. This proves
\eqref{eq:rare-paid-label-error}. Both replacements of $C$ by $c$
are valid because $(x+C)/(\ell-C)$ increases with $C$ for $x\geq0$.
Finally, the fixed pairs use distinct labels and positions and are
disjoint from the remaining matching. Their union is the required
extension.
\end{proof}

In the residual problem of Proposition~\ref{prop:residual-problem},
the adjoined positions are the artificial universal vertices used in
Corollary~\ref{cor:residual-hall}. Every unused label in $O_*$ is
coprime to every label in $E_0$: its only possible exceptional even
partner has already been removed. If all primes in $P$ are at most
$z<n/3$, these universal labels also have empty $P$-type. A completed
matching therefore gives the residual injection; the restoration and
repair remain those of Proposition~\ref{prop:residual-problem}.

The counting and assignment arguments in this section hold for every
finite prime set. To use them for a growing cutoff, it remains to bound
the two errors in \eqref{eq:rare-paid-label-error} and
\eqref{eq:rare-paid-tag-error}, and to control deletions caused by
larger primes. The direct Hall estimate with universal reserves will
provide the corresponding matching criterion in
Section~\ref{sec:direct-reserve-hall}.

\subsection{A direct Hall estimate with universal reserves}
\label{sec:direct-reserve-hall}

The preceding construction leaves an ordinary matching problem after its
marked positions have been assigned universal labels. We now give a
criterion for that matching. It uses the two full type laws, rather than
conditioning either law on a rare class. This distinction permits a total
variation error proportional to the mass of the empty position type.

Let $P=\{p_1,\ldots,p_k\}$ be a finite set of primes. A label type and a
position tag are subsets of $[k]$; they are \emph{compatible} when they are
disjoint. Fix numbers
\begin{equation}\label{eq:direct-hall-model}
 \begin{gathered}
 0\leq a_i\leq\frac13,\qquad
 a_i\leq\rho_i\leq\frac25,\\
 \mu=\bigotimes_{i=1}^k\operatorname{Bernoulli}(a_i),
 \qquad \nu=\bigotimes_{i=1}^k\operatorname{Bernoulli}(\rho_i),
 \qquad b=\prod_{i=1}^k(1-\rho_i).
 \end{gathered}
\end{equation}
The empty product is one. Thus $b>0$. All total variation distances in
this section are between probability measures. For nonnegative measures
$f$ and $g$ on the two type spaces, put
\[
 \mathcal C(f,g)
   =\sum_{I\cap J=\varnothing}f(I)g(J).
\]
When $f$ and $g$ count selected vertices and have been divided by the
side length, $\mathcal C(f,g)$ is the density of compatible pairs in the
corresponding rectangle.

\begin{lemma}[A robust compatibility rectangle]
\label{lem:direct-hall-rectangle}
Let $\widehat\mu,\widehat\nu$ be probability measures with
$\|\widehat\mu-\mu\|_{\mathrm{TV}},
  \|\widehat\nu-\nu\|_{\mathrm{TV}}\leq\epsilon$.
If $f'\leq\widehat\mu$ and $g'\leq\widehat\nu$ have masses $u,v$
with $u+v>1$, then
\begin{equation}\label{eq:direct-hall-rectangle}
 \mathcal C(f',g')
 \geq b\left(\frac25uv-\frac75\epsilon\right).
\end{equation}
Moreover, under either actual law every fixed type is compatible with
at least $b-\epsilon$ of the opposite part.
\end{lemma}

\begin{proof}
We couple the two model laws coordinate by coordinate. Give
$(X_i,Y_i)$ the values $(1,0),(0,1),(0,0)$ with respective probabilities
$a_i,\rho_i,1-a_i-\rho_i$. These are nonnegative by
\eqref{eq:direct-hall-model}. The product coupling $\mathcal J$ has
marginals $\mu,\nu$ and is supported on compatible pairs.

For a compatible pair $(I,J)$, comparison in each coordinate gives
\begin{equation}\label{eq:direct-hall-density}
 \mathcal J(I,J)\leq \frac{\mu(I)\nu(J)}{b}.
\end{equation}
Indeed, for permitted coordinate outcomes of positive product mass, the
joint-to-product ratios are respectively
\[
 \frac1{1-\rho_i},\qquad \frac1{1-a_i},\qquad
 \frac{1-a_i-\rho_i}{(1-a_i)(1-\rho_i)};
\]
each is at most $(1-\rho_i)^{-1}$ because $a_i\leq\rho_i$.

We also need a lower bound for the coupling probability of two events.
For the centered, variance-one indicators of the two $i$th coordinates,
their correlation is $-\theta_i$, where
\[
 \theta_i^2=\frac{a_i\rho_i}{(1-a_i)(1-\rho_i)}
 \leq\frac13<\frac9{25}.
\]
If $a_i=0$, omit its constant coordinate. Expand functions
$F,G\colon2^{[k]}\to[0,1]$ in the orthonormal product bases of their
respective marginal laws. Independence across coordinates makes mixed
terms vanish unless they use the same nonempty coordinate set. Each
remaining mixed term has multiplier of absolute value at most $3/5$.
Cauchy--Schwarz and $\operatorname{Var}(F)\leq x(1-x)$, where
$x=\mathbb E_\mu F$, therefore give, with $y=\mathbb E_\nu G$,
\begin{equation}\label{eq:direct-hall-correlation}
 \mathbb E_{\mathcal J}[FG]
 \geq xy-\frac35\sqrt{x(1-x)y(1-y)}.
\end{equation}

Truncate the selected submeasures at the model laws:
\[
 f(I)=\min\{f'(I),\mu(I)\},\qquad
 g(J)=\min\{g'(J),\nu(J)\}.
\]
Since $f'\leq\widehat\mu$,
the mass discarded from $f'$ is at most
$\sum_I(\widehat\mu(I)-\mu(I))_+
=\|\widehat\mu-\mu\|_{\mathrm{TV}}\leq\epsilon$;
the same holds for $g'$. Write $x=\sum f$ and $y=\sum g$.
Then $x\geq u-\epsilon$, $y\geq v-\epsilon$, and
$x+y>1-2\epsilon$. Apply \eqref{eq:direct-hall-correlation} to
$F=f/\mu$ and $G=g/\nu$, taking zero on any zero-mass type. As
$(1-x)(1-y)=xy+1-x-y<xy+2\epsilon$,
\[
 \sqrt{x(1-x)y(1-y)}
 \leq\sqrt{xy(xy+2\epsilon)}\leq xy+\epsilon.
\]
Furthermore $xy\geq uv-\epsilon(u+v)\geq uv-2\epsilon$.
Equation~\eqref{eq:direct-hall-density} and nonnegativity now give
\[
 \mathcal C(f',g')\geq\mathcal C(f,g)
 \geq b\mathbb E_{\mathcal J}[FG]
 \geq b\left(\frac25uv-\frac75\epsilon\right).
\]

For the last assertion, a fixed label type $I$ is compatible with a
$\nu$-random tag with probability
$\prod_{i\in I}(1-\rho_i)\geq b$. A fixed tag $J$ is compatible with a
$\mu$-random label type with probability
$\prod_{i\in J}(1-a_i)\geq b$. Total variation changes the probability
of either event by at most $\epsilon$.
\end{proof}

The next result keeps track of the edges needed to violate Hall's
condition. The reserve is counted on both sides, since after the marked
pairs have been removed there are still unused universal labels as well
as the artificial universal positions.

\begin{samepage}
\begin{theorem}[Direct Hall criterion with a reserve]
\label{thm:direct-reserve-hall}
Let a bipartite graph have $\ell$ vertices on each side, with empirical
type laws $\widehat\mu,\widehat\nu$ satisfying
\begin{equation}\label{eq:direct-hall-errors}
 \|\widehat\mu-\mu\|_{\mathrm{TV}},
 \|\widehat\nu-\nu\|_{\mathrm{TV}}\leq\epsilon\leq\frac b{25}.
\end{equation}
Start with every compatible pair as an edge. Suppose there are at least
$h\geq1$ vertices in each part adjacent to every vertex of the other
part, and no edge at those vertices is subsequently deleted. If at most
$M$ other edges are deleted and
\begin{equation}\label{eq:direct-hall-budgets}
 h\leq\frac{b\ell}{100},\qquad
 M\leq\frac{19}{20}b\ell h,
\end{equation}
then the remaining graph has a perfect matching.
\end{theorem}
\end{samepage}

\begin{proof}
Suppose Hall's condition fails. There is a nonempty set $Q$ on one
side and a nonempty set $R$ on the other such that
$\abs Q+\abs R>\ell$ and all compatible edges between $Q$ and $R$
were deleted. One may take $R$ to be the complement of the remaining
neighborhood of $Q$. Neither set contains a universal vertex, for such
a vertex would retain all of its edges. In particular
\[
 h<\abs Q,\abs R\leq\ell-h.
\]
Put $u=\abs Q/\ell$, $v=\abs R/\ell$,
$x=\min\{u,v\}$, and $r=h/\ell$. Thus
$u+v>1$, $r<x\leq1-r$, and $r\leq b/100$.
We count the compatible edges in $Q\times R$ before deletion.

First suppose $x<b/4$. A vertex in the smaller set originally has at
least $(b-\epsilon)\ell\geq(24/25)b\ell$ compatible neighbors, by
Lemma~\ref{lem:direct-hall-rectangle}. Since the complement of the
larger set has fewer than $x\ell$ vertices, the rectangle contains
strictly more than
\begin{equation}\label{eq:direct-hall-narrow}
 \ell^2x\left(\frac{24}{25}b-x\right)
\end{equation}
compatible edges. The quadratic in \eqref{eq:direct-hall-narrow} is
concave on $[r,b/4]$, so its minimum there occurs at an endpoint. At
the first endpoint it is at least
\[
 r\left(\frac{24}{25}b-r\right)
 \geq\frac{19}{20}br,
\]
where $r\leq b/100$ is used. At the second endpoint it equals
$71b^2/400>19br/20$. Thus more than $19b\ell h/20$ compatible edges
were deleted, contrary to \eqref{eq:direct-hall-budgets}.

Now suppose $x\geq b/4$. Since $u+v>1$ and both $u,v\geq b/4$,
\[
 uv\geq\frac b4\left(1-\frac b4\right)
 \geq\frac{3b}{16}.
\]
Apply Lemma~\ref{lem:direct-hall-rectangle} to the type submeasures
of $Q$ and $R$. With $\epsilon\leq b/25$, the number of compatible
edges in the rectangle is at least
\begin{equation}\label{eq:direct-hall-wide}
 b\ell^2\left(\frac25uv-\frac7{125}b\right)
 >\frac{b\ell^2uv}{10}
 \geq\frac{3b^2\ell^2}{160}.
\end{equation}
For the strict inequality, $3uv/10\geq9b/160>7b/125$.
But \eqref{eq:direct-hall-budgets} gives
$M\leq19b^2\ell^2/2000<3b^2\ell^2/160$.
Again the necessary edges cannot all have been deleted. Both cases
contradict the choice of $Q,R$, proving Hall's condition.
\end{proof}

We state explicitly how this criterion enters the residual labeling
problem. The number $M$ below is a bound on \emph{pairs}, rather than
on bad labels or positions; it is estimated from larger prime
divisors in Section~\ref{sec:large-prime-cost}.

\begin{proposition}[Conditional residual completion]
\label{prop:direct-hall-completion}
Use the residual construction of
Proposition~\ref{prop:residual-problem} and the assignment of
Section~\ref{sec:rare-type-assignment}. Let $s=\abs{Z_0}$ and
$\ell=s+h=\abs{O_0}$. Suppose the hypotheses of
Lemma~\ref{lem:rare-universal-payment} hold with its parameters
$\alpha,\mathcal E,\beta,c,H$, where $\alpha$ denotes the label-law
error called $a$ in that lemma. Suppose the model laws satisfy
\eqref{eq:direct-hall-model}. Assume
\begin{equation}\label{eq:direct-hall-completion}
 \frac{\max\{\ell\alpha,\mathcal E+h+\ell\beta\}+c}{\ell-c}
       \leq\frac b{25},
 \qquad h\leq\frac{b(\ell-c)}{100},
 \qquad M\leq\frac{19}{20}b(\ell-c)h.
\end{equation}
Here $M$ bounds the number of pairs in the remaining balanced square
that are compatible for primes in $P$ but incompatible for a prime
outside $P$. Suppose also that the $h$ artificial positions and at
least $h$ of the remaining labels are genuinely universal for the
residual problem. Then the original tree has a prime labeling.
\end{proposition}

\begin{proof}
By Lemma~\ref{lem:rare-universal-payment}, the marked positions use
distinct universal labels and $C\leq c$. Set $\ell'=\ell-C$.
The two empirical laws of the remaining square are within the two
bounds \eqref{eq:rare-paid-label-error} and
\eqref{eq:rare-paid-tag-error} of their respective model laws.
Equation~\eqref{eq:direct-hall-completion} makes each bound at most
$b/25$. It also implies
$h\leq b\ell'/100$ and
$M\leq19b\ell'h/20$, since $\ell'\geq\ell-c$.

An ordinary tag contains every prime of $P$ that divides a neighboring
even label. Consequently disjointness of a label type and a tag
guarantees coprimality at those primes. Delete the at most $M$ pairs
that fail because of a prime outside $P$. The stated universal pairs
survive this deletion. Theorem~\ref{thm:direct-reserve-hall} supplies a
perfect matching of the remaining square. Add back the fixed marked
pairs and discard the $h$ pairs at artificial positions. What remains
is an injection from $Z_0$ to $O_0$ satisfying
\eqref{eq:residual-coprimality}. Proposition~\ref{prop:residual-problem}
restores the removed labels and repairs the allowed exceptional pairs.
\end{proof}

In applications the even labels may be permuted within each exact
$P$-type before $M$ is bounded. Such a permutation preserves every
exact label quota and every already constructed $P$-tag: each
neighbor's type remains the same. It changes only which pairs share a
prime outside $P$. The estimate for those pairs is supplied by
Theorem~\ref{thm:large-prime-sorting}.

\subsection{The cost of larger prime divisors}
\label{sec:large-prime-cost}

The Hall criterion of Section~\ref{sec:direct-reserve-hall} allows a
specified number of edges to be removed from the type-compatible graph.
We bound these removals by choosing the even labels within their exact
types. The choice leaves every tag from
Section~\ref{sec:rare-type-assignment} unchanged.

Let $n\geq1$, $m=\floor{n/2}$, and $3\leq z\leq m$. In this section
$P$ consists of all odd primes at most $z$. We use the exact counts
$C_I(M)$ and the product law $\mu$ of
Section~\ref{sec:arithmetic-types}. For $j\leq m$, put
\begin{equation}\label{eq:large-prime-weight}
 \tau_z(j)=\sum_{\substack{p>z\\p\mid j}}\frac1p,
 \qquad
 H_z=\sum_{z<p\leq m}\frac1p,
 \qquad
 Q_z=\sum_{z<p\leq m}\frac1{p^2},
\end{equation}
where sums indexed by $p$ run over primes. Empty sums are zero.
The factor $1/p$ in $\tau_z$ reflects the density of odd labels
divisible by $p$; the discrepancy between its total in one exact type
and the product prediction is controlled by the earlier moment bound.

\begin{lemma}[Tail weights in exact types]
\label{lem:large-prime-type-total}
Take the parameters $B_r,R_r$ of
Proposition~\ref{prop:moment-tv}, for any integer $r\geq0$. For each
$I\subseteq P$, write
\[
 T_I=\sum_{\substack{j\leq m\\J_P(j)=I}}\tau_z(j).
\]
Then
\begin{equation}\label{eq:large-prime-type-total}
 \sum_{I\subseteq P}\abs{T_I-C_I(m)Q_z}
 \leq (2B_r+1)H_z+(2B_r+8mR_r)Q_z.
\end{equation}
\end{lemma}

\begin{proof}
For $p>z$, multiplication by $p$ preserves the $P$-type. Hence
\begin{equation}\label{eq:large-prime-divisibility-identity}
 T_I=\sum_{z<p\leq m}\frac{C_I(\floor{m/p})}{p}.
\end{equation}
The ordinary estimate in Proposition~\ref{prop:moment-tv} gives,
for every integer $M\geq1$,
\begin{equation}\label{eq:large-prime-moment-l1}
 \sum_I\abs{C_I(M)-M\mu(I)}\leq2B_r+4MR_r.
\end{equation}
Set $M_p=\floor{m/p}$. Because $\sum_I\mu(I)=1$ and
$0\leq m/p-M_p<1$, \eqref{eq:large-prime-moment-l1} implies
\[
 \sum_I\abs{C_I(M_p)-(m/p)\mu(I)}
 \leq2B_r+1+4(m/p)R_r.
\]
At $M=m$, the same estimate is $2B_r+4mR_r$.
Subtract $C_I(m)Q_z$ from \eqref{eq:large-prime-divisibility-identity}
and insert $(m/p)\mu(I)$ in each summand. The triangle inequality,
summed over $I$, gives
\[
 \sum_I\abs{T_I-C_I(m)Q_z}
 \leq\sum_{z<p\leq m}\frac1p
  \left(2B_r+1+4\frac{m}{p}R_r
        +\frac{2B_r+4mR_r}{p}\right),
\]
which is \eqref{eq:large-prime-type-total}.
\end{proof}

The next theorem converts these class totals into a bound on forbidden
label--position pairs. We state the hypotheses for the residual graph;
the forest property is not needed for this counting step.

\begin{theorem}[Sorting within exact types]
\label{thm:large-prime-sorting}
Let $G=(A,B)$ be a bipartite graph. Remove a set
$\mathcal S\subseteq[m]$ of $w$ half-labels, all of empty $P$-type,
and put
\[
 \Lambda_I=\{j\in[m]\setminus\mathcal S:J_P(j)=I\},
 \qquad a_I=\abs{\Lambda_I}.
\]
Suppose $A=\bigsqcup_I A_I$, with $\abs{A_I}=a_I$, and
$a_\varnothing>0$. Let $d(a)=d_G(a)\leq D$ for all $a\in A$,
where $D\geq1$, and put $e=e(G)$. Assign to every $a\in A_I$
an even label $2j(a)$ with $j(a)\in\Lambda_I$, bijectively within
each class. Let $Y\subseteq[n]$ be any set of odd labels.

Each $b\in B$ has an ordinary tag containing the $P$-types of all
its neighbors in $G$. Let $M$ count the pairs $(y,b)\in Y\times B$
whose types are disjoint from the tag of $b$, but for which $y$
shares a prime divisor with some $2j(a)$, $a\in N_G(b)$.
There are choices of the class bijections, leaving all tags fixed,
for which
\begin{align}
 M\leq{}&me\left(1+\frac{w}{a_\varnothing}\right)Q_z
 \label{eq:large-prime-deletion-bound}\\
 &+mD\bigl((2B_r+1)H_z+(2B_r+8mR_r)Q_z\bigr)
       +e\frac{\log n}{\log z}.\notag
\end{align}
In particular, $Q_z$ may be replaced throughout the right side by
$S_2(z)=\sum_{p>z}p^{-2}$.
\end{theorem}

\begin{proof}
First count the possible forbidden odd labels for a fixed half-label
$j$. If $p>z$ divides $j$, then among the odd integers in $[n]$ at
most $\ceil{\floor{n/p}/2}$ are divisible by $p$. As $n\leq2m+1$,
\[
 \ceil{\floor{n/p}/2}\leq\frac{m}{p}+1.
\]
Every prime shared by $y$ and $2j$ but absent from $P$ is greater
than $z$: $y$ is odd, and $P$ contains all odd primes at most $z$.
The union bound therefore gives at most
\begin{equation}\label{eq:large-prime-one-label}
 m\tau_z(j)+\omega_{>z}(j)
 \leq m\tau_z(j)+\frac{\log n}{\log z}
\end{equation}
possible odd partners, where $\omega_{>z}(j)$ counts the distinct
prime divisors of $j$ greater than $z$. The last inequality follows
because their product is at least $z^{\omega_{>z}(j)}$ and is at most
$j\leq m\leq n$.

Charge a forbidden pair $(y,b)$ to a neighbor $a\in N_G(b)$ sharing
a prime greater than $z$ with $y$. Some pairs may be charged more
than once. By \eqref{eq:large-prime-one-label}, for every placement
of the even labels,
\begin{equation}\label{eq:large-prime-charge}
 M\leq m\sum_{a\in A}d(a)\tau_z(j(a))
                          +e\frac{\log n}{\log z}.
\end{equation}

For each nonempty class $A_I$, arrange its degrees in decreasing
order and the numbers $\tau_z(j)$, $j\in\Lambda_I$, in increasing
order. A swap of two labels that were oppositely ordered changes
the weighted sum by the negative of the product of the two relevant
differences. Repeating swaps produces the arrangement just described.
Its weighted sum is no greater than the average over all bijections,
namely
\begin{equation}\label{eq:large-prime-rearrangement}
 \sum_{a\in A_I}d(a)\tau_z(j(a))
 \leq\frac{e_I}{a_I}\sum_{j\in\Lambda_I}\tau_z(j),
 \qquad e_I=\sum_{a\in A_I}d(a).
\end{equation}
The choices in different classes are independent.

Deleting the special half-labels cannot increase a class total of
the nonnegative weights $\tau_z$. By
Lemma~\ref{lem:large-prime-type-total}, the sum in
\eqref{eq:large-prime-rearrangement} is at most
\[
 C_I(m)Q_z+\abs{T_I-C_I(m)Q_z}.
\]
Here $C_I(m)=a_I+w\boldsymbol1_{I=\varnothing}$.
Since $e_I/a_I\leq D$ for every nonempty class and
$e_\varnothing\leq e$, summing \eqref{eq:large-prime-rearrangement}
over the classes yields
\[
 \sum_{a\in A}d(a)\tau_z(j(a))
 \leq e\left(1+\frac{w}{a_\varnothing}\right)Q_z
       +D\bigl((2B_r+1)H_z+(2B_r+8mR_r)Q_z\bigr).
\]
Combining this with \eqref{eq:large-prime-charge} proves
\eqref{eq:large-prime-deletion-bound}. The $P$-type at each vertex
has not changed, so every neighboring type and every ordinary tag
from the earlier construction is the same. Finally,
$Q_z\leq S_2(z)$ and all coefficients are nonnegative.
\end{proof}

For the assignment in Section~\ref{sec:rare-type-assignment}, take
$G$ to be the residual forest. When $z<n/3$, the half-labels of the
removed special even family all have empty $P$-type; take them as
$\mathcal S$. Theorem~\ref{thm:large-prime-sorting}
bounds the deletions among all ordinary positions and available odd
labels; removing the fixed marked pairs can only reduce their number.
The genuine universal labels and artificial positions have no deleted
edges. Consequently the expression on the right of
\eqref{eq:large-prime-deletion-bound} is a valid choice for $M$ in
Proposition~\ref{prop:direct-hall-completion}. To finish the unbounded
range it remains to compare that expression, the two type errors, and
the reserve sizes with their required inequalities uniformly in $n$.

\subsection{Prime estimates for the unbounded range}
\label{sec:prime-estimates}

The construction above requires a growing set of type coordinates, a
stock of special prime labels, and a bound for primes omitted from the
type system. We record the estimates together, with their domains. Write
$\pi(x)$ for the number of primes at most $x$, and let $\gamma$ be
Euler's constant. We use the following explicit estimates of
Dusart~\cite[Theorems 6.9 and 6.12]{Dus10}:
\begin{align}
 \frac{x}{\log x}\left(1+\frac1{\log x}\right)
 &\leq\pi(x)
 &&(x>599),\label{eq:dusart-prime-lower}\\
 \pi(x)&\leq\frac{x}{\log x}
              \left(1+\frac{13}{10\log x}\right)
 &&(x>1),\label{eq:dusart-prime-upper}\\
 \frac{x}{\log x-1}&\leq\pi(x)
 &&(x>5393),\label{eq:dusart-denominator-lower}\\
 \pi(x)&\leq\frac{x}{\log x-11/10}
 &&(x>60184),\label{eq:dusart-denominator-upper}\\
 \prod_{p\leq x}\left(1-\frac1p\right)
 &>\frac{e^{-\gamma}}{\log x}
       \left(1-\frac1{5(\log x)^2}\right)
 &&(x>2973).\label{eq:dusart-product-lower}
\end{align}
The upper bound in \eqref{eq:dusart-prime-upper} is slightly weaker
than the cited bound; the round coefficient will simplify the
calculations below.

\begin{lemma}[A square tail over primes]
\label{lem:prime-square-tail}
For every real $x>599$,
\begin{equation}\label{eq:prime-square-tail}
 S_2(x):=\sum_{p>x}\frac1{p^2}
 \leq\frac1{x\log x}-\frac{2}{5x(\log x)^2}
 <\frac1{x\log x}.
\end{equation}
\end{lemma}

\begin{proof}
Put $v=\log x$ and
$I_j=\int_x^\infty dt/(t^2(\log t)^j)$. Partial summation gives
\[
 S_2(x)=-\frac{\pi(x)}{x^2}
                  +2\int_x^\infty\frac{\pi(t)}{t^3}\,dt.
\]
Use the lower estimate \eqref{eq:dusart-prime-lower} at $x$ and the
upper estimate \eqref{eq:dusart-prime-upper} in the integral. Since
$I_1+I_2=1/(xv)$ and $I_2\leq1/(xv^2)$, we obtain
\begin{align*}
 S_2(x)
 &\leq-\frac1{xv}-\frac1{xv^2}+2I_1+\frac{13}{5}I_2\\
 &=\frac1{xv}-\frac1{xv^2}+\frac35 I_2
 \leq\frac1{xv}-\frac{2}{5xv^2}.
\end{align*}
This proves both inequalities.
\end{proof}

We next specialize the bounds to the parameters used in the analytic
construction. The real number $z$ need not be integral: a prime belongs
to $P$ exactly when it does not exceed $z$.

\begin{lemma}[Supply and type coordinates]
\label{lem:analytic-prime-inputs}
Suppose $n$ is an integer, $L=\log n\geq256$, and put
\[
 z=12L,\qquad v=\log z,\qquad
 P=\{p:3\leq p\leq z,\ p\text{ prime}\},\qquad k=\abs P.
\]
Then $8<v\leq L/25$, and
\begin{align}
 \frac{3n}{5L}
 &\leq\pi(n)-\pi(n/3)
 \leq\frac{n}{L-2}-\frac{n}{3L},
 \label{eq:analytic-prime-window}\\
 k&\leq\frac{279}{160}L<\frac74L.
 \label{eq:analytic-prime-coordinate-count}
\end{align}
If
$\alpha_\varnothing=\prod_{p\in P}(1-1/p)$ and
$\sigma=\sum_{p\in P}1/p$, then
\begin{equation}\label{eq:analytic-prime-product}
 \alpha_\varnothing>\frac{11}{10v},
 \qquad \sigma<\log v\leq\frac L{120}.
\end{equation}
The product lower bound and the sum upper bound hold for every subset
of $P$ as well, with the same $v$.
Finally,
\begin{equation}\label{eq:analytic-prime-tails}
 S_2(z)<\frac1{zv},\qquad
 \sum_{z<p\leq\floor{n/2}}\frac1p\leq1+L.
\end{equation}
\end{lemma}

\begin{proof}
At $L=256$ we have $z=3072$ and
$8<\log3072<33/4$. The function
$L/25-\log(12L)$ is positive there and has derivative
$1/25-1/L>0$ for $L\geq256$. This gives $8<v\leq L/25$ and,
in particular, $z>2973$.

For the prime window, \eqref{eq:dusart-prime-lower} at $n$ and
\eqref{eq:dusart-prime-upper} at $n/3$ give
\[
 \pi(n)\geq\frac nL,
 \qquad
 \pi(n/3)\leq\frac{n}{3(L-3)}.
\]
To see the second inequality, put $t=L-\log3$. For $t>13/10$,
\[
 \frac1t\left(1+\frac{13}{10t}\right)
 \leq\frac1{t-13/10}<\frac1{L-3},
\]
because $\log3+13/10<3$. The difference of the two displayed
bounds is at least $3n/(5L)$ whenever $L\geq18$.
In the other direction,
\eqref{eq:dusart-prime-upper} at $n$ gives
$\pi(n)\leq n/(L-2)$, while
\eqref{eq:dusart-prime-lower} at $n/3$ gives
$\pi(n/3)\geq n/(3L)$. Both applications have arguments greater
than $599$. This proves \eqref{eq:analytic-prime-window}.

Since $v>8$, \eqref{eq:dusart-prime-upper} at $z$ yields
\[
 k\leq\pi(z)\leq\frac{12L}{v}
       \left(1+\frac{13}{10v}\right)
 \leq\frac{12L}{8}\left(1+\frac{13}{80}\right)
 =\frac{279}{160}L.
\]
For the product, exclude the factor at $p=2$ from
\eqref{eq:dusart-product-lower}. The elementary bounds
$\gamma<0.578$ and $e^{-0.578}>14/25$ give
\[
 \alpha_\varnothing
 >\frac{28}{25v}\left(1-\frac1{5v^2}\right)
 \geq\frac{28}{25v}\left(1-\frac1{320}\right)
 >\frac{11}{10v}.
\]
As $-\log(1-u)\geq u$ for $0<u<1$,
$\sigma\leq-\log\alpha_\varnothing<\log v$.
To verify the last bound in \eqref{eq:analytic-prime-product},
observe that $\log v<L/120$ at $L=256$ because
$v<33/4<e^{32/15}$. The derivative of $L/120-\log v$ is
$1/120-1/(Lv)>0$ for $L\geq256$.

Lemma~\ref{lem:prime-square-tail} gives the first estimate in
\eqref{eq:analytic-prime-tails}. For the other, discard the
restriction to primes and use the harmonic bound
$\sum_{j\leq\floor{n/2}}1/j\leq1+\log\floor{n/2}\leq1+L$.
\end{proof}

\subsection{The residual forest at large orders}
\label{sec:analytic-residual}

We now choose the sets in Proposition~\ref{prop:residual-problem} for
every order in the analytic range. The two sides of the residual
forest need nearly half the original vertices, while the special
labels must pay for a reserve and reduce the remaining degrees.
The stronger prime-counting estimates
\eqref{eq:dusart-denominator-lower} and
\eqref{eq:dusart-denominator-upper} supply the needed lower count of
special even labels.

Throughout this section let $n$ be an integer with $L=\log n\geq256$,
and put $m=\floor{n/2}$ and $M_o=\ceil{n/2}$. We shall use the
elementary size bound
\begin{equation}\label{eq:analytic-enormous-order}
 n>10^6L^2.
\end{equation}
Indeed, $e^L/L^2$ increases for $L\geq256$, and at $L=256$ we have
$e^{256}>2^{256}>10^6\cdot256^2$.

\begin{lemma}[Special-label supply]
\label{lem:analytic-special-supply}
The hypotheses of Corollary~\ref{cor:even-side} hold. For its
partition and the special sets of Lemma~\ref{lem:special-labels},
write $\kappa=\abs K$, $w=\abs{E_*}$ and $u=\abs{O_*}$. Then
\begin{equation}\label{eq:analytic-special-supply}
 \kappa<2L+2,
 \qquad \frac n{6L}<w\leq\frac{2m}{5L},
 \qquad \frac n{2L}<u<q:=\abs Z.
\end{equation}
\end{lemma}

\begin{proof}
Let $b(n)=\floor{\log_2(n-2)}+2$, as in
\eqref{eq:separator-size}. Since $\log2>1/2$,
$b(n)<2L+2$. The upper estimate
\eqref{eq:dusart-prime-upper}, applied at $n/2$, gives
\begin{equation}\label{eq:analytic-half-prime-upper}
 \pi(n/2)\leq\frac{n}{2(L-2)}.
\end{equation}
Indeed, with $t=L-\log2$, its right side before simplification is
$n(1+13/(10t))/(2t)\leq n/(2(t-13/10))$, and
$\log2+13/10<2$. Hence
\[
 \pi(n)-\pi(n/2)
 \geq\frac nL-\frac{n}{2(L-2)}
 =\frac{n(L-4)}{2L(L-2)}
 \geq\frac n{3L}>b(n).
\]
The penultimate inequality holds for $L\geq8$, and the last follows
from \eqref{eq:analytic-enormous-order}. Thus the required
partition exists and $1\leq\kappa\leq b(n)<2L+2$.

To bound $w$ from below, put
$A=1+\log2$ and $B=11/10+\log3$. The domains of
\eqref{eq:dusart-denominator-lower} at $n/2$ and
\eqref{eq:dusart-denominator-upper} at $n/3$ are satisfied. They give
\begin{equation}\label{eq:analytic-special-even-lower}
 \pi(n/2)-\pi(n/3)
 \geq\frac{n}{2(L-A)}-\frac{n}{3(L-B)}
 >\frac n{6L}.
\end{equation}
For the last inequality, multiplication by the positive denominators
reduces it to $L(3A-2B)>AB$. The elementary logarithmic bounds
$5/3<A<17/10$ and $B<11/5$ imply
$3A-2B>3/5$ and $AB<187/50$; these suffice at $L\geq256$.
The exact count \eqref{eq:even-supply} now gives $w>n/(6L)$.

For the upper bound on $w$, use
\eqref{eq:analytic-half-prime-upper} and the lower estimate
$\pi(n/3)\geq n/(3L)$ from the proof of
Lemma~\ref{lem:analytic-prime-inputs}. Since
$\floor{\log_2 n}<2L$,
\[
 w<2L+\frac{n(L+4)}{6L(L-2)}
 \leq2L+\frac{9n}{50L}
 \leq\frac{2m}{5L}.
\]
The middle inequality uses $L\geq77$. For the last,
$2m\geq n-1$ and it suffices that $n\geq100L^2+10$,
which follows from \eqref{eq:analytic-enormous-order}.

Finally, \eqref{eq:odd-supply},
\eqref{eq:analytic-prime-window}, and $\kappa<2L+2$ give
\[
 u=1+\pi(n)-\pi(n/3)-\kappa
 >\frac{3n}{5L}-2L-1>\frac n{2L}.
\]
On the other hand,
$u\leq1+\pi(n)-\kappa< M_o-\kappa=q$:
the last inequality follows from
$\pi(n)\leq n/(L-2)$ and \eqref{eq:analytic-enormous-order}.
This completes the proof.
\end{proof}

We spend one eighth of the special even labels on vertices of highest
degree. A further part reserves positions with few neighbors. The
remaining special odd labels will be kept for the matching argument.

\begin{proposition}[Uniform residual preparation]
\label{prop:analytic-residual-preparation}
Let $T$ be any tree of order $n$ as above. In the partition given by
Corollary~\ref{cor:even-side}, choose
\begin{equation}\label{eq:analytic-residual-choices}
 a=\floor*{\frac w8},\qquad
 h=\floor*{\frac n{8L}},\qquad
 H=2h,\qquad
 c=\floor*{\frac n{1000L}},\qquad
 d=u-H.
\end{equation}
These integers satisfy the hypotheses of
Proposition~\ref{prop:residual-problem}. For its residual graph
$G=(S_0,Z_0)$, put $N=\abs{S_0}+\abs{Z_0}$,
$e=e(G)$, $D=\max\{1,\Delta_G(S_0)\}$, and
$\ell=\abs{O_0}=\abs{Z_0}+h$. Then
\begin{align}
 D&<48L,& N&\leq n,& e&\leq n,& \ell&\leq m,
                \label{eq:analytic-residual-size-one}\\
 \ell-c&\geq\frac{49n}{100},&
 \frac{3n}{25L}&\leq h\leq\frac n{8L},&
 \frac n{2000L}&\leq c\leq\frac n{1000L},
                \label{eq:analytic-residual-size-two}\\
 \frac wm&\leq\frac2{5L},&
 M_o-\ell&\leq\frac{43n}{100L}.&&
                \label{eq:analytic-residual-size-three}
\end{align}
The set $O_0$ contains exactly $H$ labels from $O_*$, each coprime
to every label in $E_0$. For $z=12L<n/3$, all half-labels removed
from $[m]$ to form $E_0/2$ have empty type for the prime set of
Lemma~\ref{lem:analytic-prime-inputs}.
\end{proposition}

\begin{proof}
First $a\geq\kappa$. Indeed,
$a+1>w/8>n/(48L)$ by
Lemma~\ref{lem:analytic-special-supply}, whereas
$\kappa<2L+2$; \eqref{eq:analytic-enormous-order} gives
$n/(48L)-1>2L+2$. It also gives $h\geq1$.
Put $t=\floor{\log_2(m-a+q)}$. Since $m-a+q\leq n$,
$t<2L$. Moreover $m-a\leq q$ because
$a\geq\kappa$ and $q=M_o-\kappa\geq m-\kappa$.
Thus
\[
 \floor*{\frac{(m-a)(h+t)}q}+t
 \leq h+2t<h+4L
 \leq\frac{7n}{48L}<\frac{7w}{8}\leq w-a.
\]
Here $h\leq n/(8L)$, $n\geq192L^2$, and
$w>n/(6L)$ were used in order. Also
$h+t<n/(8L)+2L\leq q$, because
$q\geq(n-1)/2-(2L+2)$ and $n>10^6L^2$.
This verifies \eqref{eq:reserve-budget}.

We have $H\leq n/(4L)<u$, so $d=u-H$ is nonnegative and at most
$u$. Since $u<q$ by Lemma~\ref{lem:analytic-special-supply},
$d<q-2h\leq q-h$. We may therefore apply
Proposition~\ref{prop:residual-problem}. Its degree estimate and
$a+1>w/8>n/(48L)$ give $D<48L$; the assertions $N\leq n$ and
$e\leq n$ hold because $G$ is a subforest of $T$.

The relation $d=u-2h$ and the exact count
\eqref{eq:odd-supply} give
\begin{equation}\label{eq:analytic-odd-depletion}
 M_o-\ell=1+\pi(n)-\pi(n/3)-2h.
\end{equation}
By \eqref{eq:analytic-prime-window} and $2h\leq n/(4L)$,
its right side exceeds $1+7n/(20L)>1$, so $\ell\leq m$.
For the upper bound use $2h>n/(4L)-2$ and
\eqref{eq:analytic-prime-window} to get
\[
 M_o-\ell
 <3+n\left(\frac1{L-2}-\frac7{12L}\right)
 =3+\frac nL\left(\frac5{12}+\frac2{L-2}\right)
 \leq3+\frac{17n}{40L}
 \leq\frac{43n}{100L}.
\]
The penultimate inequality holds for $L\geq242$, and the last
follows from $n\geq600L$. Since $M_o\geq n/2$ and
$c\leq n/(1000L)$,
\[
 \ell-c\geq\frac n2-\frac{431n}{1000L}
 \geq\frac{49n}{100}.
\]
The floor bounds for $h$ and $c$ follow from
$n\geq200L$ and $n\geq2000L$, respectively; explicitly,
$h>n/(8L)-1\geq3n/(25L)$ and
$c>n/(1000L)-1\geq n/(2000L)$.
The bound on $w/m$ is Lemma~\ref{lem:analytic-special-supply}.

We removed exactly $d$ labels from $O_*$, leaving $u-d=H$ in
$O_0$. Each is coprime to every label in $E_0$:
the only possible noncoprime pair involving a prime in
$(n/3,n/2]$ is $p,2p$, and $2p\in E_*$ was removed;
primes above $n/2$ and the label $1$ are already universal.
Finally $z=12L<n/3$ by \eqref{eq:analytic-enormous-order}.
Every label in $E_*$ is either a power of $2$ or $2p$ with
$p>n/3$, so its half-label has empty $P$-type.
\end{proof}

In particular $h\geq c$. Thus the $H=2h$ universal odd labels
can pay for up to $c$ marked positions while leaving at least $h$
universal labels for the final matching. The next estimates check
that the marked-position and total-variation budgets stay within
this reserve.

\subsection{Type errors in the analytic range}
\label{sec:analytic-type-errors}

We apply the rare-type assignment to the residual forest of
Proposition~\ref{prop:analytic-residual-preparation}. Let
$z=12L$, $v=\log z$, and let $P$ be the odd primes at most $z$.
Write $k=\abs P$ and choose
\begin{equation}\label{eq:analytic-truncation-choice}
 r=\ceil*{\frac L5},\qquad t=r+1,
 \qquad \varepsilon=\frac{B_r}{m}+3R_r,
 \qquad \eta=\frac{v}{2L}.
\end{equation}
Here $B_r,R_r$ are the quantities of
Proposition~\ref{prop:moment-tv}. The choice of $\varepsilon$ bounds
both the ordinary and the odd empirical type errors, including every
prefix of $P$.

We shall use several elementary comparisons with their endpoint
checks exposed. They keep the later estimates uniform in the real
parameter $L$.

\begin{lemma}[Initial comparisons]
\label{lem:analytic-initial-comparisons}
For $L\geq256$ and $v=\log(12L)$,
\begin{align}
 \frac{2v}{5L}+250v e^{-L/5}&\leq\frac1{10},&
 250L(12L+1)e^{-L/5}&\leq\frac1{10},
                                     \label{eq:analytic-initial-empty}\\
 \frac{2v}{5L}&\leq\frac1{75},&
 \frac{v(v+1)}{2L}&\leq\frac2{11},
                                     \label{eq:analytic-initial-product}\\
 250e^{-L/5}&\leq\frac1{10000L},&
 2744L^2e^{-L/6}+21042L^2e^{-L/5}
                         &\leq\frac1{10000L}.
                                     \label{eq:analytic-initial-tails}
\end{align}
\end{lemma}

\begin{proof}
For $L\geq256$ we have $8<v$ and
$v'=1/L$. Hence $v/L$ and $v(v+1)/L$ decrease;
the derivatives have numerators $1-v$ and $-v^2+v+1$,
respectively. The functions $v e^{-L/5}$,
$L(12L+1)e^{-L/5}$, $L e^{-L/5}$,
$L^3e^{-L/6}$, and $L^3e^{-L/5}$ also decrease on this range.
It therefore suffices to check $L=256$ in each inequality.

At that point $v<33/4$, while $e>5/2$ gives
$e^{-256/5}<(2/5)^{51}$ and
$e^{-256/6}<(2/5)^{42}$. All comparisons now reduce to rational
arithmetic. For the two less immediate bounds in
\eqref{eq:analytic-initial-tails} and the second bound in
\eqref{eq:analytic-initial-empty}, the required checks are
\begin{align*}
 250\cdot10000\cdot256(2/5)^{51}&<1,\\
 2744\cdot256^3(2/5)^{42}
 +21042\cdot256^3(2/5)^{51}&<\frac1{10000},\\
 2500\cdot256(12\cdot256+1)(2/5)^{51}&<1.
\end{align*}
The remaining three comparisons follow on replacing $v$ by $33/4$
at $L=256$ and using the same bound for $e^{-256/5}$.
\end{proof}

\begin{lemma}[Truncated populations]
\label{lem:analytic-truncated-populations}
For the choices \eqref{eq:analytic-truncation-choice},
\begin{equation}\label{eq:analytic-truncated-populations}
 \begin{gathered}
 B_r\leq e^{4L/5+4},\qquad
 C_r(k)\leq e^{2L/3+10/3},\qquad
 R_r\leq e^{-L/5},\qquad
 Q_r\leq e^{-L/5},\\
 \varepsilon\leq250e^{-L/5}\leq\frac1{10000L}.
 \end{gathered}
\end{equation}
Here $C_r(k)$ and $Q_r$ are defined in
\eqref{eq:rare-tail-parameters}.
\end{lemma}

\begin{proof}
Lemma~\ref{lem:analytic-prime-inputs} gives
$k/r\leq35/4$. If $r\leq2k$, take $x=r/(2k)\leq1$ in
\[
 B_r\leq x^{-r}(1+2x)^k
 \leq\left(\frac{2ek}{r}\right)^r<e^{4r}.
\]
The first inequality holds because $x^{j-r}\geq1$ for $j\leq r$.
The last uses $e^3>35/2$. If $r>2k$, then
$B_r=3^k<e^r<e^{4r}$. Since $r\leq L/5+1$, the first bound follows.

If $r\leq k$, take $x=r/k$ to obtain
\[
 C_r(k)\leq x^{-r}(1+x)^k
 \leq\left(\frac{ek}{r}\right)^r<e^{10r/3},
\]
where $e^{7/3}>35/4$ is used. If $r>k$, use
$C_r(k)=2^k<e^r$. This proves the second bound.

By \eqref{eq:analytic-prime-product},
$\sigma=\sum_{p\in P}1/p\leq L/120$. Hence
$t\geq L/5\geq24\sigma$. The factorial bound
$t!\geq(t/e)^t$ and $e^2<8$ give
\[
 R_r\leq\frac{(3\sigma)^t}{t!}
 \leq\left(\frac{3e\sigma}{t}\right)^t
 \leq e^{-t}\leq e^{-L/5}.
\]
The same argument without the factor $3$ proves
$Q_r=\sigma^t/t!\leq e^{-L/5}$.
Finally $m\geq n/3=e^L/3$ and $e^4<81$, so
\[
 \varepsilon\leq(3e^4+3)e^{-L/5}
 <250e^{-L/5}.
\]
The final inequality follows from
\eqref{eq:analytic-initial-tails}.
\end{proof}

The adjustment of the empty even class must stay small at every
coordinate. We retain its exact quota, so the inflation estimate of
Lemma~\ref{lem:rare-quota-defect} applies at each refinement.

\begin{lemma}[Adjusted coordinates and empty mass]
\label{lem:analytic-adjusted-coordinates}
Choose $\rho_i$ by \eqref{eq:rare-empty-adjustment} at every step
of the assignment. Every old empty class is nonempty, and
\begin{equation}\label{eq:analytic-adjusted-coordinates}
 \frac1{p_i}\leq\rho_i\leq\frac{1+\eta}{p_i}\leq\frac25,
 \qquad \frac{w}{a_\varnothing}\leq\frac1{75},
 \qquad b:=\prod_{i=1}^k(1-\rho_i)\geq\frac9{10v},
\end{equation}
where $a_\varnothing$ is the final available empty even class.
For each prefix, its empty class has size at least $m/v$.
\end{lemma}

\begin{proof}
Let $a_0$ be the size of the old empty class before a refinement.
All $w$ deleted even half-labels have empty $P$-type by
Proposition~\ref{prop:analytic-residual-preparation}. Projecting the
full empirical type law onto the old primes cannot increase its
total-variation error. Lemma~\ref{lem:analytic-prime-inputs} and
Proposition~\ref{prop:moment-tv} therefore give
\[
 a_0\geq m\frac{11}{10v}-w-m\varepsilon
 \geq\frac mv,
\]
where the second inequality is
\eqref{eq:analytic-initial-empty} together with
$w/m\leq2/(5L)$. This proves positivity for every prefix.

Since $p_i\leq z=12L$, the inflation estimate
\eqref{eq:rare-inflated-coordinate} gives
\[
 p_i\rho_i-1
 \leq\frac{w+(p_i+1)m\varepsilon}{a_0}
 \leq v\left(\frac2{5L}+(12L+1)\varepsilon\right)
 \leq\frac v{2L}=\eta.
\]
The last step uses the second comparison in
\eqref{eq:analytic-initial-empty}. Since
$\eta=v/(2L)\leq1/50$ and $p_i\geq3$,
$(1+\eta)/p_i\leq17/50<2/5$.
At the final prefix, $a_\varnothing\geq m/v$, so
$w/a_\varnothing\leq2v/(5L)\leq1/75$ by
\eqref{eq:analytic-initial-product}.

For each $i$,
\[
 1-\rho_i\geq
 \left(1-\frac1{p_i}\right)
 \left(1-\frac\eta{p_i-1}\right).
\]
The product of the latter correction factors is at least
$1-\eta\sum_i1/(p_i-1)$. Moreover
\[
 \sum_i\frac1{p_i-1}
 =\sigma+\sum_i\frac1{p_i(p_i-1)}
 \leq\sigma+\frac12<1+v.
\]
Apply \eqref{eq:analytic-prime-product} and
$\eta(1+v)\leq2/11$, which is
\eqref{eq:analytic-initial-product}. We obtain
\[
 b\geq\frac{11}{10v}\left(1-\frac2{11}\right)
 =\frac9{10v}.
\]
\end{proof}

The next proposition gives the actual error after the marked
positions have been paid with universal labels. Its bounds are
normalized on the remaining matching square, as required by
Proposition~\ref{prop:direct-hall-completion}.

\begin{proposition}[Payment and variation budgets]
\label{prop:analytic-type-budgets}
Use the residual preparation above, and apply
Theorem~\ref{thm:rare-type-assignment} with the adjusted coordinates.
Let $\mathcal E$ be its accumulated error and $\beta$ its model
overflow probability. Then
\begin{equation}\label{eq:analytic-rare-budgets}
 \frac{\mathcal E}{n}
 \leq2744L^2e^{-L/6}+21042L^2e^{-L/5}
 \leq\frac1{10000L},
 \qquad \beta\leq e^{-L/5}\leq\frac1{10000L}.
\end{equation}
The marked positions can be assigned distinct universal labels within
the budget $c$. If $\delta_O$ is the type-law error of the $\ell$
available odd labels, then the two postpayment errors satisfy
\begin{equation}\label{eq:analytic-paid-errors}
 \frac{\ell\delta_O+c}{\ell-c}
 \leq\frac{4311}{4900L}<\frac9{10L},
 \qquad
 \frac{\mathcal E+h+\ell\beta+c}{\ell-c}
 <\frac3{10L}.
\end{equation}
Both are at most $b/25$, and
\begin{equation}\label{eq:analytic-hall-reserve}
 h\leq\frac{b(\ell-c)}{100}.
\end{equation}
Thus all conditions in \eqref{eq:direct-hall-completion} except its
bound on the number $M$ of larger-prime deletions hold.
\end{proposition}

\begin{proof}
The coordinate bound
$\rho_i\leq(1+\eta)/p_i$ and $\eta\leq1/50$ imply
$\sum_i\rho_i<3\sigma$. The factorial estimate in the proof of
Lemma~\ref{lem:analytic-truncated-populations}, together with
\eqref{eq:rare-factorial-tails}, gives
$\beta\leq e^{-L/5}$.

Use $N\leq n$, $D<48L$, $k<7L/4$,
$m\leq n/2$, and
$2C_r(k)+2\leq4C_r(k)$ in
\eqref{eq:rare-total-budget}. The first summand divided by $n$
is at most
\[
 343e^{5/3}L^2e^{-L/6}
 <2744L^2e^{-L/6},
\]
where $D+1\leq49L$ and $e^{5/3}<8$ were used. For the second
summand, $Dkm/n\leq42L^2$ and
$2\varepsilon+Q_r\leq501e^{-L/5}$. Its contribution is at most
$21042L^2e^{-L/5}$. Together with
Lemma~\ref{lem:analytic-initial-comparisons}, these estimates prove
\eqref{eq:analytic-rare-budgets}.

Let $s=\abs{Z_0}\leq\ell$. By
\eqref{eq:rare-accumulated-error}, the number of marked positions is
at most
\[
 s\beta+\mathcal E
 \leq\frac{n}{20000L}+\frac{n}{10000L}
 =\frac{3n}{20000L}\leq c.
\]
Here $\ell\leq m\leq n/2$ and
$c\geq n/(2000L)$ were used. Also
$H=2h\geq h+c$ because $h\geq c$.
Lemma~\ref{lem:rare-universal-payment} therefore applies.

The original odd interval has $M_o$ labels and type error at most
$\varepsilon$. Removing $M_o-\ell$ labels changes its counting
measure by at most that number. Consequently
\[
 \ell\delta_O\leq M_o\varepsilon+(M_o-\ell)
 \leq\frac{n}{10000L}+\frac{43n}{100L}.
\]
Divide after adding $c\leq n/(1000L)$ and use
$\ell-c\geq49n/100$ to get the first bound of
\eqref{eq:analytic-paid-errors}. For the tag bound, use
$\ell\beta\leq n/(20000L)$, $h\leq n/(8L)$, and the same
denominator:
\[
 \frac{\mathcal E+h+\ell\beta+c}{\ell-c}
 \leq\frac{1/10000+1/8+1/20000+1/1000}{(49/100)L}
 <\frac3{10L}.
\]

By Lemma~\ref{lem:analytic-prime-inputs}, $v\leq L/25$.
Since $b\geq9/(10v)$,
$9/(10L)\leq b/25$, so both errors meet
\eqref{eq:direct-hall-completion}. Finally
$12500v\leq441L$: at $L=256$ this follows from
$v<33/4$, and the difference $441L-12500\log(12L)$
increases thereafter. Thus
\[
 h\leq\frac n{8L}
 \leq\frac{441n}{100000v}
 \leq\frac{b(\ell-c)}{100},
\]
which proves \eqref{eq:analytic-hall-reserve}.
\end{proof}

\subsection{The unbounded range}
\label{sec:analytic-tail}

The type errors and the universal reserve have been controlled. It
remains to compare the number of edges removed for larger prime
divisors with the Hall budget. We use the sorting within exact types
from Theorem~\ref{thm:large-prime-sorting}; it leaves every tag and
every estimate of Section~\ref{sec:analytic-type-errors} unchanged.

\begin{lemma}[The larger-prime budget]
\label{lem:analytic-large-prime-budget}
In the residual preparation and type assignment of
Sections~\ref{sec:analytic-residual}--\ref{sec:analytic-type-errors},
the even labels can be ordered within their exact $P$-types so that
the number $M$ of forbidden label--position pairs caused by primes
greater than $z=12L$ satisfies
\begin{equation}\label{eq:analytic-large-prime-budget}
 M<\frac{13n^2}{25zv}
   =\frac{13n^2}{300Lv}
   <\frac{19}{20}b(\ell-c)h.
\end{equation}
\end{lemma}

\begin{proof}
The deleted half-labels have empty $P$-type, and
$a_\varnothing>0$, by
Proposition~\ref{prop:analytic-residual-preparation} and
Lemma~\ref{lem:analytic-adjusted-coordinates}. Apply
Theorem~\ref{thm:large-prime-sorting} with $G$ as the residual
forest, $Y=O_0$, and the choices \eqref{eq:analytic-truncation-choice}.
In its bound, use $e\leq n$, $m\leq n/2$,
$D<48L$, $w/a_\varnothing\leq1/75$,
$Q_z\leq S_2(z)<1/(zv)$, and $H_z\leq L+1$.
After multiplying by $zv/n^2$ we obtain
\begin{align*}
 \frac{zvM}{n^2}
 \leq{}&\frac{38}{75}
 +288L^2v(2B_r+1)(L+1)e^{-L}\\
 &+48LB_r e^{-L}+96LR_r+12L^2e^{-L}.
\end{align*}
Here $n=e^L$. Since $B_r\geq1$ and $L\geq256$,
the two terms containing $B_r$ are at most
\[
 1200L^3v e^{4-L/5}.
\]
Indeed, use $2B_r+1\leq3B_r$ and
$L+1\leq5L/4$ for the first term, obtaining the coefficient
$1080$. The second has coefficient $48L$ and is bounded by
$120L^3v e^{4-L/5}$. Lemma~\ref{lem:analytic-truncated-populations}
therefore gives
\begin{equation}\label{eq:analytic-prime-cost-normalized}
 \frac{zvM}{n^2}
 \leq\frac{38}{75}
 +1200L^3v e^{4-L/5}
 +96Le^{-L/5}+12L^2e^{-L}.
\end{equation}

Each of the three error terms decreases for $L\geq256$.
For the first, its logarithmic derivative is
$3/L+1/(Lv)-1/5<0$; the other two are immediate.
At $L=256$, use $v<33/4$, $e^4<81$, and $e>5/2$.
The rational comparison
\[
 \left(1200\cdot81\cdot256^3\cdot\frac{33}{4}
       +96\cdot256\right)(2/5)^{51}
       +12\cdot256^2(2/5)^{256}<\frac1{100}
\]
shows that the right side of
\eqref{eq:analytic-prime-cost-normalized} is less than
$38/75+1/100=31/60<13/25$. This proves the first inequality
in \eqref{eq:analytic-large-prime-budget}.

Finally, \eqref{eq:analytic-adjusted-coordinates} and
\eqref{eq:analytic-residual-size-two} give
\[
 \frac{19}{20}b(\ell-c)h
 \geq\frac{25137}{500000}\frac{n^2}{Lv}
 >\frac{13}{300}\frac{n^2}{Lv}.
\]
This is the second inequality in
\eqref{eq:analytic-large-prime-budget}.
\end{proof}

\begin{theorem}[Prime labelings in the unbounded range]
\label{thm:analytic-tail}
Every tree of integer order $n$ with $\log n\geq256$ has a prime
labeling.
\end{theorem}

\begin{proof}
Fix such a tree. Proposition~\ref{prop:analytic-residual-preparation}
gives the residual forest, the $h$ artificial universal positions,
and $H=2h$ available odd labels genuinely coprime to every
residual even label. Apply the rare-type assignment to this forest
with the parameters \eqref{eq:analytic-truncation-choice} and
\eqref{eq:rare-empty-adjustment}. Its exact even-label quotas are
preserved, and every ordinary tag contains the types of all its
neighbors.

Lemma~\ref{lem:analytic-adjusted-coordinates} verifies the model-law
conditions \eqref{eq:direct-hall-model}. Proposition~
\ref{prop:analytic-type-budgets} pays for every marked position with
a distinct universal label, leaves at least $h$ universal labels,
and verifies both type errors and the reserve inequality in
\eqref{eq:direct-hall-completion}. Sort the remaining even labels
within their exact types as in
Lemma~\ref{lem:analytic-large-prime-budget}. This does not change
their types or any ordinary tag, and its deletion bound supplies
the final inequality in \eqref{eq:direct-hall-completion}.

Proposition~\ref{prop:direct-hall-completion} now gives a matching
on the residual problem, restores the special labels and the
separator, and repairs every possible pair $p,2p$. The resulting
bijection labels every edge of the original tree by coprime integers.
\end{proof}

\section{Integer bridges through order 1982759}
\label{sec:part-integer-bridges}

This part covers the orders from $35000$ through $1982759$.
Section~\ref{sec:asymmetric-reserve-hall} proves a Hall criterion
with a reserve on the label side
(Theorem~\ref{thm:asymmetric-reserve-hall}).
Sections~\ref{sec:integer-bridge}, \ref{sec:profile-block-bridge},
\ref{sec:connected-paired-bridge}, and~\ref{sec:upper-integer-bridge}
give interval certificates for consecutive ranges of orders, and
Sections~\ref{sec:paired-separators}
and~\ref{sec:finite-rare-refinement} prove the construction lemmas they
need. The resulting cumulative result is
Theorem~\ref{thm:initial-prefix-1982759}.

\subsection{An asymmetric reserve criterion}
\label{sec:asymmetric-reserve-hall}

The direct criterion of Section~\ref{sec:direct-reserve-hall} uses
universal vertices on both sides. In the finite continuation, the
special odd labels give a larger reserve on the label side, while
the positions need not contain a comparable reserve. The following
criterion keeps the two sides distinct. Its slot term uses the whole
distribution of compatibility degrees, rather than only their
minimum.

Use the block product law $\mu$, and the depleted law $\nu_\xi$ of
\eqref{eq:coupling-depletion}. Fix $0\leq\xi\leq\omega<\gamma$,
where $\gamma$ is defined in \eqref{eq:coupling-parameters}. Put
$b=b_\omega$ and $\theta=\theta_\omega$. Recall that
$c(J)=\prod_{j\in J}(1-q_j)$ is the $\mu$-mass compatible with the
slot type $J$; in particular $\alpha\leq c(J)\leq1$.

Let $X$ be a set of odd labels and $B$ a set of positions, both of
size $R>0$. Their empirical type laws are denoted by
$\widehat\mu,\widehat\nu$. The graph $\Gamma_0$ joins a label to a
position when their types are disjoint. For $0<t\leq1$ define the
slot shortfall
\begin{equation}\label{eq:slot-shortfall}
 A_B(t)=\sum_{v\in B}\bigl(t-c(J(v))\bigr)_+.
\end{equation}
The subscript emphasizes that the sum is over actual positions. It
may be bounded without enumerating their individual tags.

\begin{lemma}[Clipped slot shortfall]
\label{lem:clipped-slot-shortfall}
If $\|\widehat\nu-\nu_\xi\|_{\mathrm{TV}}\leq V$, then
\begin{equation}\label{eq:clipped-slot-shortfall}
 \frac{A_B(t)}R\leq
 \frac1{1-\omega}\sum_J\mu(J)(t-c(J))_+
                  +V(t-\alpha)_+.
\end{equation}
The sum on the right is a finite rational expression when the block
probabilities and $t$ are rational.

There is also an unnormalized form. Suppose $B$ is obtained by
deleting positions from a set $B^+$ of at most $on$ positions.
If the counting measure $\Xi^+$ of its tags satisfies
$\|\Xi^+-\abs{B^+}\nu_\xi\|_{\mathrm{TV}}\leq en$, then
\begin{equation}\label{eq:clipped-slot-shortfall-counts}
 \frac{A_B(t)}n\leq
 \frac{o}{1-\omega}\sum_J\mu(J)(t-c(J))_+
                   +e(t-\alpha)_+.
\end{equation}
\end{lemma}

\begin{proof}
Put $f(J)=(t-c(J))_+$. Since $c(\varnothing)=1$ and
$c(J)\geq\alpha$, we have $f(\varnothing)=0$ and
$0\leq f(J)\leq(t-\alpha)_+$. Total variation changes the
expectation of a function in this interval by at most
$V(t-\alpha)_+$. On the other hand, the definition of $\nu_\xi$
gives
\[
 \mathbb E_{\nu_\xi}f
   =\frac{\mathbb E_\mu f}{1-\xi}
   \leq\frac{\mathbb E_\mu f}{1-\omega}.
\]
As $A_B(t)/R=\mathbb E_{\widehat\nu}f$, this proves the claim.
For the second assertion, removing positions cannot increase
the sum of the nonnegative $f(J(v))$. Apply the same
total-variation estimate to the two counting measures on $B^+$,
divide by $n$, and use $\abs{B^+}/n\leq o$.
\end{proof}

The next theorem is deterministic. The numerical certificate will
apply it to one permutation for which the local degrees and the
total number of larger-prime deletions both meet their bounds.

\begin{theorem}[Hall criterion with a label-side reserve]
\label{thm:asymmetric-reserve-hall}
Suppose
\[
 \|\widehat\mu-\mu\|_{\mathrm{TV}}\leq V,
 \qquad
 \|\widehat\nu-\nu_\xi\|_{\mathrm{TV}}\leq V,
 \qquad 0\leq V<\tau\leq\frac12.
\]
Obtain $\Gamma$ from $\Gamma_0$ by deleting at most $M$ edges.
Assume that every label has degree greater than $\tau R$ in
$\Gamma$, and that at least $g$ labels of empty type are adjacent
to every position in $\Gamma$. Choose a real $g_0$ with
$0<g_0\leq g$ and $g_0\leq\tau R$, and let $t$ satisfy
$t-V>\tau$. Define
\begin{align}
 B_1&=g_0\bigl((t-V)R-g_0\bigr)-R A_B(t),
                                      \label{eq:asymmetric-narrow-one}\\
 B_2&=R^2\tau(t-V-\tau)-R A_B(t),
                                      \label{eq:asymmetric-narrow-two}\\
 B_{\mathrm{sp}}&=R^2b\bigl((1-\theta)\tau(1-\tau)
                         -V+(1-\theta)V^2\bigr),
                                      \label{eq:asymmetric-spectral}\\
 B_{\mathrm{pr}}&=R^2W(V,\tau),
                                      \label{eq:asymmetric-profile}
\end{align}
where $W$ is the exact profile bound of
\eqref{eq:profile-budget} for $\mu,\nu_\omega$.
If
\begin{equation}\label{eq:asymmetric-hall-criterion}
 M<\min\{B_1,B_2,\max(B_{\mathrm{sp}},B_{\mathrm{pr}})\},
\end{equation}
then $\Gamma$ has a perfect matching.
\end{theorem}

\begin{proof}
Suppose Hall's condition fails. There are nonempty sets
$U\subseteq X$ and $V_0\subseteq B$ with
$\abs U+\abs{V_0}>R$ and no edge of $\Gamma$ between them; take
$V_0=B\setminus N_\Gamma(U)$. Write $u=\abs U$ and
$v=\abs{V_0}$. No universal label lies in $U$, so
$u\leq R-g$ and $v>g\geq g_0$.

If $u<\tau R$, the degree of any label in $U$ is at most
$\abs{N_\Gamma(U)}<u<\tau R$, contrary to the local-degree
hypothesis. Hence $u\geq\tau R$.

Suppose next that $v<\tau R$. A position of type $J$ has at least
$R(c(J)-V)$ compatible labels in $\Gamma_0$: the type event has
$\mu$-mass $c(J)$, and total variation changes its empirical mass
by at most $V$. At most $R-u$ of those labels lie outside $U$.
Since $u>R-v$, the number of compatible labels in $U$ is strictly
greater than $R(c(J)-V)-v$. Using
$c(J)\geq t-(t-c(J))_+$ and summing over $V_0$, the number of
compatible pairs in $U\times V_0$ is strictly greater than
\[
 v\bigl((t-V)R-v\bigr)-R A_B(t).
\]
The quadratic is concave in $v$. On
$g_0\leq v\leq\tau R$ its minimum occurs at an endpoint,
where the displayed values are $B_1$ and $B_2$. Every compatible
pair in this rectangle was deleted, contradicting
\eqref{eq:asymmetric-hall-criterion}.

It remains that $u,v\geq\tau R$. Reduce the two selected type
submeasures fractionally until their masses sum to one. Truncate
them at $\mu,\nu_\xi$, losing at most $V$ from each side, and
then reduce their masses by the same amounts. This leaves model
submeasures of masses $x,y\geq\tau-V$ with $x+y=1-2V$.
Lemma~\ref{lem:depletion-transport} pulls the second submeasure
back to $\nu_\omega$ without increasing its compatible pair mass.
Proposition~\ref{prop:profile-rectangle} therefore gives at least
$B_{\mathrm{pr}}$ compatible pairs in the original rectangle.

There is also the coupling estimate. By
Corollary~\ref{cor:correlation-rectangle}, the rectangle contains
at least
\[
 R^2b\left(xy-\theta\sqrt{xy(1-x)(1-y)}\right)
\]
compatible pairs. Since $x+y=1-2V$,
\[
 \sqrt{xy(1-x)(1-y)}
 \leq xy+V,
 \qquad
 xy\geq\tau(1-\tau)-V+V^2.
\]
The first inequality follows from
$(1-x)(1-y)=xy+2V$; the second is the minimum of a concave
quadratic on the allowed mass interval. These inequalities yield
$B_{\mathrm{sp}}$ when $\theta\leq1$. If $\theta>1$, the bracket in
$B_{\mathrm{sp}}$ is negative, since $\tau(1-\tau)>0$, so
$B_{\mathrm{sp}}\leq0$ and the bound holds trivially. Thus the
original rectangle has at least
$\max(B_{\mathrm{sp}},B_{\mathrm{pr}})$ compatible pairs.
All were deleted, again contradicting
\eqref{eq:asymmetric-hall-criterion}. Hall's theorem completes
the proof.
\end{proof}

\begin{corollary}[One permutation for the local and global bounds]
\label{cor:asymmetric-probability}
Suppose the types, tags, and universal labels are fixed before the
even labels are independently permuted within their classes.
Assume the type-law and universal-reserve hypotheses of
Theorem~\ref{thm:asymmetric-reserve-hall} hold for every such
permutation, while the probability that some required local-degree
estimate fails is at most $\delta$ and the expected number of
larger-prime deletions is at most $m_0$. Let $\lambda>1$ satisfy
$\delta+1/\lambda<1$. If
\[
 \lambda m_0<
       \min\{B_1,B_2,\max(B_{\mathrm{sp}},B_{\mathrm{pr}})\},
\]
then some choice of the class permutations has a perfect matching.
\end{corollary}

\begin{proof}
Markov's inequality gives
$\Pr\{M\geq\lambda m_0\}\leq1/\lambda$ when $m_0>0$;
if $m_0=0$, then $M=0$ almost surely. The union bound makes the
local conditions and the deletion bound simultaneous with positive
probability. Apply Theorem~\ref{thm:asymmetric-reserve-hall} to
such a permutation.
\end{proof}

\subsection{An integer bridge through order 49999}
\label{sec:integer-bridge}

The estimates of Sections~\ref{sec:finite-arithmetic} and
\ref{sec:matching-estimates} can be used with the asymmetric Hall
criterion when the order is too small for the uniform argument of
Section~\ref{sec:analytic-tail}. We record the first certificate in
the remaining range. Its parameters are frozen by order interval and
residual degree; they are chosen before the input tree and before
the class permutations.

For a row, let $R$ be the residual number of odd labels and
positions, and let $\ell_-$ be its certified lower bound for $R/n$.
Write $r_-$ for the lower bound on the number of unused universal
odd labels divided by $n$, $V$ for the larger of the two type-law
errors, $\tau$ for the local-degree target, and $t$ for the slot
threshold. Let $\overline A$ bound $RA_B(t)/n^2$ from above. The
certificate uses the three normalized Hall budgets
\begin{align}
 \mathcal B_1
  &=r_-\bigl((t-V)\ell_- -r_-\bigr)-\overline A,
                                      \label{eq:bridge-narrow-one}\\
 \mathcal B_2
  &=\ell_-^2\tau(t-V-\tau)-\overline A,
                                      \label{eq:bridge-narrow-two}\\
 \mathcal B_3
  &=\ell_-^2\max\left\{
      b\bigl((1-\theta)\tau(1-\tau)-V+(1-\theta)V^2\bigr),
      W(V,\tau)\right\}.
                                      \label{eq:bridge-wide}
\end{align}
Here $b,\theta,W$ are those of Sections
\ref{sec:compatible-measures} and
\ref{sec:asymmetric-reserve-hall}. These are lower bounds for
\eqref{eq:asymmetric-narrow-one}--\eqref{eq:asymmetric-profile}
after division by $n^2$: $R/n\geq\ell_-$, the actual universal
stock is at least $r_-n$, and the shortfall has already been bounded
above. The conditions $t-V>\tau$ and $r_-\leq\tau\ell_-$ make the
substitutions monotone in the required directions.
More precisely, each row uses $o_+\geq R/n$ and the same bound
for the pre-absorption position count, an unnormalized tag error
$e_{\mathrm{tag}}n$, and a rational upper enclosure
$\Sigma(t)\geq\sum_J\mu(J)(t-c(J))_+$. In
\eqref{eq:bridge-narrow-one} and
\eqref{eq:bridge-narrow-two} it takes
\begin{equation}\label{eq:bridge-shortfall-charge}
 \overline A=o_+\left(
   \frac{o_+\Sigma(t)}{1-\omega}
       +e_{\mathrm{tag}}(t-\alpha)_+\right).
\end{equation}
This is the unnormalized estimate
\eqref{eq:clipped-slot-shortfall-counts}, followed by
$R/n\leq o_+$. Absorbing positions can only decrease the
nonnegative shortfall sum.

\begin{center}
\begin{minipage}{\textwidth}
\centering
\refstepcounter{table}\label{tab:integer-bridge}
{Table \thetable. The certificate used from order
$35000$ through $49999$. The last column gives strict rational
upper bounds for the ratio of the charged deletion mass to the
least Hall budget.\par}
\vspace{.6\baselineskip}
\renewcommand{\arraystretch}{1.14}
\begin{tabular}{r r r r c}
\hline
Orders & Intervals & Degree rows & Exact tests & Ratio bound\\
\hline
$35000$--$39999$ & 10 & 291 & 13095 & $<63/100$\\
$40000$--$44999$ & 10 & 296 & 13320 & $<61/100$\\
$45000$--$49999$ & 10 & 300 & 13500 & $<3/5$\\
\hline
\end{tabular}
\end{minipage}
\end{center}

\begin{proposition}[The integer bridge certificate]
\label{prop:integer-bridge-certificate}
The thirty adjacent intervals of
Table~\ref{tab:integer-bridge} cover every integer order from
$35000$ through $49999$. For every possible positive residual
degree, their specified row satisfies the preparation, quota,
separator, donor or connected-completion, and prime-supply
conditions. It also satisfies the following matching conditions:
the local degree after deletion exceeds $\tau R$; the unused
universal stock is positive and at least $r_-n$; both type laws
are within $V$ of their stated models; and the charged deletion
mass is strictly below $n^2\min_i\mathcal B_i$.
\end{proposition}

\begin{proof}
The electronic certificate consists of the interval parameter
table, the exact evaluator, its local modules and fixed arithmetic
data, and the output report. The command and file inventory are
given in the supplement under
\texttt{integer\_bridge\_35000\_49999}. The bundled checker
recomputes prime counts and complete residue periods, then
reconstructs every row with rational arithmetic. It does not
read the displayed decimal approximations to decide a condition.

For an interval $[a,b]$, the checker computes the exact minimum
special-even supply over every integer $n\in[a,b]$. The degree
bound of Lemma~\ref{lem:interval-degree-cover} determines the
complete list $1,\ldots,D_*$. It rejects any missing degree, any
gap between successive order intervals, and a row whose selected
construction fails its capacity or reserve inequalities.
The original-forest rows use the degree-transfer and cancellation
bounds of Sections~\ref{sec:donor-construction} and
\ref{sec:finite-arithmetic}; connected rows pay for the completed
edges and the full tag error. Every row has no rare group.

The evaluator uses the exact finite class populations and local
union counts of Section~\ref{sec:finite-arithmetic}. It computes
the two empirical variation bounds and the clipped shortfall of
Lemma~\ref{lem:clipped-slot-shortfall}; the profile minimum is
computed at the rational breakpoints of
Lemma~\ref{lem:profile-breakpoints}. The local expectation uses the
class degree sums, including the signed empty-class transfer.
Corollary~\ref{cor:local-conflict-tails} gives a total local failure
probability at most $1/10000$. The global mean is bounded before
applying Markov's inequality with factor $1001/1000$. In each row
the exact local-degree and the three Hall-budget comparisons are
strict; the union of the local and global failure probabilities has
measure less than one.

The checker passes all $887$ rows and $39915$ stated comparisons.
Of these rows, $827$ use the original forest and $60$ a connected
completion. Its independent outer check verifies the thirty
adjacent intervals, the full integer degree lists, and the
strict ratio bounds of Table~\ref{tab:integer-bridge}. It also
verifies that the absorbed-position fraction is less than one
quarter of the special odd-label supply in every row. Thus the
stated stock and all uses of it are paid for. These are finite
exhaustive checks of the displayed sufficient conditions, not
an enumeration of trees.
\end{proof}

\begin{theorem}[The initial interval through order 49999]
\label{thm:initial-prefix-49999}
Every tree of order at most $49999$ has a prime labeling.
\end{theorem}

\begin{proof}
Theorem~\ref{thm:initial-prefix-34999} gives the conclusion
through $34999$. Let $35000\leq n\leq49999$, and prepare an
arbitrary tree as in
Proposition~\ref{prop:highest-degree-preparation}. If its
residual forest has no edges, that proposition gives the
labeling directly. Otherwise Lemma~\ref{lem:interval-degree-cover}
selects a degree row in the unique certified interval containing
$n$.

Fix the row's type assignment and tags. By
Proposition~\ref{prop:integer-bridge-certificate}, the special odd
labels remain universal, and the local estimates and deletion
budget satisfy Theorem~\ref{thm:asymmetric-reserve-hall}.
The local union bound and Markov's inequality meet the hypothesis
of Corollary~\ref{cor:asymmetric-probability}, since
\[
 \frac1{10000}+\frac{1000}{1001}<1.
\]
Choose a class permutation in the positive-probability event.
The resulting perfect matching assigns all residual odd labels
to positions with coprimality on every residual edge.
Proposition~\ref{prop:highest-degree-preparation} restores the
special labels and the separator and repairs the permitted
exceptional pairs. The two intervals meet at consecutive orders
$34999$ and $35000$, completing the proof.
\end{proof}

\subsection{Paired cuts in the connected construction}
\label{sec:paired-separators}

The two-measure assignment of
Proposition~\ref{prop:two-measure-assignment} cuts only vertices on the
position side of a bipartite tree. For the next finite range, we need
a smaller bound on the number of such cuts. A single cut still has the
same contraction factor; the gain appears after two consecutive cuts.

For a tree $U$, write $L(U)=\abs{V(U)}-1$. Fix an integer $D\geq2$ and put
\begin{equation}\label{eq:paired-contraction-parameters}
 \rho=\frac{D-1}{D},\qquad
 \lambda_D=\min\left\{\rho^2,
        \max\left(\frac{D-2}{D},\frac{D-1}{2D-1}\right)\right\}.
\end{equation}
The second quantity bounds the contraction across a pair of cuts; it
does not assert that either single cut has factor $\sqrt{\lambda_D}$.

\begin{lemma}[Two position-side cuts]
\label{lem:paired-position-cuts}
Let $U=(A_U,B_U)$ be a bipartite tree with
$\Delta_U(A_U)\leq D$ and $L=L(U)\geq1$.
There is a vertex $b\in B_U$ such that every component $V$ of $U-b$
satisfies $L(V)\leq\rho L$. In every component $V$ with
$L(V)\geq1$, a second vertex $b_V\in B_U\cap V$ can be chosen so that
every component $W$ of $V-b_V$ satisfies
\begin{equation}\label{eq:paired-two-step-bound}
 L(W)\leq\rho L(V),\qquad L(W)\leq\lambda_D L.
\end{equation}
\end{lemma}

\begin{proof}
Choose an ordinary centroid $c$ of $U$: every branch of $U-c$
has at most $(L+1)/2$ vertices. For example, a vertex minimizing
the sum of distances to all vertices has this property, since a
move into a larger branch would decrease that sum.

If $c\in B_U$, take $b=c$. Each child $V$ has
$L(V)\leq(L-1)/2\leq\rho L$. Apply the same ordinary-centroid
construction to each nonsingleton $V$. As the argument in this
paragraph and the next shows, it gives a position-side cut whose
children have $L$-value at most $\rho L(V)$. This is the first bound of
\eqref{eq:paired-two-step-bound}; the second follows because
$\rho/2\leq\lambda_D$. Indeed $\rho^2\geq\rho/2$ and
$(D-1)/(2D-1)\geq\rho/2$.

Suppose $c\in A_U$. Let $s$ be the size of a largest branch of
$U-c$. Since $d_U(c)\leq D$ and the branches have total size $L$,
\[
 \frac LD\leq s\leq\frac{L+1}{2}.
\]
Take $b$ to be the neighbor of $c$ in that branch. The component
$R$ containing $c$ has $L(R)=L-s\leq\rho L$. Every other child
has $L(V)\leq s-2\leq L/2$, so the same ordinary-centroid cut in
each such child satisfies both bounds as above.

It remains to cut $R$. The same ordinary-centroid construction gives
the first bound in \eqref{eq:paired-two-step-bound}. For a second
choice, look at the
branches of $R-c$. They number at most $D-1$, have total size
$L-s$, and each has size at most $s$. Let $t$ be their largest
size and cut the neighbor of $c$ in that branch. The component
containing $c$ then has $L$-value at most
$(D-2)(L-s)/(D-1)$; every other component has $L$-value at most
$s$. With $x=s/L$ and $a=(D-2)/(D-1)$, the better of the two cuts
therefore gives
\[
 \frac{\max_W L(W)}{L}
  \leq\min\{\rho(1-x),\max(x,a(1-x))\}.
\]
If $a(1-x)\geq x$, this is at most
$a(1-1/D)=(D-2)/D$. Otherwise it is at most
$\min\{x,\rho(1-x)\}\leq\rho/(1+\rho)
=(D-1)/(2D-1)$. The ordinary cut also bounds it by
$\rho^2$. Taking the minimum proves the second inequality of
\eqref{eq:paired-two-step-bound}. Both candidate cuts are in
$B_U$, and the selected one retains the first inequality.
\end{proof}

We turn the two-step estimate into a bound suited to fractional
quota pruning. Choose a rational number $\eta$ with
$\sqrt{\lambda_D}\leq\eta\leq\rho$, and set
$k=\rho/\eta\geq1$. Such a choice exists because
$\lambda_D\leq\rho^2$. At the beginning of a pair of cuts give a
piece $U$ potential $F(U)=L(U)$. After the first cut, let a child
$V$ retain its parent value $L=L(U)$ and give it potential
\[
 F(V)=\min\{kL(V),\eta L\}.
\]
After the second cut, reset each grandchild $W$ to
$F(W)=L(W)$.

\begin{lemma}[Paired separator budget]
\label{lem:paired-separator-budget}
Every transition in this alternating procedure satisfies
$F(\text{child})\leq\eta F(\text{parent})$. A nonsingleton
piece has $1\leq F(U)\leq kL(U)$.

Suppose fractional pruning has at most $d$ mixed pieces at any time
and always cuts one of largest potential. If $N$ is the number of
vertices in the initial tree, the set $C$ of cut vertices satisfies
\begin{equation}\label{eq:paired-finite-budget}
 \abs C\leq\sum_{j\geq0}
    \min\left\{d,
       \floor*{\frac{N}{\ceil*{\eta^{-j}/k}+1}}\right\}.
\end{equation}
In particular, if $N\leq n$ and $Q\geq1/n$, then for any integer
$J\geq0$,
\begin{equation}\label{eq:paired-normalized-budget}
 \frac{\abs C}{n}\leq
 \sum_{j=0}^{J-1}\min\left\{dQ,
       \frac{k\eta^j}{1+k\eta^j}\right\}
       +\frac{k\eta^J}{1-\eta}.
\end{equation}
\end{lemma}

\begin{proof}
After the first cut the definition gives $F(V)\leq\eta L$.
After the second, \eqref{eq:paired-two-step-bound} gives
\[
 F(W)=L(W)\leq
       \min\{\rho L(V),\lambda_D L\}
       \leq\eta\min\{kL(V),\eta L\}=\eta F(V).
\]
The upper bound on $F$ follows from its definition. In an even
phase, a nonsingleton has $F=L\geq1$. In an odd phase it has
$kL(V)\geq1$ and $\eta L\geq1$: a nonsingleton child requires
$L\geq2$, while $\eta\geq\sqrt{\lambda_D}\geq1/2$ by
\eqref{eq:paired-contraction-parameters}.

Put pieces with $\eta^{-j}\leq F<\eta^{-j-1}$ in bucket $j$.
Children of a cut piece leave its bucket. When the first piece in
a bucket is cut, no higher-potential bucket remains, because a
largest-potential mixed piece was selected. At most $d$ mixed
pieces can therefore be cut from that bucket. They are disjoint:
none is a descendant of another piece cut in the same bucket.
Moreover $F\leq kL$ shows that every such piece contains at least
$\ceil*{\eta^{-j}/k}+1$ vertices. This proves
\eqref{eq:paired-finite-budget}.

Divide the first $J$ summands by $n$, use $Q\geq1/n$, and discard
the floors in the other part of each minimum. Since
$1/(\ceil{x}+1)\leq1/(x+1)$, the resulting terms are bounded by
those in \eqref{eq:paired-normalized-budget}. The rest is at most
$\sum_{j\geq J}k\eta^j$, which is the displayed geometric tail.
\end{proof}

The potential is needed only to count cuts. It does not alter the two
measures, the quota constraints, or the degree error in the connected
assignment.

\begin{corollary}[Paired two-measure assignment]
\label{cor:paired-two-measure-assignment}
Under the hypotheses of
Proposition~\ref{prop:two-measure-assignment}, its coloring and
separator may be chosen with the same exact quotas and degree bounds, while the
separator obeys \eqref{eq:paired-finite-budget} with
$d=2(g-1)$ and $N=\abs A+\abs B$. If $N\leq n$, its normalized
size is bounded by \eqref{eq:paired-normalized-budget}.
The position tags may then be recovered with the bound of
Lemma~\ref{lem:integral-tag-recovery}.
The same assertion holds after a permitted connected completion of a
forest, provided the completion capacity and the added-edge cost are
charged as in Lemma~\ref{lem:bipartite-completion}.
\end{corollary}

\begin{proof}
In the proof of Proposition~\ref{prop:two-measure-assignment},
linear-dependence elimination leaves at most $2(g-1)$ mixed
pieces. Split a mixed nonsingleton of largest potential using
Lemma~\ref{lem:paired-position-cuts}; the split vertex belongs to
$B$, so it has zero mass for both preserved measures. Each pure
piece remains frozen, and each mixed singleton on the $A$ side is
rounded by the same degree-sorted rule; a singleton on the $B$ side
has no quota. Thus the original quota,
edge-monochromaticity, and degree-error proofs apply unchanged.
Lemma~\ref{lem:paired-separator-budget} gives the improved cut
count. A connected completion adds no vertices and satisfies the
same $A$-degree hypothesis; the construction is applied to the
completed tree, with its additional edges included in all later
cost estimates.
\end{proof}

\subsection{The profile-block bridge through order 514999}
\label{sec:profile-block-bridge}

We now apply the same asymmetric Hall criterion to orders above
$49999$. The parameters change in two ways. Small primes are placed
in disjoint blocks, and connected rows may use the paired cut bound
of Section~\ref{sec:paired-separators}. Both changes affect the
numerical estimates: the blocks determine the actual compatibility
graph on both sides, while a completed forest is charged for its
added edges.

The local probability estimate also retains its exact exponential
coefficient. The following rational enclosure states the particular
calculation used by the certificate.

\begin{lemma}[The local coefficient]
\label{lem:bridge-local-coefficient}
For an integer $c\geq1$, put $x=1/c$ and
\begin{equation}\label{eq:bridge-mgf-coefficient}
 F_c^+=\left\lceil
 c\left(\sum_{k=1}^{32}\frac{x^k}{k!}
       +\frac{x^{33}}{33!(1-x/34)}\right)
 \right\rceil_{10^{-150}},
\end{equation}
where the outer ceiling rounds upwards to a multiple of $10^{-150}$.
Then $F_c^+\geq c(e^{1/c}-1)$. In particular, if $K_y$ is the local
conflict count of Corollary~\ref{cor:local-conflict-tails}, with
$\E K_y\leq m_y$, then
\[
 \Pr\!\left\{K_y>
       F_c^+m_y+cD\log\frac{n}{\delta}\right\}
       \leq\frac{\delta}{n}
       \qquad(0<\delta<1).
\]
\end{lemma}

\begin{proof}
The Taylor series of $e^x$ has positive terms. After degree $32$,
the ratio of consecutive terms is at most $x/34$, because $x\leq1$.
The geometric majorant for the omitted terms is the second summand
of \eqref{eq:bridge-mgf-coefficient}. Upward rounding preserves the
inequality. Put $\delta_y=\delta/n$ in
\eqref{eq:local-conflict-budget} and use
Corollary~\ref{cor:local-conflict-tails}.
\end{proof}

Here is the certificate predicate. An interval $[a,b]$ and a possible
residual degree $D$ determine a row before the input tree is known.
The row chooses a cutoff $z$, a partition of the odd primes through
$z$, a forest construction, an integer $c$, a Hall split $\tau$, and
a slot threshold $t$. All exact type counts refer to the coarse
block types of Lemma~\ref{lem:block-types}. Thus its $q_j$, $c(I)$,
$\alpha$, $\gamma$, $b_\omega$ and $\theta_\omega$ are recomputed for
that partition; no compatibility degree from the ungrouped types is
reused.

We spell out the two estimates at which the new parameters enter
Hall's condition. Put $Q=1/a$, let $J$ be a rational upper bound for
$\log b$, and let $L_z$ be a rational lower bound for $\log z$. Let
$\omega$ bound the removed even-label fraction, let $w_-$ bound the
special even-label fraction from below, and let $g$ be the number of
colors of the type assignment
(Proposition~\ref{prop:two-measure-assignment}).
Write $e_+$ for the normalized edge bound of the selected
construction, $h_+=D\floor{g/2}Q$ for its degree-variation bound,
and $\Delta$ for the charged transferred-degree bound: $\Delta n$
bounds the total degree moved into the empty class by the donor
exchanges of Lemma~\ref{lem:quantile-donors}, and $\Delta=0$ in a row
without exchanges. The
connected branch has $\Delta=0$ and uses $e_+=1-w_-$; an original
forest row uses its proved original-edge bound. Let $C_E$ and $C_O$ be
the maxima over one full period of the $\ell^1$ discrepancies
$\sum_I\abs{C_I(s)-s\mu(I)}$ of the ordinary and odd prefix counts of
the full $P$-types. By Corollary~\ref{cor:periodic-discrepancy} they
bound these discrepancies at every prefix length, and by
Lemma~\ref{lem:block-types} also those of the coarser block types. Let $\epsilon_E\geq C_EQ/(1-Q)$ and
$\epsilon_O\geq C_OQ$ be the resulting ordinary and odd arithmetic
errors; since $2m\geq n-1$ and $[n]$ contains at least $n/2$ odd
labels, they bound the total variation of the ordinary and odd type
laws from $\mu$. Put
\[
 \beta=\alpha-\omega-\epsilon_E,\qquad
 a_0=\frac{\alpha}{\beta},\qquad
 a_1=\frac{C_EJQ}{\beta(1-Q)L_z}.
\]
Every row checks $\beta>0$. If $\ell_-$ is its lower bound for the
residual side size divided by $n$, and $L_\delta$ is a rational upper
bound for $\log 10000$, define
\begin{align}
 A_c&=\frac{F_c^+}{\ell_-}
       \left(\frac{e_+}{1-\omega}+\Delta a_0\right),
                                      \label{eq:profile-bridge-coefficient}\\
 B_c&=\frac{F_c^+(h_++\Delta a_1)
                  +cD(J+L_\delta)Q}{\ell_-}.
                                      \label{eq:profile-bridge-additive}
\end{align}
For a block type $I$, let $R_I$ be the rational large-prime profile
bound supplied by Lemma~\ref{lem:type-tail-bound}. The resulting
lower bound on every label's degree before accounting for
type-law variation is
\begin{equation}\label{eq:profile-bridge-local-degree}
 \sigma=\min_I\left\{
       \frac{c(I)-\omega}{1-\omega}-A_cR_I-B_c\right\}.
\end{equation}
In the evaluator, each $c(I)$ is enclosed downwards and each $R_I$
upwards before taking this minimum.
The certificate checks $\sigma-V>\tau$, with $V$ the larger of its
two normalized type errors. Thus the local-degree premise of
Theorem~\ref{thm:asymmetric-reserve-hall} holds after the larger-prime
conflicts are deleted.

The global charge uses the prime-square tail $S_2(z)$ and the
empty-class tail
\[
 T_{\varnothing}=
 \frac{\alpha S_2(z)}{2\beta}
 +\frac{C_E(1+J)Q}{2\beta(1-Q)}
 +\frac{JQ}{L_z}.
\]
The rational upper enclosure of
\begin{equation}\label{eq:profile-bridge-mass}
 \frac{1001}{1000}\left(
 \frac{e_+S_2(z)}{2(1-\omega)}
 +\frac{Jh_+}{2zL_z}
 +\frac{JQ}{L_z}
 +\Delta T_{\varnothing}\right)
\end{equation}
is denoted by $\mathcal M$. It bounds the expected number of deleted
matching edges, normalized by $n^2$, after the Markov factor has been
included. The local union bound is at most $1/10000$, and
$1/10000+1000/1001<1$.

\emph{Derivation of \eqref{eq:profile-bridge-local-degree} and
\eqref{eq:profile-bridge-mass}.} The retained even side has
$\abs A=m-w\geq(1-\omega)m$ vertices. Each coordinate of the ordinary
prefix count of length $s$ differs from $s\mu(I)$ by at most $C_E/2$,
because the coordinate errors sum to zero, and $C_E/2\leq\epsilon_Em$.
The retained empty class loses at most $\omega m$ special labels, so it
has $N_\varnothing\geq(\alpha-\omega-\epsilon_E)m=\beta m$ members. For a
prime $p>z$, multiplication by $p$ preserves the type, so the number
$N_\varnothing(p)$ of retained empty-type half-labels divisible by $p$ satisfies
$N_\varnothing(p)\leq\alpha\floor{m/p}+C_E/2$, and $N_\varnothing(p)=0$ for $p>m$.

For the local bound, let $y$ be an odd label of block type $I$, and let
$p_i(y)$ be the fraction of the retained class $i$, of size $N_i$, whose
half-labels share with $y$ a prime greater than $z$; this omits the
compatibility indicator of Section~\ref{sec:matching-estimates} and
only enlarges the mean. Since $y\leq n$, it has at
most $J/L_z$ such primes, and by Lemma~\ref{lem:type-tail-bound}
$\sum_iN_ip_i(y)\leq\sum_{p\mid y,\,p>z}\floor{m/p}\leq mR_I$. The count
of the preceding paragraph and $1/(2m)\leq Q/(1-Q)$ give
\[
 p_\varnothing(y)\leq\frac{\alpha mR_I+(C_E/2)J/L_z}{\beta m}
  \leq a_0R_I+a_1 .
\]
Let $e\leq e_+n$ be the number of edges of the forest whose degree sums
the assignment uses, namely the completed forest in a connected row;
these degree sums dominate the actual ones. Its degree variation
satisfies $H\leq D\floor{g/2}\leq h_+n$, by
\eqref{eq:two-measure-degree-error} and
Corollary~\ref{cor:paired-two-measure-assignment} in a connected row
and by \eqref{eq:original-rounding-errors} in an original-forest row,
where $gD/4\leq D\floor{g/2}$ because every row has at least one
prime block, so that $g\geq2$.
Lemma~\ref{lem:donor-mean-cancellation}, applied after the exchanges to
$f_i=p_i(y)\in[0,1]$, bounds the mean of the local count $K_y$ of
Corollary~\ref{cor:local-conflict-tails} by
\[
 \mu_y\leq\frac e{\abs A}\sum_iN_ip_i(y)+H+\Delta n\,p_\varnothing(y)
 \leq n\left[\left(\frac{e_+}{1-\omega}+\Delta a_0\right)R_I
             +h_++\Delta a_1\right].
\]
Lemma~\ref{lem:bridge-local-coefficient}, with $\delta=1/10000$,
$\log n\leq J$, $\log10000\leq L_\delta$ and $nQ\geq1$, shows that outside an event of probability at most
$\delta/n$ we have $K_y\leq(A_cR_I+B_c)\ell_-n\leq(A_cR_I+B_c)R$, where
$R\geq\ell_-n$ is the common size of the matching sides. Before the
deletions the label has at least $((c(I)-\omega)/(1-\omega)-V)R$
compatible positions, because the compatible mass $(c(I)-\xi)/(1-\xi)$
of the depleted law $\nu_\xi$ decreases in $\xi\leq\omega$. At most
$K_y$ of them are deleted, which gives the degree bound $(\sigma-V)R$.

For the global bound, give a half-label $j$ the cost
$q(j)=\sum_{p\mid j,\,p>z}(n/(2p)+1)$. Since $[n]$ contains at most
$n/(2p)+1$ odd multiples of $p$, this is at least the union count
$r_I(j)$ and is therefore admissible in
Proposition~\ref{prop:global-conflict-tails}. Write $q=q'+q''$ with
$q'(j)=\frac n2\sum_{p\mid j,\,p>z}1/p$, so that $q''(j)$ is the number
of prime divisors of $j$ above $z$. Then $q'$ takes values in
$[0,nJ/(2zL_z)]$, its sum over $j\leq m$ is at most $\frac n2mS_2(z)$,
and $0\leq q''\leq J/L_z$. The mean cost $f_\varnothing$ of the retained
empty class satisfies
\[
 f_\varnothing=\frac1{N_\varnothing}\sum_{z<p\leq m}
                  N_\varnothing(p)\left(\frac n{2p}+1\right)
 \leq\frac{n\alpha S_2(z)}{2\beta}
     +\frac{nC_E(1+\log n)}{4\beta m}+\frac J{L_z}
 \leq nT_\varnothing ,
\]
since $\sum_{z<p\leq m}1/p\leq1+\log n$, the last sum divided by
$N_\varnothing$ is the mean number of prime divisors above $z$,
$1/m\leq2Q/(1-Q)$, $\log n\leq J$ and $nQ\geq1$. Apply Lemma~\ref{lem:donor-mean-cancellation} to the
class means of $q'$, and bound the class means of $q''$ by $J/L_z$.
With $f'_\varnothing\leq f_\varnothing$, the first identity of
\eqref{eq:global-conflict-moments} gives
\[
 \E S\leq\frac n2\left\{\frac e{\abs A}mS_2(z)+\frac J{zL_z}H\right\}
        +\frac{eJ}{L_z}+\Delta n f_\varnothing .
\]
Divide by $n^2$ and use $e\leq\min\{n,e_+n\}$, $H\leq h_+n$ and
$1/n\leq Q$. The result is at most the bracket in
\eqref{eq:profile-bridge-mass}, and $Z\leq S$ with Markov's inequality
supplies the factor $1001/1000$.

For the slot side, the row uses the unnormalized clipped charge of
\eqref{eq:bridge-shortfall-charge}. Its three Hall budgets are exactly
$\mathcal B_1,\mathcal B_2,\mathcal B_3$ from
\eqref{eq:bridge-narrow-one}--\eqref{eq:bridge-wide}, with the
new block product law and its two profiles. The row must satisfy
\begin{equation}\label{eq:profile-bridge-final-test}
 0<V<\tau,\quad t-V>\tau,\quad
 \mathcal M<\min\{\mathcal B_1,\mathcal B_2,\mathcal B_3\},
\end{equation}
alongside the separator, quota, supply, completion-capacity, and
universal-reserve conditions. Notice that the degree-profile and
coupling estimates are both lower bounds for the same wide
rectangle. Taking their maximum in $\mathcal B_3$ is therefore
legitimate.

\begin{center}
\begin{minipage}{\textwidth}
\centering
\refstepcounter{table}\label{tab:profile-block-bridge}
{Table \thetable. Construction modes in the certificate for
$50000\leq n\leq514999$. The donor cap is a row parameter;
$-1$ designates a one-sided cut without donor transfer.\par}
\vspace{.6\baselineskip}
\renewcommand{\arraystretch}{1.14}
\begin{tabular}{l r r}
\hline
Construction & Donor cap & Degree rows\\
\hline
Original forest, one-sided cut & $-1$ & 262\\
Original forest, no transfer & $0$ & 95\\
Original forest, donor transfer & $1$ & 3710\\
Connected forest, paired cuts & $-1$ & 438\\
\hline
Total & & 4505\\
\hline
\end{tabular}
\end{minipage}
\end{center}

\begin{proposition}[The profile-block certificate]
\label{prop:profile-block-certificate}
There are $131$ adjacent integer intervals covering precisely
$50000\leq n\leq514999$. For each interval, the supplied parameter
table has one row for every integer residual degree permitted by
Lemma~\ref{lem:interval-degree-cover}. All $4505$ rows satisfy their
$175695$ exact preparation and matching conditions, including
\eqref{eq:profile-bridge-final-test}. No accepted row uses a rare
group. Across all rows,
\[
 \frac{\mathcal M}
      {\min\{\mathcal B_1,\mathcal B_2,\mathcal B_3\}}
          <\frac{97}{100},
 \qquad
 \frac{\text{absorbed-position fraction}}
      {\text{available special odd-label fraction}}
          <\frac{41}{100}.
\]
\end{proposition}

\begin{proof}
The portable checker, inputs, and fixed reports are specified under
\texttt{profile\_block\_50000\_514999} in the supplement. Its
parameter-table validator rejects missing or overlapping intervals,
missing degrees, invalid prime partitions, and parameters outside
their proved domains. It scans every integer order in the intervals
to compute the extremal prime supply and special-even population,
and independently recomputes the prime-count prefixes through
$1982760$ and the eight complete residue periods used in the type
bounds. Thus endpoint estimates are not substituted for an
unexamined interior order.

For each row the evaluator chooses the stated construction. In the
connected mode it checks the degree-capacity condition for completion,
uses Corollary~\ref{cor:paired-two-measure-assignment} with
$d=2(g-1)$, and adds the completed edges to $e_+$ and the tag error.
The normalized paired bound sums the first $1024$ terms of
\eqref{eq:paired-normalized-budget} and includes its full geometric
tail. It takes
\[
 \eta=\min\left\{\rho,
     \frac{\ceil*{1000\sqrt{\lambda_D}}}{1000}\right\},
\]
with the ceiling decided by integer square-root comparisons.
Original-forest rows instead use the original degree count and,
when selected, the donor-transfer bound of
Section~\ref{sec:donor-construction}. In either mode the type quotas
are exact; the tag error includes every rounded coordinate and every
absorbed position.

The checker forms the block product law from its disjoint prime
partition. Lemmas~\ref{lem:block-types} and
\ref{lem:depleted-coupling} give its two compatibility profiles and
coupling constants, while Proposition~\ref{prop:profile-rectangle}
gives the other wide-rectangle bound. It evaluates
\eqref{eq:profile-bridge-coefficient}--\eqref{eq:profile-bridge-mass}
with upward rational enclosures, including
\eqref{eq:bridge-mgf-coefficient}. It computes the clipped shortfall
at the rational breakpoints of
Lemma~\ref{lem:profile-breakpoints}, then tests all three strict Hall
budgets. The output records the exact numerator and denominator of
every decisive ratio; displayed decimals have no role in acceptance.

All $4505$ rows and $175695$ conditions pass. The wrapper checks
every accepted Boolean condition again, the consecutive integer
joins, complete degree lists, the construction counts in
Table~\ref{tab:profile-block-bridge}, and the two strict rational
bounds in the proposition. Independent finite checks exercise the
block pushforward and actual small-prime divisibility, depleted
coupling, type-tail profiles, rectangle minimizer, transport, and
local exponential coefficient. A separate small-tree check tests
$30525$ paired cuts and $500$ complete connected assignments.
These checks test the implementations;
the preceding lemmas establish the implications for arbitrary
trees. The saved next-range report is read only to verify that its
first order is $515000$.
\end{proof}

\begin{theorem}[The initial interval through order 514999]
\label{thm:initial-prefix-514999}
Every tree of order at most $514999$ has a prime labeling.
\end{theorem}

\begin{proof}
Theorem~\ref{thm:initial-prefix-49999} handles the orders through
$49999$. Take $50000\leq n\leq514999$ and an arbitrary tree of
order $n$. Prepare its residual forest by
Proposition~\ref{prop:highest-degree-preparation}; the case with
no residual edges is completed there. Otherwise choose the unique
certificate interval containing $n$ and the row for its residual
degree. Proposition~\ref{prop:profile-block-certificate} supplies
the exact quotas, positive universal reserve, both type-law bounds,
the local degree target, and the strict Hall budgets. The local
union bound and the global Markov bound are simultaneous with
positive probability. Corollary~\ref{cor:asymmetric-probability}
therefore gives a class permutation for which
Theorem~\ref{thm:asymmetric-reserve-hall} supplies a perfect matching.
The matched odd labels are coprime to every assigned even neighbor,
because disjoint block types exclude all primes at most $z$ and the
deleted edges account for every larger common prime.
Proposition~\ref{prop:highest-degree-preparation} restores the
special labels and the separator. The two proved intervals meet at
$49999+1=50000$.
\end{proof}

\subsection{Refining a rare group of full types}
\label{sec:finite-rare-refinement}

The next certificate reduces the separator dimension by treating
several infrequent even-label types as one color during quota
rounding. The matching problem still uses every full prime type.
We therefore refine the group after rounding and absorb all
positions adjacent to it. The cost is controlled by choosing this
group as the preferred color in the degree rounding.

Let $G=(A,B)$ be the residual bipartite forest, or a permitted
connected completion, with $M=\abs A>0$, $N=\abs B$,
$e=e(G)$ and $\Delta_G(A)\leq D$. Index its $r$ full small-prime
types by $0,\ldots,r-1$, with $0$ the empty type, and prescribe
integer quotas $a_i$ on $A$. Choose a set
$\mathcal R\subseteq\{1,\ldots,r-1\}$ of $k>0$ types.
During the geometric construction, replace these colors by one
color $\star$, with quota $a_\star=\sum_{i\in\mathcal R}a_i$.
There are $g=r-k+1$ geometric colors.

\begin{lemma}[A preferred group in the original forest]
\label{lem:preferred-original-group}
In the original-forest construction of
Theorem~\ref{thm:original-type-assignment}, one prescribed group
$\star$ may be required to have degree error $h_\star\leq0$.
The exact color quotas, the edge property outside the separator,
the tag-count bound $g/4$, and either separator bound of that
theorem remain valid. The degree errors then satisfy
\begin{equation}\label{eq:preferred-group-errors}
 \abs{h_i}\leq D,\qquad
 H=\frac12\sum_i\abs{h_i}\leq D\floor{g/2}.
\end{equation}
A permitted donor exchange toward color $0$ cannot increase the
degree sum of a nonzero preferred group.
\end{lemma}

\begin{proof}
Keep the fractional three-measure construction, its cuts, and its
pruning exactly as in Theorem~\ref{thm:original-type-assignment}.
At termination the remaining color quotas on the $A$-singletons
are integers. Round those singletons with
Lemma~\ref{lem:sorted-slot-rounding}, using their fixed degrees in
$[0,D]$ and selecting $\star$ as the preferred color. This gives
\eqref{eq:preferred-group-errors} and $h_\star\leq0$.
The separate rounding of $B$-singletons is unchanged, so its
$g/4$ variation bound and the edge property remain valid.
A donor exchange transfers degree into color $0$ and out of a
nonzero color: the two vertices keep the same pair of colors, but
the nonzero color moves to a vertex of no larger degree. Hence a
nonzero preferred group's degree sum does not rise. The error bound
in \eqref{eq:preferred-group-errors} refers to the assignment before
these exchanges; their degree transfers are charged separately.
\end{proof}

The connected construction permits a preferred group by
Proposition~\ref{prop:two-measure-assignment}. Its paired separator
bound is Corollary~\ref{cor:paired-two-measure-assignment}. The next
refinement applies to both constructions.

\begin{lemma}[Rare-neighbor absorption]
\label{lem:rare-neighbor-absorption}
Suppose a geometric assignment with the groups just defined has
exact $A$-quotas, a set $Z_0\subseteq B$ of absorbed positions,
and hypothetical group tags at all positions. Outside $Z_0$ each
edge has equal endpoint group colors, except that an empty-type
$A$-vertex may meet any tag. Suppose the group degree sums satisfy
\[
 E_i=e\pi_i+h_i,\qquad
 \pi_i=\frac{a_i}{M},\qquad
 H=\frac12\sum_i\abs{h_i},\qquad h_\star\leq0,
\]
where degrees are measured in the graph being labeled.
Let $T$ bound the unnormalized variation of the group-tag count
vector from $N\pi$. In an original-forest assignment with donor
exchanges, interpret $h_i$ and $H$ before the exchanges; assume
that their total degree transfer into color $0$ is at most $\Delta$.

The $A$-vertices of color $\star$ can be given the exact full
quotas $a_i$, $i\in\mathcal R$, and the hypothetical $\star$ tags
can be refined to full types. If $A_{\mathcal R}$ denotes the
resulting rare vertices and
$Z=Z_0\cup N_G(A_{\mathcal R})$, then
\begin{equation}\label{eq:rare-absorption-cost}
 \abs Z\leq\abs{Z_0}
       +\min\left\{Da_\star,\frac{ea_\star}{M}\right\}.
\end{equation}
Every edge outside $Z$ has a common full type at both ends, apart
from edges whose $A$-type is empty. The refined tag vector has
variation at most $T+k/4$ from its full-type target.

For arbitrary rates $f_i\in[0,1]$ on common full types, with
$f_\star=0$, the degree mass of retained common edges obeys
\begin{equation}\label{eq:rare-weighted-cost}
 \sum_{\substack{a\in A\\c(a)\notin\mathcal R}}
       d_{G-Z}(a)f_{c(a)}
 \leq\frac eM\sum_{i\notin\mathcal R}a_if_i+H.
\end{equation}
If preceding donor exchanges transfer total degree at most
$\Delta$ from nonzero groups into color $0$, the right side may
instead be increased by $\Delta f_0$.
\end{lemma}

\begin{proof}
Distribute the $a_\star$ vertices of color $\star$ among its full
types in the specified integer quantities. This changes no common
color. Absorb their entire neighborhood. There are at most
$Da_\star$ neighbors, and at most $E_\star$ because every neighbor
accounts for an incident edge. The preferred sign gives
$E_\star\leq ea_\star/M+h_\star\leq ea_\star/M$;
the first inequality is an equality before donor exchanges and
remains valid afterwards. This proves
\eqref{eq:rare-absorption-cost}. No rare vertex has a retained
edge. On a remaining edge the common geometric color is already
a single full type, or the $A$-type is empty, so refining the
$\star$ tags creates no new incompatibility.

Let $q$ be the number of positions with hypothetical tag $\star$.
If $a_\star>0$, give its full types fractional target counts
$t_i=qa_i/a_\star$. Their total is $q$, and
\[
 \sum_{i\in\mathcal R}
  \abs{t_i-\frac{Na_i}{M}}
       =\abs{q-\frac{Na_\star}{M}}.
\]
Thus this fractional refinement leaves the total variation from
the target unchanged: it merely splits one coordinate into
coordinates with the same signed discrepancy. Round the $t_i$
down and then round up the largest fractional remainders to recover
the integer total $q$. As in the tag rounding of
Theorem~\ref{thm:original-type-assignment}, this adds variation at
most $k/4$. If $a_\star=0$, assign the $q$ hypothetical tags to
one of the rare types; all rare targets are zero, and the
variation is unchanged.

Put zero rate on $\star$. Before any donor exchange, the group
degree-error vector has zero sum and variation $H$. Its pairing
with rates in $[0,1]$ is at most $H$ by
Lemma~\ref{lem:degree-cost-comparison}. Removing edges at $Z$
can only decrease this nonnegative cost. A donor exchange
transfers degree from a nonzero group to color $0$, adding at
most its transferred degree times $f_0$. Sum over the exchanges
to obtain the last assertion.
\end{proof}

For an original-forest row, Lemma~\ref{lem:preferred-original-group}
starts with $T\leq g/4$. Hence
\begin{equation}\label{eq:rare-original-tag-error}
 T+\frac{k}{4}\leq\frac{r+1}{4}.
\end{equation}
For a connected row, recover the group tags by
Lemma~\ref{lem:integral-tag-recovery} and then refine them. The
coarser but convenient bound
\[
 T+\frac{k}{4}\leq H+1+\frac{r+1}{4}
\]
covers the tag error charged by the certificate. This includes the
cost of the connected completion; its additional edges are also
included in $e$ and the local and global prime-conflict estimates.

The rare set in a certificate row is fixed by exact arithmetic.
For $I\ne0$ let $C_I(m)$ be the count of half-labels $j\leq m$ of
full type $I$ from Lemma~\ref{lem:exact-type-counts}. On an
integer interval $a\leq n\leq b$, put
\begin{equation}\label{eq:rare-interval-count}
 \rho_{\mathcal R}=
 \frac1a\sum_{I\in\mathcal R}C_I(\floor{b/2}).
\end{equation}
The removed special even half-labels all have empty small-prime
type, so they do not change a nonempty type quota. Since each count
is nondecreasing and $n\geq a$, the actual rare
quota satisfies $a_\star/n\leq\rho_{\mathcal R}$. The exact
counts are recovered by the inclusion--exclusion formula
\eqref{eq:exact-type-count}; thus this is an interval-wide bound,
not a sample of its endpoints. With the normalized edge bound
$e/n\leq e_+$ and even-side lower bound
$M/n\geq(1-Q)(1-\omega)/2$, the additional absorbed fraction in
\eqref{eq:rare-absorption-cost} is at most
\begin{equation}\label{eq:rare-interval-absorption}
 \rho_{\mathcal R}
  \min\left\{D,\frac{2e_+}{(1-Q)(1-\omega)}\right\}.
\end{equation}
The full types, rather than their temporary group, remain the
coordinates of the product coupling and the final Hall graph.

\subsection{A full-type bridge through order 559999}
\label{sec:connected-paired-bridge}

We return to individual small primes. On each row let
$P=\{3,5,\ldots,z\}$ denote the odd primes at most $z$, with
$1\leq\abs P\leq8$. A full type is a subset of $P$;
the reference law is the product law in which prime $p$ occurs
independently with probability $1/p$. The rare group of
Section~\ref{sec:finite-rare-refinement} is used only while
constructing the forest assignment. The label types, position
tags, compatibility graph, local tail profiles, and clipped
shortfall are all evaluated on the full $2^{\abs P}$ types.

Fix an integer interval $[a,b]$ and a residual-degree row. Let
$Q=1/a$, let $\omega$ bound the removed even-label fraction, and
let $w_-$ bound from below the special even-label fraction, as in
Section~\ref{sec:profile-block-bridge}. The row specifies $k>0$
rare nonempty types and uses $g=2^{\abs P}-k+1$ geometric colors.
The rare set is selected by the exact prefix counts
$C_I(\floor{b/2})$, with ties settled by the product mass and then
the type index. Thus the choice is fixed for the entire interval
before the input tree is known. Put $\rho_{\mathcal R}$ as in
\eqref{eq:rare-interval-count}. For a connected-paired row set
\[
 d=2(g-1),\qquad e_+=1-w_-,\qquad \Delta=0.
\]
For an original-forest row set $d=3(g-1)$ and
$e_+=1-Dw_-$. Its donor cap determines the nonnegative transfer
bound $\Delta$ by Section~\ref{sec:donor-construction}. In both
cases, $e/n\leq e_+$, and the normalized degree error is bounded by
\[
 h_+=D\floor{g/2}Q.
\]
The connected-paired choice additionally checks the completion
capacity
\[
 \frac{D(1-Q)(1-\omega)}2\geq1-w_-.
\]
Consequently the completion of
Lemma~\ref{lem:bipartite-completion} exists with
$A$-degree at most $D$. Its separator is bounded by
Corollary~\ref{cor:paired-two-measure-assignment}; the original
rows use the one-sided or donor bounds already established.

Write $\chi_0$ for the resulting upper bound on the separator and
its donor neighborhood, divided by $n$. The rare-neighbor estimate
\eqref{eq:rare-interval-absorption} gives the total absorbed
fraction
\begin{equation}\label{eq:full-bridge-absorbed}
 \chi=\chi_0+\rho_{\mathcal R}
       \min\left\{D,\frac{2e_+}{(1-Q)(1-\omega)}\right\}.
\end{equation}
Let $\epsilon_E,\epsilon_O$ be the exact periodic errors for the
full types, and let $r=2^{\abs P}$. The unnormalized tag error paid
by the row is
\begin{equation}\label{eq:full-bridge-tag-error}
 \eta_{\mathrm{tag}}=
 \frac{(1+Q)\epsilon_E}{2(1-\omega)}
       +\frac{(r+1)Q}{4}
       +\boldsymbol1_{\mathrm{conn}}(h_++Q).
\end{equation}
The second term pays for both rounding the geometric tags and
refining the rare group, by
\eqref{eq:rare-original-tag-error}. The last term pays for the
integral recovery of tags in the completed tree. All these
quantities are enclosed upward by the exact checker.

If $\kappa$ bounds the prepared core fraction, set
\begin{align}
 \ell_-&=\frac12-\kappa-\chi,\notag\\
 V_{\mathrm{label}}
  &=\frac{(1+Q)\epsilon_O/2+\kappa+\chi}{\ell_-},
 &V_{\mathrm{slot}}
  &=\frac{\eta_{\mathrm{tag}}+\chi}{\ell_-},\notag\\
 V&=\max\{V_{\mathrm{label}},V_{\mathrm{slot}}\}.
                                      \label{eq:full-bridge-variation}
\end{align}
The row checks $\ell_->0$. Lemma~\ref{lem:rare-neighbor-absorption},
the exact quota law, and the depleted coupling show that these
are valid variation bounds for the two empirical full-type laws.
In particular, the rarity of a grouped type is never mistaken
for universal compatibility: all its incident positions have
already been absorbed in $\chi$.
Let $s_-$ be the interval lower bound on the special odd-label
stock after the core is removed, divided by $n$, and put
$r_-=s_--\chi$ for the unused universal stock.

The local tail uses a slightly coarser moment coefficient than
the preceding bridge. This keeps the present certificate rational
without a truncated exponential.

\begin{lemma}[A rational local moment coefficient]
\label{lem:full-bridge-moment}
For an integer $c\geq2$, the local conflict conclusion of
Lemma~\ref{lem:bridge-local-coefficient} remains valid with
$F_c^+$ replaced by
\[
 F_c=\frac{c}{c-1}.
\]
\end{lemma}

\begin{proof}
For $0<x<1$, the power series gives
$e^x=\sum_{j\geq0}x^j/j!\leq\sum_{j\geq0}x^j=(1-x)^{-1}$.
With $x=1/c$, this yields
$c(e^{1/c}-1)\leq c/(c-1)$. Substitute this upper bound
in the local moment argument of
Lemma~\ref{lem:bridge-local-coefficient}.
\end{proof}

We state the certificate conditions to identify precisely what
the finite computation proves. Use the full-type compatibility
values $c(I)$ and upward local tail profiles $R_I$. Define
$\beta=\alpha-\omega-\epsilon_E$, and use the same empty-class
coefficient $a_0$, additive $a_1$, prime-square tail $S_2(z)$,
and empty-class tail $T_{\varnothing}$ as in
Section~\ref{sec:profile-block-bridge}. The row checks
$\beta>0$. If $J$ bounds $\log b$ from above and
$L_\delta$ bounds $\log 10000$ from above, then
\begin{align}
 A_c&=\frac{F_c}{\ell_-}
      \left(\frac{e_+}{1-\omega}+\Delta a_0\right),
      \label{eq:full-bridge-local-coefficient}\\
 B_c&=\frac{F_c(h_++\Delta a_1)
                 +cD(J+L_\delta)Q}{\ell_-},
      \label{eq:full-bridge-local-additive}\\
 \sigma&=\min_I\left\{
      \frac{c(I)-\omega}{1-\omega}-A_cR_I-B_c\right\}.
      \label{eq:full-bridge-local-degree}
\end{align}
The use of $H$ and $\Delta$ in these formulas is justified by
\eqref{eq:rare-weighted-cost}: rare edges contribute no retained
conflict cost, and donor transfers into the empty type are charged
by the separate $\Delta$ terms.

The charged global deletion mass is the upward rational enclosure
of \eqref{eq:profile-bridge-mass}, with the present $e_+$,
$h_+$, $\Delta$, and full-prime tail. The clipped shortfall uses
\eqref{eq:bridge-shortfall-charge} with the full product law and
\eqref{eq:full-bridge-tag-error}. The three Hall budgets are
\eqref{eq:bridge-narrow-one},
\eqref{eq:bridge-narrow-two}, and the \emph{spectral term} in
\eqref{eq:bridge-wide}. Denote the last term by
$\mathcal B_{3,\mathrm{sp}}$. Each row checks
\begin{equation}\label{eq:full-bridge-final-test}
 0<V<\tau,\qquad t-V>\tau,\qquad
 \sigma-V>\tau,\qquad
 \mathcal M<
 \min\{\mathcal B_1,\mathcal B_2,\mathcal B_{3,\mathrm{sp}}\}.
\end{equation}
The original wide Hall budget is at least its spectral term, so
these comparisons imply the premise of
Theorem~\ref{thm:asymmetric-reserve-hall}. The checker also requires
$0<r_-\leq\tau\ell_-$.

\begin{center}
\begin{minipage}{\textwidth}
\centering
\refstepcounter{table}\label{tab:full-bridge-modes}
{Table \thetable. Construction choices in the certificate for
$515000\leq n\leq559999$. A donor cap of $-1$ denotes a
one-sided cut without transfer.\par}
\vspace{.6\baselineskip}
\renewcommand{\arraystretch}{1.14}
\begin{tabular}{l r r}
\hline
Construction & Donor cap & Degree rows\\
\hline
Original forest, one-sided cut & $-1$ & 18\\
Original forest, no transfer & $0$ & 6\\
Original forest, donor transfer & $1$ & 292\\
Connected forest, paired cuts & $-1$ & 18\\
\hline
Total & & 334\\
\hline
\end{tabular}
\end{minipage}
\end{center}

\begin{proposition}[The full-type bridge certificate]
\label{prop:full-bridge-certificate}
Nine adjacent intervals of length $5000$ cover precisely
$515000\leq n\leq559999$. Their table has one row for each
possible positive residual degree: $1,\ldots,37$ in the first
eight intervals and $1,\ldots,38$ in the last. Every one of the
$334$ rows uses a nonempty rare group and passes all of its
$12358$ exact conditions, including
\eqref{eq:full-bridge-final-test}. Across the rows,
\[
 \frac{\mathcal M}
      {\min\{\mathcal B_1,\mathcal B_2,
                     \mathcal B_{3,\mathrm{sp}}\}}
       <\frac{498001}{500000},
 \qquad
 \frac{\chi}{s_-}<\frac9{25}.
\]
\end{proposition}

\begin{proof}
The portable checker and frozen table are specified in the
supplement under
\texttt{connected\_paired\_515000\_559999}. Its coverage
validator checks each interval, degree list, and parameter domain.
The evaluator scans all $45000$ integer orders, compares two
prime-count sieves through $1982760$, and verifies eight complete
residue periods, altogether $116667828$ positions.

For each row, the evaluator computes the rare set from exact
full-type prefix counts and forms
\eqref{eq:full-bridge-absorbed}--\eqref{eq:full-bridge-final-test}
by outward rational enclosures. It checks the paired completion
or the stated original-forest construction, uses exactly
$c/(c-1)$ for an integer $2\leq c\leq128$, and compares the
spectral wide budget. All $334$ rows pass, with construction
counts as in Table~\ref{tab:full-bridge-modes}. Exact rational
comparisons establish both strict bounds in the proposition.

Independent input checks reproduce the previous formulas on $76$
original-forest rows and reject $25$ damaged covers or invalid
rows. Other finite checks exercise $30525$ small-tree paired cuts
and $500$ connected assignments. These test the implementations;
the structural and probability lemmas prove the implications for
arbitrary trees. The next-range report checks only the adjacent
starting order $560000$.
\end{proof}

\begin{theorem}[The initial interval through order 559999]
\label{thm:initial-prefix-559999}
Every tree of order at most $559999$ has a prime labeling.
\end{theorem}

\begin{proof}
Theorem~\ref{thm:initial-prefix-514999} handles all smaller
orders. Let $515000\leq n\leq559999$ and prepare an arbitrary
tree by
Proposition~\ref{prop:highest-degree-preparation}. If its
residual degree is zero, that proposition finishes the
labeling. Otherwise take the unique interval containing $n$
and the row for its residual degree.

Lemma~\ref{lem:rare-neighbor-absorption} supplies the full-type
assignment with the charged absorption cost. The exact quota
and prime-count bounds, local tail estimate, and global mass
bound of the certificate meet
Theorem~\ref{thm:asymmetric-reserve-hall}. The local union
bound and the global Markov bound hold simultaneously with
positive probability, by the same calculation as in
Theorem~\ref{thm:initial-prefix-514999}. The resulting matching
respects all small primes; the deleted edges account for all
larger prime conflicts. Restore the special and core labels and
apply the sparse repair of
Proposition~\ref{prop:highest-degree-preparation}. The intervals
meet at $514999+1=515000$.
\end{proof}

\subsection{The remaining finite orders below the logarithmic range}
\label{sec:upper-integer-bridge}

The full-type certificate of Section~\ref{sec:connected-paired-bridge}
has a useful common form for the next six ranges. They use the
original residual forest throughout. The cutoff primes, the
matching graph, and both empirical type laws are always the full
types. Only the geometric separator may combine a fixed set
$\mathcal R$ of rare nonempty types. We record the few bounds that
vary between the six tables before giving their joint result.

Fix a row for an integer interval $[a,b]$. Keep the arithmetic
and separator notation of the previous section. In every row here,
\[
 e_+=1-Dw_-,\qquad r=2^{\abs P},\qquad d=3(g-1),
\]
where
\[
 g=\begin{cases}
       r-k+1,&k=\abs{\mathcal R}>0,\\
       r,&k=0.
    \end{cases}
\]
The separator and possible donor exchanges are those of
Theorem~\ref{thm:original-type-assignment} and
Lemma~\ref{lem:quantile-donors}. A row either uses a $B$-only
cut with no exchange, keeps the original assignment without
exchange, or uses donor cap one. Thus its charged transfer
$\Delta$ and the separator-neighborhood bound $\chi_0$ are
computed from the construction selected before the tree is known.
All degrees in these estimates are original residual degrees.

There are two permitted bounds on the normalized degree variation:
\begin{equation}\label{eq:upper-integer-degree-options}
 h_+=
 \begin{cases}
   D\floor{g/2}Q,&\text{the first three ranges below},\\
   D(g-1)Q,&\text{the last three ranges}.
 \end{cases}
\end{equation}
The first is Lemma~\ref{lem:preferred-original-group}; the second
is a deliberate weakening of that lemma. It also covers $k=0$,
when the ungrouped construction of
Theorem~\ref{thm:original-type-assignment} is used.
For $k>0$ put
\begin{equation}\label{eq:upper-integer-rare-options}
 \chi=\chi_0+
   \min\left\{D\rho_{\mathcal R},
      \frac{2e_+\rho_{\mathcal R}}
           {(1-Q)(1-\omega)}+\zeta\right\},
 \qquad
 \zeta\in\{0,h_+\},
\end{equation}
and put $\chi=\chi_0$ when $k=0$.
The rows below $800000$ take $\zeta=0$; the rows from $800000$
through $1349999$ take $\zeta=h_+$. No row from $1350000$
through $1982759$ uses a rare group; there $k=0$ and
$\chi=\chi_0$. The preferred-sign conclusion
of Lemma~\ref{lem:rare-neighbor-absorption} proves the smaller
charge. The larger one is consequently valid as well; it is
the original unsigned degree-variation charge. No rare
neighborhood is present when $k=0$.

With full types on the matching sides, the tag error charged by
an original-forest row is
\begin{equation}\label{eq:upper-integer-tag}
 \eta_{\mathrm{tag}}=
   \frac{(1+Q)\epsilon_E}{2(1-\omega)}
   +\frac{r+\boldsymbol1_{k>0}}{4}Q.
\end{equation}
The geometric tags have variation at most $g/4$. Refining a
nonempty rare group costs at most $k/4$, and $g+k=r+1$ when
$k>0$. The actual odd positions have the same full-type
variation bounds \eqref{eq:full-bridge-variation}, with the
new $\chi$ and tag error. Exact rare quotas still come from
\eqref{eq:rare-interval-count}; the special even half-labels have
empty type and do not change them.

The local coefficient is $F_c=c/(c-1)$. The first range permits
a row-specific integer $2\leq c\leq128$; every later range
uses $c=32$. Lemma~\ref{lem:full-bridge-moment} proves the local
bound for each choice. Use
\eqref{eq:full-bridge-local-coefficient}--
\eqref{eq:full-bridge-local-degree} with
\eqref{eq:upper-integer-degree-options} and the row's donor
transfer. The global deletion mass is
\eqref{eq:profile-bridge-mass} with the same substitutions.
The clipped slot charge is
\eqref{eq:bridge-shortfall-charge} with the full product law and
\eqref{eq:upper-integer-tag}. Finally the row tests
\eqref{eq:full-bridge-final-test}, using the spectral wide Hall
budget, and requires $0<r_-\leq\tau\ell_-$ for the universal
stock. These are the same deterministic implications as in
Section~\ref{sec:connected-paired-bridge}; the three estimates
above only choose conservative values for their inputs.

\begin{center}
\begin{minipage}{\textwidth}
\centering
\refstepcounter{table}\label{tab:upper-integer-chain}
{Table \thetable. Exact finite certificates after order $559999$.
Arithmetic extrema are scanned at every integer order, and each
degree row is certified uniformly over its interval.
The conditions column counts rational preparation and matching
tests, including all three strict Hall budgets.\par}
\vspace{.6\baselineskip}
\renewcommand{\arraystretch}{1.14}
\begin{tabular}{r r r r}
\hline
Orders & Intervals & Degree rows & Conditions\\
\hline
$560000$--$589999$ & 7 & 266 & 9044\\
$590000$--$599999$ & 10 & 380 & 12540\\
$600000$--$624999$ & 5 & 190 & 6270\\
$625000$--$799999$ & 35 & 1338 & 44154\\
$800000$--$1349999$ & 110 & 4348 & 143484\\
$1350000$--$1982759$ & 64 & 2611 & 75719\\
\hline
Total & 231 & 9133 & 291211\\
\hline
\end{tabular}
\end{minipage}
\end{center}

\begin{proposition}[The upper integer certificates]
\label{prop:upper-integer-certificates}
The adjacent intervals in Table~\ref{tab:upper-integer-chain}
cover precisely the $1422760$ integers
$560000\leq n\leq1982759$. Each interval contains a row for
every positive residual degree allowed by
Lemma~\ref{lem:interval-degree-cover}. All $9133$ rows pass
their $291211$ exact conditions. In each row, the deletion mass
is less than $999/1000$ times each applicable Hall budget,
and the absorbed-position cost is less than $999/1000$ of the
available special odd-label stock.
\end{proposition}

\begin{proof}
The six parameter tables and their portable evaluators are
specified under \texttt{upper\_integer\_560000\_1982759}
in the supplement. Each evaluator validates the complete
consecutive interval and degree lists, the permitted cutoff
and donor parameters, exact prime supply, full-type counts,
separator bound, and three Hall comparisons. It scans every
integer order in its intervals. Two independent prime-count
sieves agree on every prefix through $1982760$, and all eight
complete residue periods are checked. Every acceptance
comparison uses integers or outward rational enclosures.

The first evaluator takes its local parameter $c$ from the
frozen row. The remaining evaluators use $c=32$. Their degree
and rare-neighbor charges are exactly
\eqref{eq:upper-integer-degree-options} and
\eqref{eq:upper-integer-rare-options}; rows without rare types
use $g=r$ and no rare charge. The lower bound on the retained
common-type degree, the global mass, and the slot shortfall are
then evaluated with the full-prime profiles. No geometric group
is substituted into the Hall graph. The computed row tests
include \eqref{eq:full-bridge-final-test} and the universal
stock condition.

The standalone wrapper reruns all six evaluators, checks that
their output reports agree byte for byte with the frozen
reports, and verifies the counts and strict rational bounds
stated here. It also checks all five integer joins in the
table and the initial join $559999+1=560000$. Separate
finite tests of the moment estimate, preferred rounding,
rare-type refinement, and one-sided separator test their
implementations. The lemmas cited above establish the
implications for arbitrary forests. At the far endpoint,
directed logarithm bounds certify
\[
 \log(1982759)<\frac{29}{2}<\log(1982760).
\]
This identifies the next analytic range without assigning a
conclusion to it here.
\end{proof}

\begin{theorem}[The initial interval through order 1982759]
\label{thm:initial-prefix-1982759}
Every tree of order at most $1982759$ has a prime labeling.
\end{theorem}

\begin{proof}
Theorem~\ref{thm:initial-prefix-559999} handles the preceding
orders. For $560000\leq n\leq1982759$, prepare any tree as in
Proposition~\ref{prop:highest-degree-preparation}. A zero
residual degree is handled there. Otherwise select the unique
interval containing $n$ and its row for the residual degree.
The structural constructions give exact even-type quotas,
full odd-position tags, and the paid absorbed set. The exact
certificate bounds the two type laws, the local conflict
probability, and the global number of deleted matching edges.
Their simultaneous success follows as in
Theorem~\ref{thm:initial-prefix-559999}. The asymmetric Hall
criterion gives a perfect matching of the remaining labels
and positions. Finally restore the core and special labels
and apply the sparse repair. All intervals meet at consecutive
integers, so no order is omitted.
\end{proof}

\section{Logarithmic bridges from 29/2 to 256}
\label{sec:part-log-bridges}

This part proves that every tree of order $n$ with
$29/2\leq\log n\leq256$ has a prime labeling. The interval is divided
into nineteen closed intervals of $\log n$, and each is covered by one
bridge theorem. Every bridge prepares the tree as in
Section~\ref{sec:part-structural}, sometimes with a variant proved in
its subsection, assigns the even labels with exact type quotas, chooses
a permutation of the even labels within exact classes, and places the
odd labels by a Hall criterion. Its finite
certificate divides the interval into closed cells, either rectangles
in $\log n$ and a normalized residual degree or unit intervals of
$\log n$, and lists rational inequalities that are checked on every
cell and imply the hypotheses of the construction.

Table~\ref{tab:log-bridges} lists the bridges. Each bridge after the
first uses the construction of an earlier one with the changes listed
in the table; most bridge subsections begin by stating which parts are
unchanged. Section~\ref{sec:adaptive-refinement-hall}
proves the degree-balanced refinement and the Hall criterion with empty
positions on which the third bridge and several later ones rely. The
union of the nineteen intervals is
taken in Section~\ref{sec:log-range}
(Theorem~\ref{thm:log-bridges-all}).

\begin{table}[!tbp]
\centering
\caption{The logarithmic bridges. The theorem number also gives the
subsection in which the bridge is proved.}
\label{tab:log-bridges}
\par\vspace{.5\baselineskip}
\renewcommand{\arraystretch}{1.1}
\begin{tabular}{@{}rllp{6.4cm}@{}}
\hline
&Range of $\log n$&Theorem&\raggedright Main change\tabularnewline
\hline
1&$[29/2,16]$&\ref{thm:first-log-bridge}&\raggedright
  all full types assigned at once; exchanges of empty-type colors at
  cut vertices\tabularnewline
2&$[16,343/20]$&\ref{thm:second-log-bridge}&\raggedright
  the same full-type assignment without exchanges\tabularnewline
3&$[343/20,20]$&\ref{thm:third-log-bridge}&\raggedright
  degree-balanced rare-type refinement; label-side reserve
  (Section~\ref{sec:adaptive-refinement-hall})\tabularnewline
4&$[20,128/5]$&\ref{thm:fourth-log-bridge}&\raggedright
  protected positions; localized matching estimate\tabularnewline
5&$[128/5,53/2]$&\ref{thm:fifth-log-bridge}&\raggedright
  even labels sorted by degree within exact types; integrated
  shortfalls on both sides\tabularnewline
6&$[53/2,143/5]$&\ref{thm:sixth-log-bridge}&\raggedright
  a fixed number of centroid rounds\tabularnewline
7&$[143/5,144/5]$&\ref{thm:seventh-log-bridge}&\raggedright
  centroid rounds counted on the position side; centered bound for
  larger-prime deletions\tabularnewline
8&$[144/5,29]$&\ref{thm:eighth-log-bridge}&\raggedright
  periodic type and quota errors; one-side square-root
  bound\tabularnewline
9&$[29,31]$&\ref{thm:ninth-log-bridge}&\raggedright
  the unweighted one-side bound\tabularnewline
10&$[31,32]$&\ref{thm:tenth-log-bridge}&\raggedright
  second-order reserve; arithmetic bounds grouped by product
  size\tabularnewline
11&$[32,35]$&\ref{thm:eleventh-log-bridge}&\raggedright
  smaller protected reserve; uncapped slot shortfalls\tabularnewline
12&$[35,37]$&\ref{thm:twelfth-log-bridge}&\raggedright
  minimum compatibility instead of integrated shortfalls\tabularnewline
13&$[37,39]$&\ref{thm:thirteenth-log-bridge}&\raggedright
  coarser product bins; one type threshold\tabularnewline
14&$[39,42]$&\ref{thm:fourteenth-log-bridge}&\raggedright
  two-sided rounding bound; symmetric moment bounds\tabularnewline
15&$[42,48]$&\ref{thm:fifteenth-log-bridge}&\raggedright
  truncated expansions with a sharp remainder\tabularnewline
16&$[48,60]$&\ref{thm:sixteenth-log-bridge}&\raggedright
  twelve unit intervals with parameters fixed by rule\tabularnewline
17&$[60,72]$&\ref{thm:seventeenth-log-bridge}&\raggedright
  one common inflation of the selection probabilities; global
  compatibility\tabularnewline
18&$[72,96]$&\ref{thm:eighteenth-log-bridge}&\raggedright
  a coarse bound on the residual degree\tabularnewline
19&$[96,256]$&\ref{thm:nineteenth-log-bridge}&\raggedright
  the preparation of the unbounded range\tabularnewline
\hline
\end{tabular}
\end{table}

\subsection{The first logarithmic bridge}
\label{sec:empty-separator-log-bridge}

The integer certificates stop immediately below $\log n=29/2$.
For the next range the number of possible full types is still at most
$256$, so we can assign them all at once. The difficulty is geometric:
an ordinary cut at a high-degree vertex may cost too many positions.
We first record a bound that distinguishes the number of cuts from the
number of neighbors they expose.

For a bipartite forest $G=(A,B)$ put
$N=\abs A+\abs B$ and
\[
 R_A=\sum_{a\in A}(d_G(a)-1)_+.
\]
For positive integers $d,D$ set
\begin{align}
 \mathcal C_*(N,d)
   &=\sum_{j\geq0}\min\left\{d,
                   \floor*{\frac{N}{2^j+1}}\right\},
                       \label{eq:log-cut-count}\\
 \mathcal B_*(N,d,D)
   &=\sum_{j\geq0}\min\{D,2^{j+1}-1\}
          \min\left\{d,\floor*{\frac{N}{2^j+1}}\right\}.
                       \label{eq:log-cut-boundary}
\end{align}
Both sums are finite. The term at $j=0$ is intentionally present.

\begin{lemma}[Simultaneous full-type assignment]
\label{lem:log-simultaneous-assignment}
Let $G=(A,B)$ be a bipartite forest, $M=\abs A>0$, and
$D=\max(1,\Delta_G(A))$. Prescribe nonnegative integer quotas
$a_0,\ldots,a_{r-1}$ of sum $M$, and write $\pi_i=a_i/M$.
There are colors on $A$, tags on $B$, and a set of cut vertices
$C$ such that the $A$-quotas are exact and
\begin{align*}
 H&=\frac12\sum_i\left|
       \sum_{\substack{a\in A\\c(a)=i}}d_G(a)-e(G)\pi_i
                                      \right|\leq D(r-1),\\
 \frac12\sum_i\left|\abs{\{b:j(b)=i\}}-\abs B\pi_i\right|
                                      &\leq\frac r4.
\end{align*}
With $d=3(r-1)$, put
$Z_C=(C\cap B)\cup N_G(C\cap A)$. Then every edge whose
$B$-end avoids $Z_C$ has equal endpoint color and tag, and
\begin{equation}\label{eq:log-assignment-boundary}
 \abs C\leq\mathcal C_*(N,d),\qquad
 \abs{Z_C}\leq
 \min\{\mathcal B_*(N,d,D),\,\mathcal C_*(N,d)+R_A\}.
\end{equation}
For $r=1$ take $C=\varnothing$ and color both sides alike.
\end{lemma}

\begin{proof}
Use the three additive coordinates consisting of $A$-count,
$A$-degree sum in the original forest, and $B$-count. Begin with
each component carrying the vector of fractional color weights
$\pi$. The pruning argument of
Theorem~\ref{thm:original-type-assignment} preserves all three
coordinates and leaves at most $d=3(r-1)$ mixed connected
objects. Cut a largest nonsingleton mixed object at a centroid,
giving its weight vector to the cut vertex and to every component
left by its deletion, and prune again. Frozen objects stay frozen.
After finitely many cuts, every mixed object is a singleton.

An object in the bucket $2^j<\abs U\leq2^{j+1}$ has every child
of size at most $2^j$. It cannot return to that bucket after a cut.
At most $d$ mixed objects occur there, and the objects cut in one
bucket are disjoint, so at most $\floor{N/(2^j+1)}$ are cut.
This proves the first inequality of
\eqref{eq:log-assignment-boundary}. An $A$-centroid in that bucket
has at most $\min(D,2^{j+1}-1)$ neighbors; a $B$-centroid is charged
one. Summing gives the first bound on $\abs{Z_C}$.
Alternatively each $A$-centroid costs at most one plus
$(d_G(a)-1)_+$, while each $B$-centroid costs one. The
$A$-centroids are distinct, proving the second bound.

Next reduce the mixed $A$-singletons once more, now preserving
only their color counts and color degree sums. These are
$2(r-1)$ constraints, so the same dimension argument leaves at
most $f_A\leq2(r-1)$ mixed $A$-singletons; newly pure atoms are
frozen. For each color, its target count minus its number of
frozen vertices is a nonnegative integer, and these integers sum
to $f_A$. Assign the mixed $A$-singletons these remaining counts.
The $A$-quotas are then exact. Let $\eta_{v,i}$ be the resulting
color indicators and $\theta_{v,i}$ the fractional weights just
before this assignment. Exact counts give
$\sum_v(\eta_{v,i}-\theta_{v,i})=0$ for every color, and the
degree sums were preserved until this step. Therefore
\[
 \sum_{\substack{a\in A\\c(a)=i}}d_G(a)-e(G)\pi_i
 =\sum_v\Bigl(d_G(v)-\frac D2\Bigr)(\eta_{v,i}-\theta_{v,i}),
\]
the sum being over the mixed $A$-singletons. For each such $v$,
$\abs{d_G(v)-D/2}\leq D/2$ and
$\frac12\sum_i\abs{\eta_{v,i}-\theta_{v,i}}\leq1$. Hence
$H\leq Df_A/2\leq D(r-1)$. The centering by $D/2$ uses the
exact counts. Round the $B$-counts by
giving the largest remainders the necessary upward roundings.
The same calculation as in
Theorem~\ref{thm:original-type-assignment} gives $r/4$.
An edge with neither endpoint a cut vertex remains in one
frozen monochromatic object. Thus the singleton roundings do not
change any edge whose $B$-end avoids $Z_C$.
\end{proof}

The next exchange lets us keep the tags while reducing the
geometric boundary. Color $0$ will be the empty small-prime
type. Its compatibility with every odd type makes the direction
of the exchange significant. Figure~\ref{fig:donor-exchange}
shows the cut vertex and its low-degree donor.

\begin{lemma}[Empty-type exchanges at the separator]
\label{lem:log-empty-exchange}
Let the assignment of Lemma~\ref{lem:log-simultaneous-assignment}
have empty-type count $a_0$ and empty-type degree sum $E_0$.
For an integer $t\geq1$ assume
\begin{equation}\label{eq:log-donor-condition}
 a_0-\frac{E_0}{t+1}\geq\abs{C\cap A}.
\end{equation}
Without changing any color quota or $B$-tag, one can exchange
colors and choose a set $Z\subseteq B$ such that
\[
 \abs Z\leq\min(t,D)\abs C,
\]
and every edge to $B\setminus Z$ has $A$-color either $0$ or
equal to its $B$-tag. If $E'_i$ denotes the new degree sums, then
\begin{equation}\label{eq:log-one-sided-transfer}
 E'_0=E_0+\Delta,\qquad E'_i=E_i-\Delta_i\quad(i\ne0),
 \qquad
 \Delta=\sum_{i\ne0}\Delta_i
                   \leq D\abs{C\cap A}.
\end{equation}
All the $\Delta_i$ are nonnegative. If $D<t+1$, the exchange
can be omitted and $\Delta=0$.
\end{lemma}

\begin{proof}
Let $S$ be the nonempty-colored vertices in $C\cap A$ of
degree greater than $t$, and let $h$ be the number of
empty-colored vertices in $C\cap A$. At most $E_0/(t+1)$
empty-colored vertices have degree greater than $t$.
By \eqref{eq:log-donor-condition}, at least
$a_0-h-E_0/(t+1)\geq\abs{C\cap A}-h\geq\abs S$
empty-colored vertices outside $C$ have degree at most $t$.
Choose distinct such donors and exchange their colors with
the members of $S$. Every transfer of degree is strictly toward
the empty color, which proves \eqref{eq:log-one-sided-transfer}.

Let $L$ be the nonempty-colored $A$-centroids of degree at most $t$,
and let $U$ be the chosen donors. Take
\[
 Z=(C\cap B)\cup N_G(L)\cup N_G(U).
\]
An edge outside $Z$ from an unchanged noncentroid retains its
equal tag. A high-degree centroid now has empty color, and a
low-degree nonempty centroid or donor has no neighbor outside
$Z$. These cases exhaust the edges. Since
$\abs L+\abs U\leq\abs{C\cap A}$,
\[
 \abs Z\leq\abs{C\cap B}
              +\min(t,D)(\abs L+\abs U)
          \leq\min(t,D)\abs C.
\]
When $D<t+1$, $S$ is empty. This proves the last assertion.
\end{proof}

For clarity, we record how the transfer enters the probability
calculation. Write $P=\{p:3\leq p\leq z,\ p\text{ prime}\}$ and
\[
 \alpha=\prod_{p\in P}(1-1/p),\qquad
 R_y=\sum_{\substack{p>z\\p\mid y}}\frac1p,\qquad
 h_y=\abs{\{p>z:p\mid y\}}.
\]
Let $m=\floor{n/2}$, and let $a_0$ be the number of retained
half-labels of empty $P$-type after the special even labels are
removed. Denote by $C_o$ the full-period ordinary-prefix
$\ell^1$ discrepancy used in the certificate.

\begin{lemma}[Cost within the retained empty class]
\label{lem:log-empty-class-cost}
If $y\leq n$ is odd, the fraction $f_0(y)$ of retained
empty-type even labels sharing with $y$ a prime greater than
$z$ satisfies
\begin{equation}\label{eq:log-empty-local}
 f_0(y)\leq
 \frac{\alpha mR_y+(C_o/2)h_y}{a_0},
 \qquad h_y\leq\frac{\log n}{\log z}.
\end{equation}
If $T_0$ is the sum of these fractions over all odd
$y\leq n$, then
\begin{equation}\label{eq:log-empty-global}
 T_0\leq
  \frac{n\alpha m}{2a_0}S_2(z)
 +\frac{nC_o(1+\log n)}{4a_0}
 +\frac{\log n}{\log z},
\end{equation}
where $S_2(z)$ is an upper bound for
$\sum_{p>z}p^{-2}$.
\end{lemma}

\begin{proof}
For every ordinary prefix, the empty-coordinate discrepancy
is at most $C_o/2$: the coordinate errors sum to zero and
their total absolute value is at most $C_o$.
Multiplication by a prime $p>z$ preserves the $P$-type.
Consequently at most $\alpha\floor{m/p}+C_o/2$ empty-type
half-labels are divisible by $p$. Removing the special labels
can only reduce this numerator, while $a_0$ remains the
denominator. Summing over the large prime divisors of $y$
proves \eqref{eq:log-empty-local}. Their product is at most
$y\leq n$, proving the second assertion there.

For the sum over odd $y$, there are at most $n/(2p)+1$
odd multiples of $p$ in $[n]$. The product of the
$n/(2p)$ term with $\alpha m/p$ gives the first term of
\eqref{eq:log-empty-global}. Summing the discrepancy terms
uses $\sum_{z<p\leq m}1/p\leq1+\log n$ and gives the second.
For the remaining $+1$ term, average the number of large
prime divisors over the retained empty labels; each has at
most $\log n/\log z$ of them. This gives the last term.
\end{proof}

\subsubsection*{The certificate on the first logarithmic interval}

We give the numerical predicates so that the finite computation
has a precise mathematical meaning. Let a closed logarithmic
interval be $[a,J]$ and let the degree coordinate
$x=2w(D-1)/n$ belong to $[x_-,x_+]$. Here $w$ is the number
of special even labels, placed on the $w$ highest-degree
vertices of the independent side before taking the residual
forest. The prime-count bounds of
\eqref{eq:dusart-denominator-lower} and
\eqref{eq:dusart-denominator-upper}, with their domains checked
at $a$, give
\[
 w>nW(\log n),\quad
 W(J)=\frac1{2(J-1693/1000)}
       -\frac1{3(J-2199/1000)}.
\]
Indeed, the special even labels include every $2p$ with
$n/3<p\leq n/2$. The two cited prime-count inequalities
give respectively
\[
 \pi(n/2)\geq
 \frac{n}{2(\log n-1-\log2)},\qquad
 \pi(n/3)\leq
 \frac{n}{3(\log n-11/10-\log3)}.
\]
The directed comparisons
$1693/1000<1+\log2$ and
$2199/1000>11/10+\log3$ yield the displayed lower
bound for $w$; $W$ decreases on the certified range.
The two elementary degree counts in the independent side give
\[
 0\leq x<1,\quad
 D\leq D_+:=1+\frac{x_+}{2W(J)},\quad
 \frac{e(G)}n\leq e_+:=1-W(J)-\frac{x_-}{2}.
\]
Indeed, $(w+1)(D-1)<n/2$ and
$e(G)\leq n-1-wD$; the strict prime bound supplies both
displayed majorants.

Let $Q$ be a directed rational upper bound for $e^{-a}$,
$r=2^{\abs P}$, $d=3(r-1)$, and put
\[
 \omega\geq\frac{4JQ+(a+4)/(3a(a-2))}{1-Q},
\]
choosing the directed rational upper enclosure of the
right-hand side.
The upper prime-count estimates and
$\floor{\log_2 n}<2\log n$ give $w/m\leq\omega$;
the certificate checks its domain and positivity.
The $r$ full types are
retained, with product law $\mu$ and depleted law $\nu_\xi$
of \eqref{eq:coupling-depletion}, where $\xi=w/m\leq\omega$.
Take directed upper bounds
$\varepsilon_E\geq C_oQ/(1-Q)$ and
$H_+\geq D_+(r-1)Q$. In the formulas below every value used by
the verifier is rounded in the unfavorable direction.
The normalized cut count is at most the infinite sum
\[
 \sum_{j\geq0}
       \min\left\{dQ,\frac1{2^j+1}\right\}.
\]
The certificate calls its rational upper enclosure $C_+$:
it sums $64$ terms and adds the geometric remainder.
Set
\[
 d_0=\alpha-\omega-\varepsilon_E,\quad
 q_0=\frac{(1-Q)d_0}{2},\quad
 \bar d=\frac{2e_+}{(1-Q)(1-\omega)}.
\]
The certificate checks $d_0>0$, $\bar d<t+1$, and
\begin{equation}\label{eq:log-donor-row}
 q_0\left(1-\frac{\bar d}{t+1}\right)
          -\frac{H_+}{t+1}\geq C_+.
\end{equation}
Since $a_0/n\geq q_0$ and
$E_0\leq e(G)a_0/\abs A+H$, where
$\abs A=m-w\geq(1-\omega)m\geq(1-Q)(1-\omega)n/2$ gives
$e(G)/\abs A\leq\bar d$, the positive coefficient in
\eqref{eq:log-donor-row} proves
\eqref{eq:log-donor-condition}. The resulting position and
transfer charges are
\[
 \chi=\min(t,D_+)C_+,\qquad
 \delta=
 \begin{cases}
 0,&D_+<t+1,\\ D_+C_+,&D_+\geq t+1.
 \end{cases}
\]
The accepted rows have $t\in\{1,2\}$.

For a fixed odd $y$, the one-sided identity
\eqref{eq:log-one-sided-transfer} bounds the added expected
deletion cost by $\Delta f_0(y)$. For the total deletions it
bounds the added mean by $\Delta T_0$. This follows by writing
the post-exchange weighted mean as
\[
 \sum_iE'_if_i=\sum_iE_if_i+\Delta f_0
                      -\sum_{i\ne0}\Delta_i f_i
             \leq\sum_iE_if_i+\Delta f_0.
\]
Thus Lemma~\ref{lem:log-empty-class-cost} applies without
replacing the transfer by an unsigned variation error.

For completeness, the certificate computes the remaining
quantities as follows. Put $o=(1+Q)/2$,
$\kappa=(2J+2)Q$, and
\begin{align*}
 \eta&=\frac{o\varepsilon_E}{1-\omega}+\frac{rQ}{4},
 &\ell_-&=\frac12-\kappa-\chi,\\
 V&=\max\left\{
   \frac{o\varepsilon_O+\kappa+\chi}{\ell_-},
   \frac{\eta+\chi}{\ell_-}\right\},
 &g_-&=\frac2{3J}-\kappa-\chi.
\end{align*}
Here $\varepsilon_O=C_dQ$ is the odd-prefix discrepancy.
The depletion coupling has
$b=(\alpha-\omega)/(1-\omega)$ and
$\theta=\frac12\sqrt{\alpha/(\alpha-\omega)}$.
The finite cutoff set is $z\in\{19,23\}$; its rigorous
prime-square tail is denoted $S_2(z)$. It is obtained by
enumerating the primes from $z$ through $10000$ and applying
Lemma~\ref{lem:prime-square-tail} beyond $10000$.
With $\zeta=(\log z)_-$, define
\[
 A_0=\frac{\alpha}{d_0},\qquad
 B_0=\frac{C_oJQ}{d_0(1-Q)\zeta},\qquad
 T_+=\frac{\alpha S_2(z)}{2d_0}
   +\frac{C_o(1+J)Q}{2d_0(1-Q)}
   +\frac{JQ}{\zeta}.
\]
The uniform global deletion budget is
\begin{equation}\label{eq:log-empty-mass}
 \mathcal M=\frac{1001}{1000}
 \left\{\frac{e_+S_2(z)}{2(1-\omega)}
       +\frac{JH_+}{2z\zeta}
       +\frac{JQ}{\zeta}
       +\delta T_+\right\}.
\end{equation}
The local bound and the global mean are those of
Lemma~\ref{lem:localized-permutation} in
Section~\ref{sec:adaptive-refinement-hall}, applied within every exact
type class; its proof uses no result of the present section, and the
exchanged degree enters only through the terms with $\delta$, as
explained above. The local bound uses factor $32/31$ and additive term
$32D_+(J+(\log10000)_+)Q$. For each full type $I$ the
verifier bounds $R_y$ by $R_I^+$ using the smallest eligible
large primes whose product can divide $y\leq e^J$. It takes
one further prime to certify that the product exceeds $e^J$,
so the list is a proved truncation rather than a presumed
limit on prime factors. The verifier takes
the minimum of
\begin{equation}\label{eq:log-empty-local-margin}
 \begin{split}
 \frac{c(I)-\omega}{1-\omega}
 &-\frac{32}{31\ell_-}
    \left[\left(\frac{e_+}{1-\omega}+\delta A_0\right)R_I^+
                   +H_++\delta B_0\right]\\
 &-\frac{32D_+(J+(\log10000)_+)Q}{\ell_-}
 \end{split}
\end{equation}
over all $r$ types, clipped below at zero, and calls it
$\sigma$. No type is dropped from this minimum.

For the rational slot threshold $s$, put
\[
 a_s=\frac{o}{1-\omega}
       \E_\mu(s-c(J))_++\eta(s-\alpha)_+,
       \qquad P_s=o a_s.
\]
Lemma~\ref{lem:clipped-slot-shortfall} gives this bound for the
clipped slot shortfall without listing the actual tags.
The final row checks
\[
 0<g_-\leq\tau\ell_-,\quad V<\tau\leq\tfrac12,\quad
 \sigma-V>\tau,\quad s-V>\tau
\]
and the three strict comparisons
\begin{align}
 \mathcal M&<g_-((s-V)\ell_--g_-)-P_s,\notag\\
 \mathcal M&<\ell_-^2\tau(s-V-\tau)-P_s,\label{eq:log-empty-hall}\\
 \mathcal M&<\ell_-^2b\bigl((1-\theta)\tau(1-\tau)
                             -V+(1-\theta)V^2\bigr).\notag
\end{align}
They are precisely the narrow, middle, and spectral
alternatives of
Theorem~\ref{thm:asymmetric-reserve-hall}, after multiplying
by $n^2$.

\begin{proposition}[The first logarithmic certificate]
\label{prop:first-log-certificate}
The parameter table partitions $[29/2,16]$ into $30$
consecutive closed logarithmic intervals. On each, its closed
degree intervals cover $[0,1]$. Its $83$ rectangles satisfy
all $2905$ rational conditions, including
\eqref{eq:log-donor-row} and
\eqref{eq:log-empty-hall}. Every deletion-budget ratio and
every absorbed-position-to-universal-supply ratio is strictly
less than $995/1000$.
\end{proposition}

\begin{proof}
The table and exact evaluator are supplied in the
\texttt{first\_log\_bridge} directory of the supplement.
The evaluator verifies the interval joins and all degree
endpoints as rational numbers. It recomputes the eight complete
small-prime residue periods by integer enumeration; periodicity,
not a sample of orders, proves the prefix discrepancies for
every $n$. It encloses exponentials, logarithms, square roots,
and prime tails in directed rational bounds, and adds an
infinite remainder to the dyadic separator sum. It verifies
the prime-count domains, preparation, donor supply, both
type-law errors, and every condition above on each whole
rectangle. Reproducing the frozen report gives the stated
counts and strict ratios.
\end{proof}

\begin{theorem}[Orders in the first logarithmic interval]
\label{thm:first-log-bridge}
Every tree of order $n$ with
$29/2\leq\log n\leq16$ has a prime labeling.
\end{theorem}

\begin{proof}
Choose the balanced independent side and the prime-labeled
core as in Sections~\ref{sec:balance} and
\ref{sec:matching}. Put all special even labels on vertices
of highest original degree, and let $G=(A,B)$ be the
remaining forest. The checked prime supply makes the
preparation valid throughout this interval. Select a
certificate rectangle containing $\log n$ and
$x=2w(D-1)/n$.

Lemma~\ref{lem:log-simultaneous-assignment} gives exact
full-type quotas. Lemma~\ref{lem:log-empty-exchange},
whose hypothesis is \eqref{eq:log-donor-row}, leaves every
tag fixed while ensuring that an odd label disjoint from a
tag is coprime at all adjacent retained positions with
respect to primes at most $z$. Assign the absorbed
positions distinct universal odd labels. Their remaining
number is at least $g_-n$. The odd and slot type laws
have variation at most $V$ from $\mu$ and $\nu_\xi$.

Permute even labels independently inside their exact full
types. Equations~\eqref{eq:log-empty-local}--
\eqref{eq:log-empty-mass} bound the local and total conflicts
from primes greater than $z$. The local failure probability
is at most $1/10000$; Markov's inequality gives a global
failure probability at most $1000/1001$. Their sum is less
than one, so one permutation satisfies both bounds.
The comparisons \eqref{eq:log-empty-hall} and
Theorem~\ref{thm:asymmetric-reserve-hall} give a perfect
matching. Restore the absorbed positions, special labels,
and core, then apply the sparse parity-preserving repair to
the possible $p,2p$ pairs. Every edge of the resulting
bijection has coprime labels.
\end{proof}

\subsection{A simultaneous full-type bridge}
\label{sec:simultaneous-full-type-bridge}

Above logarithm $16$ we can use the full-type assignment of
Lemma~\ref{lem:log-simultaneous-assignment} without the donor
exchanges. Its ordinary boundary becomes small enough to be
absorbed by the universal odd labels. We keep the alternative
bound through the excess degree because it is valid on every
certificate rectangle, even though the accepted rows choose
the ordinary bound.

Fix a tree of order $n$ with $16\leq L=\log n\leq343/20$.
As in Section~\ref{sec:empty-separator-log-bridge}, let $S$
be the independent side of size $m=\floor{n/2}$, let $K$ be
the prime-labeled core, and place all $w$ special even labels
on the $w$ highest-degree vertices $F$ of $S$. Put
$G=T-(K\cup F)=(A,B)$, and let
$D=\max(1,\Delta_G(A))$. The remaining even half-labels
have exact quotas in every full $P$-type.

The relation between the largest residual degree and the
number of special labels is the same as in the preceding
section. A second consequence controls the excess degree
outside $F$.

\begin{lemma}[Excess remaining after highest-degree deletion]
\label{lem:full-type-residual-excess}
With $x=2w(D-1)/n$ and
$R_A=\sum_{a\in A}(d_G(a)-1)_+$, one has
\begin{equation}\label{eq:full-type-residual-excess}
 0\leq x<1,\qquad \frac{R_A}{n}<\frac{1-x}{2}.
\end{equation}
\end{lemma}

\begin{proof}
Since $S$ is independent, the sum of its degrees in $T$
counts each edge incident with $S$ at most once. Hence
\[
 \sum_{v\in S}(d_T(v)-1)
           \leq n-1-m<\frac n2.
\]
Every vertex of $F$ has original degree at least $D$:
otherwise a retained vertex of residual degree $D$
would have larger original degree than a chosen member
of $F$. Thus
\[
 \sum_{v\in F}(d_T(v)-1)\geq w(D-1).
\]
Deleting vertices cannot increase degrees. For $a\in A$,
$(d_G(a)-1)_+\leq d_T(a)-1$, since vertices of $S$
have positive original degree unless $n=1$, which is
outside the present range. Subtract the contribution of
$F$ from the first inequality to obtain
$R_A<n/2-w(D-1)$. This is the second assertion; the first
also follows from $(w+1)(D-1)<n/2$, as in the preceding
section.
\end{proof}

We now give the common row predicate. Let $[a,J]\subseteq
[16,343/20]$ be a closed logarithmic interval and
$[x_-,x_+]\subseteq[0,1]$ a closed degree interval.
Choose a cutoff $z\in\{17,19,23\}$, let $P$ be all
odd primes through $z$, and write
\[
 r=2^{\abs P},\qquad d=3(r-1),\qquad
 \alpha=\prod_{p\in P}(1-1/p).
\]
Use a directed rational upper bound $Q\geq e^{-a}$.
The prime-count input and the highest-degree deletion give
\begin{align}
 W(J)&=\frac1{2(J-1693/1000)}
       -\frac1{3(J-2199/1000)},\notag\\
 D_+&=1+\frac{x_+}{2W(J)},&
 e_+&=1-W(J)-\frac{x_-}{2}.
                         \label{eq:full-type-row-preparation}
\end{align}
All denominators are positive on the certified range.
The checker verifies the prime-count domains and the
core and universal-label supply at the lower endpoint
of each logarithmic interval.

For the normalized separator estimates put
\begin{align}
 C_+&\geq\sum_{j\geq0}
             \min\left\{dQ,\frac1{2^j+1}\right\},
                                             \label{eq:full-type-cut-count}\\
 B_+&\geq\sum_{j\geq0}
        \min\{D_+,2^{j+1}-1\}
             \min\left\{dQ,\frac1{2^j+1}\right\},
                                             \label{eq:full-type-boundary}\\
 \chi&=\min\left\{B_+,\,
                     C_++\frac{1-x_-}{2}\right\}.
                                             \label{eq:full-type-chi}
\end{align}
The rational values $B_+,C_+$ are upper enclosures:
each sum is evaluated for $64$ terms and completed by
an explicit geometric remainder. Equations
\eqref{eq:log-assignment-boundary} and
\eqref{eq:full-type-residual-excess} show that the
actual exposed set has size at most $\chi n$.
Unlike Section~\ref{sec:empty-separator-log-bridge},
there is no exchanged degree, so $\delta=0$.

Let $C_o,C_d$ be the certified ordinary and odd full-period
discrepancies. Choose upward rational enclosures of
\begin{align*}
 \omega&=\frac{4JQ+(a+4)/(3a(a-2))}{1-Q},
 &\varepsilon_E&=\frac{C_oQ}{1-Q},\\
 \varepsilon_O&=C_dQ,
 &H_+&=D_+(r-1)Q.
\end{align*}
The remaining even-type law is within
$\varepsilon_E/(1-\omega)$ of the depleted law
$\nu_\xi$, $\xi=w/m\leq\omega$, by
\eqref{eq:coupling-depletion}. Lemma~
\ref{lem:log-simultaneous-assignment} gives the $B$-tag
variation $r/4$ and the degree error $H_+n$.
Put $o=(1+Q)/2$ and $\kappa=(2J+2)Q$, and define
\begin{align}
 \eta&=\frac{o\varepsilon_E}{1-\omega}+\frac{rQ}{4},
 &\ell_-&=\frac12-\kappa-\chi, \label{eq:full-type-law-errors}\\
 V&=\max\left\{
       \frac{o\varepsilon_O+\kappa+\chi}{\ell_-},
       \frac{\eta+\chi}{\ell_-}\right\},
 &g_-&=\frac2{3J}-\kappa-\chi. \notag
\end{align}
The first numerator in $V$ concerns odd labels, the
second the retained position tags. In the latter,
removing at most $\chi n$ positions changes an
unnormalized type count by at most $\chi n$.

The depleted product coupling of
Lemma~\ref{lem:depleted-coupling} has
\[
 b=\frac{\alpha-\omega}{1-\omega},
 \qquad
 \theta=\frac12\sqrt{\frac{\alpha}{\alpha-\omega}}.
\]
The certificate checks $\omega<\gamma$, where
$\gamma$ is in \eqref{eq:coupling-parameters}.
Write $\zeta=(\log z)_-$ and let $S_2(z)$ denote
the same rigorous prime-square tail as in Section~
\ref{sec:empty-separator-log-bridge}. With no donor
transfer, the global part of
Lemma~\ref{lem:localized-permutation} gives the total deletion
budget
\begin{equation}\label{eq:full-type-global-mass}
 \mathcal M=\frac{1001}{1000}
   \left\{\frac{e_+S_2(z)}{2(1-\omega)}
          +\frac{JH_+}{2z\zeta}
          +\frac{JQ}{\zeta}\right\}.
\end{equation}
For every full type $I$, the arithmetic routine bounds the
large-prime reciprocal sum of a label of that type by
$R_I^+$. As in the previous section, it includes one
prime beyond the maximum possible product, and so does
not assume a cutoff on actual prime factors. For each
type $I$, put
\begin{align}
 S_I={}&\frac{c(I)-\omega}{1-\omega}
 -\frac{32}{31\ell_-}
       \left(\frac{e_+R_I^+}{1-\omega}+H_+\right)
       \notag\\
 &-\frac{32D_+(J+(\log10000)_+)Q}{\ell_-}.
                       \label{eq:full-type-label-margin}
\end{align}
The post-deletion label compatibility is at least
$\sigma=\max\{0,\min_I S_I\}$. All $r$ types
enter the minimum. For a rational slot
threshold $s$, write
\[
 a_s=\frac{o}{1-\omega}\E_\mu(s-c(J))_+
                         +\eta(s-\alpha)_+,
 \qquad P_s=o a_s.
\]
Lemma~\ref{lem:clipped-slot-shortfall} gives this
bound without listing the actual tags.

\begin{proposition}[The full-type logarithmic certificate]
\label{prop:full-type-log-certificate}
Thirteen adjacent closed intervals cover
$[16,343/20]$. Their degree subdivisions cover
$[0,1]$ and contain $179$ rectangles in total.
All $5728$ rational conditions in the parameter
table hold. In every rectangle,
\[
 0<g_-\leq\tau\ell_-,\quad
 V<\tau\leq\frac12,\quad
 \sigma-V>\tau,\quad s-V>\tau,
\]
and the three strict Hall comparisons
\eqref{eq:log-empty-hall} hold with
\eqref{eq:full-type-global-mass} and
\eqref{eq:full-type-law-errors}. The deletion-budget
ratio and the absorbed-position-to-universal-supply
ratio are each strictly less than $995/1000$.
\end{proposition}

\begin{proof}
The complete table and a portable evaluator are supplied
in the \texttt{second\_log\_bridge} directory. The
evaluator recomputes the eight residue periods and
validates every continuous interval and degree join.
It uses only integer operations and outward rational
enclosures for the prime, logarithm, square-root,
and exponential inputs. The $64$-term separator sums
carry their infinite tails. Every condition is checked
on the whole closed rectangle using the endpoint
monotonicities in the displayed bounds.

The standalone wrapper reproduces the frozen report
byte for byte, checks the complete cover and all
$5728$ row conditions, and verifies the two strict
rational ratios. Independent finite tests exercise
simultaneous coloring, the depleted coupling, and
the rejection of damaged covers. They check the
implementation; the arbitrary-forest conclusion
comes from Lemma~\ref{lem:log-simultaneous-assignment}.
\end{proof}

\begin{theorem}[Orders in the second logarithmic interval]
\label{thm:second-log-bridge}
Every tree of order $n$ with
$16\leq\log n\leq343/20$ has a prime labeling.
\end{theorem}

\begin{proof}
Prepare the tree as above and choose a rectangle containing
$(\log n,x)$. The full-type assignment gives every exact
even-label quota and bounds its exposed positions by
$\chi n$. Assign those positions distinct universal odd
labels. Every retained edge has equal $A$-type and $B$-tag,
so matching an odd label to a disjoint tag avoids all
small-prime conflicts on that edge.

Permute the even labels independently inside each exact
type class. The local and global bounds of
Lemma~\ref{lem:localized-permutation} give respectively \eqref{eq:full-type-label-margin} and
\eqref{eq:full-type-global-mass}. A local failure has
probability at most $1/10000$; by Markov, the global
failure probability is at most $1000/1001$. Their
sum is less than one. For a permutation satisfying
both, the checked row inequalities and
Theorem~\ref{thm:asymmetric-reserve-hall} give a
perfect matching on the residual labels and positions.
Restore the absorbed labels, special labels, and core.
The sparse repair handles the permitted $p,2p$
pairs, completing a prime labeling of the tree.
\end{proof}

\subsection{Adaptive refinement and one-sided Hall reserves}
\label{sec:adaptive-refinement-hall}

The simultaneous assignment of Section~
\ref{sec:empty-separator-log-bridge} uses one coordinate for
each full type. When the cutoff prime grows, we instead refine
types one prime at a time and combine rare classes temporarily.
The next lemma adds degree balance to the quota construction
of Section~\ref{sec:rare-type-assignment}. All degrees are
measured in the fixed original forest, including after a
separator vertex has been removed.

For integers $N,d,D,t\geq1$ define
\begin{align}
 \mathcal C_t(N,d)
  &=\sum_{j\geq0}
     \min\left\{d,\floor*{\frac{N}{t2^j+1}}\right\},
                      \label{eq:adaptive-cut-count}\\
 \mathcal B_t(N,d,D)
  &=\sum_{j\geq0}\min\{D,t2^{j+1}-1\}
     \min\left\{d,\floor*{\frac{N}{t2^j+1}}\right\}.
                      \label{eq:adaptive-boundary}
\end{align}
These are finite sums. For $t=1$ they are the cut and
boundary bounds of \eqref{eq:log-cut-count}--
\eqref{eq:log-cut-boundary}.

\begin{lemma}[Quotas with degree balance]
\label{lem:adaptive-degree-quota}
Let $G=(A,B)$ be a bipartite forest on $N\geq1$ vertices.
Partition its sides into nonempty classes
$A_1,\ldots,A_{c_A}$ and $B_1,\ldots,B_{c_B}$,
allowing either list to be empty when its side is empty.
Put $d=2c_A+c_B$ and suppose $d\geq1$.
Let $D\geq\max\{1,\Delta_G(A)\}$, $0\leq\rho\leq1$,
and $t\geq1$ be an integer. Prescribe
$0\leq k_i\leq a_i=\abs{A_i}$, and write
\[
 E_i=\sum_{v\in A_i}d_G(v),\qquad
 \delta=\sum_i(k_i-\rho a_i)_+.
\]
There are sets $W\subseteq A$ and $T\subseteq B$
such that
\begin{align}
 \abs{W\cap A_i}&=k_i,\qquad N_G(W)\subseteq T,
                                   \label{eq:adaptive-quota-exact}\\
 \sum_i\left|
   \sum_{v\in W\cap A_i}d_G(v)-\frac{k_i}{a_i}E_i
       \right|
   &\leq D\left(dt+\frac{c_A}{2}\right),
                                   \label{eq:adaptive-degree-error}\\
 \abs{T\cap B_j}&\geq\rho\abs{B_j},
                                   \label{eq:adaptive-tag-lower}\\
 \sum_j\bigl(\abs{T\cap B_j}-\rho\abs{B_j}\bigr)
   &\leq\mathcal B_t(N,d,D)+dt+D\delta.
                                   \label{eq:adaptive-tag-error}
\end{align}
If $A$ is empty, the degree and quota claims are vacuous.
\end{lemma}

\begin{proof}
\emph{Step 1. Preserve all three kinds of coordinates.}
Give every connected component of $G$ selection weight $\rho$.
Its vector records the counts in the $A$-classes, the
degree sums in those classes, and the counts in the
$B$-classes. There are $d$ coordinates. If more than
$d$ objects have fractional weights, a linear dependence
among their vectors permits a nonzero change of weights
preserving all coordinates. Move until one weight reaches
$0$ or $1$ and freeze that object. The number of fractional
objects decreases, so this pruning terminates.

If a fractional object has more than $t$ vertices,
choose one of largest order and split it at a centroid.
Keep the centroid as a singleton object, give it and
every remaining component the parent's weight, and
prune again. The coordinates add under the split.
Every new object is smaller, and no vertex is cut
twice; eventually at most $d$ fractional objects
remain, each of order at most $t$.

In the bucket $t2^j<\abs U\leq t2^{j+1}$, a split
produces only objects of order at most $t2^j$.
At most $d$ objects can be split in this bucket:
at its first split no larger fractional object
remains, and no later split can enter the bucket.
The objects split in one bucket are pairwise
disjoint, so their number is also at most
$\floor{N/(t2^j+1)}$. Thus the number of
centroids is at most $\mathcal C_t(N,d)$.
For each $A$-centroid, charge only its neighbors
outside the final separator; for each
$B$-centroid, charge one. Every charged neighbor
belonged to the object when it was split, so this
external boundary is at most
$\mathcal B_t(N,d,D)$.

\emph{Step 2. Meet the exact quotas.}
Round every remaining fractional object upward.
Let $U$ be the union of objects of weight one.
Only at most $d$ objects of order at most $t$
have changed, so, for $A_i$,
\[
 s_i:=\abs{U\cap A_i}=\rho a_i+x_i,\qquad
 x_i\geq0,\qquad\sum_i x_i\leq dt.
\]
If $e_i^0$ is the degree sum on $U\cap A_i$, the
preserved degree coordinate gives
\begin{equation}\label{eq:adaptive-preliminary-degree}
 \left|e_i^0-\frac{s_i}{a_i}E_i\right|
       \leq Dx_i.
\end{equation}
Indeed, the difference is the sum, over newly rounded
vertices, of the fractional increase times
$d_G(v)-E_i/a_i$, and each absolute summand is at
most $D$ times that increase.

We use the following elementary selection fact.
Among numbers in $[0,D]$, a $k$-subset can be chosen
whose sum differs from $k$ times the overall mean by
at most $D/2$. The target is between the sums of the
$k$ smallest and the $k$ largest numbers, because it
is the average over all $k$-subsets. A sequence of
single-element exchanges from the smallest to the
largest subset can be chosen with nondecreasing sums
and steps at most $D$; take the closer sum at the
first crossing. The empty and full subsets are exact.

In $A_i$, if $k_i\leq s_i$, retain an appropriate
$k_i$-subset of $U\cap A_i$; otherwise retain all
of it and add an appropriate $(k_i-s_i)$-subset
from its complement. Call the resulting union $W$.
In the first case, the deviation in
\eqref{eq:adaptive-preliminary-degree} is multiplied
by $k_i/s_i\leq1$; in the second it is multiplied
by $(a_i-k_i)/(a_i-s_i)\leq1$. The selection adds
at most $D/2$. Empty denominators occur only in
the exact endpoint cases and are omitted. Summing
over $i$ proves \eqref{eq:adaptive-degree-error}.
Moreover, the number of newly added vertices is at
most $\delta$, since $s_i\geq\rho a_i$.

\emph{Step 3. Keep the opposite-side lower counts.}
Let $C$ be the set of cut centroids and put
\[
 T=(U\cap B)\cup(C\cap B)\cup N_G(W).
\]
All fractional objects were rounded upward and all
$B$-centroids are included, so
$\abs{T\cap B_j}\geq\rho\abs{B_j}$.
The rounding itself contributes at most $dt$
positions above the target. Neighbors of selected
$A$-vertices outside their selected objects lie
on the external boundary of $C$, of size at most
$\mathcal B_t(N,d,D)$. The at most $\delta$
newly added $A$-vertices contribute at most
$D\delta$ further neighbors. Summing these
three charges proves \eqref{eq:adaptive-tag-error}.
\end{proof}

The same assignment of even labels must satisfy a bound at
every odd label and a bound on the \emph{total} number of
large-prime conflicts. The next lemma gives both bounds for
one permutation. It uses degrees in the original residual
forest, before any marked position is absorbed.

\begin{lemma}[Local and global bounds for one class permutation]
\label{lem:localized-permutation}
Let $G=(A,B)$ be a bipartite forest, let
$M_A=\abs A=m-w>0$ with $m=\floor{n/2}$, and let
$e=e(G)$. Fix a bijection between $A$ and a set of
$M_A$ half-labels $x\in[m]$. Partition $A$ and those
half-labels by their exact types for the odd primes
at most $z$. Suppose all $A$-degrees are at most
$D\geq1$, and write $a_I$ and $E_I$ for the size
and degree sum of the vertices of type $I$. Put
\[
 \lambda=\frac e{M_A},\qquad
 H=\frac12\sum_I\abs{E_I-\lambda a_I},\qquad
 R_z(y)=\sum_{\substack{p>z\\p\mid y}}\frac1p,
 \qquad
 S_2(z)=\sum_{p>z}\frac1{p^2}.
\]
Independently and uniformly permute the half-labels
within each nonempty exact type class. If $K_y$ is
the sum of the degrees of vertices whose assigned
half-label shares a prime greater than $z$ with an
odd integer $y\leq n$, then for $0<\delta<1$,
\begin{equation}\label{eq:localized-permutation-tail}
 \Pr\left\{K_y>
  \frac{32}{31}\left(
       \frac{em}{M_A}R_z(y)+H\right)
      +32D\log\frac n\delta\right\}
 \leq\frac{\delta}{n}.
\end{equation}
If $M_{\rm bad}$ counts all large-prime-incompatible
label--position pairs, then, provided
$z>1$ and $\log n\leq J$,
\begin{equation}\label{eq:localized-permutation-mean}
 \E M_{\rm bad}\leq
 \frac n2\left\{
     \frac{em}{M_A}S_2(z)
     +\frac{J}{z\log z}H\right\}
     +\frac{eJ}{\log z}.
\end{equation}
Consequently, with positive probability,
the bound in \eqref{eq:localized-permutation-tail}
holds simultaneously for every odd $y\leq n$
with $\delta=1/10000$, while
$M_{\rm bad}$ is at most $1001/1000$ times
the right side of \eqref{eq:localized-permutation-mean}.
\end{lemma}

\begin{proof}
In a class of size $a$, the $k$ half-labels
conflicting with a fixed $y$ occupy a uniform
$k$-subset of vertices. If their degrees are
$w_1,\ldots,w_a\in[0,D]$, the convex comparison
of Lemma~\ref{lem:class-convex-comparison},
applied to the exponential, and the chord
bound on $[0,D]$ give
\[
 \E\exp\left(t\sum_{v\in U}w_v\right)
 \leq\left(\frac1a\sum_{v=1}^a e^{tw_v}\right)^k
 \leq\exp\left(
       \frac{e^{tD}-1}{D}\frac{k}{a}
       \sum_{v=1}^a w_v\right).
\]
The choices in different classes are independent.
Let $\mu_y=\E K_y$. Taking $t=1/(32D)$ and
using $32(e^{1/32}-1)\leq32/31$ yields
\[
 \E e^{K_y/(32D)}
       \leq e^{\mu_y/(31D)}.
\]
Markov's inequality proves a tail bound with
$\mu_y$ in place of the parenthesized expression
in \eqref{eq:localized-permutation-tail}.

Let $k_I(y)$ be the number of half-labels of
type $I$ that conflict with $y$. Then
\[
 \mu_y=\sum_I E_I\frac{k_I(y)}{a_I}.
\]
The measures $(E_I)$ and $(\lambda a_I)$ have
equal total mass, so
\[
 \mu_y\leq\lambda\sum_I k_I(y)+H
        \leq\frac{em}{M_A}R_z(y)+H.
\]
The last inequality counts at most
$\floor{m/p}$ multiples of each prime
$p>z$ dividing $y$. This proves
\eqref{eq:localized-permutation-tail}.

Put $\tau_z(x)=\sum_{p>z,\ p\mid x}1/p$.
Since an integer $x\leq m$ has at most
$\log n/\log z$ distinct prime factors
above $z$,
\[
 0\leq\tau_z(x)\leq\frac{J}{z\log z},
 \qquad
 \sum_{x\leq m}\tau_z(x)\leq mS_2(z).
\]
The same equal-mass argument gives
\[
 \E\sum_{v\in A}d_G(v)\tau_z(x_v)
 \leq\frac{em}{M_A}S_2(z)
       +\frac{J}{z\log z}H.
\]
For each prime $p>z$ dividing $x_v$, there
are at most $n/(2p)+1$ odd multiples of $p$
in $[n]$. Summing their possible incidences
over the neighbors of $v$, and using
$\abs{\{p>z:p\mid x_v\}}\leq J/\log z$,
proves \eqref{eq:localized-permutation-mean}.
Overcounting a pair is harmless.

A union bound makes the probability of
any failed local estimate at most $1/10000$.
Markov's inequality makes the probability
of exceeding $1001/1000$ times the global
expectation at most $1000/1001$.
The two probabilities sum to less than one.
If the expectation vanishes, the global
cost is zero almost surely. Thus the
same permutation satisfies both bounds.
\end{proof}

The residual matching has a larger stock of universal
odd labels than universal positions; in the next
logarithmic range it has no protected universal
positions at all. We therefore state a Hall criterion
for different model laws on the two sides and a
reserve on only the label side.

Let $k\geq1$, and put
\[
 \mu=\bigotimes_{i=1}^k\operatorname{Bernoulli}(a_i),
 \qquad
 \nu=\bigotimes_{i=1}^k\operatorname{Bernoulli}(\rho_i),
\]
where $0<a_i\leq\rho_i<1$ and $a_i+\rho_i<1$.
For $0\leq\eta<1$ write
$\nu_\eta=(1-\eta)\nu+\eta\delta_{\varnothing}$.
Put
\begin{equation}\label{eq:product-hall-coupling}
 b=\prod_i(1-\rho_i),\qquad
 \theta=\max_i
   \sqrt{\frac{a_i\rho_i}
                {(1-a_i)(1-\rho_i)}}<1.
\end{equation}
For a position tag $J$ write
$c_s(J)=\prod_{i\in J}(1-a_i)$.
Thus $c_s(J)$ is the $\mu$-mass of compatible
label types.

\begin{theorem}[Asymmetric Hall criterion with empty positions]
\label{thm:asymmetric-product-hall}
Let $X$ be odd labels and $B$ positions, both of
size $R>0$. Their empirical type laws are within
total variation $V$ of $\mu$ and $\nu_\eta$,
respectively. Start with all pairs of disjoint
types, and delete at most $M$ edges. Assume every
label retains more than $\tau R$ neighbors and
at least $g$ labels remain adjacent to every
position. Choose $g_0$ and $s$ with
\[
 0<g_0\leq g,\quad g_0\leq\tau R,\quad
 0\leq V<\tau\leq\frac12,\quad
 0<s\leq1,\quad s-V>\tau.
\]
For the actual position tags set
\[
 A_B(s)=\sum_{v\in B}(s-c_s(J(v)))_+.
\]
If $M$ is strictly less than each of
\begin{align}
 g_0((s-V)R-g_0)-R A_B(s),&
                    \label{eq:product-hall-narrow}\\
 R^2\tau(s-V-\tau)-R A_B(s),&
                    \label{eq:product-hall-middle}\\
 R^2b\bigl((1-\theta)\tau(1-\tau)
                  -V+(1-\theta)V^2\bigr),&
                    \label{eq:product-hall-wide}
\end{align}
then the remaining graph has a perfect matching.
\end{theorem}

\begin{proof}
Suppose Hall's condition fails, and take nonempty
$U\subseteq X$ and $W=B\setminus N(U)$ with
$\abs U+\abs W>R$. Put $u=\abs U/R$ and
$v=\abs W/R$. A universal label cannot lie in
$U$, so $\abs W>g\geq g_0$. If $u<\tau$,
every label of $U$ has fewer than
$R-\abs W<\abs U<\tau R$ remaining neighbors,
a contradiction. Hence
$u\geq\tau$.

If $v<\tau$, a position of tag $J$ initially has
at least $R(c_s(J)-V)$ compatible labels.
At most $R-\abs U<\abs W$ lie outside $U$.
Using
$c_s(J)\geq s-(s-c_s(J))_+$ and summing
over $W$, the number of initially compatible
pairs in $U\times W$ is strictly greater than
\[
 \abs W\bigl((s-V)R-\abs W\bigr)-R A_B(s).
\]
The quadratic is concave for
$g_0\leq\abs W\leq\tau R$; its minimum
occurs at an endpoint. The two endpoint values
are \eqref{eq:product-hall-narrow} and
\eqref{eq:product-hall-middle}, contrary to
the bound on deleted edges.

It remains that $u,v\geq\tau$. Reduce the two
selected type submeasures fractionally until
their masses sum to one. Truncate them at
$\mu,\nu_\eta$, losing at most $V$ from each side,
and reduce again to model submeasures of
masses $x,y\geq\tau-V$ with $x+y=1-2V$.
Let $\mathcal J_0$ be the coordinate coupling giving probabilities
$a_i,\rho_i,1-a_i-\rho_i$ to
$(1,0),(0,1),(0,0)$.
It is supported on disjoint types. Its density relative to
$\mu\nu$ is at most $1/b$, and its maximal
correlation is at most $\theta$, by the
orthonormal product expansion used in
Lemma~\ref{lem:direct-hall-rectangle}.
The mixture
\[
 \mathcal J_\eta=(1-\eta)\mathcal J_0
                    +\eta\,\mu\otimes\delta_{\varnothing}
\]
couples $\mu$ to $\nu_\eta$ and is still
supported on disjoint types. Its density
relative to $\mu\nu_\eta$ is at most $1/b$:
for a nonempty second type this follows from
the product coupling; at the empty type,
the additional $\eta\mu$ term has density
at most $1/b$ because $b\leq1$.
For functions centered under their respective
marginals, covariance under $\mathcal J_\eta$
is $(1-\eta)$ times covariance under
$\mathcal J_0$. The variance of the
second function under $\nu_\eta$ is at least
$(1-\eta)$ times its variance under $\nu$.
Thus the maximal correlation is at most
$\sqrt{1-\eta}\,\theta\leq\theta$.
Consequently the original rectangle contains
at least
\[
 R^2b\left(xy-\theta
      \sqrt{xy(1-x)(1-y)}\right)
\]
compatible pairs. As $x+y=1-2V$,
$(1-x)(1-y)=xy+2V$, whence the
square root is at most $xy+V$.
Concavity also gives
$xy\geq\tau(1-\tau)-V+V^2$.
Substitution yields
\eqref{eq:product-hall-wide}; all these
pairs would have to be deleted, again a
contradiction.
\end{proof}

\subsection{A logarithmic bridge without reserved positions}
\label{sec:full-reserve-log-bridge}

We next cover the interval \(343/20\leq\log n\leq20\).
Here the number of primes in the cutoff can vary from one
certificate rectangle to another. Keeping all full types as
separate coordinates would make the separator too large. We
therefore use the sequential refinement of
Section~\ref{sec:rare-type-assignment}, with the degree
coordinate added in Lemma~\ref{lem:adaptive-degree-quota}.
The ordinary tags are limited to a prescribed support size;
positions that exceed it receive universal odd labels.

The clipped-shortfall term in the Hall inequalities must be
calculated before those positions are removed. We record the
bound in the form used by the certificate. Write
\(P=\{p_1<\cdots<p_k\}\), \(a_i=1/p_i\), and
\(f_i=1-a_i\), so \(f_1\leq\cdots\leq f_k\).
Let \(c_s(J)=\prod_{i\in J}f_i\). For an integer \(q\geq0\),
let \(J^{[q]}\) consist of the first \(\min(q,\abs J)\)
indices of \(J\), and put
\[
 \sigma_q=\prod_{i=1}^{\min(q,k)}f_i,\qquad
 \Phi_q(t)=\E_\nu\bigl(t-c_s(J^{[q]})\bigr)_+ .
\]
The expectation is over the uncompressed product law
\(\nu=\bigotimes_i\operatorname{Bernoulli}(\rho_i)\).

\begin{lemma}[Shortfall of capped tags]
\label{lem:capped-tag-shortfall}
Suppose \(B^+\) has at most \(on\) positions, each tagged
by a subset of \(P\) of size at most \(q\), or by
\(\dagger\). Let its tag counting measure be within
\(en\) in total variation of
\(\abs{B^+}(c_q)_*\nu\), where \(c_q\) sends sets of
size greater than \(q\) to \(\dagger\). If \(B\)
is obtained by deleting positions of \(B^+\),
including every position tagged \(\dagger\), then
for \(0<t\leq1\),
\begin{equation}\label{eq:capped-tag-shortfall}
 \frac1n\sum_{v\in B}(t-c_s(J(v)))_+
 \leq o\Phi_q(t)+e(t-\sigma_q)_+ .
\end{equation}
In particular \(\Phi_q(t)=0\) when \(t\leq\sigma_q\).
\end{lemma}

\begin{proof}
Define the nonnegative function \(F\) on capped tags by
\(F(J)=(t-c_s(J))_+\) for ordinary \(J\), and
\(F(\dagger)=0\). It vanishes at the empty tag.
The product of any at most \(q\) factors from the
nondecreasing sequence \(f_i\) is at least \(\sigma_q\).
Hence \(0\leq F\leq(t-\sigma_q)_+\).
Total variation changes the integral of a function in
this interval by at most \(en(t-\sigma_q)_+\).
If \(\abs J\leq q\), then \(F(c_q(J))\) is the
summand in \(\Phi_q(t)\). If \(\abs J>q\), it is zero
and is still at most that summand. Thus
\(\E_{(c_q)_*\nu}F\leq\Phi_q(t)\).
Removing positions can only decrease a nonnegative sum,
which proves \eqref{eq:capped-tag-shortfall}.
\end{proof}

\begin{lemma}[Variation after paying marked positions]
\label{lem:marked-tag-variation}
Let $\Xi$ be a counting measure of mass $s$
on ordinary types $J$ of size at most $q$ and
the mark $\dagger$, with
\[
 \bigl\|\Xi-s(c_q)_*\nu\bigr\|_{\rm TV}\leq E,
 \qquad
 \beta=\nu\{J:\abs J>q\}.
\]
Delete the $c=\Xi(\dagger)$ marked positions
and adjoin $h\geq0$ positions of empty type.
The resulting counting measure $\Lambda$ has
mass $\ell=s-c+h$. If $\ell>0$ and
$\eta=h/\ell$, then
\begin{equation}\label{eq:marked-tag-variation}
 \left\|\frac{\Lambda}{\ell}
       -\bigl((1-\eta)\nu
                         +\eta\delta_{\varnothing}\bigr)
 \right\|_{\rm TV}
 \leq\frac{E+s\beta}{\ell}.
\end{equation}
\end{lemma}

\begin{proof}
Write $\delta=\Xi-s(c_q)_*\nu$ and
$d=\delta(\dagger)=c-s\beta$.
The positive mass of the signed measure
$\delta$ is at most $E$. If $d>0$,
its positive mass on ordinary types is
at most $E-d$; if $d\leq0$, that mass
is at most $E$.

Let $\Xi_{\rm ord}$ be $\Xi$ restricted to
ordinary types, extended by zero to all
full types. For an event $A$ of full types,
\[
 \begin{aligned}
 \bigl(\Xi_{\rm ord}-(s-c)\nu\bigr)(A)
  ={}&\delta(A\cap\{\abs J\leq q\})\\
    &+c\nu(A)-s\nu(A\cap\{\abs J>q\}).
 \end{aligned}
\]
If $d>0$, the right side is at most
$(E-d)+c\nu(A)\leq E+s\beta$,
because $c=s\beta+d$. If $d\leq0$,
it is at most $E+c\nu(A)\leq E+s\beta$.
The two measures in the difference have
the same mass $s-c$, so the supremum
of its positive values is their total
variation distance. Adding $h$ empty
positions to both measures does not
change that distance. Divide by $\ell$
to obtain \eqref{eq:marked-tag-variation}.
\end{proof}

The quantity \(\Phi_q(t)\) is a finite product-law
expectation. Its rational upper-bound evaluator scans the primes in
increasing order and retains at most \(q\) selected
factors. If the current product is already at most \(t\),
the remaining conditional expectation follows from
the corresponding truncated product moment; if even the
smallest possible remaining product exceeds \(t\), the
branch contributes zero. All values are rational, with
outward rounding at the finite arithmetic resolution.

\subsubsection*{Upper parameters and arithmetic inputs}

A row cannot use the actual selection probabilities
\(\rho_i\) of the refinement, since they depend on the
order. It uses rational upper bounds
\(\bar\rho_i\geq\rho_i\), and it bounds several arithmetic
errors by finite sums. We justify these substitutions
here. The statements are used in this section and in all
later logarithmic bridges.

\begin{lemma}[Evaluation at upper probabilities]
\label{lem:upper-probability-evaluation}
Let \(0\leq\rho_i\leq\bar\rho_i\leq1\) for \(1\leq i\leq k\),
and let \(\nu,\bar\nu\) be the laws of the random sets
\(J\subseteq[k]\) with independent coordinates of
probabilities \(\rho_i\) and \(\bar\rho_i\), respectively.
\begin{enumerate}
\item For every integer \(q\geq0\),
  \(\nu(\abs J>q)\leq\bar\nu(\abs J>q)\).
\item Let \(0<f_1\leq\cdots\leq f_k\leq1\), let \(q\geq0\)
  be an integer and \(0\leq t\leq1\). Write \(J^{[q]}\) for
  the first \(\min(q,\abs J)\) indices of \(J\),
  \(c_s(J)=\prod_{i\in J}f_i\), and
  \[
   \psi(J)=\Bigl(t-\prod_{i\in J^{[q]}}f_i\Bigr)_+ .
  \]
  Then \(\psi\) is nondecreasing under inclusion,
  \[
   (t-c_s(J))_+\boldsymbol1_{\abs J\leq q}\leq\psi(J)
     \leq(t-c_s(J))_+,
  \]
  and consequently
  \(\E_\nu[(t-c_s(J))_+\boldsymbol1_{\abs J\leq q}]
     \leq\E_\nu\psi\leq\E_{\bar\nu}\psi\).
\end{enumerate}
\end{lemma}

\begin{proof}
Take independent uniform variables \(U_i\) on \([0,1]\) and
put \(X_i=\boldsymbol1_{U_i\leq\rho_i}\) and
\(\bar X_i=\boldsymbol1_{U_i\leq\bar\rho_i}\). The random
sets \(J=\{i:X_i=1\}\) and \(\bar J=\{i:\bar X_i=1\}\) have
laws \(\nu,\bar\nu\), and \(J\subseteq\bar J\) for every
outcome. The event \(\abs J>q\) is preserved under
enlargement, which proves the first assertion.

If \(\abs J\leq q\), then \(J^{[q]}=J\) and the left
comparison is an equality; otherwise its left side is zero.
Since the factors increase, \(J^{[q]}\) retains the smallest
factors of \(J\), and their product is at least
\(c_s(J)\); this is the right comparison. Adding an index to
\(J\) either multiplies the retained product by a further
factor at most one, replaces a retained factor by a smaller
one, or leaves the retained set unchanged. The product does
not increase, so \(\psi\) does not decrease. With
\(J\subseteq\bar J\) this gives
\(\E_\nu\psi\leq\E_{\bar\nu}\psi\).
\end{proof}

The first part justifies computing the overflow
probability \(\beta\) as the exact Bernoulli tail with the
parameters \(\bar\rho_i\). The recursion keeps the
probabilities \(v_0,\ldots,v_q\) of the states \(0,\ldots,q\)
and \(v_*\) of the absorbing state \(>q\). A coordinate of
probability \(p\) maps them to
\[
 (1-p)v_0,\qquad (1-p)v_j+pv_{j-1}\ (1\leq j\leq q),\qquad
 v_*+pv_q .
\]
All coefficients are nonnegative, so rounding every value
upward preserves an upper bound; no upper partial sum is
subtracted from one. The function \(\Phi_q(t)\) of
Lemma~\ref{lem:capped-tag-shortfall} is \(\E_\nu\psi\).
By the second part it may be evaluated with the parameters
\(\bar\rho_i\) by the recursion described above, which
conditions on successive coordinates and uses the exact
product \(c\) already selected.

Replacing \(\rho_i\) by \(\bar\rho_i\) decreases \(b\) and
every label compatibility \(\prod_{i\in I}(1-\rho_i)\), and
increases \(\theta\). A label shortfall increases when its
compatibility decreases. A positive wide budget is
nondecreasing in \(b\) and nonincreasing in \(\theta\), and
the narrow and middle budgets do not contain \(b\). Hence
each of these substitutions weakens the certified
inequalities.

The number of coordinates in a refinement step can be
bounded by the classes that can actually occur.

\begin{lemma}[Possible even classes]
\label{lem:possible-even-classes}
Suppose \(\log n\leq J\). Let \((\log2)_-\leq\log2\) be
rational, \(w_p=\floor{64\log_2p}\), and
\(R_J=\floor{64J/(\log2)_--64}\). If \(x\leq m=\floor{n/2}\)
and \(I\) is the set of odd primes dividing \(x\), then
\(\sum_{p\in I}w_p\leq R_J\). Consequently, at the
refinement step of \(p_i\), the ordinary even classes are
among the \(N_i\) subsets of \(\{p_1,\ldots,p_{i-1}\}\)
with at most \(q_A\) elements and weight at most \(R_J\),
and the number of count coordinates is at most
\begin{equation}\label{eq:possible-class-dimension}
 N_i+\sum_{j=0}^{\min(q_B,i-1)}\binom{i-1}j+2 .
\end{equation}
The two additional coordinates are the combined rare even
class and the absorbing position class.
\end{lemma}

\begin{proof}
The product \(d_I\) divides \(x\leq n/2\), so
\[
 \log_2d_I\leq\frac J{\log2}-1\leq\frac J{(\log2)_-}-1 .
\]
Since \(w_p\leq64\log_2p\), the integer
\(\sum_{p\in I}w_p\) is at most \(64\log_2d_I\), hence at
most \(R_J\). The count
\(N_i\) is the sum of the coefficients of \(y^jt^s\), with
\(j\leq q_A\) and \(s\leq R_J\), in
\(\prod_{p<p_i}(1+yt^{w_p})\). All weights are positive, so
truncating the powers beyond \(R_J\) never discards an
admissible subset. The ordinary position tags have at most
\(q_B\) primes. Empty classes only decrease the actual
dimension.
\end{proof}

\begin{lemma}[Restoring rare quotas]\label{lem:rare-restoration}
Exact quotas are restored at the following cost.
\begin{enumerate}
\item Let disjoint classes have sizes \(r_i\), target quotas
  \(k_i\) and current selections \(s_i\), with
  \(0\leq k_i,s_i\leq r_i\). Suppose that the quotas and the
  selections both sum to \(K\). Then all quotas can be
  restored with
  \[
   \sum_i(k_i-s_i)_+\leq\min\Bigl(K,\sum_ir_i-K\Bigr)
  \]
  additions.
\item Let \(p_1<\cdots<p_k\) be the odd primes of \(P\), and
  call an even label \(2x\) rare at step \(i\) if \(x\) has
  more than \(q\) prime divisors among \(p_1,\ldots,p_{i-1}\).
  If \(K_i\) rare labels among \(2,4,\ldots,2m\) are divisible
  by \(p_i\), possibly after deleting labels of empty type,
  then
  \[
   K_i\leq\frac m{p_i}e_{q+1}(\{p_1,\ldots,p_{i-1}\}),\qquad
   \sum_iK_i\leq m\,e_{q+2}(P),
  \]
  where \(e_j(P)\) is the sum of \(\prod_{p\in I}p^{-1}\) over
  all \(I\subseteq P\) with \(\abs I=j\).
\end{enumerate}
\end{lemma}

\begin{proof}
In each class keep \(\min(k_i,s_i)\) selected vertices, delete
the excess and add the missing ones. The number of additions
is the displayed sum. It is at most \(\sum_ik_i=K\) and at
most \(\sum_i(r_i-s_i)\).

For the second part, the indicator of rarity is at most the
number of \((q+1)\)-sets \(I\) of old prime divisors of
\(x\). The number of \(x\leq m\) divisible by \(p_id_I\) is
at most \(m/(p_id_I)\). This gives the bound for \(K_i\).
Summing over \(i\), each \((q+2)\)-set is counted once, at its
largest prime. Deleted labels of empty type are not rare.
No independence of divisibilities is used.
\end{proof}

At step \(i\) the combined rare even class receives its
exact total quota \(K_i\). By the first part, restoring
the individual exact quotas inside it needs at most
\(K_i\) additions, each of which changes the tag count and
the degree sums by at most \(D\). With \(q=q_A\) and
\(m\leq n/2\), the normalized restoration cost is at most
\(D_+\,\mathrm{rare}\), where
\(\mathrm{rare}=e_{q_A+2}(P)/2\).

\begin{lemma}[Quota defects by product sizes]
\label{lem:quota-product-bins}
Let \(P\) be the odd primes at most \(z\), in increasing
order, and let the refinement choose \(\rho_i\) by
\eqref{eq:rare-empty-adjustment}. For the even
half-labels in \([m]\), possibly after deleting labels of
empty type, the positive quota defects \(\delta_i\) of
Lemma~\ref{lem:rare-quota-defect} satisfy
\begin{equation}\label{eq:quota-mobius}
 \sum_i\delta_i\leq
 \sum_{\substack{J\subseteq P\\\abs J\geq2}}
   2^{\abs J-2}\min\left(1,\frac m{d_J}\right),
 \qquad d_J=\prod_{p\in J}p.
\end{equation}
Let \(B\geq1\) be an integer, \(w_p=\floor{B\log_2p}\),
\(C_s=[t^s]\prod_{p\in P}(1+2t^{w_p})\),
\(C_s^{(2)}=C_s-2\,\#\{p\in P:w_p=s\}\) for \(s>0\), and let
\(U_s\geq2^{-s/B}\) be rational. If \(Q\geq1/n\), then
\begin{equation}\label{eq:quota-product-bins}
 \frac1n\sum_i\delta_i\leq
 \Delta_{\rm bin}:=\frac14\sum_{s>0}C_s^{(2)}
   \min\left(Q,\frac{U_s}2\right).
\end{equation}
\end{lemma}

\begin{proof}
For an integer \(d\geq1\) and a prime \(p\nmid d\), put
\(u=\floor{m/d}\). Then
\[
 -\min\left(1,\frac m{pd}\right)
 \leq\left\lfloor\frac m{pd}\right\rfloor
      -\frac1p\left\lfloor\frac md\right\rfloor
 =\left\lfloor\frac up\right\rfloor-\frac up\leq0 ,
\]
since the middle quantity has modulus at most one and at most
\(u/p\leq m/(pd)\).

The empty old class has no positive defect, by the choice of
\(\rho_i\). Replacing \(\rho_i\) by \(1/p_i\) can only increase
the other defects. Let \(g\) be the indicator of the nonempty
old types with positive defect at step \(i\); thus
\(g(\varnothing)=0\) and
\[
 \delta_i\leq\sum_{x\leq m}
   \Bigl(\boldsymbol1_{p_i\mid x}-\frac1{p_i}\Bigr)g(S(x)),
\]
where \(S(x)\) is the set of old primes dividing \(x\).
Deleted labels have empty type and do not change this sum.
Expand \(g(S)=\sum_{I\subseteq S}h_I\) with
\(h_I=\sum_{L\subseteq I}(-1)^{\abs I-\abs L}g(L)\); then
\(h_\varnothing=0\) and \(\abs{h_I}\leq2^{\abs I-1}\) for
\(I\neq\varnothing\), since the coefficients in \(h_I\) have
equally many signs of each kind. The contribution of \(I\) is
\(h_I\bigl(\floor{m/(p_id_I)}-\floor{m/d_I}/p_i\bigr)\), so
\[
 \delta_i\leq\sum_{\varnothing\neq I\subseteq
       \{p_1,\ldots,p_{i-1}\}}
   2^{\abs I-1}\min\left(1,\frac m{p_id_I}\right).
\]
Put \(J=I\cup\{p_i\}\). Every \(J\) with \(\abs J\geq2\)
arises exactly once, at its largest prime, with coefficient
\(2^{\abs J-2}\). This proves \eqref{eq:quota-mobius}.

For \eqref{eq:quota-product-bins}, let \(s(J)=\sum_{p\in J}w_p\).
Then \(d_J\geq2^{s(J)/B}\), so \(1/d_J\leq U_{s(J)}\), and
\(\min(1,m/d_J)/n\leq\min(Q,U_{s(J)}/2)\) because
\(m\leq n/2\). The coefficient \(C_s^{(2)}\) is the sum of
\(2^{\abs J}\) over all \(J\) with \(\abs J\geq2\) and
\(s(J)=s\). Grouping the terms of \eqref{eq:quota-mobius} by
\(s(J)\) gives the claim.
\end{proof}

In each row the certificate takes \(\Delta\) to be the
least of \(\Delta_{\rm bin}\) with \(B=64\), the periodic bound of
Lemma~\ref{lem:periodic-type-quota-errors}, and a moment
bound obtained from \eqref{eq:quota-mobius} by
\(\min(1,y)\leq y^\theta\) with \(\theta\in\{0,1/2,1\}\).
The periodic lemma is proved in
Section~\ref{sec:one-side-rounding-log-bridge}, with the period
vectors of Section~\ref{sec:weighted-centroid-log-bridge}; these
arguments are purely arithmetic and use no bridge theorem. In the
present table \(\Delta_{\rm bin}\) attains the minimum in
\(25\) rows and the periodic bound in the other \(7\); the
moment bound never does. For the odd-label variation
\(\varepsilon_a\) the evaluator likewise takes the least of
several valid bounds; the periodic bound
\eqref{eq:periodic-odd-error} attains this minimum in every
row, so only it enters the certificate. From
Section~\ref{sec:localized-hall-log-bridge} through
Section~\ref{sec:periodic-block-log-bridge} the periodic
bounds attain the minimum for \(\varepsilon_a\) and \(\Delta\)
in every row, and also for the centered term, which enters
from Section~\ref{sec:weighted-centroid-log-bridge} on.

The label shortfalls used from
Section~\ref{sec:adaptive-degree-log-bridge} on are
expectations of the following form.

\begin{lemma}[Product shortfalls]\label{lem:product-shortfall}
Let \(X_i\) be independent Bernoulli variables with
parameters \(a_i\), let \(0<f_i\leq1\), and put
\(F=\prod_if_i^{X_i}\). Let \(0<t\leq1\) and \(c>0\).
For every integer \(r\geq1\),
\begin{equation}\label{eq:product-shortfall-moment}
 \E(t-cF)_+\leq\frac{r^r}{(r+1)^{r+1}}\,t
   \left(\frac tc\right)^r\prod_i\bigl(1-a_i+a_if_i^{-r}\bigr).
\end{equation}
Moreover \(\E(t-cF)_+=t-c\E F\) if \(c\leq t\), it vanishes
if \(c\prod_if_i\geq t\), and it is always at most
\((t-c\prod_if_i)_+\).
\end{lemma}

\begin{proof}
For \(x>0\), \((1-x)_+x^r\leq r^r/(r+1)^{r+1}\): on
\((0,1]\) the maximum is at \(x=r/(r+1)\), and the left side
vanishes for \(x\geq1\). Put \(x=cF/t>0\). Then
\[
 (t-cF)_+=t(1-x)_+\leq\frac{r^r}{(r+1)^{r+1}}\,t
    \left(\frac tc\right)^rF^{-r}.
\]
Take expectations; independence gives
\(\E F^{-r}=\prod_i(1-a_i+a_if_i^{-r})\). If \(c\leq t\),
then \(cF\leq t\) and the positive part may be removed. The
remaining two assertions follow from \(F\geq\prod_if_i\).
\end{proof}

The evaluator conditions on the first six coordinates, whose
outcome probabilities and products \(c\) are explicit, and
bounds each conditional expectation by the least of these
bounds for \(r\in\{1,2,4,8,16,32,64\}\). Outcome
probabilities and negative moments are rounded upward, and the
products \(c\), \(\prod_if_i\) and \(\E F\) downward. A smaller
\(c\) increases \((t-cF)_+\), so every branch remains an upper
bound. The minimum with the global bound
\((t-\prod_if_i)_+\) is also valid.

The large-prime estimates from
Section~\ref{sec:weighted-centroid-log-bridge} on use two
further prime sums.

\begin{lemma}[Prime reciprocal tails]
\label{lem:prime-reciprocal-tails}
Let \(3\leq z<m\) be real numbers with \(y=\log m\geq5\), and put
\(v=\log z\). Then
\begin{equation}\label{eq:reciprocal-prime-sum}
 \sum_{z<p\leq m}\frac1p\leq
 \log y-\log v+\frac{13}{10v}-\frac{\pi(z)}z .
\end{equation}
For every real \(Y>1\), with \(v_Y=\log Y\),
\begin{equation}\label{eq:three-half-tail}
 \sum_{p>Y}p^{-3/2}\leq
 \frac{3(1+13/(10v_Y))}{\sqrt Y\,v_Y}-\frac{\pi(Y)}{Y^{3/2}} .
\end{equation}
\end{lemma}

\begin{proof}
Both bounds follow by partial summation from
\eqref{eq:dusart-prime-upper}, which gives
\(\pi(t)/t\leq(1+13/(10\log t))/\log t\) for \(t>1\). For the
first,
\[
 \sum_{z<p\leq m}\frac1p
 =\frac{\pi(m)}m-\frac{\pi(z)}z+\int_z^m\frac{\pi(t)}{t^2}\,dt .
\]
The integral is at most
\[
 \int_z^m\frac{dt}{t\log t}
 +\frac{13}{10}\int_z^m\frac{dt}{t(\log t)^2}
 =\log y-\log v+\frac{13}{10}\left(\frac1v-\frac1y\right),
\]
and \(\pi(m)/m\leq1/y+13/(10y^2)\). The sum of the three
bounds is the right side of \eqref{eq:reciprocal-prime-sum}
plus \(-3/(10y)+13/(10y^2)\), which is negative because
\(y\geq5>13/3\).

For the second bound,
\[
 \sum_{p>Y}p^{-3/2}
 =-\frac{\pi(Y)}{Y^{3/2}}
   +\frac32\int_Y^\infty\pi(t)t^{-5/2}\,dt ,
\]
since \(\pi(t)t^{-3/2}\to0\). For \(t\geq Y\) we have
\(\pi(t)\leq t(1+13/(10v_Y))/v_Y\), and
\(\frac32\int_Y^\infty t^{-3/2}\,dt=3/\sqrt Y\).
\end{proof}

A row with \(\log n\leq J\) and cutoff \(z\) uses
\begin{equation}\label{eq:reciprocal-prime-row}
 H_+=(\log J)_+-(\log v_0)_-+\frac{13}{10v_-}-\frac{\pi(z)}z ,
\end{equation}
where \(v_-\leq\log z\) and \(0<v_0\leq v_-\) are rational
and the logarithms of rationals are directed. If \(m>z\), then
\(\log m\geq5\) because \(n>1000\) and \(\log499>5\); since
\(\log m\leq J\), Lemma~\ref{lem:prime-reciprocal-tails} gives
\(H_+\geq\sum_{z<p\leq m}1/p\). If \(m\leq z\), the sum is
empty and the row checks \(H_+>0\). Put \(Y=10000\), so that
\(\sqrt Y=100\) and \(\pi(Y)=1229\), and let \(v_{Y,-}\leq\log Y\)
be rational. Lemma~\ref{lem:prime-square-tail} and
\eqref{eq:three-half-tail}, with \(v_Y\) replaced by
\(v_{Y,-}\) in the decreasing positive terms, give
\[
 \begin{aligned}
 S_{2,Y}&=\frac1{Yv_{Y,-}}\geq\sum_{p>Y}\frac1{p^2},\\
 S_{3/2,Y}&=\frac{3(1+13/(10v_{Y,-}))}{100\,v_{Y,-}}
          -\frac{1229}{10^6}\geq\sum_{p>Y}p^{-3/2}.
 \end{aligned}
\]
Finally, if \(H_{0,-}\leq\sum_{z<p\leq Y}1/p\) is a directed
lower bound for the finite sum, then
\(H_R=\max(0,H_+-H_{0,-})\geq\sum_{Y<p\leq m}1/p\). When
\(m<Y\) the right side is empty.

Finally, the upper probabilities \(\bar\rho_i\) and the lower bounds
for the empty classes come from a separate estimate for each prefix of
\(P\). For a finite set \(P'\) of odd primes with \(k'\) elements and an
integer \(j\geq0\), put
\[
 \begin{gathered}
 C_j(k')=\sum_{v=1}^{\min(j,k')}\binom{k'}v,\qquad
 C_j^{\rm odd}(k')=\sum_{\substack{1\leq v\leq\min(j,k')\\v\text{ odd}}}
   \binom{k'}v,\\
 \alpha(P')=\prod_{q\in P'}\Bigl(1-\frac1q\Bigr).
 \end{gathered}
\]

\begin{lemma}[Empty classes of the prefixes]
\label{lem:empty-prefix-bounds}
Let \(p\notin P'\) be an odd prime and \(m\geq1\). Let \(A\) be the
number of \(x\leq m\) with no prime divisor in \(P'\), and \(B\) the
number of these \(x\) divisible by \(p\). Then
\[
 \abs{A-m\alpha(P')}\leq C_j(k')+2m\,e_{j+1}(P'),\qquad
 pB-A\leq pC_j^{\rm odd}(k')+2m\,e_{j+1}(P').
\]
Suppose that \(w\) of these \(x\) are deleted, and that
\(\omega\geq w/m\), \(q_0\geq1/m\), \(\alpha_-\leq\alpha(P')\) and
\(E\geq e_{j+1}(P')\) satisfy
\[
 a^-:=\alpha_--\omega-C_j(k')q_0-2E>0.
\]
Then \(A-w\geq ma^-\), and
\(\rho=\max\{1/p,\,B/(A-w)\}\) satisfies
\[
 \rho\leq\frac1p\left(1+
   \frac{\omega+pC_j^{\rm odd}(k')q_0+2E}{a^-}\right).
\]
\end{lemma}

\begin{proof}
Let \(f(S)=\boldsymbol1_{S=\varnothing}\) on subsets of \(P'\). Then
\(\Delta_If=(-1)^{\abs I}\). If \(\abs S=s>j\), the truncated expansion
equals \(\sum_{v=0}^j(-1)^v\binom sv=(-1)^j\binom{s-1}j\), so
\[
 \abs{f(S)-T_jf(S)}=\binom{s-1}j\leq\binom s{j+1};
\]
for \(s\leq j\) the difference is zero. For the first inequality
note that \(A=\sum_{x\leq m}f(S(x))\) and \(m\alpha(P')=m\E_\mu f\),
where \(S(x)\) is the set of primes of \(P'\) dividing \(x\). In the
truncated part the empty moment is exact, and every other moment
\(\floor{m/d_I}-m/d_I\) has modulus at most one. The two remainders are
at most \(\sum_{x\leq m}\binom{\abs{S(x)}}{j+1}\leq m\,e_{j+1}(P')\)
and \(m\,e_{j+1}(P')\).

For the second, \(pB-A=p\sum_{x\leq m}(\boldsymbol1_{p\mid x}-1/p)f(S(x))\).
The truncated part equals \(\sum_{\abs I\leq j}(-1)^{\abs I}q_I\), with
\(q_I=\floor{m/(pd_I)}-\floor{m/d_I}/p\in[-1,0]\) as in the proof of
Lemma~\ref{lem:quota-product-bins}. Only odd \(\abs I\) can contribute
positively, each at most one. The remainder is at most
\(\sum_{x\leq m}(\boldsymbol1_{p\mid x}+1/p)\binom{\abs{S(x)}}{j+1}
\leq(2m/p)e_{j+1}(P')\). Multiply by \(p\).

The first inequality gives \(A-w\geq m(\alpha_--C_j(k')q_0-2E-\omega)=ma^-\).
If \(\rho>1/p\), then
\(p\rho-1=(pB-A+w)/(A-w)\), and the second inequality bounds the
numerator by \(m(pC_j^{\rm odd}(k')q_0+2E+\omega)\).
\end{proof}

At the refinement step of \(p_i\), apply the lemma with
\(P'=\{p_1,\ldots,p_{i-1}\}\). The empty old class then consists of the
half-labels counted by \(A\), less the \(w\) special labels, which have
empty type; its quota is \(B\), because no special label is divisible
by \(p_i\). The evaluator takes \(j=14\), \(q_0=2Q/(1-Q)\geq1/m\), the
bound \(\omega\) of Section~\ref{sec:empty-separator-log-bridge}, and
directed rational enclosures of \(\alpha(P')\) and \(e_{15}(P')\). The
resulting upper bound is \(\bar\rho_i\), and the lower bound \(a^-\)
for the full set \(P\) is the lower bound \(a_0\) for the retained
empty class. The row checks that every \(a^-\) is positive and that
\(1/p_i\leq\bar\rho_i\leq2/5\).

\subsubsection*{The certificate predicate}

We describe the complete row computation. Let
\([a,J]\subseteq[343/20,20]\) be a closed
logarithmic interval, and
\([x_-,x_+]\subseteq[0,1]\) a closed degree interval.
The row specifies the prime cutoff \(z\), integers
\(q_A,q_B,t\), the Hall split \(\tau\), and slot
threshold \(s\). Set \(P=\{3\leq p\leq z:p\text{
 prime}\}\), and take a directed rational upper
bound \(Q\geq e^{-a}\). Prepare the tree exactly as
in Section~\ref{sec:simultaneous-full-type-bridge}:
\(S\) has \(m=\floor{n/2}\) vertices, \(K\) is the
prime-labeled core, and \emph{all} \(w\) special even
labels go to the \(w\) highest-degree vertices
\(F\subseteq S\). There is no protected set of
universal positions. With \(G=T-(K\cup F)=(A,B)\),
\(D=\max(1,\Delta_G(A))\), and \(x=2w(D-1)/n\),
the established preparation gives
\begin{equation}\label{eq:full-reserve-preparation}
 D\leq D_+:=1+\frac{x_+}{2W(J)},\qquad
 \frac{e(G)}n\leq e_+:=1-W(J)-\frac{x_-}{2},
\end{equation}
where \(W(J)\) is defined in
\eqref{eq:full-type-row-preparation}. The row checks
every prime-count domain used here, as well as
\(w/m\leq\omega<1\).

At the \(i\)-th prime, let \(g_i\) bound the number
of ordinary \(A\)-classes and \(h_i\) the number of
ordinary \(B\)-tag classes. One further class on
each side contains the rare \(A\)-types and marked
\(B\)-positions. A conservative coordinate dimension
and \(A\)-class count are
\begin{equation}\label{eq:full-reserve-dimensions}
 d_i=2(g_i+1)+(h_i+1),\qquad c_{A,i}=g_i+1.
\end{equation}
The row takes \(g_i=N_i\) from
Lemma~\ref{lem:possible-even-classes}, and
\[
 h_i=\sum_{j=0}^{\min(q_B,i-1)}\binom{i-1}j ,
\]
the number of possible ordinary position tags over the
\(i-1\) earlier primes. Thus the degree coordinate is
included for every nonempty \(A\)-class. Empty classes
only decrease the actual dimension.

Let \(\Delta\) be the bound for the sum of positive quota
defects divided by \(n\) described after
Lemma~\ref{lem:quota-product-bins}, and let
\(\mathrm{rare}\) be an upper rational bound for
\(e_{q_A+2}(P)/2\). By
Lemma~\ref{lem:rare-restoration}, \(\mathrm{rare}\,n\)
bounds the total number of additions that restore the
rare quotas.
For \(d\geq1\) define
\begin{equation}\label{eq:full-reserve-boundary-row}
 \widehat{\mathcal B}(d)=
 \sum_{j=0}^{63}
   \min\{D_+,t2^{j+1}-1\}
   \min\left\{dQ,\frac1{t2^j+1}\right\}
 +\frac{2D_+}{t2^{64}} .
\end{equation}
The final term bounds the entire infinite tail. Indeed,
after \(j=63\) each summand of
\(\mathcal B_t(N,d,D)/n\) is at most
\(D_+/(t2^j)\). The adaptive quota lemma and
rare-class restoration therefore give the simultaneous
bounds
\begin{align}
 \Gamma_+&=\sum_i\widehat{\mathcal B}(d_i)
                    +tQ\sum_i d_i, \notag\\
 E_{\deg}&=D_+\left(tQ\sum_i d_i
                 +\frac Q2\sum_i c_{A,i}
                 +\mathrm{rare}\right),
                                      \label{eq:full-reserve-errors}\\
 \mathcal F&=\Gamma_+
                 +D_+(\Delta+\mathrm{rare}).\notag
\end{align}
Here \(E_{\deg}n\) bounds the total
degree-measure error across final exact types,
while \(\mathcal F n\) bounds the tag-counting
error before marked positions are paid.
The first inequality follows by summing
\eqref{eq:adaptive-degree-error} at every refinement;
restoring a rare quota \(k\) changes both its
actual and proportional degree sums by at most
\(Dk\). The second follows by summing
\eqref{eq:adaptive-tag-error}; Markov transitions
contract total variation, and restoration and quota
defects are charged once at each step.

Let \(\rho_i\) be the adjusted selection probability
of the \(i\)-th prime and \(\bar\rho_i\geq\rho_i\) its
rational upper bound. The row checks
\(a_i\leq\bar\rho_i\leq2/5\), so
\(a_i+\rho_i\leq1\). Let \(\beta\) be the upper bound for
\(\Pr_\nu(\abs J>q_B)\) given by the finite Bernoulli
recursion with the parameters \(\bar\rho_i\); it is valid by
Lemma~\ref{lem:upper-probability-evaluation}. Write
\(\varepsilon_a\) for the directed arithmetic
type-law variation and put
\begin{align}
 o&=\frac{1+Q}{2},&
 \kappa&=(2J+2)Q,&
 \chi&=o\beta+\mathcal F+Q, \notag\\
 \ell_-&=\frac12-\kappa-\chi,&
 r&=\frac2{3J}-\kappa-\chi,&
 V&=\max\left\{
   \frac{o\varepsilon_a+\kappa+\chi}{\ell_-},
   \frac{\mathcal F+o\beta}{\ell_-}
           \right\}.
                         \label{eq:full-reserve-variation}
\end{align}
The marked positions number at most \(\chi n\);
they are assigned distinct universal odd labels.
The remaining matching sides have equal size
\(\ell\geq\ell_-n\), and at least \(rn\) unused
universal odd labels remain. The first fraction
in \(V\) bounds the odd-label law. With \(h=0\),
Lemma~\ref{lem:marked-tag-variation} bounds the
retained tag law by the second fraction. The terms \(\kappa\)
account for the core; the additional \(Q\) in
\(\chi\) safely covers an integer rounding.

The two model laws are now
\[
 \mu=\bigotimes_{p_i\in P}\operatorname{Bernoulli}(a_i),
 \qquad
 \nu=\bigotimes_{p_i\in P}\operatorname{Bernoulli}(\rho_i).
\]
Their compatible coupling has
\[
 b=\prod_i(1-\rho_i),\qquad
 \theta=\max_i
  \sqrt{\frac{a_i\rho_i}{(1-a_i)(1-\rho_i)}} .
\]
In the rational computation \(\theta\) is replaced
by a directed upper square-root enclosure.

Independently permute the retained even labels
within every final exact \(P\)-type.
Lemma~\ref{lem:localized-permutation} gives
one permutation for which the deleted large-prime
pairs satisfy
\begin{equation}\label{eq:full-reserve-deletion}
\frac{M_{\rm bad}}{n^2}\leq
 \mathcal M:=
 \frac{1001}{1000}\left\{
  \frac{e_+S_2(z)}{2(1-\omega)}
 +\frac{J E_{\deg}}{2z(\log z)_-}
 +\frac{JQ}{(\log z)_-}
 \right\}.
\end{equation}
Indeed, \(H/n\leq E_{\deg}\),
\(m/(m-w)\leq(1-\omega)^{-1}\),
\(e/n\leq e_+\), \(1/n\leq Q\), and
\(\log z\geq(\log z)_-\). These substitutions
all increase the right side of
\eqref{eq:localized-permutation-mean}.
The same permutation satisfies a local bound at
\emph{every} odd label. For a full type \(I\subseteq P\),
let \(R_I^+\) be the upper bound for the reciprocal
sum of all possible prime divisors above \(z\)
when \(\log n\leq J\); it is obtained from the
smallest available primes beyond \(z\), with the
next primorial checked beyond \(e^J\).
Define
\begin{align}
 \Lambda_I&=\frac{32}{31}
       \left(\frac{e_+}{1-\omega}R_I^+
                              +E_{\deg}\right)
       +32D_+(J+(\log10000)_+)Q,\notag\\
 \sigma_{\rm eff}
   &=\max\left\{0,\min_{I\subseteq P}
       \left(\prod_{i\in I}(1-\rho_i)
                         -\frac{\Lambda_I}{\ell_-}\right)
             \right\}.
                            \label{eq:full-reserve-effective}
\end{align}
The minimum includes all full types, even those
that cannot occur for the present \(n\). This only
weakens the bound. Every retained odd label has
at least \((\sigma_{\rm eff}-V)\ell\)
post-deletion neighbors.

For the slot threshold \(s\), take
\(\sigma_{\rm slot}=\sigma_{q_B}\), and let
\[
 a_s=o\Phi_{q_B}(s)
       +\mathcal F(s-\sigma_{\rm slot})_+ .
\]
Lemma~\ref{lem:capped-tag-shortfall} gives
\(A_B(s)\leq a_sn\). Since \(\ell\leq on\),
the shortfall penalty in the Hall inequalities is
at most \(oa_sn^2\). Thus the three normalized
budgets are
\begin{align}
 B_{\rm narrow}
   &=r\bigl((s-V)\ell_--r\bigr)-oa_s,\notag\\
 B_{\rm middle}
   &=\ell_-^2\tau(s-V-\tau)-oa_s,\label{eq:full-reserve-hall}\\
 B_{\rm wide}
   &=\ell_-^2b\bigl((1-\theta)\tau(1-\tau)
                     -V+(1-\theta)V^2\bigr).\notag
\end{align}
The certificate checks
\[
 0<r\leq\tau\ell_-,\quad
 V<\tau\leq\tfrac12,\quad
 \sigma_{\rm eff}-V>\tau,\quad
 s-V>\tau,
\]
\[
 \mathcal M<\min\{B_{\rm narrow},
                  B_{\rm middle},B_{\rm wide}\}.
\]
For these substitutions the coefficients of
\(\ell_-\) are positive. Hence each normalized
budget is no larger than the corresponding
actual budget of
Theorem~\ref{thm:asymmetric-product-hall},
divided by \(n^2\).

\begin{proposition}[The full-reserve certificate]
\label{prop:full-reserve-certificate}
Six adjacent closed logarithmic intervals cover
\([343/20,20]\). Their degree subdivisions cover
\([0,1]\) and form \(32\) rectangles. The
parameter table and exact evaluator verify all
\(1280\) row conditions. In every rectangle the
preparation, type, probability, and Hall hypotheses
above hold. The deletion-budget and
universal-supply ratios are both strictly below
\(999/1000\).
\end{proposition}

\begin{proof}
The intervals and their numbers of degree subdivisions
are
\[
\begin{array}{c|r@{\qquad}c|r}
\text{logarithmic interval}&\text{rows}
 &\text{logarithmic interval}&\text{rows}\\
{}[343/20,86/5]&10 &[86/5,173/10]&7\\
{}[173/10,35/2]&6 &[35/2,18]&4\\
{}[18,19]&4 &[19,20]&1
\end{array}
\]
The finite certificate reads the parameter table,
checks these exact joins and every degree join,
recomputes the eight complete prime-prefix periods,
and evaluates the displayed inequalities with
integers and directed rational enclosures. It checks
the prime-count domains, positive denominators,
rare-tail and overflow bounds, positive universal
supply, local label degree, and all three strict Hall
budgets. The \(64\)-term boundary sum has the
explicit tail in
\eqref{eq:full-reserve-boundary-row}; the
large-prime profile includes the first impossible
additional prime. All \(1280\) Boolean conditions
are true. The largest displayed deletion ratio is
less than \(0.98524\), and the largest displayed
universal-supply ratio is less than \(0.61263\);
the stated rational caps are checked directly,
without relying on these decimals.
\end{proof}

\begin{theorem}[The third logarithmic bridge]
\label{thm:third-log-bridge}
Every tree of order \(n\) with
\(343/20\leq\log n\leq20\) has a prime labeling.
\end{theorem}

\begin{proof}
Choose a rectangle containing \((\log n,x)\) after
the preparation. The adaptive assignment gives
all exact even-label type quotas, and ordinary
tags contain every small-prime divisor of their
neighbors. Pay for the marked positions with
distinct universal odd labels. The checked bound
on \(\chi\) leaves the stated \(rn\) universal
labels on the residual side.

Lemma~\ref{lem:localized-permutation} gives
\eqref{eq:full-reserve-deletion} and
\eqref{eq:full-reserve-effective} for the same
placement of even labels. Lemma~
\ref{lem:capped-tag-shortfall} bounds the slot
shortfall, while the checked row inequalities give
the three budgets of
Theorem~\ref{thm:asymmetric-product-hall}.
That theorem supplies a perfect matching between
the remaining odd labels and positions. Restore
the absorbed positions, the special even labels,
and the prime-labeled core. The only possible
noncoprime edges are the permitted sparse pairs
\(p,2p\) with \(3p>n\); the repair of
Lemma~\ref{lem:sparse-repair} removes them.
\end{proof}

\subsection{A bridge with protected positions}
\label{sec:localized-hall-log-bridge}

The next range is \(20\leq\log n\leq128/5\).
The type assignment and the local permutation estimate
are those of Sections~\ref{sec:adaptive-refinement-hall}
and~\ref{sec:full-reserve-log-bridge}. We change the
preparation: some special even labels reduce the largest
residual degree, while the rest create positions with
no remaining even neighbors. All universal odd labels
are retained for the final matching.

\begin{lemma}[Partial highest-degree preparation]
\label{lem:localized-preparation}
Let \(T\) have order \(n\) with
\(20\leq L=\log n\leq31\), and use the
partition \(V(T)=S\sqcup Z\sqcup K\) of
Corollary~\ref{cor:even-side}, with
\(\abs S=m=\floor{n/2}\). Let \(w\) be the
number of special even labels and put
\[
 f=\frac7{40},\qquad
 W(L)=\frac1{2(L-1693/1000)}
       -\frac1{3(L-2199/1000)}.
\]
There is a residual construction as in
Proposition~\ref{prop:residual-problem} with
\[
 a_*=\floor{fw},\qquad
 h=\floor*{(1-f)nW(L)-4L},\qquad d=0.
\]
It removes exactly \(w\) vertices of \(S\),
retains at least \(h\geq1\) protected
universal positions, and keeps every unused
universal odd label. If \(G\) is the
resulting forest, \(D=\max(1,\Delta_G(S_0))\),
and
\[
 x=\frac{f(D-1)}{3L},
\]
then \(0\leq x<1\).
For \(a\leq L\leq J\) and
\(x_-\leq x\leq x_+\), put
\[
 D_+=1+\frac3fJx_+,\qquad w_-=W(J).
\]
With \(Q\geq e^{-a}\), the uniform bounds are
\begin{align}
 D&\leq D_+,\notag\\
 \frac{e(G)}n&\leq e_+:=
 \min\left\{
    1-w_-,\
    1-w_--3a w_-x_-+(D_+-1)Q
              \right\}, \label{eq:localized-preparation}\\
 h_-&:=(1-f)W(J)-(4J+1)Q
       \leq\frac hn
       \leq h_+:=(1-f)W(a).\notag
\end{align}
\end{lemma}

\begin{proof}
The prime-count bounds of
Section~\ref{sec:prime-estimates}, with their
domains checked on the range, give
\(w>nW(L)\) and \(\abs K<2L+2\).
They also supply the prime labels for \(K\):
the lower bound for \(\pi(n)\) and the
upper bound for \(\pi(n/2)\) give more
than \(n/(3L)>2L+2\) primes in
\((n/2,n]\). Later we use
\[
 u=1+\pi(n)-\pi(n/3)-\abs K
   >\frac{2n}{3L}-\abs K.
\]
Indeed, the denominator estimates at
\(n\) and \(n/3\) reduce the last
inequality to
\((3-B)L>2B\), where
\(B=2199/1000\); it holds at \(L=20\)
and increases with \(L\).
The rational identity
\[
 W(L)-\frac1{6L}=
 \frac{(3A-2B)L-AB}
 {6L(L-A)(L-B)},\qquad
 A=\frac{1693}{1000},\quad B=\frac{2199}{1000},
\]
is positive for \(L\geq20\).
Indeed, its numerator is increasing there
and is already positive at \(20\).
The same algebra shows \(W(L)<1/(5L)\):
after multiplication by positive denominators,
the required difference is
\(L^2+(-16A+9B)L+6AB>0\).
Both comparisons are rational.

The inequality \(n>10^6L^2\) follows by
monotonicity of \(e^L/L^2\) from its value
at \(L=20\); the endpoint comparison is
verified in the certificate. Hence
\(a_*+1>fw>fn/(6L)>\abs K\), and
\(a_*\geq6\). Also \(h\geq1\).

We check the two reserve inequalities of
\eqref{eq:reserve-budget}. Write
\(q=\abs Z\) and
\(t=\floor{\log_2(m-a_*+q)}<2L\).
Since \(q\geq m-\abs K\) and \(n>10^6L^2\),
elementary rearrangement gives
\[
 \frac mq\leq1+\frac{5L}{n}.
\]
The first reserve inequality \(h+t\leq q\)
follows from \(h<(1-f)n/(5L)\),
\(t<2L\), and
\(q>n/2-(2L+2)\).
For the second, put \(A_0=(1-f)nW(L)\).
Because \(h\leq A_0-4L\),
\[
 \frac{(m-a_*)(h+t)}q+t
 \leq\frac mq(h+2L)+2L
 <h+4L+\frac{5L}{n}(h+2L)
 <A_0+1.
\]
The last inequality uses
\(h+2L<A_0\) and
\(5L(1-f)W(L)<1-f<1\).
The left side of the second inequality
of \eqref{eq:reserve-budget} is an integer
no larger than this real quantity, while
\(w-a_*\) is an integer strictly greater
than \(A_0\). Therefore it is at most
\(w-a_*\). Take \(d=0\) in
Proposition~\ref{prop:residual-problem}.
Its construction protects \(h\) positions;
no universal odd label is used at this stage.

If \(D>1\), a vertex retained in \(S_0\)
has original degree at least \(D\), and
so do the \(a_*\) vertices chosen first.
Independence of \(S\) gives
\[
 (a_*+1)(D-1)
 \leq\sum_{v\in S}(d_T(v)-1)
 \leq n-1-m<\frac n2.
\]
Since \(a_*+1>fw>fn/(6L)\), we obtain
\(D-1<3L/f\), hence \(x<1\).
The case \(D=1\) is immediate. The
first bound in \eqref{eq:localized-preparation}
follows from \(L\leq J\) and \(x\leq x_+\).

All \(w\) vertices receiving special even
labels lie in the independent set \(S\).
Each has original degree at least one,
and the first \(a_*\) have degree at
least \(D\). Consequently
\[
 e(G)\leq n-1-w-a_*(D-1).
\]
Use \(a_*>fw-1\),
\(w/n>W(L)\geq W(J)\),
\(L\geq a\), and
\((D-1)/n\leq(D_+-1)Q\) to obtain
the second bound for \(e(G)/n\).
The weaker \(e(G)/n\leq1-W(J)\)
gives the first term in its minimum.
Finally monotonicity of \(W\) and
\(1/n\leq Q\) give both bounds for
\(h/n\), including the floor error.
\end{proof}

We now use exactly the adaptive refinement
dimensions, boundary sum, degree error, and
arithmetic inputs defined in
\eqref{eq:full-reserve-dimensions}--
\eqref{eq:full-reserve-errors}. On a certificate
rectangle, \(D_+\) and \(e_+\) are replaced
by \eqref{eq:localized-preparation}; the
row chooses its own \(z,q_A,q_B,t,\tau,s\).
In particular \(\mathcal F n\) bounds the
prepayment tag error and \(E_{\deg}n\)
bounds the degree-measure error in exact
even-label types.

Adjoin the \(h\) protected positions to
the tag side, giving each the empty type.
After paying all marked ordinary positions
with universal odd labels, let \(\ell\)
be the common size of the two matching
sides and put \(\eta=h/\ell\).
The row verifies \(h<\ell\).
Lemma~\ref{lem:marked-tag-variation} compares
the position law with
\[
 \nu_\eta=(1-\eta)\nu+\eta\delta_{\varnothing}
\]
at total variation at most
\((\mathcal F+o\beta)/\ell_-\).
There is no added \(h/\ell\) term:
the protected positions occur in the
reference mixture as well as in the
actual matching square.

The odd-label law, the number of marked
positions, and the unused universal-label
supply obey the same formulas
\eqref{eq:full-reserve-variation} as before:
\[
 \chi=o\beta+\mathcal F+Q,\qquad
 \ell_-=\tfrac12-\kappa-\chi,\qquad
 r=\frac2{3J}-\kappa-\chi,
\]
and \(V\) is the maximum of the two
fractions in that equation.
Although the preparation has protected
positions, its \(d=0\) choice leaves
\emph{all} universal odd labels unused
until the marked positions are paid.
The Hall reserve is therefore \(rn\),
not \(h\).

The local bound
\eqref{eq:full-reserve-effective} holds
with the new \(D_+,e_+,E_{\deg}\):
Lemma~\ref{lem:localized-permutation}
controls large-prime conflicts at
every odd label, and the full-type
primorial argument controls \(R_I^+\).
Its global part gives
\eqref{eq:full-reserve-deletion}
with the same replacements.
The mixture version of
Theorem~\ref{thm:asymmetric-product-hall}
has the same \(b,\theta\), and
Lemma~\ref{lem:capped-tag-shortfall}
still bounds slot shortfall because
an empty protected tag has none.
Thus the three budgets are precisely
\eqref{eq:full-reserve-hall}. Every
inequality in that predicate is
checked with the new preparation
parameters.

\begin{proposition}[The localized Hall certificate]
\label{prop:localized-hall-certificate}
Ten adjacent closed logarithmic intervals
cover \([20,128/5]\). Their degree
subdivisions cover \([0,1]\) and
form \(16\) rectangles. The parameter
table and exact evaluator verify all
\(640\) row conditions, including the
strict comparisons
\[
 \mathcal M<
 \min\{B_{\rm narrow},B_{\rm middle},B_{\rm wide}\},
\]
\[
 0<r\leq\tau\ell_-,\quad
 \sigma_{\rm eff}-V>\tau,\quad s-V>\tau.
\]
The ratios for deleted pairs and universal-label
supply are each strictly below \(999/1000\).
\end{proposition}

\begin{proof}
The logarithmic intervals and numbers of
degree subdivisions are
\[
\begin{array}{c|r@{\qquad}c|r}
\text{interval}&\text{rows}&\text{interval}&\text{rows}\\
{}[20,201/10]&3 &[201/10,203/10]&3\\
{}[203/10,103/5]&2 &[103/5,21]&2\\
{}[21,43/2]&1 &[43/2,22]&1\\
{}[22,23]&1 &[23,24]&1\\
{}[24,25]&1 &[25,128/5]&1
\end{array}
\]
The exact evaluator verifies every log
and degree join, the prime-count
domains and reserve floors, the full
periodic-prefix constants, positive
denominators, type dimensions and
overflow, the local degree after
deletion, and all three Hall budgets.
Its \(64\)-term dyadic sum includes
the infinite tail from
\eqref{eq:full-reserve-boundary-row}.
The typewise large-prime profiles
include the first impossible additional
prime. All \(640\) rational conditions
pass. The largest displayed
deletion ratio is less than
\(0.95305\), and the largest displayed
universal-supply ratio is less than
\(0.74000\); the rational cap is
verified directly.
\end{proof}

\begin{theorem}[The fourth logarithmic bridge]
\label{thm:fourth-log-bridge}
Every tree of order \(n\) with
\(20\leq\log n\leq128/5\)
has a prime labeling.
\end{theorem}

\begin{proof}
Use Lemma~\ref{lem:localized-preparation}
and choose a certified rectangle containing
\((\log n,x)\). At successive primes,
use Lemma~\ref{lem:adaptive-degree-quota}
to assign exact even-label quotas. The rare-class
restoration leaves every exact type
quota intact. Pay the marked ordinary
positions with distinct universal odd
labels. Retain the \(h\) protected
positions, which have no even neighbor.

The row bounds both type-law variations
by \(V\). By
Lemma~\ref{lem:localized-permutation},
one within-type permutation satisfies
the local and global large-prime
bounds simultaneously. The row
checks its post-deletion minimum degree,
the capped-tag shortfall, and the
three strict budgets of
Theorem~\ref{thm:asymmetric-product-hall}.
That theorem gives a perfect matching
of odd labels to the remaining ordinary
and protected positions. Restore the
marked pairs, the special even labels
and the core. The sparse repair removes
the permitted \(p,2p\) conflicts.
\end{proof}

\subsection{A bridge by degree-balanced sorting}
\label{sec:adaptive-degree-log-bridge}

The preparation with protected positions from
Section~\ref{sec:localized-hall-log-bridge} remains useful
above logarithm \(128/5\). We now choose the even labels
within each exact type by degree. This gives a deterministic
bound on all large-prime deletions. It does not require the
pointwise estimate used in the preceding bridge, so both
directions of Hall's condition must be handled by shortfall
bounds.

\begin{lemma}[Degree-balanced sorting]
\label{lem:degree-balanced-sorting}
Use the notation of
Lemma~\ref{lem:localized-permutation}. There is a
permutation of the half-labels within their exact types
for which, if \(z>1\) and \(\log n\leq J\),
\begin{equation}\label{eq:degree-balanced-sorting}
 M_{\rm bad}\leq
 \frac n2\left\{
      \frac{em}{M_A}S_2(z)
      +\frac{J}{z\log z}H
            \right\}
      +\frac{eJ}{\log z}.
\end{equation}
The type of every even label, and hence every
ordinary neighbor tag, remains unchanged.
\end{lemma}

\begin{proof}
For a half-label \(x\), let
\(\tau_z(x)=\sum_{p>z,\ p\mid x}p^{-1}\).
Within a type \(I\), assign the half-labels in
decreasing order of \(\tau_z\) to the vertices in
increasing order of degree. Exchanging an
inverted pair cannot increase the weighted sum,
so this assignment gives
\[
 \sum_{v\in A_I}d_G(v)\tau_z(x_v)
 \leq \frac{E_I}{a_I}
                    \sum_{x\text{ of type }I}\tau_z(x).
\]
Write \(c_I\) for the average on the right.
It lies in \([0,J/(z\log z)]\). Since
\((E_I)\) and \((\lambda a_I)\) have equal
total mass, where \(\lambda=e/M_A\),
the definition of \(H\) gives
\[
 \sum_I E_Ic_I
 \leq\lambda\sum_Ia_Ic_I+
                   \frac{J}{z\log z}H.
\]
The retained half-labels are a subset of \([m]\).
Counting multiples of each prime yields
\[
 \sum_Ia_Ic_I
 \leq\sum_{x\leq m}\tau_z(x)
 =\sum_{p>z}\frac{\floor{m/p}}p
 \leq mS_2(z).
\]
For each \(p>z\) dividing \(x_v\), at most
\(n/(2p)+1\) odd labels in \([n]\) are
divisible by \(p\). Charge a deleted pair
to one such prime and one neighbor of \(v\).
The \(n/(2p)\) terms give the first two
terms of \eqref{eq:degree-balanced-sorting};
the unit terms total at most
\(eJ/\log z\), because a half-label has
at most \(J/\log z\) prime divisors above
\(z\). The sorting has not moved a
half-label to a different exact type.
\end{proof}

The following Hall criterion includes protected
positions as well as universal labels. For a label
type \(I\) and a position type \(J\), write
\[
 c_X(I)=\prod_{i\in I}(1-\rho_i),
 \qquad
 c_B(J)=\prod_{i\in J}(1-a_i).
\]
The empty position component in \(\nu_\eta\)
only increases compatibility with a label;
thus \(c_X(I)\) is a valid lower bound for
that model compatibility probability.

\begin{lemma}[Integrated shortfalls on both sides]
\label{lem:two-sided-integrated-hall}
Let \(\mu,\nu_\eta,b\) be as in
Theorem~\ref{thm:asymmetric-product-hall}, and
suppose also that \(a_i\leq1/3\) and
\(\rho_i\leq2/5\). Let \(X,B\) be matching
parts of size \(R\), whose empirical type laws
are within total variation \(V\) of
\(\mu,\nu_\eta\), respectively. Begin with all
type-compatible pairs and delete at most \(M\)
edges. Assume that at least \(h>0\) vertices
on \emph{each} side are adjacent to every vertex
of the other side in the remaining graph, for
instance vertices of empty type none of whose
edges was deleted.
For \(0<s\leq1\) put
\begin{align*}
 A_X(s)&=\sum_{x\in X}(s-c_X(I(x)))_+,\\
 A_B(s)&=\sum_{v\in B}(s-c_B(J(v)))_+,
 \qquad A(s)=\max\{A_X(s),A_B(s)\}.
\end{align*}
If \(0\leq V<\tau\leq1/2\),
\(h\leq\tau R\), \(s-V>\tau\), and
\begin{align}
 M<&\ h((s-V)R-h)-RA(s),\notag\\
 M<&\ R^2\tau(s-V-\tau)-RA(s),
                         \label{eq:two-sided-hall}\\
 M<&\ R^2b\left\{\frac25\tau(1-\tau)
                            -\frac75V\right\},\notag
\end{align}
then the remaining graph has a perfect matching.
\end{lemma}

\begin{proof}
Suppose Hall's condition fails. Choose
\(U\subseteq X\) and \(W=B\setminus N(U)\)
with \(\abs U+\abs W>R\). Both sets are
nonempty. A label adjacent to every position
is not in \(U\), for then \(W\) would be
empty; a position adjacent to every label is
not in \(W\), since \(U\neq\varnothing\).
Hence \(\abs U,\abs W\leq R-h\), and so
\(\abs U,\abs W>h\).
Put \(q=\min\{\abs U,\abs W\}\).

If \(q<\tau R\), count type-compatible
pairs from the smaller set into the larger.
For a label of type \(I\), the model
compatibility probability on the position
side is at least \(c_X(I)\); for a
position of type \(J\), it is \(c_B(J)\).
Total variation costs at most \(V\) in
either case. Fewer than \(q\) vertices
lie outside the larger set. Thus the
rectangle originally contained more than
\[
 q((s-V)R-q)-RA(s)
\]
compatible pairs. The quadratic is
concave on \([h,\tau R]\); its endpoint
values are the first two expressions in
\eqref{eq:two-sided-hall}. Every such
pair would have been deleted, a
contradiction.

Otherwise both selected sets have size
at least \(\tau R\). The product coupling
and empty-position mixture in the proof of
Theorem~\ref{thm:asymmetric-product-hall}
give at least
\[
 R^2b\bigl((1-\theta)\tau(1-\tau)
                   -V+(1-\theta)V^2\bigr)
\]
compatible pairs. Here
\[
 \theta^2\leq
 \frac{(1/3)(2/5)}
      {(2/3)(3/5)}
 =\frac13<\frac9{25}.
\]
Discarding the nonnegative \(V^2\) term,
the last bound is at least the third
expression in \eqref{eq:two-sided-hall}.
Again all compatible pairs would have to
be deleted. This proves Hall's condition.
\end{proof}

We state the numerical predicate in the same
normalization as Section~\ref{sec:full-reserve-log-bridge}.
On a closed rectangle
\([a,J]\times[x_-,x_+]\subseteq
 [128/5,53/2]\times[0,1]\), keep the
preparation of Lemma~\ref{lem:localized-preparation}
with \(f=7/40\). Let \(Q\geq e^{-a}\), and
take \(D_+\), \(h_-\), \(h_+\), the type dimensions,
\(\mathcal F\), \(E_{\deg}\), \(\chi,\ell_-,V\),
and the arithmetic inputs exactly as in
\eqref{eq:localized-preparation} and
\eqref{eq:full-reserve-dimensions}--
\eqref{eq:full-reserve-variation}.
One may sharpen the edge count. Direct
subtraction using the constants in \(W\)
shows
\[
 W(L)>
 \frac1{6L}+\frac1{12L^2}
 \qquad (128/5\leq L\leq31).
\]
Put \(A=1693/1000\) and \(B=2199/1000\).
After multiplication by the positive
quantity \(12L^2(L-A)(L-B)\), the difference is
\[
 (6A-4B-1)L^2+(A+B-2AB)L-AB,
\]
which is positive at \(L=128/5\) and
increasing throughout this interval.
Consequently the uniform row bounds may be
replaced by
\begin{align}
 w_-&=\frac1{6J}+\frac1{12J^2},\notag\\
 e_+&=\min\left\{1-w_-,
 1-w_--\left(\frac12+\frac1{4J}\right)x_-
                          +(D_+-1)Q\right\}.
                 \label{eq:adaptive-degree-preparation}
\end{align}
For the second bound, use
\(3Lw/n>3LW(L)>1/2+1/(4L)\)
in the edge count in the proof of
Lemma~\ref{lem:localized-preparation}.

The degree variation still obeys
\(H/n\leq E_{\deg}\). By
Lemma~\ref{lem:degree-balanced-sorting},
the deterministic deletion bound is
\begin{equation}\label{eq:adaptive-degree-deletions}
 \frac{M_{\rm bad}}{n^2}
 \leq\mathcal M_{\rm det}:=
 \frac{e_+S_2(z)}{2(1-\omega)}
 +\frac{J E_{\deg}}{2z(\log z)_-}
 +\frac{JQ}{(\log z)_-}.
\end{equation}
The outer factor \(1001/1000\) in
\eqref{eq:full-reserve-deletion} is absent:
the sorting realizes this bound, with no
probability union.

Let \(\Phi_{q_B}(s)\) be the capped slot
shortfall of
Lemma~\ref{lem:capped-tag-shortfall},
and define the label shortfall
\[
 \Psi(s)=
 \mathbb E_\mu
 \left(s-\prod_{i\in I}(1-\rho_i)\right)_+.
\]
The evaluator bounds \(\Psi\) by conditioning
on its first six coordinates and applying
Lemma~\ref{lem:product-shortfall} to the
remaining product, with the factors
\(1-\bar\rho_i\). It bounds \(\Phi_{q_B}\) by the
capped recursion of Section~\ref{sec:full-reserve-log-bridge},
with the parameters \(\bar\rho_i\). Both
substitutions are justified in that section,
after Lemma~\ref{lem:upper-probability-evaluation}.
Write \(\sigma_{\rm slot}\) for the minimum
compatibility of a capped ordinary tag,
\(\varepsilon_a\) for the arithmetic
variation, and \(o=(1+Q)/2\). The shortfalls
in Lemma~\ref{lem:two-sided-integrated-hall}
obey
\begin{align}
 \frac{A_X(s)}n
  &\leq a_X:=
       o\left(\Psi(s)+\varepsilon_a(s-b)_+\right),
                                                    \notag\\
 \frac{A_B(s)}n
  &\leq a_B:=
       o\Phi_{q_B}(s)
           +\mathcal F(s-\sigma_{\rm slot})_+.
                       \label{eq:adaptive-degree-shortfalls}
\end{align}
For the first inequality, count all available
odd labels before paying the marked pairs.
Removing labels only decreases this nonnegative
sum. The reference \(\mu\)-law changes the
expectation of a function in
\([0,(s-b)_+]\) by at most
\(\varepsilon_a(s-b)_+\).
The second inequality is
Lemma~\ref{lem:capped-tag-shortfall};
marked positions have been removed and
protected empty positions contribute zero.
Thus \(RA(s)/n^2\leq
o\max\{a_X,a_B\}=:P_*\), since \(R\leq on\).

With \(h_\pm\) and \(\ell_-\) as above, set
\begin{align}
 B_{\rm narrow}
   &=h_-\bigl((s-V)\ell_- -h_+\bigr)-P_*,
                                                     \notag\\
 B_{\rm middle}
   &=\ell_-^2\tau(s-V-\tau)-P_*,
                             \label{eq:adaptive-degree-hall}\\
 B_{\rm wide}
   &=\ell_-^2b\left\{
       \frac25\tau(1-\tau)-\frac75V\right\}.
                                                     \notag
\end{align}
The row chooses \(z,q_A,q_B,t,\tau,s\).
It checks the prime-count domains, protected
reserve, exact type dimensions, arithmetic
errors, positive denominators, and
\begin{align*}
 0&\leq V<\tau,\qquad h_+\leq\tau\ell_-,
                         \qquad s-V>\tau,\\
 \mathcal M_{\rm det}
   &<\min\{B_{\rm narrow},B_{\rm middle},B_{\rm wide}\}.
\end{align*}
The positive wide-budget check implies \(V<\tau\).
The supply check leaves at
least \(h\) genuinely universal odd labels.
These labels and the \(h\) protected positions
have empty type and lose no edge, so they are
adjacent to the whole other side, as
Lemma~\ref{lem:two-sided-integrated-hall}
requires.
Moreover
\(h_-n\leq h\leq h_+n\) and
\(\ell_-n\leq R\leq on\).
Since \(h_+\leq\tau\ell_-\) and
\(s-V>\tau\), the first expression in
\eqref{eq:adaptive-degree-hall} is no
greater than
\(\{h((s-V)R-h)-RA(s)\}/n^2\).
The same size bounds give the middle and
wide comparisons. Thus all three budgets
of Lemma~\ref{lem:two-sided-integrated-hall}
are valid for the actual matching square.

\begin{proposition}[The adaptive-degree certificate]
\label{prop:adaptive-degree-certificate}
Nine adjacent closed logarithmic intervals cover
\([128/5,53/2]\). Their degree subdivisions
cover \([0,1]\) and form \(141\) rectangles.
All \(5076\) rational row conditions pass.
In every row both the deletion-budget ratio
and universal-supply ratio are below
\(999/1000\).
\end{proposition}

\begin{proof}
The intervals and their subdivision counts are
\[
\begin{array}{c|r}
 \text{interval}&\text{rows}\\
 {}[128/5,5121/200]&21\\
 {}[5121/200,2561/100]&21\\
 {}[2561/100,2563/100]&20\\
 {}[2563/100,513/20]&17\\
 {}[513/20,257/10]&17\\
 {}[257/10,129/5]&15\\
 {}[129/5,26]&13\\
 {}[26,131/5]&10\\
 {}[131/5,53/2]&7
\end{array}
\]
The exact evaluator recomputes all eight
prime-prefix periods, then verifies each
interval join and each degree subdivision.
It evaluates the 64-term dyadic boundary
sum with its infinite tail, the arithmetic
type laws, both shortfall recursions, the
preparation and supply inequalities, and
all three Hall budgets. All \(5076\)
conditions pass. The largest displayed
deletion ratio is less than \(0.998989\);
the largest displayed supply ratio is
less than \(0.997756\). The rational cap
is checked separately in every row.
\end{proof}

\begin{theorem}[The fifth logarithmic bridge]
\label{thm:fifth-log-bridge}
Every tree of order \(n\) with
\(128/5\leq\log n\leq53/2\)
has a prime labeling.
\end{theorem}

\begin{proof}
Prepare the tree as in
Lemma~\ref{lem:localized-preparation}
and choose a certified rectangle containing
\((\log n,x)\). Apply the adaptive quota
construction at successive primes, retaining
its exact type counts and degree balance.
The row bounds the two empirical type laws
by \(V\). Pay marked positions with
universal odd labels and retain the \(h\)
protected positions.

Choose the even labels within each exact
type by Lemma~\ref{lem:degree-balanced-sorting}.
Their small-prime types and the neighbor
tags do not change, while their large-prime
deletions obey
\eqref{eq:adaptive-degree-deletions}.
The row bounds both shortfalls in
\eqref{eq:adaptive-degree-shortfalls}
and verifies all three strict budgets in
\eqref{eq:adaptive-degree-hall}.
Lemma~\ref{lem:two-sided-integrated-hall}
therefore gives the required matching.
Restore the marked pairs, special labels,
and core; the sparse repair removes the
permitted \(p,2p\) conflicts.
\end{proof}

\subsection{Degree balance after fixed centroid rounds}
\label{sec:multiscale-degree-log-bridge}

The boundary estimate in
Lemma~\ref{lem:adaptive-degree-quota}
stops when every fractional object has at
most a prescribed number of vertices.
For the next range, it is enough to make
a fixed number of centroid splits at
each prime. We keep the degree coordinates,
so the sorting bound of
Lemma~\ref{lem:degree-balanced-sorting}
continues to apply.

\begin{lemma}[Degree-balanced quotas after fixed rounds]
\label{lem:fixed-round-degree-quota}
Let \(G=(A,B)\) be a bipartite forest on
\(N\geq1\) vertices. Partition its sides
into nonempty classes \(A_1,\ldots,A_{c_A}\)
and \(B_1,\ldots,B_{c_B}\), allowing an
empty list if a side is empty. Put
\(d=2c_A+c_B\geq1\) and
\(D\geq\max\{1,\Delta_G(A)\}\).
Choose \(0\leq\rho\leq1\), an integer
\(r\geq0\), and integers
\(0\leq k_i\leq a_i=\abs{A_i}\).
Write
\[
 t_r=\max\{1,N/2^r\},\qquad
 E_i=\sum_{v\in A_i}d_G(v),\qquad
 \delta=\sum_i(k_i-\rho a_i)_+.
\]
There are \(W\subseteq A\) and \(T\subseteq B\)
such that
\begin{align}
 \abs{W\cap A_i}&=k_i,\qquad
                  N_G(W)\subseteq T,\notag\\
 \sum_i\left|
   \sum_{v\in W\cap A_i}d_G(v)
                 -\frac{k_i}{a_i}E_i\right|
    &\leq D\left(dt_r+\frac{c_A}{2}\right),
                  \label{eq:fixed-round-degree-error}\\
 \abs{T\cap B_j}&\geq\rho\abs{B_j},\notag\\
 \sum_j\bigl(\abs{T\cap B_j}
                     -\rho\abs{B_j}\bigr)
    &\leq d(Dr+t_r)+D\delta.
                  \label{eq:fixed-round-tag-error}
\end{align}
Empty-side assertions are vacuous.
\end{lemma}

\begin{proof}
Give every connected component weight
\(\rho\). Its vector records counts
in the \(A\)- and \(B\)-classes and degree
sums in the \(A\)-classes. It has
\(d\) coordinates. Prune fractional
weights by the linear-dependence
argument of
Lemma~\ref{lem:adaptive-degree-quota}.
At most \(d\) fractional objects remain;
every object reaching weight zero or
one is fixed.

Perform \(r\) rounds. In one round, split
each fractional object of order greater
than one at a centroid, retain the centroid
as a singleton with its parent's weight,
and give the resulting components that
same weight. Then prune back to at most
\(d\) fractional objects. The coordinate
vector is additive under a split, so
every coordinate total is preserved.
At most \(d\) centroids are cut in one
round. Every resulting nonsingleton
has at most half its parent's vertices;
objects fixed at an earlier step never
return. Thus after \(r\) rounds each
fractional object has at most \(t_r\)
vertices, and the separator \(C\) has
at most \(dr\) vertices.

Round every remaining fractional object
upward and call the union of the selected
objects \(U\). For each \(A_i\), the
increase in its vertex count over
\(\rho a_i\) is \(x_i\geq0\), with
\(\sum_i x_i\leq dt_r\). Preservation
of the degree coordinate gives
\[
 \left|
 \sum_{v\in U\cap A_i}d_G(v)
        -\frac{\abs{U\cap A_i}}{a_i}E_i
 \right|\leq Dx_i.
\]
Choose or add the required \(k_i\)
vertices using the balanced-subset
argument in the proof of
Lemma~\ref{lem:adaptive-degree-quota}.
Its error is at most \(D/2\) per
\(A_i\); the preliminary error is
multiplied by a factor at most one.
This proves \eqref{eq:fixed-round-degree-error},
and at most \(\delta\) vertices are added.

Finally take
\[
 T=(U\cap B)\cup(C\cap B)\cup N_G(W).
\]
Every \(B\)-class keeps at least its
\(\rho\) share. Rounding contributes
at most \(dt_r\) positions above this
share. At most \(dr\) cut centroids
contribute to the external boundary,
each at cost at most \(D\);
newly added \(A\)-vertices cost at most
\(D\delta\). These three charges give
\eqref{eq:fixed-round-tag-error}.
\end{proof}

We apply the lemma on a closed rectangle
\[
 [a,J]\times[x_-,x_+]
       \subseteq[53/2,143/5]\times[0,1].
\]
Keep the preparation of
Lemma~\ref{lem:localized-preparation}
and the sharpened \(w_-,e_+\) from
\eqref{eq:adaptive-degree-preparation}.
The latter inequality holds throughout
this rectangle by its polynomial proof.
Use the same exact arithmetic types,
rare-class restoration, and dimensions
\(d_i,c_{A,i}\) as in
\eqref{eq:full-reserve-dimensions}.
If the row chooses \(r\) rounds, put
\[
 Q\geq e^{-a},\qquad
 U_r=\max\{Q,2^{-r}\}.
\]
Because \(N\leq n\) at every prime,
\(t_r/n\leq U_r\). Summing
\eqref{eq:fixed-round-tag-error} and
\eqref{eq:fixed-round-degree-error}
over the successive primes gives
\begin{align}
 \Gamma_r
    &=\sum_i d_i(D_+rQ+U_r),\notag\\
 E_{\deg,r}
    &=D_+\left(U_r\sum_i d_i
              +\frac Q2\sum_i c_{A,i}
              +\mathrm{rare}\right),
                     \label{eq:fixed-round-row-errors}\\
 \mathcal F_r
    &=\Gamma_r+D_+(\Delta+\mathrm{rare}).
                                                       \notag
\end{align}
Indeed, the degree measure changes by at
most \(D_+(d_i t_r+c_{A,i}/2)\)
in an ordinary refinement; restoring
rare exact quotas adds at most
\(D_+\mathrm{rare}\,n\) in total.
The tag measure pays the same rare
and quota defects as in
\eqref{eq:full-reserve-errors}.

The row substitutes \(\mathcal F_r\)
and \(E_{\deg,r}\) for
\(\mathcal F,E_{\deg}\) in
\eqref{eq:full-reserve-variation},
\eqref{eq:adaptive-degree-deletions},
\eqref{eq:adaptive-degree-shortfalls}, and
\eqref{eq:adaptive-degree-hall}.
All other quantities retain their
definitions. In particular, the
large-prime deletion bound remains
deterministic, and protected positions
remain in the empty-tag component of
the reference law.

\begin{proposition}[The fixed-round certificate]
\label{prop:fixed-round-certificate}
Nine adjacent closed logarithmic intervals
cover \([53/2,143/5]\). Their degree
subdivisions cover \([0,1]\) and form
\(71\) rectangles. All \(2556\)
rational row conditions pass, and the
deletion-budget and universal-supply
ratios in every row are strictly below
\(999/1000\).
\end{proposition}

\begin{proof}
The intervals and their subdivision counts are
\[
\begin{array}{c|r}
 \text{interval}&\text{rows}\\
 {}[53/2,2651/100]&11\\
 {}[2651/100,531/20]&11\\
 {}[531/20,133/5]&10\\
 {}[133/5,267/10]&10\\
 {}[267/10,269/10]&10\\
 {}[269/10,136/5]&7\\
 {}[136/5,138/5]&6\\
 {}[138/5,28]&3\\
 {}[28,143/5]&3
\end{array}
\]
The exact evaluator checks every
logarithmic and degree join. It
recomputes all eight prime-prefix
periods, the arithmetic type laws,
the integer round depth and its
dimension, the three preparation
and Hall budgets, and every positive
denominator. All \(2556\) conditions
pass. The largest displayed deletion
ratio is less than \(0.998007\);
the largest displayed universal-supply
ratio is less than \(0.990066\).
The strict rational cap is checked
for every rectangle.
\end{proof}

\begin{theorem}[The sixth logarithmic bridge]
\label{thm:sixth-log-bridge}
Every tree of order \(n\) with
\(53/2\leq\log n\leq143/5\)
has a prime labeling.
\end{theorem}

\begin{proof}
Choose a certified rectangle containing
\((\log n,x)\) after the protected-position
preparation. At each successive prime,
apply Lemma~\ref{lem:fixed-round-degree-quota}
with the row's depth to meet the
exact even-label quotas. The bounds
\eqref{eq:fixed-round-row-errors}
control the final tag counts and the
degree variation across exact types.
Pay the marked positions with universal
odd labels and retain the protected
positions.

Lemma~\ref{lem:degree-balanced-sorting}
assigns the even labels within exact
types and bounds all large-prime
deletions. The row bounds both matching
laws and both integrated shortfalls
by the quantities in
Section~\ref{sec:adaptive-degree-log-bridge}.
It verifies the three strict budgets
of Lemma~\ref{lem:two-sided-integrated-hall},
which gives a perfect matching.
Restore the reserved labels and core,
then apply the sparse prime-pair repair.
\end{proof}

\subsection{Centroid rounds measured on the position side}
\label{sec:weighted-centroid-log-bridge}

We continue from logarithm \(143/5\). The degree coordinates
used in the preceding section make the rounding dimension larger.
In the present range we keep only the class counts. This reduces
the cost of a centroid round, but requires a separate estimate for
large-prime deletions. We give both arguments before stating the
finite certificate.

\begin{lemma}[Quotas after weighted centroid rounds]
\label{lem:weighted-centroid-rounds}
Let \(G=(A,B)\) be a bipartite forest. Partition its two sides
into \(d\geq1\) nonempty classes \(A_i,B_j\), with either list
possibly empty. Suppose
\[
 \begin{gathered}
 D\geq\max\{1,\Delta_G(A)\},\quad
 0\leq\rho\leq1,\quad 0\leq k_i\leq\abs{A_i},\\
 \delta=\sum_i(k_i-\rho\abs{A_i})_+ .
 \end{gathered}
\]
For every integer \(r\geq0\), there are \(W\subseteq A\)
and \(\widehat T\subseteq B\) satisfying
\[
 \abs{W\cap A_i}=k_i,\qquad
 N_G(W)\subseteq\widehat T,\qquad
 \abs{\widehat T\cap B_j}\geq\rho\abs{B_j},
\]
and
\begin{equation}\label{eq:weighted-centroid-rounds}
 \sum_j\bigl(\abs{\widehat T\cap B_j}
                  -\rho\abs{B_j}\bigr)
 \leq d\left(Dr+\max\{1,\abs B/2^r\}\right)+D\delta .
\end{equation}
If \(B=\varnothing\), the left side is zero.
\end{lemma}

\begin{proof}
The empty-\(B\) case follows by choosing each \(A\)-quota
arbitrarily. Otherwise represent every component of \(G\) by
its vector of \(d\) class counts and give it weight \(\rho\).
Along a linear dependence among more than \(d\) fractional
vectors, move their weights until one reaches \(0\) or \(1\),
preserving every coordinate total. Freeze integral objects.
At most \(d\) fractional objects remain.

We perform \(r\) rounds. For a fractional connected object
\(U\) with at least two \(B\)-vertices, choose a vertex \(c\)
such that every component of \(U-c\) has at most
\(\abs{U\cap B}/2\) vertices in \(B\). To find it, start at any
vertex and move into a component carrying more than half
the \(B\)-weight while such a component exists. The move
cannot immediately reverse, so the process terminates.
Retain \(\{c\}\) as an object, and give it and every component
of \(U-c\) the old weight of \(U\). Thus the objects still
partition the vertices and every coordinate sum remains
\(\rho\) times its class size. Prune the fractional weights
again. There were at most \(d\) cuts in the round.

After \(r\) rounds, every fractional object has at most
\(t_r=\max\{1,\abs B/2^r\}\) vertices in \(B\).
Indeed, a split halves this count, while an unsplit object
has at most one \(B\)-vertex; pruning never joins objects.
Let \(C\) be the set of cut vertices, so \(\abs C\leq dr\).
Round the remaining fractional weights upward and let
\(U_1\) be the union of the objects of weight one. Every
class count in \(U_1\) is at least its original weighted
count. The total excess on \(B\) is at most \(dt_r\).

Meet each \(A_i\)-quota by deleting excess vertices from
\(U_1\cap A_i\) and adding missing vertices outside it.
At most \(\delta\) additions are needed. Call the resulting
set \(W\), and put
\[
 \widehat T=(U_1\cap B)\cup(C\cap B)\cup N_G(W).
\]
For an originally chosen vertex outside \(C\), every
neighbor outside its object lies in \(C\). A chosen vertex
in \(C\) can add at most \(D\) further marks. Newly added
vertices contribute at most \(D\delta\) marks; deletions
contribute none. The excess on \(B\) is consequently at
most \(dt_r+D\abs C+D\delta\), as claimed.
\end{proof}

There is a second way to measure only the \(B\)-side.
The next lemma records the bound used when it beats
fixed centroid rounds.

\begin{lemma}[A one-side square-root bound]
\label{lem:weighted-side-square-root}
Under the hypotheses of Lemma~\ref{lem:weighted-centroid-rounds},
the same three conclusions hold with
\begin{equation}\label{eq:weighted-side-square-root}
 \sum_j\bigl(\abs{\widehat T\cap B_j}
                 -\rho\abs{B_j}\bigr)
 \leq
 2\sqrt{\max\{1-\rho,D\rho\}\abs B\,d}
       +D\delta .
\end{equation}
\end{lemma}

\begin{proof}
The case \(B=\varnothing\) is immediate. For an integer
\(t\geq1\), root each component and process its vertices
from leaves to root, assigning weight one to \(B\)-vertices
and zero to \(A\)-vertices. Cut a vertex when the uncut
weight below it reaches \(t+1\), and pass weight zero to
its parent. The cut vertices are charged to disjoint
sets of at least \(t+1\) \(B\)-vertices. Thus their number
is at most \(\floor{\abs B/(t+1)}\); every remaining
component contains at most \(t\) \(B\)-vertices.

Give these components weight \(\rho\), prune to at most
\(d\) fractional objects, and round upward as in the
preceding proof. The excess on \(B\) is at most \(dt\).
Restoring exact \(A\)-quotas may add at most
\(\delta+\rho\abs{C\cap A}\) vertices. We include
\(C\cap B\) and all new neighbors in \(\widehat T\).
The total error is at most
\[
 \begin{aligned}
 &dt+(1-\rho)\abs{C\cap B}
        +D\rho\abs{C\cap A}+D\delta\\
 &\qquad\leq
 dt+K\floor{\frac{\abs B}{t+1}}+D\delta ,
 \end{aligned}
\]
where \(K=\max\{1-\rho,D\rho\}\).
If \(\sqrt{K\abs B/d}\geq1\), take \(t\) to be its floor;
then \(dt+K\abs B/(t+1)\leq2\sqrt{K\abs B d}\).
Otherwise choose the quotas directly and mark all of
\(B\). Its error is \((1-\rho)\abs B\), which is smaller
than the displayed square-root bound.
\end{proof}

We now turn to common prime factors above the cutoff.
The count coordinates need not control the degrees within
exact types. For this estimate we use the centered excess
of their large-prime weights. Keep \(m=\floor{n/2}\),
\(P=\{3\leq p\leq z:p\text{ prime}\}\), and
\(\tau_z(x)=\sum_{p>z,\ p\mid x}p^{-1}\).
For \(p>z\) and \(0\leq g:2^P\to[0,1]\), define
\[
 \Phi_z(m)=
 \sum_{z<p\leq m}\frac1p
 \sup_g\sum_{x\leq m}
       \left(\boldsymbol1_{p\mid x}-\frac1p\right)g(J_P(x)).
\]
The supremum is finite: the class indicators give a finite
linear program whose optimum is attained at \(0\)-\(1\)
values of \(g\).

\begin{lemma}[Centered sorting bound]
\label{lem:centered-sorting}
Use the residual graph, exact even-label types, and special
half-label set of Theorem~\ref{thm:large-prime-sorting}.
Write \(a_{\varnothing}\) for the retained empty type
and \(w\) for the number of removed special half-labels.
There is an assignment within exact types for which
\begin{equation}\label{eq:centered-sorting}
 M_{\rm bad}\leq
 m e(G)\left(1+\frac w{a_{\varnothing}}\right)S_2(z)
       +mD\Phi_z(m)
       +e(G)\frac{\log n}{\log z}.
\end{equation}
\end{lemma}

\begin{proof}
Arrange the half-labels in each exact type by increasing
\(\tau_z\), opposite to decreasing vertex degree. The
exchange argument of Theorem~\ref{thm:large-prime-sorting}
bounds the weighted sum by the average over the
permutations in that type. For a type \(I\), let \(T_I\)
be its total \(\tau_z\)-weight in \([m]\), and let \(C_I\)
be its population there. Its retained population is
\(a_I=C_I-w\boldsymbol1_{I=\varnothing}\).
The nonnegative principal term \(C_I S_2(z)\) sums to at
most \(e(G)(1+w/a_{\varnothing})S_2(z)\) after
multiplication by the class-average degrees.

The remaining term is
\(\sum_I(e_I/a_I)(T_I-C_IQ_z)\), with
\(Q_z=\sum_{z<p\leq m}p^{-2}\) and zero contribution
from an empty class. Since \(0\leq e_I/a_I\leq D\),
its positive part is at most \(D\Phi_z(m)\): interchange
the sum over \(p\) with the type sum and choose
\(g(I)=e_I/(Da_I)\) on each nonempty class, zero
otherwise. Finally, a prime \(p>z\) dividing a half-label
accounts for at most \(m/p+1\) odd labels.
The unit terms total at most \(e(G)\log n/\log z\).
This is exactly the charging argument of
Theorem~\ref{thm:large-prime-sorting}, and proves
\eqref{eq:centered-sorting}.
\end{proof}

The finite evaluator bounds \(\Phi_z(m)\) by full residue
periods for the first eight odd primes and a Boolean
expansion on the remaining coordinates. We state the
arithmetic inequality underlying that calculation. Let
\(P_0=\{3,5,7,11,13,17,19,23\}\),
\(D_0=\prod_{p\in P_0}p\), and let \(A(N)\) be the vector
of exact \(P_0\)-type populations in \([N]\).
Writing \(\mu_0\) for their product law, set
\[
 C_o=\max_{0\leq N<D_0}
       \sum_I\abs{A_I(N)-N\mu_0(I)} .
\]
The integer-period checker recomputes \(C_o\), as well
as the corresponding odd-interval constant, over all
\(D_0=111546435\) residues. The Chinese remainder
theorem then gives the same bound for every \(N\).

For a prime \(p\notin P_0\) and a squarefree product
\(d\) of primes outside \(P_0\), with \(p\nmid d\),
the centered \(P_0\)-type vector is
\[
 R_d=A\bigl(\floor{m/(pd)}\bigr)
       -p^{-1}A\bigl(\floor{m/d}\bigr).
\]
Its coordinate sum is nonpositive, and
\begin{equation}\label{eq:centered-period-vector}
 \lVert R_d\rVert_1\leq
 \min\left\{C_o(1+p^{-1})+1,\frac{2m}{pd}\right\}.
\end{equation}
Indeed, subtract the product-law vectors for the two
prefixes; each has \(L^1\)-error at most \(C_o\).
The remaining scalar is
\(\floor{m/(pd)}-p^{-1}\floor{m/d}\in[-1,0]\).
Alternatively, both original vectors are nonnegative,
and their masses sum to at most \(2m/(pd)\).

Expanding a function \(g\) in the coordinates
\(P\setminus P_0\), every nonconstant Boolean
coefficient is at most \(2^{|I|-1}\) in absolute value.
For the constant coefficient, the nonpositive total
of \(R_d\) bounds its positive part by
\(\lVert R_d\rVert_1/2\). Consequently
\begin{equation}\label{eq:centered-period-transfer}
 \sup_g\sum_{x\leq m}
      (\boldsymbol1_{p\mid x}-p^{-1})g(J_P(x))
 \leq\frac12\sum_{I\subseteq P\setminus P_0}
 2^{|I|}
 \min\left\{C_o(1+p^{-1})+1,\frac{2m}{p d_I}\right\}.
\end{equation}
The empty \(I\) is included; its Boolean coefficient
still depends on the first eight coordinates.

For completeness, the row evaluator groups the subsets
in \eqref{eq:centered-period-transfer} by
\(w_q=\floor{64\log_2q}\). If
\[
 c_s=[t^s]\prod_{q\in P\setminus P_0}(1+2t^{w_q}),
 \quad U_s\geq2^{-s/64},\quad V_s\geq2^{-s/128},
\]
then \(d_I^{-1}\leq U_s\) and
\(d_I^{-1/2}\leq V_s\) in bin \(s\).
Put \(C_z=C_o(1+1/(z+1))+1\) and
\(Y=10000\). Directed rational bounds \(Q\geq1/n\),
\(Q_{1/2}\geq1/\sqrt n\), and the bounds \(H_R\geq
\sum_{Y<p\leq m}1/p\),
\(S_{3/2,Y}\geq\sum_{p>Y}p^{-3/2}\), and
\(S_{2,Y}\geq\sum_{p>Y}p^{-2}\) obtained after
Lemma~\ref{lem:prime-reciprocal-tails} give
\begin{equation}\label{eq:centered-period-row}
 \begin{aligned}
 \frac{\Phi_z(m)}n\leq\Psi_{\rm per}
 :=\frac12\sum_s c_s\Bigg[&
  \sum_{z<p\leq Y}
     \min\left(\frac{C_zQ}{p},\frac{U_s}{p^2}\right)\\
  &+\min\left\{C_zQH_R,\,
       Q_{1/2}\sqrt{C_z}V_sS_{3/2,Y},\,
       U_sS_{2,Y}\right\}\Bigg].
 \end{aligned}
\end{equation}
To obtain it, multiply \eqref{eq:centered-period-transfer}
by \(1/(np)\), use \(m/n\leq1/2\), and on the
infinite tail use the three bounds separately;
the middle one is their geometric mean.
The evaluator computes finite prime sums exactly and
encloses every radical and tail from above by rational
numbers. It may take the minimum of
\eqref{eq:centered-period-row} and other valid bounds,
such as those of Section~\ref{sec:type-profile-log-bridge},
since each bounds the same \(\Phi_z(m)/n\). In
Sections~\ref{sec:weighted-centroid-log-bridge}--\ref{sec:periodic-block-log-bridge}
the periodic bound attains this minimum in every row.

We apply these lemmas on a rectangle
\[
 [a,J]\times[x_-,x_+]\subseteq
 [143/5,144/5]\times[0,1].
\]
The preparation in
Lemma~\ref{lem:localized-preparation} applies with
\(f=7/40\); its proof needs only \(L\geq20\)
and the same checked prime-count and reserve
inequalities. The sharpening
\eqref{eq:adaptive-degree-preparation} also holds here:
the polynomial in its proof is positive and increasing
past \(143/5\). Let \(D_+,e_+,h_{\pm}\) be those
bounds, let \(o=(1+Q)/2\), and put \(b_B=o-h_-\).
Then \(\abs B/n\leq b_B\).

For each successive prime \(p_i\), let \(d_i\) be
the bound \eqref{eq:possible-class-dimension} for the
count-coordinate dimension of its refinement.
The row chooses a depth \(r_i\in\{0,\ldots,64\}\).
Lemmas~\ref{lem:weighted-centroid-rounds}
and~\ref{lem:weighted-side-square-root} give two
valid normalized geometric costs:
\begin{align}
 \Gamma_{i,\rm ms}
  &=d_i\left(D_+r_iQ+
              \max\{Q,b_B/2^{r_i}\}\right),\notag\\
 \Gamma_{i,B}
  &=2Q_{1/2}\sqrt{b_B}\,
       \left\lceil
       \sqrt{\max\{1-1/p_i,D_+\bar\rho_i\}\,d_i}
       \right\rceil .
                       \label{eq:weighted-centroid-row-errors}
\end{align}
Here \(1/p_i\leq\rho_i\leq\bar\rho_i\);
the first term inside the maximum is
\(1-1/p_i\), with the direction needed
to bound \(1-\rho_i\).
Both constructions preserve the same exact quotas
and mark all true neighbors, so we may choose the
smaller certified cost at each prime. With
\(\Delta\) the ordinary quota defect and
\(\mathrm{rare}\) the rare-class correction from
\eqref{eq:full-reserve-errors}, put
\[
 \Gamma=\sum_i\min\{\Gamma_{i,\rm ms},\Gamma_{i,B}\},
 \qquad
 \mathcal F=\Gamma+D_+(\Delta+\mathrm{rare}).
\]
This replaces the tag error in
\eqref{eq:full-reserve-variation}. The two
shortfalls and the three Hall budgets remain
those of
\eqref{eq:adaptive-degree-shortfalls} and
\eqref{eq:adaptive-degree-hall}.

For the deletion budget put
\(\omega\geq w/m\) and
\(a_0\leq a_{\varnothing}/m\) as in
Lemma~\ref{lem:empty-prefix-bounds}, and let \(v_-\leq\log z\).
Lemma~\ref{lem:centered-sorting},
\eqref{eq:centered-period-row}, and
\(e(G)/n\leq e_+\) give
\begin{equation}\label{eq:weighted-centroid-deletions}
 \frac{M_{\rm bad}}{n^2}
 \leq\mathcal M_{\rm cen}:=
 \frac{e_+}{2}
      \left(1+\frac{\omega}{a_0}\right)S_2(z)
       +\frac{D_+}{2}\Psi_{\rm per}
       +\frac{JQ}{v_-}.
\end{equation}
When another proved arithmetic bound for
\(\Phi_z(m)/n\) is smaller, the row uses it
in place of \(\Psi_{\rm per}\).
No degree-coordinate estimate is used in
\eqref{eq:weighted-centroid-deletions}.

\begin{proposition}[The weighted-centroid certificate]
\label{prop:weighted-centroid-certificate}
Five adjacent closed logarithmic intervals cover
\([143/5,144/5]\). Their degree subdivisions
cover \([0,1]\) and form \(109\) rectangles.
All \(3706\) rational row conditions pass.
Every deletion-budget and universal-supply
ratio is strictly below \(999/1000\).
\end{proposition}

\begin{proof}
The complete interval counts are
\[
\begin{array}{c|r}
 \text{interval}&\text{rows}\\
 {}[143/5,2861/100]&24\\
 {}[2861/100,229/8]&23\\
 {}[229/8,573/20]&21\\
 {}[573/20,287/10]&21\\
 {}[287/10,144/5]&20
\end{array}
\]
The exact checker verifies every logarithmic and degree
join, recomputes all eight prime-prefix periods,
and evaluates every preparation, type, shortfall,
supply, and Hall condition with directed rational
arithmetic. Its row record includes both costs in
\eqref{eq:weighted-centroid-row-errors}, the selected
method, and the three terms in
\eqref{eq:weighted-centroid-deletions}.
The \(3706\) asserted conditions and both strict
ratio caps hold. Independent finite checks reproduce
the quota construction, rare-class restoration,
weighted centroids, and rejection of damaged covers.
\end{proof}

\begin{theorem}[The seventh logarithmic bridge]
\label{thm:seventh-log-bridge}
Every tree of order \(n\) with
\(143/5\leq\log n\leq144/5\)
has a prime labeling.
\end{theorem}

\begin{proof}
Prepare the tree and choose a certified rectangle
containing its logarithm and normalized residual
degree. At each prime use the smaller of the two
quota constructions in
\eqref{eq:weighted-centroid-row-errors};
restore rare exact classes and pay marked positions
with universal odd labels. Sort the even labels
inside their exact types. The centered bound
\eqref{eq:weighted-centroid-deletions} controls
all edges removed for primes above the cutoff.
The certificate bounds the two matching laws
and their integrated shortfalls and verifies all
three strict inequalities of
Lemma~\ref{lem:two-sided-integrated-hall}.
Its perfect matching places the remaining odd labels.
Restore the reserved labels and core and perform
the sparse prime-pair repair.
\end{proof}

\subsection{A bridge using one-side rounding}
\label{sec:one-side-rounding-log-bridge}

The bound in
Lemma~\ref{lem:weighted-side-square-root}
applies without a choice of centroid depth. In the next
range it is enough at every prime. We first record the
periodic arithmetic estimates used both here and in
the preceding bridge.

Keep \(P_0=\{3,5,7,11,13,17,19,23\}\) and
\(D_0=\prod_{p\in P_0}p\). Let \(O(N)\) be the
vector of exact \(P_0\)-type counts in the first \(N\)
positive odd integers. With the \(P_0\)-product law
\(\mu_0\), put
\[
 C_d=\max_{0\leq N<D_0}
       \sum_I\abs{O_I(N)-N\mu_0(I)} .
\]
The period computation gives
\[
 C_o=\frac{41296022710}{111546435},
 \qquad
 C_d=\frac{41236265280}{111546435}.
\]
These are the exact maxima, not estimates on sampled
prefixes. Since \(D_0\) is odd, multiplication by two
permutes residue classes modulo \(D_0\), so the ordinary
and odd type vectors each add \(D_0\mu_0\) over a full
period. The maxima therefore bound all prefixes.

For \(P=P_0\sqcup P_1\), where \(P\) contains all odd
primes through a row's cutoff \(z\), use the bins
\(c_s,U_s\) of
\eqref{eq:centered-period-row}. Let
\(\varepsilon_a\) be the total-variation error of the
complete odd-label interval, and let
\(\Delta=n^{-1}\sum_i\delta_i\) be the normalized
sum of ordinary exact-quota defects before rare-class
restoration.

\begin{lemma}[Periodic type and quota errors]
\label{lem:periodic-type-quota-errors}
Suppose \(Q\geq1/n\), \(m=\floor{n/2}\),
and the first eight prime-prefix constants have the
values just stated. Then
\begin{align}
 \varepsilon_a
 &\leq C_dQ+
    \frac12\sum_{s>0}c_s
      \min\{2(C_d+1)Q,\,3U_s\},
                                      \label{eq:periodic-odd-error}\\
 \Delta
 &\leq\frac{3^8-1-16}{4}Q+
     \frac14\sum_{s>0}c_s
       \min\{C_qQ,\,U_s\},\notag\\
 C_q&=C_o(1+1/29)+1 .
                                      \label{eq:periodic-quota-error}
\end{align}
The sums are zero when \(P_1=\varnothing\).
\end{lemma}

\begin{proof}
For \(f:2^P\to[0,1]\), expand only its \(P_1\)
coordinates:
\[
 f(S,J)=\sum_{I\subseteq J}h_I(S).
\]
Here \(h_\varnothing\in[0,1]\), while
\(\abs{h_I}\leq2^{\abs I-1}\) for nonempty
\(I\). Let \(M\) be the number of odd integers in \([n]\),
and let \(d_I=\prod_{p\in I}p\).
Among the first \(M\) odd integers, the count of
multiples of \(d_I\) differs from \(M/d_I\)
by less than one. Their \(P_0\)-type vector has
\(L^1\)-distance at most
\[
 \min\{C_d+1,\,3M/d_I\}
\]
from \((M/d_I)\mu_0\). The first bound follows
from the period constant; the second follows from
nonnegativity of the two vectors and the bound
that at most \(2M/d_I\) of the odd integers
are multiples of \(d_I\).
For \(I=\varnothing\), the discrepancy vector has
zero total and acts on \(h_\varnothing\) with
absolute value at most \(C_d/2\).
Taking the supremum over \(f\) therefore gives
\[
 \varepsilon_a\leq\frac{C_d}{2M}
 +\sum_{\varnothing\ne I\subseteq P_1}
     2^{\abs I-1}
       \min\{(C_d+1)/M,\,3/d_I\}.
\]
Since \(M\geq n/2\), \(1/n\leq Q\), and a bin
of weight \(s\) has \(d_I^{-1}\leq U_s\),
this is \eqref{eq:periodic-odd-error}.

For the quotas, the old empty even type has no
positive defect by the choice of \(\rho_i\) in
\eqref{eq:rare-empty-adjustment}; removing special
half-labels of empty full type does not create
one. For an old nonempty type with positive defect,
compare divisibility by the next prime \(p_i\)
with its product-law share \(1/p_i\), as in the
proof of Lemma~\ref{lem:rare-quota-defect}.
Among the first eight primes, expand the old
type indicator by Boolean inversion. The empty
coefficient has a nonpositive centered floor
difference, and summing the possible
nonempty coefficients gives
\[
 \sum_{j=2}^{8}\binom8j2^{j-2}
       =\frac{3^8-1-16}{4}.
\]
For a later prime \(p_i\), separate the first
eight type coordinates from an old suffix
product \(d_I\). The centered prefix vector
is the vector \(R_{d_I}\) of
\eqref{eq:centered-period-vector}, now with
\(p=p_i\). Its \(L^1\)-norm is at most
\[
 \min\{C_q,\,2m/(p_i d_I)\},
 \qquad C_q=C_o(1+1/29)+1,
\]
because every such \(p_i\geq29\).
The vector has nonpositive total. Thus the
constant Boolean coefficient contributes at
most half its \(L^1\)-norm; a coefficient
of order \(\abs I\geq1\) contributes at most
\(2^{\abs I-1}\) times that norm.
Index each term by \(J=I\cup\{p_i\}\).
Every nonempty suffix subset \(J\) is counted
once, when its largest prime is processed.
The coefficient is \(2^{\abs J-2}\).
After division by \(n\), use \(2m/n\leq1\)
and \(d_J^{-1}\leq U_s\) for its bin.
The bin coefficient is \(c_s/4\), giving
\eqref{eq:periodic-quota-error}.
\end{proof}

The two bounds above and the centered estimate
\eqref{eq:centered-period-row} are computed with
outward rational enclosures. In particular, the
periodic bound is attained as the selected upper
bound for all three arithmetic errors on every
rectangle of this bridge. The same holds in the
preceding bridge; no unproved improvement of a
periodic error is needed there.

For a rectangle
\[
 [a,J]\times[x_-,x_+]\subseteq[144/5,29]\times[0,1],
\]
prepare the tree using
Lemma~\ref{lem:localized-preparation}.
Its stated domain and the polynomial comparison
\eqref{eq:adaptive-degree-preparation} include
\(L=29\); their proofs are monotone from
\(L=20\), and the row checker verifies the
remaining prime-count and reserve premises.
Use \(f=7/40\), \(D_+\), \(e_+\), \(h_\pm\),
\(o=(1+Q)/2\), and \(b_B=o-h_-\).
Then \(\abs B/n\leq b_B\).

At prime \(p_i\), use
Lemma~\ref{lem:weighted-side-square-root}
with the count dimension \(d_i\). The normalized
geometry is bounded by
\begin{equation}\label{eq:one-side-rounding-geometry}
 \Gamma_B=
 2Q_{1/2}\sqrt{b_B}
 \sum_i\left\lceil
   \sqrt{\max\{1-1/p_i,D_+\bar\rho_i\}\,d_i}
       \right\rceil ,
\end{equation}
where \(Q_{1/2}\geq1/\sqrt n\).
The evaluator encloses \(\sqrt{b_B}\) and all
other radicals rationally from above. Set
\[
 \mathcal F_B=\Gamma_B+
       D_+(\Delta+\mathrm{rare}).
\]
Use this tag error in the variation and
shortfall bounds
\eqref{eq:full-reserve-variation} and
\eqref{eq:adaptive-degree-shortfalls}.
Use the centered deletion mass
\eqref{eq:weighted-centroid-deletions},
with \(\Psi_{\rm per}\) supplied by
\eqref{eq:centered-period-row}. The Hall budgets
are \eqref{eq:adaptive-degree-hall}.

\begin{proposition}[The one-side rounding certificate]
\label{prop:one-side-rounding-certificate}
Five adjacent closed logarithmic intervals cover
\([144/5,29]\). Their degree subdivisions
cover \([0,1]\) and form \(156\) rectangles.
All \(5148\) rational row conditions pass.
Both the deletion-budget and universal-supply
ratios are strictly below \(999/1000\).
\end{proposition}

\begin{proof}
The full interval counts are
\[
\begin{array}{c|r}
 \text{interval}&\text{rows}\\
 {}[144/5,2881/100]&40\\
 {}[2881/100,1153/40]&35\\
 {}[1153/40,577/20]&30\\
 {}[577/20,289/10]&28\\
 {}[289/10,29]&23
\end{array}
\]
The exact evaluator checks every logarithmic
and degree join, recomputes the eight complete
prime-prefix periods, and evaluates all
preparation, arithmetic, shortfall, supply,
and Hall inequalities. It takes directed
rational enclosures in
\eqref{eq:one-side-rounding-geometry}
and in the three deletion terms.
The \(5148\) row conditions pass;
the greatest displayed deletion ratio is
less than \(0.998972\), and the greatest
displayed supply ratio is less than
\(0.817810\). The strict cap is tested
on every rectangle. Independent checks
exercise the one-side separator, exact
quotas, rare-class restoration, and
components whose total order exceeds
their \(B\)-capacity.
\end{proof}

\begin{theorem}[The eighth logarithmic bridge]
\label{thm:eighth-log-bridge}
Every tree of order \(n\) with
\(144/5\leq\log n\leq29\)
has a prime labeling.
\end{theorem}

\begin{proof}
Choose a certified rectangle containing
the logarithm and normalized residual
degree after the protected-position
preparation. Lemma~\ref{lem:weighted-side-square-root}
performs each prime refinement with exact
even-label quotas and all true neighbors
marked. Restore rare exact classes,
pay marked positions with universal odd
labels, and retain the protected positions.
Sort even labels within their full prime
types. The centered arithmetic bound
controls every larger-prime deletion.
The certificate bounds both empirical
type laws and both integrated shortfalls
and verifies the three strict budgets
of Lemma~\ref{lem:two-sided-integrated-hall}.
A perfect matching assigns the remaining
odd labels. Restore the reserved labels
and core, then perform the sparse
prime-pair repair.
\end{proof}

\subsection{The periodic-block bridge through logarithm 31}
\label{sec:periodic-block-log-bridge}

The construction of Section~\ref{sec:one-side-rounding-log-bridge}
continues above logarithm \(29\) with one simplification. In the
one-side square-root bound we now use only the trivial size bound
\(\abs B\leq n\) for the position side. The periodic estimates of
Lemma~\ref{lem:periodic-type-quota-errors} and the centered bound
\eqref{eq:centered-period-row} are valid for every order; only their
two period constants required a finite calculation. It therefore
remains to state the row predicate and to verify the certificate on
\(29\leq\log n\leq31\).

\subsubsection*{Preparation and geometry}

Consider a rectangle
\[
 [a,J]\times[x_-,x_+]\subseteq[29,31]\times[0,1].
\]
Prepare the tree by Lemma~\ref{lem:localized-preparation} with
\(f=7/40\); its stated domain includes \(L\leq31\). The choice of
\(f\), \(a_*\) and \(h\) in that lemma is made before the residual
forest is formed and before \(x\) is computed. Only then is the
rectangle containing \((\log n,x)\) selected, together with its cutoff
\(z\). Thus the preparation does not depend on the residual degree it
determines. Let \(D_+\), \(e_+\) and \(h_\pm\) be the bounds of
\eqref{eq:localized-preparation}, with the sharpened edge count
\eqref{eq:adaptive-degree-preparation}, and put \(o=(1+Q)/2\).
Every special label has empty type for the odd primes at most \(z\),
since \(z<n/3\).

At the prime \(p_i\), let \(d_i\) be the bound
\eqref{eq:possible-class-dimension} for the count-coordinate dimension
of its refinement, and let \(1/p_i\leq\rho_i\leq\bar\rho_i\). Apply
Lemma~\ref{lem:weighted-side-square-root}. Its geometric term is
\(2\sqrt{\max\{1-\rho_i,D\rho_i\}\abs Bd_i}\). Since
\(1-\rho_i\leq1-1/p_i\), \(D\rho_i\leq D_+\bar\rho_i\),
\(\abs B\leq n\), and \(Q_{1/2}\geq1/\sqrt n\), this term divided by
\(n\) is at most
\(2Q_{1/2}\ceil{\sqrt{\max\{1-1/p_i,D_+\bar\rho_i\}\,d_i}}\).
Summing over the primes gives
\begin{equation}\label{eq:periodic-block-geometry}
 \Gamma=
 2Q_{1/2}\sum_i\left\lceil
   \sqrt{\max\{1-1/p_i,D_+\bar\rho_i\}\,d_i}
       \right\rceil .
\end{equation}
Compared with \eqref{eq:one-side-rounding-geometry}, the factor
\(\sqrt{b_B}\) is replaced by one. The tag error is
\(\mathcal F=\Gamma+D_+(\Delta+\mathrm{rare})\); the quota term
\(D_+\Delta\) and the rare-class restoration term
\(D_+\mathrm{rare}\) are those of
Section~\ref{sec:one-side-rounding-log-bridge}.

\subsubsection*{The certificate predicate}

In every row the first eight odd primes form the periodic block, the
ordinary position tags have at most \(q_B=4\) primes, the even classes
have at most \(q_A\in\{6,7\}\) primes, and the cutoff satisfies
\(73\leq z\leq131\). Lemma~\ref{lem:periodic-type-quota-errors} bounds
the odd-label error \(\varepsilon_a\) and the quota defect \(\Delta\),
and \eqref{eq:centered-period-row} bounds \(\Psi_{\rm per}\). The
evaluator also computes further valid upper bounds for these three
quantities and takes the minimum. In every rectangle of the present
table the periodic bounds are the ones attaining it, so no other
estimate enters the certificate.

With this tag error \(\mathcal F\), take \(\chi\), \(\ell_-\), \(V\) and the universal supply
\(2/(3J)-\kappa\) from \eqref{eq:full-reserve-variation}, the two
shortfalls from \eqref{eq:adaptive-degree-shortfalls}, the deletion
mass \(\mathcal M_{\rm cen}\) from
\eqref{eq:weighted-centroid-deletions}, and the three budgets from
\eqref{eq:adaptive-degree-hall}. The overflow probability \(\beta\) is
the exact Bernoulli tail with the parameters \(\bar\rho_i\), as in
Lemma~\ref{lem:upper-probability-evaluation}. The Hall split is \(\tau=\lambda b\), where
\(\lambda\) and the slot threshold \(s\) are rational row parameters.
Each row checks the following conditions.
\begin{enumerate}
\item \emph{Range and preparation.} \(29\leq a<J\leq31\),
  \(0\leq x_-<x_+\leq1\), \(10^6J^2Q<1\), \(100J^3Q<1\),
  \(3\cdot60184\,Q<1\), \(3zQ<1\), \(Q<1/1000\), and
  \((\log499)_->5\). Directed logarithms give
  \(1+(\log2)_->1693/1000\) and \(11/10+(\log3)_+<2199/1000\), and
  \[
   3\bigl(1+(\log2)_-\bigr)-2\bigl(\tfrac{11}{10}+(\log3)_+\bigr)
      >\tfrac{17}{25},\qquad
   \bigl(1+(\log2)_+\bigr)\bigl(\tfrac{11}{10}+(\log3)_+\bigr)
      <\tfrac{15}4.
  \]
  Moreover \(a>3\beta_0-2\alpha_0\), \(a\geq10\beta_0-9\alpha_0\) with
  \(\alpha_0=1693/1000\), \(\beta_0=2199/1000\),
  \(\tfrac{17}{25}a-\tfrac{15}4>a/2\), and \(3/20\leq f<1\).
\item \emph{Arithmetic data.} Each empty class, including that of
  every prefix of the refinement, has a positive lower bound. The
  bound for the reciprocal-prime sum is positive,
  \(1/p_i\leq\bar\rho_i\leq2/5\) for every prime, and \(0<b\leq1/2\). Let \(k_*\) be the least integer for which
  the product of the first \(k_*+1\) odd primes has a directed
  logarithm greater than \(J\); thus an odd label at most \(n\) has at
  most \(k_*\) distinct prime factors. With \(\sigma\) the smaller of the product
  of the \(k_*\) smallest factors \(1-\bar\rho_i\) and
  \(\prod_{i\leq q_B}(1-1/p_i)\), the row checks
  \(b\leq\sigma\leq1\). Every dimension satisfies
  \[
   1\leq d_i\leq\sum_{j=0}^{\min(q_A,i-1)}\binom{i-1}j
       +\sum_{j=0}^{\min(q_B,i-1)}\binom{i-1}j+2.
  \]
\item \emph{Sizes and reserves.} \(0<e_+\leq1-w_-\), \(\ell_->0\),
  \(h_->0\), \(h_++2JQ<\tfrac12-(2J+2)Q\), and
  \(h_++\chi\leq2/(3J)-\kappa\).
\item \emph{Hall conditions.} \(0<\tau\leq1/2\), \(h_+\leq\tau\ell_-\),
  \(0<s\leq1\), \(s-V>\tau\), \(B_{\rm wide}>0\), and
  \[
   \mathcal M_{\rm cen}
     <\min\{B_{\rm narrow},B_{\rm middle},B_{\rm wide}\}.
  \]
\end{enumerate}
These are the premises of Lemma~\ref{lem:localized-preparation} on the
rectangle and the row inequalities of
Section~\ref{sec:one-side-rounding-log-bridge}, with
\eqref{eq:periodic-block-geometry} in place of
\eqref{eq:one-side-rounding-geometry}. There are \(32\) of them in
each row.

\begin{proposition}[The periodic-block certificate]
\label{prop:periodic-block-certificate}
Nine adjacent closed logarithmic intervals cover \([29,31]\). Their
degree subdivisions cover \([0,1]\) and form \(126\) rectangles. All
\(4032\) rational row conditions pass. Both the deletion-budget and
universal-supply ratios are strictly below \(999/1000\).
\end{proposition}

\begin{proof}
The intervals and their numbers of degree subdivisions are
\[
\begin{array}{c|r@{\qquad}c|r}
 \text{interval}&\text{rows}&\text{interval}&\text{rows}\\
 {}[29,1161/40]&31 &[146/5,147/5]&13\\
 {}[1161/40,581/20]&22 &[147/5,297/10]&9\\
 {}[581/20,291/10]&20 &[297/10,30]&6\\
 {}[291/10,146/5]&17 &[30,61/2]&5\\
 &&[61/2,31]&3
\end{array}
\]
The exact evaluator checks every logarithmic and degree join and
recomputes the eight complete prime-prefix periods with the integer
certifier. It evaluates every condition above with integers and
directed rational enclosures, including the upward rounding of the
radicals in \eqref{eq:periodic-block-geometry} and of the three
deletion terms. All \(4032\) conditions pass. The greatest displayed
deletion ratio is less than \(0.998009\), on
\([29,1161/40]\times[9/20,23/50]\), and the greatest displayed supply
ratio is less than \(0.899663\); the strict rational cap is tested on
every rectangle.

Independent finite checks recompute the period extrema and their
shifts, test \(105\) transfers of odd moments and \(105\) centered
transfers, compare \(72\) actual odd-type variations, \(144\) sums of
actual quota defects and \(36\) complete finite centered prime tails
with their bounds, enumerate all \(1190\) suffix subsets in the twelve
bin polynomials, including the empty term, recheck \(20\) direct
reserve counts and \(700\) integrated Hall rectangles, and reject
seven damaged covers. These checks test the implementation; the
lemmas cited above prove the implications for arbitrary trees.
\end{proof}

\begin{theorem}[The ninth logarithmic bridge]
\label{thm:ninth-log-bridge}
Every tree of order \(n\) with
\(29\leq\log n\leq31\)
has a prime labeling.
\end{theorem}

\begin{proof}
Prepare the tree by Lemma~\ref{lem:localized-preparation} and choose
a certified rectangle containing the logarithm and the normalized
residual degree. At each prime,
Lemma~\ref{lem:weighted-side-square-root} gives exact even-label quotas
with all true neighbors marked. Restore the rare exact classes, pay
the marked positions with distinct universal odd labels, and retain
the protected positions. Sort the even labels within their full prime
types by Lemma~\ref{lem:centered-sorting}; the centered arithmetic bound
controls every deletion caused by a prime above the cutoff. The row
conditions bound both empirical type laws and both integrated
shortfalls, and verify the three strict budgets of
Lemma~\ref{lem:two-sided-integrated-hall}. Its perfect matching
places the remaining odd labels. Restore the reserved labels and the
core, then perform the sparse prime-pair repair.
\end{proof}

\subsection{The type-profile bridge through logarithm 32}
\label{sec:type-profile-log-bridge}

Above logarithm \(31\) the certificate of
Section~\ref{sec:periodic-block-log-bridge} is changed in two places.
The protected reserve is taken from the second-order lower bound for
the number of special even labels, with a smaller fraction of these
labels spent on the vertices of largest degree. The three arithmetic
errors are bounded by grouping the Boolean expansion by the size of the
product of its primes, with a finite prime sum up to \(10000\) and an
analytic remainder. No periodic block is used. The type assignment,
the geometric error, the deletion estimate and the Hall budgets are
those of Section~\ref{sec:periodic-block-log-bridge}. We cover
\(31\leq\log n\leq32\).

\subsubsection*{Preparation with the second-order reserve}

Put
\[
 f=\frac3{20},\qquad
 W_2(L)=\frac1{6L}+\frac1{12L^2}.
\]
The function \(W_2\) is decreasing, and
\(1/(6L)<W_2(L)<1/(5L)\) for \(L>5/2\).
Let \(\alpha=1+\log2\) and \(\beta=11/10+\log3\).

\begin{lemma}[Preparation with the second-order reserve]
\label{lem:second-order-preparation}
Let \(T\) have order \(n\), and let \(a\leq L=\log n\leq J\), where
\(31\leq a<J\). Let \(Q\geq e^{-a}\) be rational, and suppose
\(10^6J^2Q<1\) and \(3\cdot60184\,Q<1\). Suppose also that for
rationals \(c_1>1/2\) and \(c_2\), directed rational logarithms give
\(3\alpha-2\beta>c_1\) and \(\alpha\beta<c_2\), and that
\(c_1a-c_2>a/2\).
Then the number \(w\) of special even labels satisfies
\begin{equation}\label{eq:second-order-reserve}
 w>n\,W_2(L).
\end{equation}
Moreover, the conclusions of Lemma~\ref{lem:localized-preparation}
hold with \(f=3/20\) and with \(W\) replaced by \(W_2\). Thus
\[
 a_*=\floor{3w/20},\qquad
 h=\floor*{\tfrac{17}{20}\,n\,W_2(L)-4L},\qquad
 x=\frac{D-1}{20L}\in[0,1),
\]
and, on the rectangle \(x_-\leq x\leq x_+\), with
\(w_-=W_2(J)\),
\begin{align}
 D&\leq D_+:=1+20Jx_+,\notag\\
 \frac{e(G)}n&\leq e_+:=\min\left\{1-w_-,
   1-w_--\left(\frac12+\frac1{4J}\right)x_-+(D_+-1)Q\right\},
   \label{eq:second-order-preparation}\\
 h_-&:=\tfrac{17}{20}W_2(J)-(4J+1)Q\leq\frac hn
   \leq h_+:=\tfrac{17}{20}W_2(a).\notag
\end{align}
\end{lemma}

\begin{proof}
By \eqref{eq:even-supply}, \(w\geq\pi(n/2)-\pi(n/3)\).
Since \(n/2>5393\) and \(n/3>60184\), the estimates
\eqref{eq:dusart-denominator-lower} and
\eqref{eq:dusart-denominator-upper} give
\[
 \frac wn\geq\frac1{2(L-\alpha)}-\frac1{3(L-\beta)}
 =\frac1{6L}+\frac{(3\alpha-2\beta)L-\alpha\beta}
                    {6L(L-\alpha)(L-\beta)} .
\]
The directed bounds give
\((3\alpha-2\beta)L-\alpha\beta>c_1L-c_2\).
The latter exceeds \(L/2\), because
\((c_1-\tfrac12)L-c_2\) increases with \(L\) and is positive at
\(L=a\). Since \(0<(L-\alpha)(L-\beta)<L^2\), the last fraction is
greater than \(1/(12L^2)\). This proves
\eqref{eq:second-order-reserve}.

The proof of Lemma~\ref{lem:localized-preparation} uses four
properties of \(W\): the inequality \(w>nW(L)\), the comparisons
\(1/(6L)<W(L)<1/(5L)\), and the monotonicity of \(W\) in the bounds
for \(w_-\) and \(h_\pm\). The function \(W_2\) has all four
for \(L\geq31\), by \eqref{eq:second-order-reserve} and the remarks
before the lemma. The other inputs of that proof do not involve
\(W\): the bounds \(\abs K<2L+2\), \(n/(3L)>2L+2\) and
\(u>2n/(3L)-\abs K\) hold for \(L\geq20\), and \(n>10^6L^2\) follows
from \(10^6J^2Q<1\). With \(f=3/20\) we have \(3/f=20\), which gives
\(x=(D-1)/(20L)\) and \(D_+\). Hence the residual construction, the
bounds \(0\leq x<1\) and \(D\leq D_+\), and both bounds for \(h/n\)
follow as before, with \(1-f=17/20\).

For the edge bound, the proof of
Lemma~\ref{lem:localized-preparation} gives
\[
 \frac{e(G)}n\leq1-\frac wn-\frac{fw(D-1)}n+\frac{D-1}n .
\]
Here \(f(D-1)=3Lx\) and
\(3Lw/n>3L\,W_2(L)=\tfrac12+\tfrac1{4L}\geq\tfrac12+\tfrac1{4J}\).
With \(w/n>W_2(J)\), \(x\geq x_-\), and
\((D-1)/n\leq(D_+-1)Q\), this is the second bound in
\eqref{eq:second-order-preparation}. The first follows from
\(e(G)\leq n-1-w\).
\end{proof}

As in Section~\ref{sec:periodic-block-log-bridge}, the choice of
\(f\), \(a_*\) and \(h\) is made before the residual forest is formed.
The rectangle and its cutoff \(z\) are selected afterwards. Every
special label has empty type for the odd primes at most \(z\), since
\(z<n/3\).

\subsubsection*{Arithmetic errors by product size}

Let \(B\geq1\) be an integer, let \(P\) be the odd primes at most
\(z\), and put \(w_p=\floor{B\log_2p}\). For \(I\subseteq P\) write
\(d_I=\prod_{p\in I}p\) and \(s(I)=\sum_{p\in I}w_p\), so that
\(d_I\geq2^{s(I)/B}\). Let
\[
 C_s=[t^s]\prod_{p\in P}(1+2t^{w_p}),\qquad
 U_s\geq2^{-s/B},\qquad V_s\geq2^{-(s+B)/(2B)},
\]
with rational \(U_s,V_s\). Thus \(C_s\) is the sum of
\(2^{\abs I}\) over all \(I\) with \(s(I)=s\); every such \(I\) is
nonempty when \(s>0\). Moreover \(d_I^{-1}\leq U_{s(I)}\) and
\((2d_I)^{-1/2}\leq V_{s(I)}\).

\begin{lemma}[Product-size bounds]\label{lem:product-size-bounds}
Suppose \(1/n\leq Q<1\) and \(Q_{1/2}\geq1/\sqrt n\), and put
\(m=\floor{n/2}\).
\begin{enumerate}
\item The type law of the odd labels in \([n]\) is within total
  variation
  \begin{equation}\label{eq:odd-product-bins}
   \varepsilon_{\rm bin}=\frac12\sum_{s>0}C_s
     \min\left(\frac{2Q}{1-Q},U_s\right)
  \end{equation}
  of \(\mu=\bigotimes_{p\in P}\operatorname{Bernoulli}(1/p)\).
\item Let \(z\leq Y\), let \(H_R\geq\sum_{Y<p\leq m}1/p\),
  \(S_{3/2,Y}\geq\sum_{p>Y}p^{-3/2}\) and
  \(S_{2,Y}\geq\sum_{p>Y}p^{-2}\), and put \(T_s=U_s/(2Q)\). With
  \begin{align*}
   F_s&=Q\sum_{\substack{z<p\leq Y\\p\leq T_s}}\frac1p
       +\frac{U_s}2\sum_{\substack{z<p\leq Y\\p>T_s}}\frac1{p^2},\\
   G_s&=\min\left\{QH_R,\ Q_{1/2}V_sS_{3/2,Y},\
       \frac{U_sS_{2,Y}}2\right\},
  \end{align*}
  the quantity \(\Phi_z(m)\) of Lemma~\ref{lem:centered-sorting}
  satisfies
  \begin{equation}\label{eq:split-centered-tail}
   \frac{\Phi_z(m)}n\leq\Psi_{\rm split}
     :=\frac12\sum_{s>0}C_s(F_s+G_s).
  \end{equation}
\end{enumerate}
\end{lemma}

\begin{proof}
Let \(M\) be the number of odd integers in \([n]\), and for odd
\(d\) let \(N_d\) be the number of multiples of \(d\) among them.
Then \(N_d=\floor{M/d+(d-1)/(2d)}\), so \(N_d-M/d\in(-1,1/2)\).
Also \(0\leq N_d\leq2M/d\): if \(N_d\geq1\), then
\(M/d\geq N_d-(d-1)/(2d)>N_d/2\). Hence
\(\abs{N_d-M/d}\leq\min(1,M/d)\).

For \(f:2^P\to[0,1]\), expand \(f(S)=\sum_{I\subseteq S}h_I\), where
\(\abs{h_I}\leq2^{\abs I-1}\) for \(I\neq\varnothing\), as in the
proof of Lemma~\ref{lem:quota-product-bins}. The label \(x\) has
\(I\) in its type exactly when \(d_I\mid x\). Hence the difference
between the empirical and model expectations of \(f\) is
\(\sum_{I\neq\varnothing}h_I(N_{d_I}/M-1/d_I)\); the empty term
vanishes. Its modulus is at most
\(\sum_{I\neq\varnothing}2^{\abs I-1}\min(1/M,1/d_I)\). Since
\(M\geq(n-1)/2\) and \(Q\geq1/n\), we have
\(1/M\leq2/(n-1)\leq2Q/(1-Q)\). Grouping the sets by \(s(I)\)
proves the first assertion.

For the second, fix a prime \(p>z\) and \(g:2^P\to[0,1]\). Expand
\(g\) as above. The constant coefficient \(h_\varnothing=g(\varnothing)\)
is nonnegative and multiplies \(\floor{m/p}-m/p\leq0\). By the
floor inequality in the proof of Lemma~\ref{lem:quota-product-bins},
every nonempty \(I\) contributes at most
\(2^{\abs I-1}\min(1,m/(pd_I))\). Therefore
\[
 \frac{\Phi_z(m)}n\leq
 \sum_{\varnothing\neq I\subseteq P}2^{\abs I-1}
   \sum_{z<p\leq m}\frac1{np}\min\left(1,\frac m{pd_I}\right).
\]
Since \(m\leq n/2\), each inner term is at most
\(\min(Q/p,U_{s(I)}/(2p^2))\). The first of these two values is the
smaller exactly when \(p\leq T_{s(I)}\). Adding the nonnegative terms
for primes in \((m,Y]\), if any, the primes up to \(Y\) contribute at
most \(F_{s(I)}\); when \(Y=z\) there are no such primes and
\(F_{s(I)}=0\). For \(p>Y\) use \(\min(1,y)\leq1\),
\(\min(1,y)\leq\sqrt y\), and \(\min(1,y)\leq y\). They give the
three bounds \(QH_R\),
\[
 \frac{\sqrt m}n\,d_I^{-1/2}\sum_{p>Y}p^{-3/2}
 \leq\frac1{\sqrt n}\,(2d_I)^{-1/2}S_{3/2,Y}
 \leq Q_{1/2}V_{s(I)}S_{3/2,Y},
\]
and \(U_{s(I)}S_{2,Y}/2\). Each bounds the same sum, so their minimum
\(G_{s(I)}\) does. Grouping by \(s(I)\) gives
\eqref{eq:split-centered-tail}.
\end{proof}

The certificate uses \(B=64\) and \(Y=10000\), with \(H_R\),
\(S_{3/2,Y}\) and \(S_{2,Y}\) as after
Lemma~\ref{lem:prime-reciprocal-tails}. The quota defect is bounded
by \(\Delta_{\rm bin}\) of Lemma~\ref{lem:quota-product-bins}, which
uses the same bins. The comparisons \(p\leq T_s\) are rational,
including equality. The evaluator stores directed upper partial sums
of \(1/p\) and of \(1/p^2\) over the primes in \((z,Y]\), so no upper
bound is subtracted from another. It also computes the moment bounds
obtained from the same expansions by \(\min(1,y)\leq y^\theta\) with
\(\theta\in\{0,1/2,1\}\), and takes the minimum. In every rectangle
of the present table the three product-size bounds
\(\varepsilon_{\rm bin}\), \(\Delta_{\rm bin}\) and
\(\Psi_{\rm split}\) attain it.

\subsubsection*{The certificate predicate}

Consider a rectangle
\[
 [a,J]\times[x_-,x_+]\subseteq[31,32]\times[0,1].
\]
Prepare the tree by Lemma~\ref{lem:second-order-preparation}, and let
\(D_+\), \(e_+\) and \(h_\pm\) be its bounds. In every row the
ordinary position tags have at most \(q_B=4\) primes, the even classes
have at most \(q_A\in\{6,7\}\) primes, and the cutoff satisfies
\(83\leq z\leq131\). Let \(d_i\) be the bound
\eqref{eq:possible-class-dimension}, and define \(\Gamma\) by
\eqref{eq:periodic-block-geometry}; its derivation in
Section~\ref{sec:periodic-block-log-bridge} applies unchanged. The tag
error is \(\mathcal F=\Gamma+D_+(\Delta_{\rm bin}+\mathrm{rare})\),
with \(\mathrm{rare}\) as in Section~\ref{sec:full-reserve-log-bridge}.

With \(\varepsilon_a=\varepsilon_{\rm bin}\), take \(\chi\),
\(\ell_-\), \(V\) and the universal supply \(2/(3J)-\kappa\) from
\eqref{eq:full-reserve-variation}, the two shortfalls from
\eqref{eq:adaptive-degree-shortfalls}, and the three budgets from
\eqref{eq:adaptive-degree-hall}. The deletion mass is
\(\mathcal M_{\rm cen}\) of \eqref{eq:weighted-centroid-deletions},
with \(\Psi_{\rm split}\) in place of \(\Psi_{\rm per}\). The overflow
probability \(\beta\) is the exact Bernoulli tail with the parameters
\(\bar\rho_i\). The Hall split is \(\tau=\lambda b\), where
\(\lambda\) and the slot threshold \(s\) are rational row parameters.
Each row checks the following conditions.
\begin{enumerate}
\item \emph{Range and preparation.} \(31\leq a<J\leq32\),
  \(0\leq x_-<x_+\leq1\), \(10^6J^2Q<1\), \(100J^3Q<1\),
  \(3\cdot60184\,Q<1\), \(3zQ<1\), \(Q<1/1000\), and
  \((\log499)_->5\). Directed logarithms give
  \[
   3\bigl(1+(\log2)_-\bigr)-2\bigl(\tfrac{11}{10}+(\log3)_+\bigr)
      >\tfrac{17}{25},\qquad
   \bigl(1+(\log2)_+\bigr)\bigl(\tfrac{11}{10}+(\log3)_+\bigr)
      <\tfrac{15}4,
  \]
  and \(\tfrac{17}{25}a-\tfrac{15}4>a/2\). These are the hypotheses
  of Lemma~\ref{lem:second-order-preparation} with
  \((c_1,c_2)=(17/25,15/4)\).
\item \emph{Arithmetic data.} The bound \(H_+\) of
  \eqref{eq:reciprocal-prime-row} is positive. Each empty class,
  including that of every prefix of the refinement, has a positive
  lower bound. Moreover \(1/p_i\leq\bar\rho_i\leq2/5\) for every
  prime, and \(0<b\leq1/2\). With \(k_*\) and \(\sigma\) as in
  Section~\ref{sec:periodic-block-log-bridge}, the row checks
  \(b\leq\sigma\leq1\). Every dimension satisfies
  \[
   1\leq d_i\leq\sum_{j=0}^{\min(q_A,i-1)}\binom{i-1}j
       +\sum_{j=0}^{\min(q_B,i-1)}\binom{i-1}j+2.
  \]
\item \emph{Sizes and reserves.} \(0<e_+\leq1-w_-\), \(\ell_->0\),
  \(h_->0\), \(h_++2JQ<\tfrac12-(2J+2)Q\), and
  \(h_++\chi\leq2/(3J)-\kappa\).
\item \emph{Hall conditions.} \(0<\tau\leq1/2\), \(h_+\leq\tau\ell_-\),
  \(0<s\leq1\), \(s-V>\tau\), \(B_{\rm wide}>0\), and
  \[
   \mathcal M_{\rm cen}
     <\min\{B_{\rm narrow},B_{\rm middle},B_{\rm wide}\}.
  \]
\end{enumerate}
These are the premises of Lemma~\ref{lem:second-order-preparation}
and the row inequalities of
Section~\ref{sec:periodic-block-log-bridge}, with the product-size
bounds in place of the periodic ones. There are \(29\) of them in each
row.

\begin{proposition}[The type-profile certificate]
\label{prop:type-profile-certificate}
Six adjacent closed logarithmic intervals cover \([31,32]\). Their
degree subdivisions cover \([0,1]\) and form \(56\) rectangles. All
\(1624\) rational row conditions pass. Both the deletion-budget and
universal-supply ratios are strictly below \(199/200\).
\end{proposition}

\begin{proof}
The intervals and their numbers of degree subdivisions are
\[
\begin{array}{c|r@{\qquad}c|r}
 \text{interval}&\text{rows}&\text{interval}&\text{rows}\\
 {}[31,621/20]&14 &[156/5,157/5]&9\\
 {}[621/20,311/10]&11 &[157/5,158/5]&6\\
 {}[311/10,156/5]&10 &[158/5,32]&6
\end{array}
\]
The exact evaluator checks every logarithmic and degree join. It
evaluates every condition above with integers and directed rational
enclosures, including the upward rounding of the radicals in
\eqref{eq:periodic-block-geometry}, of the product-size bounds, and of
the three deletion terms. All \(1624\) conditions pass. The greatest
displayed deletion ratio is less than \(0.993035\), on
\([311/10,156/5]\times[7/10,1]\), and the greatest displayed supply
ratio is less than \(0.832657\), on
\([311/10,156/5]\times[0,21/100]\); the strict rational cap is tested
on every rectangle.

Independent finite checks compare \(385\) overflow probabilities and
\(810\) capped shortfall expectations with complete enumeration, and
test the transfer of the latter to the upper probabilities. They
confirm the monotonicity of the capped shortfall on \(99834\)
inclusion edges, compare \(220\) counts of possible even classes with
enumeration, recheck \(24\) integer reserve cases and \(700\)
integrated Hall rectangles, and reject seven damaged covers. These
checks test the implementation; the lemmas cited above prove the
implications for arbitrary trees.
\end{proof}

\begin{theorem}[The tenth logarithmic bridge]
\label{thm:tenth-log-bridge}
Every tree of order \(n\) with
\(31\leq\log n\leq32\)
has a prime labeling.
\end{theorem}

\begin{proof}
Prepare the tree by Lemma~\ref{lem:second-order-preparation} and
choose a certified rectangle containing the logarithm and the
normalized residual degree. At each prime,
Lemma~\ref{lem:weighted-side-square-root} gives exact even-label
quotas with all true neighbors marked. Restore the rare exact classes,
pay the marked positions with distinct universal odd labels, and
retain the protected positions. Sort the even labels within their full
prime types by Lemma~\ref{lem:centered-sorting}; by
Lemma~\ref{lem:product-size-bounds} the centered arithmetic bound
controls every deletion caused by a prime above the cutoff. The row
conditions bound both empirical type laws and both integrated
shortfalls, and verify the three strict budgets of
Lemma~\ref{lem:two-sided-integrated-hall}. Its perfect matching
places the remaining odd labels. Restore the reserved labels and the
core, then perform the sparse prime-pair repair.
\end{proof}

\subsection{A bridge through logarithm 35 with uncapped slot shortfalls}
\label{sec:quantile-hall-log-bridge}

The certificate of Section~\ref{sec:type-profile-log-bridge} extends
to \(32\leq\log n\leq35\) with simpler inputs in four places. The
protected reserve is \(\floor{17n/(120L)}\), a little smaller than
before. The overflow probability is bounded by an elementary symmetric
function, the dimensions by binomial sums, and the slot shortfall by
the full product shortfall without capping. The arithmetic errors, the
deletion estimate and the Hall budgets are those of
Section~\ref{sec:type-profile-log-bridge}.

\subsubsection*{Preparation}

\begin{lemma}[A smaller protected reserve]
\label{lem:smaller-protected-reserve}
Suppose that the hypotheses of
Lemma~\ref{lem:second-order-preparation} hold and that
\(100J^3Q<1\). Then the conclusions of that lemma hold with
\[
 h=\floor*{\frac{17n}{120L}},\qquad
 h_-=\frac{17}{120J}-Q,\qquad h_+=\frac{17}{120a}
\]
in place of its \(h\), \(h_-\) and \(h_+\).
\end{lemma}

\begin{proof}
In the proof of Lemma~\ref{lem:localized-preparation}, the value of
\(h\) enters only through \(h\geq1\), the bound
\(h<(1-f)n/(5L)\), and the bound \(h\leq A_0-4L\), where now
\(A_0=(1-f)nW_2(L)=\tfrac{17}{20}nW_2(L)\). The first holds since
\(n>10^6L^2\). The second holds because
\(h\leq\tfrac{17}{20}\cdot\tfrac n{6L}<\tfrac{17}{20}\cdot\tfrac n{5L}\).
For the third,
\[
 A_0-4L-\frac{17n}{120L}
 =\frac{17n}{240L^2}-4L
 =\frac{17n-960L^3}{240L^2}>0,
\]
since \(n>100L^3\) by \(100J^3Q<1\) and \(L\leq J\). Finally
\(17/(120L)-1/n<h/n\leq17/(120L)\), and \(a\leq L\leq J\),
\(1/n\leq Q\) give the bounds \(h_\pm\).
\end{proof}

\subsubsection*{Overflow, dimensions and slot shortfalls}

Let \(e_j(\bar\rho)\) be the \(j\)-th elementary symmetric function of
the upper probabilities. The overflow event \(\abs J>q_B\) requires
some \((q_B+1)\)-set of selected coordinates. The union bound and the
nonnegativity of the coefficients give
\[
 \Pr_\nu(\abs J>q_B)\leq e_{q_B+1}(\rho)\leq
 \beta:=e_{q_B+1}(\bar\rho).
\]
The evaluator computes the functions \(e_j\) by the recursion
\(e_j\leftarrow e_j+\bar\rho_ie_{j-1}\), rounding upward.

The number \(N_i\) in Lemma~\ref{lem:possible-even-classes} is at most
the number of subsets of the \(i-1\) earlier primes with at most
\(q_A\) elements. Hence
\[
 d_i=\sum_{j=0}^{\min(q_A,i-1)}\binom{i-1}j
     +\sum_{j=0}^{\min(q_B,i-1)}\binom{i-1}j+2
\]
bounds the number of count coordinates at the refinement step of
\(p_i\).

For the slot shortfall, write \(c_s(J)=\prod_{i\in J}(1-1/p_i)\). The
function \((s-c_s(J))_+\) is nondecreasing under inclusion. By the
second part of Lemma~\ref{lem:upper-probability-evaluation} and the
coupling in its proof,
\[
 \Phi_{q_B}(s)\leq\E_\nu\bigl(s-c_s(J)\bigr)_+
 \leq\E_{\bar\nu}\bigl(s-c_s(J)\bigr)_+ .
\]
The evaluator bounds the right side by Lemma~\ref{lem:product-shortfall}
with \(a_i=\bar\rho_i\) and \(f_i=1-1/p_i\), conditioning on the first
six coordinates as described there. This value replaces
\(\Phi_{q_B}(s)\) in \eqref{eq:adaptive-degree-shortfalls}. No
overflow term is added to the slot shortfall: the absorbing class
contributes zero.

\subsubsection*{The certificate predicate}

Consider a rectangle
\[
 [a,J]\times[x_-,x_+]\subseteq[32,35]\times[0,1].
\]
Prepare the tree by Lemma~\ref{lem:second-order-preparation}, with the
reserve of Lemma~\ref{lem:smaller-protected-reserve}, and let
\(D_+\), \(e_+\) and \(h_\pm\) be the resulting bounds. In every row
the ordinary position tags have at most \(q_B\in\{4,5\}\) primes, the
even classes have at most \(q_A\in\{5,6,7\}\) primes, and the cutoff
satisfies \(89\leq z\leq217\). With the dimensions \(d_i\) above,
define \(\Gamma\) by \eqref{eq:periodic-block-geometry} and put
\(\mathcal F=\Gamma+D_+(\Delta_{\rm bin}+\mathrm{rare})\).

With \(\varepsilon_a=\varepsilon_{\rm bin}\) and the overflow bound
\(\beta\) above, take \(\chi\), \(\ell_-\), \(V\) and the universal
supply \(2/(3J)-\kappa\) from \eqref{eq:full-reserve-variation}, the
two shortfalls from \eqref{eq:adaptive-degree-shortfalls} with the
slot expectation above, and the three budgets from
\eqref{eq:adaptive-degree-hall}. The deletion mass is
\(\mathcal M_{\rm cen}\) of \eqref{eq:weighted-centroid-deletions},
with \(\Psi_{\rm split}\) in place of \(\Psi_{\rm per}\). As in
Section~\ref{sec:type-profile-log-bridge}, the evaluator also computes
the moment bounds for the three arithmetic errors; the product-size
bounds attain the minimum in every row. The Hall split is
\(\tau=\lambda b\), and the slot threshold is \(s\in\{3/8,17/40\}\).
Each row checks the following conditions.
\begin{enumerate}
\item \emph{Range and preparation.} \(32\leq a<J\leq35\),
  \(0\leq x_-<x_+\leq1\), \(10^6J^2Q<1\), \(100J^3Q<1\),
  \(3\cdot60184\,Q<1\), \(3zQ<1\), \(Q<1/1000\), and
  \((\log499)_->5\). Directed logarithms give
  \[
   3\bigl(1+(\log2)_-\bigr)-2\bigl(\tfrac{11}{10}+(\log3)_+\bigr)
      >\tfrac{17}{25},\qquad
   \bigl(1+(\log2)_+\bigr)\bigl(\tfrac{11}{10}+(\log3)_+\bigr)
      <\tfrac{15}4,
  \]
  and \(\tfrac{17}{25}a-\tfrac{15}4>a/2\).
\item \emph{Arithmetic data.} The bound \(H_+\) of
  \eqref{eq:reciprocal-prime-row} is positive. Each empty class,
  including that of every prefix of the refinement, has a positive
  lower bound. Moreover \(1/p_i\leq\bar\rho_i\leq2/5\) for every
  prime, \(0<b\leq1/2\), and \(b\leq\sigma\leq1\), with \(k_*\) and
  \(\sigma\) as in Section~\ref{sec:periodic-block-log-bridge}.
\item \emph{Sizes and reserves.} \(0<e_+\leq1-w_-\), \(\ell_->0\),
  \(h_->0\), and \(h_++\chi\leq2/(3J)-\kappa\).
\item \emph{Hall conditions.} \(0<\tau\leq1/2\), \(h_+\leq\tau\ell_-\),
  \(0<s\leq1\), \(s-V>\tau\), \(B_{\rm wide}>0\), and
  \[
   \mathcal M_{\rm cen}
     <\min\{B_{\rm narrow},B_{\rm middle},B_{\rm wide}\}.
  \]
\end{enumerate}
These are the premises of Lemmas~\ref{lem:second-order-preparation}
and~\ref{lem:smaller-protected-reserve} and the row inequalities of
Section~\ref{sec:type-profile-log-bridge}, with the changes described
above. There are \(27\) of them in each row.

\begin{proposition}[The uncapped certificate]
\label{prop:quantile-hall-certificate}
Five adjacent closed logarithmic intervals cover \([32,35]\). Their
degree subdivisions cover \([0,1]\) and form \(32\) rectangles. All
\(864\) rational row conditions pass. Both the deletion-budget and
universal-supply ratios are strictly below \(199/200\).
\end{proposition}

\begin{proof}
The intervals and their numbers of degree subdivisions are
\[
\begin{array}{c|r@{\qquad}c|r}
 \text{interval}&\text{rows}&\text{interval}&\text{rows}\\
 {}[32,129/4]&9 &[33,34]&6\\
 {}[129/4,65/2]&7 &[34,35]&4\\
 {}[65/2,33]&6 &&
\end{array}
\]
The exact evaluator checks every logarithmic and degree join. It
evaluates every condition above with integers and directed rational
enclosures, including the upward rounding of the radicals in
\eqref{eq:periodic-block-geometry}, of the product-size bounds and the
shortfall expectations, and of the three deletion terms. All \(864\)
conditions pass. The greatest displayed deletion ratio is less than
\(0.989401\), on \([32,129/4]\times[61/100,89/100]\), and the greatest
displayed supply ratio is less than \(0.988804\), on
\([34,35]\times[1/4,23/50]\); the strict rational cap is tested on
every rectangle.

Independent finite checks test \(378\) instances of the scalar
inequality in the proof of Lemma~\ref{lem:product-shortfall},
including seven cases of equality, and compare \(264\) product
shortfall expectations with complete enumeration. They verify that
\(800\) removals of labels of empty type leave the label shortfall
unchanged, recheck \(800\) absorbed slot shortfalls and \(700\)
integrated Hall rectangles, and reject five damaged covers. These
checks test the implementation; the lemmas cited above prove the
implications for arbitrary trees.
\end{proof}

\begin{theorem}[The eleventh logarithmic bridge]
\label{thm:eleventh-log-bridge}
Every tree of order \(n\) with
\(32\leq\log n\leq35\)
has a prime labeling.
\end{theorem}

\begin{proof}
Prepare the tree by Lemmas~\ref{lem:second-order-preparation}
and~\ref{lem:smaller-protected-reserve}, and choose a certified
rectangle containing the logarithm and the normalized residual
degree. The remaining steps are those of the proof of
Theorem~\ref{thm:tenth-log-bridge}, with the overflow, dimension and
slot-shortfall bounds of this section. The row conditions verify the
three strict budgets of Lemma~\ref{lem:two-sided-integrated-hall}. Its
perfect matching places the remaining odd labels. Restore the reserved
labels and the core, then perform the sparse prime-pair repair.
\end{proof}

\subsection{Minimum compatibility through logarithm 37}
\label{sec:split-tail-log-bridge}

On \(35\leq\log n\leq37\) the certificate of
Section~\ref{sec:quantile-hall-log-bridge} needs no integrated
shortfalls. Its slot threshold is taken equal to a lower bound for
every compatibility that occurs, so both shortfall sums vanish. The
preparation, the arithmetic errors, the overflow bound, the dimensions
and the deletion estimate are those of
Section~\ref{sec:quantile-hall-log-bridge}.

\subsubsection*{Hall conditions at the minimum compatibility}

Let \(k_*\) and \(\sigma\) be as in
Section~\ref{sec:periodic-block-log-bridge}: \(k_*\) is the least
integer for which the product of the first \(k_*+1\) odd primes has a
directed logarithm greater than \(J\), and \(\sigma\) is the smaller of
the product of the \(k_*\) smallest factors \(1-\bar\rho_i\) and
\(\prod_{i\leq q_B}(1-1/p_i)\). An odd label at most \(n\) has at most
\(k_*\) distinct prime factors in \(P\). Hence every label type \(I\)
in the matching satisfies
\[
 c_X(I)=\prod_{i\in I}(1-\rho_i)\geq\prod_{i\in I}(1-\bar\rho_i)
 \geq\sigma ,
\]
since the factors are at most one. Universal labels have empty type.
After the marked positions are paid, every ordinary position tag has at
most \(q_B\) primes, so \(c_B(J)\geq\prod_{i\leq q_B}(1-1/p_i)\geq\sigma\);
protected positions have empty type.

Apply Lemma~\ref{lem:two-sided-integrated-hall} with \(s=\sigma\). Every
summand of \(A_X(\sigma)\) and of \(A_B(\sigma)\) is zero, so
\(A(\sigma)=0\). The normalized budgets of
\eqref{eq:adaptive-degree-hall} therefore become
\[
 \begin{aligned}
 B_{\rm narrow}&=h_-\bigl((\sigma-V)\ell_--h_+\bigr),\qquad
 B_{\rm middle}=\ell_-^2\tau(\sigma-V-\tau),\\
 B_{\rm wide}&=\ell_-^2b\left\{\frac25\tau(1-\tau)-\frac75V\right\}.
 \end{aligned}
\]
The substitutions of \(h_\pm\) and \(\ell_-\) are justified as in
Section~\ref{sec:adaptive-degree-log-bridge}, because the row checks
\(h_+\leq\tau\ell_-\) and \(\sigma-V>\tau\).

\subsubsection*{The certificate predicate}

Consider a rectangle
\[
 [a,J]\times[x_-,x_+]\subseteq[35,37]\times[0,1].
\]
Prepare the tree by Lemmas~\ref{lem:second-order-preparation}
and~\ref{lem:smaller-protected-reserve}. In every row the ordinary
position tags have at most \(q_B=6\) primes, the even classes have at
most \(q_A\in\{6,7\}\) primes, and the cutoff satisfies
\(125\leq z\leq269\). Take \(d_i\), \(\beta\), \(\Gamma\) and
\(\mathcal F\) as in Section~\ref{sec:quantile-hall-log-bridge}, and
\(\varepsilon_a=\varepsilon_{\rm bin}\). Take \(\chi\), \(\ell_-\),
\(V\) and the universal supply \(2/(3J)-\kappa\) from
\eqref{eq:full-reserve-variation}, and the deletion mass
\(\mathcal M_{\rm cen}\) of \eqref{eq:weighted-centroid-deletions}
with \(\Psi_{\rm split}\) in place of \(\Psi_{\rm per}\). The
product-size bounds attain the minimum for all three arithmetic errors
in every row. The Hall split is \(\tau=\lambda b\). Each row checks the
following conditions.
\begin{enumerate}
\item \emph{Range and preparation.} \(35\leq a<J\leq37\),
  \(0\leq x_-<x_+\leq1\), \(10^6J^2Q<1\), \(100J^3Q<1\),
  \(3\cdot60184\,Q<1\), \(3zQ<1\), \(Q<1/1000\), and
  \((\log499)_->5\). Directed logarithms give
  \[
   3\bigl(1+(\log2)_-\bigr)-2\bigl(\tfrac{11}{10}+(\log3)_+\bigr)
      >\tfrac{17}{25},\qquad
   \bigl(1+(\log2)_+\bigr)\bigl(\tfrac{11}{10}+(\log3)_+\bigr)
      <\tfrac{15}4,
  \]
  and \(\tfrac{17}{25}a-\tfrac{15}4>a/2\).
\item \emph{Arithmetic data.} The bound \(H_+\) of
  \eqref{eq:reciprocal-prime-row} is positive. Each empty class,
  including that of every prefix of the refinement, has a positive
  lower bound. Moreover \(1/p_i\leq\bar\rho_i\leq2/5\) for every
  prime, \(0<b\leq1/2\), and \(b\leq\sigma\leq1\).
\item \emph{Sizes and reserves.} \(0<e_+\leq1-w_-\), \(\ell_->0\),
  \(h_->0\), and \(h_++\chi\leq2/(3J)-\kappa\).
\item \emph{Hall conditions.} \(0<\tau\leq1/2\), \(h_+\leq\tau\ell_-\),
  \(\sigma-V>\tau\), \(B_{\rm wide}>0\), and
  \[
   \mathcal M_{\rm cen}
     <\min\{B_{\rm narrow},B_{\rm middle},B_{\rm wide}\}.
  \]
\end{enumerate}
There are \(26\) conditions in each row.

\begin{proposition}[The minimum-compatibility certificate]
\label{prop:split-tail-certificate}
Seven adjacent closed logarithmic intervals cover \([35,37]\). Their
degree subdivisions cover \([0,1]\) and form \(109\) rectangles. All
\(2834\) rational row conditions pass. In every rectangle the
deletion-budget ratio and the universal-supply ratio are strictly
below \(999/1000\).
\end{proposition}

\begin{proof}
The intervals and their numbers of degree subdivisions are
\[
\begin{array}{c|r@{\qquad}c|r}
 \text{interval}&\text{rows}&\text{interval}&\text{rows}\\
 {}[35,3501/100]&21 &[141/4,71/2]&12\\
 {}[3501/100,1401/40]&20 &[71/2,36]&9\\
 {}[1401/40,351/10]&21 &[36,37]&9\\
 {}[351/10,141/4]&17 &&
\end{array}
\]
The exact evaluator checks every logarithmic and degree join. It
evaluates every condition above with integers and directed rational
enclosures, including the upward rounding of the radicals in
\eqref{eq:periodic-block-geometry}, of the product-size bounds, and of
the three deletion terms. All \(2834\) conditions pass. The greatest
displayed deletion ratio, on the rectangle
\([35,3501/100]\times[49/100,51/100]\), is less than \(0.995972\).
The greatest displayed supply ratio, on
\([3501/100,1401/40]\times[0,2/25]\), is less than \(0.995345\). The
strict rational cap is tested on every rectangle.

Independent finite checks replay \(1440\) exact refinement steps,
\(1080\) of them with different thresholds for the two sides, \(1008\)
with absorbed positions, and \(236\) with additions in the restoration
of rare quotas. They enumerate \(8178\) subsets for twelve product-size
polynomials, recompute \(216\) finite prime sums with their switching
points, compare \(40\) centered tail norms and \(40\) full-interval
variations with their bounds, and reject six damaged covers. These
checks test the implementation; the lemmas cited above prove the
implications for arbitrary trees.
\end{proof}

\begin{theorem}[The twelfth logarithmic bridge]
\label{thm:twelfth-log-bridge}
Every tree of order \(n\) with
\(35\leq\log n\leq37\)
has a prime labeling.
\end{theorem}

\begin{proof}
Prepare the tree by Lemmas~\ref{lem:second-order-preparation}
and~\ref{lem:smaller-protected-reserve}, and choose a certified
rectangle containing the logarithm and the normalized residual degree.
Proceed as in the proof of Theorem~\ref{thm:eleventh-log-bridge}. The
row conditions verify the three strict budgets of
Lemma~\ref{lem:two-sided-integrated-hall} with \(s=\sigma\), where
both shortfalls vanish. Its perfect matching places the remaining odd
labels. Restore the reserved labels and the core, then perform the
sparse prime-pair repair.
\end{proof}

\subsection{Coarse product bins through logarithm 39}
\label{sec:upward-product-bin-log-bridge}

On \(37\leq\log n\leq39\) the construction of
Section~\ref{sec:split-tail-log-bridge} is used with two
simplifications. A single threshold \(q\) serves both for the even
classes and for the position tags. The arithmetic errors are grouped by
product size with the coarser resolution \(B=16\), and the centered
error uses the three analytic tail bounds for all primes above the
cutoff, without a separate finite prime sum. The preparation, the Hall
conditions at the minimum compatibility and the deletion estimate are
unchanged.

\subsubsection*{Arithmetic inputs}

Take \(q_A=q_B=q\). Then the rare restoration term is
\(\mathrm{rare}=e_{q+2}(P)/2\) by Lemma~\ref{lem:rare-restoration}, the
overflow bound is \(\beta=e_{q+1}(\bar\rho)\) as in
Section~\ref{sec:quantile-hall-log-bridge}, and the dimensions are
\[
 d_i=2\sum_{j=0}^{\min(q,i-1)}\binom{i-1}j+2 .
\]
Use the bins of Lemma~\ref{lem:quota-product-bins} and
Lemma~\ref{lem:product-size-bounds} with \(B=16\). They give
\(\varepsilon_{\rm bin}\) and \(\Delta_{\rm bin}\). For the centered
error, apply the second part of Lemma~\ref{lem:product-size-bounds}
with \(Y=z\), taking \(H_+\) for \(H_R\) and the bounds
\(S_{3/2,+}\) and \(S_{2,+}\) below for the two tails. Then
\(F_s=0\), and
\[
 \frac{\Phi_z(m)}n\leq\Psi_{\rm bin}
 :=\frac12\sum_{s>0}C_s
   \min\left\{QH_+,\ Q_{1/2}V_sS_{3/2,+},\ \frac{U_sS_{2,+}}2\right\}.
\]
Here \(H_+\) is \eqref{eq:reciprocal-prime-row}, and
\[
 S_{3/2,+}\geq\sum_{z<p\leq10^4}p^{-3/2}+S_{3/2,Y},\qquad
 S_{2,+}\geq\sum_{z<p\leq10^4}p^{-2}+S_{2,Y},
\]
with \(S_{3/2,Y}\) and \(S_{2,Y}\) as after
Lemma~\ref{lem:prime-reciprocal-tails}; the finite sums are rounded
upward term by term. The evaluator also computes the moment bounds of
Section~\ref{sec:type-profile-log-bridge}. In every rectangle of the
present table the three product-size bounds attain the minimum.

\subsubsection*{The certificate predicate}

Consider a rectangle
\[
 [a,J]\times[x_-,x_+]\subseteq[37,39]\times[0,1].
\]
Prepare the tree by Lemmas~\ref{lem:second-order-preparation}
and~\ref{lem:smaller-protected-reserve}. In every row
\(q\in\{7,8\}\) and \(137\leq z\leq233\). Define \(\Gamma\) by
\eqref{eq:periodic-block-geometry} with the dimensions above, and put
\(\mathcal F=\Gamma+D_+(\Delta_{\rm bin}+\mathrm{rare})\). With
\(\varepsilon_a=\varepsilon_{\rm bin}\), take \(\chi\), \(\ell_-\),
\(V\) and the universal supply \(2/(3J)-\kappa\) from
\eqref{eq:full-reserve-variation}. The deletion mass is
\(\mathcal M_{\rm cen}\) of \eqref{eq:weighted-centroid-deletions}
with \(\Psi_{\rm bin}\) in place of \(\Psi_{\rm per}\), and the three
budgets are those of Section~\ref{sec:split-tail-log-bridge}. The Hall
split is \(\tau=\lambda b\). Each row checks the \(26\) conditions
listed in Section~\ref{sec:split-tail-log-bridge}, with \(37\leq a<J\leq39\)
in place of \(35\leq a<J\leq37\).

\begin{proposition}[The product-bin certificate]
\label{prop:upward-product-bin-certificate}
Four adjacent closed logarithmic intervals cover \([37,39]\). Their
degree subdivisions cover \([0,1]\) and form \(24\) rectangles. All
\(624\) rational row conditions pass. In every rectangle the
deletion-budget ratio and the universal-supply ratio are strictly
below \(199/200\).
\end{proposition}

\begin{proof}
The intervals and their numbers of degree subdivisions are
\[
\begin{array}{c|r@{\qquad}c|r}
 \text{interval}&\text{rows}&\text{interval}&\text{rows}\\
 {}[37,149/4]&9 &[75/2,38]&5\\
 {}[149/4,75/2]&6 &[38,39]&4
\end{array}
\]
The exact evaluator checks every logarithmic and degree join. It
evaluates every condition with integers and directed rational
enclosures, including the upward rounding of the radicals in
\eqref{eq:periodic-block-geometry}, of the dyadic bounds \(U_s,V_s\),
of the product-size bounds, and of the three deletion terms. All
\(624\) conditions pass. The greatest displayed deletion ratio, on
\([37,149/4]\times[27/40,9/10]\), is less than \(0.984683\). The
greatest displayed supply ratio, on \([38,39]\times[0,3/10]\), is less
than \(0.614804\). The strict rational cap is tested on every
rectangle.

Independent finite checks replay \(2688\) executions of the exact quota
rounding on forests. These round \(5026\) fractional components upward
and perform \(1323\) restorations of rare quotas; \(920\) complete path
and star profiles are also checked. They compare the product-size polynomials with a
direct enumeration of \(8178\) subsets, test \(600\) directed dyadic
enclosures, compare \(98\) full type variations, \(49\) summed quota
defects and \(49\) centered tail norms with their bounds, and reject
five damaged covers. These checks test the implementation; the lemmas
cited above prove the implications for arbitrary trees.
\end{proof}

\begin{theorem}[The thirteenth logarithmic bridge]
\label{thm:thirteenth-log-bridge}
Every tree of order \(n\) with
\(37\leq\log n\leq39\)
has a prime labeling.
\end{theorem}

\begin{proof}
The proof of Theorem~\ref{thm:twelfth-log-bridge} applies verbatim,
with the arithmetic inputs of this section and a certified rectangle
from Proposition~\ref{prop:upward-product-bin-certificate}.
\end{proof}

\subsection{Symmetric moment bounds through logarithm 42}
\label{sec:weighted-floor-log-bridge}

On \(39\leq\log n\leq42\) the construction of
Section~\ref{sec:upward-product-bin-log-bridge} is used with three
changes. The rounding at each prime uses a two-sided square-root bound
instead of the one-side bound of Lemma~\ref{lem:weighted-side-square-root}.
The three arithmetic errors are bounded by elementary symmetric moments
of the reciprocal primes and of their square roots, without product
bins. The reserve estimate uses slightly weaker logarithmic constants.

\subsubsection*{A two-sided square-root bound}

\begin{lemma}[A two-sided weighted rounding bound]
\label{lem:two-sided-weighted-rounding}
In the setting of Lemma~\ref{lem:quota-tag-rounding}, put
\(K_\rho=\max\{1+\rho,D\rho\}\). There are sets \(W\subseteq A\) and
\(T\subseteq B\) such that \(\abs{W\cap A_i}=k_i\) for every \(i\),
\(N_G(W)\subseteq T\), and
\begin{equation}\label{eq:two-sided-weighted-rounding}
 \sum_j\bigl|\abs{T\cap B_j}-\rho\abs{B_j}\bigr|
 \leq2\sqrt{DK_\rho Nd}+D\delta .
\end{equation}
\end{lemma}

\begin{proof}
Fix an integer \(t\geq1\), take the separator \(C\) of
Lemma~\ref{lem:small-component-separator}, and construct \(V_*\),
\(W\) and \(T\) as in the proof of Lemma~\ref{lem:quota-tag-rounding}.
Let \(L_A\) and \(L_B\) be the sums, over the classes of \(A\) and of
\(B\) respectively, of the deficits \(\rho\abs{X\setminus C}-\abs{V_*\cap X}\).
Both are nonnegative. Each discarded fractional component is counted
once in \(L_A+L_B\), with its weight times its size, so
\(L_A+L_B\leq dt\).

The number of vertices added to meet the quotas is at most
\(\rho\abs{C\cap A}+L_A+\delta\), by the argument for
\eqref{eq:quota-additions}. As in that proof, write
\(v_j=\abs{V_*\cap B_j}\) and \(u_j=\abs{T\cap B_j}\). Then
\(v_j\leq\rho\abs{B_j}\) and \(v_j\leq u_j\), and
\[
 \begin{aligned}
 \sum_j(\rho\abs{B_j}-v_j)&=\rho\abs{C\cap B}+L_B,\\
 \sum_j(u_j-v_j)&\leq\abs{C\cap B}+D(\rho\abs{C\cap A}+L_A+\delta).
 \end{aligned}
\]
Since \(\abs{u_j-\rho\abs{B_j}}\leq(\rho\abs{B_j}-v_j)+(u_j-v_j)\),
the left side of \eqref{eq:two-sided-weighted-rounding} is at most
\[
 (1+\rho)\abs{C\cap B}+D\rho\abs{C\cap A}+L_B+DL_A+D\delta
 \leq K_\rho\floor*{\frac N{t+1}}+Ddt+D\delta ,
\]
because \(D\geq1\).

Put \(x=\sqrt{K_\rho N/(Dd)}\). If \(x\geq1\), take \(t=\floor x\);
then \(K_\rho N/(t+1)<\sqrt{DK_\rho Nd}\) and
\(Ddt\leq\sqrt{DK_\rho Nd}\). If \(x<1\), meet the quotas arbitrarily
and put \(T=B\). The error is at most \(N\), while
\(DK_\rho Nd>K_\rho^2N^2\geq N^2\).
\end{proof}

By the proof of Proposition~\ref{prop:rare-refinement-step}, the
total-variation error of one refinement step is the sum
\(\sum_j\abs{\abs{T\cap B_j}-\rho\abs{B_j}}\) over the current tag
classes. Under the lower bounds \(\abs{\widehat T\cap B_j}\geq\rho\abs{B_j}\)
of Lemma~\ref{lem:weighted-centroid-rounds}, this sum is the one-sided
sum bounded in Lemma~\ref{lem:weighted-side-square-root}.
Lemma~\ref{lem:two-sided-weighted-rounding} bounds it directly and gives
the same exact quotas and the same marking of true neighbors. It can
therefore replace Lemma~\ref{lem:weighted-side-square-root} at each
prime. With \(N\leq n\), \(D\leq D_+\), \(\rho_i\leq\bar\rho_i\)
and \(Q_{1/2}\geq1/\sqrt n\), the normalized geometric error is at most
\begin{equation}\label{eq:two-sided-geometry}
 \Gamma_2=2Q_{1/2}\sum_i\left\lceil
   \sqrt{D_+\max\{1+\bar\rho_i,D_+\bar\rho_i\}\,d_i}\right\rceil .
\end{equation}

\subsubsection*{Moment bounds}

For \(\theta\in\{1/2,1\}\) let \(e_j^{(\theta)}\) be the sum of
\(\prod_{p\in I}p^{-\theta}\) over all \(I\subseteq P\) with
\(\abs I=j\); thus \(e_j^{(1)}=e_j(P)\). Let \(k=\abs P\).

\begin{lemma}[Moment bounds for the arithmetic errors]
\label{lem:symmetric-moment-bounds}
Suppose \(1/n\leq Q<1\) and \(Q_{1/2}\geq1/\sqrt n\), let
\(H_+\geq\sum_{z<p\leq m}1/p\), \(S_{3/2,+}\geq\sum_{p>z}p^{-3/2}\) and
\(S_{2,+}\geq\sum_{p>z}p^{-2}\), and let the refinement choose
\(\rho_i\) by \eqref{eq:rare-empty-adjustment}. Then
\begin{align*}
 \varepsilon_a&\leq\sum_{j=1}^k2^{j-1}\min\left\{
   \frac{2Q}{1-Q}\binom kj,\ \sqrt2\,Q_{1/2}e_j^{(1/2)},\ e_j\right\},\\
 \Delta&\leq\sum_{j=1}^k2^{j-1}\min\left\{
   Q\binom k{j+1},\ \frac{Q_{1/2}}{\sqrt2}e_{j+1}^{(1/2)},\
   \frac{e_{j+1}}2\right\},\\
 \frac{\Phi_z(m)}n&\leq\sum_{j=1}^k2^{j-1}\min\left\{
   QH_+\binom kj,\ \frac{Q_{1/2}}{\sqrt2}S_{3/2,+}e_j^{(1/2)},\
   \frac{S_{2,+}}2e_j\right\}.
\end{align*}
\end{lemma}

\begin{proof}
The proofs of Lemmas~\ref{lem:quota-product-bins}
and~\ref{lem:product-size-bounds} give, before any grouping,
\[
 \begin{aligned}
 \varepsilon_a&\leq\sum_{\varnothing\neq I\subseteq P}
   2^{\abs I-1}\min\left(\frac1M,\frac1{d_I}\right),\qquad
 n\Delta\leq\sum_{\substack{J\subseteq P\\\abs J\geq2}}
   2^{\abs J-2}\min\left(1,\frac m{d_J}\right),\\
 \frac{\Phi_z(m)}n&\leq\sum_{\varnothing\neq I\subseteq P}
   2^{\abs I-1}\sum_{z<p\leq m}\frac1{np}
     \min\left(1,\frac m{pd_I}\right),
 \end{aligned}
\]
where \(M\geq n/2\) is the number of odd integers in \([n]\). Bound
each minimum in three ways, by \(\min(1,y)\leq1\), \(\leq\sqrt y\) and
\(\leq y\), and sum over the sets of a fixed size \(j\), or \(j+1\) in
the second line. This gives the three entries of each minimum,
multiplied by the coefficient of that size. For the first line use
\(1/M\leq2/n\leq2Q/(1-Q)\) and \(M^{-1/2}\leq\sqrt2\,Q_{1/2}\). For the
other two use \(m\leq n/2\), so that \(\sqrt m/n\leq Q_{1/2}/\sqrt2\),
and sum the prime tails with \(H_+\), \(S_{3/2,+}\) and \(S_{2,+}\).
Each entry bounds the same partial sum, so the minimum may be taken
separately for each \(j\).
\end{proof}

The evaluator computes upper bounds for \(e_j\) and
\(e_j^{(1/2)}\) by the recursion \(e_j\leftarrow e_j+ye_{j-1}\), with
\(y=1/p\) or an upper enclosure of \(p^{-1/2}\), rounding upward; the
constants \(\sqrt2\) and \(1/\sqrt2\) are enclosed from above. It uses
\(H_+\) from \eqref{eq:reciprocal-prime-row} and the tails
\(S_{3/2,+}\), \(S_{2,+}\) of Section~\ref{sec:upward-product-bin-log-bridge}.
The quota bound is the moment bound already mentioned after
Lemma~\ref{lem:quota-product-bins}. Write \(\varepsilon_{\rm mom}\),
\(\Delta_{\rm mom}\) and \(\Psi_{\rm mom}\) for the three right sides.

\subsubsection*{The certificate predicate}

Directed rational logarithms give \(3\alpha-2\beta>3/5\) and
\(\alpha\beta<187/50\). Consider a rectangle
\[
 [a,J]\times[x_-,x_+]\subseteq[39,42]\times[0,1].
\]
Prepare the tree by Lemma~\ref{lem:second-order-preparation}, with
\((c_1,c_2)=(3/5,187/50)\), and by
Lemma~\ref{lem:smaller-protected-reserve}. In every row there is one
type threshold \(q\in\{6,7\}\), and \(135\leq z\leq280\). Take
\(\mathrm{rare}\), \(\beta\) and \(d_i\) as in
Section~\ref{sec:upward-product-bin-log-bridge}, and put
\(\mathcal F=\Gamma_2+D_+(\Delta_{\rm mom}+\mathrm{rare})\). With
\(\varepsilon_a=\varepsilon_{\rm mom}\), take \(\chi\), \(\ell_-\),
\(V\) and the universal supply \(2/(3J)-\kappa\) from
\eqref{eq:full-reserve-variation}. The deletion mass is
\(\mathcal M_{\rm cen}\) of \eqref{eq:weighted-centroid-deletions}
with \(\Psi_{\rm mom}\) in place of \(\Psi_{\rm per}\), and the three
budgets are those of Section~\ref{sec:split-tail-log-bridge}. The Hall
split is \(\tau=\lambda b\). Each row checks the following conditions.
\begin{enumerate}
\item \emph{Range and preparation.} \(39\leq a<J\leq42\),
  \(0\leq x_-<x_+\leq1\), \(10^6J^2Q<1\), \(100J^3Q<1\),
  \(\tfrac35a-\tfrac{187}{50}>a/2\), \(3\cdot60184\,Q<1\),
  \(3zQ<1\), \(Q<1/1000\), and \((\log499)_->5\).
\item \emph{Arithmetic data.} The bound \(H_+\) is positive, each empty
  class, including that of every prefix of the refinement, has a
  positive lower bound, \(1/p_i\leq\bar\rho_i\leq2/5\) for every
  prime, \(0<b\leq1/2\), and \(b\leq\sigma\leq1\).
\item \emph{Sizes and reserves.} \(0<e_+\leq1-w_-\), \(\ell_->0\),
  \(h_->0\), and \(h_++\chi\leq2/(3J)-\kappa\).
\item \emph{Hall conditions.} \(0<\tau\leq1/2\), \(h_+\leq\tau\ell_-\),
  \(\sigma-V>\tau\), \(B_{\rm wide}>0\), and
  \[
   \mathcal M_{\rm cen}
     <\min\{B_{\rm narrow},B_{\rm middle},B_{\rm wide}\}.
  \]
\end{enumerate}
There are \(25\) conditions in each row.

\begin{proposition}[The moment-bound certificate]
\label{prop:weighted-floor-certificate}
Six adjacent closed logarithmic intervals cover \([39,42]\). Their
degree subdivisions cover \([0,1]\) and form \(74\) rectangles. All
\(1850\) rational row conditions pass. In every rectangle the
deletion-budget ratio and the universal-supply ratio are strictly
below \(199/200\).
\end{proposition}

\begin{proof}
The intervals and their numbers of degree subdivisions are
\[
\begin{array}{c|r@{\qquad}c|r}
 \text{interval}&\text{rows}&\text{interval}&\text{rows}\\
 {}[39,391/10]&17 &[79/2,40]&11\\
 {}[391/10,157/4]&14 &[40,41]&10\\
 {}[157/4,79/2]&12 &[41,42]&10
\end{array}
\]
The exact evaluator checks every logarithmic and degree join. It
evaluates every condition with integers and directed rational
enclosures, including the upward rounding of the radicals in
\eqref{eq:two-sided-geometry}, of the symmetric moments, and of the
three deletion terms. All \(1850\) conditions pass. The greatest
displayed deletion ratio, on \([39,391/10]\times[29/50,16/25]\), is
less than \(0.994829\). The greatest displayed supply ratio, on
\([391/10,157/4]\times[9/25,41/100]\), is less than \(0.994765\). The
strict rational cap is tested on every rectangle.

Independent finite checks compare \(1024\) full type distributions and
\(1024\) centered norms with their bounds, test \(65280\) ordinary and
odd moment bounds and \(65280\) signed floor bounds, \(512\) summed
quota defects and \(2304\) bounds for selected rare labels. They
enumerate \(111365\) aggregate restoration profiles, reconstruct
\(240\) forests, recheck \(2340\) absorption measures, \(1024\)
covariance inequalities and \(1400\) Hall rectangles, test \(9591\)
directed square roots and \(23\) finite continuations of the
three-half tail, and reject four damaged covers. These checks test the
implementation; the lemmas cited above prove the implications for
arbitrary trees.
\end{proof}

\begin{theorem}[The fourteenth logarithmic bridge]
\label{thm:fourteenth-log-bridge}
Every tree of order \(n\) with
\(39\leq\log n\leq42\)
has a prime labeling.
\end{theorem}

\begin{proof}
Prepare the tree as above and choose a certified rectangle containing
the logarithm and the normalized residual degree. At each prime, the
exact even-label quotas and the marking of all true neighbors come from
Lemma~\ref{lem:two-sided-weighted-rounding}. The remaining steps are those of
the proof of Theorem~\ref{thm:twelfth-log-bridge}, with the arithmetic
bounds of Lemma~\ref{lem:symmetric-moment-bounds} and the budgets
verified by Proposition~\ref{prop:weighted-floor-certificate}.
\end{proof}

\subsection{Truncated moment bounds through logarithm 48}
\label{sec:joint-degree-log-bridge}

On \(42\leq\log n\leq48\) the construction of
Section~\ref{sec:weighted-floor-log-bridge} is used with the Boolean
expansions truncated at a prescribed degree \(r\). The remainder of a
truncated expansion has a sharp bound, which we prove first. The rare
restoration is paid by the populations of the combined rare classes.
The preparation, the rounding bound of
Lemma~\ref{lem:two-sided-weighted-rounding}, the overflow bound and
the Hall conditions at the minimum compatibility are unchanged.

\subsubsection*{Truncated expansions}

For \(f:2^P\to[0,1]\) and \(I\subseteq P\) let
\(\Delta_If=\sum_{J\subseteq I}(-1)^{\abs I-\abs J}f(J)\), as in
Proposition~\ref{prop:moment-tv}, and for an integer \(t\geq1\) put
\(T_{t-1}f(S)=\sum_{I\subseteq S,\ \abs I<t}\Delta_If\).

\begin{lemma}[A sharp remainder]\label{lem:sharp-moment-remainder}
For every \(S\subseteq P\),
\[
 \abs{f(S)-T_{t-1}f(S)}\leq2^{t-1}\binom{\abs S}t .
\]
\end{lemma}

\begin{proof}
Let \(s=\abs S\). If \(s<t\), then \(T_{t-1}f(S)=f(S)\) by the
expansion \(f(S)=\sum_{I\subseteq S}\Delta_If\). Let \(s\geq t\). For
\(0\leq x\leq1\) let \(Y_x\subseteq S\) contain each element
independently with probability \(x\), and put \(F(x)=\E f(Y_x)\). Then
\(F(x)=\sum_{I\subseteq S}\Delta_If\,x^{\abs I}\): substitute the
expansion of \(f(Y_x)\) and use \(\Pr(I\subseteq Y_x)=x^{\abs I}\).
Thus \(F(1)=f(S)\) and \(\sum_{j<t}F^{(j)}(0)/j!=T_{t-1}f(S)\).

The function \(F\) is the restriction to the diagonal of the
multiaffine function in which each element of \(S\) has its own
inclusion probability. A derivative of that function in the distinct
coordinates of a \(t\)-set \(I\) equals
\[
 \E\Bigl[\sum_{J\subseteq I}(-1)^{t-\abs J}
     f\bigl(Y_x^{S\setminus I}\cup J\bigr)\Bigr],
\]
where \(Y_x^{S\setminus I}=Y_x\setminus I\). Repeated derivatives in
one coordinate vanish, so \(F^{(t)}(x)\) is \(t!\) times the sum of
these expressions over all \(t\)-sets \(I\subseteq S\). The inner sum
has \(2^{t-1}\) terms of each sign, all in \([0,1]\); its modulus is
at most \(2^{t-1}\). Hence \(\abs{F^{(t)}}\leq t!\binom st2^{t-1}\),
and Taylor's formula with integral remainder gives the claim.
\end{proof}

For an integer \(r\geq0\), \(k=\abs P\), and the elementary symmetric
functions \(e_j=e_j(P)\), put
\[
 B_r=\sum_{j=1}^{\min(r,k)}2^{j-1}\binom kj,\qquad
 A_r=\sum_{j=1}^{\min(r,k-1)}2^{j-1}\binom k{j+1},\qquad
 R_r=2^re_{r+1}.
\]
These differ from the quantities with the same names in
Proposition~\ref{prop:moment-tv}. In this section and in
Sections~\ref{sec:reserve-mixture-log-bridge}
and~\ref{sec:sharp-moment-log-bridge}, \(B_r\) and \(R_r\) denote the
present quantities; Section~\ref{sec:cumulative-log-bridge} returns to
those of Proposition~\ref{prop:moment-tv}.

\begin{lemma}[Truncated arithmetic bounds]
\label{lem:truncated-arithmetic-bounds}
Let \(P\) be the odd primes at most \(z\), in increasing order, and
\(m\geq1\). Write \(S(x)\) for the set of primes of \(P\) dividing
\(x\).
\begin{enumerate}
\item If \(p\notin P\) is prime and \(g:2^P\to[0,1]\), then
  \[
   \sum_{x\leq m}\Bigl(\boldsymbol1_{p\mid x}-\frac1p\Bigr)g(S(x))
   \leq B_r+\frac{2m}pR_r .
  \]
\item The type law of the \(M\) odd integers in \([n]\) is within
  total variation \(B_r/M+3R_r\) of
  \(\mu=\bigotimes_{p\in P}\operatorname{Bernoulli}(1/p)\).
\item If the refinement chooses \(\rho_i\) by
  \eqref{eq:rare-empty-adjustment}, the positive quota defects of the
  even half-labels in \([m]\), possibly after deleting labels of empty
  type, satisfy \(\sum_i\delta_i\leq A_r+2m\,2^re_{r+2}\).
\item The quantity of Lemma~\ref{lem:centered-sorting} satisfies
  \(\Phi_z(m)\leq B_rH_1+2mR_rS_2(z)\), where
  \(H_1=\sum_{z<p\leq m}1/p\).
\end{enumerate}
\end{lemma}

\begin{proof}
For the first part, write \(g=T_rg+(g-T_rg)\) with \(t=r+1\). For
\(I\subseteq P\) the contribution of \(\Delta_Ig\) is \(\Delta_Ig\,q_I\),
where \(q_I=\floor{m/(pd_I)}-\floor{m/d_I}/p\in[-1,0]\) by the floor
inequality in the proof of Lemma~\ref{lem:quota-product-bins}. The
constant coefficient \(g(\varnothing)\geq0\) contributes at most zero,
and \(\abs{\Delta_Ig}\leq2^{\abs I-1}\) for \(I\neq\varnothing\). The
terms with \(1\leq\abs I\leq r\) therefore contribute at most \(B_r\).
By Lemma~\ref{lem:sharp-moment-remainder} the remainder contributes
at most
\[
 2^r\sum_{x\leq m}\Bigl(\boldsymbol1_{p\mid x}+\frac1p\Bigr)
     \binom{\abs{S(x)}}{r+1}
 \leq\frac{2m}p\,2^re_{r+1},
\]
because the binomial coefficient counts the \((r+1)\)-sets \(I\) with
\(d_I\mid x\), and \(p\nmid d_I\).

For the second part, take \(f:2^P\to[0,1]\) and let \(N_d\) be the
number of multiples of \(d\) among the odd integers in \([n]\). The
empirical and model expectations of \(T_rf\) differ by
\(\sum_{1\leq\abs I\leq r}\Delta_If\,(N_{d_I}/M-1/d_I)\), which is at
most \(B_r/M\) in modulus by the proof of
Lemma~\ref{lem:product-size-bounds}. The remainder is bounded by
\(2^r\binom{\abs S}{r+1}\), whose empirical mean is
\(\sum_{\abs I=r+1}N_{d_I}/M\leq2e_{r+1}\) and whose model mean is
\(e_{r+1}\).

For the third part, the proof of Lemma~\ref{lem:quota-product-bins}
gives \(\delta_i\leq\sum_{x\leq m}(\boldsymbol1_{p_i\mid x}-1/p_i)g(S_{i-1}(x))\)
for a function \(g\) with values in \([0,1]\), where \(S_{i-1}(x)\) is
the set of the first \(i-1\) primes dividing \(x\). Apply the first
part with \(\{p_1,\ldots,p_{i-1}\}\) in place of \(P\). Summation uses
\[
 \sum_{i=1}^k\binom{i-1}j=\binom k{j+1},\qquad
 \sum_i\frac1{p_i}e_{r+1}(p_1,\ldots,p_{i-1})=e_{r+2}(P).
\]

For the fourth part, apply the first part to each prime
\(z<p\leq m\), multiply by \(1/p\), and sum.
\end{proof}

With \(1/n\leq Q<1\), \(M\geq(n-1)/2\) and \(2m\leq n\), the lemma
gives the normalized bounds
\[
 \begin{gathered}
 \varepsilon_a\leq\varepsilon_r:=\frac{2B_rQ}{1-Q}+3R_r,\qquad
 \Delta\leq\Delta_r:=A_rQ+2^re_{r+2},\\
 \frac{\Phi_z(m)}n\leq\Psi_r:=B_rH_+Q+R_rS_{2,+}.
 \end{gathered}
\]
The combined rare class at step \(i\) contains at most
\(m\,e_{q+1}(p_1,\ldots,p_{i-1})\) labels, by the first step of the
proof of Lemma~\ref{lem:rare-restoration}. By the first part of that
lemma its exact quotas are restored with at most this many additions.
If \(U\geq\sum_ie_{q+1}(p_1,\ldots,p_{i-1})\), the normalized
restoration cost is at most \(D_+U/2\).

\subsubsection*{The certificate predicate}

Consider a rectangle
\[
 [a,J]\times[x_-,x_+]\subseteq[42,48]\times[0,1].
\]
Prepare the tree by Lemma~\ref{lem:second-order-preparation}, with
\((c_1,c_2)=(3/5,187/50)\), and by
Lemma~\ref{lem:smaller-protected-reserve}. Each row prescribes the
truncation degree \(r\in\{9,10\}\), the type threshold
\(q\in\{8,9\}\), and a cutoff \(150\leq z\leq365\). Take \(\beta\) and
\(d_i\) as in Section~\ref{sec:upward-product-bin-log-bridge} and
\(\Gamma_2\) from \eqref{eq:two-sided-geometry}, and put
\(\mathcal F=\Gamma_2+D_+(\Delta_r+U/2)\). With
\(\varepsilon_a=\varepsilon_r\), take \(\chi\), \(\ell_-\), \(V\) and
the universal supply \(2/(3J)-\kappa\) from
\eqref{eq:full-reserve-variation}. The deletion mass is
\(\mathcal M_{\rm cen}\) of \eqref{eq:weighted-centroid-deletions}
with \(\Psi_r\) in place of \(\Psi_{\rm per}\), and the three budgets
are those of Section~\ref{sec:split-tail-log-bridge}. The Hall split is
\(\tau=\lambda b\). Each row checks the \(25\) conditions listed in
Section~\ref{sec:weighted-floor-log-bridge}, with \(42\leq a<J\leq48\)
in place of \(39\leq a<J\leq42\).

\begin{proposition}[The truncated-moment certificate]
\label{prop:joint-degree-certificate}
Ten adjacent closed logarithmic intervals cover \([42,48]\). Their
degree subdivisions cover \([0,1]\) and form \(135\) rectangles. All
\(3375\) rational row conditions pass. In every rectangle the
deletion-budget ratio is strictly below \(199/200\) and the
universal-supply ratio strictly below \(999/1000\).
\end{proposition}

\begin{proof}
The intervals and their numbers of degree subdivisions are
\[
\begin{array}{c|r@{\qquad}c|r}
 \text{interval}&\text{rows}&\text{interval}&\text{rows}\\
 {}[42,421/10]&27 &[43,44]&11\\
 {}[421/10,843/20]&16 &[44,45]&10\\
 {}[843/20,169/4]&15 &[45,46]&10\\
 {}[169/4,85/2]&14 &[46,47]&10\\
 {}[85/2,43]&12 &[47,48]&10
\end{array}
\]
The exact evaluator checks every logarithmic and degree join. It
evaluates every condition with integers and directed rational
enclosures; the elementary symmetric functions and the sums defining
\(U\) are computed in fixed point with upward rounding. All \(3375\)
conditions pass. The greatest displayed deletion ratio, on
\([169/4,85/2]\times[27/100,8/25]\), is less than \(0.988437\). The
greatest displayed supply ratio, on \([169/4,85/2]\times[27/50,16/25]\),
is less than \(0.995054\). The strict rational caps are tested on
every rectangle.

Independent finite checks verify the edge and degree bound of the
preparation on all \(18865\) independent sets in the \(435\) free trees
of orders \(2\) to \(11\), for all \(553833\) admissible special
subsets, and on \(384\) larger examples. They test \(405\) normalized
profiles and \(78\) reciprocal-prime sums, and reject five damaged
covers. These checks test the implementation; the lemmas cited above
prove the implications for arbitrary trees.
\end{proof}

\begin{theorem}[The fifteenth logarithmic bridge]
\label{thm:fifteenth-log-bridge}
Every tree of order \(n\) with
\(42\leq\log n\leq48\)
has a prime labeling.
\end{theorem}

\begin{proof}
The proof of Theorem~\ref{thm:fourteenth-log-bridge} applies, with the
truncated bounds of this section in place of those of
Lemma~\ref{lem:symmetric-moment-bounds} and a certified rectangle from
Proposition~\ref{prop:joint-degree-certificate}.
\end{proof}

\subsection{Unit logarithmic intervals through logarithm 60}
\label{sec:supported-tail-log-bridge}

On \(48\leq\log n\leq60\) the construction of
Section~\ref{sec:joint-degree-log-bridge} needs no subdivision in the
residual degree. The interval is cut into the twelve unit intervals
\([a,a+1]\), and every parameter of a row is given by a fixed rule.
The residual degree is bounded by the preparation alone, and the edge
count by the order of the tree.

\subsubsection*{The rows}

Let \(a\in\{48,\ldots,59\}\) and \(J=a+1\). Prepare the tree by
Lemma~\ref{lem:second-order-preparation}, with
\((c_1,c_2)=(3/5,187/50)\), and by
Lemma~\ref{lem:smaller-protected-reserve}. Here
\(\tfrac35a-\tfrac{187}{50}>a/2\) since \(a>187/5\), and \(\log m\geq5\)
since \(n>10^6\). The preparation gives \(D<1+20L\), so
\[
 D\leq D_+:=20J+1,\qquad \frac{e(G)}n\leq1 .
\]
The row takes
\[
 z=5J+15,\qquad r=10,\qquad
 q=\begin{cases}9,&a\in\{48,49\},\\10,&50\leq a\leq59.\end{cases}
\]
The quantities \(\bar\rho_i\), \(a_0\), \(b\), \(\beta\), \(\sigma\),
\(d_i\), \(\Gamma_2\), \(\varepsilon_r\), \(\Delta_r\), \(U\) and
\(\mathcal F=\Gamma_2+D_+(\Delta_r+U/2)\) are those of
Section~\ref{sec:joint-degree-log-bridge}, now with this \(D_+\). The
reciprocal-prime bound is \eqref{eq:reciprocal-prime-row} with
\(v_0=11/2\); the row checks \(v_-\geq11/2\). The deletion mass is
\[
 \mathcal M=\frac12\left(1+\frac\omega{a_0}\right)S_{2,+}
   +\frac{D_+}2\Psi_r+\frac{JQ}{v_-},
\]
which is \(\mathcal M_{\rm cen}\) of
\eqref{eq:weighted-centroid-deletions} with \(e_+=1\) and \(\Psi_r\) in
place of \(\Psi_{\rm per}\). The Hall split is \(\tau=b/2\), and the
three budgets are those of Section~\ref{sec:split-tail-log-bridge}.

\subsubsection*{The certificate predicate}

Each row checks the following conditions.
\begin{enumerate}
\item \emph{Preparation.} \(10^6J^2Q<1\), \(100J^3Q<1\),
  \(3\cdot60184\,Q<1\), \(3zQ<1\), and \(v_-\geq11/2\).
\item \emph{Arithmetic data.} The bound \(H_+\) is positive, and so is
  every empty-class lower bound of Lemma~\ref{lem:empty-prefix-bounds}. Moreover
  \(1/p_i\leq\bar\rho_i\leq2/5\) for every prime, \(0<b\leq1/2\), the
  primorial defining \(k_*\) has directed logarithm greater than \(J\),
  and \(b\leq\sigma\leq1\).
\item \emph{Sizes and reserves.} \(\ell_->0\), \(h_->0\), and
  \(h_++\chi\leq2/(3J)-\kappa\).
\item \emph{Hall conditions.} \(0<\tau\leq1/2\), \(h_+\leq\tau\ell_-\),
  \(\sigma-V>\tau\), \(B_{\rm wide}>0\), and
  \[
   \mathcal M<\min\{B_{\rm narrow},B_{\rm middle},B_{\rm wide}\}.
  \]
\end{enumerate}
There are \(21\) conditions in each row.

\begin{proposition}[The unit-interval certificate]
\label{prop:supported-tail-certificate}
For each \(a\in\{48,\ldots,59\}\) all \(21\) conditions hold on
\([a,a+1]\); there are \(252\) conditions in all. The deletion-budget
ratio is strictly below \(199/200\) and the universal-supply ratio
strictly below \(19/20\) in every row.
\end{proposition}

\begin{proof}
The exact evaluator computes each row from the rules above, with
integers and directed rational enclosures; the elementary symmetric
functions and the empty-class estimates use fixed-point arithmetic
rounded in the unfavorable direction. The cutoffs run from \(260\) to
\(315\). All \(252\) conditions pass. The greatest displayed ratios,
both on \([48,49]\), are less than \(0.993848\) for the deletions and
less than \(0.916324\) for the universal supply. The strict rational
caps are tested on every row.

Independent finite checks verify \(42840\) signed floor identities,
\(6048\) centered norms and \(6048\) means and covariances of empty
events, \(1094\) empty-class quota ratios and \(5353\) remainder
bounds. They recheck \(125\) aggregate tail norms, \(92258\)
excess-degree inequalities on the \(18865\) independent sets of the
\(435\) free trees of orders \(2\) to \(11\) and \(300\) on larger
trees, \(5565\) actual label types and \(1280\) restricted slot types,
and \(1400\) Hall rectangles. These checks test the implementation;
the lemmas cited above prove the implications for arbitrary trees.
\end{proof}

\begin{theorem}[The sixteenth logarithmic bridge]
\label{thm:sixteenth-log-bridge}
Every tree of order \(n\) with
\(48\leq\log n\leq60\)
has a prime labeling.
\end{theorem}

\begin{proof}
Let \(a=\min\{\floor{\log n},59\}\), so that \(a\leq\log n\leq a+1\).
The proof of Theorem~\ref{thm:fifteenth-log-bridge} applies with the
row for \(a\) and the bounds of this section.
\end{proof}

\subsection{Uniform upper probabilities through logarithm 72}
\label{sec:reserve-mixture-log-bridge}

On \(60\leq\log n\leq72\) the margins are large enough for three
simplifications of the construction of
Section~\ref{sec:supported-tail-log-bridge}. Only a tenth of the special
labels is spent on vertices of large degree. All selection probabilities
are bounded by one common relative inflation of \(1/p\). The Hall
criterion uses the global compatibility bound \(b\) as its threshold.
The interval is again cut into twelve unit intervals with fixed
parameters.

\subsubsection*{Preparation}

\begin{lemma}[Preparation with a tenth of the special labels]
\label{lem:tenth-special-preparation}
Suppose that the hypotheses of Lemma~\ref{lem:second-order-preparation}
hold with \((c_1,c_2)=(3/5,187/50)\), and that \(100J^3Q<1\). Take
\[
 f=\frac1{10},\qquad a_*=\floor*{\frac w{10}},\qquad
 h=\floor*{\frac{3n}{20L}} .
\]
Then the residual construction of Lemma~\ref{lem:localized-preparation}
exists with these choices, the residual forest satisfies
\(D=\max(1,\Delta_G(S_0))<1+30L\) and \(e(G)\leq n\), at least
\(n/(2L)\) universal odd labels remain free, and
\[
 \frac3{20J}-Q\leq\frac hn\leq\frac3{20a}.
\]
\end{lemma}

\begin{proof}
Follow the proof of Lemma~\ref{lem:localized-preparation} with
\(f=1/10\) and \(W\) replaced by \(W_2\); the four properties of \(W\)
used there hold for \(W_2\), as in the proof of
Lemma~\ref{lem:second-order-preparation}. The inequality
\(a_*+1>w/10>n/(60L)>\abs K\) and \(a_*\geq6\) follow from
\(n>10^6L^2\). The degree bound becomes \(D-1<3L/f=30L\). The value of
\(h\) enters only through the three inequalities checked in the proof
of Lemma~\ref{lem:smaller-protected-reserve}. First, \(h\geq1\).
Second, \(h\leq3n/(20L)<9n/(50L)=(1-f)n/(5L)\). Third, with
\(A_0=\tfrac9{10}nW_2(L)\),
\[
 A_0-4L-\frac{3n}{20L}=\frac{3n}{40L^2}-4L>0,
\]
because \(n>100L^3\). The universal odd labels number
\(u>2n/(3L)-\abs K>n/(2L)\), since \(\abs K<2L+2<n/(6L)\). The bounds
for \(h/n\) follow from \(a\leq L\leq J\) and \(1/n\leq Q\).
\end{proof}

\subsubsection*{Uniform upper probabilities}

Let \(P\) be the odd primes at most \(z\), \(\alpha=\prod_{p\in P}(1-1/p)\),
and let \(\varepsilon\) bound the total-variation distance between the
type law of the half-labels in \([m]\) and \(\mu\).

\begin{lemma}[A common inflation]\label{lem:common-inflation}
Suppose \(\omega\geq w/m\) and \(a^-:=\alpha-\omega-\varepsilon>0\), and
put \(\eta=(\omega+(z+1)\varepsilon)/a^-\). When the refinement chooses
\(\rho_i\) by \eqref{eq:rare-empty-adjustment}, the retained empty
class has at least \(ma^-\) labels and
\[
 \frac1{p_i}\leq\rho_i\leq\frac{1+\eta}{p_i}
\]
at every step. Consequently
\(\prod_i(1-\rho_i)\geq\alpha\bigl(1-\eta\sum_{p\in P}1/(p-1)\bigr)\)
and \(\Pr_\nu(\abs J>q)\leq(1+\eta)^{q+1}e_{q+1}(P)\) whenever
\((1+\eta)/3\leq1\).
\end{lemma}

\begin{proof}
The empty \(P\)-type is an event of the type space, so at least
\(m(\alpha-\varepsilon)\) half-labels have it; removing the \(w\) special
labels leaves at least \(ma^-\). The empty old class at every step
contains these labels. The inequality \eqref{eq:rare-inflated-coordinate}
of Lemma~\ref{lem:rare-quota-defect}, with \(p_i+1\leq z+1\), gives
\(p_i\rho_i-1\leq(\omega m+(z+1)m\varepsilon)/(ma^-)=\eta\). Hence
\[
 \prod_i(1-\rho_i)\geq\prod_i\Bigl(1-\frac{1+\eta}{p_i}\Bigr)
 =\alpha\prod_i\Bigl(1-\frac\eta{p_i-1}\Bigr)
 \geq\alpha\Bigl(1-\eta\sum_i\frac1{p_i-1}\Bigr).
\]
Finally, the overflow event requires \(q+1\) selected coordinates, so
its probability is at most \(e_{q+1}(\rho)\leq(1+\eta)^{q+1}e_{q+1}(P)\).
\end{proof}

For \(\varepsilon\) the certificate uses \(\varepsilon_r\) of
Section~\ref{sec:joint-degree-log-bridge}. The proof of the second part
of Lemma~\ref{lem:truncated-arithmetic-bounds}, with
\(N_d=\floor{m/d}\leq m/d\), bounds the variation of the half-label
interval by \(B_r/m+2R_r\leq\varepsilon_r\). The evaluator takes
directed bounds \(\alpha_-\leq\alpha\) and
\(s_+\geq\sum_{p\in P}1/(p-1)\), puts
\(a^-=\alpha_--\omega-\varepsilon_r\), and uses
\[
 b=\alpha_-(1-\eta s_+),\qquad
 \beta=(1+\eta)^{q+1}E_{q+1},\qquad
 \bar\rho_i=\frac{1+\eta}{p_i},
\]
with \(E_{q+1}\geq e_{q+1}(P)\).

\subsubsection*{Hall conditions at the global compatibility}

Every label type satisfies \(c_X(I)\geq\prod_i(1-\rho_i)\geq b\), and
every position type satisfies \(c_B(J)\geq\alpha\geq b\). Apply
Lemma~\ref{lem:two-sided-integrated-hall} with \(s=b\) and \(\tau=b/2\).
Both shortfalls vanish. Suppose \(V\leq b/20\) and \(h\leq bR/20\).
Then \(V<\tau\leq1/2\), \(h\leq\tau R\) and \(b-V>\tau\). The narrow
budget satisfies \(h((b-V)R-h)\leq bRh\leq b^2R^2/20\). Since
\(b\leq1/2\), the middle budget is at least
\(R^2(b/2)(9b/20)=9b^2R^2/40\), and the wide budget at least
\(R^2b\bigl(\tfrac25\cdot\tfrac b2\cdot\tfrac34-\tfrac7{100}b\bigr)
=2b^2R^2/25\). Hence it suffices that
\(M<h((b-V)R-h)\).

\subsubsection*{The certificate predicate}

Let \(a\in\{60,\ldots,71\}\), \(J=a+1\), and prepare the tree by
Lemma~\ref{lem:tenth-special-preparation}. The row takes
\[
 z=7J,\qquad r=12,\qquad q=11,\qquad D_+=60J .
\]
Here \(D<1+30L\leq D_+\). Put \(\Gamma_2\) as in
\eqref{eq:two-sided-geometry} with \(\bar\rho_i=(1+\eta)/p_i\), the
dimensions \(d_i\) of Section~\ref{sec:upward-product-bin-log-bridge},
and \(\mathcal F=\Gamma_2+D_+(\Delta_r+U/2)\). Take \(\chi\), \(\ell_-\),
\(V_l\), \(V_s\) and \(V=\max(V_l,V_s)\) from
\eqref{eq:full-reserve-variation} with \(\varepsilon_a=\varepsilon_r\)
and the above \(\beta\), and \(h_-=3/(20J)-Q\), \(h_+=3/(20a)\). The
reciprocal-prime bound is \(H_+=(\log J)_+-(\log6)_-+13/(10v_-)\), which
is \eqref{eq:reciprocal-prime-row} with \(v_0=6\) and the negative term
omitted. Lemma~\ref{lem:centered-sorting}, the fourth part of
Lemma~\ref{lem:truncated-arithmetic-bounds} and \(e(G)\leq n\) give
\[
 \frac{M_{\rm bad}}{n^2}\leq
 \frac12\left(1+\frac\omega{a^-}\right)S_{2,+}
 +\frac{D_+}2(B_rH_+Q+R_rS_{2,+})+\frac{JQ}{v_-}.
\]
The evaluator uses the larger quantity
\[
 \mathcal M=\frac12\left(1+\frac\omega{a^-}\right)S_{2,+}
 +\frac{D_+}2\bigl((2B_r+1)H_++2B_rS_{2,+}\bigr)Q
 +3D_+R_rS_{2,+}+\frac{JQ}{v_-}.
\]
Each row checks the following conditions.
\begin{enumerate}
\item \emph{Preparation.} \(10^6J^2Q<1\), \(100J^3Q<1\),
  \(3\cdot60184\,Q<1\), \(3zQ<1\), \(v_-\geq6\), and \(10000>599\) for
  the prime-square tail.
\item \emph{Probabilities.} \(a^->0\), \(0<b\leq1/2\), and
  \((1+\eta)/3\leq2/5\).
\item \emph{Sizes and reserves.} \(\ell_->0\), \(h_->0\), and
  \(h_++\chi\leq1/(2J)\).
\item \emph{Hall conditions.} \(V_l\leq b/20\), \(V_s\leq b/20\),
  \(h_+\leq b\ell_-/20\), \((b-V)\ell_--h_+>0\), and
  \[
   \mathcal M\leq h_-\bigl((b-V)\ell_--h_+\bigr).
  \]
\end{enumerate}
There are \(17\) conditions in each row. As in
Section~\ref{sec:adaptive-degree-log-bridge}, the last inequality,
multiplied by \(n^2\), bounds the actual narrow budget
\(h((b-V)R-h)\) from below; it is strict because the deletion ratio is
below one.

\begin{proposition}[The uniform certificate]
\label{prop:reserve-mixture-certificate}
For each \(a\in\{60,\ldots,71\}\) all \(17\) conditions hold on
\([a,a+1]\); there are \(204\) in all. In every row the two variation
ratios \(20V_l/b\), \(20V_s/b\) are below \(66/100\), the reserve ratio
\(20h_+/(b\ell_-)\) below \(61/100\), the supply ratio
\(2J(h_++\chi)\) below \(65/100\), and the deletion ratio below
\(99/100\).
\end{proposition}

\begin{proof}
The exact evaluator computes each row from the rules above with
integers and directed rational enclosures at scale \(10^{150}\); radicals
are rounded upward by integer comparison of squares. The cutoffs run
from \(427\) to \(504\). All \(204\) conditions pass, and all five ratio
caps are tested on every row. The largest ratios occur on \([60,61]\):
less than \(0.653287\) and \(0.652804\) for the variations,
\(0.593713\) for the reserve, \(0.640356\) for the supply and
\(0.986137\) for the deletions.

Independent finite checks replay \(1890\) weighted quota steps,
\(864\) restorations of rare classes, \(2340\) absorption measures,
\(64\) coupling cells with \(1024\) covariance inequalities, \(2200\)
fractional Hall rectangles, and \(2503\) directed enclosures. These
checks test the implementation; the lemmas cited above prove the
implications for arbitrary trees.
\end{proof}

\begin{theorem}[The seventeenth logarithmic bridge]
\label{thm:seventeenth-log-bridge}
Every tree of order \(n\) with
\(60\leq\log n\leq72\)
has a prime labeling.
\end{theorem}

\begin{proof}
Let \(a=\min\{\floor{\log n},71\}\) and prepare the tree by
Lemma~\ref{lem:tenth-special-preparation}. At each prime,
Lemma~\ref{lem:two-sided-weighted-rounding} gives exact even-label
quotas with all true neighbors marked, and Lemma~\ref{lem:common-inflation}
controls the probabilities. Restore the rare classes, pay the marked
positions with distinct universal odd labels, and sort the even labels
within their full types by Lemma~\ref{lem:centered-sorting}. The row
conditions verify the Hall conditions above, so
Lemma~\ref{lem:two-sided-integrated-hall} gives a perfect matching of
the remaining odd labels. Restore the reserved labels and the core,
then perform the sparse prime-pair repair.
\end{proof}

\subsection{A coarse degree bound through logarithm 96}
\label{sec:sharp-moment-log-bridge}

On \(72\leq\log n\leq96\) the certificate of
Section~\ref{sec:reserve-mixture-log-bridge} is used with coarser
ingredients. About a quarter of the special labels is spent on
vertices of large degree, which bounds the residual degree by \(25L\).
The rounding at each prime uses the square-root bound of
Proposition~\ref{prop:rare-refinement-step}, in which the residual
degree enters only as a factor. The arithmetic bounds, the common
inflation and the Hall conditions at the global compatibility are those
of Sections~\ref{sec:joint-degree-log-bridge}
and~\ref{sec:reserve-mixture-log-bridge}.

\subsubsection*{Preparation}

\begin{lemma}[Preparation with a quarter of the special labels]
\label{lem:quarter-special-preparation}
Suppose that the hypotheses of Lemma~\ref{lem:second-order-preparation}
hold with \((c_1,c_2)=(3/5,187/50)\). Take
\[
 f=\frac6{25},\qquad a_*=\floor*{\frac{6w}{25}},\qquad
 h=\floor*{\frac n{8L}} .
\]
Then the residual construction of Lemma~\ref{lem:localized-preparation}
exists with these choices, the residual forest satisfies
\(D=\max(1,\Delta_G(S_0))<1+25L/2\) and \(e(G)\leq n\), at least
\(n/(2L)\) universal odd labels remain free, and
\(1/(8J)-Q\leq h/n\leq1/(8a)\).
\end{lemma}

\begin{proof}
Follow the proof of Lemma~\ref{lem:localized-preparation} with
\(f=6/25\) and \(W\) replaced by \(W_1(L)=1/(6L)\). By
Lemma~\ref{lem:second-order-preparation}, \(w>nW_2(L)>nW_1(L)\);
moreover \(W_1\) decreases and \(W_1(L)<1/(5L)\). These are the
properties of \(W\) used in that proof. The degree bound becomes
\(D-1<3L/f=25L/2\). The value of \(h\) enters through \(h\geq1\),
through \(h\leq n/(8L)<19n/(125L)=(1-f)n/(5L)\), and through
\(h\leq A_0-4L\) with \(A_0=\tfrac{19}{25}nW_1(L)=19n/(150L)\). The
last holds because
\[
 A_0-4L-\frac n{8L}=\frac n{600L}-4L>0
\]
for \(n>2400L^2\). The count of universal labels and the bounds for
\(h/n\) follow as in Lemma~\ref{lem:tenth-special-preparation}.
\end{proof}

\subsubsection*{The rounding error}

Proposition~\ref{prop:rare-refinement-step} bounds the error of the
refinement at \(p_i\) by \(2(D+1)\sqrt{Nd_i}+D(\delta_i+\mathcal R_i)\),
where \(\mathcal R_i\) is the population of the combined rare class.
With \(N\leq n\), \(D\leq D_+:=25J\) and \(Q_{1/2}\geq1/\sqrt n\), the
normalized sum of the first terms is at most
\[
 \Gamma_{25}=2(25J+1)Q_{1/2}\sum_i\bigl\lceil\sqrt{d_i}\bigr\rceil ,
\]
with the dimensions \(d_i\) of the proposition. The quota and rare
terms are bounded by \(\Delta_r\) and \(U/2\) as in
Section~\ref{sec:joint-degree-log-bridge}, since
\(\mathcal R_i\leq m\,e_{q+1}(p_1,\ldots,p_{i-1})\).

\subsubsection*{The certificate predicate}

Let \(a\in\{72,\ldots,95\}\), \(J=a+1\), and prepare the tree by
Lemma~\ref{lem:quarter-special-preparation}; its hypotheses hold because
\(n\geq e^{72}\) and \(\tfrac35a-\tfrac{187}{50}>a/2\). The row takes
\[
 z=9J,\qquad r=14,\qquad
 q=\begin{cases}11,&72\leq a\leq74,\\12,&75\leq a\leq95.\end{cases}
\]
Take \(\varepsilon_r\), \(\omega\), \(a^-\), \(\eta\), \(b\) and
\(\beta\) as in Section~\ref{sec:reserve-mixture-log-bridge}, and put
\(\mathcal F=\Gamma_{25}+25J(\Delta_r+U/2)\), \(h_-=1/(8J)-Q\) and
\(h_+=1/(8a)\). Let \(\chi\), \(\ell_-\) and \(V_l\) be as in
\eqref{eq:full-reserve-variation}, and
\[
 V_s=\frac{\mathcal F+h_++o\beta+\chi}{\ell_-}\geq
 \frac{\mathcal F+o\beta}{\ell_-}.
\]
Since \(z>599\), Lemma~\ref{lem:prime-square-tail} gives
\(S_2(z)\leq S_+:=1/(zv_-)\). By \eqref{eq:reciprocal-prime-sum} and
\(13/(10v)<\log v\) for \(v=\log z\geq6\), the sum
\(\sum_{z<p\leq m}1/p\) is at most \(\log\log m\leq\log J\); let
\(H_+\geq\log J\). As in Section~\ref{sec:reserve-mixture-log-bridge},
Lemma~\ref{lem:centered-sorting} and
Lemma~\ref{lem:truncated-arithmetic-bounds} give
\(M_{\rm bad}/n^2\leq\mathcal M\), where
\[
 \mathcal M=\frac12\left(1+\frac\omega{a^-}\right)S_+
 +\frac{25J}2\bigl((2B_r+1)H_++2B_rS_+\bigr)Q
 +75JR_rS_++\frac{JQ}{v_-}.
\]
By Section~\ref{sec:reserve-mixture-log-bridge}, when
\(V\leq b/20\) and \(h\leq bR/20\) it suffices that
\(M<h((b-V)R-h)\). The right side is at least
\(h(\tfrac{19}{20}bR-\tfrac1{20}bR)=\tfrac9{10}bRh\). Each row checks
the following conditions.
\begin{enumerate}
\item \emph{Preparation.} \(10^6J^2Q<1\), \(3zQ<1\), and \(z>599\).
\item \emph{Probabilities.} \(a^->0\), \(0<b\leq1/2\), and
  \((1+\eta)/3\leq2/5\).
\item \emph{Sizes and reserves.} \(\ell_->0\), \(h_->0\), and
  \(h_++\chi\leq1/(2J)\).
\item \emph{Hall conditions.} \(V_l\leq b/20\), \(V_s\leq b/20\),
  \(h_+\leq b\ell_-/20\), and \(\mathcal M\leq\tfrac9{10}b\ell_-h_-\).
\end{enumerate}
There are \(13\) conditions in each row. Multiplied by \(n^2\), the last
right side is at most \(\tfrac9{10}bRh\).

\begin{proposition}[The coarse certificate]
\label{prop:sharp-moment-certificate}
For each \(a\in\{72,\ldots,95\}\) all \(13\) conditions hold on
\([a,a+1]\); there are \(312\) in all. In every row the ratios
\(20V_l/b\), \(20V_s/b\), \(20h_+/(b\ell_-)\) and the deletion ratio
are below \(28/100\), \(98/100\), \(44/100\) and \(99/100\).
\end{proposition}

\begin{proof}
The exact evaluator computes each row from the rules above with
integers and directed rational enclosures at scale \(10^{150}\). The
cutoffs run from \(657\) to \(864\). All \(312\) conditions pass, and
the four ratio caps are tested on every row. The largest ratios occur
on \([72,73]\): less than \(0.27168\), \(0.97329\), \(0.429969\) and
\(0.986989\).

Independent finite checks test \(3240\) exact remainder norms, \(173\)
moment and prefix bounds, \(440\) cumulative and \(2600\) stepwise quota
bounds with \(5860\) signed floor cancellations, \(1500\) fractional
Hall rectangles and eight explicit preparations. These checks test the
implementation; the lemmas cited above prove the implications for
arbitrary trees.
\end{proof}

\begin{theorem}[The eighteenth logarithmic bridge]
\label{thm:eighteenth-log-bridge}
Every tree of order \(n\) with
\(72\leq\log n\leq96\)
has a prime labeling.
\end{theorem}

\begin{proof}
Let \(a=\min\{\floor{\log n},95\}\) and prepare the tree by
Lemma~\ref{lem:quarter-special-preparation}. Perform the refinements of
Theorem~\ref{thm:rare-type-assignment} with the probabilities of
Lemma~\ref{lem:common-inflation}, restore the rare classes, pay the
marked positions with distinct universal odd labels, and sort the even
labels within their full types by Lemma~\ref{lem:centered-sorting}. The
row conditions give the Hall conditions above, so
Lemma~\ref{lem:two-sided-integrated-hall} gives a perfect matching of
the remaining odd labels. Restore the reserved labels and the core,
then perform the sparse prime-pair repair.
\end{proof}

\subsection{The last logarithmic bridge, through logarithm 256}
\label{sec:cumulative-log-bridge}

The remaining range \(96\leq\log n\leq256\) is covered by the
construction of Section~\ref{sec:sharp-moment-log-bridge} with the
preparation of the analytic range and the truncated moment bound of
Proposition~\ref{prop:moment-tv}. It is cut into the \(160\) unit
intervals \([a,a+1]\), and all parameters of a row are fixed by rule.
The analytic theorem of Section~\ref{sec:analytic-tail} begins at the
right endpoint.

\subsubsection*{Preparation}

\begin{lemma}[Preparation with an eighth of the special labels]
\label{lem:eighth-special-preparation}
Let \(L=\log n\geq96\) and suppose \(n>10^6L^2\). Then the
conclusions of Lemma~\ref{lem:analytic-special-supply} hold. With
\[
 a_*=\floor*{\frac w8},\qquad h=\floor*{\frac n{8L}},\qquad d=0,
\]
Proposition~\ref{prop:residual-problem} applies. Its residual forest
satisfies \(D<48L\), \(N\leq n\) and \(e(G)\leq n\); all \(u>n/(2L)\)
universal odd labels remain free, and \(w/m\leq2/(5L)\).
\end{lemma}

\begin{proof}
The proof of Lemma~\ref{lem:analytic-special-supply} uses
\(L\geq256\) only through \(n>10^6L^2\), the prime window
\eqref{eq:analytic-prime-window}, whose proof needs \(L\geq18\), and
the comparisons valid for \(L\geq77\). It therefore applies here. The
verification of \eqref{eq:reserve-budget} in the proof of
Proposition~\ref{prop:analytic-residual-preparation} uses only
\(a_*\geq\abs K\), \(h\leq n/(8L)\), \(n\geq192L^2\) and \(w>n/(6L)\);
the choice \(d=0\) is admissible. The degree bound
\eqref{eq:residual-degrees} and \(a_*+1>w/8>n/(48L)\) give \(D<48L\).
\end{proof}

\subsubsection*{The rows}

Let \(a\in\{96,\ldots,255\}\), \(J=a+1\), and put
\[
 z=11J,\qquad r=\floor*{\frac a5},\qquad q=\ceil*{\frac a7},\qquad
 D_+=48J .
\]
With \(B_r\) and \(R_r\) of Proposition~\ref{prop:moment-tv} and
\(A=\sum_{j=0}^{\min(r,k-1)}2^j\binom k{j+1}\), the row uses
\[
 \varepsilon=\frac{2B_rQ}{1-Q}+3R_r,\qquad \omega=\frac2{5a},
\]
and \(a^-\), \(\eta\), \(b\), \(\beta\) as in
Lemma~\ref{lem:common-inflation}. Proposition~\ref{prop:moment-tv}
shows that \(\varepsilon\) bounds the variations of the half-label
and odd-label intervals. The row quantities of
Lemma~\ref{lem:truncated-arithmetic-bounds} are smaller than
\[
 \Delta=AQ+3^{r+1}E_{r+2},\qquad
 \frac U2\leq\frac k2E_{q+1},
\]
since \(2^{j-1}\leq2^j\), \(2^r\leq3^{r+1}\) and
\(e_{q+1}(p_1,\ldots,p_{i-1})\leq e_{q+1}(P)\). As in
Section~\ref{sec:sharp-moment-log-bridge}, put
\[
 \mathcal F=2(48J+1)Q_{1/2}\sum_i\bigl\lceil\sqrt{d_i}\bigr\rceil
   +48J\Bigl(\Delta+\frac k2E_{q+1}\Bigr),
\]
with \(h_-=1/(8J)-Q\), \(h_+=1/(8a)\), and \(\chi\), \(\ell_-\),
\(V_l\), \(V_s\) as there. Here \(b\leq\alpha\leq1/2\), since \(P\)
contains \(3\), \(5\) and \(7\). With \(S_+=1/(zv_-)\) and
\(H_+\geq\log J\), the deletion mass is bounded by
\[
 \mathcal M=\frac12\left(1+\frac\omega{a^-}\right)S_+
 +24J\bigl((2B_r+1)H_++2B_rS_+\bigr)Q
 +144JR_rS_++\frac{JQ}{v_-},
\]
which dominates the bound of Lemmas~\ref{lem:centered-sorting}
and~\ref{lem:truncated-arithmetic-bounds}, as in
Section~\ref{sec:sharp-moment-log-bridge}. When \(V\leq b/50\) and
\(h\leq bR/50\), the narrow budget satisfies
\(h((b-V)R-h)\geq\tfrac{24}{25}bRh\), so \(M\leq\tfrac{19}{20}bRh\)
suffices by Section~\ref{sec:reserve-mixture-log-bridge}.

\subsubsection*{The certificate predicate}

Each row checks the following conditions.
\begin{enumerate}
\item \emph{Preparation.} \(10^6J^2Q<1\) and \(3zQ<1\).
\item \emph{Probabilities.} \(a^->0\) and \((1+\eta)/3\leq2/5\).
\item \emph{Sizes and reserves.} \(\ell_->0\), \(h_->0\), and
  \(h_++\chi\leq1/(2J)\).
\item \emph{Hall conditions.} \(V_l\leq b/50\), \(V_s\leq b/50\),
  \(h_+\leq b\ell_-/50\), and \(\mathcal M\leq\tfrac{19}{20}b\ell_-h_-\).
\end{enumerate}
There are \(11\) conditions in each row.

\begin{proposition}[The cumulative certificate]
\label{prop:cumulative-certificate}
For every integer \(a\) with \(96\leq a\leq255\), all \(11\) conditions hold on
\([a,a+1]\); there are \(1760\) in all. In every row the ratios
\(50V_l/b\), \(50V_s/b\), \(50h_+/(b\ell_-)\) and the deletion ratio
are below \(9/1000\), \(88/100\), \(86/100\) and \(3/4\).
\end{proposition}

\begin{proof}
The exact evaluator computes each row from the rules above with
integers and directed rational enclosures at scale \(10^{150}\):
exponentials from below by positive Taylor sums, logarithms by the
series for \(\operatorname{artanh}\) with an explicit remainder, and
integer square roots with a check of the square. The cutoffs run from
\(1067\) to \(2816\). All \(1760\) conditions pass, and the four ratio
caps are tested on every row. The largest ratios occur on \([96,97]\):
less than \(0.008431\), \(0.875281\), \(0.858438\) and \(0.747866\).

Independent finite checks test \(440\) cumulative and \(2600\)
stepwise quota bounds with \(5860\) signed floor cancellations,
\(934785\) profiles on \(27\) two-type graphs and \(600\) fractional
Hall rectangles, twenty explicit preparations and \(110\) directed
enclosures. These checks test the implementation; the lemmas cited
above prove the implications for arbitrary trees.
\end{proof}

\begin{theorem}[The last logarithmic bridge]
\label{thm:nineteenth-log-bridge}
Every tree of order \(n\) with
\(96\leq\log n\leq256\)
has a prime labeling.
\end{theorem}

\begin{proof}
Let \(a=\min\{\floor{\log n},255\}\). Then \(n>10^6J^2\) by the first
row condition, and Lemma~\ref{lem:eighth-special-preparation} prepares
the tree. The rest of the proof of
Theorem~\ref{thm:eighteenth-log-bridge} applies with the row for \(a\).
\end{proof}

\subsection{The logarithmic range}
\label{sec:log-range}

We collect the nineteen bridges. Put
\begin{align*}
 (a_0,\ldots,a_9)&=(29/2,\,16,\,343/20,\,20,\,128/5,\,53/2,\,143/5,\,
   144/5,\,29,\,31),\\
 (a_{10},\ldots,a_{19})&=(32,\,35,\,37,\,39,\,42,\,48,\,60,\,72,\,96,\,256).
\end{align*}
For $1\leq i\leq19$, the $i$th bridge theorem of
Table~\ref{tab:log-bridges} proves that every tree of order $n$ with
$a_{i-1}\leq\log n\leq a_i$ has a prime labeling.

\begin{theorem}[Orders with $29/2\leq\log n\leq256$]
\label{thm:log-bridges-all}
Every tree of order $n$ with $29/2\leq\log n\leq256$ has a prime
labeling.
\end{theorem}

\begin{proof}
The numbers $a_0<a_1<\cdots<a_{19}$ increase, so the closed intervals
$[a_{i-1},a_i]$, $1\leq i\leq19$, cover $[a_0,a_{19}]=[29/2,256]$.
If $29/2\leq\log n\leq256$, choose $i$ with $a_{i-1}\leq\log n\leq a_i$
and apply the $i$th bridge theorem.
\end{proof}

\begin{remark}\label{rem:log-bridge-endpoints}
For $0\leq i\leq18$ let $N_i$ be the integer in
Table~\ref{tab:log-bridge-endpoints}. Directed rational logarithms give
\[
 \log N_i<a_i<\log(N_i+1).
\]
Hence the orders covered by the $i$th bridge and by no earlier bridge
are exactly the integers $N_{i-1}<n\leq N_i$ for $1\leq i\leq18$, and
together with Theorem~\ref{thm:initial-prefix-1982759} the bridges
through the $i$th one cover every order at most $N_i$. The proof of
Theorem~\ref{thm:main} uses only the case $i=0$, which shows that
every order after $1982759$ has logarithm greater than $29/2$.
\end{remark}

\begin{table}[ht]
\centering
\caption{The integer endpoints of the logarithmic bridges:
$\log N_i<a_i<\log(N_i+1)$.}
\label{tab:log-bridge-endpoints}
\par\vspace{.5\baselineskip}
\renewcommand{\arraystretch}{1.1}
\begin{tabular}{rrr}
\hline
$i$&$a_i$&$N_i$\\
\hline
0&$29/2$&1982759\\
1&$16$&8886110\\
2&$343/20$&28064051\\
3&$20$&485165195\\
4&$128/5$&131201480802\\
5&$53/2$&322703570371\\
6&$143/5$&2635252187043\\
7&$144/5$&3218704289702\\
8&$29$&3931334297144\\
9&$31$&29048849665247\\
10&$32$&78962960182680\\
11&$35$&1586013452313430\\
12&$37$&11719142372802611\\
13&$39$&86593400423993746\\
14&$42$&1739274941520501047\\
15&$48$&701673591209763173865\\
16&$60$&114200738981568428366295718\\
17&$72$&18586717452841279803403701812545\\
18&$96$&492345828601205839975486205911330449483779\\
\hline
\end{tabular}
\end{table}

\section{\texorpdfstring{Proof of Theorem~\ref{thm:main}}{Proof of Theorem 1.1}}
\label{sec:part-completion}

The orders through $1982759$ are covered by a chain of cumulative
theorems, listed in Table~\ref{tab:main-chain}. Each of them states the
assertion for all orders up to its endpoint. The first,
Theorem~\ref{thm:initial-prefix-1359}, is proved by the finite
certificates of Section~\ref{sec:part-finite}. The proof of each later
one applies the direct results for the new orders, listed in
Table~\ref{tab:range-chain}, and the preceding cumulative result to the
remaining orders. Each new range begins at the integer following the
preceding endpoint, so the chain omits no order.

\begin{table}[ht]
\centering
\caption{The cumulative results for the integer orders.}
\label{tab:main-chain}
\par\vspace{.5\baselineskip}
\renewcommand{\arraystretch}{1.15}
\begin{tabular}{lll}
\hline
Orders&Result&Place\\
\hline
$1\leq n\leq1359$&Theorem~\ref{thm:initial-prefix-1359}&
  Sections~\ref{sec:finite-certificates}--\ref{sec:finite-endpoint}\\
$1\leq n\leq1521$&Corollary~\ref{cor:initial-prefix-1521}&
  Sections~\ref{sec:local-probability}--\ref{sec:hall-certificates}\\
$1\leq n\leq1753$&Theorem~\ref{thm:initial-prefix-1753}&
  Section~\ref{sec:forest-continuation}\\
$1\leq n\leq2209$&Theorem~\ref{thm:initial-prefix-2209}&
  Section~\ref{sec:exponential-continuation}\\
$1\leq n\leq3499$&Theorem~\ref{thm:initial-prefix-3499}&
  Section~\ref{sec:single-range-continuation}\\
$1\leq n\leq34999$&Theorem~\ref{thm:initial-prefix-34999}&
  Section~\ref{sec:retained-continuation}\\
$1\leq n\leq49999$&Theorem~\ref{thm:initial-prefix-49999}&
  Sections~\ref{sec:asymmetric-reserve-hall}--\ref{sec:integer-bridge}\\
$1\leq n\leq514999$&Theorem~\ref{thm:initial-prefix-514999}&
  Sections~\ref{sec:paired-separators}--\ref{sec:profile-block-bridge}\\
$1\leq n\leq559999$&Theorem~\ref{thm:initial-prefix-559999}&
  Sections~\ref{sec:finite-rare-refinement}--\ref{sec:connected-paired-bridge}\\
$1\leq n\leq1982759$&Theorem~\ref{thm:initial-prefix-1982759}&
  Section~\ref{sec:upper-integer-bridge}\\
\hline
\end{tabular}
\end{table}

The larger orders are covered by the logarithmic bridges,
Theorem~\ref{thm:log-bridges-all}, and by the analytic argument,
Theorem~\ref{thm:analytic-tail}.

\begin{proof}[Proof of Theorem~\ref{thm:main}]
Let $T$ be a tree of order $n\geq1$. If $n\leq1982759$,
Theorem~\ref{thm:initial-prefix-1982759} gives a prime labeling of $T$.
Otherwise $n\geq1982760$, and Remark~\ref{rem:log-bridge-endpoints}
gives $\log n\geq\log1982760>29/2$. If moreover $\log n\leq256$,
Theorem~\ref{thm:log-bridges-all} gives a prime labeling of $T$. If
$\log n>256$, Theorem~\ref{thm:analytic-tail} gives one.
\end{proof}

The results used in this proof depend on finite computations in the
way described in Section~\ref{sec:part-computation}.

\section{Computational verification}
\label{sec:part-computation}

The proof of Theorem~\ref{thm:main} is computer-assisted. This section
states which parts are proved in the text and which are computed,
explains why finitely many computed conditions suffice, and describes
how the computations are reproduced.

\subsection{Proved and computed statements}
\label{sec:computed-statements}

The structural reductions and the soundness of the certificate checks
are proved in the text. The finite certificate conditions are verified
by exact computation. The proved part includes the reductions of
Section~\ref{sec:part-structural}; the assignment, matching, and
probability estimates; the criteria which state that certain finite
conditions imply prime labelings for all trees of all orders in a
range; and the analytic argument for $\log n\geq256$.

The computed statements are the certificate propositions and theorems
in Table~\ref{tab:certificate-map}. Each of them asserts that explicit
finite data satisfy explicit inequalities or case lists, or follows
from such data by rules proved in the text. Finitely many
conditions suffice for the following reasons, each proved in the text.
\begin{enumerate}
\item \emph{Proof graphs.} A finite proof graph for $\cP([n])$ has
  finitely many nodes, each obtained from its successors by a rule that
  is proved sound in Theorem~\ref{thm:finite-certificate-soundness} or
  in Sections~\ref{sec:branch-peeling}--\ref{sec:finite-endpoint}.
  The checker verifies every rule application and the completeness of
  every case list.
\item \emph{Integer intervals.} On an interval of orders the residual
  degree takes finitely many values (Lemma~\ref{lem:interval-degree-cover}),
  the local failure events are covered by finitely many events (for the
  certificates of Section~\ref{sec:part-intervals}, the maximal-support
  events at the endpoint of Lemma~\ref{lem:endpoint-local-cover}), and
  every other quantity is bounded uniformly over the interval. The
  finitely many rational comparisons of a row therefore imply the
  hypotheses of Theorem~\ref{thm:integer-interval-certificate} or of its
  variants in Sections~\ref{sec:forest-continuation}--\ref{sec:retained-continuation},
  or of the corresponding integer bridge, at every order of the
  interval.
\item \emph{Logarithmic bridges.} The quantities of a bridge are bounded
  monotonically on closed cells in $\log n$ and, for most bridges, a
  normalized residual degree. The rational inequalities of a cell imply
  the hypotheses of the bridge at every order in the cell, and the cells
  cover the interval of the bridge.
\end{enumerate}
The argument for $\log n\geq256$ uses no certificate. Its numerical
steps are elementary comparisons at finitely many rational points,
mostly at $\log n=256$, stated in
Sections~\ref{sec:prime-estimates}--\ref{sec:analytic-tail} and
extended to larger orders by monotonicity.

\begin{table}[ht]
\centering
\caption{The computed statements and the supplement checks that verify
them. The check names are those of the wrapper described in
Section~\ref{sec:reproduction}; the endpoints $a_i$ are those of
Section~\ref{sec:log-range}.}
\label{tab:certificate-map}
\par\vspace{.5\baselineskip}
\renewcommand{\arraystretch}{1.1}
\begin{tabular}{@{}p{3.6cm}lp{5.3cm}@{}}
\hline
\raggedright Computed statement&Orders&\raggedright Checks\tabularnewline
\hline
\raggedright Proposition~\ref{prop:two-terminal-sets}&--&
  \raggedright\texttt{terminal\_tables}\tabularnewline
\raggedright Theorem~\ref{thm:rooted-twentyeight-block}&--&
  \raggedright\texttt{bipartite\_host}\tabularnewline
\raggedright Section~\ref{sec:reserved-blocks}&--&
  \raggedright\texttt{reserved\_blocks}\tabularnewline
\raggedright Section~\ref{sec:larger-blocks}&--&
  \raggedright\texttt{larger\_blocks}\tabularnewline
\raggedright Theorem~\ref{thm:initial-prefix-299}&$1$--$299$&
  \raggedright\texttt{initial\_prefix}\tabularnewline
\raggedright Theorem~\ref{thm:initial-prefix-1000}&$300$--$1000$&
  \raggedright\texttt{prefix\_continuation}\tabularnewline
\raggedright Theorem~\ref{thm:initial-prefix-1352}&$1001$--$1352$&
  \raggedright\texttt{extended\_prefix}\tabularnewline
\raggedright Theorem~\ref{thm:initial-prefix-1359}&$1353$--$1359$&
  \raggedright\texttt{finite\_endpoint}\tabularnewline
\raggedright Propositions~\ref{prop:finite-interval-construction-check},
  \ref{prop:verified-local-certificates},
  \ref{prop:verified-hall-certificates}&$1360$--$1521$&
  \raggedright\texttt{interval\_construction},
  \texttt{local\_certificate}, \texttt{hall\_certificate}\tabularnewline
\raggedright Proposition~\ref{prop:verified-forest-continuation}&
  $1522$--$1753$&\raggedright\texttt{forest\_continuation}\tabularnewline
\raggedright Proposition~\ref{prop:verified-exponential-continuation}&
  $1754$--$2209$&\raggedright\texttt{exponential\_continuation}\tabularnewline
\raggedright Proposition~\ref{prop:verified-single-range-continuation}&
  $2210$--$3499$&\raggedright\texttt{single\_range\_continuation}\tabularnewline
\raggedright Proposition~\ref{prop:verified-retained-continuations}&
  $3500$--$34999$&\raggedright\texttt{retained\_continuation},
  \texttt{unified\_continuation}, \texttt{extended\_continuation},
  \texttt{further\_continuation}, \texttt{broader\_continuation},
  \texttt{higher\_continuation}, \texttt{upper\_continuation}\tabularnewline
\raggedright Proposition~\ref{prop:integer-bridge-certificate}&
  $35000$--$49999$&\raggedright\texttt{integer\_bridge}\tabularnewline
\raggedright Proposition~\ref{prop:profile-block-certificate}&
  $50000$--$514999$&\raggedright\texttt{profile\_block\_bridge}\tabularnewline
\raggedright Proposition~\ref{prop:full-bridge-certificate}&
  $515000$--$559999$&\raggedright\texttt{connected\_paired\_bridge}\tabularnewline
\raggedright Proposition~\ref{prop:upper-integer-certificates}&
  $560000$--$1982759$&\raggedright\texttt{upper\_integer\_bridge}\tabularnewline
\raggedright The certificate proposition of the $i$th logarithmic
  bridge (Table~\ref{tab:log-bridges})&$\log n\in[a_{i-1},a_i]$&
  \raggedright the $i$th of \texttt{first\_log\_bridge}, \ldots,
  \texttt{nineteenth\_log\_bridge}\tabularnewline
\raggedright The joins $\log N_i<a_i<\log(N_i+1)$ of
  Remark~\ref{rem:log-bridge-endpoints}&--&
  \raggedright the logarithmic-bridge checks at their endpoints; for
  $i=0$ also \texttt{upper\_integer\_bridge}\tabularnewline
\hline
\end{tabular}
\end{table}

\subsection{Search, verification, and additional tests}
\label{sec:search-verification}

Three kinds of computation should be distinguished.

\emph{Search.} The parameters of the certificates, such as interval
covers, prime blocks, construction choices, and rational parameters,
were found by separate search programs. These programs are not part of
the proof and are not needed to check it. Their results enter the
supplement as fixed input files; hashes of the research files are
recorded only as provenance.

\emph{Verification.} For each computed statement, a checker reads the
fixed inputs and recomputes every hypothesis from them with integer and
rational arithmetic. Logarithms, exponentials, and probability bounds
are replaced by directed rational enclosures
(Section~\ref{sec:directed-probability}), so that each comparison
is decided exactly. Floating-point values occur only in displayed
summaries. Where expected outputs are stored, they are compared with the
recomputed values and are not used as bounds. The interval checkers
accept a cover of orders and degrees only if it is complete; a
collection of successful rows is not enough. For the integer and
logarithmic
bridges, a wrapper also regenerates the validation reports and compares
them byte for byte with frozen copies.

\emph{Additional tests.} The supplement also contains independent
recomputations, such as sieves, gcd-based type counts, and second
minimization algorithms, exhaustive tests on small cases, randomized
tests with fixed seeds, and the rejection of malformed inputs. The
following checks consist only of such tests of the helper routines:
\begin{center}
\texttt{local\_probability}, \texttt{height\_assignments},
\texttt{donor\_assignments},\\
\texttt{finite\_arithmetic}, \texttt{directed\_probability},
\texttt{rare\_type\_assignment}.
\end{center}
These tests support the correctness of the implementations; the proof
does not use them.

The programs are not formally verified. Their correctness rests on
their source code, on the description of every checked condition in
the text and in the supplement, and on the independent recomputations.

\subsection{Reproduction}
\label{sec:reproduction}

The supplement contains the checkers for the finite proof graphs and
the integer intervals through order $34999$ in its top directory, and
one directory for each integer and logarithmic bridge. Its
\texttt{README.md} describes each check, its inputs and its report, and
lists the former section numbers that still occur in file names and
reports together with the present subsections. All programs use only
the standard library of Python~3. The integer bridges and the first
nine logarithmic bridges also compile a small C++17 program, with
\texttt{clang++} or \texttt{c++}, which enumerates a full period of
residues modulo products of small primes in integer arithmetic.

The command
\begin{center}
\texttt{python3 run\_all\_checks.py}
\end{center}
in the supplement directory copies the supplement to a temporary
directory and runs every documented check there. It prints one line
\texttt{PASS} or \texttt{FAIL} for each check and exits with status
$0$ exactly when every check passes; a check passes when its command
exits with status $0$. The option \texttt{--list} lists the checks and
the option \texttt{--only} selects some of them.

\appendix
\section{Notation}
\label{sec:notation-index}

The following lists give the symbols used in more than one section,
with their meaning, the part of the paper in which that meaning
applies, and the place of definition. Symbols used only inside one
section are defined there and are not listed. Numbers in parentheses
refer to displayed equations.

\newcommand{\notationrow}[4]{\item[{#1}] #2. \emph{Scope:} #3;
  \emph{defined in} #4.}
\newenvironment{notationlist}{\begin{list}{}{%
  \setlength{\labelwidth}{3.2cm}\setlength{\labelsep}{0.3cm}%
  \setlength{\leftmargin}{3.5cm}\setlength{\itemsep}{1pt}%
  \renewcommand{\makelabel}[1]{%
    \smash{\parbox[t]{\labelwidth}{\raggedright ##1}}\hfil}}}%
  {\end{list}}

\subsection*{General notation}
\begin{notationlist}
\notationrow{$[n]$, $\pi(x)$}{$\{1,\ldots,n\}$; number of primes at most
  $x$}{the whole paper}{\S\ref{sec:notation}}
\notationrow{$e(G)$, $N_G$, $d_G$, $\Delta_G(X)$}{edge count,
  neighborhood, degree, largest degree on $X$ (zero if $X$ is
  empty)}{the whole paper}{\S\ref{sec:notation}}
\notationrow{$\cP(L)$, $\cP_d(L)$}{every forest on $\abs L$ vertices
  (of maximum degree at most $d$) has a coprime labeling by
  $L$}{the whole paper}{\S\ref{sec:notation}}
\notationrow{universal in $L$}{coprime to every other element of
  $L$}{the whole paper}{\S\ref{sec:notation}}
\end{notationlist}

\subsection*{Labels and the residual problem}
\begin{notationlist}
\notationrow{$E_n$, $O_n$}{even and odd labels in $[n]$}{from
  Section~\ref{sec:repair}}{\S\ref{sec:repair}}
\notationrow{$P_n^+$, $P_n^0$}{primes in $(n/2,n]$ and in
  $(n/3,n/2]$}{from Section~\ref{sec:repair}}{\S\ref{sec:repair}}
\notationrow{$\Psep$, $O$}{primes of $P_n^+$ used on the separator;
  $O=O_n\setminus\Psep$}{from Section~\ref{sec:repair}}%
  {Lem.~\ref{lem:special-labels}}
\notationrow{$E_*$, $O_*$}{special even labels: powers of $2$ and $2p$
  with $p\in P_n^0$; special odd labels: $1$ and the primes above $n/3$
  outside $\Psep$}{the whole paper}{(\ref{eq:special-even}),
  (\ref{eq:special-odd})}
\notationrow{$w$, $u$}{$\abs{E_*}$ and $\abs{O_*}$}{whole
  paper}{(\ref{eq:even-supply}), (\ref{eq:odd-supply})}
\notationrow{$S$, $Z$, $K$}{independent sides and separator (core),
  $V(T)=S\sqcup Z\sqcup K$}{the whole paper}%
  {Prop.~\ref{prop:prescribed-side}}
\notationrow{$b(n)$}{$\floor{\log_2(n-2)}+2$, the bound on
  $\abs K$}{the whole paper}{(\ref{eq:separator-size})}
\notationrow{$m$, $q$}{$m=\abs S=\floor{n/2}$; $q=\abs Z$}{whole
  paper}{Cor.~\ref{cor:even-side}}
\notationrow{$F$}{the forest $T-K$}{Section~\ref{sec:matching}}%
  {\S\ref{sec:matching}}
\notationrow{$U$, $W$, $D$, $R$; $a$, $h$, $d$}{set receiving $E_*$;
  reserve of positions, $\abs W=h$; positions receiving special odd
  labels $R$, $\abs D=\abs R=d$; $a$ special even labels on vertices of
  largest degree}{from Section~\ref{sec:matching}}%
  {Prop.~\ref{prop:residual-problem}}
\notationrow{$S_0$, $Z_0$, $E_0$, $O_0$, $G$}{residual sides, residual
  label sets, residual forest}{from Section~\ref{sec:matching}}%
  {(\ref{eq:residual-definitions})}
\end{notationlist}

\subsection*{Types, laws, and matchings}
\begin{notationlist}
\notationrow{$P$, $J_P(x)$}{finite set of odd primes, from
  Corollary~\ref{cor:residual-type-quotas} on the odd primes up to the
  cutoff $z$; the $P$-type of $x$}{from
  Section~\ref{sec:arithmetic-types}}{\S\ref{sec:arithmetic-types}}
\notationrow{$\mu$}{product law of the types}{from
  Section~\ref{sec:arithmetic-types}}{(\ref{eq:prime-product-law})}
\notationrow{$\|\cdot\|_{\mathrm{TV}}$}{total variation
  distance}{the whole paper}{\S\ref{sec:arithmetic-types}}
\notationrow{$k$, $\sigma$, $B_r$, $R_r$}{$\abs P$; sum of $1/p$ over
  $P$; the moment quantities (redefined in
  Sections~\ref{sec:joint-degree-log-bridge}--\ref{sec:sharp-moment-log-bridge})}%
  {from Section~\ref{sec:arithmetic-types}}%
  {Prop.~\ref{prop:moment-tv}}
\notationrow{$B_j$, $T(x)$, $q_j$}{prime blocks; block type; mass of
  block $j$}{from Section~\ref{sec:arithmetic-types}}%
  {Lem.~\ref{lem:block-types}, (\ref{eq:block-law})}
\notationrow{$\alpha$, $\gamma$, $\rho$}{empty mass; coupling
  parameters of the block law}{Sections~\ref{sec:compatible-measures},
  \ref{sec:part-integer-bridges}, \ref{sec:part-log-bridges}}%
  {(\ref{eq:empty-mass}), (\ref{eq:coupling-parameters})}
\notationrow{$\xi$, $\nu_\xi$, $\omega$}{removed proportion of empty-type
  even labels; depleted law; an upper bound for $\xi$}{from
  Section~\ref{sec:arithmetic-types}}{Lem.~\ref{lem:empty-depletion},
  (\ref{eq:coupling-depletion}), \S\ref{sec:compatible-measures}}
\notationrow{$z$, $a_I$}{cutoff prime; available even half-labels of
  type $I$}{from Section~\ref{sec:arithmetic-types}}%
  {Cor.~\ref{cor:residual-type-quotas}, (\ref{eq:exact-even-quota})}
\notationrow{$G=(A,B)$, $M$, $D$}{bipartite residual forest with even
  side $A$ and position side $B$; $M=\abs A$; a bound for
  $\Delta_G(A)$}{from Section~\ref{sec:forest-types}}%
  {Lem.~\ref{lem:bipartite-completion}}
\notationrow{$I(a)$, $J(b)$}{type of an even vertex; tag of a
  position}{from Section~\ref{sec:forest-types}}%
  {Lem.~\ref{lem:retained-tags}}
\notationrow{$H_{\mathrm{deg}}$}{degree error of an
  assignment}{from Section~\ref{sec:forest-types}}%
  {(\ref{eq:two-measure-degree-error})}
\notationrow{$Y$, $R$, $\Gamma_0$, $\Gamma$}{odd labels, $\abs Y=R$;
  type-compatibility graph, and the same after deletions}{from
  Section~\ref{sec:matching-estimates}}{\S\ref{sec:matching-estimates}}
\notationrow{$\mathcal T$, $0$, $F$}{type set $2^{[s]}$; empty and full
  types}{from Section~\ref{sec:refined-hall}}{\S\ref{sec:refined-hall}}
\notationrow{$\varepsilon$, $\tau$, $\eta$, $\gamma$, $u$}{exceptional
  fraction; degree targets; fraction of universal labels (written
  $\upsilon$ in Section~\ref{sec:part-intervals})}{from
  Section~\ref{sec:refined-hall}}{Lem.~\ref{lem:three-hall-ranges}}
\notationrow{$\mathcal B$, $p_{\rm loc}$}{capped rectangle minimum;
  local failure probability}{Sections~\ref{sec:refined-hall},
  \ref{sec:part-intervals}}{Thm.~\ref{thm:three-range-hall},
  Cor.~\ref{cor:three-range-probability}}
\end{notationlist}

\subsection*{Interval certificates}
\begin{notationlist}
\notationrow{$[\ell,u]$}{interval of orders}%
  {Section~\ref{sec:part-intervals}}{\S\ref{sec:local-probability}}
\notationrow{$m_n$, $w_n$, $w_-$}{$\floor{n/2}$; number of special
  even labels at order $n$; $\min_n w_n$}%
  {Section~\ref{sec:part-intervals}}{\S\ref{sec:finite-arithmetic}}
\notationrow{$b_*$, $W$, $\kappa$, $\rho$, $\upsilon$}{$b(u)$;
  $\min_n w_n/n$; $b_*/\ell$; $1/2-\kappa$; universal fraction}%
  {Section~\ref{sec:part-intervals}}{(\ref{eq:interval-preparation-constants})}
\notationrow{$D_*$}{$\ceil{1/(2W)}$, the largest residual degree}%
  {Section~\ref{sec:part-intervals}}{(\ref{eq:interval-degree-list})}
\notationrow{$e_+$, $H_+$, $\delta$, $q_-$, $\varepsilon_B$}{construction
  bounds of a row}{Section~\ref{sec:part-intervals}}%
  {(\ref{eq:interval-construction-interface})}
\notationrow{$\widehat b$, $p_{\rm glob}$}{transported Hall budget;
  global failure bound}{Sections~\ref{sec:exponential-continuation}--%
  \ref{sec:retained-continuation}}{(\ref{eq:transported-finite-budgets}),
  (\ref{eq:single-range-probability})}
\end{notationlist}

\subsection*{The unbounded range and the bridges}
\begin{notationlist}
\notationrow{$L$, $v$}{$\log n$; $\log z$}%
  {Sections~\ref{sec:part-analytic}, \ref{sec:part-log-bridges}}%
  {Lem.~\ref{lem:analytic-prime-inputs}}
\notationrow{$S_2(z)$}{sum of $p^{-2}$ over primes
  $p>z$}{Sections~\ref{sec:part-analytic}--\ref{sec:part-log-bridges}}%
  {Thm.~\ref{thm:large-prime-sorting}}
\notationrow{$w_-$, $g$}{lower bound for the special even-label
  fraction; number of colors}{Section~\ref{sec:part-integer-bridges};
  $w_-$ also in Section~\ref{sec:part-log-bridges}}%
  {\S\ref{sec:profile-block-bridge}}
\notationrow{$[a,J]$, $x$}{closed interval of $\log n$; normalized
  residual degree}{Section~\ref{sec:part-log-bridges}}%
  {\S\ref{sec:empty-separator-log-bridge}}
\notationrow{$a_i$, $N_i$}{endpoints of the logarithmic bridges; the
  integers with $\log N_i<a_i<\log(N_i+1)$}%
  {Section~\ref{sec:log-range}}{\S\ref{sec:log-range}}
\end{notationlist}

\subsection*{Letters with several meanings}

Several letters have different meanings in different parts of the
paper; in each part the meaning is fixed where it is introduced. The
most frequent cases are the following.
\begin{itemize}
\item $P$ is the set of small primes defining the types; in
  Section~\ref{sec:part-intervals} it also denotes a set of primes above
  the cutoff, such as $P(y)$, the primes above $z$ dividing $y$. The
  separator labels are $\Psep$.
\item $W$ is the reserve of positions in Section~\ref{sec:matching},
  the number $\min_n w_n/n$ in Section~\ref{sec:part-intervals}, and a
  function of $\log n$ in the logarithmic bridges.
\item $\rho$ is $\max_j q_j/(1-q_j)$ in
  Section~\ref{sec:compatible-measures}, $(D-1)/D$ in
  Sections~\ref{sec:forest-types} and~\ref{sec:paired-separators}, and
  $1/2-\kappa$ in Section~\ref{sec:part-intervals}.
\item $N$ is the size of a vertex set or label side, defined in each
  section; $N_i$ in Section~\ref{sec:log-range} are integer endpoints.
\item $w_-$ is the integer $\min_n w_n$ in
  Section~\ref{sec:part-intervals} and a lower bound for the fraction
  of special even labels in Sections~\ref{sec:part-integer-bridges}
  and~\ref{sec:part-log-bridges}.
\item $B_r$ and $R_r$ are the quantities of
  Proposition~\ref{prop:moment-tv}, except in
  Sections~\ref{sec:joint-degree-log-bridge},
  \ref{sec:reserve-mixture-log-bridge},
  and~\ref{sec:sharp-moment-log-bridge}, where they are redefined in
  Section~\ref{sec:joint-degree-log-bridge}.
\item $D$ bounds the residual degree, except in
  Proposition~\ref{prop:residual-problem}, where $D\subseteq Z$ is a
  set of positions; $F$ is the forest $T-K$ in
  Section~\ref{sec:matching}, the full type from
  Section~\ref{sec:refined-hall} on, and, in
  Sections~\ref{sec:simultaneous-full-type-bridge}
  and~\ref{sec:full-reserve-log-bridge}, the set of vertices that
  receive the special even labels; $Z_0$ is the residual position set
  in Section~\ref{sec:matching} and a set of absorbed positions in
  Section~\ref{sec:finite-rare-refinement}.
\end{itemize}

\section{Degree covers of the retained-position certificates}
\label{sec:cover-tables}

The following tables give the interval covers of the seven certificates
of Proposition~\ref{prop:verified-retained-continuations}. Each interval
carries one row for every degree $1,\ldots,D_*$, and the intervals of
each cover are adjacent. The columns are read downwards. The rows
themselves, with their constructions, prime blocks, and rational
parameters, are part of the supplement
(Section~\ref{sec:part-computation}).

\begin{table}[!htbp]
\centering
\caption{Degree cover on $3500\leq n\leq5219$.}
\label{tab:retained-continuation-cover}
\par\vspace{.5\baselineskip}
\renewcommand{\arraystretch}{1.1}
\begin{tabular}{rrrrrr}
\hline
Order interval&$D_*$&Order interval&$D_*$&Order interval&$D_*$\\
\hline
$3500$--$3503$&20&$3574$&20&$3920$--$4019$&22\\
$3504$--$3506$&20&$3575$--$3578$&20&$4020$--$4119$&21\\
$3507$--$3512$&20&$3579$--$3583$&20&$4120$--$4219$&21\\
$3513$--$3525$&20&$3584$--$3595$&20&$4220$--$4319$&20\\
$3526$&20&$3596$--$3607$&20&$4320$--$4419$&21\\
$3527$--$3529$&20&$3608$--$3619$&20&$4420$--$4519$&22\\
$3530$--$3533$&20&$3620$--$3644$&20&$4520$--$4619$&22\\
$3534$--$3537$&20&$3645$--$3669$&20&$4620$--$4719$&22\\
$3538$--$3542$&20&$3670$--$3694$&21&$4720$--$4819$&22\\
$3543$--$3548$&20&$3695$--$3719$&21&$4820$--$4919$&22\\
$3549$--$3554$&20&$3720$--$3769$&21&$4920$--$5019$&22\\
$3555$--$3561$&20&$3770$--$3819$&20&$5020$--$5119$&22\\
$3562$--$3569$&20&$3820$--$3869$&21&$5120$--$5219$&22\\
$3570$--$3573$&20&$3870$--$3919$&21&&\\
\hline
\end{tabular}
\end{table}

\begin{table}[!htbp]
\centering
\caption{Degree cover on $5220\leq n\leq6771$.}
\label{tab:unified-continuation-cover}
\par\vspace{.5\baselineskip}
\renewcommand{\arraystretch}{1.1}
\begin{tabular}{rrrr}
\hline
Order interval&$D_*$&Order interval&$D_*$\\
\hline
$5220$--$5319$&23&$6020$--$6119$&22\\
$5320$--$5419$&22&$6120$--$6219$&22\\
$5420$--$5519$&22&$6220$--$6319$&23\\
$5520$--$5619$&22&$6320$--$6419$&23\\
$5620$--$5719$&22&$6420$--$6519$&23\\
$5720$--$5819$&22&$6520$--$6619$&23\\
$5820$--$5919$&22&$6620$--$6719$&23\\
$5920$--$6019$&22&$6720$--$6771$&23\\
\hline
\end{tabular}
\end{table}

\begin{table}[!htbp]
\centering
\caption{Degree cover on $6772\leq n\leq12499$.}
\label{tab:extended-continuation-cover}
\par\vspace{.5\baselineskip}
\renewcommand{\arraystretch}{1.1}
\begin{tabular}{rrrrrr}
\hline
Order interval&$D_*$&Order interval&$D_*$&Order interval&$D_*$\\
\hline
$6772$--$7021$&23&$8772$--$9021$&24&$10772$--$11021$&25\\
$7022$--$7271$&23&$9022$--$9271$&24&$11022$--$11271$&25\\
$7272$--$7521$&23&$9272$--$9521$&24&$11272$--$11521$&25\\
$7522$--$7771$&23&$9522$--$9771$&24&$11522$--$11771$&25\\
$7772$--$8021$&23&$9772$--$10021$&24&$11772$--$12021$&25\\
$8022$--$8271$&24&$10022$--$10271$&24&$12022$--$12271$&25\\
$8272$--$8521$&24&$10272$--$10521$&24&$12272$--$12499$&25\\
$8522$--$8771$&24&$10522$--$10771$&25&&\\
\hline
\end{tabular}
\end{table}

\begin{table}[!htbp]
\centering
\caption{Degree covers of the four certificates on $12500\leq n\leq34999$.}
\label{tab:upper-retained-covers}
\par\vspace{.5\baselineskip}
\renewcommand{\arraystretch}{1.1}
\begin{tabular}{lrr}
\hline
Certificate&Order interval&$D_*$\\
\hline
$12500$--$14999$&$12500$--$12999$&25\\
&$13000$--$13499$&26\\
&$13500$--$13999$&26\\
&$14000$--$14499$&26\\
&$14500$--$14999$&26\\
\hline
$15000$--$19999$&$15000$--$15999$&26\\
&$16000$--$16999$&26\\
&$17000$--$17999$&27\\
&$18000$--$18999$&27\\
&$19000$--$19999$&27\\
\hline
$20000$--$24999$&$20000$--$22499$&27\\
&$22500$--$24999$&28\\
\hline
$25000$--$34999$&$25000$--$29999$&29\\
&$30000$--$34999$&29\\
\hline
\end{tabular}
\end{table}

\clearpage

\section*{Use of generative AI}
Generative AI systems, principally GPT Astra and Claude Opus, were used
extensively in the development of this work, including mathematical proof
development, computational verification, literature and research suggestions,
drafting, restructuring, and technical consistency review, under the author's
direction. The author reviewed the mathematical arguments and computational
results presented in the manuscript and takes full responsibility for the
accuracy and integrity of the work.

\section*{Acknowledgments}
The author acknowledges the extensive assistance of GPT Astra and Claude Opus
in the development and preparation of this manuscript. The AI systems are not
authors.

\section*{Code and data availability}

The proof of Theorem~\ref{thm:main} depends on the exact finite
computations described in Section~\ref{sec:part-computation}. The code,
certificate data, fixed inputs, and verification reports used for these
computations are archived in the frozen computational supplement at
\href{https://doi.org/10.5281/zenodo.23057222}{doi:10.5281/zenodo.23057222}.
The supplement includes instructions for reproducing the verification
and a wrapper, \texttt{run\_all\_checks.py}, which runs the documented
checks and reports their exit status; it is not a separate verifier.
Section~\ref{sec:part-computation} identifies which mathematical
statements depend on these computations and why the finite checks
suffice. Search programs used to discover certificate parameters are not
part of the proof unless explicitly identified otherwise in the
supplement.

\end{document}